\documentclass[reqno]{amsart}
\usepackage{amsthm,amsfonts,amssymb,euscript}
\usepackage{latexsym, multicol, fancybox}
\usepackage{graphicx}
\usepackage{color}
\usepackage{amsmath, amsthm, amssymb, bm}
\usepackage{epstopdf}
\usepackage{caption}
\usepackage{psfrag}

\usepackage{mathrsfs}
 \usepackage{xcolor}
 \usepackage[citebordercolor={green}]{hyperref}
 \usepackage{tikz}
 \usepackage{bbm}
 \usepackage{dsfont}
 \usepackage{bbold}
 
\numberwithin{equation}{section}

 \usepackage{tikz}

\newtheorem{theorem}{Theorem}[section]
\newtheorem{lemma}[theorem]{Lemma}
\newtheorem{proposition}[theorem]{Proposition}
\newtheorem{corollary}[theorem]{Corollary}
\newtheorem{definition}[theorem]{Definition}
\newtheorem{remark}[theorem]{Remark}
\newtheorem{conjecture}[theorem]{Conjecture}

\newcommand{\bea}{\begin{eqnarray}}
\newcommand{\eea}{\end{eqnarray}}
\def\beaa{\begin{eqnarray*}}
\def\eeaa{\end{eqnarray*}}
\def\ba{\begin{array}}
\def\ea{\end{array}}
\def\be#1{\begin{equation} \label{#1}}
\def \eeq{\end{equation}}

\newcommand{\bsub}{\begin{subequations}}
\newcommand{\esub}{\end{subequations}}

\newcommand{\nn}{\nonumber}

\newcommand{\ds}{\displaystyle}

\def\a{{\alpha}}

\def\b{{\beta}}
\def\be{{\beta}}
\def\ga{\gamma}
\def\Ga{\Gamma}
\def\de{\delta}
\def\De{\Delta}
\def\ep{\epsilon}

\def\la{\lambda}
\def\La{\Lambda}

\def\Si{\Sigma}
\def\om{\omega}
\def\Om{\Omega}

\def\vphi{\varphi}

\def\th{\theta}

\def\ze{\zeta}

\def\nab{\nabla}

\def\Up{\Upsilon}

\def\pr{{\partial}}
\def\les{\lesssim}
\def\c{\cdot}

\def\AA{{math\cal A}}

\def\MM{{\mathcal M}}
\def\NN{{\mathcal N}}

\def\LL{{\mathcal L}}
\def\II{{\mathcal I}}

\def\HH{{\mathcal H}}

\def\TT{{\mathcal T}}

\def\JJ{{\mathcal J}}

\def\Lie{{\mathcal L}}

\def\DD{{\mathcal D}}

\def\AA{{\mathcal A}}

\def\HH{{\mathcal H}}

\def\Lie{{\mathcal L}}

\def\lap{{\triangle}}

\def\B{{\bf B}}

\def\D{{\bf D}}
\def\E{{\bf E}}
\def\F{{\bf F}}

\def\M{{\bf M}}
\def\L{{\bf L}}
\def\O{{\bf O}}
\def\Q{{\bf Q}}
\def\R{{\bf R}}

\def\T{T}
\def\Z{Z}

\def\g{{\bf g}}

\def\RRR{{\Bbb R}}

\def\CCC{{\Bbb C}}
\def\f12{{\frac 1 2}}

\def\dual{{\,\,^*}}
\def\div{{\mbox div\,}}
\def\curl{{\mbox curl\,}}

\def\Hb{\,\underline{H}}

\def\Xh{\,^{(h)}X}

\def\trch{{\mbox tr}\chi}
\def\chih{{\widehat \chi}}
\def\chib{{\underline \chi}}
\def\chibh{{\underline{\chih}}}

\def\etab{{\underline \eta}}
\def\omb{{\underline{\om}}}
\def\bb{{\underline{\b}}}
\def\aa{\protect\underline{\a}}
\def\xib{{\underline \xi}}

\def\Xib{\underline{\Xi}}

\def\Ab{\protect\underline{A}}
\def\Bb{\protect\underline{B}}
\def\Xb{\protect\underline{X}}

\def\Xh{\widehat{X}}
\def\Xbh{\widehat{\Xb}}

\def\tr{\mbox{tr}}
\def\atr{\,^{(a)}\mbox{tr}}

\def\trchb{{\tr\chib}}

\def\Div{\mbox{Div}}

\def\atrch{\atr\chi}
\def\atrchb{\atr\chib}

\def\hot{\widehat{\otimes}}
\def\rhod{\,\dual\hspace{-2pt}\rho}

\def\fb{\protect\underline{f}}
\def\err{{\mbox{Err}}}
\def\ov{\overline}

\def\f12{\frac 1 2}
\def\lab{\label}
\def\nabc{\,^{(c)}\nab}

\def\bsplit{\begin{split}}

\newcommand{\Mext}{{\,{}^{(ext)}\mathcal{M}}}

\newcommand{\Mint}{{\protect \,{}^{(int)}\mathcal{M}}}

\def\Mint{{\protect \, ^{(int)}\MM}}

\def\qf{\mathfrak{q}}
\def\qfb{\protect\underline{\qf}}

\def\Nk{\mathfrak{N}}

\def\Rk{\mathfrak{R}}

\def\Ik{\mathfrak{I}}
\def\Jk{\mathfrak{J}}

\def\sk{\mathfrak{s}}
\def\dkb{ \, \mathfrak{d}     \mkern-9mu /}
\def\dk{\mathfrak{d}}

\DeclareFontFamily{U}{mathx}{\hyphenchar\font45}
\DeclareFontShape{U}{mathx}{m}{n}{
      <5> <6> <7> <8> <9> <10>
      <10.95> <12> <14.4> <17.28> <20.74> <24.88>
      mathx10
      }{}
\DeclareSymbolFont{mathx}{U}{mathx}{m}{n}
\DeclareFontSubstitution{U}{mathx}{m}{n}
\DeclareMathAccent{\widecheck}{0}{mathx}{"71}

\def\Zc{\widecheck{Z}}
\def\Hc{\widecheck{H}}
\def\Hbc{\widecheck{\Hb}}
\def\trXc{\widecheck{\tr X}}
\def\trXbc{\widecheck{\tr\Xb}}
\def\Pc{\widecheck{P}}
\def\ombc{\widecheck{\omb}}
\def\Gac{\widecheck{\Ga}}
\def\Rc{\widecheck R}

\def\rhoc{\widecheck{\rho}}

\def\omc{\widecheck \omega}
\def\ombc{\underline{\widecheck{\omega}}}

\def\trchc{\widecheck{\tr\chi}}
\def\trchbc{\widecheck{\tr\chib}}

\def\DDc{\,^{(c)} \DD}

\newcommand{\deh}{\delta_{\mathcal{H}}}
\newcommand{\dec}{\delta_{dec}}
\newcommand{\dee}{\delta_{extra}}
\newcommand{\dt}{\delta_B}

\def\That{{{ \widehat T}}}

\newcommand{\Lieb}{\Lie \mkern-10mu /\,}

\def\DDov{\ov{\DD}}

\def\DDs{ \, \DD \hspace{-2.4pt}\dual    \mkern-20mu /}
\def\DDd{ \, \DD \hspace{-2.4pt}    \mkern-8mu /}

\def\Bdot{\dot{B}}

\def\Rdot{\dot{\R}}

\def\DDc{\,^{(c)} \DD}

\def\DDb{\ov{\DD}}
\def\DDbc{\ov{\DDc}}

\def\Lied{\dot{\Lie}}

\def\DDov{\ov{\DD}}

\def\Ddot{\dot{\D}}
\def\squared{\dot{\square}}

 \def\Ft{\widetilde{F}}

\def\eS{\,^{S} \hspace{-1.5pt} e}

\def\Ddot{\dot{\D}}

\def\Ddot{\dot{\D}}

\def\Kh{\,^{(h)}K}

\def\Rhat{{\widehat{R}}}

\def\nabc{\,^{(c)}\nab}

\DeclareFontFamily{U}{mathx}{\hyphenchar\font45}
\DeclareFontShape{U}{mathx}{m}{n}{
      <5> <6> <7> <8> <9> <10>
      <10.95> <12> <14.4> <17.28> <20.74> <24.88>
      mathx10
      }{}
\DeclareSymbolFont{mathx}{U}{mathx}{m}{n}
\DeclareFontSubstitution{U}{mathx}{m}{n}
\DeclareMathAccent{\widecheck}{0}{mathx}{"71}

\def\trXc{\widecheck{\tr X}}
\def\trXbc{\widecheck{\tr\Xb}}
\def\Hc{\widecheck{H}}
\def\Zc{\widecheck{Z}}
\def\Pc{\widecheck{P}}
\def\ombc{\widecheck{\omb}}

      \def\ntrap{trap\mkern-18 mu\big/\,}
          \def\Mtrap{\,\MM_{trap}}
\def\Mntrap{{\MM_{\ntrap}}}

    \def\DDc{\,^{(c)} \DD}

\def\Bdot{\dot{B}}

\def\Pdot{\dot{P}}
\def\Bbdot{\dot{\Bb}}

      \def\ntrap{trap\mkern-18 mu\big/\,}
\def\Mntrap{{\MM_{\ntrap}}}

\def\Sk{\mathfrak G}

\def\Rkext{\,^{(ext)} \Rk}
\def\Skext{\,^{(ext)} \Sk}
\def\Skint{\,^{(int)} \Sk}
\def\Rkint{\,^{(int)} \Rk}

\def\Pdot{\dot{P}}
\def\Bdot{\dot{B}}
\def\Bbdot{\dot{\Bb}}

\def\Bt{\widetilde{B}}

\def \Bbt{\widetilde{\Bb}}
\def\At{\widetilde{A}}
\def \Abt{\widetilde{\Ab}}

     \def\kl{{k_L}}
    
    \def\DDc{\,^{(c)} \DD}

\def\Psib{\und{\Psi}}

\def\Pdot{\dot{P}}
\def\Bdot{\dot{B}}
\def\Bbdot{\dot{\Bb}}

\def\und{\underline}

\def\N{\mathbf{N}}

\newcommand{\phis}[1]{\pmb\phi_{s}^{(#1)}}

\def\gam{\g_{a,m}}
\def\dhor{\delta_{\HH}}
\def\dbl{\delta_{\textbf{BL}}}
\def\tmod{t_{\text{mod}}}
\def\phimod{\phi_{\text{mod}}}
\def\tt{\tau}
\def\ttt{\widetilde{\tau}}
\def\tphi{\varphi}

\def\Xcal{\mathcal{X}}

\def\dred{\delta_{\text{red}}}

\def\Ytau{Y_\tau}

\def\tx{\widetilde{x}}
\def\tJk{\widetilde{\Jk}}
\def\tq{\widetilde{q}}
\def\rr{\widetilde{r}}

\def\EMF{{\bf EMF}}

\def\EF{{\bf EF}}

\def\BEF{{\bf BEF}}
\def\EF{{\bf EF}}
\def\BF{{\bf BF}}
\def\BE{{\bf BE}}

\newcommand{\Nmic}{N_0}

\def\qs{|q|^2}

\def\reg{\mathbf{k}}

\def\kst{k_*}

\begin{document}

\title{Kerr stability in the full subextremal range}
\author{J\'{e}r\'{e}mie Szeftel}

\begin{abstract}
In the present  paper, we extend the results in the sequence of works \cite{KS-GCM1} \cite{KS-GCM2} \cite{KS:Kerr} by Sergiu Klainerman and the author, \cite{GKS22} by Elena  Giorgi, Sergiu Klainerman and the author, and \cite{Shen} by Dawei Shen from the slowly rotating regime $|a|\ll m$ to the full subextremal range $|a|<m$, thereby completing the proof of the Kerr stability conjecture. The restriction $|a|\ll m$ enters in fact only through the estimates for wave equations and for the Bianchi system established in \cite{GKS22}. To remove this restriction, we crucially rely on the two companion papers \cite{MaSz24} \cite{MaSz26} by Siyuan Ma and the author, in which we prove energy-Morawetz estimates respectively for the inhomogeneous scalar wave equation and for inhomogeneous Teukolsky equations on perturbations of Kerr with $|a|<m$. Additionally, some of the results of the present paper are proved in our companion paper \cite{Sze}.
\end{abstract}

\maketitle

\tableofcontents

%%%%%%%%%%%%%%%%%%%%

\section{Introduction}

%%%%%%%%%%%%%%%%%%%%%

%%%%%%%%%%%%%%%%%%%%%%%%%%%%%%%%%%%%%%%%

\subsection{Kerr stability conjecture and rough statement of the main result}
  
%%%%%%%%%%%%%%%%%%%%%%%%%%%%%%%%%%%%%%%% 

%%%%%%%%%%%%%%%%%%%%%%%%%%%

\subsubsection{Einstein vacuum equations}

%%%%%%%%%%%%%%%%%%%%%%%%%%%

The Einstein vacuum equations (EVE) in a {Lorentzian manifold} $(\MM, \g)$ take the form
\bea\lab{eq:EVE:intro}
\mathbf{R}_{\a\b}=0,
\eea
where $\mathbf{R}_{\a\b}$ {denotes} the Ricci curvature tensor of the metric $\g$. Foundational contributions by Choquet-Bruhat \cite{CB52}, and by  
Choquet-Bruhat and Geroch \cite{CBG69}, formulate EVE as an evolution problem of hyperbolic type and associate to any suitable sufficiently regular initial data set a unique (up to diffeomorphisms) maximal Cauchy development.

%%%%%%%%%%%%%%%%%%%%%%%%%%%

\subsubsection{Kerr solution}

%%%%%%%%%%%%%%%%%%%%%%%%%%%

The EVE admit a family of explicit solutions, found by Kerr \cite{Kerr63} in 1963, which describe asymptotically flat, stationary, axially symmetric black hole spacetimes. The metrics of Kerr spacetimes are parameterized by an angular momentum per unit mass $a$ and a mass $m$,  satisfying $|a|\leq m$, and take the following form in the Boyer--Lindquist coordinates $(t,r,\th, \phi)$
\bea\lab{eq:expressionofKerrmetricinBLcoordinates:intro}
\gam=-\frac{\Delta \qs}{\Sigma^2} dt^2 + \frac{\sin^2\th\Sigma^2 }{\qs}\bigg(d\phi - \frac{2amr}{\Sigma^2} dt\bigg)^2 +\frac{\qs}{\Delta} dr^2 + \qs d\th^2,
\eea
where
\bea
\Delta = r^2 - 2mr +a^2, \quad \qs=r^2+a^2\cos^2\th, \quad \Sigma^2=(r^2+a^2)^2 - a^2\sin^2\th \Delta.
\eea
Note that the particular case $a=0$ with $m>0$ corresponds to the family of Schwarzschild spacetimes, introduced by Schwarzschild \cite{Sch16} in 1916.

We consider in this work the family of \textit{subextremal} Kerr spacetimes, in which the two parameters $(a,m)$ satisfy the strict inequality $|a|<m$. Such a subextremal Kerr spacetime contains a black hole $\{r<r_+\}$ with a nondegenerate event horizon located at $\{r=r_+\}$ where $r_+:=m+\sqrt{m^2-a^2}$ is the larger root of $\Delta=\De(r)$, see Figure \ref{fig:penrosediagramofKerr} for the corresponding Penrose diagram.

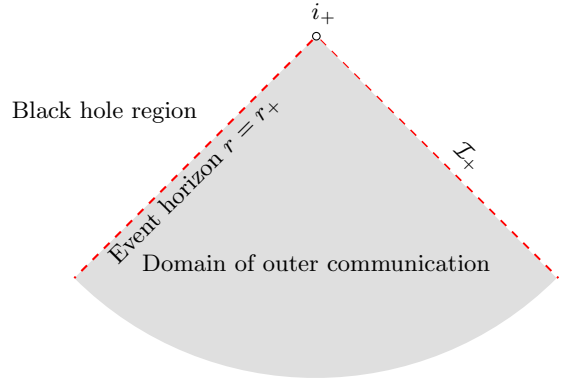
\begin{figure}[htbp]
  \begin{center}
\begin{tikzpicture}[scale=1]
\tikzstyle{every node}=[font=\small]
      \draw[dashed, color=red, thick] (0.05,3.95) -- (3.2,0.8);
  \fill[lightgray!50] (-0.05,3.95)--(-3.2,0.8) arc(225:315: 4.51 and 4.51) -- (0.05,3.95);
   \draw[dashed, color=red, thick] (-0.05,3.95)--(-3.2,0.8) ;
     \node at (-2.8,3) {Black hole region};
     \node at (0.1,4.3) {$i_+$};
       \node[rotate=315]  at (2.0,2.4) {$\II_+$};
       \node[rotate=45] at (-1.6,2.1) {Event horizon $r=r_+$};
         \draw[] (0,4) circle (0.05);
         \node at (0,1) {Domain of outer communication};  
\end{tikzpicture}
\end{center}
\caption{\footnotesize{Penrose diagram of subextremal Kerr spacetimes.}}
\lab{fig:penrosediagramofKerr}
\end{figure}

%%%%%%%%%%%%%%%%%%%%%%%%%%%

\subsubsection{Kerr stability conjecture}

%%%%%%%%%%%%%%%%%%%%%%%%%%%

The \textit{black hole stability conjecture} is one of the central open problems in general relativity. We provide a rough statement below.

\begin{conjecture}[Kerr stability conjecture]
The maximal Cauchy development of any initial data set for EVE, that is sufficiently close to a subextremal Kerr initial data set in a suitable sense, has a complete future null infinity and a domain of outer communication\footnote{The domain of outer communication is the complement of the black hole region.} which is asymptotic to a nearby member of the subextremal Kerr family. 
\end{conjecture}

%%%%%%%%%%%%%%%%%%%%%%%%%%
  
\subsubsection{First version of the main result}

%%%%%%%%%%%%%%%%%%%%%%%%%%
 
The following theorem proves the Kerr stability conjecture.
   \begin{theorem}[Main theorem, rough version]
\lab{MainThm-firstversion}
Let $a_0$ and $m_0$ be real constants such that $|a_0|<m_0$. The future globally hyperbolic   development  of  a general,   asymptotically  flat, initial data set, sufficiently close to a   $Kerr(a_0, m_0) $   initial data set in the sense that 
\bea\lab{eq:rmk:remarkonroughtopologyinitialdatamainTh}
\g=\g_{a_0,m_0}+O(\ep_0r^{-\frac{3}{2}-\de})\quad\textrm{as}\quad r\to +\infty\quad\textrm{along the spacelike initial hypersurface}\quad\Si_0,
\eea
for $\de>0$ and for $\ep_0>0$ small enough, has a complete    future null infinity  $\II^+$ and converges in  its causal past  $\JJ^{-1}(\II^{+})$  to another  nearby Kerr spacetime $Kerr(a_f, m_f)$ with parameters    $(a_f, m_f)$ close to the initial ones $(a_0, m_0)$.
 \end{theorem}

 \begin{figure}[ht!]
 \lab{fig0-introd}
\centering
\includegraphics[scale=0.35]{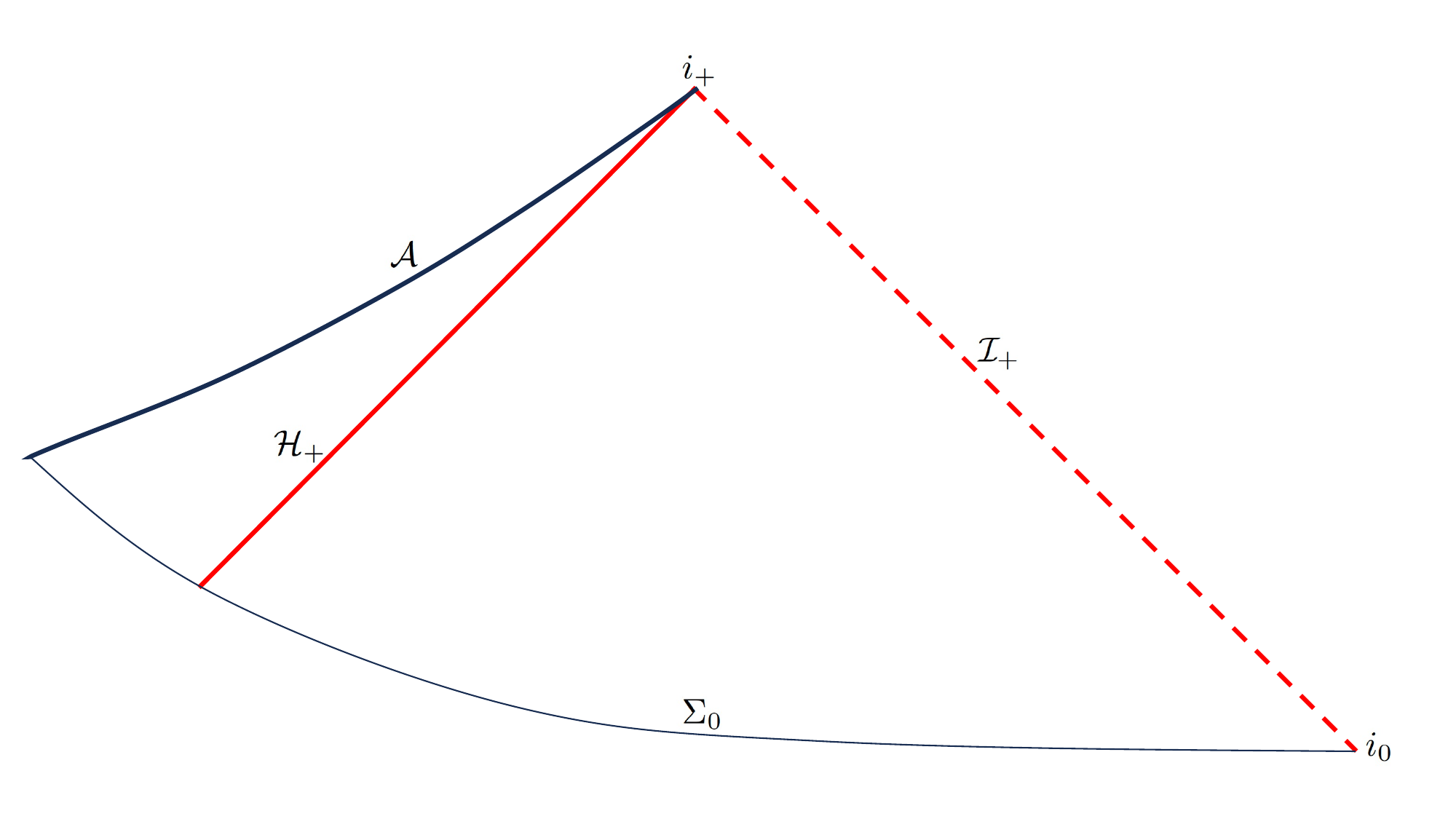}
\caption{\footnotesize{The Penrose diagram of the  final   space-time in  Theorem \ref{MainThm-firstversion} with initial spacelike hypersurface $\Si_0$,  future space-like boundary $\AA$, and $\II^+$ the complete future null infinity.  The hypersurface $\HH_+$} denotes the future event horizon.}
\end{figure} 

Prior to this work, the Kerr stability conjecture has been addressed in the case of Schwarzschild (i.e., for $a=0$) in polarized symmetry by Klainerman-Szeftel \cite{KS}, in the case of Schwarzschild for a codimension 3 subset of initial data by Dafermos-Holzegel-Rodnianski-Taylor \cite{DHRT}, and in the case of Kerr for $|a|\ll m$ in the sequence of works \cite{KS-GCM1} \cite{KS-GCM2} \cite{KS:Kerr} by Klainerman and the author, \cite{GKS22} by Giorgi, Klainerman and the author, and \cite{Shen} by Shen, where the initial data in \cite{KS} and \cite{GKS22} \cite{KS-GCM1} \cite{KS-GCM2} \cite{KS:Kerr} \cite{Shen} are as in \eqref{eq:rmk:remarkonroughtopologyinitialdatamainTh}, while the ones in \cite{DHRT} fall off like $O(\ep_0r^{-\frac{5}{2}})$ as $r\to +\infty$ along $\Si_0$. Recently, it has also been proved for all $|a|<m$ by Hintz \cite{Hi2} \cite{Hi3} \cite{Hi1} for a class of well-prepared initial data satisfying
\bea\lab{eq:intro:wellpreparedinitialdatainresultbyPeterHintz:0}
\g-\g_{a_0,m_0} = \textrm{finite polyhomogeneous expansion} + O(\ep_0r^{-3-\de})\quad\textrm{as}\quad r\to +\infty,
\eea
for $\de>0$ and for $\ep_0>0$ small enough, which allows for a reduction to finite-dimensional gauge modifications\footnote{In general, infinite-dimensional gauge modifications are needed, as expected from the infinite-dimensionality of the kernel of the linearized operator generated by general covariance, see Section \ref{sec:comparisionwithresultbyHintz} and Remark 1.3 in \cite{Hi1} for more details.}. We refer to Section \ref{sec:previousresultsonKerrstabilityconjecture} for a detailed review of the literature, and in particular to Section \ref{sec:comparisionwithresultbyHintz} for a comparison of \cite{Hi1} with our main result.

%%%%%%%%%%%%%%%%%%%%%%%%%%%%%%%%%%%%
 
\subsection{Strategies to approach the Kerr stability conjecture}
\lab{sec:intro:strategytoappraochKerrstability}

%%%%%%%%%%%%%%%%%%%%%%%%%%%%%%%%%%%%

To put the review on advances related to the Kerr stability conjecture, presented in Section \ref{sec:previousresultsonKerrstabilityconjecture}, into context, we begin with a broad overview of the two approaches that have been developed to tackle this problem.

%%%%%%%%%%%%%%%%%%%%%%%%%%%%%%%%%%%%%%%%%%%
  
\subsubsection{Main difficulties}
\lab{sec:intro:maindifficultiesoftheproofifKerrstab}

%%%%%%%%%%%%%%%%%%%%%%%%%%%%%%%%%%%%%%%%%%%

We start with a short summary of the proof of the stability of Minkowski before outlining the main new difficulties arising in the context of the Kerr stability conjecture.

%%%%%%%%%%%%%%%%%%%%%%%%%%%%%%%%%%%%%%%%%%%
  
\paragraph{\textit{The stability of Minkowski}.}

%%%%%%%%%%%%%%%%%%%%%%%%%%%%%%%%%%%%%%%%%%%

The stability of Minkowski has been proved by D. Christodoulou and S. Klainerman in \cite{Ch-Kl}, see \cite{SmuMinkreview} and \cite{Shen2} for a review of the literature in the wake of this groundbreaking result. The main  innovations introduced in \cite{Ch-Kl} are as follows:
\begin{enumerate}
\item A robust method to derive quantitive decay for linear fields on perturbations of Minkowski based on the energy method, using generators of approximate symmetries of Minkowski as multipliers and commutators, in conjunction with Klainerman-Sobolev inequalities. 

\item The construction of geometric integrable\footnote{This refers to the fact that the orthogonal space $\{e_3,e_4\}^\perp$ to the null pair $(e_3, e_4)$ is integrable in the sense of Frobenius.} well-posed gauges connected to the distinguished integrable null pair $(e_4=\pr_t+\pr_r, e_3=\pr_t-\pr_r)$ of Minkowski.

\item The choice of gauge is used to exhibit a special structure of the nonlinear terms of \eqref{eq:EVE:intro}, known as the null structure, which then needs to be exploited. 
\end{enumerate}

%%%%%%%%%%%%%%%%%%%%%%%%%%%%%%%%%%%%%%%%%%%%%%
 
\paragraph{\textit{Main new difficulties in connection to the Kerr stability conjecture}.}

%%%%%%%%%%%%%%%%%%%%%%%%%%%%%%%%%%%%%%%%%%%%%%

Compared to the stability of Minkowski, the following main new difficulties need to be addressed, see the reviews in \cite{Kl-review} and \cite{Ch1} for more context:
\begin{enumerate}
\item \textit{Trapping and superradiance.} Kerr black holes possess much less symmetries than Minkowski, and their more involved geometry displays a priori obstructions to derive quantitive decay for linear fields on Kerr in the form of trapping and superradiance.

\item \textit{Large coupling of the Linearized Gravity System.} The linearized EVE around Kerr, known as the Linearized Gravity System (LGS), form a large coupled system. It can be derived for metric perturbations or curvature perturbations, see Sections \ref{sec:intro:strategycurvarureperturbation} and \ref{sec:intro:strategymetricpertubations}.

\item \textit{Modulation in infinite dimension.} LGS has a large kernel which translates at the nonlinear level in the need to track dynamically the change in the parameters $a$ and $m$, as well as the change in the gauge during evolution. Such a procedure is known as modulation, and tracking the change in the gauge requires to do modulation in infinite dimension.  

\item \textit{Non-integrability of the principal null directions of Kerr.} Kerr possesses a natural pair of null directions, known as principal null directions, with the principal null pair regular across the future event horizon given in Boyer-Lindquist coordinates by
\bea
\lab{def:e3e4inKerr:intro}
 e_4 = \frac{r^2+a^2}{|q|^2} \pr_t +\frac{\De}{|q|^2} \pr_r +\frac{a}{|q|^2} \pr_{\phi}, \qquad 
 e_3=\frac{r^2+a^2}{\De} \pr_t -\pr_r +\frac{a}{\De} \pr_{\phi}.
\eea
This null pair diagonalizes the curvature tensor and thus plays a distinguished role in dealing with the large coupling of LGS. Unlike their analog in Minkowski, these null directions are non-integrable whenever $a\neq 0$, requiring to extend the construction of geometric well-posed gauges in the context of the stability of Minkowski from the integrable setting to the non-integrable setting.
\end{enumerate} 

To address the difficulty related to the coupling of LGS, two strategies have emerged, based respectively on formulating the problem w.r.t. curvature perturbations or metric perturbations. We review them below, starting with the one used in the proof of Theorem \ref{MainThm-firstversion}.

%%%%%%%%%%%%%%%%%%%%%%%%%%%%%%%%%%%%%%%%%%%
  
\subsubsection{Strategy based on curvature perturbations}
\lab{sec:intro:strategycurvarureperturbation}

%%%%%%%%%%%%%%%%%%%%%%%%%%%%%%%%%%%%%%%%%%%

%%%%%%%%%%%%%%%%%%%%%%%%%%%%%%%%%%%%%%%%%%%
  
\paragraph{\textit{LGS in curvature perturbations and Teukolsky scalars}.}

%%%%%%%%%%%%%%%%%%%%%%%%%%%%%%%%%%%%%%%%%%%

Curvature perturbations refer to a formulation of \eqref{eq:EVE:intro} leading to a coupled system involving curvature components and Ricci coefficients computed w.r.t. a null frame (or null tetrad in the context of Newman-Penrose formalism \cite{NP}), where: 
\begin{itemize}
\item the curvature components solve the Bianchi system (resulting from a combination of \eqref{eq:EVE:intro} and the Bianchi identities), which is a first order hyperbolic system,

\item the Ricci coefficients solve the null structure equations, which, for a suitable choice of gauge, correspond to a coupled system of transport equations and elliptic equations of Hodge type. 
\end{itemize}
As mentioned in Section \ref{sec:intro:maindifficultiesoftheproofifKerrstab}, at the level of LGS, this results in a highly coupled system of tensorial equations. The main breakthrough to disentangle this system is due to Teukolsky in \cite{Teuk}, who identified the so-called \textit{Teukolsky scalars}, which are specific contractions of the curvature tensor w.r.t. the principal null frame satisfying, at the linearized level, decoupled wave-type equations  known as \textit{Teukolsky equations}.

%%%%%%%%%%%%%%%%%%%%%%%%%%%%%%%%%%%%%%%%%%%%%%
  
\paragraph{\textit{Strategy based on curvature perturbations and the energy method}.}

%%%%%%%%%%%%%%%%%%%%%%%%%%%%%%%%%%%%%%%%%%%%%%

The strategy based on curvature perturbations and the energy method builds first on the derivation of decay estimates for the simplest model problem, namely the scalar wave equation on Kerr. Such decay estimates for solutions $\psi$ to the scalar wave equation $\square_{\gam}\psi=0$ are derived based on the following strategy, see Section \ref{sec:stateoftheartformodelproblems} for more references:
\begin{enumerate}
\item The first and most fundamental step,  which appeared first in Blue-Soffer \cite{BS03}, is  the derivation of  an \textit{energy-Morawetz} estimate which deals with issues due to trapping\footnote{For $a\neq 0$, energy-Morawetz estimates need also to deal with superradiance which was first addressed in \cite{DaRo2} \cite{Ta-Toh} \cite{A-B}.} and provides a weak, albeit unconditional, integrated decay for $\psi$.

\item Then, using the \textit{redshift} vectorfield introduced in Dafermos-Rodnianski \cite{DaRo1} allows to remove the degeneracy near the horizon.

\item Finally, using in addition \textit{$r^p$-weighted} estimates\footnote{These estimates go back to Morawetz \cite{Mor1} \cite{Mor2} in the context of the exterior problem for the scalar wave equation in Minkowski  and were first used in the context of the scalar wave equation on Schwarzschild in \cite{Da-Ro3}.} in the far $r$ region in conjunction with a simple mean value argument, as introduced in Dafermos-Rodnianski \cite{Da-Ro3}, allows to obtain pointwise decay estimates for $\psi$.
\end{enumerate}

The strategy based on curvature perturbations and the energy method then proceeds as follows, see Section \ref{sec:previousresultsonKerrstabilityconjecture} for more references:
\begin{enumerate}
\item \textit{Decay estimates for Teukolsky scalars.} The strategy for the scalar wave equation outlined above cannot immediately extend to solutions to Teukolsky equations as they do not derive from a Lagrangian and are hence not directly amenable to the energy method. Instead, one applies a transformation, introduced in the work of Chandrasekhar \cite{Chand2} and first used in the context of decay estimates for Teukolsky scalars in Dafermos-Holzegel-Rodnianski \cite{D-H-R}, which yields a coupled tensorial wave-transport system for which the above strategy for the scalar wave can then be adapted, ultimately leading to quantitative decay estimates for Teukolsky scalars.

\item \textit{Decay estimates and triangular structure of the main equations.} Starting from the decay of Teukolsky scalars, one then recovers decay for the rest of the curvature components and the Ricci coefficients by exhibiting a triangular structure in the main equations. 

\item \textit{Modulation.} While Teukolsky scalars depend only quadratically on change of gauges, the rest of the system depends linearly on such changes. In particular, the previous step cannot be achieved without first designing a modulation procedure to deal with the instabilities coming from the kernel of LGS. In practice, the construction of preferred co-dimension 2 spheres, the General Covariant Modulated (GCM) spheres in Klainerman-Szeftel \cite{KS-GCM1} \cite{KS-GCM2}, tracks the change of gauge (the infinite dimensional part of the modulation scheme), and a suitable definition in \cite{KS-GCM2} of the mass and angular momentum on such GCM spheres deals with the remaining finite dimensional part of the modulation scheme.

\item \textit{Higher-order energy estimates.} In addition to decay estimates, energy-type estimates are derived for all components up to top order without losing derivatives. This requires in particular to rely on the hyperbolic character of the Bianchi system, see for instance Part III in \cite{GKS22} by Giorgi-Klainerman-Szeftel.

\item\textit{Non-integrable formalism and non-integrable gauges.}  As mentioned in Section \ref{sec:intro:maindifficultiesoftheproofifKerrstab}, the principal null pair of Kerr, which lies at the core of the definition of Teukolsky scalars and of the triangular structure of LGS, is non-integrable. This requires to extend the formalism and the well-posed geometric gauges in \cite{Ch-Kl} to the non-integrable setting as done in Part I of  \cite{GKS22} and Klainerman-Szeftel \cite{KS:Kerr} respectively.

\item \textit{Null structure.} Finally, as in the proof of \cite{Ch-Kl}, the null structure of the nonlinear terms plays an essential role, particularly in the derivation of decay estimates for the Teukolsky scalars. 
\end{enumerate}

%%%%%%%%%%%%%%%%%%%%%%%%%%%%%%%%%%%%%%%%%%%
  
\subsubsection{Strategy based on metric perturbations}
\lab{sec:intro:strategymetricpertubations}

%%%%%%%%%%%%%%%%%%%%%%%%%%%%%%%%%%%%%%%%%%%

Metric perturbations refers to a formulation of \eqref{eq:EVE:intro} in generalized wave coordinates which results in a coupled system of quasilinear wave equations for the metric coefficients. One can try to use the energy method together with a decoupling of LGS in metric perturbations, in the spirit of the approach outlined in Section \ref{sec:intro:strategycurvarureperturbation}. So far,  this only works for LGS in metric perturbations around Schwarzschild, see Hung-Keller-Wang \cite{HKW} and Johnson \cite{Johnson}. We present below another approach which relies instead on spectral methods.

%%%%%%%%%%%%%%%%%%%%%%%%%%%%%%%%%%%%%%%%%%%
  
\paragraph{\textit{Strategy based on metric perturbations and spectral methods}.}

%%%%%%%%%%%%%%%%%%%%%%%%%%%%%%%%%%%%%%%%%%%

Metric perturbations, in conjunction with spectral methods, have been first used in the seminal work of Hintz-Vasy in \cite{HVas} on the nonlinear stability of Kerr-de Sitter black holes. Roughly speaking, it follows the following steps (see the review in \cite{HafnerReview} for more details concerning the application of spectral methods in general relativity):
\begin{enumerate}
\item First, moving to the RHS all nonlinear terms in \eqref{eq:EVE:intro} written in generalized wave coordinates, including the quasilinear ones, reduces the problem to a coupled system of wave equations on Kerr with nonlinear forcing terms on the RHS.

\item To estimate the resulting coupled system of wave equations on Kerr, one then takes Laplace-Fourier transform in time which transforms it into a spectral problem for a second order system. 

\item The goal is then to derive estimates for the corresponding resolvent, the hardest part being the analysis of the regularity of the resolvant near the 0 time frequency.

\item  Next, going back to physical space variables using Plancherel allows to deduce quantitative decay estimates up to non-decaying modes.

\item A suitable mode stability result\footnote{The fact that mode stability holds for LGS on Kerr in metric perturbations is due to Andersson-H\"afner-Whiting \cite{AHW}.} then allows to show that all non-decaying modes can in fact be re-interpreted as being due to the dynamical change in mass, angular momentum and gauge so that they can be incorporated in a modulation procedure.

\item Also, as the formulation described in the first step does not preserve the constraints, it is necessary to use a constraint damping procedure ensuring that no artificial non-decaying modes are generated.  

\item Finally, the control of nonlinear terms uses the following two ingredients:
\begin{enumerate}
\item the weak null structure of EVE in harmonic gauge, exhibited in the new proof of the stability of Minkowski by Lindblad-Rodnianski \cite{LiRo1} \cite{LiRo2}, is exploited, 

\item the fact that quasilinear terms were moved to the RHS in the first step above induces a loss of derivatives which is overcome using a Nash-Moser iteration scheme.  
\end{enumerate}
\end{enumerate}

%%%%%%%%%%%%%%%%%%%%%%%%%%%%%%%%%%%%%%%%%%%%%%%%%%%%%%%%
 
\paragraph{\textit{Comparison with the strategy based on curvature perturbations and the energy method}.}

%%%%%%%%%%%%%%%%%%%%%%%%%%%%%%%%%%%%%%%%%%%%%%%%%%%%%%%%

We now compare the strategies outlined in Section \ref{sec:intro:strategycurvarureperturbation} and in the present section:
\begin{itemize}
\item \textit{Advantages of the spectral approach vs the energy method.} The main advantage of the spectral approach over the energy method is that it requires few properties concerning the structure of LGS in metric perturbations. By contrast, the energy method is highly dependent on the disentanglement of LGS in curvature perturbations using Teukolsky equations\footnote{Note however that Teukolsky equations play a central role in the mode stability result \cite{AHW} used in the spectral approach.} and on the transformation of Teukolsky into a wave-transport system amenable to the energy method.

\item \textit{Advantages of the energy method vs the spectral approach.} An advantage of the energy method is that it provides direct access to physically relevant quantities, such as Bondi's mass-loss formula, angular momentum and center of mass, which are more naturally connected to curvature perturbations. Also, unlike the energy method, the spectral approach necessarily uses a Nash-Moser iteration scheme to overcome the fact that it relies on a formulation which loses derivatives. More importantly, contrary to the energy method, the spectral approach does not allow so far to deal with infinite dimensional gauge modifications\footnote{Such modifications are in general needed, as expected from the infinite-dimensionality of the kernel of LGS due to general covariance, see the discussion in Section \ref{sec:intro:maindifficultiesoftheproofifKerrstab}.} which in turn restricts its application to well-prepared initial data, see also the discussion in Section \ref{sec:comparisionwithresultbyHintz}.
\end{itemize}

%%%%%%%%%%%%%%%%%%%%%%%%%%%%%%%%%%%%
 
\subsection{Previous results concerning the Kerr stability conjecture}
\lab{sec:previousresultsonKerrstabilityconjecture}

%%%%%%%%%%%%%%%%%%%%%%%%%%%%%%%%%%%%

In the wake of the initial contributions by Regge-Wheeler \cite{RW57}, Zerilli \cite{Ze}, Vishveshwara \cite{Vishev}, Teukolsky \cite{Teuk}, Chandrasekhar \cite{Chand2}, Wald \cite{Wald} and Whiting \cite{Whit}, a substantial literature has been motivated by the Kerr stability conjecture. We first review results concerning model problems before turning to results on the proof of the Kerr stability conjecture in particular cases.

%%%%%%%%%%%%%%%%%%%%%%%%%%%%%%%
 
\subsubsection{State of the art concerning model problems}
\lab{sec:stateoftheartformodelproblems}

%%%%%%%%%%%%%%%%%%%%%%%%%%%%%%%

We have outlined the general strategy allowing to derive pointwise decay for the scalar wave equation and Teukolsky equations in Section \ref{sec:intro:strategycurvarureperturbation} based in particular on a combination of energy-Morawetz estimates, redshift estimates and $r^p$-weighted estimates. We start with a short review of the derivation of energy-Morawetz estimates for these two models and refer to the introductions in  \cite{MaSz24} \cite{MaSz26} for more context.

%%%%%%%%%%%%%%%%%%%%%%%%%%%%%%%%%%%%%%%%%%%%%
 
\paragraph{\textit{Energy-Morawetz estimates for the scalar wave equation on Kerr}.}

%%%%%%%%%%%%%%%%%%%%%%%%%%%%%%%%%%%%%%%%%%%%%

In Schwarzschild, i.e., $a=0$, a Morawetz estimate, together with a degenerate conserved energy estimate based on the causal Killing vectorfield $\pr_t$, is first derived by Blue-Soffer in \cite{BS03} using a radial vectorfield linearly degenerating at the trapping radius $r=3m$. It is used in \cite{BS03} and later in \cite{BS07,BS06} to study the decay of semilinear wave equations in Schwarzschild spacetime. Dafermos-Rodnianski \cite{DaRo1} later introduced so-called redshift estimates which exploit the positivity of the surface gravity on the event horizon to remove the degeneracy near the event horizon. Subsequently, Marzuola-Metcalfe-Tataru-Tohaneanu \cite{MMTT} proved Strichartz estimates.

In the slowly rotating case, i.e., $|a|\ll m$, Dafermos-Rodnianski \cite{DaRo2} obtained a uniform boundedness result using a frequency decomposition into modes. Soon after, the first energy-Morawetz estimates, as well as weak decay estimates, for scalar fields in slowly rotating Kerr were proved by Tataru-Tohaneanu \cite{Ta-Toh} using microlocal multipliers. Andersson-Blue \cite{A-B} then obtained a new proof using a purely physical space approach based on Carter's Killing $2$-tensor. 

The proof of energy-Morawetz estimates in the full subextremal range has been obtained by Dafermos-Rodnianski-Shlapentokh-Rothman \cite{mDiRySR2014} using mode decomposition based on the full separability of the wave equation in Kerr. Recently, He-Klainerman \cite{HeKl1} \cite{HeKl2} have recovered energy-Morawetz estimates in the range $|a|\leq 0.75m$ relying almost entirely on physical space methods\footnote{Section 1.7 in \cite{HeKl1} contains the claim that the authors can push their methods all the way to $|a|<m$.}.

%%%%%%%%%%%%%%%%%%%%%%%%%%%%%%%%%%%%%%%%%
 
\paragraph{\textit{Energy-Morawetz estimates for Teukolsky equations on Kerr}.}

%%%%%%%%%%%%%%%%%%%%%%%%%%%%%%%%%%%%%%%%%

In order to derive energy-Morawetz estimates, one must first address the question of mode stability for solutions to Teukolsky equations in Kerr spacetimes. The absence of exponentially growing mode solutions was proved in the seminal work of Whiting \cite{Whit}. It was later extended in \cite{AMPW17}, see also \cite{TdC20}, to show the absence of non-trivial mode solutions with real frequencies. Recently, an unconditional bound for the horizon flux using a generalized physical space version of Whiting's transform has been derived in \cite{HeKl2}.

In Schwarzschild,  energy-Morawetz estimates for the Teukolsky equations were first obtained by Dafermos-Holzegel-Rodnianski \cite{D-H-R}. The proof relies on a physical-space analog of the Chandrasekhar's transformation \cite{Chand2} that converts the Teukolsky equations into a Regge-Wheeler type wave equation \cite{RW57}, to which the techniques developed for the scalar wave equation can be directly applied. Generalizations to Kerr spacetimes were achieved in the slowly rotating case by Ma \cite{Ma} and Dafermos-Holzegel-Rodnianski \cite{D-H-R-Kerr}, and for the full subextremal range by Millet \cite{Millet}\footnote{While \cite{Millet} derives sharp decay estimates for solutions to Teukolsky equations (for any half-integer spin) by relying on the spectral approach discussed in Section \ref{sec:intro:strategymetricpertubations}, one can easily adapt the methodology in that paper to derive a weak Morawetz estimate, though with a loss of several derivatives, see Section 12 in \cite{MaSz26}.} and Shlapentokh-Rothman-Teixeira da Costa \cite{SRTdC20, SRTdC23}.

%%%%%%%%%%%%%%%%%%%%%%%%%%%%%
 
\paragraph{\textit{Decay for LGS on Kerr}.}

%%%%%%%%%%%%%%%%%%%%%%%%%%%%%

A  first  quantitative proof of the  linear stability\footnote{Linear stability refers to the derivation of some decay estimates for solutions of LGS, even when they are not compatible with nonlinear applications. This is in fact the case for all the results on the decay of LGS on Kerr reviewed in this paragraph.} of Schwarzschild  spacetime  was established by  Dafermos-Holzegel-Rodnianski in \cite{D-H-R} relying on an analog in the linearized setting of the strategy based on curvature perturbations outlined in Section \ref{sec:intro:strategycurvarureperturbation} which starts with their analysis of the Teukolsky equation on Schwarzschild, see also Hung-Keller-Wang \cite{HKW} and Johnson   \cite{Johnson} for the linear stability of Schwarzschild in metric perturbations. The linear stability of Kerr for $|a|\ll m$ was obtained in Anderson-B\"ackdahl-Blue-Ma \cite{ABBMa2019} using curvature perturbations and starting from the analysis in \cite{Ma}, and in H\"afner-Hintz-Vasy \cite{HHV} using metric perturbations by following an analog at the linearized level of the strategy outlined in Section \ref{sec:intro:strategymetricpertubations}. The latter result was recently extended to a proof of the  linear stability of Kerr for all $|a|<m$ in H\"afner-Hintz-Vasy \cite{HHV1}. Note also that the results in \cite{ABBMa2019} extend to $|a|<m$ as well when combined with the energy-Morawetz estimates in \cite{SRTdC20, SRTdC23}.

%%%%%%%%%%%%%%%%%%%%%%%%%%%%%%%%%%%%%%%%%%%%%%%%%%%%%%%
 
\paragraph{\textit{Energy-Morawetz estimates for (non-)linear scalar waves on perturbations of Kerr}.}

%%%%%%%%%%%%%%%%%%%%%%%%%%%%%%%%%%%%%%%%%%%%%%%%%%%%%%%

To address the nonlinear stability of Kerr, it is important to extend the energy-Morawetz estimates for the scalar wave equation and Teukolsky equations on Kerr reviewed above to perturbations of Kerr. We start first with results concerning energy-Morawetz estimates for the scalar wave equation on perturbations of Kerr. Lindblad-Tohaneanu proved in \cite{LT18,LT20} global existence for small data solutions to a quasilinear wave equation in small perturbations of, first Schwarzschild and then Kerr with $|a|/m\ll 1$, using a microlocal approach adapted to the constant time level sets as in \cite{Ta-Toh}. Dafermos-Holzegel-Rodnianski-Taylor \cite{DHRT22} later  generalized those results to allow the presence of quadratic semilinear terms satisfying the null condition. The first result proving energy-Morawetz estimates for the scalar wave equation in perturbations of any subextremal Kerr background, which plays an important role in Section \ref{sec:energyMorawetzforPc}, is due to Ma and the author \cite{MaSz24}, see also \cite{DHRT24} by Dafermos-Holzegel-Rodnianski-Taylor for the proof of global existence of small data quasilinear waves in the full subextremal range on metrics coinciding with Kerr for $r$ large enough.

%%%%%%%%%%%%%%%%%%%%%%%%%%%%%%%%%%%%%%%%%%%%%%%%%%%%%%%
 
\paragraph{\textit{Energy-Morawetz estimates for Teukolsky equations on perturbations of Kerr}.}

%%%%%%%%%%%%%%%%%%%%%%%%%%%%%%%%%%%%%%%%%%%%%%%%%%%%%%%

Next, we discuss results concerning energy-Morawetz estimates for Teukolsky equations on perturbations of Kerr. This has been achieved in the context of the recent proofs of the nonlinear stability of Schwarzschild and of Kerr spacetimes for $|a|\ll m$: see Chapter 10 of \cite{KS} in the context of the nonlinear stability of Schwarzschild under polarized axisymmetry, Chapters 12 and 13 of \cite{DHRT} in the context of the nonlinear stability of Schwarzschild spacetimes for a codimension-3 set of initial data, and Chapter 9 of \cite{GKS22} in the context of the nonlinear stability of slowly rotating Kerr, i.e., with $|a|\ll m$. The proof of energy-Morawetz estimates for Teukolsky equations on perturbations of Kerr for $|a|<m$ has been recently obtained by Ma and the author \cite{MaSz26}, building in particular on \cite{MaSz24}, and plays a crucial role in  the present paper.

%%%%%%%%%%%%%%%%%%%%%%%%%%%%%%%%%%%%
 
\subsubsection{Stability of Schwarzschild}
\lab{sec:reviewproofstabSchwarzschild}

%%%%%%%%%%%%%%%%%%%%%%%%%%%%%%%%%%%%

Generic perturbations of Schwarzschild metrics converge asymptotically to Kerr metrics with $|a|\ll m$ so that the stability of Schwarzschild can only hold in a restricted class of perturbations. In particular, the first  nonlinear stability result of the Schwarzschild space  was established by Klainerman and the author in  \cite{KS} for axial polarized perturbations which is the  simplest  assumption ensuring that the final state is itself Schwarzschild. This paper developed a number of ideas that subsequently influenced both the proof of Kerr stability for small angular momentum, reviewed below, and the proof of the Kerr stability conjecture in the present paper, particularly with regard to the overall structure of the proof, as well as the procedure of modulation in infinite dimension based on GCM spheres already mentioned in Section \ref{sec:intro:strategycurvarureperturbation}. 

The result in \cite{KS} was later extended by Dafermos-Holzegel-Rodnianski-Taylor \cite{DHRT} to a co-dimension 3 subset of the initial data such that the final state is  still Schwarzschild, based on an additional   three  dimensional modulation.

%%%%%%%%%%%%%%%%%%%%%%%%%%%%%%%%%%%%
 
\subsubsection{Kerr stability for small angular momentum}
\lab{sec:reviewproofstabKerrforsmalla}

%%%%%%%%%%%%%%%%%%%%%%%%%%%%%%%%%%%%

Kerr stability for $|a|\ll m$ has been proved in the sequence of works \cite{KS-GCM1} \cite{KS-GCM2} \cite{KS:Kerr} by Klainerman and the author, \cite{GKS22} by Giorgi, Klainerman and the author, and \cite{Shen} by Shen. More precisely:
\begin{itemize}
\item The papers \cite{KS-GCM1} \cite{KS-GCM2} and \cite{Shen} address the modulation procedure in infinite dimensions. The construction of GCM spheres in \cite{KS-GCM1} \cite{KS-GCM2} is used both to construct the last sphere $S_*$ of the bootstrap spacetime (the one most towards the future), and for the construction of preferred spacelike hypersufaces in \cite{Shen}. In turn, \cite{Shen} is used to construct the spacelike future boundary $\Si_*$ intialized from the last sphere $S_*$. This construction of $S_*$ and $\Si_*$ corresponds to the infinite dimensional part of the modulation scheme. Finally, we choose the mass and angular momentum of the bootstrap spacetime respectively as the Hawking mass of $S_*$ and the definition of the angular momentum in \cite{KS-GCM2} on $S_*$ to deal with the remaining finite dimensional part of the modulation scheme. The results in \cite{KS-GCM1} \cite{KS-GCM2} and \cite{Shen} generalize the constructions in Chapter 9 of \cite{KS} and hold for general perturbations of Kerr for all angular momenta.

\item The paper \cite{KS:Kerr} contains in particular the following:
\begin{itemize}
\item The construction of two non-integrable\footnote{Recall from the discussion in Section \ref{sec:intro:maindifficultiesoftheproofifKerrstab} that the principal null pair of Kerr in \eqref{def:e3e4inKerr:intro} is non-integrable for $a\neq 0$ which justifies in turn the need for non-integrable gauges.} gauges, the Principal Geodesic (PG) gauge and the Principal Temporal (PT) gauge, where the PG gauge is used for decay estimates and loses derivatives while the PT gauge is well posed and used for high-order energy estimates.  

\item The statement of Kerr stability for small angular momentum and its decomposition in nine main intermediary steps, Theorems M0 to M8, in Section 3.7 of \cite{KS:Kerr}.

\item The derivation of decay estimates for all linearized curvature and Ricci coefficients except the Teukolsky scalars, and the derivation of higher order energy estimates for the Ricci coefficients.  
\end{itemize}
The results in \cite{KS:Kerr} hold for perturbations of Kerr for all $|a|<m$.

\item The paper \cite{GKS22} contains in particular the following:
\begin{itemize}
\item The extension of the formalism in \cite{Ch-Kl} to the non-integrable setting.

\item The derivation of decay estimates for the Teukolsky scalars.

\item The derivation of higher order energy estimates for the curvature components. 
\end{itemize}
The limitation to $|a|\ll m$ in \cite{GKS22} is only used in the proof of decay estimates for the Teukolsky scalars and of higher order energy estimates for the curvature components. 
\end{itemize}

In particular, the results in \cite{KS-GCM1} \cite{KS-GCM2} \cite{KS:Kerr} \cite{Shen} already hold for $|a|<m$, and to prove the main result in this paper, it suffices to extend the validity of the proof of decay estimates for the Teukolsky scalars and of higher order energy estimates for the curvature components in \cite{GKS22} to all subextremal angular momenta.

%%%%%%%%%%%%%%%%%%%%%%%%%%%%%%%%%%%%
 
\subsubsection{Kerr stability for well-prepared initial data}
\lab{sec:comparisionwithresultbyHintz}

%%%%%%%%%%%%%%%%%%%%%%%%%%%%%%%%%%%%

In the sequence of papers \cite{Hi2} \cite{Hi3} \cite{Hi1}, Hintz has recently proved Kerr stability in the full subextremal range for initial data perturbations satisfying 
\bea\lab{eq:intro:wellpreparedinitialdatainresultbyPeterHintz}
\g-\g_{a_0,m_0} = \textrm{finite polyhomogeneous expansion} + O(\ep_0r^{-3-\de})\quad\textrm{as}\quad r\to +\infty,
\eea
for $\de>0$ and for $\ep_0>0$ small enough. The proof follows the spectral approach in metric perturbations outlined in Section \ref{sec:intro:strategymetricpertubations}, and relies in particular on a Nash-Moser iteration scheme to close nonlinear estimates, on b-calculus to derive spectral estimates, and on the mode stability result in \cite{AHW}.

%%%%%%%%%%%%%%%%%%%%%%%%%%%%%%%%%%%%%%%%%%%%%%%%%%%%%%%%%%%%

\paragraph{\textit{The choice of initial data in \eqref{eq:intro:wellpreparedinitialdatainresultbyPeterHintz}}.}

%%%%%%%%%%%%%%%%%%%%%%%%%%%%%%%%%%%%%%%%%%%%%%%%%%%%%%%%%%%%

We provide further details on the specific choice of initial data \eqref{eq:intro:wellpreparedinitialdatainresultbyPeterHintz} in\footnote{We are grateful to Peter Hintz for helpful explanations regarding the choice of initial data in \cite{Hi1}, on which the discussion in this paragraph is based.}  \cite{Hi1}. To this end, we note that the proof in \cite{Hi1} relies on an iterative scheme in which, at each step: 
\begin{enumerate}
\item The linearized Einstein vacuum equations in wave gauge on a perturbation of Kerr are solved.

\item The resulting solution is decomposed into the sum of a well-decaying part and non-decaying or slowly decaying parts.

\item The non-decaying and slowly decaying parts are then identified with, and incorporated into, a linearized change of Kerr parameters and a linearized change of gauge for the next iteration.
\end{enumerate} 
In particular, carrying out the last step above requires proving that the non-decaying and slowly decaying parts, that do not arise from a change of Kerr parameters, lie in the image of the gauge modification map. For the linearized operator on Kerr, this follows from an explicit algebraic computation, see \cite{HHV1}. On the other hand, the gauge modification procedure in \cite{Hi1} extends to perturbations of Kerr only when the non-decaying and slowly decaying parts belong to a finite dimensional space. Now, as expected from the fact that the kernel of LGS due to general covariance is infinite dimensional, the slowly decaying part belongs, in general, to an infinite-dimensional space, thereby preventing the above strategy from being applied directly. Nevertheless, one can proceed in the following two special cases:
\begin{enumerate}
\item \textit{Finite polyhomogeneous expansion.} If the initial data admits a finite polyhomogeneous expansion, then the slowly decaying part lies in a finite-dimensional space.

\item \textit{Sufficiently fast decay.} If the initial data decays sufficiently fast, namely, $O(r^{-3-\de})$ as $r\to +\infty$, then the slowly decaying part already decays fast enough to close the nonlinear argument, without the need to identify it with, and incorporate it into, a linearized change of  gauge. 
\end{enumerate} 
The well-prepared choice of initial data in \eqref{eq:intro:wellpreparedinitialdatainresultbyPeterHintz} is in fact a combination of these two special cases.

%%%%%%%%%%%%%%%%%%%%%%%%%%%%%%%%%%%%%%

\paragraph{\textit{Comparison with Theorem \ref{MainThm-firstversion}}.}

%%%%%%%%%%%%%%%%%%%%%%%%%%%%%%%%%%%%%%

We now compare the result in \cite{Hi1} with our main Theorem \ref{MainThm-firstversion}. To this end, we first note that the finite polyhomogeneous expansion in \eqref{eq:intro:wellpreparedinitialdatainresultbyPeterHintz} has a very specific structure, whereas addressing the stability of a solution to an evolution PDE requires considering a structureless set of initial data. We therefore restrict our comparison to the part of \eqref{eq:intro:wellpreparedinitialdatainresultbyPeterHintz} relevant to the Kerr stability conjecture, namely, the $O(r^{-3-\de})$ structureless decay. In particular:
\begin{itemize}
\item Recall from \eqref{eq:rmk:remarkonroughtopologyinitialdatamainTh} that the initial data perturbations in Theorem \ref{MainThm-firstversion} correspond to an $O(r^{-\frac{3}{2}-\de})$ structureless decay. Thus, the structureless decay in \eqref{eq:intro:wellpreparedinitialdatainresultbyPeterHintz} falls off significantly faster compared to the initial data perturbation in Theorem \ref{MainThm-firstversion}. In fact, \eqref{eq:intro:wellpreparedinitialdatainresultbyPeterHintz} also fails to cover a substantial part of the slowly rotating initial Kerr perturbations corresponding to the stability result reviewed in Section \ref{sec:reviewproofstabKerrforsmalla}, which likewise exhibit an $O(r^{-\frac{3}{2}-\de})$ structureless decay.

\item Moreover, the structureless decay in \eqref{eq:intro:wellpreparedinitialdatainresultbyPeterHintz} also falls off significantly faster than $O(r^{-2})$, which is believed to be a threshold for a broad class of slowly decaying, physically relevant initial data, see for example Damour \cite{Damour}, Christodoulou \cite{Ch} and Kehrberger \cite{Kehr}.
\end{itemize}

%%%%%%%%%%%%%%%%%%%%%%%%%%%%%%%%%%%%%%

\paragraph{\textit{Analogy with the stability of Minkowski}.}

%%%%%%%%%%%%%%%%%%%%%%%%%%%%%%%%%%%%%%

To place the results in \cite{Hi1} and Theorem \ref{MainThm-firstversion} in context, it is useful to draw an analogy with results on the stability of Minkowski space. Stability of Minkowski was first proved by Friedrich \cite{Fri} for a restricted class of initial data that coincide with Schwarzschild outside of a compact set, with a framework which does not allow to extend the result beyond this setting. A proof for a much broader class of initial data was subsequently given by Christodoulou-Klainerman in \cite{Ch-Kl}. Their approach introduced a number of influential ideas that proved sufficiently flexible to allow for numerous further extensions, see \cite{SmuMinkreview} \cite{Shen2} and references therein. Another example of a flexible approach is the later proof of the stability of Minkowski in harmonic coordinates by Lindblad-Rodnianski, first proved for restricted data that coincide with Schwarzschild outside of a compact set in \cite{LiRo1}, and later extended to a much broader class of initial data in \cite{LiRo2}.

Now, while the class of initial data in \cite{Hi1} is broader than the one in \cite{Fri} and the two approaches differ in nature, they nevertheless share the common feature that neither can be extended beyond their respective class of initial data without introducing new ideas\footnote{In particular, extending the Kerr stability result in \cite{Hi1} would require new ideas to deal with infinite dimensional gauge modifications in the setting of that proof, see Remark 1.3 in \cite{Hi1}.}. By contrast, we expect the flexibility of the proof of Theorem \ref{MainThm-firstversion} to allow for further extensions beyond the class of initial data in \eqref{eq:rmk:remarkonroughtopologyinitialdatamainTh}, as was the case of the result of  Christodoulou-Klainerman.

%%%%%%%%%%%%%%%%%%%%%%%%%%%%%%%%%%%%
 
\subsubsection{Stability of other families of black holes}

%%%%%%%%%%%%%%%%%%%%%%%%%%%%%%%%%%%%

We end this section by mentioning results concerning back hole stability in other settings. In the context of Einstein equations with a positive cosmological constant, we have already mentioned the proof of the stability of Kerr-de Sitter for $|a|\ll m$ by Hintz-Vasy in \cite{HVas} which has been reproved by Fang in \cite{Fang1} \cite{Fang2} and recently extended to the full subextremal range conditionally on mode stability by Hintz-Petersen-Vasy \cite{HPV}. Note that these results concern the stationary region of Kerr-de Sitter, while the stability of the expanding region of Kerr-de Sitter has been obtained by Fournodavlos-Schlue in \cite{FoSc} for any angular momentum, see also the new proof by Hintz-Vasy \cite{HVas1} in harmonic coordinates. 

We also mention advances concerning the stability of charged black holes. The linear stability of Reissner-Nordstr\"om has been proved by Giorgi in \cite{Gio}, the linear stability of Kerr-Newman black holes  for small charges and small angular momenta has been proved by He in \cite{He}, and GCM constructions in the context of the nonlinear stability of Kerr-Newman have been recently obtained by Fang-Giorgi-Wan \cite{FGW1} \cite{FGW2}.

%%%%%%%%%%%%%%%%%%%%%%%%%%%%%%%%%%%%%%%%%%%%

\subsection{Precise statement of the main result and strategy of the proof}

%%%%%%%%%%%%%%%%%%%%%%%%%%%%%%%%%%%%%%%%%%%%

%%%%%%%%%%%%%%%%%%%%%%%%%%%%%%%%%%%%%%%%%%%%

\subsubsection{Precise version of the main result}
\lab{sec:statemeentofthemainresultonKerrstabilitycomplete}

%%%%%%%%%%%%%%%%%%%%%%%%%%%%%%%%%%%%%%%%%%%%

Our main result extends the validity of the main theorem of \cite{KS:Kerr} from initial perturbations of Kerr with parameters $|a_0|\ll m_0$ to the full subextremal range, i.e., $|a_0|<m_0$, which completes the proof of the Kerr stability conjecture. Reproducing the full statement of the main theorem of \cite{KS:Kerr}, stated in Section 3.4.3 of \cite{KS:Kerr}, and extended here to the range $|a_0|<m_0$, would require to repeat Sections 3.1 to 3.4 of \cite{KS:Kerr}, a substantial part of the material in these sections being irrelevant for the rest of this paper. We will thus provide a short statement mentioning where all relevant definitions can be found in Sections 3.1 to 3.4 of \cite{KS:Kerr}.

%%%%%%%%%%%%%%%%%%%%%%%%%%%%%%%%%%%%%%%%%%%%

\paragraph{\textit{Smallness constants}.}

%%%%%%%%%%%%%%%%%%%%%%%%%%%%%%%%%%%%%%%%%%%%

The following constants are involved in the statement of the main theorem:
\begin{itemize}
\item The constants $m_0>0$ and $a_0$, with $|a_0|<m_0$, are the mass and the angular momentum per unit mass of the Kerr solution relative to which our initial perturbation is measured. 

\item The integer $k_{large}$ corresponds to the maximum number of derivatives of the solution.

\item The size of the initial data layer norm is measured by $\ep_0>0$. 

\item $r_0>0$ will be used, given a coordinate $r$, to separate the spacetime in $\{r\leq r_0\}$ and $\{r\geq r_0\}$. 

\item $\deh>0$ will be tied to the definition $\{r=r_+(1-\deh)\}$ of the left causal boundary of the spacetime. 

\item $\dec$ is tied to decay rates  in $\tau$ where $\tau$ is a global time function on $\MM$.

\item $\dt$ is a small parameter involved in the $r$-power measuring the initial data perturbation.
\end{itemize}

  In what follows, $m_0$ and $a_0$ are fixed constants with $|a_0|<m_0$, $\deh$, $\dec$ and $\dt$ are fixed, sufficiently small, universal constants, and $r_0$ and $k_{large}$ are  fixed, sufficiently large, universal constants, satisfying 
\bea\lab{eq:constraintsonthemainsmallconstantsepanddelta:intro}
0<\deh\ll 1-\frac{|a_0|}{m_0},\qquad 0<\dec<\frac{\dt}{2}\ll 1, \qquad r_0\gg \max\{m_0,1\},\qquad k_{large}\gg \frac{1}{\dt}.
\eea
Then, $\ep_0$ is chosen such that
\bea\lab{eq:constraintsonthemainsmallconstantsepanddelta:bis:intro}
0<\ep_0\ll \min\left\{\deh, \dec, \dt, \frac{1}{r_0}, \frac{1}{k_{large}}, 1-\frac{|a_0|}{m_0}\right\}, \qquad \ep_0\ll \frac{|a_0|}{m_0}\quad \textrm{ in the case }a_0\neq 0.
\eea
Also, we introduce the integer $k_{small}$ which corresponds to the  number of derivatives  for which the solution satisfies decay estimates. It is related to $k_{large}$ by
\bea\lab{eq:choiceksmallmaintheorem}
k_{small}=\left \lfloor\frac 1 2 k_{large}\right \rfloor +1.
\eea

From now on, in the rest of the paper, $\lesssim$ means bounded by a constant depending only on geometric universal constants (such as Sobolev embeddings, elliptic estimates,...) as well as the constants 
$$m_0,\, a_0, \, \deh, \,\dec,\, \dt, \, r_0, \, k_{large}$$
\textit{but not on} $\ep_0$.

%%%%%%%%%%%%%%%%%%%%%%%%%%%%%%%%%%%%%%%%%%%%

\paragraph{\textit{Initial data layer and admissible future null complete spacetimes}.}

%%%%%%%%%%%%%%%%%%%%%%%%%%%%%%%%%%%%%%%%%%%%

\begin{figure}[ht!]
\centering
\includegraphics[scale=0.35]{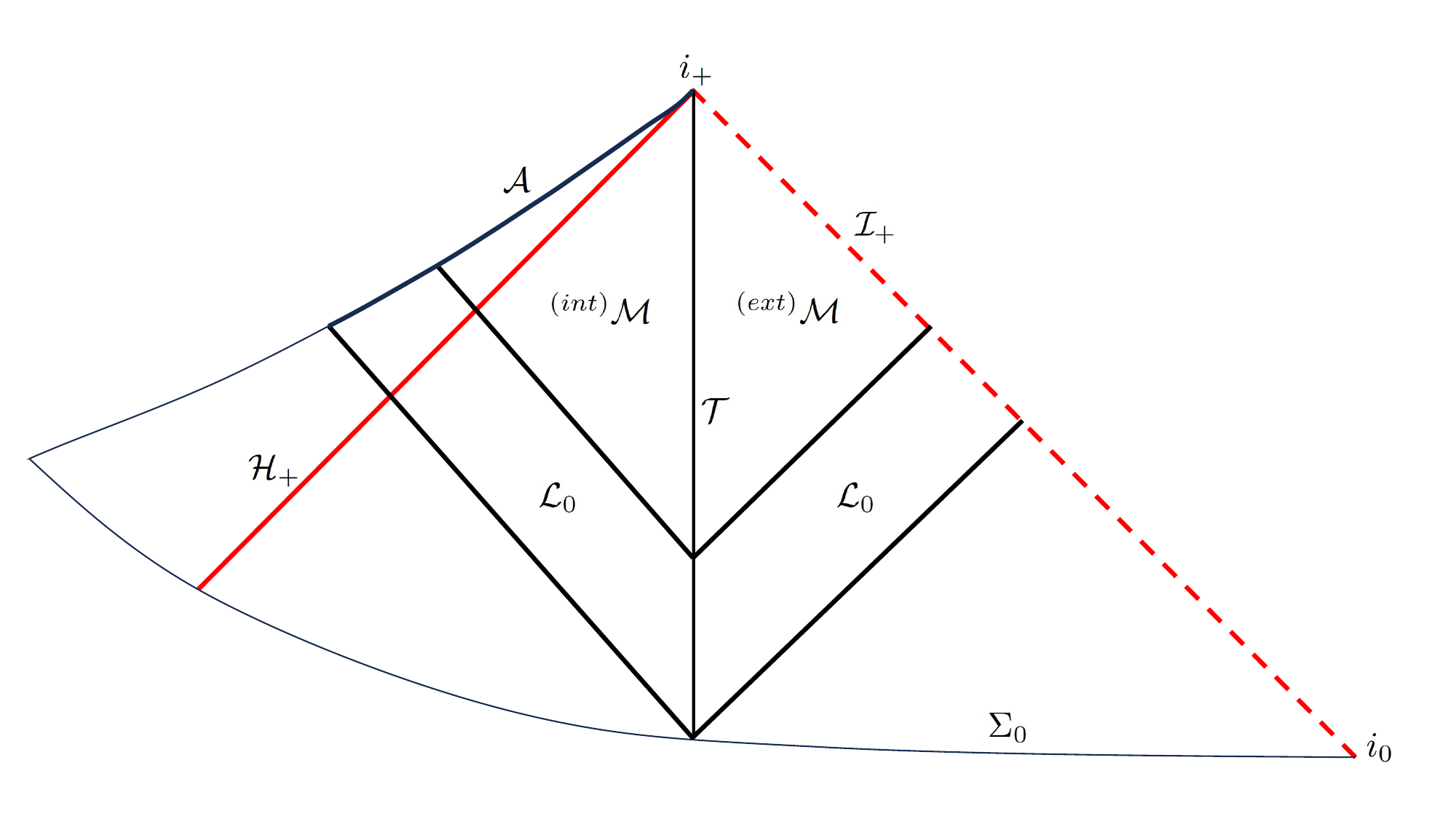}
\caption{\footnotesize{Penrose diagram of  an admissible future complete spacetime}}
\label{fig:penrosediagramfuturecomplete}
\end{figure}

We now recall the choice of initial data layer and the notion of admissible future null complete spacetimes introduced in Section 3 of \cite{KS:Kerr}, see also the Penrose diagram in Figure \ref{fig:penrosediagramfuturecomplete}.

We start with the notion of $(\ep_0, k)$-admissible initial data layer introduced in Definition 3.4.2 of \cite{KS:Kerr}:
\begin{itemize}
\item The initial data layer, denoted by $\LL_0$, is defined in Section 3.1 of \cite{KS:Kerr}, see also Figure \ref{fig:penrosediagramfuturecomplete}.

\item In addition, following Definition 3.4.2 in \cite{KS:Kerr}, we say that an initial data layer $\LL_0=\LL_0(a_0, m_0)$, defined  as in  Section 3.1 of \cite{KS:Kerr}, is admissible   if it lies in the future of an asymptotically flat initial data set,  supported  on a spacelike hypersurface  $\Si_0$, of ADM mass $ m_0  $ and angular momentum per unit mass $a_0$, see the lower part of Figure \ref{fig:penrosediagramfuturecomplete}.

\item Also,  we say that   $\LL_0$ is   $(\ep_0, k)$-admissible if it verifies the bound in Definition 3.4.2 in \cite{KS:Kerr} measuring the size of the initial data perturbation w.r.t. the initial layer norm $\Ik_k$ defined as in Section 3.3.6 of \cite{KS:Kerr}, with $k$ denoting the number of derivatives in $\Ik_k$ and $\ep_0$ measuring its size.
\end{itemize}

\begin{remark}
\lab{remark:resultsKl-Ni}
We can immediately identify a large class of initial data sets  on $\Si_0$  which generate  $(\ep_0, k)$-admissible  initial data layers by:
\begin{enumerate}
\item considering on $\Si_0$ the initial data sets for Kerr stability constructed in Theorem 3.3 of \cite{FST} by Fang, Touati and the author, with the constants $(q, \de)$ in that result chosen as\footnote{This corresponds to initial data sets measured in weighted Sobolev spaces implying in particular the pointwise estimate $\g=\gam+O(\ep_0r^{-\frac{3}{2}-\dt})$ as $r\to +\infty$ along the spacelike initial hypersurface $\Si_0$, as well as corresponding estimates for higher order weighted derivatives.} $q=1$ and $\de=\frac{1}{2}+\dt$, 

\item propagating the above initial data sets all the way to the initial data layer by a combination of local existence in the left region below $\LL_0$ in Figure \ref{fig:penrosediagramfuturecomplete} and of Theorem 1.7  in \cite{Shen1} by Shen, with the constant $s$ in that result chosen as $s=4+2\dt$ in the external region on the right below $\LL_0$. 
\end{enumerate}
\end{remark}

We also need the notion of an admissible future null complete spacetime introduced in Definition 3.4.5 of \cite{KS:Kerr} which is an asymptotically flat vacuum spacetime  $\MM$, as in Figure \ref{fig:penrosediagramfuturecomplete}, satisfying in particular the following properties:
\begin{itemize}
\item  It is   a  future  development   of  an admissible initial data layer  set $\LL_0(a_0, m_0)$, in the sense  of  Definition 3.4.2 in \cite{KS:Kerr} recalled above.

 \item The future null infinity  $\II^{+}$ of   $\MM$  is  complete.  The other future boundary of $\MM$ is given by the spacelike hypersurface $\AA=\{r=r_+(1-\dhor)\}$, where $r$ is a global coordinate function on $\MM$, and belongs to the complement of $\JJ^{-}(\II^{+})$.

 \item $\MM=\Mext\cup\Mint$, where $\Mext=\{r\geq r_0\}$ and $\Mint=\{r\leq r_0\}$.
 
 \item $\TT=\Mext\cap\Mint=\{r=r_0\}$ is  timelike.
\end{itemize}

Finally, given an  admissible  future null  complete spacetime $\MM$ as in Definition 3.4.5 of \cite{KS:Kerr} recalled above, and constants $(a_\infty, m_\infty)$, $|a_\infty|<m_\infty$, we define  the norms $\Nk^{(Sup)}_{k_{large}}(a_\infty, m_\infty )$, $\Nk^{(Dec)}_{k_{small}}(a_\infty, m_\infty )$ as in  Section 3.3.5 of \cite{KS:Kerr}, where\footnote{Note that the first and the second line of \eqref{eq:assumptionsonMMforpartII} correspond essentially to the estimates $\Nk^{(Sup)}_{k_{large}}(a_\infty, m_\infty )\leq \ep$ and $\Nk^{(Dec)}_{k_{small}}(a_\infty, m_\infty )\leq \ep$ with $(a_\infty, m_\infty )\to (a, m)$ and $k_{large}\to\kl$.} $\Nk^{(Sup)}_{k_{large}}(a_\infty, m_\infty )$ denotes $r$-weighted supremum norms for $k_{large}$ weighted derivatives of the main linearized quantities, and where, given a suitable global time function $\tau$ on $\MM$ in addition to the coordinate $r$, $\Nk^{(Dec)}_{k_{small}}(a_\infty, m_\infty )$ denotes $(r, \tau)$-weighted supremum norms for $k_{small}$ weighted derivatives of the main linearized quantities, with $k_{small}$ as in \eqref{eq:choiceksmallmaintheorem}.

%%%%%%%%%%%%%%%%%%%%%%%%%%%%%%%%%%%%%%%%%%%%%%%%%%%%

\paragraph{\textit{Statement of the main theorem}.} 

%%%%%%%%%%%%%%%%%%%%%%%%%%%%%%%%%%%%%%%%%%%%%%%%%%%%

We are now ready to give the following precise version of our main theorem, which exactly corresponds to the extension of the validity of the main theorem of \cite{KS:Kerr}, stated in Section 3.4.3 of \cite{KS:Kerr},  from initial perturbations of Kerr with parameters $|a_0|\ll m_0$ to the full subextremal range, i.e., $|a_0|<m_0$.

\begin{theorem}[Main theorem, precise version]
\lab{thm:mainKerrstab:preciseversion}
Let $\LL=\LL(a_0, m_0) $  be an     $(\ep_0, k_{large}+10)$-admissible initial data layer as in Definition 3.4.2 of \cite{KS:Kerr},  with  $|a_0|<m_0$,  $ k_{large}$ sufficiently 
large,    and   $\ep_0>0$  sufficiently small. In particular, we assume the bound (3.4.7) in \cite{KS:Kerr} for the initial layer norm measuring the perturbation. Then   $\LL=\LL(a_0, m_0) $   possesses an    admissible future complete development  $\MM_\infty$  as in Definition 3.4.5 of \cite{KS:Kerr}. 

Moreover, there exist constants $(a_\infty, m_\infty)$,  $|a_\infty| < m_\infty$,  such that the following  estimates hold true,      relative to the norms     $\Nk^{(Sup)}_{k_{large}}=\Nk^{(Sup)}_{k_{large}}(a_\infty, m_\infty ), \, \Nk^{(Dec)}_{k_{small}}=\Nk^{(Dec)}_{k_{small}}(a_\infty, m_\infty ) $ introduced in Section 3.3.5 of \cite{KS:Kerr}, 
\bea\lab{def:bootstrapasumptionsglobalnorms}
\Nk^{(Sup)}_{k_{large}}+\Nk^{(Dec)}_{k_{small}} +| a_\infty- a_0| +|m_\infty- m_0| \le C\ep_0
\eea
where $C$ is a large universal constant and $k_{small}=\left \lfloor\frac 1 2 k_{large}\right \rfloor +1$.

Finally, all the other conclusions in the main theorem of \cite{KS:Kerr} also hold in the case $|a_0|<m_0$. In particular, $m_\infty$ coincides with the final Bondi mass, a Bondi mass law formula holds, and $a_\infty m_\infty$ is related to the limit towards $i_+$ of a notion of  angular momentum on $\II_+$, see Section 3.4.3 of \cite{KS:Kerr} for more details.
\end{theorem}

%%%%%%%%%%%%%%%%%%%%%%%%%%

\subsubsection{Main intermediary results}
\lab{sec:mainintermediaryresults:intro}

%%%%%%%%%%%%%%%%%%%%%%%%%%

The proof of Theorem \ref{thm:mainKerrstab:preciseversion}, as the one of the main theorem in Section 3.4 of \cite{KS:Kerr}, is divided in a sequence of nine intermediary steps, called Theorem M0--M8, see Section 3.7 in \cite{KS:Kerr}. Also, recall from Section \ref{sec:reviewproofstabKerrforsmalla} that the results in \cite{KS-GCM1} \cite{KS-GCM2} \cite{KS:Kerr} \cite{Shen} already hold for all subextremal angular momenta so that the proof of Theorem M0, M3--M7, and the control of Ricci coefficients in Theorem M8 already hold for all subextremal angular momenta. To conclude the proof of Theorem \ref{thm:mainKerrstab:preciseversion}, it thus suffices to extend the validity of the decay estimates for Teukolsky scalars in Theorem M1 and Theorem M2 and the curvature estimates in Theorem M8 in \cite{GKS22} to all subextremal angular momenta, i.e., to prove Theorems \ref{theoremM1:intro}, \ref{theoremM2:intro} and \ref{THEOREMM8:INTRO} stated below\footnote{See respectively Theorems \ref{theoremM1:Chap11}, \ref{thm:restatementofTheoremM2} and \ref{prop:rpweightedestimatesiterationassupmtionThM8} for the precise version of Theorems \ref{theoremM1:intro}, \ref{theoremM2:intro} and \ref{THEOREMM8:INTRO}.}.

\begin{theorem}[Extension of Theorem M1 in \cite{KS:Kerr} to the full subextremal range, rough version]
\lab{theoremM1:intro}
Assume that the spacetime\footnote{In Theorems \ref{theoremM1:intro}, \ref{theoremM2:intro} and \ref{THEOREMM8:INTRO}, $\MM$ is a bootstrap spacetime, see Section \ref{section:SpacetimeMM-chap6} for the precise geometric set-up.} $\MM$  verifies the assumptions {stated in Section \ref{section:assumptionsneededforstatementofTheoremM1}}, as well as  the assumption \eqref{eq:controlofinitialdataforThM8-intro} on the  initial data. Then, if $\ep_0>0$ is sufficiently small, there exists $\dee>\dec$ such that the curvature component $A$ introduced in Definition \ref{def:complexRicciandcurvaturecoefficients} satisfies in $\MM$, 
\beaa
\sup_{\MM}\Big(r^2\tau^{1+\dee}+r^3(2r+\tau)^{\frac{1}{2}+\dee}\Big)|A|+\sup_{\MM}\Big(r^3\tau^{1+\dee}+r^4\tau^{\frac{1}{2}+\dee}+r^{\frac{9}{2}+\dec}\Big)|\nab_3A|\\
+\sup_{\MM}\Big(r^4\tau^{1+\dee}+r^{\frac{9}{2}+\dec}\Big)|\nab_3^2A| &\les& \ep_0,
\eeaa
as well as corresponding estimates for a suitable number of higher order weighted derivatives.
\end{theorem}

\begin{theorem}[Extension of Theorem M2 in \cite{KS:Kerr} to the full subextremal range, rough version]
\lab{theoremM2:intro}
Assume that the global frame of $\MM$ satisfies the assumptions of {Section} \ref{sec:assumptionsontheframe:Chapter12}, as well as  the assumption \eqref{eq:controlofinitialdataforThM8-intro} on the  initial data. Then, we  have the following decay estimates for the curvature component $\Ab$ introduced in Definition \ref{def:complexRicciandcurvaturecoefficients} along the spacelike hypersurface\footnote{The spacelike hypersurface $\Si_*$ is a part of the future boundary of the bootstrap spacetime $\MM$, see Section \ref{section:SpacetimeMM-chap6}.} $\Si_*$ of $\MM$ introduced in Section \ref{section:SpacetimeMM-chap6}, 
\beaa
\int_{\Sigma_*}\tau^{2+2\dec}|\Ab|^2 \les \ep_0^2,
\eeaa
as well as corresponding estimates for a suitable number of higher order weighted derivatives.
\end{theorem}

Next, we extend the curvature estimates in Theorem M8, i.e., the estimates for curvature stated in Theorem {9.65} of \cite{KS:Kerr}, to the full subextremal range. To this end, we introduce weighed $L^2$ type norms $\Rk_k$ and $\Sk_k$ respectively for curvature and Ricci coefficients, and decompose $\Rk_k$ and $\Sk_k$ into their restrictions $\Rkint$, $\Skint$ to $\Mint$ and $\Rkext$, $\Skext$ to $\Mext$, see Section \ref{subsection:MainNormsM8} for the precise definitions. 

\begin{theorem}[Extension of Theorem {9.65} in \cite{KS:Kerr} to the full subextremal range, rough version]
\lab{THEOREMM8:INTRO}
Assume that the spacetime $\MM$  as defined in Section \ref{sec:geometricsetupspacetimeMMforsections789curvatureestimatesTheoremM8}  
verifies the assumptions in Section \ref{section:MainAss.ThmM8}. Then,  we  have the following boundedness estimates for all components of curvature, for $k_{small}-1\leq J\leq k_{large}+6$,
\bea
   \Rkint_{J+1}^2&\les&  r_0^{18}\Big(\ep_J^{\frac{1}{6}}(\Sk_{J+1}+\Rk_{J+1})^{\frac{11}{6}}+\ep_0^2\Big),\\
\Rkext_{J+1}^2 &\les& r_0^{3+\de_B}\Rkint^2_{J+1}+ r_0^{-\de_B} \Skext^2_{J+1} +\ep_J^2 +\ep_0^2 {+\ep(\Rk_{J+1}+\Sk_{J+1})^2},
\eea
where the constants in $\les$ are independent of $r_0$ and $\ep_J$  is such that $\Sk_J+\Rk_J\leq \ep_J$.
\end{theorem}

%%%%%%%%%%%%%%%%%%%%%%%%%%%%%%%%%%%%%%%%%%%%

\subsubsection{Main strategy of the proof}

%%%%%%%%%%%%%%%%%%%%%%%%%%%%%%%%%%%%%%%%%%%%

%%%%%%%%%%%%%%%%%%%%%%%%%%%%%%%%%%%%%%%%%%%%%%%%%%%%

\paragraph{\textit{Strategy of the proof of Theorems \ref{theoremM1:intro} and \ref{theoremM2:intro}}.}

%%%%%%%%%%%%%%%%%%%%%%%%%%%%%%%%%%%%%%%%%%%%%%%%%%%%

The core of the proof of Theorems \ref{theoremM1:intro} and \ref{theoremM2:intro} is the derivation of energy-Morawetz estimates for Teukolsky equations which follows from the following steps:
\begin{enumerate}
\item First, we exhibit the structure of the Teukolsky wave-transport system in perturbations of Kerr in Section \ref{sec:sectionontheteukolskywavetransportsystem} and prove in particular that it can be put in the form needed to apply the main energy-Morawetz estimates of \cite{MaSz26}. The derivation of the Teukolsky wave-transport system is given in the companion paper \cite{Sze}.

\item Then, we need additional constructions in order to fulfill the assumptions in \cite{MaSz26}. More precisely, we construct in Section \ref{sec:prelimnariesfortheproofofdecayesetimatesforTeuk}:
\begin{itemize} 
\item a coordinate system with additional properties compared to the ones in \cite{KS:Kerr},

\item complex 1-forms satisfying specific algebraic identities,

\item a triplet of 1-forms satisfying specific algebraic identities needed for the procedure of regular scalarization of tensorial wave equations introduced in \cite{MaSz26}.
\end{itemize}
\item With these constructions, we can adapt in Section \ref{sec:energyMorawetzesitmatesforTeukoslkyonMM:upto15derivatives} the energy-Morawetz estimates for Teukolsky equations in \cite{MaSz26} to our setting. The difference compared to \cite{MaSz26} is that the spacetime $\MM$ does not extend all the way to $\II_+$ but instead stops at the spacelike future boundary $\Si_*$ which requires to revisit a number of the constructions in \cite{MaSz26}. In the end, this yields energy-Morawetz estimates for up to 14 derivatives of solutions to Teukolsky equations on perturbations of Kerr, see Theorem \ref{thm:main:MaSz26}. The proof of Theorem \ref{thm:main:MaSz26} is given in the companion paper \cite{Sze}.

\item On then upgrades this energy-Morawetz estimates for solutions to Teukolsky equations to higher order derivatives in Section \ref{sec:energyMorawetzesitmatesforTeukoslkyonMM:higherorderderivatives} using a  somewhat standard procedure which relies as a first step on the commutation with Lie derivatives of approximate Killing fields. 
\end{enumerate}
Once these energy-Morawetz estimates for solutions to the Teukolsky equations on perturbations of Kerr have been established for an arbitrary number of derivatives, the remainder of the proofs of Theorems \ref{theoremM1:intro} and \ref{theoremM2:intro} follows readily. Indeed, the remaining steps, such as $r^p$-weighted estimates, were already proven for $|a|<m$ in Chapters 10 to 12 of \cite{GKS22}.

%%%%%%%%%%%%%%%%%%%%%%%%%%%%%%%%%%%%%%%%%%%%%%%%%

\paragraph{\textit{Strategy of the proof of Theorem \ref{THEOREMM8:INTRO}}.}

%%%%%%%%%%%%%%%%%%%%%%%%%%%%%%%%%%%%%%%%%%%%%%%%%

The goal of Theorem \ref{THEOREMM8:INTRO} is to estimate the renormalized curvature components $(A, B, \Pc, \Bb, \Ab)$, where the curvature components $(A, B, P, \Bb, \Ab)$ are introduced in Definition \ref{def:complexRicciandcurvaturecoefficients} and $\Pc=P+\frac{2m}{q^3}$ with $q=r+ia\cos\th$. The proof is part of the iteration argument set up in Section 9.4.7 of \cite{KS:Kerr} on the number of derivatives of the solution. As a consequence, we may assume that we control $J$ derivatives, and we would like to recover $J+1$ derivatives of $(A, B, \Pc, \Bb, \Ab)$ for $k_{small}-1\leq J\leq k_{large}+6$. 

The first step, carried out in Section \ref{sec:energyMorawetzforPc}, is to derive energy-Morawetz estimates for $r^2\Pc$. As in Chapter 14 of \cite{GKS22}, the proof of such energy-Morawetz estimates for $r^2\Pc$ relies on the following structural ingredients:
\begin{enumerate}
\item the fact that $P$ satisfies, as a consequence of the Bianchi system, a scalar wave equation which, unlike the case of Teukolsky equations, couples to the rest of the system even at the linearized level,

\item the fact that we have proved energy-Morawetz estimates for the solutions to Teukolsky equations which do not lose derivative and hence also hold in the context of Theorem \ref{THEOREMM8:INTRO},

\item the existence of an identity relating angular derivatives of $\Pc$ to Teukolsky scalars as a consequence of the Bianchi system.
\end{enumerate}

Based on these structural facts, the first step in Chapter 14 of \cite{GKS22} is to derive energy-Morawetz estimates for first order derivatives of $P$ w.r.t. approximate Killing fields, see Proposition 14.1.8 in \cite{GKS22}. This crucially uses the smallness of $a$ to absorb one of the coupling terms on the RHS of the corresponding wave equation, namely the one generated by the potential $W$ in (14.2.6) of \cite{GKS22} which satisfies $\Im(W)=O(mar^{-4})$. In order to derive energy-Morawetz estimates for $r^2\Pc$ for all $|a|<m$, we thus must proceed differently:
\begin{enumerate}
\item We commute the wave equation for $P$ with $\nab_{\Rhat}^2$, where $\Rhat$ coincides with $\frac{\De}{r^2+a^2}\pr_r$ in Kerr, and obtain, after linearization, renormalization and using again the wave equation for $P$, an energy-Morawetz for $r^2\nab_{\Rhat}^2\Pc$ with angular derivatives of $\Pc$ on the RHS.

\item Together with the control of angular derivatives that follows from a combination of energy-Morawetz estimates for the solutions to Teukolsky equations and the above-mentioned identity relating angular derivatives of $\Pc$ to Teukolsky scalars, we obtain an unconditional energy-Morawetz estimate for $r^2\nab_{\Rhat}^2\Pc$ and angular derivatives of $r^2\Pc$. 

\item We then use the wave equation for $P$ to recover the remaining derivatives, and conclude with the use of $r^p$ estimates and redshift estimates to recover weighted derivatives and remove degeneracies near the horizon.\end{enumerate}

 Next, as in Chapter 15 of \cite{GKS22}, we rely in Section \ref{sec:energyMorawetzforABBbAb} on the energy-Morawetz estimates for $r^2\Pc$ of Section \ref{sec:energyMorawetzforPc} to derive energy-Morawetz estimates for the remaining curvature components $A$, $B$, $\Bb$ and $\Ab$. Now, the restriction to $|a|\ll m$ in the derivation of energy-Morawetz estimates for $A$, $B$, $\Bb$ and $\Ab$ in Chapter 15 of \cite{GKS22} appears when recovering  angular derivatives from the non-integrable Hodge operators that appear in the Bianchi system which in turn relies on the non-integrable Hodge estimates of Corollary 13.4.1 in \cite{GKS22}. The main new ingredient to derive energy-Morawetz estimates for the curvature components $A$, $B$, $\Bb$ and $\Ab$ in the full subextremal range consists in removing the condition on the size of $|a|$ in non-integrable Hodge estimates by proceeding as follows:
 \begin{enumerate}
\item First, carefully revisiting the proof in Chapter 15 of \cite{GKS22} reveals that the restriction on the size of $|a|$ in Chapter 15 of \cite{GKS22} can be more precisely attributed to the use of non-integrable Hodge estimates on the level hypersurfaces $\Si(\tau)$ of the time function $\tau$ and on the future spacelike boundaries $\AA$ and $\Si_*$ of the spacetime $\MM$ when recovering the control of the energy and the fluxes of angular derivatives of $A$, $B$, $\Bb$ and $\Ab$.

\item We thus need to improve non-integrable Hodge estimates on $\Si(\tau)$, $\AA$ and $\Si_*$. To this end, we first notice that the size of $|a|$ is not needed when $r$ is sufficiently large which already takes care of the fluxes on $\Si_*$ and $\Si_{r\geq r_1}(\tau)$ for $r_1\gg m$.

\item Then, we show a second version of non-integrable Hodge estimates which allow to remove the condition on the size of $|a|$ on $\Si(\tau)\cap\{r_+\leq r\leq r_1\}$ provided we have:
\begin{itemize}
\item for $A$ and $B$, an a priori control of their $e_3$ derivatives,

\item for $\Ab$ and $\Bb$, an a priori control of their $e_4$ derivatives,
\end{itemize}
which turns out to be the case by using in particular the control of $\Pc$ in Section \ref{sec:energyMorawetzforPc}.

\item A last step, relying in particular on redshift estimates applied to the commutation w.r.t. Lie derivatives of approximate Killing fields of the Bianchi system allows to treat the remaining region $\Si_{r\leq r_+}(\tau)$ as well as the flux on $\AA$.   
\end{enumerate}

Finally, to conclude the proof of Theorem \ref{THEOREMM8:INTRO}, it remains to combine the energy-Morawetz estimates for $(A, B, \Pc, \Bb, \Ab)$ derived in Sections \ref{sec:energyMorawetzforPc} and \ref{sec:energyMorawetzforABBbAb} with the $r^p$-estimates for the Bianchi system in Chapter 16 of \cite{GKS22} which are already valid without restriction on $a$.

%%%%%%%%%%%%%%%%%%%%%%%%%%%%%%%%%%%%

\subsection{Organization of the paper}

%%%%%%%%%%%%%%%%%%%%%%%%%%%%%%%%%%%%

In Section \ref{sec:nonintergrableformalism}, we review the non-integrable formalism introduced in \cite{GKS20} \cite{GKS22}. Next, Section \ref{sec:sectionontheteukolskywavetransportsystem} is dedicated to the derivation of the Teukolsky wave-transport system in perturbations of Kerr. Section \ref{sec:prelimnariesfortheproofofdecayesetimatesforTeuk} provides the geometric set-up as well as geometric constructions needed for the derivation of energy-Morawetz estimates for Teukolsky equations. Then, we prove Theorems \ref{theoremM1:intro} and \ref{theoremM2:intro} in Section \ref{sec:decayestimatesforAandAb:analogPartIIGKS22}. Finally, Section \ref{sec:higherordercurvatureestimatesforprovingThM8largea:00} is dedicated to the proof of Theorem \ref{THEOREMM8:INTRO}.

%%%%%%%%%%%%%%%%%%%%%%%%%%%%%%%%%%%%

\subsection{Acknowledgements}

%%%%%%%%%%%%%%%%%%%%%%%%%%%%%%%%%%%%

The author would like to emphasize that this work particularly owes to the contributions of his collaborators Elena Giorgi, Sergiu Klainerman and Siyuan Ma, as well as the ones of his former student Dawei Shen. 
In particular, he wishes to express his deepest gratitude to Sergiu Klainerman without whom this program would have been unthinkable, and he would also like to acknowledge the essential contribution of Siyuan Ma in  managing the last (but not least) hurdles. The author is supported by the ERC grant ERC-2023 AdG 101141855 BlaHSt.

%%%%%%%%%%%%%%%%%%%%

\section{Non-integrable formalism}
\lab{sec:nonintergrableformalism}

%%%%%%%%%%%%%%%%%%%%%

In this section, we review the formalism for non-integrable structures introduced in \cite{GKS20} \cite{GKS22}. This will be used throughout the paper.

%%%%%%%%%%%%%%%%%%%%%%%%%%%%

\subsection{Definition of the non-integrable formalism}

%%%%%%%%%%%%%%%%%%%%%%%%%%%%

%%%%%%%%%%%%%%%%%%%%%%%%%%

\subsubsection{Null pairs and horizontal structures}
\lab{subsection:review-horiz.structures}

%%%%%%%%%%%%%%%%%%%%%%%%%%

Consider a fixed null pair $e_3, e_4$, i.e., $\g(e_3, e_3)=\g(e_4, e_4)=0$,   $\g(e_3, e_4)=-2,$ and
 denote  by  $\O(\MM)$ the vectorspace  of horizontal vectorfields $X$  on $\MM$, i.e., $\g(e_3, X)= \g(e_4, X)=0$.
  Given a fixed   orientation  on $\MM$,  with corresponding  volume form  $\in$,  we define  the induced 
 volume form on   $\O(\MM)$ by,  
 \beaa
  \in(X, Y):=\frac 1 2\in(X, Y, e_3, e_4). 
  \eeaa
 A null  frame on $\MM$ consists of a choice of horizontal vectorfields  $e_1, e_2$, such that\footnote{We use greek 
 indices $\a, \b, \ga$ for $1,2,3,4$ and latin indices $a,b$ for $1,2$.}
 \beaa
 \g(e_a, e_b)=\de_{ab}\qquad  a, b=1,2.
 \eeaa  
 The commutator $[X,Y]$ of two horizontal vectorfields
may fail however to be horizontal. We say that the pair $(e_3, e_4 )$ is integrable if   $\O(\MM)$  forms an integrable distribution, i.e., $X, Y\in\O(\MM) $ implies that $[X,Y]\in\O(\MM)$. As is well-known,  the  principal null pair in Kerr fails to be integrable.
Given an arbitrary vectorfield $X$, we denote by $^{(h)}X$
its  horizontal projection, 
\beaa
{}^{(h)}X := X+ \frac 1 2 \g(X,e_3)e_4+ \frac 1 2   \g(X,e_4) e_3. 
\eeaa
A  $k$-covariant tensor-field $U$ is said to be horizontal,  $U\in \O_k(\MM)$,
if  for any $X_1,\ldots, X_k$ we have 
\beaa
U(X_1,\ldots, X_k)=U( ^{(h)} X_1,\ldots, {}^{(h)}X_k).
\eeaa

\begin{definition}\label{definition-SS-real}
We denote by $\sk_0=\sk_0(\MM, \mathbb{R})$ the set of pairs of real scalar functions on $\MM$,  by $\sk_1=\sk_1(\MM, \mathbb{R})$ the  set of real horizontal $1$-forms  on $\MM$   and by $\sk_2=\sk_2(\MM, \mathbb{R})$ the set of symmetric traceless   horizontal real $2$-tensors on $\MM$.
\end{definition}

\begin{definition}\label{definition-hodge-duals}
We define the dual of $\xi\in\sk_1$ and $U\in\sk_2$ by
\beaa
\dual \xi_{a}:=\in_{ab}\xi_b,\qquad \dual U_{ab}:=\in_{ac} U_{cb}.
\eeaa
\end{definition}

Note that given $\xi, \eta\in\sk_1$ and $U\in\sk_2$, we have 
\beaa
\dual(\dual \xi)=-\xi, \qquad \dual (\dual U)=-U,\qquad \dual\xi \c  \eta=-\xi\c\dual\eta.
\eeaa
 Also, given  $\xi, \eta\in\sk_1 $,  $U, V\in \sk_2$  we denote
\beaa
\xi\c \eta&:=&\de^{ab} \xi_a\eta_b,\qquad 
\xi\wedge\eta:=\in^{ab} \xi_a\eta_b=\xi\c\dual \eta,\qquad 
(\xi\hot \eta)_{ab}:=\xi_a \eta_b +\xi_b \eta_a-\de_{ab} \xi\c \eta,\\
(\xi\c U)_a&:=&\de^{bc} \xi_b U_{ac}, \qquad  (U\wedge V)_{ab} := \in^{ab}U_{ac}V_{cb}.
\eeaa

For any $ X, Y\in \O(\MM)$ we define  the induced metric $g(X, Y):=\g(X, Y)$ and the null second fundamental forms
\bea
\chi(X,Y):=\g(\D_Xe_4 ,Y), \qquad \chib(X,Y):=\g(\D_Xe_3,Y).
\eea
Observe that  $\chi$ and $\chib$  are  symmetric if and only if   the horizontal structure is integrable. Indeed this  follows easily from the following formulas
 \beaa
 \chi(X,Y)-\chi(Y,X)&=&\g(\D_X e_4, Y)-\g(\D_Ye_4,X)=-\g(e_4, [X,Y]),\\
 \chib(X,Y)-\chib(Y,X)&=&\g(\D_X e_3, Y)-\g(\D_Ye_3,X)=-\g(e_3, [X,Y]).
\eeaa
  Note  that  we  can view  $\chi$ and $\chib$ as horizontal 2-covariant tensor-fields
 by extending their definition to arbitrary vectorfields  $X, Y$  by setting  $\chi(X, Y)= \chi( ^{(h)}X, ^{(h)}Y)$,  $\chib(X, Y)= \chib( ^{(h)}X, ^{(h)}Y)$.
 Given an horizontal 2-tensor $U$  we define its trace $\tr U$  and anti-trace $\atr U$
\beaa
\tr (U):=\de^{ab}U_{ab}, \qquad \atr U:=\in^{ab} U_{ab}.
\eeaa
Accordingly we  decompose $\chi, \chib$ as follows,
\beaa
\chi_{ab}=\chih_{ab} +\frac 1 2 \de_{ab} \trch+\frac 1 2 \in_{ab}\atrch,\qquad \chib_{ab}=\chibh_{ab} +\frac 1 2 \de_{ab} \trchb+\frac 1 2 \in_{ab}\atrchb,
\eeaa
where $\chih$ and $\chibh$ denote respectively the symmetric traceless part of $\chi$ and $\chib$.

We define the horizontal covariant operator $\nab$ as follows. Given $X, Y\in \O(\MM)$
 \bea
 \nab_X Y&:=&^{(h)}(\D_XY)=\D_XY- \frac 1 2 \chib(X,Y)e_4 -  \frac 1 2 \chi(X,Y) e_3.
 \eea
 In particular, for  all  $X,Y, Z\in \O(\MM)$,
 \beaa
 Z g (X,Y)=g(\nab_Z X, Y)+ g(X, \nab_ZY).
 \eeaa

In the integrable case, $\nab$ coincides with the Levi-Civita connection
 of the metric induced on the integral surfaces of   $\O(\MM)$.  
 Given $X$ horizontal, $\D_4X$ and $\D_3 X$ are in general not horizontal.
 We define $\nab_4 X$ and $\nab_3 X$  to be the horizontal projections
 of the former.  More precisely,
 \beaa
 \nab_4 X&:=&^{(h)}(\D_4 X)=\D_4 X- \frac 1 2 \g(X, \D_4 e_3 ) e_4- \frac 1 2  \g(X, \D_4 e_4)  e_3 ,\\
 \nab_3 X&:=&^{(h)}(\D_3 X)=\D_3 X-   \frac 1 2 \g(X, \D_3e_3) e_3 - \frac 1 2   \g(X, \D_3 e_4 ) e_3. 
 \eeaa
The definition can be easily extended to arbitrary  $  \O_k(\MM) $ tensor-fields  $U$ 
\beaa
 \nab_4U(X_1,\ldots, X_k)&:=&e_4 (U(X_1,\ldots, X_k))- \sum_i U( X_1,\ldots, \nab_4 X_i, \ldots, X_k),\\
  \nab_3 U(X_1,\ldots, X_k)&:=&e_3 (U(X_1,\ldots, X_k)) -\sum_i U( X_1,\ldots, \nab_3 X_i, \ldots, X_k).
 \eeaa 
Finally, for a given horizontal   1-form $\xi$, we  define the frame independent   operators
\bea\lab{eq:defintiondivcurlandnabhot}
\div\xi:=\de^{ab}\nab_b\xi_a,\qquad 
\curl\xi:=\in^{ab}\nab_a\xi_b,\qquad 
(\nab\hot \xi)_{ba}:=\nab_b\xi_a+\nab_a  \xi_b-\de_{ab}( \div \xi).
\eea
  
%%%%%%%%%%%%%%%%%%%%%%%%%

\subsubsection{Ricci and curvature  coefficients}

%%%%%%%%%%%%%%%%%%%%%%%%%

Given a null frame $(e_1, e_2, e_3, e_4)$ we define the following connection coefficients,
 \bea
 \begin{split}
\chib_{ab}&:=\g(\D_ae_3, e_b),\qquad \,\,\,\,\,\,\,\chi_{ab}:=\g(\D_ae_4, e_b),\\
\xib_a&:=\frac 1 2 \g(\D_3 e_3 , e_a),\qquad\,\,\,\,\, \xi_a:=\frac 1 2 \g(\D_4 e_4, e_a),\\
\omb&:=\frac 1 4 \g(\D_3e_3 , e_4),\qquad\,\,\,\,\,\,\, \om:=\frac 1 4 \g(\D_4 e_4, e_3),\qquad \\
\etab_a&:=\frac 1 2\g(\D_4 e_3, e_a),\qquad \quad \eta_a:=\frac 1 2 \g(\D_3 e_4, e_a),\qquad\\
 \ze_a&:=\frac 1 2 \g(\D_{e_a}e_4,  e_3),
 \end{split}
\eea
which account for all the  connection coefficients except $\g(\D_{e_\mu} e_b, e_a)$, $\mu=1,2,3,4$, $a, b=1,2$.

We also define the curvature components 
\bea
\a_{ab}:={\bf R}_{a4b4},\quad \b_a:=\frac 12 {\bf R}_{a434}, \quad \rho:=\frac 1 4 {\bf R}_{3434}, \quad\rhod:=\frac 1 4 \dual{\bf R}_{3434},\quad \bb_a :=\frac 1 2 {\bf R}_{a334}, \quad \aa_{ab}:={\bf R}_{a3b3},
\eea
where $\dual{\bf R}$ denotes the Hodge dual of the curvature tensor ${\bf R}$.

%%%%%%%%%%%%%%%%%%%%%%%%%%%%%%%%%%%%

\subsubsection{The tensorial wave operator}

%%%%%%%%%%%%%%%%%%%%%%%%%%%%%%%%%%%%

In order to define the tensorial wave operator in a covariant way, we first introduce the covariant derivative $\Ddot$ acting on mixed tensors of the type $\T_k (\MM)\otimes   \O_l (\MM)$, i.e., tensors  of the form  $U_{\nu_1\ldots \nu_k,  a_1\ldots a_l}$, 
for which we define
\beaa
\Ddot_\mu U_{\nu_1\ldots \nu_k,  a_1\ldots a_l}&:=& e_\mu(U_{\nu_1\ldots \nu_k,  a_1\ldots a_l}) -U_{\D_\mu e_{\nu_1}\ldots \nu_k,  a_1\ldots a_l}-\ldots- U_{\nu_1\ldots  \D_\mu e_{\nu_k},  a_1\ldots a_l}\\
&-& U_{\nu_1\ldots \nu_k,   ^{(h)}(\D_\mu e_{a_1})\ldots a_l}-  U_{\nu_1\ldots \nu_k,   a_1 \ldots ^{(h)}(\D_\mu e_{a_l})}.
\eeaa

\begin{proposition}
\lab{Proposition:commutehorizderivatives}
For a tensor $\Psi\in \O_1 (\MM)$, we   have  the following formula
 \bea
( \Ddot _\mu\Ddot_\nu -\Ddot_\nu\Ddot _\mu)\Psi_a=\Rdot_{a b  \mu\nu}\Psi^b
 \eea
with an immediate generalization to tensors $\Psi\in \O_l (\MM)$, where, with $(\La_\a)_{\b\ga}= \g(\D_\a e_\ga, e_\b)$,
 \bea
 \lab{eq:DefineRdot}
 \bsplit
 \Rdot_{ab   \mu\nu}&:= {\bf R}_{ab    \mu\nu}+ \frac 1 2  \B_{ab   \mu\nu},\\
  \B_{ab   \mu\nu} &:=  (\La_\mu)_{3a} (\La_\nu)_{b4}+  (\La_\mu)_{4a} (\La_\nu)_{b3}- (\La_\nu)_{3a} (\La_\mu)_{b4}-  (\La_\nu)_{4a} (\La_\mu)_{b3}.
  \end{split}
 \eea 
 \end{proposition}

\begin{proof}
See Proposition 2.1.27 in \cite{GKS22}.
\end{proof}

\begin{proposition}
\lab{proposition:componentsofB}
The components of $\B$   are given   by the following formulas:
\bea
\begin{split}
\B_{ a   b  c 3}&=     -  \trchb  \big( \de_{ca}\eta_b-  \de_{cb} \eta_a\big)  -  \atrchb \big( \in_{ca}  \eta_b -  \in_{cb}  \eta_a\big) \\
&+ 2 \big(- \chibh_{ca}  \eta_b + \chibh_{cb} \eta_a-  \chi_{ca} \xib_b+  \chi_{cb} \xib_a\big),\\
\B_{ a   b  c 4}&=     -  \trch  \big( \de_{ca}\etab_b-  \de_{cb} \etab_a\big)  -  \atrch \big( \in_{ca}  \etab_b -  \in_{cb}  \etab_a\big) \\
&+ 2 \big(- \chih_{ca}  \etab_b + \chih_{cb} \etab_a-  \chib_{ca} \xi_b+  \chib_{cb} \xi_a\big),\\
\B_{ a   b  3 4} &=4\big(-\eta_a \etab_b+\etab_a\eta_b -\xib_a \xi_b+\xi_a \xib_b\big),\\
\B_{abcd} &= \left(- \frac 12  \trch \trchb-\frac 1 2 \atrch \atrchb+\chih \c \chibh\right)\in_{ab}\in_{cd}.
\end{split}
\eea
\end{proposition}

\begin{proof}
See Proposition 2.2.4 in \cite{GKS22}.
\end{proof}

Then, we  define  the wave operator for $\psi \in \sk_k(\mathbb{C})$, $k=0,1,2$, to be, see Definition 2.3.1 in \cite{GKS22},  
 \bea\label{eq:def=squared-2}
 \squared_k\psi:= \g^{\mu\nu} \Ddot_\mu\Ddot_ \nu \psi.
\eea

%%%%%%%%%%%%%%%%%%%%%%%%%

\subsubsection{Commutation formulas}
\lab{sec:generalcommutationformulasrealcase}

%%%%%%%%%%%%%%%%%%%%%%%%%

We start with the following general commutation formulas.
\begin{lemma}
   \lab{LEMMA:COMM-GEN-B}
Let $U_{A}= U_{a_1\ldots a_k} $ be a general $k$-horizontal  tensorfield, and let $\B_{ab\mu\nu}$ be given by \eqref{eq:DefineRdot}. Then:
\begin{enumerate}
\item  We have
\bea
\,[\nab_3, \nab_b] U_A =- \chib_{bc} \nab_c U_A+( \eta_b-\ze_b) \nab_3 U_A +\xib_b \nab_4 U_A +\sum_{i=1}^k\Big(-\in_{a_i c} \dual\bb_b +  \frac 1 2 \B_{a_i c 3b} \Big) U_{a_1\ldots }\,^ c \,_{\ldots a_k}. 
\eea

\item We have
\bea
\,[\nab_4, \nab_b] U_A =- \chi_{bc} \nab_c U_a+( \etab_b+\ze_b) \nab_4 U_a +\xi_b \nab_3 U_a +\sum_{i=1}^k\Big(\in_{a_i c} \dual\b_b +  \frac 1 2 \B_{a_i c 4b} \Big) U_{a_1\ldots }\,^ c \,_{\ldots a_k}.
\eea

\item We have
\bea
\, [\nab_4, \nab_3] U_A= 2(\etab_b-\eta_b ) \nab_b U_A + 2 \om \nab_3 U_A -2\omb \nab_4 U_A+ \sum_{i=1}^k\Big( - \in_{a_i b}\dual \rho+  \frac 1 2 \B_{a_i b 43} \Big) U_{a_1\ldots}\,^b\,_{\ldots a_k}. 
\eea
\end{enumerate}
\end{lemma}

\begin{proof}
See the proof of Lemma 2.2.7 in \cite{GKS22}.
\end{proof}

Using the values of $\B_{ab\mu\nu}$ given by Proposition \ref{proposition:componentsofB}, we specialize  to the case of $\sk_0$, $\sk_1$ and $\sk_2$.

 \begin{lemma}
   \lab{lemma:comm}
   The following commutation formulas hold true:
   \begin{enumerate}
\item Given   $f \in \sk_0$, we have
       \bea\label{eq:comm-nab3-nab4-naba-f-general}
       \begin{split}
        \,[\nab_3, \nab_a] f &=-\frac 1 2 \left(\trchb \nab_a f+\atrchb \dual \nab_a f\right)+(\eta_a-\ze_a) \nab_3 f-\chibh_{ab}\nab_b f  +\xib_a \nab_4 f,\\
         \,[\nab_4, \nab_a] f &=-\frac 1 2 \left(\trch \nab_a f+\atrch \dual \nab_a f\right)+(\etab_a+\ze_a) \nab_4 f-\chih_{ab}\nab_b f  +\xi_a \nab_3 f, \\
         \, [\nab_4, \nab_3] f&= 2(\etab-\eta ) \c \nab f + 2 \om \nab_3 f -2\omb \nab_4 f. 
         \end{split}
       \eea

  \item   Given  $u\in \sk_1$, we have
    \bea\label{commutator-3-a-u-b}\label{commutator-u-in-SS1}
         \bsplit            
\,  [\nab_3,\nab_a] u_b    &=-\frac 1 2 \trchb \big( \nab_a u_b+\eta_b u_a-\de_{ab} \eta \c u \big) -\frac 1 2 \atrchb \big( \dual \nab_a u_b+\eta_b \dual u_a-\in_{ab} \eta\c u\big) \\
&+(\eta-\ze)_a \nab_3 u_b+\err_{3ab}[u],\\
  \err_{3ab}[u] &=-\dual \bb_a\dual u_b+\xib_a\nab_4 u_b-\xib_b \chi_{ac} u_c+\chi_{ab} \,\xib\c u-\chibh_{ac}\nab_c u_b-\eta_b\chibh_{ac}u_c+\chibh_{ab}\eta\c u,
   \end{split}
   \eea
   \bea\label{commutator-4-a-u-b}
   \bsplit
\,  [\nab_4,\nab_a] u_b    &=-\frac 1 2 \trch \big( \nab_a u_b+\etab_b u_a-\de_{ab} \etab \c u \big) -\frac 1 2 \atrch \big( \dual \nab_a u_b+\etab_b \dual u_a-\in_{ab} \etab\c u\big)\\
&+(\etab +\ze)_a \nab_4 u_b +\err_{4ab}[u],\\
   \err_{4ab}[u]&=\dual \b_a\dual u_b+\xi_a\nab_3 u_b-\xi_b \chib_{ac} u_c+\chib_{ab} \,\xi\c u-\chih_{ac}\nab_c u_b-\etab_b\chih_{ac}u_c+\chih_{ab}\etab\c u, 
      \end{split}
   \eea
   \bea
   \bsplit
 \, [\nab_4, \nab_3] u_a&=2 \om \nab_3 u_a -2\omb \nab_4 u_a+ 2(\etab_b-\eta_b ) \nab_b u_a +2(\etab \c u ) \eta_{a} -2 (\eta \c u )\etab_{a}\\
 &  -2 \dual \rho \dual u_a +\err_{43a}[u],\\
 \err_{43a}[u]&= 2 \big( \xib_{a}  \xi_b- \xi_{a}  \xib_b )u^b.
\end{split}
\eea

\item  Given  $u\in \sk_2$, we have 
    \bea\label{commutator-u-in-SS2}\label{commutator-3-a-u-bc}
         \bsplit            
\,  [\nab_3,\nab_a] u_{bc}    &=-\frac 1 2 \trchb\, (\nab_a u_{bc}+\eta_bu_{ac}+\eta_c u_{ab}-\de_{a b}(\eta \c u)_c-\de_{a c}(\eta \c u)_b )\\
&-\frac 1 2 \atrchb\, (\dual \nab_a u_{bc} +\eta_b\dual u_{ac}+\eta_c\dual u_{ab}- \in_{a b}(\eta \c u)_c- \in_{a c}(\eta \c u)_b )\\
&+(\eta_a-\ze_a)\nab_3 u_{bc}+\err_{3abc}[u],\\
\err_{3abc}[u]&= -2\dual \bb_a \dual u_{bc}+\xib_a \nab_4 u_{bc} -\xib_b\chi_{ad}u_{dc} -\xib_c\chi_{ad}u_{bd}+\chi_{ab}\xib_d u_{dc} \\
&+\chi_{ac}\xib_d u_{bd}-\chibh_{ad} \nab_d u_{bc} -\eta_b\chibh_{ad}u_{dc} - \eta_c\chibh_{ad}u_{bd}+\chibh_{ab}\eta_du_{dc} +\chibh_{ac}\eta_du_{bd},
   \end{split}
   \eea
   \bea\label{commutator-4-a-u-bc}
   \bsplit
\,  [\nab_4,\nab_a] u_{bc}    &=-\frac 1 2 \trch\, (\nab_a u_{bc}+\etab_bu_{ac}+\etab_c u_{ab}-\de_{a b}(\etab \c u)_c-\de_{a c}(\etab \c u)_b )\\
&-\frac 1 2 \atrch\, (\dual \nab_a u_{bc} +\etab_b\dual u_{ac}+\etab_c\dual u_{ab}- \in_{a b}(\etab \c u)_c- \in_{a c}(\etab \c u)_b )\\
&+(\etab_a+\ze_a)\nab_4 u_{bc}+\err_{4abc}[u],\\
\err_{4abc}[u]&= 2\dual \b_a \dual u_{bc}+\xi_a \nab_3 u_{bc} -\xi_b\chib_{ad}u_{dc} -\xi_c\chib_{ad}u_{bd}+\chib_{ab}\xi_d u_{dc} +\chib_{ac}\xi_d u_{bd}\\
& -\chih_{ad} \nab_d u_{bc} -\etab_b\chih_{ad}u_{dc} - \etab_c\chih_{ad}u_{bd}+\chih_{ab}\etab_du_{dc} +\chih_{ac}\etab_du_{bd}, 
     \end{split}
   \eea
   \bea\label{commutator-4-3-u-bc}
   \bsplit
   \, [\nab_4, \nab_3] u_{ab} &=2 \om \nab_3 u_{ab} -2\omb \nab_4 u_{ab} + 2(\etab_c-\eta_c ) \nab_c u_{ab} + 4 \eta \hot (\etab \c u)  \\
   &-4 \etab \hot (\eta \c u)-4 \dual \rho \dual u_{ab}+\err_{43ab}[u],\\
\err_{43ab}[u]&= 2 \big( \xib_{a}  \xi_c- \xi_{a}  \xib_c )u^c\,_{b}+2 \big( \xib_{b}  \xi_c- \xi_{b}  \xib_c )u_{a} \,^c.
   \end{split}
\eea
       \end{enumerate}
 \end{lemma}

\begin{proof}
See the proof of Lemma 2.2.8 in \cite{GKS22}.
\end{proof}

%%%%%%%%%%%%%%%%%%%%%%%%%%%%%

\subsubsection{Horizontal Lie derivatives}
\label{subsect:horizontalliederivatives}

%%%%%%%%%%%%%%%%%%%%%%%%%%%%%

Recall that the Lie derivative of a $k$-covariant tensor $U$ relative to a vectorfield  $X$ is given by
\beaa
\Lie_X{U}\big(e_{\a_1}, \ldots , e_{\a_k}\big) = X\big(U_{\a_1\ldots\a_k}\big) -  U\big(\Lie_Xe_{\a_1}, \ldots,e_{\a_k}\big) - U\big(e_{\a_1}, \ldots, \Lie_Xe_{\a_k}\big),
\eeaa
where $\Lie_X Y=[X, Y]$. We define horizontal Lie derivatives as follows, see Definition 2.2.12 in \cite{GKS22}.

\begin{definition}[Horizontal Lie derivatives]\label{definition:hor-Lie-derivative}
Given   vectorfields  $X$, $Y$,  the horizontal Lie  derivative  $\Lieb_XY$ is given by
 \beaa
 \Lieb_X Y :=\Lie_X Y+ \frac 1 2 \g(\Lie_XY, e_3) e_4+  \frac 1 2 \g(\Lie_XY, e_4) e_3.
 \eeaa
 Given  a horizontal covariant k-tensor $U$,  the horizontal  Lie derivative $\Lieb_X U $ is defined   to be the projection of $\Lie_X U$ to the  horizontal space, i.e.,
  \beaa
 \Lieb_X U\big(e_{a_1}, \ldots,e_{a_k}\big) := X\big(U_{a_1\ldots a_k}\big)- U\big(\Lieb_Xe_{a_1},\ldots, e_{a_k} \big)-\ldots -U\big(e_{a_1},\ldots,  \Lieb_Xe_{a_k} \big).
  \eeaa 

Also, given a mixed tensor $U$ of the type $\T_k (\MM)\otimes   \O_l (\MM)$, we define the general horizontal derivative $\Lied_XU$ as follows 
\beaa
&& \Lied_XU\big(e_{\a_1},\ldots, e_{\a_k},  e_{a_1},\ldots, e_{a_l}\big) \\
&:=& X\big(U_{\a_1\ldots \a_k,  a_1\ldots a_l}\big) -U\big(\Lie_Xe_{\a_1},\ldots, e_{\a_k},  e_{a_1}, \ldots, e_{a_l}\big) -\ldots  - U\big(e_{\a_1},\ldots.  \Lie_Xe_{\a_k},  e_{a_1},\ldots, e_{a_l}\big)\\
&& -U\big(e_{\a_1}, \ldots, e_{\a_k},  \Lieb_Xe_{a_1},\ldots, e_{a_l}\big) - \ldots -  U\big(e_{\a_1},\ldots, e_{\a_k},   e_{a_1}, \ldots, \Lieb_X e_{a_l}\big).
\eeaa
\end{definition}

%%%%%%%%%%%%%%%%%%%%%%%%%%%%%%%%%%%%%%%%%%%%    
    
\subsection{Main equations using complex notations}\label{section:complex-notations}

%%%%%%%%%%%%%%%%%%%%%%%%%%%%%%%%%%%%%%%%%%%%

%%%%%%%%%%%%%%%%%%%%%%%%

\subsubsection{Complex notations}
\lab{sec:complexnotationsRicciandcurvature}

%%%%%%%%%%%%%%%%%%%%%%%%

Recall Definition \ref{definition-SS-real} of the set of real horizontal tensors $\sk_k=\sk_k(\MM, \mathbb{R})$ on $\MM$ for $k=0,1,2$. We now define the corresponding complexified versions.

\begin{definition} 
\lab{def:skC:horizontaltensors}
We denote by $\sk_k(\mathbb{C})$, $k=0,1,2$, the following set of  horizontal tensors on $\MM$: 
\beaa
 a+ i b \in \sk_0(\mathbb{C})\,\,\,\textrm{if}\,\,\, (a, b) \in \sk_0, \quad F= f+ i \dual f  \in \sk_1(\mathbb{C})\,\,\,\textrm{if}\,\,\, f \in \sk_1, \quad U=u + i \dual u \in \sk_2(\mathbb{C})\,\,\,\textrm{if}\,\,\, u \in \sk_2,
\eeaa
where $F\in\sk_1(\mathbb{C})$ and $U\in\sk_2(\mathbb{C})$ are anti-self dual, i.e., $\dual F=-iF$ and $\dual U=-iU$. 
\end{definition}

\begin{definition}\lab{def:complexRicciandcurvaturecoefficients}
We define the following complexified curvature components 
\beaa
A:=\a+i\dual\a, \quad B:=\b+i\dual\b, \quad P:=\rho+i\dual\rho,\quad \Bb:=\bb+i\dual\bb, \quad \Ab:=\aa+i\dual\aa,
\eeaa      
with $A, \Ab\in\sk_2(\mathbb{C})$, $B, \Bb\in\sk_1(\mathbb{C})$, $P\in\sk_0(\mathbb{C})$, and the following complexified Ricci coefficients     
\beaa
&& X:=\chi+i\dual\chi, \quad \Xb:=\chib+i\dual\chib, \quad H:=\eta+i\dual \eta, \quad \Hb:=\etab+i\dual \etab,  \\ 
&& Z:=\ze+i\dual\ze, \quad \Xi:=\xi+i\dual\xi, \quad \Xib:=\xib+i\dual\xib,
\eeaa    
with $\widehat{X}, \Xbh\in\sk_2(\mathbb{C})$, $H, \Hb, Z, \Xi, \Xib\in\sk_1(\mathbb{C})$, where $\widehat{X}, \Xbh$, as well as $\tr X, \tr\Xb$ are given by
\beaa
\tr X := \trch-i\atrch, \quad \widehat{X}:=\chih+i\dual\chih, \quad \tr\Xb:=\trchb -i\atrchb, \quad \Xbh:=\chibh+i\dual\chibh.
\eeaa
\end{definition}

\begin{definition}
We define derivatives of complex quantities as follows
\begin{itemize}
\item For two scalar functions $a$ and $b$, we define
\beaa
\DD(a+ib) &:=& (\nabla+i\dual\nabla)(a+ib).
\eeaa

\item For a 1-form $f$, we define
\beaa
\ov{\DD}\c(f+i\dual f) &:=& (\nabla-i\dual\nabla)\c(f+i\dual f)
\eeaa
and  
\beaa
\DD\hot(f+i\dual f) &:=& (\nabla+i\dual\nabla)\hot(f+i\dual f).
\eeaa

\item For a symmetric traceless 2-tensor $u$, we define
\beaa
\ov{\DD}\c(u+i\dual u) &:=& (\nabla-i\dual\nabla)\c(u+i\dual u).
\eeaa
\end{itemize}
\end{definition}

%%%%%%%%%%%%%%%%%%%%%%%%%%%%%

\subsubsection{Main equations  in complex form}

%%%%%%%%%%%%%%%%%%%%%%%%%%%%%

The complex notations allow us to rewrite the Ricci equations in  a more compact  form. 
\begin{proposition}
\label{prop-nullstr:complex} {We have}
\beaa
\nab_3\tr\Xb +\frac{1}{2}(\tr\Xb)^2+2\omb\,\tr\Xb &=& \DD\c\ov{\Xib}+\Xib\c\ov{\Hb}+\ov{\Xib}\c(H-2Z)-\frac{1}{2}\Xbh\c\ov{\Xbh},\\
\nab_3\Xbh+\Re(\tr\Xb) \Xbh+ 2\omb\,\Xbh&=&\frac 1 2  \DD\hot \Xib+\frac 1 2   \Xib\hot(H+\Hb-2Z)-\Ab,
\eeaa
\beaa
\nab_3\tr X +\frac{1}{2}\tr\Xb\tr X-2\omb\tr X &=& \DD\c\ov{H}+H\c\ov{H}+2P+\Xib\c\ov{\Xi}-\frac{1}{2}\Xbh\c\ov{\Xh},\\
\nab_3\widehat{X} +\frac{1}{2}\tr\Xb\, \widehat{X} -2\omb\widehat{X} &=&\frac 1 2  \DD\hot H  +\frac 1 2 H\hot H -\frac{1}{2}\ov{\tr X} \widehat{\Xb}+\frac{1}{4}\Xib\hot\Xi,
\eeaa
\beaa
\nab_4\tr\Xb +\frac{1}{2}\tr X\tr\Xb -2\om\tr\Xb &=& \DD\c\ov{\Hb}+\Hb\c\ov{\Hb}+2\ov{P}+\Xi\c\ov{\Xib}-\frac{1}{2}\Xh\c\ov{\Xbh},\\
\nab_4\widehat{\Xb} +\frac{1}{2}\tr X\, \widehat{\Xb} -2\om\widehat{\Xb} &=&\frac  12  \DD\hot\Hb  +\frac 1 2 \Hb\hot\Hb -\frac{1}{2}\ov{\tr\Xb} \widehat{X}+\frac{1}{4}\Xi\hot\Xib,
\eeaa
\beaa
\nab_4\tr X +\frac{1}{2}(\tr X)^2+2\om\tr X &=& \DD\c\ov{\Xi}+\Xi\c\ov{H}+\ov{\Xi}\c({\Hb}+2Z)-\frac{1}{2}\Xh\c\ov{\Xh},\\
\nab_4\Xh+\Re(\tr X)\Xh+ 2\om\Xh&=&\frac 1 2  \DD\hot \Xi+\frac 12   \Xi\hot(\Hb+H+2Z)-A.
\eeaa
Also,
\beaa
\nab_3Z +\frac{1}{2}\tr\Xb(Z+H)-2\omb(Z-H) &=& -2\DD\omb -\frac{1}{2}\widehat{\Xb}\c(\ov{Z}+\ov{H})+\frac{1}{2}\tr X\Xib+2\om\Xib -\Bb+\frac{1}{2}\ov{\Xib}\c\Xh,\\
\nab_4Z +\frac{1}{2}\tr X(Z-\Hb)-2\om(Z+\Hb) &=& 2\DD\om +\frac{1}{2}\widehat{X}\c(-\ov{Z}+\ov{\Hb})-\frac{1}{2}\tr\Xb\Xi-2\omb\Xi -B-\frac{1}{2}\ov{\Xi}\c\Xbh,\\
\nab_3\Hb -\nab_4\Xib &=&  -\frac{1}{2}\ov{\tr\Xb}(\Hb-H) -\frac{1}{2}\Xbh\c(\ov{\Hb}-\ov{H}) -4\om\Xib+\Bb,\\
\nab_4H -\nab_3\Xi &=&  -\frac{1}{2}\ov{\tr X}(H-\Hb) -\frac{1}{2}\Xh\c(\ov{H}-\ov{\Hb}) -4\omb\Xi-B,
\eeaa
and
\beaa
\nab_3\om+\nab_4\omb -4\om\omb -\xi\c \xib -(\eta-\etab)\c\ze +\eta\c\etab&=&   \rho.
\eeaa
Also,
\beaa
\frac{1}{2}\ov{\DD}\c\Xh +\frac{1}{2}\Xh\c\ov{Z} &=& \frac{1}{2}\DD\ov{\tr X}+\frac{1}{2}\ov{\tr X}Z-i\Im(\tr X)H-i\Im(\tr \Xb)\Xi-B,\\
\frac{1}{2}\ov{\DD}\c\Xbh -\frac{1}{2}\Xbh\c\ov{Z} &=& \frac{1}{2}\DD\ov{\tr\Xb}-\frac{1}{2}\ov{\tr\Xb}Z-i\Im(\tr\Xb)\Hb-i\Im(\tr X)\Xib+\Bb,
\eeaa
and,
\beaa
\curl\ze&=&-\frac 1 2 \chih\wedge\chibh   +\frac 1 4 \big(  \trch\atrchb-\trchb\atrch   \big)+\om \atrchb -\omb\atrch+\dual \rho.
\eeaa
\end{proposition}

The complex notations allow us to rewrite the Bianchi identities as follows.  
  \begin{proposition}\label{prop:bianchi:complex} 
    We have,
 \beaa
 \nab_3A -\frac 1 2 \DD\hot B &=& -\frac{1}{2}\tr\Xb A+4\omb A +\frac 1 2 (Z+4H)\hot B -3\ov{P}\Xh,\\
\nab_4B -\frac{1}{2}\ov{\DD}\c A &=& -2\ov{\tr X} B -2\om B +\frac{1}{2}A\c  (\ov{2Z +\Hb})+3\ov{P} \,\Xi,\\
\nab_3B-\DD\ov{P} &=& -\tr\Xb B+2\omb B+\ov{\Bb}\c \Xh+3\ov{P}H +\frac{1}{2}A\c\ov{\Xib},\\
\nab_4P -\frac{1}{2}\DD\c \ov{B} &=& -\frac{3}{2}\tr X P +\frac{1}{2}(2\Hb+Z)\c\ov{B} -\ov{\Xi}\c\Bb -\frac{1}{4}\Xbh\c \ov{A}, \\
\nab_3P +\frac{1}{2}\ov{\DD}\c\Bb &=& -\frac{3}{2}\ov{\tr\Xb} P -\frac{1}{2}(\ov{2H-Z})\c\Bb +\Xib\c \ov{B} -\frac{1}{4}\ov{\Xh}\c\Ab, \\
\nab_4\Bb+\DD P &=& -\tr X\Bb+2\om\Bb+\ov{B}\c \Xbh-3P\Hb -\frac{1}{2}\Ab\c\ov{\Xi},\\
\nab_3\Bb +\frac{1}{2}\ov{\DD}\c\Ab &=& -2\ov{\tr\Xb}\,\Bb -2\omb\,\Bb -\frac{1}{2}\Ab\c (\ov{-2Z +H})-3P \,\Xib,\\
\nab_4\Ab +\frac{1}{2}\DD\hot\Bb &=& -\frac{1}{2}\tr X \Ab+4\om\Ab +\frac{1}{2}(Z-4\Hb)\hot \Bb -3P\Xbh.
\eeaa
    \end{proposition} 
    
    \begin{proof} 
    See Proposition 2.4.12 in \cite{GKS22}.
    \end{proof}

%%%%%%%%%%%%%%%%%%%%%%%%%%%%%%%%%%%%%

\subsubsection{Main complex equations using conformal derivatives}

%%%%%%%%%%%%%%%%%%%%%%%%%%%%%%%%%%%%%

Following Sections 2.2.9 and 2.4.3 in \cite{GKS22}, we introduce conformal derivatives. Consider  frame transformations of the form
\beaa
e_3'=\la^{-1} e_3, \qquad  e'_4 = \la e_4 , \qquad e_a'= e_a.
\eeaa
Note that  under   the   above  mentioned  frame transformation we have
\beaa
 \Xb'&=&\la^{-1} \Xb, \quad X'=\la X, \quad \Xi'= \la^2\Xi, \quad   H'=H, \quad \Hb'=\Hb,  \quad \Xib'=\la^{-2}\Xib,\\
   A'&=&\la^2A,\quad B'=\la B,    \quad P'=P,    \quad  \Bb'=\la^{-1} \Bb,\quad  \Ab'=\la^{-2} \Ab,
   \eeaa
   and
   \beaa
 \omb'&=& \la^{-1}\left(\omb +\frac{1}{2} e_3(\log \la)\right), \quad \om'= \la\left(\om -\frac{1}{2} e_4(\log \la)\right), \quad
 Z'= Z - \DD(\log \la).
\eeaa

\begin{definition}[$s$-conformally invariants]
\lab{def:sconformalinvariants}
We say that  a horizontal tensor $f$ is $s$-conformally invariant  if, under the  conformal  frame transformation above,  it changes as $f'=\la^s f $. 
\end{definition}

\begin{remark}
If $f$ $s$-conformal invariant, then  $\nab_3 f, \nab_4 f, \nab_a f$ are not conformal invariant.
\end{remark} 

 We correct the lack of being conformal invariant by making the following  definition.  

\begin{lemma}\label{lemma:definition-conformal-derivatives}
If $f$ is $s$-conformal invariant, then: 
 \begin{enumerate}
 \item $\nabc_3 f:= \nab_3f-2 s \omb f$ is $(s-1)$-conformally invariant.
 
 \item $\nabc_4 f:= \nab_4f+2 s \om f$ is $(s+1)$-conformally invariant.
 
 \item $\nabc_a f:= \nab_af+ s \ze_a f$ is $s$-conformally invariant. 
  \end{enumerate}
\end{lemma}

\begin{proof}
Immediate verification.
\end{proof}

\begin{remark} 
Note that $s$ is precisely what   in \cite{Ch-Kl} is called       the  signature of the tensor. In GHP formalism \cite{GHP}, the signature is related to the boost weights of the complex scalars.
\end{remark}

\begin{definition}
 We define the following conformal angular derivatives in the complex notation:
 \begin{itemize}
\item For $a+i b \in \sk_0(\mathbb{C}) $  we define
\beaa
\DDc(a+ib) &:=& \big(\nabc +i\dual \nabc\big)(a+ib).
\eeaa

\item For  $f+i \dual f \in\sk_1(\mathbb{C}) $ we define
\beaa
\DDc(f+i\dual f) &:=& \big(\nabc+i\dual \nabc\big) \c (f+i\dual f),
\\
\DDc \hot(f+i\dual f) &:=& (\nabc+i\dual\nabc )\hot(f+i\dual f).
\eeaa
\item For $u+ i \dual u \in \sk_2(\mathbb{C})$ we define
\beaa
\DDc \c(u+i\dual u) &:=& \big(\nabc +i\dual\nabc \big)\c(u+i\dual u).
\eeaa
\item In all the above cases we set
\beaa
\DDbc&:=&\nabc-i\nabc.
\eeaa
\end{itemize}
\end{definition}

Using these definitions we rewrite the main equations as follows.
\begin{proposition}
\label{prop-nullstr:complex-conf}
We have
\beaa
\nabc_3\tr\Xb +\frac{1}{2}(\tr\Xb)^2 &=& \DDc\c\ov{\Xib}+\Xib\c\ov{\Hb}+\ov{\Xib}\c H-\frac{1}{2}\Xbh\c\ov{\Xbh},\\
\nabc_3\Xbh+\Re(\tr\Xb) \Xbh&=&\frac 1 2  \DDc\hot \Xib+  \frac 1 2  \Xib\hot(H+\Hb)-\Ab,
\eeaa
\beaa
\nabc_3\tr X +\frac{1}{2}\tr\Xb\tr X &=& \DDc\c\ov{H}+H\c\ov{H}+2P+\Xib\c\ov{\Xi}-\frac{1}{2}\Xbh\c\ov{\Xh},\\
\nabc_3\widehat{X} +\frac{1}{2}\tr\Xb\, \widehat{X} &=&\frac 1 2 \DDc\hot H  +\frac 1 2 H\hot H -\frac{1}{2}\ov{\tr X} \widehat{\Xb}+\frac 1 4 \Xib\hot\Xi,
\eeaa
\beaa
\nabc_4\tr\Xb +\frac{1}{2}\tr X\tr\Xb &=& \DDc\c\ov{\Hb}+\Hb\c\ov{\Hb}+2\ov{P}+\Xi\c\ov{\Xib}-\frac{1}{2}\Xh\c\ov{\Xbh},\\
\nabc_4\widehat{\Xb} +\frac{1}{2}\tr X\, \widehat{\Xb} &=&\frac 1 2  \DDc\hot\Hb  +\frac 1 2 \Hb\hot\Hb -\frac{1}{2}\ov{\tr\Xb} \widehat{X}+\frac 1 4 \Xi\hot\Xib,
\eeaa
\beaa
\nabc_4\tr X +\frac{1}{2}(\tr X)^2 &=& \DDc\c\ov{\Xi}+\Xi\c\ov{H}+\ov{\Xi}\c\Hb-\frac{1}{2}\Xh\c\ov{\Xh},\\
\nabc_4\Xh+\Re(\tr X)\Xh&=&\frac 1 2  \DDc\hot \Xi+\frac 1 2   \Xi\hot(\Hb+H)-A,
\eeaa
\beaa
\nabc_3\Hb -\nabc_4\Xib &=&  -\frac{1}{2}\ov{\tr\Xb}(\Hb-H) -\frac{1}{2}\Xbh\c(\ov{\Hb}-\ov{H}) +\Bb,\\
\nabc_4H -\nabc_3\Xi &=&  -\frac{1}{2}\ov{\tr X}(H-\Hb) -\frac{1}{2}\Xh\c(\ov{H}-\ov{\Hb}) -B.
\eeaa
Also,
\beaa
\frac{1}{2}\ov{\DDc}\c\Xh &=& \frac{1}{2}\DDc\ov{\tr X}-i\Im(\tr X)H-i\Im(\tr \Xb)\Xi-B,\\
\frac{1}{2}\ov{\DDc}\c\Xbh &=& \frac{1}{2}\DDc\ov{\tr\Xb}-i\Im(\tr\Xb)\Hb-i\Im(\tr X)\Xib+\Bb.
\eeaa
\end{proposition}

\begin{proof}
See Proposition 2.4.14 in \cite{GKS22}.
\end{proof}
    
  \begin{proposition}\label{prop:bianchi:complex-conf} 
    We have
 \beaa
 \nabc_3A -\frac 1 2 \DDc\hot B &=& -\frac{1}{2}\tr\Xb A + 2 H   \hot B -3\ov{P}\Xh,\\
\nabc_4B -\frac{1}{2} \DDbc \c A &=& -2\ov{\tr X} B +\frac{1}{2}A\c \ov{\Hb}+3\ov{P} \,\Xi,\\
\nabc_3B-\DDc\ov{P} &=& -\tr\Xb B+\ov{\Bb}\c \Xh+3\ov{P}H +\frac{1}{2}A\c\ov{\Xib},\\
\nabc_4P -\frac{1}{2}\DDc\c \ov{B} &=& -\frac{3}{2}\tr X P + \Hb \c\ov{B} -\ov{\Xi}\c\Bb -\frac{1}{4}\Xbh\c \ov{A}, 
\eeaa
\beaa
\nabc_3P +\frac{1}{2}\DDbc \c\Bb &=& -\frac{3}{2}\ov{\tr\Xb} P - \ov{H} \c\Bb +\Xib\c \ov{B} -\frac{1}{4}\ov{\Xh}\c\Ab, \\
\nabc_4\Bb+\DDc P &=& -\tr X\Bb+\ov{B}\c \Xbh-3P\Hb -\frac{1}{2}\Ab\c\ov{\Xi},\\
\nabc_3\Bb +\frac{1}{2}\DDbc \c\Ab &=& -2\ov{\tr\Xb}\,\Bb  -\frac 1 2  \Ab\c \ov{H}-3P \,\Xib,\\
\nabc_4\Ab +\frac 1 2 \DDc\hot\Bb &=& -\frac{1}{2}\tr X \Ab - 2 \Hb\hot \Bb -3P\Xbh.
\eeaa
    \end{proposition} 
    
\begin{proof}
See Proposition 2.4.15 in \cite{GKS22}.
\end{proof}

%%%%%%%%%%%%%%%%%%%%%%%%%%%%%%%

\subsection{Kerr values}

%%%%%%%%%%%%%%%%%%%%%%%%%%%%%%%

%%%%%%%%%%%%%%%%%%%%%%%%%%%%%%%

\subsubsection{Normalized coordinates in Kerr spacetimes}
\label{subsect:normalizedcoords}

%%%%%%%%%%%%%%%%%%%%%%%%%%%%%%%

The Kerr metric in Boyer--Lindquist coordinates $(t,r,\th,\phi)$ is given by
\begin{align}
\gam={}& \g_{tt}dt^2 +\g_{rr}dr^2+(\g_{t\phi}+\g_{\phi t})dtd\phi +\g_{\phi\phi}d\phi^2 +\g_{\th\th}d\th^2,
\end{align}
where
\begin{equation}
\begin{split}
\g_{tt}={}&-\frac{\Delta-a^2\sin^2\theta}{|q|^2}, \quad \g_{t\phi}={}\g_{\phi t}=-\frac{2amr\sin\theta}{|q|^2}, \quad \g_{rr}={}\frac{|q|^2}{\Delta},\\
\g_{\phi\phi}={}&\frac{(r^2+a^2)^2-a^2\sin^2\theta\Delta}{|q|^2}\sin^2\theta, \quad \g_{\th\th}={}|q|^2,
\end{split}
\end{equation}
with 
\bea
\Delta:=r^2-2mr+a^2,\qquad |q|^2:=r^2+a^2\cos^2\th.
\eea
{In particular, $\partial_{t}$ and $\partial_{\phi}$ are Killing vectorfields and the larger root 
\begin{align}
r_+:=m + \sqrt{m^2 -a^2}
\end{align}
of $\Delta=\De(r)$ corresponds to the location of the event horizon.

It is well-known that the metric is singular on the event horizon in both the Boyer--Lindquist and the tortoise coordinates. To extend the Kerr metric beyond the future event horizon, we define the ingoing Eddington--Finkelstein coordinates $(v_+, r,\th,\phi_+)$ by 
\bea\lab{eq:definitionofingoingEFcoordiantesvplusandphiplus}
dv_+=dt+\frac{r^2+a^2}{\De}dr, \quad d\phi_+=d\phi+\frac{a}{\Delta}dr \,\,\, \text{mod } 2\pi.
\eea
The Kerr metric in this coordinate system is 
\begin{align}\lab{eq:KerrmetriciningoingEddigtonFinkelstein}
\gam={}&-\bigg(1-\frac{2mr}{|q|^2}\bigg)dv_+^2 +2dr dv_+ -\frac{4amr\sin^2\th}{|q|^2}dv_+d\phi_+ -2a\sin^2\th dr d\phi_+\nn\\
&+|q|^2 d\th^2+\frac{(r^2+a^2)^2-a^2\sin^2\theta\Delta}{|q|^2}\sin^2\th d\phi_+^2.
\end{align}

In the following lemma, we introduce coordinates systems, referred to as \textit{normalized coordinates}, used throughout the paper. 

\begin{lemma}[Normalized coordinates]
\label{lem:specificchoice:normalizedcoord}
We fix constants $\dhor$ and $\dbl$ such that 
$$0<\dhor\ll \dbl\ll 1-\frac{|a|}{m}.$$ 
There exists a choice of smooth functions $\tmod=\tmod(r)$ and $\phimod=\phimod(r)$ such that the coordinate systems $(\tt, r, x^1_0, x^2_0)$ and $(\tt, r, x^1_p, x^2_p)$, defined respectively on $\th\neq 0, \pi$ and $\th\neq\frac{\pi}{2}$, with 
\bea\lab{eq:definitionof:specificchoice:normalizedcoord}
\tau=v_+-\tmod, \quad \tphi=\phi_+ -  \phimod, \quad x^1_0=\th, \quad x^2_0=\tphi, \quad x^1_p=\sin\th\cos\tphi, \quad x^2_p=\sin\th\sin\tphi,
\eea
satisfy the following properties:
\begin{enumerate}
\item defining the causal spacetime region $\MM$ and corresponding spacelike boundary $\AA$ by 
\beaa
\bsplit
{}\qquad\MM&:=\big(\{(\tt, r, x^1_0, x^2_0),\,\, \th\neq 0, \pi\}\cup\{(\tt, r, x^1_p, x^2_p),\,\, \th\neq \pi/2\}\big)\cap\{r\geq r_+(1-\dhor)\}, \\ 
{}\qquad\AA&:=\pr\MM=\big(\{(\tt, r, x^1_0, x^2_0),\,\, \th\neq 0, \pi\}\cup\{(\tt, r, x^1_p, x^2_p),\,\, \th\neq \pi/2\}\big)\cap\{r=r_+(1-\dhor)\},
\end{split}
\eeaa
$\MM$ is covered by $(\tt, r, x^1_0, x^2_0)$ and $(\tt, r, x^1_p, x^2_p)$ with the metric components and  inverse metric components being smooth on their respective coordinate patch, 

\item $(\tt, r, x^1_0, x^2_0)$ coincides with Boyer-Lindquist coordinates\footnote{In particular, we have
\beaa
\tmod'(r)=\frac{r^2+a^2}{\De}, \qquad \phimod'(r)=\frac{a}{\De}\quad\textrm{on}\quad r\in [r_+(1+2\dbl), 12m].
\eeaa}
 in $r\in [r_+(1+2\dbl),  12m]$, 

\item for $r\notin (r_+(1+\dbl), 13m)$, we choose 
\beaa
\begin{split}
\tmod'(r) &=\frac{m^2}{r^2}, \qquad \phimod'(r)=0\quad\textrm{on}\quad r\leq r_+(1+\dbl),\\
\tmod'(r)&=\frac{2(r^2+a^2)}{\De} -\frac{m^2}{r^2}, \qquad \phimod'(r)=\frac{2a}{\De} \quad\textrm{on}\quad r\geq 13m,
\end{split}
\eeaa

\item the level sets of $\tt$ in $\MM$ are globally spacelike, transverse to the future event horizon $\HH_+$ and the spacelike boundary $\AA$, and asymptotically null to future null infinity $\II_+$.
\end{enumerate}

Furthermore, the nontrivial inverse metric components in the coordinate system $(\tt, r, \th, \tphi)$ are
\begin{align}
\label{eq:inverse:hypercoord}
\gam^{\tt \tt}={}&\frac{a^2\sin^2\th}{|q|^2} -\frac{2(r^2+a^2)}{|q|^2}\tmod'+\frac{\Delta}{|q|^2}(\tmod')^2, \qquad \gam^{rr}=\frac{\Delta}{|q|^2},\nn\\ 
\gam^{\tt r}={}&\gam^{r \tt}=\frac{r^2+a^2}{|q|^2}\left(1-\frac{\De}{r^2+a^2}\tmod'\right), \qquad \gam^{r\tphi} =\gam^{\tphi r} =\frac{a}{|q|^2}-\frac{\Delta}{|q|^2}\phimod',\nn\\
\gam^{\tt\tphi}={}& \gam^{\tphi \tt}= \frac{a}{|q|^2} (1-\tmod')-\phimod'\frac{r^2+a^2}{|q|^2}\left(1-\frac{\De}{r^2+a^2}\tmod'\right),\nn\\
\gam^{\th\th}={}&\frac{1}{|q|^2}, \qquad \gam^{\tphi\tphi}=\frac{1}{|q|^2\sin^2\th}-\frac{2a}{|q|^2}\phimod' +\frac{\Delta}{|q|^2}(\phimod')^2.
\end{align}
\end{lemma}

\begin{remark}\lab{rmk:phimodprimeisproportionaltoa!!}
Additionally, we may choose $\phimod$ such that 
\beaa
\phimod'(r)=a\phi_{\textrm{mod},0}'(r), \qquad \phi_{\textrm{mod},0}'(r)\geq 0\quad \forall r\in(r_+(1-\dhor), +\infty),
\eeaa
so that $\phimod'(r)$ has the same sign as $a$. From now on, we will assume that our choice of $\phimod$ satisfies this property. In view of \eqref{eq:inverse:hypercoord}, it implies that the inverse metric coefficients $\gam^{\a\b}$ in the normalized coordinates system $(\tau, r, x^1, x^2)$ are invariant under the change $(a, \tphi)\to (-a, -\tphi)$. 
\end{remark}

\begin{proof}
See the proof of Lemma 2.1 in \cite{MaSz24}.
\end{proof}

%%%%%%%%%%%%%%%%%%%%%%%%%%%%%%%

\subsubsection{Principal null pair in Kerr}
\label{subsect:principalnullpairinKerr}

%%%%%%%%%%%%%%%%%%%%%%%%%%%%%%%

We consider the principal null pair of Kerr which is regular across the future event horizon, i.e., in Boyer-Lindquist coordinates, 
\bea
\lab{def:e3e4inKerr}
 e_4 = \frac{r^2+a^2}{|q|^2} \pr_t +\frac{\De}{|q|^2} \pr_r +\frac{a}{|q|^2} \pr_{\phi}, \qquad 
 e_3=\frac{r^2+a^2}{\De} \pr_t -\pr_r +\frac{a}{\De} \pr_{\phi}.
\eea
Also, we consider its associated horizontal bundle $\{e_3, e_4\}^\perp$, which, for $\th\neq 0, \pi$, is spanned by 
\bea
\lab{def:e1e2inKerr}
 e_1=\frac{1}{|q|}\pr_\th,\quad e_2=\frac{a\sin\th}{|q|}\pr_t+\frac{1}{|q|\sin\th}\pr_\phi,
\eea
and we define the complex-valued scalar $q$ and the complex horizontal  $1$-forms $\Jk$ and $\Jk_\pm$ as 
\bea\lab{eq:def:Jkandq}
\bsplit
q=& r+ia\cos\th, \qquad \Jk=j+i\dual j, \qquad \Jk_\pm=j_\pm+i\dual j_\pm, \qquad j_1=0, \quad j_2=\frac{\sin\th}{|q|}, \\
(j_+)_1=&\frac{1}{|q|} \cos\th\cos\tphi, \,\,\,\, (j_+)_2=-\frac{1}{|q|} \sin\tphi, \,\,\,\, (j_-)_1 =\frac{1}{|q|} \cos\th\sin\tphi, \,\,\,\, (j_-)_2=\frac{1}{|q|}  \cos\tphi,
\end{split}
 \eea
where $\Jk$ and $\Jk_\pm$ are regular (even at the axis) as well as anti-self dual, i.e., $\Jk,\, \Jk_\pm\in\sk_1(\mathbb{C})$, and where the coordinate $\tphi$ involved in the definition of $j_\pm$ has been introduced in \eqref{eq:definitionof:specificchoice:normalizedcoord}. Note in particular the following identities, with $(x^1_p, x^2_p)$ introduced in \eqref{eq:definitionof:specificchoice:normalizedcoord}, 
\bea\lab{eq:usefulalgebraicidentitiesinvolvingscalarproductsReJkReJkpm}
\bsplit
\Re(\Jk)\c\Re(\Jk) &=\frac{(\sin\th)^2}{|q|^2}, \qquad \Re(\Jk)\c\Re(\Jk_+)=-\frac{x^2_p}{|q|^2},\qquad \Re(\Jk)\c\Re(\Jk_-)=\frac{x^1_p}{|q|^2},\\
\dual(\Re(\Jk))\c\Re(\Jk_+) &=\frac{\cos\th x^1_p}{|q|^2},\qquad \dual(\Re(\Jk))\c\Re(\Jk_-)=\frac{\cos\th x^2_p}{|q|^2}.
\end{split}
\eea

The complexified Ricci coefficients w.r.t. this principal null pair are given by  
\bea\lab{eq:KerrvaluesofcomplexifiedRicci}
\begin{split}
&\Xh=\Xbh=\Xi=\Xib=\omb=0, \qquad  \tr X=\frac{2\De \ov{q}}{|q|^4}, \qquad \tr\Xb=-\frac{2}{\ov{q}}, \qquad \om=- \frac 12\pr_r\left(\frac{\Delta}{|q|^2}\right),\\
&H=Z=\frac{a}{\ov{q}}\Jk=\frac{aq}{|q|^2}\Jk, \qquad \Hb=-\frac{a}{q}\Jk= -\frac{a\ov{q}}{|q|^2}\Jk.
\end{split}
\eea

The complexified curvature components are given by
\bea
A=B=\Bb=\Ab=0,\quad P=-\frac{2m}{q^3}.
\eea

Also, the principal null frame acts on the the normalized coordinates of Lemma \ref{lem:specificchoice:normalizedcoord} as follows 
\bea\lab{eq:actionofingoingprincipalnullframeonnormalizedcoordinates}
\bsplit
e_3(r)&=-1, \qquad\quad e_4(r)=\frac{\De}{|q|^2}, \qquad\qquad\qquad\qquad\quad\, e_1(r)=0, \quad\,\,\,\,\, e_2(r)=0,\\
e_3(\tau)&=\tmod'(r), \quad e_4(\tau)=\frac{2(r^2+a^2) - \De\tmod'(r)}{|q|^2}, \quad e_1(\tau)=0, \quad\,\,\,\, e_2(\tau)=\frac{a\sin\th}{|q|},\\
e_3(\th)&=0, \qquad\quad\,\,\,\,\, e_4(\th)=0, \qquad\qquad\qquad\qquad\qquad\,\, e_1(\th)=\frac{1}{|q|}, \quad e_2(\th)=0,\\
e_3(\tphi)&=\phimod'(r), \quad e_4(\tphi)=\frac{2a -\De\phimod'(r)}{|q|^2},\qquad\quad\,\,\,\,\, e_1(\tphi)=0, \quad\,\,\,\, e_2(\tphi)=\frac{1}{|q|\sin\th}.
\end{split}
\eea
In particular, recalling that $x^1_p=\sin\th\cos\tphi$ and $x^2_p=\sin\th\sin\tphi$, we have
\bea\lab{eq:actionofingoingprincipalnullframeonnormalizedcoordinates:bis}
\bsplit
e_3(x^1_p)&=-\phimod'(r)x^2_p, \qquad\, e_4(x^1_p)=-\frac{2a -\De\phimod'(r)}{|q|^2}x^2_p,\\ 
e_3(x^2_p)&=\phimod'(r)x^1_p, \qquad\quad e_4(x^2_p)=\frac{2a -\De\phimod'(r)}{|q|^2}x^1_p,
\end{split}
\eea
and, in view of the definition of $\Jk$ and $\Jk_\pm$,  
\bea\lab{eq:actionofingoingprincipalnullframeonnormalizedcoordinates:ter}
\DD(\tau)=a\Jk, \qquad \DD(\cos\th)=i\Jk,\qquad \DD(x^1_p)=\Jk_+, \qquad \DD(x^2_p)=\Jk_-.
\eea

Moreover the derivatives of $\Jk$ and $\Jk_\pm$ w.r.t. the principal null frame satisfy the following
\bea
\bsplit
\nab_3\Jk &=\frac{1}{\ov{q}}\Jk, \qquad \nab_4\Jk =- \frac{\De \ov{q}}{|q|^4}\Jk, \qquad \nab_3\Jk_\pm =\frac{1}{\ov{q}}\Jk_\pm, \qquad \nab_4 \Jk_\pm =- \frac{\De \ov{q}}{|q|^4}\Jk_{\pm} \mp  \frac{2a}{|q|^2}\Jk_{\mp},\\
\ov{\DD}\c\Jk &=\frac{4i(r^2+a^2)\cos\th}{|q|^4},\qquad \DD\hot\Jk=0, \\ 
\ov{\DD}\c \Jk_+ &= - \frac{4r^2 }{|q|^4}x^1_p - \frac{4ia^2\cos\th}{|q|^4}x^2_p, \qquad \ov{\DD}\c \Jk_- = - \frac{4r^2 }{|q|^4}x^2_p + \frac{4ia^2\cos\th}{|q|^4}x^1_p, \qquad \DD\hot \Jk_{\pm} =0.
\end{split}
\eea

%%%%%%%%%%%%%%%%%%%%%%%%%%%%%%%%%%%

\section{Teukolsky wave-transport system}
\lab{sec:sectionontheteukolskywavetransportsystem}

%%%%%%%%%%%%%%%%%%%%%%%%%%%%%%%%%%%

The main goal of this section is to state our main result concerning the Teukolsky wave-transport system in perturbations of Kerr. To this end, we first need to introduce some notations and properties for perturbations of Kerr.

%%%%%%%%%%%%%%%%%%%%%%%%%%%%%%%%%%%

\subsection{Perturbations of Kerr}
\lab{sec:Kerrpertbasic}

%%%%%%%%%%%%%%%%%%%%%%%%%%%%%%%%%%%

We consider a given vacuum spacetime $(\MM, \g)$ together with a null pair $(e_3, e_4)$ and its corresponding horizontal structure as in Section \ref{subsection:review-horiz.structures}. We will use the complexified Ricci and curvature coefficients  of Definition \ref{def:complexRicciandcurvaturecoefficients}. Moreover, we assume that $\MM$ is endowed with a pair of constants $(a, m)$, scalar functions $(\tau, r, \th, \vphi)$ and complex horizontal 1-forms $\Jk$, $\Jk_\pm$. We start by defining linearized  quantities.

%%%%%%%%%%%%%%%%%%%%%%%%%%%%%%%%%%%%%%%%%%%%%

\subsubsection{Definition of linearized quantities}
\lab{sec:definitionoflinearizedquantities:chap4}

%%%%%%%%%%%%%%%%%%%%%%%%%%%%%%%%%%%%%%%%%%%%%

Recall from Section \ref{subsect:principalnullpairinKerr} that, w.r.t. the principal null pair regular across  the future event horizon, the following quantities vanish in Kerr
\beaa
\Xh, \quad \Xbh,\quad \Xi, \quad \Xib, \quad \omb, \quad A, \quad B, \quad \Bb, \quad \Ab,\quad \nab(r), \quad e_4(\th), \quad e_3(\th), \quad \DD\hot\Jk.
\eeaa

We renormalize below all other quantities, not vanishing in Kerr as follows.
\begin{definition}
\lab{def:renormalizationofallnonsmallquantitiesinPGstructurebyKerrvalue}
We  define  the following renormalizations.
\begin{enumerate}
\item Linearization of the complex-valued Ricci and curvature coefficients:
\beaa
\bsplit
\trXc &:= \tr X-\frac{2\ov{q}\De}{|q|^4}, \qquad\trXbc := \tr\Xb+\frac{2}{\ov{q}},\\ 
\Pc &:= P+\frac{2m}{q^3},\qquad\qquad\, \omc  := \om  + \frac{1}{2}\pr_r\left(\frac{\De}{|q|^2} \right),\\
 \Hc &:= H-\frac{aq}{|q|^2}\Jk, \qquad\quad \Hbc:=\Hb+\frac{a\ov{q}}{|q|^2}\Jk,\qquad\quad \Zc := Z-\frac{aq}{|q|^2}\Jk.
 \end{split}
\eeaa

\item Linearization of derivatives of the scalar functions\footnote{Note that in Kerr we have $e_3(q)=e_3(r)$ and $e_4(q)=e_4(r)$ so that it suffices to linearize $e_3(r)$ and $e_4(r)$.} $r$, $\cos\th$, $q$, $\tau$, $x^1_p$ and $x^2_p$:
\beaa
\bsplit
\widecheck{e_3(r)} :=& e_3(r)+1, \qquad\qquad\qquad\,\,\,\,\, \widecheck{e_4(r)} := e_4(r)-\frac{\Delta}{|q|^2},\\
\widecheck{e_3(\tau)}:=& e_3(\tau)-\tmod'(r), \qquad\qquad\, \widecheck{e_4(\tau)}:=e_4(\tau)-\frac{2(r^2+a^2) - \De\tmod'(r)}{|q|^2},\\
\widecheck{e_3(x^1_p)}:=& e_3(x^1_p)+\phimod'(r)x^2_p, \qquad \widecheck{e_4(x^1_p)}:=  e_4(x^1_p)+\frac{2a -\De\phimod'(r)}{|q|^2}x^2_p, \\ \widecheck{e_3(x^2_p)}:=& e_3(x^2_p)-\phimod'(r)x^1_p, \qquad \widecheck{e_4(x^2_p)}:=e_4(x^2_p)-\frac{2a -\De\phimod'(r)}{|q|^2}x^1_p,\\
\widecheck{\DD q} :=& \DD q+a\Jk, \qquad\qquad\, \widecheck{\DD \ov{q}} :=\DD \ov{q}-a\Jk,\qquad \,\,\,\widecheck{\DD(\cos\th)} := \DD(\cos\th) -i\Jk,\\
\widecheck{\DD(\tau)}:=& \DD(\tau)-a\Jk, \qquad  \widecheck{\DD(x^1_p)} := \DD(x^1_p) - \Jk_+, \qquad \widecheck{\DD(x^1_p)} := \DD(x^2_p) - \Jk_-.
\end{split}
\eeaa

\item Linearization of derivatives of the complex 1-forms $\Jk$ and $\Jk_{\pm}$:
\beaa
\bsplit
\widecheck{\nab_3\Jk}:=&\nab_3\Jk -\frac{1}{\ov{q}}\Jk, \qquad\qquad  \widecheck{\nab_4\Jk}:=\nab_4\Jk +\frac{\De \ov{q}}{|q|^4}\Jk,\qquad \widecheck{\ov{\DD}\c\Jk}:= \ov{\DD}\c\Jk-\frac{4i(r^2+a^2)\cos\th}{|q|^4},\\
\widecheck{\nab_3\Jk_\pm}:=&\nab_3\Jk_\pm - \frac{1}{\ov{q}}\Jk_\pm\pm\phimod'(r)\Jk_{\mp},\qquad \widecheck{\nab_4 \Jk_\pm}:=\nab_4 \Jk_\pm + \frac{\De \ov{q}}{|q|^4}\Jk_{\pm} \pm\frac{2a-\De\phimod'(r)}{|q|^2}\Jk_{\mp},\\
\widecheck{\ov{\DD}\c \Jk_+}:=&\ov{\DD}\c \Jk_+ + \frac{4r^2 }{|q|^4}x^1_p + \frac{4ia^2\cos\th}{|q|^4}x^2_p,\qquad \widecheck{\ov{\DD}\c \Jk_-}:=\ov{\DD}\c \Jk_- + \frac{4r^2 }{|q|^4}x^2_p - \frac{4ia^2\cos\th}{|q|^4}x^1_p.
\end{split}
\eeaa
 \end{enumerate}
\end{definition}

%%%%%%%%%%%%%%%%%%%%%%%%%%%%%%%%%%%%%%%%%%%%%

\subsubsection{Definition of the notations $\Ga_b$ and $\Ga_g$ for error terms}
\lab{sec:definitionofGabandGagfirsttime}

%%%%%%%%%%%%%%%%%%%%%%%%%%%%%%%%%%%%%%%%%%%%%

\begin{definition}
\lab{definition.Ga_gGa_b}
The set of all linearized quantities is of the form $\Ga_g\cup \Ga_b$ with  $\Ga_g,  \Ga_b$
 defined as follows.
 \begin{enumerate}
\item 
 The set $\Ga_g$ is given by $\Ga_g=\Ga_{g,1}\cup \Ga_{g, 2}\cup\Ga_{g,3}$   with
 \bea
 \bsplit
 \Ga_{g,1} &= \Big\{\Xi, \quad \omc, \quad\trXc,\quad  \Xh,\quad \Zc,\quad \Hbc, \quad \trXbc , \quad r\Pc, \quad  rB, \quad  rA\Big\},\\
 \Ga_{g,2} &= \Big\{\widecheck{e_4(r)}, \quad r^{-1}\nab(r), \quad \widecheck{e_4(\tau)},   \quad r^{-1}\widecheck{\DD(\tau)}, \quad e_4(\cos\th)\Big\},\\
  \Ga_{g,3} &= \Big\{r\widecheck{\nab_4\Jk}\Big\}.
 \end{split}
 \eea
 
 \item The set $\Ga_b$ is given by $\Ga_b=\Ga_{b,1}\cup \Ga_{b, 2}\cup \Ga_{b,3}$   with
 \bea
 \bsplit
 \Ga_{b,1}&= \Big\{\Hc, \quad \Xbh, \quad \omb, \quad \Xib,\quad  r\Bb, \quad \Ab\Big\},\\
  \Ga_{b, 2}&= \Big\{r^{-1}\widecheck{e_3(r)}, \quad r^{-1}\widecheck{e_3(\tau)}, \quad  \widecheck{\DD(\cos\th)}, \quad e_3(\cos\th)\Big\}, \\
   \Ga_{b,3}&=\bigg\{ r\,\widecheck{\ov{\DD}\c\Jk}, \quad r\,\DD\hot\Jk, \quad r\,\widecheck{\nab_3\Jk}\bigg\}. 
   \end{split}
 \eea
\end{enumerate}
\end{definition}

\begin{remark}
The justification for the above decompositions has to do with the expected  decay properties of the linearized  components in perturbations of Kerr. More precisely, we will consider perturbations of Kerr for which $\Ga_g$ and $\Ga_b$ satisfy the following estimates,    see Section \ref{sec:controlofGabandGagfromBA} for details, 
\bea\lab{eq:expectedbehaviorGabGag:chap2}
\bsplit
\big|\dk^{\leq s}\Ga_g|&\les \ep \min\Big\{ r^{-2 }\tau^{-1/2-\dec},  \, r^{-1}\tau^{-1-\dec} \Big\}, \\
\big|\nab_3\dk^{\leq s-1}\Ga_g| &\les \ep  r^{-2 }\tau^{-1-\dec},\\
\big|\dk^{\leq s}\Ga_b\big| &\les \ep  r^{-1 }\tau^{-1-\dec},
\end{split}
\eea
for a small constant $\dec>0$, where $\dk=\{\nab_3, r\nab_4, \dkb=r\nab \}$ denotes weighted derivatives.
\end{remark}

\begin{remark}
In view of \eqref{eq:expectedbehaviorGabGag:chap2}, we note that $\Ga_g$ satisfies the assumptions of $\Ga_b$ and that $r^{-1}\Ga_b$ satisfies the assumptions of $\Ga_g$. Thus, in the rest of the paper, we will systematically replace $\Ga_g+\Ga_b$ by $\Ga_b$ and $r^{-1}\Ga_b+\Ga_g$ by $\Ga_g$.  
\end{remark}

\begin{remark}\lab{rmk:linearizedquantitiesJkpmandxbpnotincludedindefintionGagGab}
Note that we do not include linearized quantities associated to $\Jk_{\pm}$ and first order derivatives of $x^1_p$ and $x^2_p$ in the definition of $(\Ga_b, \Ga_g)$. Indeed, these quantities will only be used in the context of energy-Morawetz estimates for few derivatives when applying the results of \cite{MaSz24} \cite{MaSz26}, see Section \ref{sec:energyMorawetzesitmatesforTeukoslkyonMM:upto15derivatives}, while we will need to estimate $(\Ga_b, \Ga_g)$ up to top order derivatives.
\end{remark}

%%%%%%%%%%%%%%%%%%%%%%%%%%

\subsubsection{Commutation formulas revisited}

%%%%%%%%%%%%%%%%%%%%%%%%%%

The following lemma extends the commutation formulas of Section \ref{sec:generalcommutationformulasrealcase} to complex derivatives and makes use of the notations $(\Ga_b, \Ga_g)$ of Definition \ref{definition.Ga_gGa_b}.

   \begin{lemma}\label{LEMMA:COMMUTATION-FORMULAS-1}
   The following commutation formulas hold true.
   \begin{enumerate}
   \item 
   Let $h \in \sk_0(\CCC)$ $s$-conformally invariant. Then 
   \bea\label{eq:comm-nab4-nab3-DD-h-precise}
   \begin{split}
 \, [\nab_4 , \DD]h  &=  -\frac{1}{2}\tr X\DD h+(\Hb+Z)\nab_4 h -\frac 1 2 \Xh \c\ov{\DD} h+\Xi \nab_3h , \\
 \, [\nab_3 , \DD]h  &=   -\frac{1}{2}\tr \Xb\DD h+(H-Z)\nab_3 h -\frac 1 2 \Xbh\c\ov{\DD} h+\Xib \nab_4h. 
\end{split}
 \eea    
    
   \item
    Let $F\in  \sk_1 (\mathbb{C})$.  Then
    \bea\label{eq:comm:nab4-nab3-DDhot-precise}
    \begin{split}
\, [\nab_4,  \DD \hot] F  &=-\frac 1 2 \tr X\left( \DD \hot  F+\Hb\hot F\right)+(\Hb+Z)\hot\nab_4 F+ \Xi \hot \nab_3 F \\
 &-B \hot F - \frac 1 2 \tr \Xb \Xi \hot  F-\frac 1 2\Xh \c \ov{\DD} F+\frac12\Xh (\ov{\Hb}\c F)+ (\Ga_b \c \Ga_g) F, \\
\, [\nab_3,  \DD \hot] F  &=-\frac 1 2 \tr \Xb\left( \DD \hot  F+H \hot F\right)+(H-Z)\hot\nab_3 F+ \Xib \hot \nab_4 F \\
 &+\Bb \hot F - \frac 1 2 \tr X \Xib \hot  F-\frac 1 2\Xbh \c \ov{\DD} F+\frac12\Xbh (\ov{H}\c F)+ (\Ga_b \c \Ga_g) F.
 \end{split}
\eea
Using the schematic structure of the error terms, the above can be written as
\bea
 \lab{commutator-nab-43-D-hot} 
\bsplit
\, [\nab_4, \mathcal{D}\hot ]F&=- \frac 1 2 \tr X( \mathcal{D}\hot F + \underline{H} \hot F)+ (\underline{H}+Z) \hot \nab_4 F+\Xi \c \nabc_3 F+ r^{-1} \Ga_g \c  \dk^{\leq 1} F, \\
\, [\nab_3, \mathcal{D}\hot] F&=- \frac 1 2 \tr \Xb( \mathcal{D}\hot F + H \hot F)+ (H-Z) \hot \nab_3 F+r^{-1}\Ga_b \c \dk^{\leq 1} F.
\end{split}
\eea

\item Let $U\in \sk_2(\mathbb{C})$. Then
\bea\label{eq:comm:nab4-nab3-DDc-precise}
\begin{split}
\,[\nab_4,  \ov{\DD} \c ] U &=-\frac 1 2\ov{\tr X} \big(  \ov{\DD}\c U - 2\ov{\Hb} \c U\big) +\ov{(\Hb + Z)}\c\nab_4 U+ \ov{\Xi} \c \nab_3 U \\
&+2\ov{B} \c U  -\frac 1 2 \ov{\tr \Xb} \ov{\Xi}\c  U  -\frac 1 2 \Xh \c \ov{\DD} U-\frac 1 2 (\ov{\Xh}\c U)\ov{\Hb}+ (\Ga_b \c \Ga_g) U, \\
\,[\nab_3,  \ov{\DD} \c ] U &=-\frac 1 2\ov{\tr \Xb} \big(  \ov{\DD}\c U - 2\ov{H} \c U\big) +\ov{(H- Z)}\c\nab_3 U+ \ov{\Xib} \c \nab_4 U \\
&-2\ov{\Bb} \c U  -\frac 1 2 \ov{\tr X} \ov{\Xib}\c  U  -\frac 1 2 \Xbh \c \ov{\DD} U-\frac 1 2 (\ov{\Xbh}\c U)\ov{H}+ (\Ga_b \c \Ga_g) U.
\end{split}
\eea
Using the schematic structure of the error terms, the above can be written as
 \bea\label{commutator-nabc-3-ov-DDc-U}
 \bsplit
\, [\nab_4, \ov{\DD}\c] U&=- \frac 1 2\ov{\tr X}\, ( \ov{\DD} \c U - 2 \ov{\Hb} \c U)+\ov{(\Hb+Z)} \c \nab_4 U+\Xi \c \nabc_3 U +r^{-1} \Ga_g  \c \dk^{\leq 1} U, \\
\, [\nab_3, \ov{\DD}\c] U&=- \frac 1 2\ov{\tr\Xb}\, ( \ov{\DD} \c U -  2 \ov{H} \c U)+\ov{(H-Z)} \c \nab_3 U + r^{-1}\Ga_b \c  \dk^{\leq 1} U.
\end{split}
\eea
Similarly, for $F\in  \sk_1 (\mathbb{C})$ we have
 \bea\label{commutator-nab4-ov-DDcF}
 \bsplit
\, [\nab_4, \ov{\DD}\c] F&=- \frac 1 2\ov{\tr X}\, ( \ov{\DD} \c F -  \ov{\Hb} \c F)+\ov{(\Hb+Z)} \c \nab_4 F +\ov{\Xi} \c \nabc_3 F+ r^{-1} \Ga_g  \c \dk^{\leq 1} F, \\
\, [\nab_3, \ov{\DD}\c] F&=- \frac 1 2\ov{\tr\Xb}\, ( \ov{\DD} \c F -   \ov{H} \c F)+\ov{(H-Z)} \c \nab_3 F +r^{-1}\Ga_b \c  \dk^{\leq 1} F.
\end{split}
\eea
\item Let $U\in \sk_2(\mathbb{C})$. Then
\bea\label{correct-commutator-1} \lab{commutator-nab-4-D-c} \label{commutator-nab-3-nab-4-Psi}\label{correct-commutator}
\,[\nab_3, \nab_4]U = - 2\om \nab_3 U+ 2\omb \nab_4 U  + 2 (\eta_c-\etab_c) \nab_c U +4i \left(- \rhod+ \eta \wedge \etab  \right) U+\left(\Ga_b  \c \Ga_g \right) U.
\eea
Also,
\bea\label{commutator-nab-3-nab-a-U}
\bsplit
\, [\nab_3, \nab_a] U_{bc}&=-\frac  1 2   \trchb\, \Big(\nab_a U_{bc}+\eta_bU_{ac}+\eta_c U_{ab}-\de_{a b}(\eta \c U)_c-\de_{a c}(\eta \c U)_b \Big)\\
&-\frac 1 2 \atrchb\, \Big(\dual \nab_a  U_{bc} +\eta_b \dual U_{ac}+\eta_c \dual U_{ab}- \in_{a b}(\eta \c  U)_c- \in_{a c}(\eta \c  U)_b \Big)\\
 &+(\eta_a-\ze_a)\nab_3 U_{bc}+r^{-1}\Ga_b \c \dk^{\leq 1} U,
 \end{split}
\eea
\end{enumerate}
where the above error terms may also contain terms which are 
quadratic   in  the perturbation and  enjoy  better  decay  properties,  or are higher order  and decay at least as good.
\end{lemma}

\begin{proof}
See  Lemma 4.2.1 in \cite{GKS22}. 
\end{proof}

%%%%%%%%%%%%%%%%%%%%%%%%%%%%%%%%%%

\subsubsection{Approximate Killing vectorfields $\T$ and $\Z$}

%%%%%%%%%%%%%%%%%%%%%%%%%%%%%%%%%%

\begin{definition}\lab{Definition:vfsTZ} 
In $\MM$, we define $\T$ and $\Z$ as follows:
\bea\lab{eq:definitionofTandPhithataretheapproximateKillingvectorifeldinKerrpert}
\T := \frac{1}{2}\left(e_4+\frac{\Delta}{|q|^2}e_3 -2a\Re(\Jk)^be_b\right),\quad \Z := \frac 1 2 \left(2(r^2+a^2)\Re(\Jk)^be_b -a(\sin\th)^2 e_4 -\frac{a(\sin\th)^2\De}{ |q|^2} e_3\right).
\eea
\end{definition}

The following lemma provides basic relations between $\Lieb_\T, \Lieb_\Z$ and $\nab_\T, \nab_\Z$.
 \begin{lemma}\lab{lemma:basicpropertiesLiebTfasdiuhakdisug:chap9}
For a horizontal covariant k-tensor $U$, we have
\beaa
\nab_\T U_{b_1\cdots b_k} &=& \Lieb_\T U_{b_1\cdots b_k} +\frac{2amr\cos\th}{|q|^4}\sum_{j=1}^k\in_{b_jc} U_{b_1\cdots c\cdots b_k}+\Ga_b \c U,\\
\nab_\Z U_{b_1\cdots b_k} &=& \Lieb_\Z U_{b_1\cdots b_k} -\frac{\cos\th((r^2+a^2)^2-a^2(\sin\th)^2\De)}{|q|^4}\sum_{j=1}^k\in_{b_jc} U_{b_1\cdots c\cdots b_k}+r\Ga_b \c U.
\eeaa
\end{lemma}

\begin{proof}
See Lemma 4.3.6 in \cite{GKS22}.
\end{proof}

\begin{lemma}\lab{lemma:commutatorbetweenLieTLieZandnabnab4nab3}
We have 
\bea\lab{commutatorbetweenLieTLieZandnabnab4nab3:1}
\bsplit
 [\nab, \Lieb_T]U=&  \Ga_b \c \nab_3U +r^{-1}\dk^{\leq 1}(\Ga_b \c U),\\ 
   [\nab_4, \Lieb_T]U  =& \Ga_g \c \nab_3U+r^{-1}\dk^{\leq 1}(\Ga_b \c U),\\
    [\nab_3, \Lieb_T]U=& \dk^{\leq 1}(\Ga_b \c U),
\end{split}
\eea
and {
\bea\lab{commutatorbetweenLieTLieZandnabnab4nab3:2}
\bsplit
 [\nab, \Lieb_Z]U =&  r\Ga_g \c \nab_3U +\dk^{\leq 1}(\Ga_b\c U),\\ 
   [\nab_4, \Lieb_Z]U=&  r \xi \c \nab_3U+\dk^{\leq 1}(\Ga_g \c U),\\
     [\nab_3, \Lieb_Z]U=     & r\Ga_b\nab_3 U +r\dk^{\leq 1}(\Ga_g \c  U).
\end{split}
\eea}
This yields in particular the following non-sharp commutators
\bea
[r\nab, \Lieb_\T]=r\Ga_b\dk^{\leq 1}, \qquad [r\nab, \Lieb_\Z]=r^2\Ga_g\dk^{\leq 1}.
\eea
\end{lemma} 

\begin{proof}
See Lemma C.5.2 and (4.3.3) in \cite{GKS22}.
\end{proof}

%%%%%%%%%%%%%%%%%%%%%%%%%%%%%%%%%%%

\subsubsection{Decomposition of the wave operator in null frames}

%%%%%%%%%%%%%%%%%%%%%%%%%%%%%%%%%%%

The following lemma provides the decomposition of $\squared_k$ in null frames. 
\begin{lemma}\label{lemma:expression-wave-operator}
The wave operator for $\psi\in {\sk_k(\mathbb{C})}$, $k=0,1,2$, is given by
\bea\label{eq:wave-squared}
\begin{split}
\squared_k \psi&=-\nab_4 \nab_3 \psi  -\frac 1 2 \trchb \nab_4\psi+\left(2\om -\frac 1 2 \trch\right) \nab_3\psi+\lap_k\psi+2\etab \c\nab \psi \\
&+ ki \left( \rhod- \eta \wedge \etab \right) \psi+(\Ga_b \c \Ga_g) \c \psi,\\
\squared_k \psi&=-\nab_3\nab_4\psi +\left(2\omb -\frac 1 2 \trchb\right)\nab_4\psi -\frac 1 2 \trch \nab_3\psi+\lap_k\psi+2\eta \c\nab \psi \\
&- ki \left( \rhod- \eta \wedge \etab \right) \psi+(\Ga_b \c \Ga_g) \c \psi,
\end{split}
\eea
where $\lap_k=\nab^a \nab_a$ denotes the horizontal Laplacian for $k$-tensors. 
\end{lemma}

\begin{proof}
See Lemma 4.7.5 in \cite{GKS22} for the first identity of \eqref{eq:wave-squared} in the case $k=2$. The proof of Lemma 4.7.5 in \cite{GKS22} immediately extends to $k=0,1,2$ and to the second identity of \eqref{eq:wave-squared}.
\end{proof}

We will also need to decompose the following analog of Lemma \ref{lemma:expression-wave-operator} using complex derivatives.
\begin{corollary}\label{corollary-wave-complex} 
We have, for $\psi\in \sk_2(\CCC)$,
\bea
\begin{split}
\squared_2 \psi&=-\nab_4 \nab_3 \psi +\frac 1 4  \DD\hot( \DDb \c \psi)+\left(2\om -\frac 1 2 \tr X\right) \nab_3\psi- \frac 1 2 \tr\Xb \nab_4\psi+2\etab \c\nab \psi \\
& +  \left( - \frac 1 4 \tr X \ov{\tr \Xb}- \frac 1 4 \tr \Xb  \ov{\tr X}  - 2\ov{P}\right) \psi- 2i \left(\eta \wedge \etab\right)  \psi+(\Ga_b \c \Ga_g) \c \psi.
\end{split}
\eea
\end{corollary}

\begin{proof}
See Corollary 4.7.9 in \cite{GKS22}.
\end{proof}

%%%%%%%%%%%%%%%%%%%%%%%%%%%%%%%%%%%%%

\subsection{Teukolsky wave-transport system in perturbations of Kerr}
\lab{sec:precisederivationTeukolskywave-transportsystem:Kerrpert}

%%%%%%%%%%%%%%%%%%%%%%%%%%%%%%%%%%%%%

%%%%%%%%%%%%%%%%%%%%%%%%%%%%%%%%%%%%%%%%%%%%%%%%%%

\subsubsection{Definition  of $\pmb\phi_{s}^{(p)}$, $s=\pm 2$, $p=0,1,2$, in perturbations of Kerr}

%%%%%%%%%%%%%%%%%%%%%%%%%%%%%%%%%%%%%%%%%%%%%%%%%%

We start with the case $s=+2$. In perturbations of Kerr, we define $\pmb\phi_{+2}^{(p)}\in\sk_2(\mathbb{C})$, $p=0,1,2$ as follows
\bea\lab{eq:definitionofthephiplus2phierarchy:perturbationofKerr}
\bsplit
\pmb\phi_{+2}^{(0)} &:=\frac{\ov{q}}{q}A,\\
\pmb\phi_{+2}^{(1)} &:= \ov{q}^2\left(\nabc_3+\frac{1}{2}\trchb -\frac{3}{2}\frac{\atrchb^2}{\trchb} -2i\atrchb\right)A,\\
\pmb\phi_{+2}^{(2)} &:=|q|^2\ov{q}^2\left( \nabc_3\nabc_3 A + C_1  \nabc_3A + C_2   A\right),
\end{split}
\eea
where the scalar function $C_1$, $C_2$ are given by\footnote{Notice that we always work with a frame in ingoing normalization, so that $\trchb=-\frac{2r}{|q|^2}+O(\frac{\ep}{r^2})\neq 0$ which justifies the fact that we may divide by $\trchb$ in \eqref{eq:C1-C2-comparison-Ma}.} 
\bea\lab{eq:C1-C2-comparison-Ma}
\begin{split}
C_1&=2\trchb -\frac{2\atrchb^2}{\trchb}  -4 i \atrchb, \\
C_2  &= \frac 1 2 \trchb^2- 4\atrchb^2+\frac 3 2\frac{\atrchb^4}{\trchb^2} +  i \left(-2\trchb\atrchb +\frac{4\atrchb^3}{\trchb}\right).
\end{split}
\eea

\begin{remark}\lab{rmk:compasisionphiplus2p=0withMaSz26andphiplus2p=2withqfGKS22}
Note that the definition of $\pmb\phi_{+2}^{(0)}$ in \eqref{eq:definitionofthephiplus2phierarchy:perturbationofKerr} coincides with the one in (5.31) of \cite{MaSz26}. Also, comparing the definition of $\pmb\phi_{+2}^{(2)}$ in \eqref{eq:definitionofthephiplus2phierarchy:perturbationofKerr} with Definition 5.2.2 in \cite{GKS22} for $\qf$, we have $\qf=\pmb\phi_{+2}^{(2)}$. 
\end{remark}

Next, we consider the case $s=-2$. In perturbations of Kerr, we define $\pmb\phi_{-2}^{(p)}\in\sk_2(\mathbb{C})$, $p=0,1,2$ as follows
\bea\lab{eq:definitionofthephiminus2phierarchy:perturbationofKerr}
\bsplit
\pmb\phi_{-2}^{(0)} &:=\frac{q}{\ov{q}}\left(\frac{\De}{|q|^2}\right)^2\Ab,\\
\pmb\phi_{-2}^{(1)} &:=q^2\frac{\De}{|q|^2}\left(\nabc_4 +\frac{1}{2}\trch  -\frac{3}{2}\left(\frac{\atrch}{\trch}\right)_{\vartheta}\atrch -2i\atrch\right)\Ab,\\
\pmb\phi_{-2}^{(2)} &:=|q|^2q^2\left( \nabc_4\nabc_4\Ab + \underline{C}_1  \nabc_4\Ab + \underline{C}_2   \Ab\right),
\end{split}
\eea
where the scalar function $\underline{C}_1$, $\underline{C}_2$ are given by
\bea\label{eq:Cb1-Cb2-comparison-Ma}
\bsplit
\und{C}_1=& 2\trch -2\left(\frac{\atrch}{\trch}\right)_{\vartheta}\atrch -4i\atrch, \\ 
\und{C}_2  =& \frac 1 2 \trch^2 -4(\atrchb)^2 +\frac 3 2\left(\left(\frac{\atrch}{\trch}\right)_{\vartheta}\right)^2\atrch^2 +i\left(-2\trch\atrch +4\left(\frac{\atrch}{\trch}\right)_{\vartheta}\atrch^2\right),
\end{split}
\eea
and where we have introduced the following 0-conformally invariant scalar
\bea\lab{eq:defintionatrchovertrchvarthetamodificationtoavoiddividingbytrchnearHHplus}
\left(\frac{\atrch}{\trch}\right)_{\vartheta}:=\vartheta(r)\frac{a\cos\th}{r}+(1-\vartheta(r))\frac{\atrch}{\trch},
\eea
with $\vartheta$ supported in $r\leq 11m$ and $\vartheta=1$ for $r\leq 10m$.

\begin{remark}
Since we will always work with the frame in ingoing normalization, so that $\trch=\frac{2\De r}{|q|^4}+O(\frac{\ep}{r^2})$, we cannot divide by $\trch$ near $r=r_+$ which justifies the definition of the 0-conformally invariant scalar in   \eqref{eq:defintionatrchovertrchvarthetamodificationtoavoiddividingbytrchnearHHplus}.
\end{remark}

\begin{remark}\lab{rmk:compasisionphiminus2p=0withMaSz26andphiminus2p=2withqfGKS22}
Note that the definition of $\pmb\phi_{-2}^{(0)}$ in \eqref{eq:definitionofthephiminus2phierarchy:perturbationofKerr} coincides with the one in (5.31) of \cite{MaSz26}. Also, comparing the definition of $\pmb\phi_{-2}^{(2)}$ in \eqref{eq:definitionofthephiminus2phierarchy:perturbationofKerr} with Definition 5.2.2 in \cite{GKS22} for $\qfb$, we have $\qfb=\pmb\phi_{-2}^{(2)}$. 
\end{remark}

%%%%%%%%%%%%%%%%%%%%%%%%%%%%%%%%%%%%%%

\subsubsection{Teukolsky wave-transport system in perturbations of Kerr}

%%%%%%%%%%%%%%%%%%%%%%%%%%%%%%%%%%%%%%

The following theorem provides the tensorial wave equations satisfied by $\pmb\phi_{s}^{(p)}$, $s=\pm 2$, $p=0,1,2$.
\begin{theorem}\lab{thm:derivationoftheTeukolskytensorialwavesystemfors=plusminus2:kerrpert:alternateformnullframeinsteadcoordvectorfield}
Let $(\Ga_g, \Ga_b)$ be given by Definition \ref{definition.Ga_gGa_b}. Then, the horizontal tensors $\pmb\phi_{s}^{(p)}\in\sk_2(\mathbb{C})$, $s=\pm 2$, $p=0,1,2$, defined in \eqref{eq:definitionofthephiplus2phierarchy:perturbationofKerr} \eqref{eq:definitionofthephiminus2phierarchy:perturbationofKerr} satisfy the following  tensorial wave equations
\bsub
\lab{eq:TensorialTeuSysandlinearterms:rescaleRHScontaine2:general:Kerrperturbation:alternateformnullframeinsteadcoordvectorfield}
\bea
\lab{eq:TensorialTeuSys:rescaleRHScontaine2:general:Kerrperturbation:alternateformnullframeinsteadcoordvectorfield}
\bigg(\squared_2 -\frac{4ia\cos\th}{|q|^2}\nab_{\T}- \frac{4-2\de_{p0}}{|q|^2}\bigg){\phis{p}} = \widetilde{\L}_{s}^{(p)}[\pmb\phi_{s}]+\widetilde{\N}_{W,s}^{(p)}, \quad s=\pm2, \quad p=0,1,2,
\eea
where the linear coupling terms $\widetilde{\L}_{s}^{(p)}[{\pmb\phi_s}]$ have the following schematic forms
\bea
\lab{eq:tensor:Lsn:onlye_2present:general:Kerrperturbation:alternateformnullframeinsteadcoordvectorfield}
\bsplit
{\widetilde{\L}_{s}^{(0)}[\pmb\phi_{s}]}={}& (2sr^{-3} +O(mr^{-4}))\phis{1}+ O(mr^{-3}) \nab_{\widehat{\mathcal{X}}_s}^{\leq 1}\phis{0},\\
{\widetilde{\L}_{s}^{(1)}[\pmb\phi_{s}]}={}& (sr^{-3} +O(mr^{-4}))\phis{2}+ O(mr^{-3}) \nab_{\widehat{\mathcal{X}}_s}^{\leq 1}  \phis{1}+O(mr^{-2})\nab_{\widehat{Z}}^{\leq 1}\phis{0},\\
{\widetilde{\L}_{s}^{(2)}[\pmb\phi_{s}]}={}&O(mr^{-3})\phis{2}+O(mr^{-2})\nab_{\Z+a\T}^{\leq 1}\phis{1}+O(m^2 r^{-2})\phis{0},
\end{split}
\eea
with $\widehat{\mathcal{X}}_s$, $s=\pm 2$, and $\widehat{Z}$ being the regular vectorfields defined respectively by
\bea\lab{eq:formofregularhorizontalvectorfieldwidetildemathcalXs:Kerrperturbation}
\widehat{\mathcal{X}}_s:=|q|^2\left(s\Re(\Jk)\c\nab -\frac{2a\cos\th}{r}\dual\Re(\Jk)\c\nab\right), \quad s=\pm 2,
\eea
and 
\bea\lab{eq:formofregularhorizontalvectorfieldwidehatZ:Kerrperturbation}
\widehat{Z}:=|q|^2\left(\Re(\Jk)\c\nab +\frac{a(\cos\th)^2}{|q|^2}\nab_{\T}\right),
\eea
with all the coefficients in  \eqref{eq:tensor:Lsn:onlye_2present:general:Kerrperturbation:alternateformnullframeinsteadcoordvectorfield} being independent of coordinates $\tau$ and 
$\tphi$, and with the coefficients in front of the terms $\phis{2}$ and $\nab_{\Z+a\T}\phis{1}$ on the RHS of equation of ${\widetilde{\L}_{s}^{(2)}[\pmb\phi_{s}]}$ in \eqref{eq:tensor:Lsn:onlye_2present:general:Kerrperturbation:alternateformnullframeinsteadcoordvectorfield} being real functions, and where the nonlinear correction terms $\widetilde{\N}_{W,s}^{(p)}$ have the following schematic form
\bea\lab{eq:schematicformofNpWsplus2} 
\bsplit
\widetilde{\N}^{(0)}_{W,+2} =& r^{-1}  \dk^{\leq 1}\big( \Ga_g \c  B\big) + \nab_3\Xi  \c B {+\big( r^{-1}\dkb^{\leq1} \Pc+ r^{-3} \Ga_b \big)\c \Xi +(\Ga_g \c \Xi) \c \Bb}\\
 & +(\Ga_b \c \Ga_g)\c A +\Ga_b\nab_3A+r^{-1}\dk^{\leq 1}(\Ga_bA),\\
 \widetilde{\N}^{(1)}_{W,+2} =& r\dk^{\leq 2}\big( \Ga_g \c  B\big)+r^2\dk^{\leq 1}\big(\nab_3\Xi  \c B) +r\dk^{\leq 1}\big((\dkb^{\leq1} \Pc+ r^{-2} \Ga_b)\c \Xi\big)\\
& +r^2\dk^{\leq 1}((\Ga_g \c \Xi) \c \Bb) +r^2\dk^{\leq 1}(\Ga_b\c\nab_3A)+r\dk^{\leq 2}(\Ga_b\c A)+r^2\dk^{\leq 1}\big(\Ga_b \c \Ga_g\c A\big),\\
 \widetilde{\N}^{(2)}_{W,+2} =& r^2 \dk^{\leq 3} (\Ga_g \c (A, B))+ {\nab_3 (r^3 \dk^{\leq 2}( \Ga_g \c (A, B)))} + \nab_3 (r^3\dk^{\leq1}\big(\Xi \c (\dkb^{\leq1} \Pc,  r^{-2}\Ga_b )\big))\\
&+\dk^{\leq 1} (\Ga_g \c \qf) {+\nab_3(r^4(\dk^{\leq 1}(\Ga_g\c\Xi)\c\Bb))},
\end{split}
\eea
\bea\lab{eq:schematicformofNpWsminus2}
   \bsplit
  \widetilde{\N}^{(0)}_{W,-2} =& r^{-1}\dk^{\leq 1}\big(  \Ga_b \c \Bb \big) + \Ga_b \c \Ga_b \c \Ga_g +r^{-1}\Ga_g\pmb\phi_{-2}^{(1)}+\Ga_g\nab_3\Ab +r^{-1}\Ga_b\dk^{\leq 1}\Ab,\\
 \widetilde{\N}^{(1)}_{W,-2} =& \dk^{\leq 2}\big(\Ga_b\c \Ga_g\big)+r\dk^{\leq 2}(\Ga_b\c\xi),\\
  \widetilde{\N}^{(2)}_{W,-2} =&  r^2 \dk^{\leq 2}(\Ga_b \c (A, B))+ \dk^{\leq 3} (\Ga_g \c \Ga_b),
\end{split}
 \eea
\esub
under the following additional properties for the global null frame of $\MM$:
\begin{itemize}
\item For \eqref{eq:schematicformofNpWsplus2} to hold, we require in addition that 
\bea\label{eq:additional-conditions-Hc-Xi-gRW-eq:0}
\Hc \in \Ga_g, \qquad \nab_3\Xi\in r^{-1}\dk^{\leq 1}\Ga_g.
\eea

\item For \eqref{eq:schematicformofNpWsminus2} to hold, we require in addition that 
\bea
\lab{eq:Xi-Hb-chapter12:0}
\Xi\in r^{-2}\Ga_g, \qquad \Hbc\in r^{-1}\Ga_g.
\eea
\end{itemize}
\end{theorem}

\begin{remark}
An alternate form of the Teukolsky wave system \eqref{eq:TensorialTeuSysandlinearterms:rescaleRHScontaine2:general:Kerrperturbation:alternateformnullframeinsteadcoordvectorfield} will be given in Corollary  \ref{cor:derivationoftheTeukolskytensorialwavesystemfors=plusminus2:kerrpert}.
\end{remark}

\begin{proof}
See Theorem 3.15 in \cite{Sze}.
\end{proof}

The following proposition provides the transport equations satisfied by $\pmb\phi_{s}^{(p)}$, $s=\pm 2$, $p=0,1,2$.
\begin{proposition}\lab{prop:transportequationphis=plusminus2p=0and1}
The horizontal tensors $\pmb\phi_{s}^{(p)}\in\sk_2(\mathbb{C})$, $s=\pm 2$, $p=0,1,2$, defined in \eqref{eq:definitionofthephiplus2phierarchy:perturbationofKerr} \eqref{eq:definitionofthephiminus2phierarchy:perturbationofKerr} satisfy the following transport equations in ingoing normalization 
\bsub\lab{eq:transportequationsins=plus2andminus2caseforp=0and1}
\bea\lab{eq:transportequationsins=+2caseforp=0and1}
\nab_3\left(\frac{r\ov{q}}{q}\left(\frac{r^2}{|q|^2}\right)^{p-2}\pmb\phi_{+2}^{(p)}\right) &=& \frac{\ov{q}}{rq}\left(\frac{r^2}{|q|^2}\right)^{p-1}\pmb\phi_{+2}^{(p+1)}+\N_{T,+2}^{(p)},
\eea
and
\bea\lab{eq:transportequationsins=-2caseforp=0and1}
\nab_4\left(\frac{rq}{\ov{q}}\left(\frac{r^2}{|q|^2}\right)^{p-2}\pmb\phi_{-2}^{(p)}\right) = \frac{q}{r\ov{q}}\left(\frac{r^2}{|q|^2}\right)^{p-1}\frac{\De}{|q|^2}\pmb\phi_{-2}^{(p+1)}+\N_{T,-2}^{(p)},
\eea
where the error terms $\N_{T,+2}^{(p)}$, $p=0,1$, are given by 
\bea\lab{eq:transportequationsins=+2caseforp=0and1:RHSNT+2p=0and1}
\N_{T,+2}^{(0)}=r\Ga_bA, \qquad \N_{T,+2}^{(1)}=r^3\Ga_b\nab_3A+r^2\dk^{\leq 1}(\Ga_b)A,
\eea
and where the error terms $\N_{T,-2}^{(p)}$, $p=0,1$, are given by 
\bea\lab{eq:transportequationsins=-2caseforp=0and1:RHSNT-2p=0and1}
\N_{T,-2}^{(0)}=r\Ga_g\c\Ga_b, \qquad \N_{T,-2}^{(1)}=r^2\dk^{\leq 1}(\Xi)\c\Ga_b+r\Ga_g\dk^{\leq 1}\Ga_b.
\eea
\esub
\end{proposition}

\begin{remark}
Note that the transport equations \eqref{eq:transportequationsins=+2caseforp=0and1} \eqref{eq:transportequationsins=-2caseforp=0and1} coincide with the ones in (5.33a) (5.33b) of \cite{MaSz26}. 
\end{remark}

%\begin{remark}
%In the frame used for the control of $\qfb$ in \cite{GKS22}, we have $\Xi\in r^{-1}\Ga_g$ (see Section 3.6 of \cite{KS:Kerr}), so that we have in fact
%\beaa
%r^2\dk^{\leq 1}(\Xi)\c\Ga_b+r\Ga_g\dk^{\leq 1}\Ga_b=r\Ga_g\dk^{\leq 1}\Ga_b.
%\eeaa
%\end{remark}

\begin{proof}
See Proposition 3.16 in \cite{Sze}.
\end{proof}

%%%%%%%%%%%%%%%%%%%%%%%%%%%%%%%%%%%%%%%%%

\subsubsection{Alternate form of the Teukolsky wave system} 
  
%%%%%%%%%%%%%%%%%%%%%%%%%%%%%%%%%%%%%%%%%

We start with the following lemma. 
\begin{lemma}
Assume that we have\footnote{Recall from Remark \ref{rmk:linearizedquantitiesJkpmandxbpnotincludedindefintionGagGab} that $\widecheck{e_4(x^b_p)}$, $\widecheck{\DD(x^b_p)}$ and ${\widecheck{e_3(x^b_p)}}$ do not appear in Definition \ref{definition.Ga_gGa_b} for $(\Ga_g, \Ga_b)$.} 
\bea\lab{eq:propoertiesoflinearizedfirstorderderivativesoftphiintermsofGabGag}
\widecheck{e_4(x^b_p)}\in\Ga_g, \quad \widecheck{\DD(x^b_p)}\in\Ga_b, \quad {\widecheck{e_3(x^b_p)}}\in\Ga_b, \,\, b=1,2.
\eea
Moreover, assume that $\Jk$ and $\Jk_{\pm}$ satisfy\footnote{In practice, the identities in \eqref{eq:identiteReJkcReJkandReJkcReJkpmuptoGagtermswhichisenough} will be exactly enforced, i.e. without $\Ga_g$ error terms, see \eqref{eq:usefulalgebraicidentitiesinvolvingscalarproductsReJkReJkpm:Kerrpert:tilde}.} 
\bea\lab{eq:identiteReJkcReJkandReJkcReJkpmuptoGagtermswhichisenough}
\Re(\Jk)\c\Re(\Jk) =\frac{(\sin\th)^2}{|q|^2}+\Ga_g, \quad \Re(\Jk)\c\Re(\Jk_+)=-\frac{x^2_p}{|q|^2}+\Ga_g, \quad \Re(\Jk)\c\Re(\Jk_-)=\frac{x^1_p}{|q|^2}+\Ga_g.
\eea
Then, the vectorfields $\T$ and $\Z$ defined by \eqref{eq:definitionofTandPhithataretheapproximateKillingvectorifeldinKerrpert} satisfy
\bea\lab{eq:expressionofTandPhiintermsofprttandprPhi}
\T=\pr_{\tt}+r\Ga_b\dk, \qquad \Z=\pr_{\tphi}+r^2\Ga_g\dk.
\eea 
\end{lemma}

\begin{proof}
The proof follows immediately by plugging the definition \eqref{eq:definitionofTandPhithataretheapproximateKillingvectorifeldinKerrpert}  of $\T$ and $\Z$ in the identities $\T=\T(x^\a)\pr_{x^\a}$ and $\Z=\Z(x^\a)\pr_{x^\a}$, and then by expressing $\widecheck{e_\a(x^\b)}$ in terms of $\Ga_b$ and $\Ga_g$ using Definition \ref{definition.Ga_gGa_b}, \eqref{eq:propoertiesoflinearizedfirstorderderivativesoftphiintermsofGabGag} and \eqref{eq:identiteReJkcReJkandReJkcReJkpmuptoGagtermswhichisenough}.
\end{proof}

 We now provide an alternate form of the Teukolsky wave-system of Theorem \ref{thm:derivationoftheTeukolskytensorialwavesystemfors=plusminus2:kerrpert:alternateformnullframeinsteadcoordvectorfield} 
using coordinates vectorfield instead of null frames. 
\begin{corollary}\lab{cor:derivationoftheTeukolskytensorialwavesystemfors=plusminus2:kerrpert}
Let $(\Ga_g, \Ga_b)$ be given by Definition \ref{definition.Ga_gGa_b}, and assume also that \eqref{eq:propoertiesoflinearizedfirstorderderivativesoftphiintermsofGabGag} and 
\eqref{eq:identiteReJkcReJkandReJkcReJkpmuptoGagtermswhichisenough} hold. Then, the horizontal tensors $\pmb\phi_{s}^{(p)}\in\sk_2(\mathbb{C})$, $s=\pm 2$, $p=0,1,2$, defined in \eqref{eq:definitionofthephiplus2phierarchy:perturbationofKerr} \eqref{eq:definitionofthephiminus2phierarchy:perturbationofKerr} satisfy the following  tensorial wave equations
\bsub
\lab{eq:TensorialTeuSysandlinearterms:rescaleRHScontaine2:general:Kerrperturbation}
\bea
\lab{eq:TensorialTeuSys:rescaleRHScontaine2:general:Kerrperturbation:cor}
\bigg(\squared_2 -\frac{4ia\cos\th}{|q|^2}\nab_{\pr_{\tt}}- \frac{4-2\de_{p0}}{|q|^2}\bigg){\phis{p}} = \L_{s}^{(p)}[\pmb\phi_{s}]+\N_{W,s}^{(p)}, \quad s=\pm2, \quad p=0,1,2,
\eea
where the linear coupling terms $\L_{s}^{(p)}[{\pmb\phi_s}]$ have the following schematic forms
\bea
\lab{eq:tensor:Lsn:onlye_2present:general:Kerrperturbation:cor}
\bsplit
{\L_{s}^{(0)}[\pmb\phi_{s}]}={}& (2sr^{-3} +O(mr^{-4}))\phis{1}+ O(mr^{-3}) \nab_{\Xcal_s}^{\leq 1}\phis{0},\\
{\L_{s}^{(1)}[\pmb\phi_{s}]}={}& (sr^{-3} +O(mr^{-4}))\phis{2}+ O(mr^{-3}) \nab_{\Xcal_s}^{\leq 1}  \phis{1}+O(mr^{-2})\nab_{\pr_{\tphi}+a\pr_{\tt}}^{\leq 1}\phis{0},\\
{\L_{s}^{(2)}[\pmb\phi_{s}]}={}&O(mr^{-3})\phis{2}+O(mr^{-2})\nab_{\pr_{\tphi}+a\pr_{\tt}}^{\leq 1}\phis{1}+O(m^2 r^{-2})\phis{0},
\end{split}
\eea
with $\Xcal_s$, $s=\pm 2$, being the regular vectorfields defined by
\bea\lab{eq:formofregularhorizontalvectorfieldmathcalXs:Kerrperturbation:cor}
\mathcal{X}_{s} &:=& s\big(\pr_{\tphi}+a(\sin\th)^2\pr_\tau\big) - \frac{2a\cos\th}{r}\sin\th\pr_\th, \quad s=\pm 2, 
\eea
with all the coefficients in  \eqref{eq:tensor:Lsn:onlye_2present:general:Kerrperturbation:cor} being independent of coordinates $\tau$ and 
$\tphi$, and with the coefficients in front of the terms $\phis{2}$ and $\nab_{\pr_{\tphi}+a\pr_{\tt}}\phis{1}$ on the RHS of equation of ${\L_{s}^{(2)}[\pmb\phi_{s}]}$ in \eqref{eq:tensor:Lsn:onlye_2present:general:Kerrperturbation:cor} being real functions, and where the nonlinear correction terms $\N_{W,s}^{(p)}$ are related to $\widetilde{\N}_{W,s}^{(p)}$ introduced in \eqref{eq:TensorialTeuSysandlinearterms:rescaleRHScontaine2:general:Kerrperturbation:alternateformnullframeinsteadcoordvectorfield} as follows, for $s=\pm 2$, 
\bea\lab{eq:schematicformofNpWs:comparisionNWandwidetildeWW}
\bsplit
\N_{W,s}^{(0)} =& \widetilde{\N}_{W,s}^{(0)} + r^{-1}\Ga_b\c\dk\pmb\phi_s^{(0)},\\
\N_{W,s}^{(1)} =& \widetilde{\N}_{W,s}^{(1)} + r^{-1}\Ga_b\c\dk\pmb\phi_s^{(1)}+r^{-1}\Ga_b\c\dk\pmb\phi_s^{(0)},\\
\N_{W,s}^{(2)} =& \widetilde{\N}_{W,s}^{(2)} + r^{-1}\Ga_b\c\dk\pmb\phi_s^{(2)}+\Ga_g\c\dk\pmb\phi_s^{(1)}.
\end{split}
\eea
\esub
\end{corollary}

\begin{proof}
We start by proving the following identities 
\bea\lab{eq:linkbetweenmathcalXsandwidetildemathcalXsaswellaslinkmathcalZwithcoordvfandlinkenabZplusanabTwithcoordvf}
\widehat{\mathcal{X}}_s=\mathcal{X}_s+r\Ga_b\c\dk, \,\,\, s=\pm 2, \quad \widehat{Z}=\nab_{\pr_{\tphi}+a\pr_{\tt}}+r\Ga_b\c\dk, \quad \nab_{\Z}+a\nab_{\T}=\nab_{\pr_{\tphi}+a\pr_{\tt}}+r^2\Ga_g\dk,
\eea
where $\mathcal{X}_s$, $\widehat{\mathcal{X}}_s$, $s=\pm 2$, and $\widehat{Z}$ are the vectorfields given by 
\eqref{eq:formofregularhorizontalvectorfieldmathcalXs:Kerrperturbation:cor}, \eqref{eq:formofregularhorizontalvectorfieldwidetildemathcalXs:Kerrperturbation} and \eqref{eq:formofregularhorizontalvectorfieldwidehatZ:Kerrperturbation}. Indeed, in view of \eqref{eq:propoertiesoflinearizedfirstorderderivativesoftphiintermsofGabGag} and 
\eqref{eq:identiteReJkcReJkandReJkcReJkpmuptoGagtermswhichisenough}, we have
\beaa
\Jk\c\nab &=& \Jk\c\nab(x^\a)\nab_{\pr_{x^\a}}\\ 
&=& \Jk\c\big(a\Re(\Jk)+\Ga_b\big)\nab_{\pr_\tau}+\Jk\c r\Ga_g\nab_{\pr_r}+\Jk\c\Big(\Re(\Jk_+)\nab_{\pr_{x^1}}+\Re(\Jk_-)\nab_{\pr_{x^2}}+r^{-1}\Ga_b\dk\Big)\\
&=& a\Jk\c\Re(\Jk)\nab_{\pr_\tau}+\Jk\c\Re(\Jk_+)\nab_{\pr_{x^1}}+\Jk\c\Re(\Jk_-)\nab_{\pr_{x^2}}+r^{-1}\Ga_b\dk\\
&=& \frac{1}{|q|^2}\Big(a(\sin\th)^2\nab_{\pr_\tau}+\nab_{-x^2\pr_{x^1}+x^1\pr_{x^2}}\Big)+\frac{i\cos\th}{|q|^2}\nab_{x^1\pr_{x^1}+x^2\pr_{x^2}}\Big)+r^{-1}\Ga_b\dk\\
&=& \frac{1}{|q|^2}\nab_{\pr_{\tphi}+a(\sin\th)^2\pr_\tau}+\frac{i}{|q|^2}\nab_{\sin\th\pr_\th}+r^{-1}\Ga_b\dk
\eeaa
so that
\bea\lab{eq:formofReJkcnabanddualReJkcnab}
\Re(\Jk)\c\nab = \frac{1}{|q|^2}\nab_{\pr_{\tphi}+a(\sin\th)^2\pr_\tau}+r^{-1}\Ga_b\dk, \qquad \dual\Re(\Jk)\c\nab &=& \frac{1}{|q|^2}\nab_{\sin\th\pr_\th}+r^{-1}\Ga_b\dk.
\eea
Then, \eqref{eq:linkbetweenmathcalXsandwidetildemathcalXsaswellaslinkmathcalZwithcoordvfandlinkenabZplusanabTwithcoordvf} follows immediately from \eqref{eq:formofReJkcnabanddualReJkcnab} and \eqref{eq:expressionofTandPhiintermsofprttandprPhi}.

Next, note that \eqref{eq:TensorialTeuSys:rescaleRHScontaine2:general:Kerrperturbation:alternateformnullframeinsteadcoordvectorfield}
\eqref{eq:tensor:Lsn:onlye_2present:general:Kerrperturbation:alternateformnullframeinsteadcoordvectorfield} immediately yields \eqref{eq:TensorialTeuSys:rescaleRHScontaine2:general:Kerrperturbation:cor}
\eqref{eq:tensor:Lsn:onlye_2present:general:Kerrperturbation:cor} with $\N_{W,s}^{(p)}$, $s=\pm 2$, $p=0,1,2$, being defined by 
\beaa
\N_{W,s}^{(0)} &=& \widetilde{\N}_{W,s}^{(0)} + O(r^{-2})\nab_{\T-\pr_\tau}\pmb\phi_s^{(0)}+O(r^{-3})\nab_{\widehat{\Xcal}_s-\Xcal_s}\pmb\phi_s^{(0)},\\
\N_{W,s}^{(1)} &=& \widetilde{\N}_{W,s}^{(1)} + O(r^{-2})\nab_{\T-\pr_\tau}\pmb\phi_s^{(1)}+O(r^{-3})\nab_{\widehat{\Xcal}_s-\Xcal_s}\pmb\phi_s^{(1)}+O(r^{-2})\big(\widehat{Z} -\nab_{\pr_{\tphi}+a\pr_\tau}\big)\pmb\phi_s^{(0)},\\
\N_{W,s}^{(2)} &=& \widetilde{\N}_{W,s}^{(2)} + O(r^{-2})\nab_{\T-\pr_\tau}\pmb\phi_s^{(2)}+O(r^{-2})\nab_{\Z-\pr_{\tphi}}\pmb\phi_s^{(1)}+O(r^{-2})\nab_{\T-\pr_\tau}\pmb\phi_s^{(1)},
\eeaa
which together with \eqref{eq:expressionofTandPhiintermsofprttandprPhi} and \eqref{eq:linkbetweenmathcalXsandwidetildemathcalXsaswellaslinkmathcalZwithcoordvfandlinkenabZplusanabTwithcoordvf} implies, for $s=\pm 2$, 
\beaa
\N_{W,s}^{(0)} &=& \widetilde{\N}_{W,s}^{(0)} + r^{-1}\Ga_b\c\dk\pmb\phi_s^{(0)},\\
\N_{W,s}^{(1)} &=& \widetilde{\N}_{W,s}^{(1)} + r^{-1}\Ga_b\c\dk\pmb\phi_s^{(1)}+r^{-1}\Ga_b\c\dk\pmb\phi_s^{(0)},\\
\N_{W,s}^{(2)} &=& \widetilde{\N}_{W,s}^{(2)} + r^{-1}\Ga_b\c\dk\pmb\phi_s^{(2)}+\Ga_g\c\dk\pmb\phi_s^{(1)},
\eeaa
as stated. This concludes the proof of Corollary \ref{cor:derivationoftheTeukolskytensorialwavesystemfors=plusminus2:kerrpert}. 
\end{proof}

%%%%%%%%%%%%%%%%%%%%%%%%%%%%%%%%%%%%%%%%%

\subsubsection{Alternate form of the Teukolsky transport system for $s=-2$} 
  
%%%%%%%%%%%%%%%%%%%%%%%%%%%%%%%%%%%%%%%%%

In this section, we provide an alternate form of the transport equations \eqref{eq:transportequationsins=-2caseforp=0and1} for $\pmb\phi_{-2}^{(p)}$, $p=0,1$, which we derive instead for $\Ab$ and the -1-conformally invariant horizontal tensor $\underline{\Psi}$ given by the following proposition.

\begin{definition}\lab{def:PsibforintegrationofAbfromqfb:chap12}
Let $\Psib\in \sk_2(\CCC)$ {be} given by 
\beaa
\und{\Psi} &:=& \frac{q^4}{r^2}\left(\nabc_4+2\tr X -\frac{3|\tr X|^2}{2\trch}\right)\Ab.
\eeaa
\end{definition}

We have the following proposition.
\begin{proposition}\lab{cor:systemoftransportequationsforPsibandAbfromqfb:chap12}
Let $\Psib$ {be} as in Definition \ref{def:PsibforintegrationofAbfromqfb:chap12}. Then, $\und{\Psi}\in \dk^{\leq 1}\Ga_b$,  and $(\und{\Psi}, \Ab)$ satisfies the following system of transport equations 
\beaa
\nabc_4(r\Psib)=\frac{q}{r\ov{q}}\pmb\phi_{-2}^{(2)}  +r\dk^{\leq 1}(\Ga_g\c\Ga_b), \qquad \nabc_4\left( \frac{q^4}{r^3}\Ab\right)= \frac{1}{r}\Psib+r\Ga_g\c\Ga_b.
\eeaa
\end{proposition}

\begin{proof}
See Corollary 12.1.4 in \cite{GKS22}, noticing also that $\qfb=\pmb\phi_{-2}^{(2)}$ in view of Remark \ref{rmk:compasisionphiminus2p=0withMaSz26andphiminus2p=2withqfGKS22}.
\end{proof}

%%%%%%%%%%%%%%%%%%%%%%%%%%%%%%%%%%%%%%%%%%%%%%%%%%%%%%%%%%%%%%%%%%%%%%%%%%%%%%%%%%%

\section{Preliminaries for Sections \ref{sec:decayestimatesforAandAb:analogPartIIGKS22} and \ref{sec:higherordercurvatureestimatesforprovingThM8largea:00}} 
\lab{sec:prelimnariesfortheproofofdecayesetimatesforTeuk}
  
%%%%%%%%%%%%%%%%%%%%%%%%%%%%%%%%%%%%%%%%%%%%%%%%%%%%%%%%%%%%%%%%%%%%%%%%%%%%%%%%%%%

We start this section by providing the geometric set-up for Sections \ref{sec:decayestimatesforAandAb:analogPartIIGKS22} and \ref{sec:higherordercurvatureestimatesforprovingThM8largea:00}. We then exhibit a suitable coordinates systems and a regular triplet on $\MM$, both being needed in order to apply the energy-Morawetz estimates for Teukolsky equations in perturbations of Kerr of \cite{MaSz26} in Section \ref{sec:energyMorawetzesitmatesforTeukoslkyonMM:upto15derivatives}.

%%%%%%%%%%%%%%%%%%%%%%%%%%%%%%%%%%%%%%%%%%%%%%%%%%%%%%%%%%%%%%%%%%%%%%%%%%%%%%%%%%%%%%

\subsection{Geometric set-up for Sections \ref{sec:decayestimatesforAandAb:analogPartIIGKS22} and \ref{sec:higherordercurvatureestimatesforprovingThM8largea:00}}
\lab{sec:geometricsetupfordecsayTeukolsky} 
  
%%%%%%%%%%%%%%%%%%%%%%%%%%%%%%%%%%%%%%%%%%%%%%%%%%%%%%%%%%%%%%%%%%%%%%%%%%%%%%%%%%%%%%

%%%%%%%%%%%%%%%%%%%%%%%

\subsubsection{Choices of constants}
\lab{sec:smallnesconstants}

%%%%%%%%%%%%%%%%%%%%%%%

The following constants  are  involved in the derivation of decay estimates for Teukolsky:
\begin{itemize}
\item The constants $m>0$ and $a$, with $|a|<m$, are the mass and the angular momentum per unit mass of the Kerr solution relative to which the perturbation to Kerr is measured. 

\item The integer $\kl$ which corresponds to the maximum number of derivatives of the solution.

\item The size of the initial data  perturbation  is measured by $\ep_0>0$. 

\item The size of the perturbation to Kerr is measured by $\ep>0$. 

\item $r_0>0$ is tied to  $\Mint\cap\Mext=\{r=r_0\}$. 

\item The constant $\dhor$ is tied to the boundary of $\MM$ given by $\pr\MM=\AA=\{r=r_+(1-\dhor)\}$. 

\item The constant $\dred$ measures the width of the redshift region.

\item The constant $\dbl$ appears in the construction of normalized coordinates, see Lemma \ref{lem:specificchoice:normalizedcoord}.

\item The constant $\dec$ is tied to decay estimates in $(r, \tau)$ of Kerr perturbations, see Section \ref{sec:controlofGabandGagfromBA}. 

\item The constant $\de$ is tied to weights in the energy-Morawetz norms and $r^p$-norms of Section \ref{subsection:basicnormsforpsi}.
\end{itemize}

These  constants are chosen such that 
\bea\lab{eq:constraintsonthemainsmallconstantsepanddelta}
0<\ep_0,\,\ep\ll\dhor\ll\dred\ll \dbl \ll 1-\frac{|a|}{m}, \quad\,\,\, \ep_0,\,\ep\ll \de\ll\dec, \quad\,\,\, \ep_0,\,\ep\ll \frac{1}{r_0},\, \frac{1}{\kl},\quad\,\,\, \frac{1}{\kl}\ll\dec,
\eea
and
\bea\lab{eq:constraintbetweenepep*ep0}
\ep=\ep_0^{\frac{2}{3}}.
\eea

From now on, in the rest of the paper, $\lesssim$ means bounded by a positive constant multiple, with this positive constant depending only on universal constants (such as constants arising from Sobolev embeddings, elliptic estimates,...) as well as the constants 
$$m,\,\, a, \,\, \dhor,\,\, \dred,\,\, \dbl, \,\, \dec,\,\,\de,\, r_0, \, \kl,$$
\textit{but not on} $\ep$ and $\ep_0$. {Also, note that the} constants $\dhor, \dred$ and $\dbl$ can be {chosen} to be only dependent on $m$ and $a$.

%%%%%%%%%%%%%%%%%%%%%%%%     
   
     \subsubsection{The spacetime $\MM$}
     \lab{section:SpacetimeMM-chap6}

%%%%%%%%%%%%%%%%%%%%%%%%

As in Section \ref{sec:Kerrpertbasic}, we consider a given vacuum spacetime $(\MM, \g)$ together with a null pair $(e_3, e_4)$ and its corresponding horizontal structure as in Section \ref{subsection:review-horiz.structures}. We will use the complexified Ricci and curvature coefficients  of Definition \ref{def:complexRicciandcurvaturecoefficients}. Moreover, we assume that $\MM$ is endowed with a pair of constants $(a, m)$, scalar functions $(\tau, r, \th, \tphi)$ and complex horizontal 1-forms $\Jk$, $\Jk_{\pm}$. 

In addition, we assume the following:
\begin{enumerate}
\item The level sets $\Si(\tau)$ of $\tau$ are spacelike, and $\tau\in[1,\tau_*]$ on $\MM$ for some arbitrary large constant $\tau_*$. 

\item The boundary of $\MM$ is given by 
\bea\lab{eq:prMM=AAcupSistarcupSi1cupSitaustar}
\pr\MM=\AA\cup\Si_*\cup\Si(1)\cup\Si(\tau_*)
\eea
where 
\bea\lab{eq:defAA=r=rplusoneminusedeh}
\AA:=\Big\{r=r_+(1-\deh), \, 1\leq\tau\leq\tau_*\Big\}, 
\eea
and $\Si_*$ is a spacelike hypersurface such that 
\bea\lab{eq:rangetauandronSigmastar}
1\leq\tau\leq\tau_*,\qquad r\geq r_*, \qquad \tau+r=c_*+h_*(r), \qquad h_*(r)=O(m^2r^{-1}),\qquad\textrm{on}\quad\Si_*,
\eea
with $r_*$, $c_*$ constants, and with $r_*$ satisfying the following dominance condition
\bea\lab{eq:dominantconditionforrstarcomparedtotaustartonS*}
r_*\simeq \ep_0^{-1}\tau_*^{1+\dec}.
\eea

\item Let $r_0$ a large enough fixed constant. We decompose $\MM$ as follows
\bea\lab{eq:definitionofMextandMint}
\Mint:=\MM\cap\{r\leq r_0\}, \qquad \Mext:=\MM\cap\{r\geq r_0\}. 
\eea
\end{enumerate} 

\begin{remark}
In view of \eqref{eq:rangetauandronSigmastar} and \eqref{eq:dominantconditionforrstarcomparedtotaustartonS*}, we have on $\Si_*$
\bea\lab{eq:controlofsizeofc*assimeqr*:onSi*}
c_*=\big(1+O(\ep_0)\big)r_*,
\eea
as well as 
\beaa
|r-r_*|\les |\tau-\tau_*|+\frac{m^2}{r_*}\les \tau_*+\frac{m^2}{r_*}\les \ep_0r_*+\frac{m^2}{r_*}\les \ep_0r_*\quad\Longrightarrow\quad r\simeq r_*
\eeaa
so that \eqref{eq:dominantconditionforrstarcomparedtotaustartonS*} in fact holds for any $r$ on $\Si_*$, i.e., 
\bea\lab{eq:dominantconditionforrstarcomparedtotaustartonS*:holdsforallronSi*}
r\simeq \ep_0^{-1}\tau_*^{1+\dec}\quad\textrm{on}\quad\Si_*.
\eea
\end{remark}

Finally, we introduce the following subregions of $\MM$
\bsub\lab{eq:defofsubregionsofMM}
 \begin{align}
 \MM(\tt_1,\tt_2):={}&\MM\cap\{\tt_1\leq \tt\leq \tt_2\}, \quad \forall\tau_1<\tau_2,\\
 \MM_{red}:={}&\MM\cap\{r\leq r_+(1+\dred)\},\\
 \Mtrap:={}&\MM_{r_+(1+2\dbl), 10m},\\ 
 \Mntrap:={}&\MM\setminus\Mtrap.
 \end{align}
 \esub

\begin{remark}\lab{rmk:choiceofmaingloblaframeinthewholepaper}
The global null frame $(e_3, e_4, e_1, e_2)$ will be chosen later as follows:
\begin{itemize}
\item In Section \ref{sec:proofofThmM1}, the global null frame $(e_3, e_4, e_1, e_2)$ is the one in Section 3.6.4 in \cite{KS:Kerr}.

\item In Section \ref{sec:proofofThmM2}, the global null frame $(e_3, e_4, e_1, e_2)$ is the one in Section 3.6.5 in \cite{KS:Kerr}. 

\item In Sections \ref{sec:higherordercurvatureestimatesforprovingThM8largea} \ref{sec:energyMorawetzforPc} \ref{sec:energyMorawetzforABBbAb}, the global null frame $(e_3, e_4, e_1, e_2)$ is the one in Section 9.6.1 in \cite{KS:Kerr}.
\end{itemize}
\end{remark}

\begin{remark}
The conditions \eqref{eq:rangetauandronSigmastar} and \eqref{eq:dominantconditionforrstarcomparedtotaustartonS*} on $\Si_*$ are consistent with (3.2.7) and (3.4.5) in \cite{KS:Kerr}.
\end{remark}

%%%%%%%%%%%%%%%%%%%%%%%%%%%%%%%%%

\subsubsection{A second global null frame on $\MM$}
\lab{sec:secondglobalnullframeonMM}

%%%%%%%%%%%%%%%%%%%%%%%%%%%%%%%%%

For the constructions in Sections \ref{sec:constructionofcoordinatessatisfyingassumptionsMaSz24} and \ref{sec:defandmainpropofregtripletOmiinKerrpert}, we will rely on the existence a second global null frame $(e_3', e_4', e_1', e_2')$ on $\MM$ with additional properties.

%%%%%%%%%%%%%%%%%%%%%%%%%%%%%%%%%

\paragraph{\textit{Transformation between two null frames}.}

%%%%%%%%%%%%%%%%%%%%%%%%%%%%%%%%%

In order to compare $(e_3', e_4', e_1', e_2')$ to the null frame of $\MM$, we rely on the following null frame transformations.

\begin{lemma}
\lab{Lemma:Generalframetransf}
The following transformation formulas hold true.
\begin{enumerate}
\item A general null transformation\footnote{In full generality, one could also rotate $e_1', e_2'$, but this would not change the horizontal structure and, as it  turns out, is in fact  not needed. The dot product and magnitude  $|\c |$ are taken with respect to the standard euclidian norm of $\RRR^2$.}  between two  null frames $(e_4, e_3, e_1, e_2)$ and $(e_4', e_3', e_1', e_2')$ on $\MM$
  can be written in   the form,
 \bea
 \lab{General-frametransformation}
 \bsplit
  e_4'&=\la\left(e_4 + f^b  e_b +\frac 1 4 |f|^2  e_3\right),\\
  e_a'&= \left(\de_a^b +\frac{1}{2}\fb_af^b\right) e_b +\frac 1 2  \fb_a  e_4 +\left(\frac 1 2 f_a +\frac{1}{8}|f|^2\fb_a\right)   e_3,\\
 e_3'&=\la^{-1}\left( \left(1+\frac{1}{2}f\c\fb  +\frac{1}{16} |f|^2  |\fb|^2\right) e_3 + \left(\fb^b+\frac 1 4 |\fb|^2f^b\right) e_b  + \frac 1 4 |\fb|^2 e_4 \right),
 \end{split}
 \eea
 where $\la$, $f_a= f^a , \fb_a= \fb^a $  are    scalar functions, called the transition coefficients of the change of frame. Note, in particular,
  \beaa
   e_a'= e_a +\frac 1 2  \fb_a \la^{-1} e_4'  +\frac 1 2 f_a e_3,\qquad  e_3' =\la^{-1}\left(  e_3 +  \fb^ae_a' -\frac 1 4 |\fb|^2\la^{-1} e_4'\right). 
 \eeaa

  \item The inverse transformation is given by the formulas
  \bea
 \lab{General-frametransformation'}
 \bsplit
  e_4&=\la'\left(e'_4 + f_b'  e'_b +\frac 1 4 |f'|^2  e'_3\right),\\
  e_a&= \left(\de_a^b  +\frac{1}{2}\fb'_a f'^b \right) e'_b +\frac 1 2  \fb'_a  e'_4 +\left(\frac 1 2 f'_a +\frac{1}{8}|f'|^2\fb'_a\right)   e'_3,\\
 e_3&=(\la')^{-1}\left( \left(1+\frac{1}{2}f'\c\fb'  +\frac{1}{16} |f'|^2  |\fb'|^2\right) e'_3 + \left(\fb'^b+\frac 1 4 |\fb'|^2f'^b\right) e'_b  + \frac 1 4 |\fb'|^2 e'_4 \right),
 \end{split}
 \eea
 where
\bea
\lab{relations:laffb-to-primes}
\bsplit
\la' &= \la^{-1} \left(1+\frac{1}{2}f\c\fb  +\frac{1}{16} |f|^2  |\fb|^2\right),\\
f_a'  &= -\frac{\la}{1+\frac{1}{2}f\c\fb  +\frac{1}{16} |f|^2  |\fb|^2}\left(f_a +\frac{1}{4}|f|^2\fb_a\right),\\
\fb_a' &= -\la^{-1}\left(\fb_a+\frac 1 4 |\fb|^2f_a\right).
\end{split}
\eea
Moreover
\bea
\lab{relations:ffbf'fb'}
\bsplit
\fb'_a f_b'=\fb_b f_a, \qquad  \la' |f'|^2 = \la |f|^2, \qquad 
(\la')^{-1}|\fb'|^2= \la^{-1} |\fb|^2.
\end{split}
\eea
 Denoting    $F=(f, \fb,  \la-1)$,   we also write, for small $|F|$,      
\bea
\lab{eq:f-f'transf-simplified}
\la' = \la^{-1} \left(1+\frac{1}{2}f\c\fb \right)+ O(|F|^3),\qquad f'_a=-\la f_a + O(|F|^3),\qquad \fb'_a=  -\la^{-1}\fb_a + O(|F|^3).
\eea
  \end{enumerate}
  \end{lemma}

\begin{proof}
See Lemma 2.2.1 in \cite{KS:Kerr}.
\end{proof}

%%%%%%%%%%%%%%%%%%%%%%%%%%%%%%%%%%%%%%%%%%%%%%

\paragraph{\textit{Definition and properties of the second global null frame on $\MM$}.}

%%%%%%%%%%%%%%%%%%%%%%%%%%%%%%%%%%%%%%%%%%%%%%

We denote by $(e_3', e_4', e_1', e_2')$ the second global null frame on $\MM$. It is obtained from the global null frame $(e_3, e_4, e_1, e_2)$ by the following null frame transformation 
 \bea
 \lab{frametransformation:fromfirsttosecondglobalframe}
 \bsplit
  e_4'&=e_4 + f^b  e_b +\frac 1 4 |f|^2  e_3,\\
  e_a'&= \left(\de_a^b +\frac{1}{2}\fb_af^b\right) e_b +\frac 1 2  \fb_a  e_4 +\left(\frac 1 2 f_a +\frac{1}{8}|f|^2\fb_a\right)   e_3,\\
 e_3'&=\left(1+\frac{1}{2}f\c\fb  +\frac{1}{16} |f|^2  |\fb|^2\right) e_3 + \left(\fb^b+\frac 1 4 |\fb|^2f^b\right) e_b  + \frac 1 4 |\fb|^2 e_4,
 \end{split}
 \eea
which corresponds to the transformation \eqref{General-frametransformation} with the choice $(f,\fb, \la=1)$, where the estimates satisfied by $(f, \fb)$ will be specified later. Then, we assume that the Ricci coefficients associated to $(e_4', e_3', e_1', e_2')$ satisfy 
\bea\lab{eq:additionalpropertiesseconddgolbalnullframe}
\xi'=0\quad\textrm{on}\quad\Mext\cap\{r\geq r_0+1\},
\eea
which is the additional property of the second global frame of $\MM$ that will be used\footnote{More precisely, the assumption \eqref{eq:additionalpropertiesseconddgolbalnullframe} is used to cancel the dangerous term $\xi_a'\nab_3$ in the commutator $[\nab_4', \nab_a']$ (for which $\xi'\in\Ga_g'$ would not suffice). This will be needed in Sections \ref{sec:constructionofcoordinatessatisfyingassumptionsMaSz24}  and \ref{sec:defandmainpropofregtripletOmiinKerrpert} to integrate various transport equations.} in Sections \ref{sec:constructionofcoordinatessatisfyingassumptionsMaSz24} and \ref{sec:defandmainpropofregtripletOmiinKerrpert}. 

\begin{remark}\lab{rmk:choiceofsecondgloblaframeinthewholepaper}
The second global null frame $(e_3', e_4', e_1', e_2')$ will be chosen as follows\footnote{Note that \eqref{eq:additionalpropertiesseconddgolbalnullframe} holds true both the global null frame of Section 3.6.5 in \cite{KS:Kerr} and the one of Section 9.6.1 in \cite{KS:Kerr}.}:
\begin{itemize}
\item In Sections \ref{sec:proofofThmM1} \ref{sec:proofofThmM2}, the second global null frame $(e_3', e_4', e_1', e_2')$ is the one of Section 3.6.5 in \cite{KS:Kerr}.

\item In Sections \ref{sec:higherordercurvatureestimatesforprovingThM8largea} \ref{sec:energyMorawetzforPc} \ref{sec:energyMorawetzforABBbAb}, the second global null frame $(e_3', e_4', e_1', e_2')$ is the one in Section 9.6.1 in \cite{KS:Kerr}.
\end{itemize}
In particular, in view of Remark \ref{rmk:choiceofmaingloblaframeinthewholepaper}, the global null frames $(e_3, e_4, e_1, e_2)$ and $(e_3', e_4', e_1', e_2')$ coincide in Sections \ref{sec:proofofThmM2} \ref{sec:higherordercurvatureestimatesforprovingThM8largea} \ref{sec:energyMorawetzforPc} \ref{sec:energyMorawetzforABBbAb},  and hence we have $f=\fb=0$ in that case for the null frame transformation \eqref{frametransformation:fromfirsttosecondglobalframe}. On the other hand, the null frame transformation \eqref{frametransformation:fromfirsttosecondglobalframe} is non trivial in Section \ref{sec:proofofThmM1} as  the global null frames $(e_3, e_4, e_1, e_2)$ and $(e_3', e_4', e_1', e_2')$ do not coincide in that case. 
\end{remark}

%%%%%%%%%%%%%%%%%%%%%%%%%%%%%%%%%

\subsubsection{Geometric set-up on $\Si_*$}
\lab{sec:geometricsetuponSigma*}

%%%%%%%%%%%%%%%%%%%%%%%%%%%%%%%%%

In this section, we provide the geometric set-up on $\Si_*$. This set-up will be needed for the constructions in Sections \ref{sec:constructionofcoordinatessatisfyingassumptionsMaSz24} and \ref{sec:defandmainpropofregtripletOmiinKerrpert}.

%%%%%%%%%%%%%%%%%%%%%%%%%%%%%%%%%

\paragraph{\textit{Integrable null frame on $\Si_*$}.}

%%%%%%%%%%%%%%%%%%%%%%%%%%%%%%%%%

We denote by $((e_*)_3, (e_*)_4, (e_*)_1, (e_*)_2)$ a null frame defined on $\Si_*$ and satisfying the following\footnote{Recall from \eqref{eq:rangetauandronSigmastar} that $\tau+r=c_*+h_*(r)$ on $\Si_*$ so that $(e_*)_a\in T\Si_*$ in fact follows from $(e_*)_a(r)=(e_*)_a(\tt)=0$.}
\bea\lab{eq:identitiesonSigma*fornullframeadaptedSigma*}
(e_*)_a\in T\Si_*, \quad a=1,2, \qquad (e_*)_a(r)=(e_*)_a(\tt)=0, \quad a=1,2,\quad\textrm{on}\quad\Si_*.
\eea
Note that \eqref{eq:identitiesonSigma*fornullframeadaptedSigma*} implies in particular that $((e_*)_3, (e_*)_4, (e_*)_1, (e_*)_2)$ is integrable, i.e.,
\bea\lab{eq:thenullframeadaptedtoSigma*isintegrable}
\atrch_*=\atrchb_*=0\quad\textrm{on}\quad\Si_*.
\eea
In addition, we impose the following transversality conditions on $\Si_*$ 
\bea\lab{eq:tranversalityconditionforthefoliationonSi*:chap5}
\xi_*=0, \qquad \om_*=0, \qquad \etab_*=-\ze_*, \qquad  (e_*)_4(r)=1, \qquad (e_*)_4(\tau)=0,
\eea
which allows us to make sense of all the Ricci coefficients in the frame of $\Si_*$.

We also introduce the vectorfield $\nu_*$ given by
\bea\lab{eq:definitionofvectorfieldnu*tangenttoSigma*inspanofe3*ande4*}
\nu_*=(e_*)_3+b_*(e_*)_4,
\eea
where the scalar function $b_*$ is uniquely chosen\footnote{We have $e_4(\tau+r-c_*-O(m^2r^{-1}))=1+O(mr^{-1})+\Ga_g$ which does not vanish for $r$ large enough, so that $e_4$ is transversal to $\Si_*$ in view of \eqref{eq:rangetauandronSigmastar}. This yields the existence and uniqueness of $b_*$ such that $\nu_*$ is tangent to $\Si_*$.} such that $\nu_*$ is tangent to $\Si_*$. In particular, $(\nu_*, (e_*)_1, (e_*)_2)$ forms an orthogonal basis of the tangent space of $\Si_*$. Also, we assume that the coordinates $(\th, \tphi)$ are propagated by $\nu_*$ along $\Si_*$, i.e.,
\bea\lab{eq:corrdinatesthandvarphiarepropagatedalongSi*bynu*}
\nu_*(\th)=0, \qquad \nu_*(\tphi)=0\quad\textrm{on}\quad \Si_*.
\eea

With respect to the null frame $((e_*)_3, (e_*)_4, (e_*)_1, (e_*)_2)$, we define the following linearized quantities 
\beaa
&&\widecheck{\trch}_* :=  \ds\trch_*-\frac{2}{r}, \qquad \widecheck{\trchb}_* := \ds\trchb_*+\frac{2\Up}{r},\qquad \ombc_* := \ds\omb_*-\frac{m}{r^2},\qquad \rhoc_* := \ds \rho_* +\frac{2m}{r^3},\\
\\
&&\widecheck{e_3(r)} := e_3(r) +\Up,\qquad \widecheck{e_3(\tau)} := e_3(\tau)-2,\qquad \widecheck{b_*}:= \ds b_*+1+\frac{2m}{r},
\eeaa
where $\Up: = 1-\frac{2m}{r}$, where we may linearize by Schwarzschild values in view of \eqref{eq:dominantconditionforrstarcomparedtotaustartonS*}. We then denote by $\Ga_{g}^*, \Ga_{b}^*$ the sets of  linearized quantities below
\bea\lab{eq:defintionofGagandGabforChapter5glsdfiuhgs}
\bsplit
\Ga_{g}^*&:=\Big\{ \widecheck{\trch}_*, \quad  \chih_*,  \quad   \ze_*,  \quad  \widecheck{\trchb}_*, \quad  r\a_*,\quad  r\b_*, \quad  r \rhoc_*, \quad  r \rhod_*\Big\},\\
\Ga_{b}^*&:=\Big\{ \eta_*,\quad  \chibh_*, \quad \widecheck{\omb}_*, \quad \xib_*,\quad r\bb_*,\quad \aa_*,  \quad  r^{-1}\widecheck{(e_*)_3(r)}, \quad r^{-1}\widecheck{(e_*)_3(\tau)},  \quad r^{-1} \widecheck{b_*}\Big\}.
\end{split}
\eea

%%%%%%%%%%%%%%%%%%%%%%%%%%%%%%%%%%%%%%%%%%%%%%%%%%%%%%%%%%%%%%

\paragraph{\textit{Uniformization on the last sphere $S_*$ of $\Si_*$ and 1-forms $(f_0, f_\pm,\Jk, \Jk_\pm)$ on $\Si_*$}.}

%%%%%%%%%%%%%%%%%%%%%%%%%%%%%%%%%%%%%%%%%%%%%%%%%%%%%%%%%%%%%%%

In view of \eqref{eq:rangetauandronSigmastar}, the last sphere $S_*$ of $\Si_*$ is defined by 
\beaa
S_*:=\{\tau=\tau_*\}\cap\{r=r_*\}, \qquad S_*\subset\Si_*.
\eeaa 
The coordinates $(\th, \vphi)$ provide coordinates on $S_*$ and we assume that the  induced metric $g_*$ on $S_*$  takes the form
 \bea\lab{eq:uniformizationforinducemetriconS*}
 g_*= r^2e^{2\phi_*}\Big( (d\th)^2+ \sin^2 \th (d\vphi)^2\Big),
 \eea
for some scalar function $\phi_*$ on $S_*$. For the constructions in Definition \ref{def:definitionoff0fplusfminus} below, we rely on a special orthonormal basis $((e_*)_1, (e_*)_2)$ of the tangent space of $S_*$ given by
\bea\lab{eq:specialorthonormalbasisofSstar}
(e_*)_1=\frac{1}{re^{\phi_*}}\pr_\th, \qquad (e_*)_2=\frac{1}{r\sin\th e^{\phi_*}}\pr_\vphi, \quad\textrm{on}\quad S_*.
\eea
We now introduce the following 1-forms, defined on $S_*$ and then extended to $\Si_*$.  

\begin{definition}\lab{def:definitionoff0fplusfminus}
Let $f_0$, $f_+$ and $f_-$ be the 1-forms defined on $S_*$ by:
\bea\lab{eq:definitionofthereal1formsf0fpmonS*}
\bsplit
&(f_0)_1 =0,  \qquad (f_0)_2 =\sin\th, \qquad (f_+)_1 =\cos\th\cos\vphi,  \qquad (f_+)_2 =-\sin\vphi,\\
& (f_-)_1 =\cos\th\sin\vphi,  \qquad (f_-)_2=\cos\vphi, \quad \textrm{on}\quad S_*,
\end{split}
\eea
in the orthonormal basis $((e_*)_1, (e_*)_2)$ of $S_*$ given by \eqref{eq:specialorthonormalbasisofSstar}, and extended to $\Si_*$ by:
\bea\lab{eq:extendionofthereal1formsf0fpmfromS*toSigma*}
(\nab_*)_{\nu_*} f_0=0, \qquad (\nab_*)_{\nu_*} f_+=0, \qquad (\nab_*)_{\nu_*} f_-=0. 
\eea
\end{definition}

We also introduce the following renormalization for angular derivatives of $f_0$ and $f_\pm$.
\begin{definition}\lab{def:renormalizationforf0fpfm}
We introduce the notations
\beaa
\widecheck{\curl_*(f_0)}:=\curl_*(f_0)-\frac{2}{r}\cos\th, \quad \widecheck{\div_*(f_\pm)}:=\div_*(f_\pm)+\frac{2}{r}J^{(\pm)}, \quad J^{(+)}:=x^1_p, \quad J^{(-)}:=x^2_p,
\eeaa
where we recall that $x^1_p=\sin\th\cos\tphi$ and $x^2_p=\sin\th\sin\tphi$, and 
\beaa
\widecheck{\nab_*\cos\th}:=\nab_*\cos\th+\frac{1}{r}\dual f_0, \qquad \widecheck{\nab_*J^{(+)}}:=\nab_*J^{(+)}-\frac{1}{r}f_+, \qquad \widecheck{\nab_*J^{(-)}}:=\nab_*J^{(-)}-\frac{1}{r}f_-.
\eeaa
\end{definition} 

Finally, the complex 1-forms $\Jk$ and $\Jk_{\pm}$ are defined as follows on $\Si_*$   
\bea\lab{eq:relationbetweenJkJkpmandf0onSigmastar}
\Jk = \frac{1}{|q|}\left(f_0+i\dual f_0\right), \qquad \Jk_\pm := \frac{1}{|q|}\left(f_\pm+i\dual f_\pm\right)\quad \textrm{on}\quad\Si_*.
\eea

%%%%%%%%%%%%%%%%%%%%%%%%%%%%%%%%%%%%%%%%%%%%%%%%%%%%%%%%%

\paragraph{\textit{Comparison of $(e_3', e_4', e_1', e_2')$ and $((e_*)_3, (e_*)_4, (e_*)_1, (e_*)_2)$ on $\Si_*$}.}

%%%%%%%%%%%%%%%%%%%%%%%%%%%%%%%%%%%%%%%%%%%%%%%%%%%%%%%%%

In order to compare on $\Si_*$ the null frame $((e_*)_3, (e_*)_4, (e_*)_1, (e_*)_2)$ to the global null frame $(e_3', e_4', e_1', e_2')$ introduced in \eqref{frametransformation:fromfirsttosecondglobalframe}, we rely on the following null frame transformations coefficients
\bea\lab{eq:changeofframecoefffromfromframeSigmastarttodoubleprimedframe}
\la_*:=\frac{\De}{|q|^2}, \qquad f_*:=\frac{a}{r}f_0, \qquad \fb_*:=\frac{a}{r}f_0+\widecheck{\fb_*},
\eea
where $f_0$ is given by Definition \ref{def:definitionoff0fplusfminus}, and where the estimates satisfied by the horizontal 1-form $\widecheck{\fb_*}$ will be provided later. On $\Si_*$, we then have the following transformation formula  between the global null frame $(e_3', e_4', e_1', e_2')$ on $\MM$ and the null frame $((e_*)_3, (e_*)_4, (e_*)_1, (e_*)_2)$ of $\Si_*$
\bea\lab{eq:changeofframefromfromframeSigmastarttodoubleprimedframe}
 \bsplit
   e_4' &=\la_*\left((e_*)_4 + f_*^b  (e_*)_b +\frac 1 4 |f_*|^2  (e_*)_3\right),\\
  e_a' &= \left(\de_a^b +\frac{1}{2}(\fb_*)_af_*^b\right)(e_*)_b +\frac 1 2  (\fb_*)_a  (e_*)_4 +\left(\frac 1 2 (f_*)_a +\frac{1}{8}|f_*|^2(\fb_*)_a\right)   (e_*)_3,\\
 e_3' &= \la_*^{-1}\left(\left(1+\frac{1}{2}f_*\c\fb_*  +\frac{1}{16} |f_*|^2  |\fb_*|^2\right) (e_*)_3 + \left(\fb_*^b+\frac 1 4 |\fb_*|^2f^b_*\right) (e_*)_b  + \frac 1 4 |\fb_*|^2 (e_*)_4\right),
 \end{split}
 \eea
 which corresponds to the transformation \eqref{General-frametransformation} with the choice $(f,\fb, \la)=(f_*, \fb_*, \la_*)$.

\begin{remark}
The geometric set-up on $\Si_*$ of this Section \ref{sec:geometricsetuponSigma*} is the one of Sections 3 and 5 of \cite{KS:Kerr}. More precisely, the setting in the above paragraph on the integrable null frame on $\Si_*$ is the one of Section 3.2.3 and Definition 3.3.2 in \cite{KS:Kerr}. Also, the setting in the above paragraph on uniformization on the last sphere $S_*$ of $\Si_*$ and on the 1-forms $(f_0, f_\pm,\Jk, \Jk_\pm)$ on $\Si_*$ is the one of (5.1.17) and (5.6.1)--(5.6.5) in \cite{KS:Kerr}. Finally, the change of frame \eqref{eq:changeofframecoefffromfromframeSigmastarttodoubleprimedframe} \eqref{eq:changeofframefromfromframeSigmastarttodoubleprimedframe} corresponds to (3.2.4)--(3.2.6) or (9.1.5) in \cite{KS:Kerr} taking the choice of the second global null frame $(e_3', e_4', e_1', e_2')$ discussed in Remark \ref{rmk:choiceofsecondgloblaframeinthewholepaper} into account\footnote{More precisely, if the second global null frame $(e_3', e_4', e_1', e_2')$ is the one of Section 3.6.5 in \cite{KS:Kerr} (which is the case in Sections \ref{sec:proofofThmM1} and \ref{sec:proofofThmM2}), then we are in the case (3.2.4)--(3.2.6) in \cite{KS:Kerr} and $\widecheck{\fb_*}=\frac{2m}{r}+\widecheck{(e_*)_3(r)}+O(m^2r^{-2})+O(r^{-2})\widecheck{b}_*$ in that case. On the other hand, if the second global null frame $(e_3', e_4', e_1', e_2')$ is the one of Section 9.6.1 in \cite{KS:Kerr} (which is the case in Sections \ref{sec:higherordercurvatureestimatesforprovingThM8largea} \ref{sec:energyMorawetzforPc} \ref{sec:energyMorawetzforABBbAb}), then we are in the case (9.1.5) in \cite{KS:Kerr} and $\widecheck{\fb_*}=0$ in that case.}.  
\end{remark}

%%%%%%%%%%%%%%%%%%%%%%%%%%%

   \subsubsection{Admissible perturbations of Kerr}
   \lab{sec:controlofGabandGagfromBA}

%%%%%%%%%%%%%%%%%%%%%%%%%%%

We assume that the global null frame $(e_3, e_4, e_1, e_2)$ satisfies the following estimates on $\MM$ 
\bea\lab{eq:assumptionsonMMforpartII}
\bsplit
r^3|\dk^{\leq k}\xi|+ r^2|\dk^{\leq k}\Ga_g|+r|\dk^{\leq k}\Ga_b| &\leq \ep, \qquad\qquad\,\,\, k\leq \kl,\\
r^2|\dk^{\leq k}\Ga_g|\leq \frac{\ep}{\tau^{\frac{1}{2}+\dec}}, \quad r|\dk^{\leq k}(\Ga_g, \Ga_b)| &\leq \frac{\ep}{\tau^{1+\dec}}, \qquad k\leq \frac{\kl}{2},
\end{split}
\eea
where the notation $(\Ga_b,\Ga_g)$ has been introduced in Definition \ref{definition.Ga_gGa_b}, where $\dk=\{\nab_3, r\nab_4, \dkb=r\nab \}$ denotes weighted derivatives, and where the small constant $\dec>0$ has been introduced in Section \ref{sec:smallnesconstants}. Moreover, we assume that the second global null frame $(e_3', e_4', e_1', e_2')$ satisfies the following estimates on $\MM$ 
\bea\lab{eq:assumptionsonMMforpartII:secondglobalframeauxassfor15derivatives}
\bsplit
 r|\dk^{\leq k}(f, \fb)|+r^2|\dk^{\leq k}\Ga_g'|+r^2|\dk^{\leq k}\widecheck{e_4'(x^b_p)}|+r^3|\dk^{\leq k}\widecheck{\nab_4'\Jk_{\pm}}|&\leq \frac{\ep}{\tau^{\frac{1}{2}+\dec}}, \quad k\leq 16,\\
 |\dk^{\leq k}(f, \fb)|+r|\dk^{\leq k}(\Ga_g', \Ga_b')|+r|\dk^{\leq k}(\widecheck{e_4'(x^b_p)}, \widecheck{e_3'(x^b_p)}, \widecheck{\nab'(x^b_p)})|\\
 +r^2|\dk^{\leq k}(\widecheck{\nab_4'\Jk_{\pm}}, \widecheck{\nab_3'\Jk_{\pm}}, \widecheck{\nab'\Jk_{\pm}})|&\leq \frac{\ep}{\tau^{1+\dec}}, \quad k\leq 16,
\end{split}
\eea 
where $(\Ga_b', \Ga_g')$ denotes the analog of\footnote{The Kerr values that are subtracted to obtain the linearized quantities appearing in $(\Ga_b', \Ga_g')$ are defined as in Section \ref{sec:definitionoflinearizedquantities:chap4}, i.e., still using $(\tau, r, x^1, x^2)$, $\Jk$ and $\Jk_\pm$.} $(\Ga_b,\Ga_g)$ for the global null frame $(e_3', e_4', e_1', e_2')$, and where the change of frame coefficients $(f, \fb)$ have been introduced in \eqref{frametransformation:fromfirsttosecondglobalframe}. Finally, we have the following estimates on the spacelike hypersurface $\Si_*$ 
\bea\lab{eq:assumptionsonSigmastarforpartII}
\bsplit
r^2\left|\dk_*^k\left[\div_*(f_0),\, \widecheck{\curl_*(f_0)},\, \nab_*\hot f_0,\, \nab_*f_0 - \frac{1}{r}\cos\th\in\right]\right| \\
+ r^2\left|\dk_*^k\left[\widecheck{\div_*(f_\pm)},\, \curl_*(f_\pm),\, \nab_*\hot f_\pm,\, \nab_*f_\pm + \frac{1}{r}J^{(\pm)}\de\right]\right|\\ 
+r^2|\dk_*^k\Ga_g^*|+r^2|\dk_*^k\widecheck{\nab_*\cos\th}|+r^2|\dk_*^k\widecheck{\nab_*J^{(\pm)}}|+r|\dk_*^k\widecheck{\fb_*}| &\leq \frac{\ep}{\tau^{\frac{1}{2}+\dec}}, \qquad k\leq 15,\\
r|\dk_*^k\Ga_b^*| &\leq \frac{\ep}{\tau^{1+\dec}}, \qquad k\leq 15,\\
\int_{1}^{\tau_*}\tau^{2+2\dec}\big(\|r\dk_*^k\eta_*\|^2_{L^\infty(S_*(\tau))}+\|\dk_*^k\widecheck{(e_*)_3(\tau)}\|^2_{L^\infty(S_*(\tau))}\big)d\tau&\leq\ep^2, \qquad k\leq 15,
\end{split}
\eea
where $\dk_*=\{(\nab_*)_\nu, \dkb_*=r\nab_*\}$ denotes weighted derivatives on $\Si_*$, where $S_*(\tau_0)=\Si_*\cap\{\tau=\tau_0\}$, and where the various quantities appearing in \eqref{eq:assumptionsonSigmastarforpartII} have been introduced in Section \ref{sec:geometricsetuponSigma*}. Additionally, the scalar function $\phi_*$ on $S_*$ appearing in \eqref{eq:uniformizationforinducemetriconS*} satisfies 
\bea\lab{eq:controloftheconformaluniformizationfactorofinducedmetricS*}
|\dkb_*^{\leq 15}\phi_*|\leq\frac{\ep}{r\tau^{\frac{1}{2}+\dec}}.
\eea

\begin{remark}
In this section $k_L$  is an unspecified  large positive  integer satisfying in particular $\kl\geq 32$.
\end{remark}

\begin{remark}
The bounds \eqref{eq:assumptionsonSigmastarforpartII} on $\Si_*$ and \eqref{eq:assumptionsonMMforpartII:secondglobalframeauxassfor15derivatives} on $\MM$ will be used in Sections \ref{sec:constructionofcoordinatessatisfyingassumptionsMaSz24} and \ref{sec:defandmainpropofregtripletOmiinKerrpert} to initialize various quantities on $\Si_*$ and then transport them inside $\MM$. In turn, the constructions in Sections \ref{sec:constructionofcoordinatessatisfyingassumptionsMaSz24} and \ref{sec:defandmainpropofregtripletOmiinKerrpert} will be used in Section \ref{sec:energyMorawetzesitmatesforTeukoslkyonMM:upto15derivatives} to apply the results in \cite{MaSz24} \cite{MaSz26}.
\end{remark}

\begin{remark}
Recall the choices of the global null frames $(e_3, e_4, e_1, e_2)$ and $(e_3', e_4', e_1', e_2')$ discussed in Remarks \ref{rmk:choiceofmaingloblaframeinthewholepaper} and \ref{rmk:choiceofsecondgloblaframeinthewholepaper}. Assumptions \eqref{eq:assumptionsonMMforpartII} hold for the global null frames on $\MM$ of Sections 3.6.4 and 3.6.5 in \cite{KS:Kerr}, respectively used in Sections \ref{sec:proofofThmM1} and \ref{sec:proofofThmM2}, and the one of Lemma 9.6.5 in \cite{KS:Kerr} used in Sections \ref{sec:higherordercurvatureestimatesforprovingThM8largea} \ref{sec:energyMorawetzforPc} \ref{sec:energyMorawetzforABBbAb}. Also, assumptions \eqref{eq:assumptionsonMMforpartII:secondglobalframeauxassfor15derivatives} and \eqref{eq:additionalpropertiesseconddgolbalnullframe} hold for the global null frame of Section 3.6.5 in \cite{KS:Kerr} and for the one in Section 9.6.1 in \cite{KS:Kerr}. For $(f_0, f_+, f_-)$, $\widecheck{\nab_*J^{(\pm)}}$ and $\widecheck{b_*}$, the bounds  \eqref{eq:assumptionsonSigmastarforpartII} on $\Si_*$ hold  in view of Proposition 5.6.4 in \cite{KS:Kerr}, while for $(\widecheck{f_*}, \widecheck{\fb_*})$, they hold in view of the properties of the global null frames on $\MM$ of Sections 3.6.4 and 3.6.5 in \cite{KS:Kerr}, and for $\eta_*$ and $\widecheck{(e_*)_3(\tau)}$, they hold in view of Proposition 5.3.1 in \cite{KS:Kerr}. Finally, the bound \eqref{eq:controloftheconformaluniformizationfactorofinducedmetricS*} for $\phi_*$ follows from effective uniformization, see Corollary 5.6.5 in \cite{KS:Kerr}.
\end{remark}
 
\begin{remark}\lab{rmk:xiisactuallybetterthanGag}
Note that the assumptions for $\xi$ in \eqref{eq:assumptionsonMMforpartII} {would follow from} $\xi\in r^{-1}\Ga_g$, while $\xi$ is a priori only in $\Ga_g$ according to Definition \ref{definition.Ga_gGa_b}. These  stronger assumptions \eqref{eq:assumptionsonMMforpartII} for $\xi$ will always hold:
\begin{itemize}
\item In Section \ref{sec:proofofThmM1}, this follows from the assumptions \eqref{eq:GlobalFrame-HcinGa_g2} and the fact that $\xi\in\Ga_g$. 
\item In Section \ref{sec:proofofThmM2}, this follows from the assumptions \eqref{eq:Xi-Hb-chapter12}. 
\item In Sections \ref{sec:higherordercurvatureestimatesforprovingThM8largea} \ref{sec:energyMorawetzforPc} \ref{sec:energyMorawetzforABBbAb}, this follows from the assumptions \eqref{eq:specialidentityforthegloablframeofMMinpartIII}.
\end{itemize}
\end{remark}

We derive the following simple consequence of \eqref{eq:assumptionsonMMforpartII}.
\begin{lemma}\lab{lemma:identitymodqsquaretensornormReJksquareminussinthsquareovermodqsquareuptoerrorterm}
The following estimates hold on $\MM$
\bea\lab{eq:thestatedestimateformodqsquaretensornormReJksquareminussinthsquareovermodqsquareonMMgeq3m}
\begin{split}
&\left|\dk^{\leq \frac{\kl}{2}}\left(|\Re(\Jk)|^2 - \frac{(\sin\th)^2}{|q|^2}\right)\right|\les\frac{\ep}{r^3\tau^{\frac{1}{2}+\dec}}, \quad\,\,\,\, \left|\dk^{\leq \frac{\kl}{2}}\left(|\Re(\Jk)|^2 - \frac{(\sin\th)^2}{|q|^2}\right)\right|\les\frac{\ep}{r^2\tau^{1+\frac{3\dec}{4}}},\\
&\left|\dk^{\leq\kl}\left(|\Re(\Jk)|^2 - \frac{(\sin\th)^2}{|q|^2}\right)\right|\les\frac{\ep}{r^3}.
\end{split}
\eea
\end{lemma}

\begin{proof}
First, note from \eqref{eq:definitionofthereal1formsf0fpmonS*} that $|f_0|^2=(\sin\th)^2$ on $S_*$ which together with \eqref{eq:extendionofthereal1formsf0fpmfromS*toSigma*} and \eqref{eq:corrdinatesthandvarphiarepropagatedalongSi*bynu*} implies that $|f_0|^2=(\sin\th)^2$ on $\Si_*$. In view of \eqref{eq:relationbetweenJkJkpmandf0onSigmastar}, this yields
\bea\lab{eq:RealpartofJkdotRealpartofJkisequaltosinthsquaredividedbymodofqsquareonSi*}
|\Re(\Jk)|^2 &=& \frac{|f_0|^2}{|q|^2}=\frac{(\sin\th)^2}{|q|^2}\quad\textrm{on}\quad\Si_*.
\eea
Then, we compute
\beaa
\nab_4\left(|\Re(\Jk)|^2 - \frac{(\sin\th)^2}{|q|^2}\right) &=& 2\Re(\Jk)\c\Re(\nab_4\Jk) +  \frac{\cos\th\nab_4(\cos\th)}{|q|^2} +\frac{(\sin\th)^2}{|q|^4}\nab_4(|q|^2)\\
&=& 2\Re(\Jk)\c\Re\left(-\frac{\De\ov{q}}{|q|^4}\Jk\right)+\frac{2r\De(\sin\th)^2}{|q|^6}+r^{-2}\Ga_g\\
&=& -\frac{2r\De}{|q|^4}\left(|\Re(\Jk)|^2 - \frac{(\sin\th)^2}{|q|^2}\right)+r^{-2}\Ga_g
\eeaa
so that 
\bea\lab{eq:transporteqfornab4h0withh0=modqsquaretensornormReJksquareminussinthsquareovermodqsquare}
\nab_4h_0 = \Ga_g, \qquad h_0:=|q|^2\left(|\Re(\Jk)|^2 - \frac{(\sin\th)^2}{|q|^2}\right).
\eea
Note that \eqref{eq:RealpartofJkdotRealpartofJkisequaltosinthsquaredividedbymodofqsquareonSi*} and \eqref{eq:transporteqfornab4h0withh0=modqsquaretensornormReJksquareminussinthsquareovermodqsquare} imply 
\bea\lab{eq:controlofweightedderivativesofh0onSi*withh0=modqsquaretensornormReJksquareminussinthsquareovermodqsquare}
|\dk^{\leq\frac{\kl}{2}}h_0|\les \frac{\ep}{r\tau^{\frac{1}{2}+\dec}}\quad\textrm{on}\quad\Si_*.
\eea

Next, we assume the following bootstrap assumptions for the scalar function $h_0$
\bea\lab{eq:bootstrapasumptionforh0withh0=modqsquaretensornormReJksquareminussinthsquareovermodqsquare}
|\dk^{\leq \frac{\kl}{2}}h_0|\leq\frac{\sqrt{\ep}}{r\tau^{\frac{1}{2}+\dec}}, \qquad |\dk^{\leq \frac{\kl}{2}}h_0|\leq\frac{\sqrt{\ep}}{\tau^{1+\frac{3\dec}{4}}}, \quad\textrm{on}\quad\MM\cap\{r\geq r_1\}
\eea
for some $r_1\geq 3m$. Relying on the commutation formulas \eqref{commutatorbetweenLieTLieZandnabnab4nab3:1} and in Lemma \ref{LEMMA:COMM-GEN-B}, we infer, for $k_1+k_2\leq\frac{\kl}{2}$, 
\beaa
\nab_4\left((r\nab)^{k_1}\Lieb_\T^{k_2}h_0\right) &=& O(r^{-2})(r\nab)^{\leq k_1}\Lieb_\T^{k_2}h_0+\dk^{k_1+k_2}\Ga_g+\dk^{\leq k_1+k_2}(\Ga_g\c h_0),
\eeaa
where $\T$ is given by \eqref{eq:definitionofTandPhithataretheapproximateKillingvectorifeldinKerrpert}. Integrating from $\Si_*$ where \eqref{eq:controlofweightedderivativesofh0onSi*withh0=modqsquaretensornormReJksquareminussinthsquareovermodqsquare} holds, and using \eqref{eq:assumptionsonMMforpartII} and \eqref{eq:bootstrapasumptionforh0withh0=modqsquaretensornormReJksquareminussinthsquareovermodqsquare}, we infer on $\MM\cap\{r\geq r_1\}$, for $k_1+k_2\leq\frac{\kl}{2}$, 
\beaa
|(r\nab)^{k_1}\Lieb_\T^{k_2}h_0| &\les& \frac{\ep}{r_*\tau^{\frac{1}{2}+\dec}}+\int_r^{r_*}\frac{\ep}{{r'}^2\tau^{\frac{1}{2}+\dec}}dr'\les\frac{\ep}{r\tau^{\frac{1}{2}+\dec}},\\
|(r\nab)^{k_1}\Lieb_\T^{k_2}h_0| &\les& \frac{\ep}{r_*\tau^{\frac{1}{2}+\dec}}+\int_r^{r_*}\left(\frac{\ep}{{r'}^2\tau^{\frac{1}{2}+\dec}}\right)^{\frac{\dec}{2}}\left(\frac{\ep}{r'\tau^{1+\dec}}\right)^{1-\frac{\dec}{2}}dr'\les\frac{\ep}{\tau^{1+\frac{3\dec}{4}}},
\eeaa
where we used also \eqref{eq:dominantconditionforrstarcomparedtotaustartonS*} in the last inequality. Together with Lemma \ref{lemma:basicpropertiesLiebTfasdiuhakdisug:chap9},  \eqref{eq:transporteqfornab4h0withh0=modqsquaretensornormReJksquareminussinthsquareovermodqsquare} and \eqref{eq:assumptionsonMMforpartII}, this yields on $\MM\cap\{r\geq r_1\}$, for $k\leq\frac{\kl}{2}$,  
\beaa
|(r\nab_4, \nab_{\T}, r\nab)^kh_0| \les \frac{\ep}{r\tau^{\frac{1}{2}+\dec}},\qquad |(r\nab_4, \nab_{\T}, r\nab)^kh_0| \les \frac{\ep}{\tau^{1+\frac{3\dec}{4}}}.
\eeaa
Together with the fact that $\nab_3$ is generated by $(\nab_{\T}, \nab_4, \nab)$, we deduce
\beaa
|\dk^{\leq \frac{\kl}{2}}h_0|\les\frac{\ep}{r\tau^{\frac{1}{2}+\dec}}, \qquad |\dk^{\leq \frac{\kl}{2}}h_0|\les\frac{\ep}{\tau^{1+\frac{3\dec}{4}}}, \quad\textrm{on}\quad\MM\cap\{r\geq r_1\}
\eeaa
which improves the bootstrap assumptions \eqref{eq:bootstrapasumptionforh0withh0=modqsquaretensornormReJksquareminussinthsquareovermodqsquare}. We deduce that $r_1=3m$ and that, in view of the definition of $h_0$ in \eqref{eq:transporteqfornab4h0withh0=modqsquaretensornormReJksquareminussinthsquareovermodqsquare}, 
\bea\lab{eq:provesthestatedestimateformodqsquaretensornormReJksquareminussinthsquareovermodqsquareonMMgeq3m}
\left|\dk^{\leq \frac{\kl}{2}}\left(|\Re(\Jk)|^2 - \frac{(\sin\th)^2}{|q|^2}\right)\right|\les\frac{\ep}{r^3\tau^{\frac{1}{2}+\dec}}, \quad\,\,\,\, \left|\dk^{\leq \frac{\kl}{2}}\left(|\Re(\Jk)|^2 - \frac{(\sin\th)^2}{|q|^2}\right)\right|\les\frac{\ep}{r^2\tau^{1+\frac{3\dec}{4}}},
\eea
on $\MM\cap\{r\geq 3m\}$. 

The estimate \eqref{eq:provesthestatedestimateformodqsquaretensornormReJksquareminussinthsquareovermodqsquareonMMgeq3m} proves the first two estimates in \eqref{eq:thestatedestimateformodqsquaretensornormReJksquareminussinthsquareovermodqsquareonMMgeq3m} on $\MM\cap\{r\geq 3m\}$. To extend these estimates to $\MM\cap\{r\leq 3m\}$, we compute
\beaa
\nab_3\left(|\Re(\Jk)|^2 - \frac{(\sin\th)^2}{|q|^2}\right) &=& 2\Re(\Jk)\c\Re(\nab_3\Jk) +  \frac{\cos\th\nab_3(\cos\th)}{|q|^2} +\frac{(\sin\th)^2}{|q|^4}\nab_3(|q|^2)\\
&=& 2\Re(\Jk)\c\Re\left(\frac{1}{\ov{q}}\Jk\right)-\frac{2r(\sin\th)^2}{|q|^4}+r^{-2}\Ga_b\\
&=& \frac{2r}{|q|^2}\left(|\Re(\Jk)|^2 - \frac{(\sin\th)^2}{|q|^2}\right)+r^{-2}\Ga_b.
\eeaa
Using this transport equation, we easily extend \eqref{eq:provesthestatedestimateformodqsquaretensornormReJksquareminussinthsquareovermodqsquareonMMgeq3m} from $r=3m$ to $\MM\cap\{r\leq 3m\}$ which concludes the proof of the first two estimates in \eqref{eq:thestatedestimateformodqsquaretensornormReJksquareminussinthsquareovermodqsquareonMMgeq3m}. Finally, the proof of the last estimate in \eqref{eq:thestatedestimateformodqsquaretensornormReJksquareminussinthsquareovermodqsquareonMMgeq3m} is very similar and left to the reader. This concludes the proof of Lemma \ref{lemma:identitymodqsquaretensornormReJksquareminussinthsquareovermodqsquareuptoerrorterm}.
\end{proof}

%%%%%%%%%%%%%%%%%%%%%%%%%%%%%%%%%%%%%%

\subsection{Suitable coordinates system and complex 1-forms on $\MM$} 
\lab{sec:constructionofcoordinatessatisfyingassumptionsMaSz24} 
  
%%%%%%%%%%%%%%%%%%%%%%%%%%%%%%%%%%%%%%

In this section, we exhibit a coordinates system and complex 1-forms satisfying the assumptions needed to apply the energy-Morawetz estimates in perturbations of Kerr of \cite{MaSz24} \cite{MaSz26}, respectively for scalar waves and Teukolsky equations, in Section \ref{sec:energyMorawetzesitmatesforTeukoslkyonMM:upto15derivatives}.

%%%%%%%%%%%%%%%%%%%%%%%%%%%%%%%%%%%%%%%%%%%%%%

\subsubsection{Construction of the coordinate $\ttt$ on $\MM$}
\lab{sec:constructionofthecoordinatetttonMM}

%%%%%%%%%%%%%%%%%%%%%%%%%%%%%%%%%%%%%%%%%%%%%%

In view of Definition \ref{definition.Ga_gGa_b}, we have $\widecheck{e_3(\tau)}=r\Ga_b$ and $\widecheck{\nab(x^b)}\in\Ga_b$, $b=1,2$. In turn, this implies that the corresponding linearized inverse metric coefficients for the coordinates $(\tau, r, x^1, x^2)$ do not satisfy the assumptions of \cite{MaSz24} \cite{MaSz26}. In this section, we construct a new coordinates $(\ttt, \tx^1, \tx^2)$ that satisfy improved estimates. To this end, we first define the following vectorfield on $\MM$.

\begin{definition}[Vectorfield $\Ytau$]
\lab{def:definitionofthevectorfildYtau}
The vectorfield $\Ytau$ is defined w.r.t. the second global frame $(e_3', e_4', e_1', e_2')$ of $\MM$ introduced in Section \ref{sec:secondglobalnullframeonMM} as follows 
\beaa
\Ytau:=\Ytau^3e_3'+\Ytau^4e_4',
\eeaa
where $\Ytau^3$ and $\Ytau^4$ are defined on $\MM$ by
\beaa
\Ytau^4:=\chi(r)\frac{e_3'(\tau)}{e_4'(\tau)} +(1-\chi(r))\frac{|q|^2}{\De}, \qquad \Ytau^3:=-\chi(r)-(1-\chi(r))\frac{e_4'(\tau)}{e_3'(\tau)}\frac{|q|^2}{\De},
\eeaa
with $\chi(r)$ a smooth cut-off function such that $\chi(r)=1$ for $r\leq r_+(1+\dbl)$ and $\chi(r)=0$ for $r\geq 13m$. 
\end{definition}

\begin{remark}\lab{rmk:YtauistangenttolevelhypersurfacesSioftauonMM}
Since $\Ytau^3$ and $\Ytau^4$ are chosen such that $e_3'(\tau)\Ytau^3=-e_4'(\tau)\Ytau^4$, we have $\Ytau(\tau)=0$ so that the vectorfield $\Ytau$ is tangent to the level hypersurface $\Si(\tau)$ of $\tau$. In particular, $\MM$ is covered by the integral curves of $\Ytau$ initialized on $\Si_*$. 
\end{remark}

\begin{remark}\lab{rmk:formofYautforrgeqr0intermsofdoubleprimedframe}
The constant $r_0$ is chosen large enough and satisfies in particular $r_0\geq 13m$. Thus, in the view the above definition of $\Ytau$, we have in $\Mext\cap\{r\geq r_0+1\}$
\beaa
\Ytau=\frac{|q|^2}{\De}e_4' -\frac{|q|^2}{\De}\frac{e_4'(\tau)}{e_3'(\tau)}e_3'=\frac{|q|^2}{\De}e_4' +\big(O(m^2r^{-2})+\Ga_g'\big)e_3',
\eeaa
where the global null frame $(e_3', e_4', e_1', e_2')$ has been introduced in \eqref{frametransformation:fromfirsttosecondglobalframe}.
\end{remark}

We start with the construction of the coordinate $\ttt$.
\begin{lemma}\lab{lemma:controloftildeuonMext}
There exists a scalar function $\ttt$ on $\MM$ such that 
\bea\lab{eq:controloftildeuonMext:ealphaoftttcheckinrGag}
r|\dk^{\leq 15}(\widecheck{e_4(\ttt)}, \widecheck{e_3(\ttt)}, \widecheck{\nab(\ttt)})|\les\frac{\ep}{\tau^{\frac{1}{2}+\dec}}, \qquad |\dk^{\leq 15}(\widecheck{e_4(\ttt)}, \widecheck{e_3(\ttt)}, \widecheck{\nab(\ttt)})|\les\frac{\ep}{\tau^{1+\dec}},
\eea
and
\bea\lab{eq:additionestimatewidechecke4tttbyeprminus2}
|\dk^{\leq 15}\widecheck{e_4(\ttt)}|\les\frac{\ep}{r^2}, 
\eea
where we have used above the following linearized quantities 
\beaa
\widecheck{e_3(\ttt)}:= e_3(\ttt)-\tmod'(r), \quad \widecheck{e_4(\ttt)}:=e_4(\ttt)-\frac{2(r^2+a^2) - \De\tmod'(r)}{|q|^2},\quad \widecheck{\DD(\ttt)}:=& \DD(\ttt)-a\Jk.
\eeaa
Also, we have
\bea\lab{eq:additionestimatetttminust}
|\ttt-\tt|\les \frac{\ep}{\tt^{\frac{1}{2}+\de_{dec}}}\quad\textrm{on}\quad\MM\cap\left\{r\geq \frac{1}{2}\tau^{\frac{1}{2}}\right\}, \qquad  \ttt=\tt\quad\textrm{on}\quad\MM\cap\left\{r\leq \frac{1}{2}\tau^{\frac{1}{2}}\right\}.
\eea
\end{lemma}

\begin{remark}
The point of introducing $\ttt$ is that $\widecheck{e_3(\ttt)}\in r\Ga_g$ while we have only $\widecheck{e_3(\tt)}\in r\Ga_b$. 
\end{remark}

\begin{proof}
We first introduce an auxiliary scalar function $\widehat{\tau}$ defined on $\Mext\cap\{r\geq r_0+1\}$ as follows 
\bea\lab{definitionofthenewcoordinatetildeu}
\widehat{\tau}=\tt_*\quad\textrm{on}\quad S_*, \qquad \nu_*(\widehat{\tau})=2\quad\textrm{on}\quad \Si_*, \qquad \Ytau(\widehat{\tau})=0\quad\textrm{on}\quad \Mext\cap\{r\geq r_0+1\},
\eea
and we start with the control of $\widehat{\tau}$ on $\Si_*$ where we assume the following bootstrap assumption\footnote{The bootstrap assumption for $(e_*)_4(\widehat{\tau})$ in \eqref{eq:bootassfordk*15widecheckealphaoftttonSi*} is purposely non-sharp.} 
\bea\lab{eq:bootassfordk*15widecheckealphaoftttonSi*}
r^{\frac{3}{2}}|\dk_*^{\leq 15}(e_*)_4(\widehat{\tau})|+r|\dk_*^{\leq 15}(\widecheck{(e_*)_3(\widehat{\tau})}, \nab_*(\widehat{\tau}))|\leq\frac{\sqrt{\ep}}{\tau^{\frac{1}{2}+\dec}}\quad\textrm{on}\quad\Si_*\cap\{\tt_{**}\leq\tau\leq\tt_*\},
\eea
for some $1\leq\tt_{**}<\tt_*$, with the linearized quantity $\widecheck{(e_*)_3(\widehat{\tau})}$ defined as follows
\beaa
\widecheck{(e_*)_3(\widehat{\tau})}:=(e_*)_3(\widehat{\tau})-2.
\eeaa

We now improve the bootstrap assumptions \eqref{eq:bootassfordk*15widecheckealphaoftttonSi*}. First, in view of \eqref{definitionofthenewcoordinatetildeu} and \eqref{eq:definitionofvectorfieldnu*tangenttoSigma*inspanofe3*ande4*}, we have on $\Si_*$
\beaa
2=\nu_*(\widehat{\tau})=(e_*)_3(\widehat{\tau})+b_*(e_*)_4(\widehat{\tau})\quad\Longrightarrow\quad \widecheck{(e_*)_3(\widehat{\tau})}=-b_*(e_*)_4(\widehat{\tau}),
\eeaa
which, together with \eqref{eq:bootassfordk*15widecheckealphaoftttonSi*}, \eqref{eq:dominantconditionforrstarcomparedtotaustartonS*} and the control of $b_*$ provided by \eqref{eq:assumptionsonSigmastarforpartII} yields 
\bea\lab{eq:improvementofeq:bootassfordk*15widecheckealphaoftttonSi*:e3}
r|\dk_*^{\leq 15}\widecheck{(e_*)_3(\widehat{\tau})}| \les r|\dk_*^{\leq 15}(e_*)_4(\widehat{\tau})| \les\frac{\sqrt{\ep}}{r^{\frac{1}{2}}\tau^{\frac{1}{2}+\dec}}\les\frac{\ep}{\tau^{\frac{1}{2}+\dec}} \quad\textrm{on}\quad\Si_*\cap\{\tt_{**}\leq\tau\leq\tt_*\},
\eea
which improves the bootstrap assumptions \eqref{eq:bootassfordk*15widecheckealphaoftttonSi*} for $\widecheck{(e_*)_3(\widehat{\tau})}$. 

Next, we improve the non-sharp bootstrap assumptions \eqref{eq:bootassfordk*15widecheckealphaoftttonSi*} for $(e_*)_4(\widehat{\tau})$. In view of the change of frame formula \eqref{eq:changeofframefromfromframeSigmastarttodoubleprimedframe}, \eqref{eq:changeofframecoefffromfromframeSigmastarttodoubleprimedframe}, we have
\beaa
 \bsplit
\big(1+O(mr^{-1})\big)e_4'(\widehat{\tau}) &= (e_*)_4(\widehat{\tau}) + O(r^{-1})\nab_*(\widehat{\tau}) +O(r^{-2})(e_*)_3(\widehat{\tau}),\\
r|\dk_*^{\leq 15}e_3'(\widehat{\tau})| &\les 1 \quad\textrm{on}\quad\Si_*\cap\{\tt_{**}\leq\tau\leq\tt_*\},
 \end{split}
 \eeaa
where we used the bootstrap assumptions \eqref{eq:bootassfordk*15widecheckealphaoftttonSi*} and  the control of $\widecheck{\fb_*}$ in  \eqref{eq:assumptionsonSigmastarforpartII} to get the non-sharp bound on the second line. On the other hand, we have $\Ytau(\widehat{\tau})=0$ in view of \eqref{definitionofthenewcoordinatetildeu} which together with Remark \ref{rmk:formofYautforrgeqr0intermsofdoubleprimedframe}, \eqref{eq:dominantconditionforrstarcomparedtotaustartonS*} and \eqref{eq:assumptionsonMMforpartII:secondglobalframeauxassfor15derivatives} implies, on $\Si_*\cap\{\tt_{**}\leq\tau\leq\tt_*\}$, 
\beaa
\left|\dk_*^{\leq 15}e_4'(\widehat{\tau})\right| \les r^{-2}\left|\dk_*^{\leq 15}e_3'(\widehat{\tau})\right|\les r^{-2}\les \sqrt{\ep_0}r^{-\frac{3}{2}},
\eeaa
where we used the above non-sharp bound for $e_3'(\widehat{\tau})$. Plugging in the above identity relating $e_4'(\widehat{\tau})$ and $(e_*)_4(\widehat{\tau})$, and using  the bootstrap assumptions \eqref{eq:bootassfordk*15widecheckealphaoftttonSi*} and \eqref{eq:dominantconditionforrstarcomparedtotaustartonS*}, we infer 
\bea\lab{eq:improvementofeq:bootassfordk*15widecheckealphaoftttonSi*:e4}
r^{\frac{3}{2}}|\dk_*^{\leq 15}(e_*)_4(\widehat{\tau})| \les r^{\frac{3}{2}}|\dk_*^{\leq 15}e_4'(\widehat{\tau})|+r^{-2}\les \sqrt{\ep_0}r^{-\frac{3}{2}} \quad\textrm{on}\quad\Si_*\cap\{\tt_{**}\leq\tau\leq\tt_*\},
\eea
which improves the non-sharp bootstrap assumptions \eqref{eq:bootassfordk*15widecheckealphaoftttonSi*} for $(e_*)_4(\widehat{\tau})$. 

Next, we improve the bootstrap assumptions \eqref{eq:bootassfordk*15widecheckealphaoftttonSi*} for $\nab_*(\widehat{\tau})$. In view of \eqref{eq:comm-nab3-nab4-naba-f-general} and \eqref{eq:thenullframeadaptedtoSigma*isintegrable}, we have on $\Si_*$
\beaa
       \begin{split}
        \,[(e_*)_3, (e_*)_a]\widehat{\tau} &=2(\eta_*)_a +\Ga_g^*+ \big(O(r^{-1})+\Ga_b^*)\nab_*(\widehat{\tau})+\Ga_b^*\widecheck{(e_*)_3(\widehat{\tau})}+\Ga_b^*(e_*)_4(\widehat{\tau}),\\
         \,[(e_*)_4, (e_*)_a]\widehat{\tau} &=\Ga_g^*+ \big(O(r^{-1})+\Ga_g^*)\nab_*(\widehat{\tau})+\Ga_g^*\widecheck{(e_*)_3(\widehat{\tau})}+\Ga_g^*(e_*)_4(\widehat{\tau}), 
         \end{split}
       \eeaa
and hence, in view of \eqref{eq:definitionofvectorfieldnu*tangenttoSigma*inspanofe3*ande4*}, we have, for $a=1,2$, 
\beaa
\,[\nu_*, (e_*)_a]\widehat{\tau} &=& 2(\eta_*)_a +\Ga_g^*+ \big(O(r^{-1})+\Ga_b^*)\nab_*(\widehat{\tau})+\Ga_b^*\widecheck{(e_*)_3(\widehat{\tau})}+\big(-e_a(b_*)+\Ga_b^*\big)(e_*)_4(\widehat{\tau}).
\eeaa
Also, we have in view of \eqref{definitionofthenewcoordinatetildeu} 
\beaa
(\nab_*)_{\nu_*}\nab_*\widehat{\tau} &=& \nab_*(\nu_*(\widehat{\tau}))+[(\nab_*)_{\nu_*}, \nab_*]\widehat{\tau} \\
&=& [(\nab_*)_{\nu_*}, \nab_*]\widehat{\tau}. 
\eeaa
Together with the above structure of $\,[\nu_*, (e_*)_a]\widehat{\tau}$, and relying on \eqref{eq:dominantconditionforrstarcomparedtotaustartonS*}, \eqref{eq:assumptionsonSigmastarforpartII} and \eqref{eq:bootassfordk*15widecheckealphaoftttonSi*}, we infer on $\Si_*\cap\{\tt_{**}\leq\tau\leq\tt_*\}$
\beaa
\int_{\tau}^{\tau_*}|r\dk_*^{\leq 15}(\nab_*)_{\nu_*}\nab_*\widehat{\tau}|d\tau' &\les& \frac{1}{\tau^{\frac{1}{2}+\dec}}\left(\int_{1}^{\tau_*}\tau^{2+2\dec}\|r\dk_*^k\eta_*\|^2_{L^\infty(S_*(\tau))}d\tau\right)^{\frac{1}{2}}\\
&&+\sqrt{\ep}\int_{\tau}^{\tau_*}\frac{d\tau'}{r{\tau'}^{\frac{1}{2}+\dec}}\\
&\les& \frac{\ep}{\tau^{\frac{1}{2}+\dec}} +\frac{\sqrt{\ep}\tau_*^{\frac{1}{2}-\dec}}{r_*}\les  \frac{\ep}{\tau^{\frac{1}{2}+\dec}}.
\eeaa
Using the general commutation formulas of Lemma \ref{LEMMA:COMM-GEN-B}, which in view of the values of $\B_{ab\mu\nu}$ given by Proposition \ref{proposition:componentsofB} yield for a $k$-tensor $U$
\beaa
\,[(\nab_*)_3, \nab_*] U &=& \big(O(r^{-1})+\Ga_b^*\big)\c\nab_*U+ \Ga_b^*(\nab_*)_3 U +\Ga_b^*(\nab_*)_4 U +r^{-1}\Ga_b^*\c U,\\
\,[(\nab_*)_4, \nab_*] U &=&  \big(O(r^{-1})+\Ga_g^*\big)\c\nab_*U+\Ga_g^* (\nab_*)_4 U +\Ga_g^*(\nab_*)_3 U +r^{-1}\Ga_g^*\c U,
\eeaa
and using also \eqref{eq:dominantconditionforrstarcomparedtotaustartonS*}, \eqref{eq:assumptionsonSigmastarforpartII} and \eqref{eq:bootassfordk*15widecheckealphaoftttonSi*}, we deduce on $\Si_*\cap\{\tt_{**}\leq\tau\leq\tt_*\}$
\beaa
\sum_{j=1}^{16}\int_{\tau}^{\tau_*}|r^j(\nab_*)_{\nu_*}\nab_*^j\widehat{\tau}|d\tau' &\les& \int_{\tau}^{\tau_*}|r\dk_*^{\leq 15}(\nab_*)_{\nu_*}\nab_*\widehat{\tau}|d\tau'\\
&&+\sum_{j_1+j_2=1}^{14}\int_{\tau}^{\tau_*}|r^{2+j_2}(r\nab_*)^{j_1}[(\nab_*)_{\nu_*}, \nab_*]\nab_*^{1+j_2}\widehat{\tau}|d\tau'\\
&\les&   \frac{\ep}{\tau^{\frac{1}{2}+\dec}}+\sqrt{\ep}\int_{\tau}^{\tau_*}\frac{d\tau'}{r{\tau'}^{\frac{1}{2}+\dec}}\\
&\les& \frac{\ep}{\tau^{\frac{1}{2}+\dec}} +\frac{\sqrt{\ep}\tau_*^{\frac{1}{2}-\dec}}{r_*}\les\frac{\ep}{\tau^{\frac{1}{2}+\dec}}. 
\eeaa
Integrating from $S_*$ where $\nab_*(\widehat{\tau})=0$ in view of \eqref{definitionofthenewcoordinatetildeu}, we obtain 
\beaa
r|(r\nab_*)^{\leq 15}\nab_*(\widehat{\tau})|\les\frac{\ep}{\tau^{\frac{1}{2}+\dec}}\quad\textrm{on}\quad\Si_*\cap\{\tt_{**}\leq\tau\leq\tt_*\}
\eeaa
which together with \eqref{eq:improvementofeq:bootassfordk*15widecheckealphaoftttonSi*:e3} and \eqref{eq:improvementofeq:bootassfordk*15widecheckealphaoftttonSi*:e4}, and using again the above commutation formulas, improves \eqref{eq:bootassfordk*15widecheckealphaoftttonSi*}. We thus deduce $\tt_{**}=1$ and 
\bea\lab{eq:bootassfordk*15widecheckealphaoftttonSi*:finalimprovement}
r^{\frac{3}{2}}|\dk_*^{\leq 15}(e_*)_4(\widehat{\tau})|+r|\dk_*^{\leq 15}(\widecheck{(e_*)_3(\widehat{\tau})}, \nab_*(\widehat{\tau}))|\les\frac{\ep}{\tau^{\frac{1}{2}+\dec}}\quad\textrm{on}\quad\Si_*.
\eea

Next, relying on the change of frame formula \eqref{eq:changeofframefromfromframeSigmastarttodoubleprimedframe}, \eqref{eq:changeofframecoefffromfromframeSigmastarttodoubleprimedframe}, as well as \eqref{eq:dominantconditionforrstarcomparedtotaustartonS*}, \eqref{eq:assumptionsonSigmastarforpartII} and \eqref{eq:bootassfordk*15widecheckealphaoftttonSi*:finalimprovement}, we obtain 
\beaa
r\left|\dk_*^{\leq 15}\left(e_3'(\widehat{\tau})-2\left(1+\frac{2m}{r}\right), \nab'(\widehat{\tau})-\frac{a}{r}f_0\right)\right|\les r^{-1}+\frac{\ep}{\tau^{\frac{1}{2}+\dec}}\les \frac{\ep}{\tau^{\frac{1}{2}+\dec}}\quad\textrm{on}\quad\Si_*.
\eeaa
Introducing the notations
\beaa
\widecheck{e_3'(\widehat{\tau})}:= e_3'(\widehat{\tau})-\tmod'(r), \qquad \widecheck{e_4'(\widehat{\tau})}:=e_4'(\widehat{\tau})-\frac{2(r^2+a^2) - \De\tmod'(r)}{|q|^2},\qquad \widecheck{\DD'(\widehat{\tau})}:=& \DD'(\widehat{\tau})-a\Jk,
\eeaa
and noticing that 
\beaa
\widecheck{e_3'(\widehat{\tau})}= e_3'(\widehat{\tau})-2\left(1+\frac{2m}{r}\right)+O(m^2r^{-2}), \quad \widecheck{e_4'(\widehat{\tau})}=e_4'(\widehat{\tau})- \frac{m^2}{r^2}+O(m^3r^{-3}),
\eeaa
we infer, using \eqref{eq:dominantconditionforrstarcomparedtotaustartonS*} and \eqref{eq:relationbetweenJkJkpmandf0onSigmastar},
\beaa
r\left|\dk_*^{\leq 15}\left(\widecheck{e_3'(\widehat{\tau})}, \widecheck{\DD'(\widehat{\tau})}\right)\right|\les r^{-1}+\frac{\ep}{\tau^{\frac{1}{2}+\dec}}\les \frac{\ep}{\tau^{\frac{1}{2}+\dec}}\quad\textrm{on}\quad\Si_*.
\eeaa
On the other hand, we have $\Ytau(\widehat{\tau})=0$ in view of \eqref{definitionofthenewcoordinatetildeu} which together with Remark \ref{rmk:formofYautforrgeqr0intermsofdoubleprimedframe} and \eqref{eq:assumptionsonMMforpartII:secondglobalframeauxassfor15derivatives} implies, on $\Si_*$
\beaa
r^2\left|\dk^{\leq 15}\widecheck{e_4'(\widehat{\tau})}\right|\les \left|\dk^{\leq 15}\widecheck{e_3'(\widehat{\tau})}\right|+r^2\left|\dk^{\leq 15}\Ga_g'\right|\les \left|\dk^{\leq 15}\widecheck{e_3'(\widehat{\tau})}\right|+\frac{\ep}{\tau^{\frac{1}{2}+\dec}}\quad\textrm{on}\quad\Si_*.
\eeaa
Furthermore, using Remark \ref{rmk:formofYautforrgeqr0intermsofdoubleprimedframe}, \eqref{eq:changeofframecoefffromfromframeSigmastarttodoubleprimedframe}, \eqref{eq:changeofframefromfromframeSigmastarttodoubleprimedframe} and \eqref{eq:definitionofvectorfieldnu*tangenttoSigma*inspanofe3*ande4*}, we have
\beaa
\Ytau &=& \Big(1+O(r^{-2})\big(1, \widecheck{\fb_*}, \widecheck{b_*}\big)+\Ga_g'\Big)(e_*)_4+\Big(O(r^{-2})\big(1, \widecheck{\fb_*}, \widecheck{b_*}\big)+\Ga_g'\Big)\dk_*,
\eeaa
which, together with \eqref{eq:assumptionsonSigmastarforpartII} and \eqref{eq:assumptionsonMMforpartII:secondglobalframeauxassfor15derivatives}, implies for a tensor $U$ defined in a neighborhood of $\Si_*$
\beaa
|\dk_*^j\dk U|\les r|\dk_*^j\nab_{(e_*)_4}U|+|\dk_*^jU|\les r|\dk_*^{\leq j}\nab_{\Ytau}U|+|\dk_*^{\leq j+1}U|, \quad j\leq 14,
\eeaa
and then
\beaa
\left|\dk_*^{j_1}\dk^{j_2}U\right| &\les& \left|\dk_*^{\leq j_1+1}\dk^{j_2-1}U\right| + r\left|\dk_*^{\leq j_1}\nab_{Y_{\tau}}\dk^{j_2-1}U\right|\\
&\les& \left|\dk_*^{\leq j_1+1}\dk^{\leq j_2-1}U\right| + r\left|\dk_*^{\leq j_1}\dk^{j_2-1}\nab_{Y_{\tau}}U\right|, \quad 1\leq j_2\leq 15-j_1.
\eeaa
The above estimates finally imply 
\bea\lab{eq:controlforwidecheckealpha'oftttonSi*}
r^2\left|\dk^{\leq 15}\widecheck{e_4'(\widehat{\tau})}\right|+r\left|\dk^{\leq 15}\left(\widecheck{e_3'(\widehat{\tau})}, \widecheck{\DD'(\widehat{\tau})}\right)\right|\les \frac{\ep}{\tau^{\frac{1}{2}+\dec}}\quad\textrm{on}\quad\Si_*.
\eea

Next, we extend \eqref{eq:controlforwidecheckealpha'oftttonSi*} to $\Mext\cap\{r\geq r_0+1\}$. To this end, we assume the following bootstrap assumptions
\bea\lab{eq:bootstrapassumptionsderivativetttonMextrgeqr0plus1doubleprimeframe}
r^2\left|\dk^{\leq 15}\widecheck{e_4'(\widehat{\tau})}\right|+r\left|\dk^{\leq 15}\left(\widecheck{e_3'(\widehat{\tau})}, \widecheck{\DD'(\widehat{\tau})}\right)\right| &\leq& \frac{\sqrt{\ep}}{\tau^{\frac{1}{2}+\dec}}\quad\textrm{on}\quad\Mext\cap\{r\geq r_1\},
\eea
for some $r_1\geq r_0+1$. We then introduce the approximate Killing vectorfield $\T'$ w.r.t. the global null frame $(e_4', e_3', e_a')$, i.e. 
\bea\lab{eq:defonT'expressedindoubleprimedframe}
\T' &=& \frac{1}{2}\left(e_4'+\frac{\Delta}{|q|^2}e_3' -2a\Re(\Jk)^be_b'\right),
\eea
which satisfies in view of Lemma C.5.1 in \cite{GKS22}
\bea\lab{eq:commutationformulasforT'expressedindoubleprimedframe}
\, [\T', e_3'] =\Ga_b'\c\nab' +\Ga_b' \c e_3'+\Ga_b'\c e_4', \qquad \, [\T', e_4'] = \Ga_b'\c\nab' +{\Ga_g'} \c e_3'+\Ga_b' \c e_4'.
\eea
In view of Remark \ref{rmk:formofYautforrgeqr0intermsofdoubleprimedframe} and \eqref{definitionofthenewcoordinatetildeu}, we deduce
\beaa
\Ytau(\T'(\widehat{\tau})) &=& \dk^{\leq 1}\Ga_g' +\dk^{\leq 1}\Ga_g'\widecheck{e_3'(\widehat{\tau})}+\Ga_b'\widecheck{\DD'(\widehat{\tau})},
\eeaa
and hence
\bea\lab{eq:transportequationYtauT'tttonMextrgeqr0plus1}
\Ytau(\widecheck{\T'(\widehat{\tau})}) &=& \dk^{\leq 1}\Ga_g' +\dk^{\leq 1}\Ga_g'\widecheck{e_3'(\widehat{\tau})}+\Ga_b'\widecheck{\DD'(\widehat{\tau})}, \qquad \widecheck{\T'(\widehat{\tau})}:=T'(\widehat{\tau})-1.
\eea
Also, we have in view of Remark \ref{rmk:formofYautforrgeqr0intermsofdoubleprimedframe}, \eqref{commutatorbetweenLieTLieZandnabnab4nab3:1}, Lemma \ref{LEMMA:COMM-GEN-B} and \eqref{eq:additionalpropertiesseconddgolbalnullframe}
\bea\lab{eq:commutatorsYtauwithT'andrnab'}
\bsplit
\,[\nab_{\Ytau}', \Lieb_{\T'}]U =&  \dk^{\leq 1}\Ga_g'\c \nab_3'U+r^{-1}\dk^{\leq 1}(\Ga_b' \c U),\\
\,[\nab_{\Ytau}', r\nab']U =& \big(O(r^{-2})+\Ga_g'\big)\c(r\nab')^{\leq 1}U+\Ga_g'\c\nab_4'U+\big(O(r^{-3})+\dk^{\leq 1}\Ga_g'\big)\nab_3'U\\
=& \big(O(r^{-2})+\dk^{\leq 1}\Ga_g'\big)\c r\nab'U+\big(O(r^{-3})+\dk^{\leq 1}\Ga_g'\big)\c\nab_{\Ytau}'U\\
&+\big(O(r^{-3})+\dk^{\leq 1}\Ga_g'\big)\Lieb_{\T'}U\quad\textrm{on}\quad\Mext\cap\{r\geq r_0+1\},
\end{split}
\eea
where we also used Lemma \ref{lemma:basicpropertiesLiebTfasdiuhakdisug:chap9} and the fact that $(\nab_4', \nab_3')$ are generated by $(\nab_{\T'}, \nab_{\Ytau}', \nab')$. Now, commuting \eqref{eq:transportequationYtauT'tttonMextrgeqr0plus1} with \eqref{eq:commutatorsYtauwithT'andrnab'} and integrating the resulting transport equation from $\Si_*$ where $\dk^k\widecheck{\T'(\widehat{\tau})}$ is under control thanks to \eqref{eq:controlforwidecheckealpha'oftttonSi*}, we infer on $\Mext\cap\{r\geq r_1\}$, using also \eqref{eq:assumptionsonMMforpartII:secondglobalframeauxassfor15derivatives} and the bootstrap assumptions \eqref{eq:bootstrapassumptionsderivativetttonMextrgeqr0plus1doubleprimeframe},  
\beaa
\left|(\Lieb_{\T'}, r\nab_{\Ytau}', r\nab')^{\leq 15}\widecheck{\T'(\widehat{\tau})}\right| &\les& \frac{\ep}{r_*\tau^{\frac{1}{2}+\dec}}+\int_r^{r_*}\frac{\ep}{{r'}^2\tau^{\frac{1}{2}+\dec}}dr'\les\frac{\ep}{r\tau^{\frac{1}{2}+\dec}}.
\eeaa
Using again Lemma \ref{lemma:basicpropertiesLiebTfasdiuhakdisug:chap9} and the fact that $(\nab_4', \nab_3')$ are generated by $(\nab_{\T'}, \nab_{\Ytau}', \nab')$, we deduce 
\bea\lab{eq:controlofwidecheckT'widehattau:improvementbootass}
\left|\dk^{\leq 15}\widecheck{\T'(\widehat{\tau})}\right|  \les\frac{\ep}{r\tau^{\frac{1}{2}+\dec}} \quad\textrm{on}\quad\Mext\cap\{r\geq r_1\}.
\eea

Next, we control $\widecheck{\DD'(\widehat{\tau})}$ in $\Mext\cap\{r\geq r_1\}$. In view of the commutation formulas \eqref{eq:comm-nab3-nab4-naba-f-general}, \eqref{eq:additionalpropertiesseconddgolbalnullframe}, Remark \ref{rmk:formofYautforrgeqr0intermsofdoubleprimedframe} and \eqref{definitionofthenewcoordinatetildeu}, we have
\beaa
\nab_{\Ytau}'(r\nab'\widehat{\tau}) &=& \big(O(r^{-2})+\Ga_g'\big)\c r\nab'\widehat{\tau} + \big(O(r^{-2})+\Ga_g'\big) re_4'(\widehat{\tau})+\big(O(r^{-2})+\dk^{\leq 1}\Ga_g'\big)e_3'(\widehat{\tau}).
\eeaa
Using the fact that $(\nab_4', \nab_3')$ are generated by $(\nab_{\T'}, \nab_{\Ytau}', \nab')$, and using again \eqref{definitionofthenewcoordinatetildeu}, we infer
\beaa
\nab_{\Ytau}'(r\nab'\widehat{\tau}) &=& \big(O(r^{-2})+\Ga_g'\big)\c r\nab'\widehat{\tau} +\big(O(r^{-2})+\dk^{\leq 1}\Ga_g'\big)\T'(\widehat{\tau}).
\eeaa
Linearizing, and using in particular the fact that $r\widecheck{\nab_4'\Jk}\in\Ga_g'$, we deduce 
\bea\lab{eq:transportequationYtaurnab'tttonMextrgeqr0plus1}
\nab_{\Ytau}'(r\widecheck{\nab'\widehat{\tau}}) = \dk^{\leq 1}\Ga_g'+\big(O(r^{-2})+\Ga_g'\big)\c r\widecheck{\nab'\widehat{\tau}} +\big(O(r^{-2})+\dk^{\leq 1}\Ga_g'\big)\widecheck{\T'(\widehat{\tau})}.
\eea
Now, commuting \eqref{eq:transportequationYtaurnab'tttonMextrgeqr0plus1} with \eqref{eq:commutatorsYtauwithT'andrnab'} and integrating the resulting transport equation from $\Si_*$ where $\dk^kr\widecheck{\nab'(\widehat{\tau})}$ is under control thanks to \eqref{eq:controlforwidecheckealpha'oftttonSi*}, we infer on $\Mext\cap\{r\geq r_1\}$, using also \eqref{eq:assumptionsonMMforpartII:secondglobalframeauxassfor15derivatives} and  \eqref{eq:controlofwidecheckT'widehattau:improvementbootass},  
\beaa
\left|(\Lieb_{\T'}, r\nab_{\Ytau}', r\nab')^{\leq 15}\widecheck{\nab'(\widehat{\tau})}\right| &\les& \frac{\ep}{r\tau^{\frac{1}{2}+\dec}}+\frac{1}{r}\int_r^{r_*}\frac{\ep}{{r'}^2\tau^{\frac{1}{2}+\dec}}dr'\les\frac{\ep}{r\tau^{\frac{1}{2}+\dec}}.
\eeaa
Using again Lemma \ref{lemma:basicpropertiesLiebTfasdiuhakdisug:chap9} and the fact that $(\nab_4', \nab_3')$ are generated by $(\nab_{\T'}, \nab_{\Ytau}', \nab')$, we deduce 
\bea\lab{eq:controlofwidechecknab'widehattau:improvementbootass}
\left|\dk^{\leq 15}\widecheck{\nab'(\widehat{\tau})}\right|  \les\frac{\ep}{r\tau^{\frac{1}{2}+\dec}} \quad\textrm{on}\quad\Mext\cap\{r\geq r_1\}.
\eea

Next, using \eqref{definitionofthenewcoordinatetildeu}, Remark \ref{rmk:formofYautforrgeqr0intermsofdoubleprimedframe}, \eqref{eq:assumptionsonMMforpartII:secondglobalframeauxassfor15derivatives},  
\eqref{eq:controlofwidecheckT'widehattau:improvementbootass},  \eqref{eq:controlofwidechecknab'widehattau:improvementbootass}, and the decomposition of $\T'$ on the null frame $(e_4', e_3', e_1', e_2')$, we infer
\beaa
r^2\left|\dk^{\leq 15}\widecheck{e_4'(\widehat{\tau})}\right| &\les& r^2\left|\dk^{\leq 15}\Ga_g'\right|+r\left|\dk^{\leq 15}\left(\widecheck{\T'(\widehat{\tau})}, \widecheck{\DD'(\widehat{\tau})}\right)\right|\\
&\les& \frac{\ep}{\tau^{\frac{1}{2}+\dec}}\quad\textrm{on}\quad\Mext\cap\{r\geq r_1\},
\eeaa
which together with \eqref{eq:controlofwidecheckT'widehattau:improvementbootass},  \eqref{eq:controlofwidechecknab'widehattau:improvementbootass}, and the decomposition of $\T'$ on the null frame $(e_4', e_3', e_1', e_2')$ implies the following improvement of the bootstrap assumptions \eqref{eq:bootstrapassumptionsderivativetttonMextrgeqr0plus1doubleprimeframe}
\beaa
r^2\left|\dk^{\leq 15}\widecheck{e_4'(\widehat{\tau})}\right|+r\left|\dk^{\leq 15}\left(\widecheck{e_3'(\widehat{\tau})}, \widecheck{\DD'(\widehat{\tau})}\right)\right| &\les& \frac{\ep}{\tau^{\frac{1}{2}+\dec}}\quad\textrm{on}\quad\Mext\cap\{r\geq r_1\}.
\eeaa
We thus deduce $r_1=r_0+1$ and 
\bea\lab{eq:bootstrapassumptionsderivativetttonMextrgeqr0plus1doubleprimeframe:finalimprovement:doubleprimeframe}
r^2\left|\dk^{\leq 15}\widecheck{e_4'(\widehat{\tau})}\right|+r\left|\dk^{\leq 15}\left(\widecheck{e_3'(\widehat{\tau})}, \widecheck{\DD'(\widehat{\tau})}\right)\right| &\les& \frac{\ep}{\tau^{\frac{1}{2}+\dec}}\quad\textrm{on}\quad\Mext\cap\{r\geq r_0+1\}.
\eea
Together with the change of frame formulas \eqref{frametransformation:fromfirsttosecondglobalframe}, and using also the control of $(f, \fb)$ in \eqref{eq:assumptionsonMMforpartII:secondglobalframeauxassfor15derivatives}, we infer the following analog of \eqref{eq:bootstrapassumptionsderivativetttonMextrgeqr0plus1doubleprimeframe:finalimprovement:doubleprimeframe} in the global null frame $(e_3, e_4, e_1, e_2)$
\bea\lab{eq:bootstrapassumptionsderivativetttonMextrgeqr0plus1doubleprimeframe:finalimprovement}
r^2\left|\dk^{\leq 15}\widecheck{e_4(\widehat{\tau})}\right|+r\left|\dk^{\leq 15}\left(\widecheck{e_3(\widehat{\tau})}, \widecheck{\DD(\widehat{\tau})}\right)\right| &\les& \frac{\ep}{\tau^{\frac{1}{2}+\dec}}\quad\textrm{on}\quad\Mext\cap\{r\geq r_0+1\}.
\eea

Note that \eqref{eq:bootstrapassumptionsderivativetttonMextrgeqr0plus1doubleprimeframe:finalimprovement} yields only the first estimate in \eqref{eq:controloftildeuonMext:ealphaoftttcheckinrGag}, but not the second one\footnote{The lack of integrable decay in $\tau$ for $\widecheck{e_\a'(\widehat{\tau})}$ on $\Mext$ is due to the initialization of the transport equation \eqref{eq:transportequationYtaurnab'tttonMextrgeqr0plus1} on $\Si_*$ which requires to rely on the control $|\dk^{\leq 15}r\widecheck{\nab'(\widehat{\tau})}|\les \ep\tau^{-\frac{1}{2}-\dec}$  provided on $\Si_*$ by \eqref{eq:controlforwidecheckealpha'oftttonSi*}.}. To obtain the desired function $\ttt$, we must thus modify $\widehat{\tau}$. To this end, we introduce the smooth cut-off function $\chi=\chi(s)$ such that $\chi=0$ for $s\leq 1$ and $\chi=1$ for $s\geq 2$, and we defined $\ttt$ as follows 
\bea\lab{eq:definitionoftttasaninterpolationoftauandwidehattau}
\ttt &=& \chi\left(\frac{\tau^{\frac{1}{2}}}{r}\right)\tau+\left(1-\chi\left(\frac{\tau^{\frac{1}{2}}}{r}\right)\right)\widehat{\tau}
\eea 
so that we have in particular 
\beaa
\ttt=\widehat{\tau}\quad\textrm{on}\quad\Mext\cap\{r\geq \tau^{\frac{1}{2}}\}, \qquad \ttt=\tau\quad\textrm{on}\quad\MM\cap\left\{r\leq \frac{1}{2}\tau^{\frac{1}{2}}\right\},
\eeaa
which proves the second part of \eqref{eq:additionestimatetttminust}.

To derive estimates for the scalar function $\ttt$ introduced in \eqref{eq:definitionoftttasaninterpolationoftauandwidehattau}, we in particular need to control $\widehat{\tau}-\tau$. Along $\Si_*$, we have
\beaa
\nu_*(\tau) &=& (e_*)_3(\tau)+b_*(e_*)_4(\tau)=2+O(m^2r^{-2})+\widecheck{(e_*)_3(\tau)}+O(r^{-2})\widecheck{b_*}\\
&=& 2+O(m^2r^{-2})+\widecheck{(e_*)_3(\tau)}+r^{-1}\Ga_b^*
\eeaa
and we infer from \eqref{eq:dominantconditionforrstarcomparedtotaustartonS*} and \eqref{eq:assumptionsonSigmastarforpartII}
\beaa
\int_1^{\tau_*}\tau^{2+2\dec}|\nu_*(\tau) - 2|^2 &\les& \int_1^{\tau_*}\tau^{2+2\dec}|\widecheck{(e_*)_3(\tau)}|^2+\frac{\tau_*^{3+2\dec}}{r_*^4}\les \ep^2. 
\eeaa
As $\nu_*(\widehat{\tau})=2$, we deduce 
\beaa
\int_{\tau}^{\tau_*}|\nu_*(\tau-\widehat{\tau})|d\tau' &=& \int_{\tau}^{\tau_*}|\nu_*(\tau)-2|d\tau'\\
&\les& \frac{1}{\tau^{\frac{1}{2}+\dec}}\left(\int_1^{\tau_*}\tau^{2+2\dec}|\nu_*(\tau) - 2|^2\right)^{\frac{1}{2}}\les \frac{\ep}{\tau^{\frac{1}{2}+\dec}}
\eeaa
and hence, since $\widehat{\tau}=\tau=\tau_*$ on $S_*$, we obtain, integrating from $S_*$ along $\Si_*$
\beaa
|\widehat{\tau}-\tau| &\les& \frac{\ep}{\tau^{\frac{1}{2}+\dec}} \quad\textrm{on}\quad\Si_*.
\eeaa
As $\Ytau(\tau)=0$ on $\MM$ in view of Remark \ref{rmk:YtauistangenttolevelhypersurfacesSioftauonMM}, and $\Ytau(\widehat{\tau})=0$ on $\Mext\cap\{r\geq r_0+1\}$ in view of \eqref{definitionofthenewcoordinatetildeu}, the above estimate for $\widehat{\tau}-\tau$ immediately extends from $\Si_*$ to $\Mext\cap\{r\geq r_0+1\}$, i.e. 
\bea\lab{eq:estimateforwidehattauminustauonMextrgeqr0plus1}
|\widehat{\tau}-\tau| &\les& \frac{\ep}{\tau^{\frac{1}{2}+\dec}} \quad\textrm{on}\quad\Mext\cap\{r\geq r_0+1\}.
\eea
In particular, in view of \eqref{eq:definitionoftttasaninterpolationoftauandwidehattau}, we immediately infer
\beaa
|\ttt-\tau| &\les& \frac{\ep}{\tau^{\frac{1}{2}+\dec}} \quad\textrm{on}\quad\MM\cap\left\{r\geq \frac{1}{2}\tau^{\frac{1}{2}}\right\},
\eeaa
as stated in the first part of \eqref{eq:additionestimatetttminust}.

We may now estimate $\widecheck{e_\a(\ttt)}$. In view of \eqref{eq:definitionoftttasaninterpolationoftauandwidehattau}, we have
\beaa
\widecheck{e_\a(\ttt)} &=& \chi\left(\frac{\tau^{\frac{1}{2}}}{r}\right)\widecheck{e_\a(\tau)}+\left(1-\chi\left(\frac{\tau^{\frac{1}{2}}}{r}\right)\right)\widecheck{e_\a(\widehat{\tau})}+O(r^{-1})\chi'\left(\frac{\tau^{\frac{1}{2}}}{r}\right)(\tau-\widehat{\tau})
\eeaa
which, together with the control for $\widecheck{e_\a(\tau)}$ provided by \eqref{eq:assumptionsonMMforpartII}, the control of $\widecheck{e_\a(\widehat{\tau})}$ provided by \eqref{eq:bootstrapassumptionsderivativetttonMextrgeqr0plus1doubleprimeframe:finalimprovement}, and \eqref{eq:estimateforwidehattauminustauonMextrgeqr0plus1} yields
\beaa
r|\dk^{\leq 15}(\widecheck{e_4(\ttt)}, \widecheck{e_3(\ttt)}, \widecheck{\nab(\ttt)})|\les\frac{\ep}{\tau^{\frac{1}{2}+\dec}}, \qquad |\dk^{\leq 15}(\widecheck{e_4(\ttt)}, \widecheck{e_3(\ttt)}, \widecheck{\nab(\ttt)})|\les\frac{\ep}{\tau^{1+\dec}},
\eeaa
as stated in \eqref{eq:controloftildeuonMext:ealphaoftttcheckinrGag}, where we used in particular the fact that $\tau^{1+\dec}\gtrsim r\tau^{\frac{1}{2}+\dec}$ on the support of $\chi(\frac{\tau^{\frac{1}{2}}}{r})$, $\tau^{1+\dec}\les r\tau^{\frac{1}{2}+\dec}$ on the support of $1-\chi(\frac{\tau^{\frac{1}{2}}}{r})$ and $\tau^{1+\dec}\simeq r\tau^{\frac{1}{2}+\dec}$ on the support of $\chi'(\frac{\tau^{\frac{1}{2}}}{r})$.

Finally, it remains to prove \eqref{eq:additionestimatewidechecke4tttbyeprminus2}. We have 
\beaa
\widecheck{e_4(\ttt)} &=& \chi\left(\frac{\tau^{\frac{1}{2}}}{r}\right)\widecheck{e_4(\tau)}+\left(1-\chi\left(\frac{\tau^{\frac{1}{2}}}{r}\right)\right)\widecheck{e_4(\widehat{\tau})}+O(r^{-1})\chi'\left(\frac{\tau^{\frac{1}{2}}}{r}\right)(\tau-\widehat{\tau})\\
&=& \chi\left(\frac{\tau^{\frac{1}{2}}}{r}\right)\Ga_g +\left(1-\chi\left(\frac{\tau^{\frac{1}{2}}}{r}\right)\right)\widecheck{e_4(\widehat{\tau})}+O(r^{-1})\chi'\left(\frac{\tau^{\frac{1}{2}}}{r}\right)(\tau-\widehat{\tau}).
\eeaa
Together with \eqref{eq:assumptionsonMMforpartII}, the control of $\widecheck{e_4(\widehat{\tau})}$ provided by \eqref{eq:bootstrapassumptionsderivativetttonMextrgeqr0plus1doubleprimeframe:finalimprovement}, and \eqref{eq:estimateforwidehattauminustauonMextrgeqr0plus1} yields
\beaa
r^2|\dk^{\leq 15}(\widecheck{e_4(\ttt)}| &\les& \ep+\frac{\ep r}{\tau^{\frac{1}{2}+\dec}}\left|\chi'\left(\frac{\tau^{\frac{1}{2}}}{r}\right)\right|\les\ep
\eeaa
as stated in \eqref{eq:additionestimatewidechecke4tttbyeprminus2}. This concludes the proof of Lemma \ref{lemma:controloftildeuonMext}.
\end{proof}

%%%%%%%%%%%%%%%%%%%%%%%%%%%%%%%%%%%%%%%%%%%%%%

\subsubsection{Construction of the coordinates $(\rr, \tx^1, \tx^2)$ on $\MM$}
\lab{sec:constructionofthecoordinatestx1tx2onMM}

%%%%%%%%%%%%%%%%%%%%%%%%%%%%%%%%%%%%%%%%%%%%%%

In addition to the new scalar function $\ttt$ introduced in Lemma \ref{lemma:controloftildeuonMext}, we define new coordinates $(\rr, \tx^1, \tx^2)$ on $\MM$ as follows 
\bea\lab{definitionofthenewcoordinatetilder}
\rr=r\quad\textrm{on}\quad \Si_*, \qquad \Ytau(\rr)=[\Ytau(\rr)]_K(\rr),\quad\textrm{on}\quad \MM,
\eea
and
\bea\lab{definitionofthenewcoordinatetildethetandtphi}
\bsplit
&\widetilde{\th}=\th, \quad \widetilde{\tphi}=\tphi,\quad\textrm{on}\quad \Si_*, \qquad \Ytau(\widetilde{\th})=0,\quad \Ytau(\widetilde{\tphi})=[\Ytau(\widetilde{\tphi})]_K(\rr),\quad\textrm{on}\quad \MM,\\
& (\tx^1_0, \tx^2_0):=(\widetilde{\th}, \widetilde{\tphi}),\qquad (\tx^1_p, \tx^2_p):=(\sin\widetilde{\th}\cos\widetilde{\tphi}, \sin\widetilde{\th}\cos\widetilde{\tphi}),
\end{split}
\eea
where the Kerr values $[\Ytau(\rr)]_K(\rr)$ and $[\Ytau(\widetilde{\tphi})]_K(\rr)$ are given by
\beaa
\,[\Ytau(\rr)]_K(\rr) &:=& [\Ytau^4|q|^{-2}]_K(\rr)\De(\rr)-[\Ytau^3]_K(\rr),\\
\,[\Ytau(\widetilde{\tphi})]_K(\rr) &:=& [\Ytau^4|q|^{-2}]_K(\rr)(2a -\De(\rr)\phimod'(\rr))+[\Ytau^3]_K(\rr)\phimod'(\rr),
\eeaa
with $[\Ytau^4|q|^{-2}]_K(\rr)$ and $[\Ytau^3]_K(\rr)$ defined by 
\bea\lab{eq:defintionofKerrvaluesofYtau4modqminus2andYtau3}
\bsplit
\,[\Ytau^4|q|^{-2}]_K(\rr) :=&\chi(\rr)\frac{\tmod'(\rr)}{2(\rr^2+a^2) - \De(\rr)\tmod'(\rr)} +(1-\chi(\rr))\frac{1}{\De(\rr)}, \\
\,[\Ytau^3]_K(\rr) :=& -\chi(\rr)-(1-\chi(\rr))\frac{2(\rr^2+a^2) - \De(\rr)\tmod'(\rr)}{\De(\rr)\tmod'(\rr)},
\end{split}
\eea
and with $\chi$ the cut-off appearing in Definition \ref{def:definitionofthevectorfildYtau}. 

We start with the control of $\rr$ and its derivatives.
\begin{lemma}\lab{lemma:controlofnewcoordinatewidetilder}
The coordinates $\rr$ introduced in \eqref{definitionofthenewcoordinatetilder} satisfies
\bea\lab{eq:esitmatefor16weightedderivativesofrrminusr}
|\dk^{\leq 16}(\rr-r)|\les\frac{\ep}{r\tau^{\frac{1}{2}+\dec}}, \qquad |\dk^{\leq 16}(\rr-r)|\les\frac{\ep}{\tau^{1+\frac{7\dec}{8}}},\quad\textrm{on}\quad\MM,
\eea
and 
\bea\lab{eq:esitmatefor15weightedderivativesoffirstorderdrivgivescheckofrr}
\bsplit
\left|\dk^{\leq 15}\widecheck{e_4(\rr)}\right|+r^{-1}|\dk^{\leq 15}\nab(\rr)| &\les \frac{\ep}{r^2\tau^{\frac{1}{2}+\dec}}, \quad b=1,2, \quad\textrm{on}\quad\MM,\\
r^{-1}|\dk^{\leq 15}(\nab(\rr), \widecheck{e_3(\rr)})|+|\dk^{\leq 15}\widecheck{e_4(\rr)}| &\les \frac{\ep}{r\tau^{1+\frac{7\dec}{8}}}, \quad b=1,2, \quad\textrm{on}\quad\MM,
\end{split}
\eea
where the linearized quantities $\widecheck{e_4(\rr)}$ and $\widecheck{e_3(\rr)}$ are given by
\beaa
\widecheck{e_3(\rr)} := e_3(\rr)+1, \qquad\widecheck{e_4(\rr)} := e_4(\rr)-\frac{\Delta}{|q|^2}.
\eeaa
\end{lemma}

\begin{proof}
In view of \eqref{definitionofthenewcoordinatetilder} and the definition of $\Ytau$, we have on $\MM$
\bea\lab{eq:transportequationYtauofrrminusronMM}
\Ytau(\rr - r) &=& O(1)\widecheck{e_4'(r)}+O(r^{-2})\widecheck{e_3'(r)}+\Ga_g'=\Ga_g'.
\eea
Integrating the transport equation \eqref{eq:transportequationYtauofrrminusronMM} from $\Si_*$ where $\rr=r$ in view of \eqref{definitionofthenewcoordinatetilder}, and using the control of $\Ga_g'$ in  \eqref{eq:assumptionsonMMforpartII:secondglobalframeauxassfor15derivatives}, we infer
\bea\lab{eq:intermediarycontrolorrrminusr:0derivatives}
\bsplit
|\rr-r| \les& \, \int_r^{r_*}\frac{\ep}{{r'}^2\tau^{\frac{1}{2}+\dec}}dr'\les \frac{\ep}{r\tau^{\frac{1}{2}+\dec}}\quad\textrm{on}\quad\MM,\\
|\rr-r| \les& \,\int_r^{r_*}\frac{\ep}{{r'}^{1+\frac{\dec}{4}}\tau^{\frac{1}{2}+\frac{7\dec}{8}}}dr'\les \frac{\ep}{\tau^{1+\frac{7\dec}{8}}}\quad\textrm{on}\quad\MM.
\end{split}
\eea

It remains to extend \eqref{eq:intermediarycontrolorrrminusr:0derivatives} to higher derivatives. We first focus on the region $\Mext\cap\{r\geq r_0+1\}$ and commute \eqref{eq:transportequationYtauofrrminusronMM} with \eqref{eq:commutatorsYtauwithT'andrnab'}. We then integrate the resulting transport equation from $\Si_*$, where $\rr=r$ in view of \eqref{definitionofthenewcoordinatetilder}, and infer, using also  \eqref{eq:assumptionsonMMforpartII:secondglobalframeauxassfor15derivatives},  
\beaa
&& |(\Lieb_{\T'}, r\nab_{\Ytau}', r\nab')^{\leq 16}(\rr-r)|\les\frac{\ep}{r\tau^{\frac{1}{2}+\dec}},\quad\textrm{on}\quad\Mext\cap\{r\geq r_0+1\},\\ 
&& |(\Lieb_{\T'}, r\nab_{\Ytau}', r\nab')^{\leq 16}(\rr-r)|\les\frac{\ep}{\tau^{1+\frac{7\dec}{8}}},\quad\textrm{on}\quad\Mext\cap\{r\geq r_0+1\}.
\eeaa
Using Lemma \ref{lemma:basicpropertiesLiebTfasdiuhakdisug:chap9} and the fact that $(\nab_4', \nab_3')$ are generated by $(\nab_{\T'}, \nab_{\Ytau}', \nab')$, we deduce 
\beaa
|\dk^{\leq 16}(\rr-r)|\les\frac{\ep}{r\tau^{\frac{1}{2}+\dec}}, \qquad |\dk^{\leq 16}(\rr-r)|\les\frac{\ep}{\tau^{1+\frac{7\dec}{8}}},\quad\textrm{on}\quad\Mext\cap\{r\geq r_0+1\}.
\eeaa
Then, we easily extend this estimate from $\{r=r_0+1\}$ to $\MM\cap\{r\leq r_0+1\}$ by proceeding similarly as above which yields 
\beaa
|\dk^{\leq 16}(\rr-r)|\les\frac{\ep}{r\tau^{\frac{1}{2}+\dec}}, \qquad |\dk^{\leq 16}(\rr-r)|\les\frac{\ep}{\tau^{1+\frac{7\dec}{8}}},\quad\textrm{on}\quad\MM,
\eeaa
as stated in \eqref{eq:esitmatefor16weightedderivativesofrrminusr}. Finally, \eqref{eq:esitmatefor15weightedderivativesoffirstorderdrivgivescheckofrr} follows from \eqref{eq:esitmatefor16weightedderivativesofrrminusr}, the fact that $\widecheck{e_4(r)}\in\Ga_g$, $\nab(r)\in\Ga_g$, and $\widecheck{e_4(r)}\in r\Ga_b$, and from \eqref{eq:assumptionsonMMforpartII}. This concludes the proof of Lemma \ref{lemma:controlofnewcoordinatewidetilder}.
\end{proof}

The control of $\tx^b$, $b=1,2$, and their derivatives is provided by the following lemma.
\begin{lemma}\lab{lemma:gainforwidechecknabJp}
The coordinates $\tx^b$, $b=1,2$, introduced in \eqref{definitionofthenewcoordinatetildethetandtphi} satisfy 
\bea\lab{eq:esitmatefor16weightedderivativesoftxbminusxb}
|\dk^{\leq 16}(\tx^b - x^b)|\les\frac{\ep}{r\tau^{\frac{1}{2}+\dec}}, \qquad |\dk^{\leq 16}(\tx^b - x^b)|\les\frac{\ep}{\tau^{1+\frac{3\dec}{4}}},\quad\textrm{on}\quad\MM,
\eea
and 
\bea\lab{eq:esitmatefor15weightedderivativesoffirstorderdrivgivescheckoftxb}
\bsplit
\left|\dk^{\leq 15}\widecheck{e_4(\tx^b)}\right|+|\dk^{\leq 15}\widecheck{\nab(\tx^b)}| &\les \frac{\ep}{r^2\tau^{\frac{1}{2}+\dec}}, \quad b=1,2, \quad\textrm{on}\quad\MM,\\
|\dk^{\leq 15}\widecheck{\nab(\tx^b)}|+|\dk^{\leq 15}\big(\widecheck{e_4(\tx^b)}, \,\widecheck{e_3(\tx^b)}\big)| &\les \frac{\ep}{r\tau^{1+\frac{3\dec}{4}}}, \quad b=1,2, \quad\textrm{on}\quad\MM,
\end{split}
\eea
where the linearized quantities $\widecheck{e_\a(\tx^b)}$, $b=1,2$, are given by
\beaa
\bsplit
\widecheck{e_3(\tx^1_p)}:=& e_3(\tx^1_p)+\phimod'(r)x^2_p, \qquad \widecheck{e_4(\tx^1_p)}:=  e_4(\tx^1_p)+\frac{2a -\De\phimod'(r)}{|q|^2}x^2_p, \\ 
\widecheck{e_3(\tx^2_p)}:=& e_3(\tx^2_p)-\phimod'(r)x^1_p, \qquad \widecheck{e_4(\tx^2_p)}:=e_4(\tx^2_p)-\frac{2a -\De\phimod'(r)}{|q|^2}x^1_p,\\
\widecheck{\DD(\cos(\widetilde{\th}))} :=& \DD(\cos(\widetilde{\th})) -i\Jk, \qquad  \widecheck{\DD(\tx^1_p)} := \DD(\tx^1_p) - \Jk_+, \qquad \widecheck{\DD(\tx^2_p)} := \DD(\tx^2_p) - \Jk_-.
\end{split}
\eeaa
\end{lemma}

\begin{proof}
We start with the estimate for $\tx^b - x^b$, $b=1,2$. Using \eqref{definitionofthenewcoordinatetildethetandtphi}, we have 
\bea\lab{eq:transportequationYtauoftxbminusxbonMM}
\nn\Ytau(\tx^b - x^b) &=& -\Ytau^3e_3'(x^b)-\Ytau^4e_4'(x^b)\\
&=& \Ga_g'+O(1)\widecheck{e_4'(x^b)}+O(r^{-2})\widecheck{e_3'(x^b)}, \quad b=1,2,\quad\textrm{on}\quad\MM.
\eea
Integrating from $\Si_*$, where $\tx^b=x^b$ in view of \eqref{definitionofthenewcoordinatetildethetandtphi}, we infer on $\MM$, in view of the control of $\Ga_g'$, $\widecheck{e_4'(x^b)}$ and $\widecheck{e_3'(x^b)}$ given by \eqref{eq:assumptionsonMMforpartII:secondglobalframeauxassfor15derivatives}, 
\bea\lab{eq:intermediarycontroloftxbminusxb:0derivatives}
\bsplit
|\tx^b - x^b| &\les \int_r^{r_*}\frac{\ep}{{r'}^2\tau^{\frac{1}{2}+\dec}}\les\frac{\ep}{r\tau^{\frac{1}{2}+\dec}}, \quad b=1,2,\\ 
|\tx^b - x^b| &\les \int_r^{r_*}\frac{\ep}{{r'}^{1+\frac{\dec}{2}}\tau^{1+\frac{3\dec}{4}}}\les\frac{\ep}{\tau^{1+\frac{3\dec}{4}}},\quad b=1,2,
\end{split}
\eea

It remains to extend \eqref{eq:intermediarycontroloftxbminusxb:0derivatives} to higher derivatives. We first focus on the region $\Mext\cap\{r\geq r_0+1\}$ and commute \eqref{eq:transportequationYtauoftxbminusxbonMM} with \eqref{eq:commutatorsYtauwithT'andrnab'}. We then integrate the resulting transport equation from $\Si_*$, where $\tx^b=x^b$ in view of \eqref{definitionofthenewcoordinatetildethetandtphi}, and infer, using also  \eqref{eq:assumptionsonMMforpartII:secondglobalframeauxassfor15derivatives},  
\beaa
&& |(\Lieb_{\T'}, r\nab_{\Ytau}', r\nab')^{\leq 16}(\tx^b-x^b)|\les\frac{\ep}{r\tau^{\frac{1}{2}+\dec}},\quad\textrm{on}\quad\Mext\cap\{r\geq r_0+1\},\\ 
&& |(\Lieb_{\T'}, r\nab_{\Ytau}', r\nab')^{\leq 16}(\tx^b-x^b)|\les\frac{\ep}{\tau^{1+\frac{3\dec}{4}}},\quad\textrm{on}\quad\Mext\cap\{r\geq r_0+1\}.
\eeaa
Using Lemma \ref{lemma:basicpropertiesLiebTfasdiuhakdisug:chap9} and the fact that $(\nab_4', \nab_3')$ are generated by $(\nab_{\T'}, \nab_{\Ytau}', \nab')$, we deduce 
\beaa
|\dk^{\leq 16}(\tx^b-x^b)|\les\frac{\ep}{r\tau^{\frac{1}{2}+\dec}}, \quad |\dk^{\leq 16}(\tx^b-x^b)|\les\frac{\ep}{\tau^{1+\frac{3\dec}{4}}},\,\,b=1,2,\,\,\,\,\textrm{on}\,\,\,\,\Mext\cap\{r\geq r_0+1\}.
\eeaa
Then, we easily extend this estimate from $\{r=r_0+1\}$ to $\MM\cap\{r\leq r_0+1\}$ by proceeding similarly as above which yields 
\beaa
|\dk^{\leq 16}(\tx^b-x^b)|\les\frac{\ep}{r\tau^{\frac{1}{2}+\dec}}, \qquad |\dk^{\leq 16}(\tx^b-x^b)|\les\frac{\ep}{\tau^{1+\frac{3\dec}{4}}},\quad b=1,2,\quad\textrm{on}\quad\MM,
\eeaa
as stated in \eqref{eq:esitmatefor16weightedderivativesoftxbminusxb}. 

Next, we prove \eqref{eq:esitmatefor15weightedderivativesoffirstorderdrivgivescheckoftxb}. The estimate \eqref{eq:esitmatefor16weightedderivativesoftxbminusxb}, together with \eqref{eq:assumptionsonMMforpartII:secondglobalframeauxassfor15derivatives} immediately yields 
\bea\lab{eq:esitmatefor15weightedderivativesoffirstorderdrivgivescheckoftxb:goodnab4andnabbutnotbab3} 
|\dk^{\leq 15}\widecheck{e_4(\tx^b)}| \les \frac{\ep}{r^2\tau^{\frac{1}{2}+\dec}}, \quad 
|\dk^{\leq 15}\widecheck{e_4(\tx^b)}| \les \frac{\ep}{r\tau^{1+\frac{3\dec}{4}}},\,\, b=1,2,\,\,\textrm{on}\,\,\MM,
\eea
and 
\bea\lab{eq:esitmatefor15weightedderivativesoffirstorderdrivgivescheckoftxb:nab3onlyrGagbutnotGab} 
|\dk^{\leq 15}\widecheck{e_3(\tx^b)}| \les \frac{\ep}{r\tau^{\frac{1}{2}+\dec}}, \quad 
r|\dk^{\leq 15}\widecheck{\nab(\tx^b)}|+|\dk^{\leq 15}\widecheck{e_3(\tx^b)}| \les \frac{\ep}{\tau^{1+\frac{3\dec}{4}}},\,\, b=1,2,\quad\textrm{on}\quad\MM.
\eea
Note that \eqref{eq:esitmatefor15weightedderivativesoffirstorderdrivgivescheckoftxb:goodnab4andnabbutnotbab3} implies the stated estimates for $\widecheck{e_4(\tx^b)}$ in \eqref{eq:esitmatefor15weightedderivativesoffirstorderdrivgivescheckoftxb}. On the other hand, we need to improve on the estimate  \eqref{eq:esitmatefor15weightedderivativesoffirstorderdrivgivescheckoftxb:nab3onlyrGagbutnotGab} for $\widecheck{e_3(\tx^b)}$ and $\widecheck{\nab(\tx^b)}$. To this end, we first notice that we have in view of \eqref{eq:corrdinatesthandvarphiarepropagatedalongSi*bynu*} and \eqref{eq:assumptionsonSigmastarforpartII},
\beaa
|\dk_*^{\leq 15}\widecheck{\nab_*x^1}|+|\dk_*^{\leq 15}\widecheck{\nab_*x^2}|  \leq \frac{\ep}{r^2\tau^{\frac{1}{2}+\dec}}, \qquad \nu_*(x^1)=\nu_*(x^2)=0, \quad\textrm{on}\quad\Si_*,
\eeaa
which together with \eqref{definitionofthenewcoordinatetildethetandtphi} implies 
\bea\lab{eq:firstestimatesforderivesoftx1andtx2onSi*:tangentialtoSi*}
|\dk_*^{\leq 15}\widecheck{\nab_*\tx^1}|+|\dk_*^{\leq 15}\widecheck{\nab_*\tx^2}|  \leq \frac{\ep}{r^2\tau^{\frac{1}{2}+\dec}}, \qquad \nu_*(\tx^1)=\nu_*(\tx^2)=0, \quad\textrm{on}\quad\Si_*.
\eea
Also, using again \eqref{definitionofthenewcoordinatetildethetandtphi}, as well as Remark \ref{rmk:formofYautforrgeqr0intermsofdoubleprimedframe}, we have
\bea\lab{eq:e4'ofxtbintermsofe3'onMextrgeqr0plus1}
e_4'(\tx^b) &=& \big(O(m^2r^{-2})+\Ga_g'\big)e_3'(\tx^b), \quad b=1,2, \quad\textrm{on}\quad\Mext\cap\{r\geq r_0+1\}.
\eea
In view of the change of frame formula \eqref{eq:changeofframefromfromframeSigmastarttodoubleprimedframe}, \eqref{eq:changeofframecoefffromfromframeSigmastarttodoubleprimedframe}, we deduce
\beaa
(e_*)_4(\tx^b) &=& O(m^2r^{-2})  +O(r^{-1})\widecheck{\nab_*\tx^b}+ \big(O(m^2r^{-2})+\Ga_g'\big)\Big((e_*)_3(\tx^b)+O(r^{-1})\widecheck{\fb_*}\Big),
\eeaa
which together with \eqref{eq:firstestimatesforderivesoftx1andtx2onSi*:tangentialtoSi*}, and in particular $(e_3)_*(\tx^b)=-b_*(e_*)_4(x^b)$, and using also \eqref{eq:assumptionsonMMforpartII:secondglobalframeauxassfor15derivatives} and \eqref{eq:assumptionsonSigmastarforpartII}, yields
\bea\lab{eq:firstestimatesforderivesoftx1andtx2onSi*:transeversalderivativestoSi*}
|\dk_*^{\leq 15}(e_*)_4(\tx^b)|+|\dk_*^{\leq 15}(e_*)_3(\tx^b)| \les r^{-2}, \quad b=1,2,\quad\textrm{on}\quad\Si_*.
\eea

Next, using again the change of frame formula \eqref{eq:changeofframefromfromframeSigmastarttodoubleprimedframe}, \eqref{eq:changeofframecoefffromfromframeSigmastarttodoubleprimedframe}, we have
\beaa
 \bsplit
  e_a'(\tx^b) =&(e_*)_a(\tx^b) +O(r^{-3})+O(r^{-2})\widecheck{\fb_*}+O(r^{-2})\widecheck{\nab_*(\tx^b)}+O(r^{-1})(e_*)_4(\tx^b)+O(r^{-1})(e_*)_3(\tx^b),\\
 e_3'(\tx^b) =& \big(1+O(mr^{-1})\big)(e_*)_3(\tx^b)+O(r^{-2})+O(r^{-1})\widecheck{\fb_*}+O(r^{-1})\widecheck{\nab_*(\tx^b)}\\
 &+O(r^{-2})(e_*)_4(\tx^b)+O(r^{-2})(e_*)_3(\tx^b),
 \end{split}
 \eeaa
which, in view of Definition \ref{def:renormalizationforf0fpfm} and \eqref{eq:relationbetweenJkJkpmandf0onSigmastar} yields
\beaa
 \bsplit
  \widecheck{\nab'(\tx^b)} =& \widecheck{\nab_*(\tx^b)} +O(r^{-3})+O(r^{-2})\widecheck{\fb_*}+O(r^{-2})\widecheck{\nab_*(\tx^b)}+O(r^{-1})(e_*)_4(\tx^b)+O(r^{-1})(e_*)_3(\tx^b),\\
 \widecheck{e_3'(\tx^b)} =& \big(1+O(mr^{-1})\big)(e_*)_3(\tx^b)+O(r^{-2})+O(r^{-1})\widecheck{\fb_*}+O(r^{-1})\widecheck{\nab_*(\tx^b)}+O(r^{-2})(e_*)_4(\tx^b)\\
 &+O(r^{-2})(e_*)_3(\tx^b).
 \end{split}
 \eeaa
Together with \eqref{eq:firstestimatesforderivesoftx1andtx2onSi*:tangentialtoSi*} and \eqref{eq:firstestimatesforderivesoftx1andtx2onSi*:transeversalderivativestoSi*}, we deduce
\beaa
|\dk_*^{\leq 15}\widecheck{\nab'(\tx^b)}| \les \frac{\ep}{r^2\tau^{\frac{1}{2}+\dec}}+\frac{1}{r^3}, \qquad |\dk_*^{\leq 15}\widecheck{e_3'(\tx^b)}|\les \frac{1}{r^2} \quad\textrm{on}\quad\Si_*
\eeaa
which together with \eqref{eq:dominantconditionforrstarcomparedtotaustartonS*} implies 
\bea\lab{eq:estimatesforderivesoftx1andtx2onSi*indoubleprimeframe:differentiationbydk*}
|\dk_*^{\leq 15}\widecheck{\nab'(\tx^b)}| \les \frac{\ep}{r^2\tau^{\frac{1}{2}+\dec}}, \qquad |\dk_*^{\leq 15}\widecheck{e_3'(\tx^b)}|\les \frac{\ep}{r\tau^{1+\dec}}, \quad b=1,2, \quad\textrm{on}\quad\Si_*.
\eea
Then, combining \eqref{eq:e4'ofxtbintermsofe3'onMextrgeqr0plus1} and \eqref{eq:estimatesforderivesoftx1andtx2onSi*indoubleprimeframe:differentiationbydk*}, and using the commutation formulas of Section \ref{sec:generalcommutationformulasrealcase}, we immediately infer 
\bea\lab{eq:estimatesforderivesoftx1andtx2onSi*indoubleprimeframe:differentiationbydksoallderivatives}
|\dk^{\leq 15}\widecheck{\nab'(\tx^b)}| \les \frac{\ep}{r^2\tau^{\frac{1}{2}+\dec}}, \qquad |\dk^{\leq 15}\widecheck{e_3'(\tx^b)}|\les \frac{\ep}{r\tau^{1+\dec}}, \quad b=1,2, \quad\textrm{on}\quad\Si_*.
\eea

Next, we extend \eqref{eq:estimatesforderivesoftx1andtx2onSi*indoubleprimeframe:differentiationbydksoallderivatives} to $\MM$ starting with the estimate for $\widecheck{e_3'(\tx^b)}$. To this end, we rely on the vectorfield $\T'$ introduced in \eqref{eq:defonT'expressedindoubleprimedframe} and commute the transport equations in \eqref{definitionofthenewcoordinatetildethetandtphi} w.r.t. $\Lieb_{\T'}$. In view of the first commutation formula in \eqref{eq:commutatorsYtauwithT'andrnab'}, we obtain on $\MM$
\bea\lab{eq:transportequationrLiebTprimetxb}
\nab_{\Ytau}'\T'(\tx^b) &=& O(r^{-3})\T'(\rr)+\dk^{\leq 1}(\Ga_g')e_3'(\tx^b)+r^{-1}\dk^{\leq 1}(\Ga_b'\tx^b).
\eea
Integrating from $\Si_*$ where \eqref{eq:estimatesforderivesoftx1andtx2onSi*indoubleprimeframe:differentiationbydksoallderivatives} holds, and using \eqref{eq:esitmatefor15weightedderivativesoffirstorderdrivgivescheckoftxb:goodnab4andnabbutnotbab3}, \eqref{eq:esitmatefor15weightedderivativesoffirstorderdrivgivescheckoftxb:nab3onlyrGagbutnotGab}, \eqref{eq:assumptionsonMMforpartII:secondglobalframeauxassfor15derivatives}, as well as \eqref{eq:esitmatefor15weightedderivativesoffirstorderdrivgivescheckofrr}, we infer
\beaa
|\T'(\tx^b)| &\les& \frac{\ep}{r_*\tau^{1+\dec}}+\int_r^{r_*}\frac{\ep}{{r'}^2\tau^{1+\dec}}\les\frac{\ep}{r\tau^{1+\dec}},\quad\textrm{on}\quad\MM.
\eeaa
Commuting the above transport equation for $T'(\tx^b)$, we then extend the above estimates to higher order derivatives to obtain 
\beaa
|\dk^{\leq 15}\T'(\tx^b)| &\les& \frac{\ep}{r\tau^{1+\frac{3\dec}{4}}},\quad\textrm{on}\quad\MM.
\eeaa
Together with \eqref{eq:esitmatefor15weightedderivativesoffirstorderdrivgivescheckoftxb:goodnab4andnabbutnotbab3}, \eqref{eq:esitmatefor15weightedderivativesoffirstorderdrivgivescheckoftxb:nab3onlyrGagbutnotGab}, the definition of $\T'$ and Lemma \ref{lemma:basicpropertiesLiebTfasdiuhakdisug:chap9}, we deduce 
\bea\lab{eq:recoveringthecorrectestimateforwidechecke3txb15weigthedderivatives}
|\dk^{\leq 15}\widecheck{e_3(\tx^b)}| &\les& \frac{\ep}{r\tau^{1+\frac{3\dec}{4}}}, \quad b=1,2, \quad\textrm{on}\quad\MM,
\eea
which is the stated estimate for $\widecheck{e_3(\tx^b)}$ in \eqref{eq:esitmatefor15weightedderivativesoffirstorderdrivgivescheckoftxb}.

It remains to recover the stated estimate for $\widecheck{\nab(\tx^b)}$ in \eqref{eq:esitmatefor15weightedderivativesoffirstorderdrivgivescheckoftxb}. We start with the region $\Mext\cap\{r\geq r_0+1\}$. In view of the second commutation formula in \eqref{eq:commutatorsYtauwithT'andrnab'} and \eqref{definitionofthenewcoordinatetildethetandtphi}, we have
\beaa
\nab_{\Ytau}'(r\nab'\tx^b) &=& \left[e_4' + \big(O(m^2r^{-2})+\Ga_g'\big)e_3', r\nab'\right](\tx^b)+r\nab'([\Ytau(\tx^b)]_K(r, \tx^1, \tx^2))\\
&=& r^{-2}\Ga_g'+\big(O(r^{-2})+\Ga_g'\big)\c r\nab'\tx^c + \big(O(r^{-2})+\Ga_g'\big) re_4'(\tx^b)+\big(O(r^{-2})+\dk^{\leq 1}\Ga_g'\big)e_3'(\tx^b).
\eeaa
Linearizing, and using in particular the fact that $r\widecheck{\nab_4'\Jk}\in\Ga_g'$, we deduce 
\bea\lab{eq:transportequationYtaurnab'txbonMextrgeqr0plus1}
\nab_{\Ytau}'(r\widecheck{\nab'\tx^b}) = \Ga_g'+\big(O(r^{-2})+\Ga_g'\big)\c r\widecheck{\nab'\tx^c} +\big(O(r^{-2})+\dk^{\leq 1}\Ga_g'\big)\big(\widecheck{e_4'(\tx^b)}, \widecheck{e_3'(\tx^b)}\big).
\eea
Now, commuting \eqref{eq:transportequationYtaurnab'txbonMextrgeqr0plus1} with \eqref{eq:commutatorsYtauwithT'andrnab'} and integrating the resulting transport equation from $\Si_*$ where $\dk^kr\widecheck{\nab'(\tx^b)}$ is under control thanks to \eqref{eq:estimatesforderivesoftx1andtx2onSi*indoubleprimeframe:differentiationbydksoallderivatives}, we infer on $\Mext\cap\{r\geq r_0+1\}$, using also \eqref{eq:assumptionsonMMforpartII:secondglobalframeauxassfor15derivatives} and \eqref{eq:esitmatefor15weightedderivativesoffirstorderdrivgivescheckoftxb:goodnab4andnabbutnotbab3} \eqref{eq:recoveringthecorrectestimateforwidechecke3txb15weigthedderivatives}, 
\beaa
\left|(\Lieb_{\T'}, r\nab_{\Ytau}', r\nab')^{\leq 15}\widecheck{\nab'(\widehat{\tau})}\right| &\les& \frac{\ep}{rr_*\tau^{\frac{1}{2}+\dec}}+\frac{1}{r}\int_r^{r_*}\frac{\ep}{{r'}^2\tau^{\frac{1}{2}+\dec}}dr'\les\frac{\ep}{r^2\tau^{\frac{1}{2}+\dec}}.
\eeaa
Using again Lemma \ref{lemma:basicpropertiesLiebTfasdiuhakdisug:chap9} and the fact that $(\nab_4', \nab_3')$ are generated by $(\nab_{\T'}, \nab_{\Ytau}', \nab')$, we deduce 
\beaa
\left|\dk^{\leq 15}\widecheck{\nab'(\tx^b)}\right|  \les\frac{\ep}{r^2\tau^{\frac{1}{2}+\dec}} \quad\textrm{on}\quad\Mext\cap\{r\geq r_0+1\}.
\eeaa
Together with the change of frame formulas \eqref{frametransformation:fromfirsttosecondglobalframe}, and using also the control of $(f, \fb)$ in \eqref{eq:assumptionsonMMforpartII:secondglobalframeauxassfor15derivatives}, we infer 
\beaa
\left|\dk^{\leq 15}\widecheck{\nab(\tx^b)}\right|  \les\frac{\ep}{r^2\tau^{\frac{1}{2}+\dec}} \quad\textrm{on}\quad\Mext\cap\{r\geq r_0+1\}.
\eeaa
Together with 
\eqref{eq:esitmatefor15weightedderivativesoffirstorderdrivgivescheckoftxb:goodnab4andnabbutnotbab3}, \eqref{eq:esitmatefor15weightedderivativesoffirstorderdrivgivescheckoftxb:nab3onlyrGagbutnotGab}  and \eqref{eq:recoveringthecorrectestimateforwidechecke3txb15weigthedderivatives}, this yields the stated estimates  \eqref{eq:esitmatefor15weightedderivativesoffirstorderdrivgivescheckoftxb} which concludes the proof of Lemma \ref{lemma:gainforwidechecknabJp}.
\end{proof}

%%%%%%%%%%%%%%%%%%%%%%%%%%%%%%%%%%%%%%%%%%%%%%

\subsubsection{Construction of complex 1-forms $\widetilde{\Jk}$, $\widetilde{\Jk}_{\pm}$}
\lab{sec:constructionoftJkandtJkpminKerrpert}
  
%%%%%%%%%%%%%%%%%%%%%%%%%%%%%%%%%%%%%%%%%%%%%%

We define new complex 1-forms on $\MM$ as follows 
\bea\lab{definitionofthenewcomplex1formstJktJkpm}
\bsplit
&\tJk=\Jk, \quad \tJk_{\pm}=\Jk_{\pm},\quad\textrm{on}\quad \Si_*, \\
& \nab_{\Ytau}'\tJk = \left(\frac{1}{\ov{\tq}}[\Ytau^3]_K(\rr) -\frac{\De(\rr) \ov{\tq}}{|\tq|^2}[\Ytau^4|q|^{-2}]_K(\rr)\right)\tJk,\quad\textrm{on}\quad \MM,\\
& \nab_{\Ytau}'\tJk_{\pm} = \left(\frac{1}{\ov{\tq}}[\Ytau^3]_K(\rr) -\frac{\De(\rr)\ov{\tq}}{|\tq|^2}[\Ytau^4|q|^{-2}]_K(\rr)\right)\tJk_{\pm}\\
&\qquad\qquad \mp  \Big(\phimod'(r)[\Ytau^3]_K(\rr)+\big(2a-\De\phimod'(r)\big)[\Ytau^4|q|^{-2}]_K(\rr)\Big)\tJk_{\mp},\quad\textrm{on}\quad \MM,
\end{split}
\eea
where $[\Ytau^3]_K(\rr)$ and $[\Ytau^4|q|^{-2}]_K(\rr)$ are given by \eqref{eq:defintionofKerrvaluesofYtau4modqminus2andYtau3} and where we have introduced the following notation 
\beaa
\tq:=\rr+ia\cos(\widetilde{\th}).
\eeaa

The control of $\tJk$ and $\tJk_\pm$ is provided by the following lemma.
\begin{lemma}\lab{lemma:controlofthecomplex1formswidetildeJkandwidetildeJkplusminus}
The complex 1-forms $\tJk$ and $\tJk_\pm$ introduced in \eqref{definitionofthenewcomplex1formstJktJkpm} satisfy 
\bea\lab{eq:esitmatefor16weightedderivativesoftJkandtJKpmminusJkandJkpm}
|\dk^{\leq 16}(\tJk-\Jk, \tJk_{\pm}-\Jk_{\pm})|\les\frac{\ep}{r^2\tau^{\frac{1}{2}+\dec}}, \qquad |\dk^{\leq 16}(\tJk-\Jk, \tJk_{\pm}-\Jk_{\pm})|\les\frac{\ep}{r\tau^{1+\frac{3\dec}{4}}},\quad\textrm{on}\quad\MM,
\eea
and 
\bea\lab{eq:esitmatefor15weightedderivativesoffirstorderdrivgivescheckoftJkandtJkpm} 
\bsplit
|\dk^{\leq 15}(\widecheck{\nab_4\tJk}, \widecheck{\nab_4\tJk_{\pm}})|&\les\frac{\ep}{r^3\tau^{\frac{1}{2}+\dec}}, \quad\textrm{on}\quad\MM,\\ 
|\dk^{\leq 15}(\widecheck{\nab_4\tJk}, \widecheck{\nab_4\tJk_{\pm}}, \widecheck{\nab\tJk_{\pm}}, \widecheck{\nab\tJk_{\pm}}, \widecheck{\nab_3\tJk}, \widecheck{\nab_3\tJk_{\pm}})| &\les\frac{\ep}{r^2\tau^{1+\frac{3\dec}{4}}},\quad\textrm{on}\quad\MM.
\end{split}
\eea
 where the linearized quantities associated to $\tJk_{\pm}$ are given by 
\beaa
\bsplit
\widecheck{\nab_3\tJk_\pm}:=&\nab_3\tJk_\pm - \frac{1}{\ov{q}}\Jk_\pm \pm\phimod'(r)\Jk_{\mp},\qquad \widecheck{\nab_4 \tJk_\pm}:=\nab_4 \tJk_\pm + \frac{\De\ov{q}}{|q|^4}\Jk_{\pm} \pm\frac{2a-\De\phimod'(r)}{|q|^2}\Jk_{\mp},\\
\widecheck{\ov{\DD}\c \tJk_+}:=&\ov{\DD}\c \tJk_+ + \frac{4r^2 }{|q|^4}x^1_p + \frac{4ia^2\cos(\th)}{|q|^4}x^2_p,\qquad \widecheck{\ov{\DD}\c \tJk_-}:=\ov{\DD}\c \tJk_- + \frac{4r^2 }{|q|^4}x^2_p - \frac{4ia^2\cos(\th)}{|q|^4}x^1_p.
\end{split}
\eeaa
 
Also, the following identities hold on $\MM$
\bea\lab{eq:usefulalgebraicidentitiesinvolvingscalarproductsReJkReJkpm:Kerrpert:tilde}
\bsplit
\Re(\tJk)\c\Re(\tJk) &=\frac{(\sin(\widetilde{\th}))^2}{|\tq|^2}, \qquad \Re(\tJk)\c\Re(\tJk_+)=-\frac{\tx^2_p}{|\tq|^2},\qquad \Re(\tJk)\c\Re(\tJk_-)=\frac{\tx^1_p}{|\tq|^2},\\
\dual(\Re(\tJk))\c\Re(\tJk_+) &=\frac{\cos(\widetilde{\th}) \tx^1_p}{|\tq|^2},\qquad \dual(\Re(\tJk))\c\Re(\tJk_-)=\frac{\cos(\widetilde{\th})\tx^2_p}{|\tq|^2}.
\end{split}
\eea
\end{lemma}

\begin{proof}
In view of \eqref{definitionofthenewcomplex1formstJktJkpm} and the definition of $\Ytau$, we have on $\MM$
\beaa
\nab_{\Ytau}'(\tJk - \Jk) &=& \left(\frac{1}{\ov{\tq}}\Ytau^3 -\frac{\De(\rr) \ov{\tq}}{|\tq|^4}\Ytau^4\right)\tJk -\left(\frac{1}{\ov{q}}\Ytau^3 -\frac{\De(r) \ov{q}}{|q|^4}\Ytau^4\right)\Jk + O(1)\widecheck{\nab_4'\Jk}+O(r^{-2})\widecheck{\nab_3'\Jk}\\
&=& -\frac{1}{r}(\tJk - \Jk)+O(mr^{-2})(\tJk - \Jk)+O(r^{-3})(\rr-r)+O(r^{-3})(\tx^b-x^b) + r^{-1}\Ga_g'
\eeaa
and hence
\bea\lab{eq:transportequationYtauofrtimestJkminusJk}
\nab_{\Ytau}'\big(r(\tJk - \Jk)\big) &=& O(mr^{-2})r(\tJk - \Jk)+O(r^{-2})(\rr-r)+O(r^{-2})(\tx^b-x^b) + \Ga_g'.
\eea
Similarly, using again, \eqref{definitionofthenewcomplex1formstJktJkpm}, we have
\bea\lab{eq:transportequationYtauofrtimestJkpmminusJkpm}
\nn\nab_{\Ytau}'\big(r(\tJk_{\pm} - \Jk_{\pm})\big) &=& O(mr^{-2})r(\tJk_{\pm} - \Jk_{\pm})+O(mr^{-2})r(\tJk_{\mp} - \Jk_{\mp})\\
&&+O(r^{-2})(\rr-r)+O(r^{-2})(\tx^b-x^b) + \Ga_g'.
\eea
Integrating the transport equations \eqref{eq:transportequationYtauofrtimestJkminusJk} and \eqref{eq:transportequationYtauofrtimestJkpmminusJkpm} from $\Si_*$ where $\tJk=\Jk$ and $\tJk_{\pm}=\Jk_{\pm}$ in view of \eqref{definitionofthenewcomplex1formstJktJkpm}, and using the control of $\Ga_g'$, $\rr-r$ and $\tx^b-x^b$ respectively in  \eqref{eq:assumptionsonMMforpartII:secondglobalframeauxassfor15derivatives}, \eqref{eq:esitmatefor16weightedderivativesofrrminusr} and \eqref{eq:esitmatefor16weightedderivativesoftxbminusxb}, we infer
\bea\lab{eq:intermediarycontrolortJkminusJk:0derivatives}
\bsplit
|\tJk-\Jk|+|\tJk_{\pm}-\Jk_{\pm}| \les& \, \frac{1}{r}\int_r^{r_*}\frac{\ep}{{r'}^2\tau^{\frac{1}{2}+\dec}}dr'\les \frac{\ep}{r^2\tau^{\frac{1}{2}+\dec}}\quad\textrm{on}\quad\MM,\\
|\tJk-\Jk|+|\tJk_{\pm}-\Jk_{\pm}|  \les& \,\frac{1}{r}\int_r^{r_*}\frac{\ep}{{r'}^{1+\frac{\dec}{2}}\tau^{\frac{1}{2}+\frac{3\dec}{4}}}dr'\les \frac{\ep}{r\tau^{1+\frac{3\dec}{4}}}\quad\textrm{on}\quad\MM.
\end{split}
\eea

It remains to extend \eqref{eq:intermediarycontrolortJkminusJk:0derivatives} to higher derivatives. We first focus on the region $\Mext\cap\{r\geq r_0+1\}$ and commute \eqref{eq:transportequationYtauofrtimestJkminusJk} and \eqref{eq:transportequationYtauofrtimestJkpmminusJkpm} with \eqref{eq:commutatorsYtauwithT'andrnab'}. We then integrate the resulting transport equation from $\Si_*$, where $\tJk=\Jk$ and $\tJk_{\pm}=\Jk_{\pm}$ in view of \eqref{definitionofthenewcomplex1formstJktJkpm}, and infer, using also  \eqref{eq:assumptionsonMMforpartII:secondglobalframeauxassfor15derivatives}, as well as \eqref{eq:esitmatefor16weightedderivativesofrrminusr} and \eqref{eq:esitmatefor16weightedderivativesoftxbminusxb}, 
\beaa
&& |(\Lieb_{\T'}, r\nab_{\Ytau}', r\nab')^{\leq 16}(\tJk-\Jk, \tJk_{\pm}-\Jk_{\pm})|\les\frac{\ep}{r^2\tau^{\frac{1}{2}+\dec}},\quad\textrm{on}\quad\Mext\cap\{r\geq r_0+1\},\\ 
&& |(\Lieb_{\T'}, r\nab_{\Ytau}', r\nab')^{\leq 16}(\tJk-\Jk, \tJk_{\pm}-\Jk_{\pm})|\les\frac{\ep}{r\tau^{1+\frac{3\dec}{4}}},\quad\textrm{on}\quad\Mext\cap\{r\geq r_0+1\}.
\eeaa
Using Lemma \ref{lemma:basicpropertiesLiebTfasdiuhakdisug:chap9} and the fact that $(\nab_4', \nab_3')$ are generated by $(\nab_{\T'}, \nab_{\Ytau}', \nab')$, we deduce 
\beaa
&& |\dk^{\leq 16}(\tJk-\Jk, \tJk_{\pm}-\Jk_{\pm})|\les\frac{\ep}{r^2\tau^{\frac{1}{2}+\dec}}, \quad\textrm{on}\quad\Mext\cap\{r\geq r_0+1\},\\ 
&&|\dk^{\leq 16}(\tJk-\Jk, \tJk_{\pm}-\Jk_{\pm})|\les\frac{\ep}{r\tau^{1+\frac{3\dec}{4}}},\quad\textrm{on}\quad\Mext\cap\{r\geq r_0+1\}.
\eeaa
Then, we easily extend this estimate from $\{r=r_0+1\}$ to $\MM\cap\{r\leq r_0+1\}$ by proceeding similarly as above which yields 
\beaa
|\dk^{\leq 16}(\tJk-\Jk, \tJk_{\pm}-\Jk_{\pm})|\les\frac{\ep}{r^2\tau^{\frac{1}{2}+\dec}}, \qquad |\dk^{\leq 16}(\tJk-\Jk, \tJk_{\pm}-\Jk_{\pm})\les\frac{\ep}{r\tau^{1+\frac{3\dec}{4}}},\quad\textrm{on}\quad\MM,
\eeaa
as stated in \eqref{eq:esitmatefor16weightedderivativesoftJkandtJKpmminusJkandJkpm}. 

Next, \eqref{eq:esitmatefor16weightedderivativesoftJkandtJKpmminusJkandJkpm}, together with the fact that $\widecheck{\nab_4\Jk}, \widecheck{\nab_4\Jk_{\pm}}\in r^{-1}\Ga_g$ and $\widecheck{\nab\Jk}, \widecheck{\nab\Jk_{\pm}}\in r^{-1}\Ga_b$ immediately yields 
\bea\lab{eq:esitmatefor15weightedderivativesoffirstorderdrivgivescheckoftJkandtJkpm:goodnab4andnabbutnotbab3} 
\bsplit
|\dk^{\leq 15}(\widecheck{\nab_4\tJk}, \widecheck{\nab_4\tJk_{\pm}})|&\les\frac{\ep}{r^3\tau^{\frac{1}{2}+\dec}}, \quad\textrm{on}\quad\MM,\\ 
|\dk^{\leq 15}(\widecheck{\nab_4\tJk}, \widecheck{\nab_4\tJk_{\pm}}, \widecheck{\nab\tJk_{\pm}}, \widecheck{\nab\tJk_{\pm}})| &\les\frac{\ep}{r^2\tau^{1+\frac{3\dec}{4}}},\quad\textrm{on}\quad\MM,
\end{split}
\eea
and 
\bea\lab{eq:esitmatefor15weightedderivativesoffirstorderdrivgivescheckoftJkandtJkpm:nab3onlyrGagbutnotGab} 
|\dk^{\leq 15}(\widecheck{\nab_3\tJk}, \widecheck{\nab_3\tJk_{\pm}})|\les\frac{\ep}{r^2\tau^{\frac{1}{2}+\dec}}, \qquad |\dk^{\leq 15}(\widecheck{\nab_3\tJk}, \widecheck{\nab_3\tJk_{\pm}})| &\les\frac{\ep}{r\tau^{1+\frac{3\dec}{4}}},\quad\textrm{on}\quad\MM.
\eea
Note that \eqref{eq:esitmatefor15weightedderivativesoffirstorderdrivgivescheckoftJkandtJkpm:goodnab4andnabbutnotbab3} implies the stated estimates for $(\widecheck{\nab_4\tJk}, \widecheck{\nab_4\tJk_{\pm}})$ and $(\widecheck{\nab\tJk}, \widecheck{\nab\tJk_{\pm}})$ in \eqref{eq:esitmatefor15weightedderivativesoffirstorderdrivgivescheckoftJkandtJkpm}. On the other hand, we need to improve on the estimate  \eqref{eq:esitmatefor15weightedderivativesoffirstorderdrivgivescheckoftJkandtJkpm:nab3onlyrGagbutnotGab} for $(\widecheck{\nab_3\tJk}, \widecheck{\nab_3\tJk_{\pm}})$. To this end, we rely on the vectorfield $\T'$ introduced in \eqref{eq:defonT'expressedindoubleprimedframe} and commute the transport equations in \eqref{definitionofthenewcomplex1formstJktJkpm} w.r.t. $\Lieb_{\T'}$. In view of the first commutator formula in  \eqref{eq:commutatorsYtauwithT'andrnab'}, we obtain on $\MM$
\beaa
\bsplit
\nab_{\Ytau}'\Lieb_{\T'}\tJk =& -\frac{1}{r}\Lieb_{\T'}\tJk+O(mr^{-2})\Lieb_{\T'}\tJk+O(r^{-3})\T'(\rr)+O(r^{-3})\T'(\tx^b)+\Ga_g'\nab_3'\tJk+r^{-1}\Ga_b'\c\dk \tJk+r^{-2}\dk\Ga_b',\\
\nab_{\Ytau}'\Lieb_{\T'}\tJk_{\pm} =&  -\frac{1}{r}\Lieb_{\T'}\tJk_{\pm}+O(mr^{-2})\Lieb_{\T'}\tJk_{\pm}+O(mr^{-2})\Lieb_{\T'}\tJk_{\mp}+O(r^{-3})\T'(\rr)+O(r^{-3})\T'(\tx^b)\\
&+\Ga_g'\nab_3'\tJk_{\pm}+r^{-1}\Ga_b'\c\dk\tJk_{\pm}+r^{-2}\dk\Ga_b' 
\end{split}
\eeaa
and hence
\beaa
\bsplit
\nab_{\Ytau}'(r\Lieb_{\T'}\tJk) =& O(mr^{-2})r\Lieb_{\T'}\tJk+O(r^{-2})\T'(\rr)+O(r^{-2})\T'(\tx^b)+\Ga_g'r\nab_3'\tJk+\Ga_b'\c\dk \tJk+r^{-1}\dk\Ga_b',\\
\nab_{\Ytau}'(r\Lieb_{\T'}\tJk_{\pm}) =&  O(mr^{-2})r\Lieb_{\T'}\tJk_{\pm}+O(mr^{-2})r\Lieb_{\T'}\tJk_{\mp}+O(r^{-2})\T'(\rr)+O(r^{-2})\T'(\tx^b)\\
&+\Ga_g'r\nab_3'\tJk_{\pm}+\Ga_b'\c\dk\tJk_{\pm}+r^{-1}\dk\Ga_b'. 
\end{split}
\eeaa
Integrating from $\Si_*$, and using \eqref{eq:esitmatefor15weightedderivativesoffirstorderdrivgivescheckoftJkandtJkpm:goodnab4andnabbutnotbab3}, \eqref{eq:esitmatefor15weightedderivativesoffirstorderdrivgivescheckoftJkandtJkpm:nab3onlyrGagbutnotGab}, \eqref{eq:assumptionsonMMforpartII:secondglobalframeauxassfor15derivatives}, as well as \eqref{eq:esitmatefor15weightedderivativesoffirstorderdrivgivescheckofrr} and \eqref{eq:esitmatefor15weightedderivativesoffirstorderdrivgivescheckoftxb}, we infer
\beaa
|(\Lieb_{\T'}\tJk, \Lieb_{\T'}\tJk_{\pm})| &\les& \frac{1}{r}\int_r^{r_*}\frac{\ep}{{r'}^2\tau^{1+\frac{3\dec}{4}}}\les\frac{\ep}{r^2\tau^{1+\frac{3\dec}{4}}},\quad\textrm{on}\quad\MM.
\eeaa
Commuting the above transport equations for $\Lieb_{\T'}\tJk$ and $\Lieb_{\T'}\tJk_{\pm}$, we then extend the above estimates to higher order derivatives to obtain 
\beaa
|\dk^{\leq 15}(\Lieb_{\T'}\tJk, \Lieb_{\T'}\tJk_{\pm})| &\les& \frac{\ep}{r^2\tau^{1+\frac{3\dec}{4}}},\quad\textrm{on}\quad\MM.
\eeaa
Together with \eqref{eq:esitmatefor15weightedderivativesoffirstorderdrivgivescheckoftJkandtJkpm:goodnab4andnabbutnotbab3}, \eqref{eq:esitmatefor15weightedderivativesoffirstorderdrivgivescheckoftJkandtJkpm:nab3onlyrGagbutnotGab}, the definition of $\T'$ and Lemma \ref{lemma:basicpropertiesLiebTfasdiuhakdisug:chap9}, we deduce 
\beaa
\bsplit
|\dk^{\leq 15}(\widecheck{\nab_4\tJk}, \widecheck{\nab_4\tJk_{\pm}})|&\les\frac{\ep}{r^3\tau^{\frac{1}{2}+\dec}}, \quad\textrm{on}\quad\MM,\\ 
|\dk^{\leq 15}(\widecheck{\nab_4\tJk}, \widecheck{\nab_4\tJk_{\pm}}, \widecheck{\nab\tJk_{\pm}}, \widecheck{\nab\tJk_{\pm}}, \widecheck{\nab_3\tJk}, \widecheck{\nab_3\tJk_{\pm}})| &\les\frac{\ep}{r^2\tau^{1+\frac{3\dec}{4}}},\quad\textrm{on}\quad\MM,
\end{split}
\eeaa
as stated in \eqref{eq:esitmatefor15weightedderivativesoffirstorderdrivgivescheckoftJkandtJkpm}. 

Finally, it remains to establish the identities \eqref{eq:usefulalgebraicidentitiesinvolvingscalarproductsReJkReJkpm:Kerrpert:tilde}. To this end, we first prove that these identities hold on $\Si_*$. We start by computing scalar products satisfied by the 1-forms $f_0$ and $f_\pm$ on $\Si_*$. First, we have in view of the definition \eqref{eq:definitionofthereal1formsf0fpmonS*} of $f_0$ and $f_\pm$ on $S_*$
\bea\lab{eq:scalarproductreal1formsf0fpmonSi*}
\bsplit
f_0\c f_0 &=(\sin\th)^2, \qquad f_0\c f_+=-\sin\th\sin\tphi, \qquad \dual f_0\c f_+=\sin\th\cos\th\cos\tphi, \\
f_0\c f_- &=\sin\th\cos\tphi, \qquad \dual f_0\c f_-=\sin\th\cos\th\sin\tphi, \\
f_+\c f_+&=(\cos\th)^2(\cos\tphi)^2+(\sin\tphi)^2,\qquad f_-\c f_- =(\cos\th)^2(\sin\tphi)^2+(\cos\tphi)^2.
\end{split}
\eea
In view of the transport equations \eqref{eq:corrdinatesthandvarphiarepropagatedalongSi*bynu*} and \eqref{eq:extendionofthereal1formsf0fpmfromS*toSigma*} along $\Si_*$ in the direction $\nu_*$, the identities also hold on $\Si_*$. Together with the definition \eqref{eq:relationbetweenJkJkpmandf0onSigmastar} of $\Jk$ and $\Jk_{\pm}$ on $\Si_*$, we immediately infer the following identities on $\Si_*$
\beaa
\bsplit
\Re(\Jk)\c\Re(\Jk) &=\frac{(\sin(\th))^2}{|q|^2}, \qquad \Re(\Jk)\c\Re(\Jk_+)=-\frac{x^2_p}{|q|^2},\qquad \Re(\Jk)\c\Re(\Jk_-)=\frac{x^1_p}{|q|^2},\\
\dual(\Re(\Jk))\c\Re(\Jk_+) &=\frac{\cos(\th)x^1_p}{|q|^2},\qquad \dual(\Re(\Jk))\c\Re(\Jk_-)=\frac{\cos(\th)x^2_p}{|q|^2}.
\end{split}
\eeaa
In view of the initialization of $\tJk$, $\tJk_{\pm}$, $\rr$ and $\tx^b$, $b=1,2$ on $\Si_*$ respectively in \eqref{definitionofthenewcomplex1formstJktJkpm}, \eqref{definitionofthenewcoordinatetilder} and
\eqref{definitionofthenewcoordinatetildethetandtphi}, this implies that the identities \eqref{eq:usefulalgebraicidentitiesinvolvingscalarproductsReJkReJkpm:Kerrpert:tilde} hold on $\Si_*$. They are then immediately transported\footnote{This is due to the fact that the identities \eqref{eq:usefulalgebraicidentitiesinvolvingscalarproductsReJkReJkpm:Kerrpert:tilde} are true in Kerr, see \eqref{eq:usefulalgebraicidentitiesinvolvingscalarproductsReJkReJkpm}, and that the transport equations \eqref{definitionofthenewcomplex1formstJktJkpm}, \eqref{definitionofthenewcoordinatetilder} and
\eqref{definitionofthenewcoordinatetildethetandtphi} in $\Ytau$ agree with the ones in Kerr.} to $\MM$ in view of the transport equations in $\Ytau$ for $\tJk$, $\tJk_{\pm}$, $\rr$ and $\tx^b$, $b=1,2$ respectively in \eqref{definitionofthenewcomplex1formstJktJkpm}, \eqref{definitionofthenewcoordinatetilder} and
\eqref{definitionofthenewcoordinatetildethetandtphi}. For instance, the first identity in \eqref{eq:usefulalgebraicidentitiesinvolvingscalarproductsReJkReJkpm:Kerrpert:tilde} satisfies the following transport equation 
\beaa
\Ytau\left(\Re(\tJk)\c\Re(\tJk)-\frac{(\sin(\widetilde{\th}))^2}{|\tq|^2}\right) &=& -\frac{2\rr}{|\tq|^2}\big(\De(\rr)[\Ytau^4|q|^{-2}]_K(\rr)-[\Ytau^3]_K(\rr)\big)\left(\Re(\tJk)\c\Re(\tJk)-\frac{(\sin(\widetilde{\th}))^2}{|\tq|^2}\right)
\eeaa
where we used the fact that $\Re(\tJk)\c\dual\Re(\tJk)=0$, while the second identity in \eqref{eq:usefulalgebraicidentitiesinvolvingscalarproductsReJkReJkpm:Kerrpert:tilde} satisfies the following transport equation  
\beaa
&&\Ytau\left(\Re(\tJk)\c\Re(\tJk_+)+\frac{\tx^2_p}{|\tq|^2}\right)\\ 
&=& -\frac{2\rr}{|\tq|^2}\big(\De(\rr)[\Ytau^4|q|^{-2}]_K(\rr) - [\Ytau^3]_K(\rr)\big)\left(\Re(\tJk)\c\Re(\tJk_+)+\frac{\tx^2_p}{|\tq|^2}\right)\\
&& - \Big(\phimod'(r)[\Ytau^3]_K(\rr)+\big(2a-\De\phimod'(r)\big)[\Ytau^4|q|^{-2}]_K(\rr)\Big)\left(\Re(\tJk)\c\Re(\tJk_-)-\frac{\tx^1_p}{|\tq|^2}\right).
\eeaa
This concludes the proof of Lemma \ref{lemma:controlofthecomplex1formswidetildeJkandwidetildeJkplusminus}.
\end{proof}

%%%%%%%%%%%%%%%%%%%%%%%%%%%%%%%%%%%%%%%%%%%%%%%%%%%%%%%%%%%

\subsubsection{Notation $(\widetilde{\Ga}_g, \widetilde{\Ga}_b)$ associated to $(\rr, \tx^1, \tx^2)$ and $(\tJk, \tJk_{\pm})$}

%%%%%%%%%%%%%%%%%%%%%%%%%%%%%%%%%%%%%%%%%%%%%%%%%%%%%%%%%%%

In this section, we introduce a new notation $(\widetilde{\Ga}_g, \widetilde{\Ga}_b)$ for error terms. To this end, we first introduce linearized quantities with Kerr values computed w.r.t. $(\rr, \tx^1, \tx^2)$ and $(\tJk, \tJk_{\pm})$.

\begin{definition}
\lab{def:renormalizationofallnonsmallquantitiesinPGstructurebyKerrvalue:widetildecase}
We  define  the following renormalizations.
\begin{enumerate}
\item Linearization of the complex-valued Ricci and curvature coefficients:
\beaa
\bsplit
\widecheck{\trXc} &:= \tr X-\frac{2\ov{\tq}\De(\rr)}{|\tq|^4}, \qquad\widecheck{\trXbc} := \tr\Xb+\frac{2}{\ov{\tq}},\\ 
\widecheck{\Pc} &:= P+\frac{2m}{\tq^3},\qquad\qquad\, \widecheck{\omc}  := \om  + \frac{1}{2}\pr_{\rr}\left(\frac{\De(\rr)}{|\tq|^2} \right),\\
\widecheck{\Hc} &:= H-\frac{a\tq}{|\tq|^2}\tJk, \qquad\quad \widecheck{\Hbc}:=\Hb+\frac{a\ov{\tq}}{|\tq|^2}\tJk,\qquad\quad \widecheck{\Zc} := Z-\frac{a\tq}{|\tq|^2}\tJk.
 \end{split}
\eeaa

\item Linearization of derivatives of the scalar functions $\rr$, $\cos(\widetilde{\th})$, $\tq$, $\ttt$, $\tx^1_p$ and $\tx^2_p$:
\beaa
\bsplit
\widecheck{\widecheck{e_3(\rr)}} :=& e_3(\rr)+1, \qquad\qquad\qquad\,\,\,\,\, \widecheck{\widecheck{e_4(\rr)}} := e_4(\rr)-\frac{\Delta(\rr)}{|\tq|^2},\\
\widecheck{\widecheck{e_3(\ttt)}}:=& e_3(\ttt)-\tmod'(\rr), \qquad\qquad\, \widecheck{\widecheck{e_4(\ttt)}}:=e_4(\ttt)-\frac{2(\rr^2+a^2) - \De(\rr)\tmod'(\rr)}{|\tq|^2},\\
\widecheck{\widecheck{e_3(\tx^1_p)}}:=& e_3(\tx^1_p)+\phimod'(\rr)\tx^2_p, \qquad \widecheck{\widecheck{e_4(x^1_p)}}:=  e_4(\tx^1_p)+\frac{2a -\De(\rr)\phimod'(\rr)}{|\tq|^2}\tx^2_p, \\ 
\widecheck{\widecheck{e_3(\tx^2_p)}}:=& e_3(\tx^2_p)-\phimod'(\rr)\tx^1_p, \qquad \widecheck{\widecheck{e_4(\tx^2_p)}}:=e_4(\tx^2_p)-\frac{2a -\De(\rr)\phimod'(\rr)}{|\tq|^2}\tx^1_p,\\
\widecheck{\widecheck{\DD\tq}} :=& \DD\tq+a\tJk, \qquad\qquad\, \widecheck{\widecheck{\DD \ov{\tq}}} :=\DD \ov{\tq}-a\tJk,\qquad \,\,\,\widecheck{\widecheck{\DD(\cos(\widetilde{\th}))}} := \DD(\cos(\widetilde{\th})) -i\tJk,\\
\widecheck{\widecheck{\DD(\ttt)}}:=& \DD(\ttt)-a\tJk, \qquad  \widecheck{\widecheck{\DD(\tx^1_p)}} := \DD(\tx^1_p) - \tJk_+, \qquad \widecheck{\widecheck{\DD(\tx^2_p)}} := \DD(\tx^2_p) - \tJk_-.
\end{split}
\eeaa

\item Linearization of derivatives of the complex 1-forms $\tJk$ and $\tJk_{\pm}$:
\beaa
\bsplit
\widecheck{\widecheck{\nab_3\tJk}}:=&\nab_3\tJk -\frac{1}{\ov{\tq}}\tJk, \qquad\qquad  \widecheck{\widecheck{\nab_4\tJk}}:=\nab_4\tJk +\frac{\De(\rr) \ov{\tq}}{|\tq|^4}\tJk,\qquad \widecheck{\widecheck{\ov{\DD}\c\tJk}}:= \ov{\DD}\c\tJk-\frac{4i(\rr^2+a^2)\cos(\widetilde{\th})}{|\tq|^4},\\
\widecheck{\widecheck{\nab_3\tJk_\pm}}:=&\nab_3\tJk_\pm - \frac{1}{\ov{\tq}}\tJk_\pm \pm\phimod'(\rr)\tJk_{\mp},\qquad \widecheck{\widecheck{\nab_4 \tJk_\pm}}:=\nab_4 \tJk_\pm + \frac{\De(\rr) \ov{\tq}}{|\tq|^4}\tJk_{\pm} \pm\frac{2a-\De(\rr)\phimod'(\rr)}{|\tq|^2}\tJk_{\mp},\\
\widecheck{\widecheck{\ov{\DD}\c \tJk_+}}:=&\ov{\DD}\c \tJk_+ + \frac{4\rr^2 }{|\tq|^4}\tx^1_p + \frac{4ia^2\cos(\widetilde{\th})}{|\tq|^4}\tx^2_p,\qquad \widecheck{\widecheck{\ov{\DD}\c \tJk_-}}:=\ov{\DD}\c \tJk_- + \frac{4\rr^2 }{|\tq|^4}\tx^2_p - \frac{4ia^2\cos(\widetilde{\th})}{|\tq|^4}\tx^1_p.
\end{split}
\eeaa
 \end{enumerate}
\end{definition}

\begin{remark}
In order to avoid conflict of notations concerning linearized quantities, we denote by one check quantities linearized using $(r, x^1, x^2)$ and $(\Jk, \Jk_{\pm})$, see Definition \ref{def:renormalizationofallnonsmallquantitiesinPGstructurebyKerrvalue}, and by two checks the ones linearized using $(\rr, \tx^1, \tx^2)$ and $(\tJk, \tJk_{\pm})$, see Definition 
\ref{def:renormalizationofallnonsmallquantitiesinPGstructurebyKerrvalue:widetildecase}.
\end{remark}

We are now ready to define the quantities $(\widetilde{\Ga}_g, \widetilde{\Ga}_b)$.
\begin{definition}
\lab{definition.Ga_gGa_b:widetildecase}
The set of all linearized quantities in Definition \ref{def:renormalizationofallnonsmallquantitiesinPGstructurebyKerrvalue:widetildecase} is of the form $\widetilde{\Ga}_g\cup \widetilde{\Ga}_b$ with  $\widetilde{\Ga}_g,  \widetilde{\Ga}_b$
 defined as follows.
 \begin{enumerate}
\item 
 The set $\widetilde{\Ga}_g$ is given by $\widetilde{\Ga}_g=\widetilde{\Ga}_{g,1}\cup \widetilde{\Ga}_{g, 2}\cup\widetilde{\Ga}_{g,3}$   with
 \bea
 \bsplit
 \widetilde{\Ga}_{g,1} &= \Big\{\Xi, \quad \widecheck{\omc}, \quad\widecheck{\trXc},\quad  \Xh,\quad \widecheck{\Zc},\quad \widecheck{\Hbc}, \quad \widecheck{\trXbc} , \quad r\widecheck{\Pc}, \quad  rB, \quad  rA\Big\},\\
 \widetilde{\Ga}_{g,2} &= \Big\{\widecheck{\widecheck{e_4(\rr)}}, \,\,\,\, r^{-1}\nab(\rr), \,\,\,\, r^{-1}\widecheck{\widecheck{e_4(\ttt)}}, \,\,\,\, r^{-1}\widecheck{\widecheck{\DD(\ttt)}}, \,\,\,\, r^{-1}\widecheck{\widecheck{e_3(\ttt)}},  \,\,\,\, e_4(\cos(\widetilde{\th})), \,\,\,\, \widecheck{\widecheck{\DD(\cos(\widetilde{\th}))}},\\
 &\qquad \widecheck{\widecheck{e_4(\tx^1_p)}}, \,\,\,\,   \widecheck{\widecheck{\DD(\tx^1_p)}}, \,\,\,\,   \widecheck{\widecheck{e_4(\tx^2_p)}},\,\,\,\,  \widecheck{\widecheck{\DD(\tx^2_p)}}\Big\},\\
  \widetilde{\Ga}_{g,3} &= \Big\{r\widecheck{\widecheck{\nab_4\tJk},} \quad r\widecheck{\widecheck{\nab_4\tJk_\pm}}\Big\}.
 \end{split}
 \eea
 
 \item The set $\widetilde{\Ga}_b$ is given by $\widetilde{\Ga}_b=\widetilde{\Ga}_{b,1}\cup \widetilde{\Ga}_{b, 2}\cup \widetilde{\Ga}_{b,3}$   with
 \bea
 \bsplit
 \widetilde{\Ga}_{b,1}&= \Big\{\widecheck{\Hc}, \quad \Xbh, \quad \omb, \quad \Xib,\quad  r\Bb, \quad \Ab\Big\},\\
  \widetilde{\Ga}_{b, 2}&= \Big\{r^{-1}\widecheck{\widecheck{e_3(\rr)}}, \quad e_3(\cos(\widetilde{\th})), \quad \widecheck{\widecheck{e_3(\tx^1_p)}}, \quad \widecheck{\widecheck{e_3(\tx^2_p)}}\Big\}, \\
   \widetilde{\Ga}_{b,3}&=\bigg\{ r\,\widecheck{\widecheck{\ov{\DD}\c\tJk}}, \quad r\,\DD\hot\tJk, \quad r\,\widecheck{\widecheck{\nab_3\tJk}}, \quad r\,\widecheck{\widecheck{\ov{\DD}\c\tJk_\pm}}, \quad r\,\DD\hot\tJk_\pm, \quad r\,\widecheck{\widecheck{\nab_3\tJk_\pm}}\bigg\}. 
   \end{split}
 \eea
\end{enumerate}
\end{definition}

 As a consequence of \eqref{eq:assumptionsonMMforpartII}, \eqref{eq:controloftildeuonMext:ealphaoftttcheckinrGag}, and Lemmas \ref{lemma:controlofnewcoordinatewidetilder}, \ref{lemma:gainforwidechecknabJp} and \ref{lemma:controlofthecomplex1formswidetildeJkandwidetildeJkplusminus}, $(\widetilde{\Ga}_{g}, \widetilde{\Ga}_{b})$ introduced in Definition \ref{definition.Ga_gGa_b:widetildecase} satisfy
\bea\lab{eq:decaypropertiesofwidehatGab1widehatGag1}
|\dk^{\leq 15}\widetilde{\Ga}_{g}|\les \ep\min\Big\{r^{-2}\tau^{-\frac{1+\dec}{2}}, \, r^{-1}\tau^{-1-\frac{3\dec}{4}}\Big\}, \quad\,\,\, |\dk^{\leq 15}\widetilde{\Ga}_{b}|\les \ep r^{-1}\tau^{-1-\frac{3\dec}{4}}\,\,\,\textrm{on}\,\,\,\MM,
\eea
and we also have
\bea\lab{eq:additionestimatewidechecke4tttbyeprminus2:withdoublewidecheck}
|\dk^{\leq 15}\xi|\les\frac{\ep}{r^3}, \qquad |\dk^{\leq 15}\widecheck{\widecheck{e_4(\ttt)}}|\les\frac{\ep}{r^2}.
\eea

\begin{remark}
We list below the additional properties of the new coordinates system $(\ttt, \rr, \tx^1, \tx^2)$ and the new complex 1-forms $(\tJk, \tJk_{\pm})$ compared to $(\tt, r, x^1, x^2)$ and $(\Jk, \Jk_{\pm})$:
\begin{itemize}
\item $\widecheck{\widecheck{e_3(\ttt)}}\in r\widetilde{\Ga}_g$, while we only have $\widecheck{e_3(\tt)}\in\Ga_b$.

\item $\widecheck{\widecheck{\nab(\tx^b)}}\in\widetilde{\Ga}_g$, while we only have $\widecheck{\nab(x^b)}\in\Ga_b$.

\item The coordinates $(\rr, \tx^1, \tx^2)$ satisfy explicit transport equations along $\Ytau$, see \eqref{definitionofthenewcoordinatetilder} \eqref{definitionofthenewcoordinatetildethetandtphi}.

\item The complex 1-forms $(\tJk, \tJk_{\pm})$ satisfy the identities \eqref{eq:usefulalgebraicidentitiesinvolvingscalarproductsReJkReJkpm:Kerrpert:tilde}.
\end{itemize}
On the other hand, $\widecheck{\widecheck{e_4(\ttt)}}\in r\widetilde{\Ga}_g$ with $\widecheck{\widecheck{e_4(\ttt)}}$ satisfying \eqref{eq:additionestimatewidechecke4tttbyeprminus2:withdoublewidecheck} which is weaker than the corresponding estimate $\widecheck{e_4(\tt)}\in \Ga_g$. 
\end{remark}

%%%%%%%%%%%%%%%%%%%%%%%%%%%%%%%%%%%%%%%%%%%%%%%%%%%%%%%%%%%

\subsubsection{Metric coefficients in the coordinates system $(\ttt, \rr, \tx^1, \tx^2)$}

%%%%%%%%%%%%%%%%%%%%%%%%%%%%%%%%%%%%%%%%%%%%%%%%%%%%%%%%%%%
 
We now consider the coordinates system $(\ttt, \rr, \tx^1, \tx^2)$ with $\ttt$, $\rr$ and $(\tx^1, \tx^2)$ introduced respectively in Lemma \ref{lemma:controloftildeuonMext}, and in \eqref{definitionofthenewcoordinatetilder} and 
\eqref{definitionofthenewcoordinatetildethetandtphi} and we linearize the components of the inverse metric in this coordinates system as follows
\bea
\widecheck{\widecheck{\g}}^{\a\b} := \g^{\a\b} - \g^{\a\b}_{Kerr},
\eea
where $\g^{\a\b}_{Kerr}$ denotes the corresponding value in Kerr which is an explicit scalar function only depending on $(\rr, \widetilde{\th})$ and $(a, m)$. Then, we have the following control for the linearized inverse metric coefficients.

\begin{lemma}[Control of the linearized inverse metric coefficients]
\lab{lemma:controloflinearizedinversemetriccoefficients}
In the coordinates system $(\ttt, \rr, \tx^1, \tx^2)$, recalling the notation $(\widetilde{\Ga}_{g}, \widetilde{\Ga}_{b})$ induced in Definition \ref{definition.Ga_gGa_b:widetildecase}, we have 
\bea\lab{eq:controloflinearizedinversemetriccoefficients}
\widecheck{\widecheck{\g}}^{\rr\rr}=r\widetilde{\Ga}_{b}, \quad \widecheck{\widecheck{\g}}^{\rr\ttt}=r\widetilde{\Ga}_{g},\quad \widecheck{\widecheck{\g}}^{\ttt\ttt}=r\widetilde{\Ga}_{g}, \quad \widecheck{\widecheck{\g}}^{\rr a}=\widetilde{\Ga}_{b}, \quad \widecheck{\widecheck{\g}}^{\ttt a}=\widetilde{\Ga}_{g},\quad \widecheck{\widecheck{\g}}^{ab}=r^{-1}\widetilde{\Ga}_{g},
\eea 
and
\bea\lab{eq:controloflinearizedinversemetriccoefficients:inversegtautau}
|\dk^{\leq 15}\widecheck{\widecheck{\g}}^{\ttt\ttt}|\les \frac{\ep}{r^2}.
\eea 
\end{lemma}

\begin{remark}
Using the coordinates system $(\ttt, r, \tx^1, \tx^2)$ instead of $(\tt, \rr, x^1, x^2)$ allows in particular a gain one power of $r$ for $\widecheck{\widecheck{\g}}^{\rr\ttt}$ and $\widecheck{\widecheck{\g}}^{ab}$ due to the improved estimates for $\widecheck{\widecheck{e_3(\ttt)}}$ and $\widecheck{\widecheck{\nab(\tx^b)}}$ compared to $\widecheck{e_3(\tt)}$ and $\widecheck{\nab(\cos\th)}$. 
\end{remark}

\begin{proof}
In a general coordinate system $(x^\a)$, we have 
\beaa
\g^{\a\b} &=& -\frac{1}{2}e_4(x^\a)e_3(x^\b)-\frac{1}{2}e_3(x^\a)e_4(x^\b)+\nab(x^\a)\c\nab(x^\b),
\eeaa
see Lemma 4.1.3 in \cite{KS:Kerr}. In the $(\ttt, \rr, \tx^1, \tx^2)$ coordinates system, since we have, in view of Definition \ref{definition.Ga_gGa_b:widetildecase} 
\beaa
e_4(\rr)=1+\widetilde{\Ga}_{g}, \quad e_4(\ttt)=O(r^{-2}),\quad e_4(\tx^b)=O(r^{-2}), \,\, b=1,2,\quad \nab(\rr)\in r\widetilde{\Ga}_{g}, 
\eeaa
and 
\beaa
e_3(\rr), \, e_3(\ttt)=O(1), \quad e_3(\tx^b)=O(r^{-1}), \,\, b=1,2, \quad \nab(\ttt)=O(r^{-1}), \quad \nab(\tx^b)=O(r^{-1}),  \,\, b=1,2,
\eeaa
we infer
\beaa
\widecheck{\widecheck{\g}}^{\rr\rr} &=& -\widecheck{\widecheck{e_3(\rr)}}+\widetilde{\Ga}_{g},\\
\widecheck{\widecheck{\g}}^{\rr\ttt} &=& -\frac{1}{2}\widecheck{\widecheck{e_3(\ttt)}}s +O(r^{-2})\widecheck{\widecheck{e_3(\rr)}}+O(1)\widecheck{\widecheck{e_4(\ttt)}}+\widetilde{\Ga}_{g},\\
\widecheck{\widecheck{\g}}^{\ttt\ttt} &=& O(r^{-2})\widecheck{\widecheck{e_3(\ttt)}}+O(1)\widecheck{\widecheck{e_4(\ttt)}}+O(r^{-1})\widecheck{\widecheck{\nab(\ttt)}}+(r\widetilde{\Ga}_{g})^2,\\
\widecheck{\widecheck{\g}}^{\rr b} &=& -\frac{1}{2}\widecheck{\widecheck{e_3(\tx^b)}}+O(1)\widecheck{\widecheck{e_4(\tx^b)}}+\widetilde{\Ga}_{g}, \quad b=1,2,\\
\widecheck{\widecheck{\g}}^{\ttt b} &=& O(r^{-2})\widecheck{\widecheck{e_3(\ttt)}}+O(r^{-2})\widecheck{\widecheck{e_3(\tx^b)}} +O(r^{-1})\widecheck{\widecheck{e_4(\ttt)}}+O(1)\widecheck{\widecheck{e_4(\tx^b)}}+O(r^{-1})\widecheck{\widecheck{\nab(\ttt)}}\\
&& +O(r^{-1})\widecheck{\widecheck{\nab(\tx^b)}}+r\widetilde{\Ga}_{g}\c\widetilde{\Ga}_{b},\quad b=1,2.\\
\widecheck{\widecheck{\g}}^{bc} &=& O(r^{-2})\widecheck{\widecheck{e_3(\tx^b)}}+O(r^{-2})\widecheck{\widecheck{e_3(\tx^c)}}+O(r^{-1})\widecheck{\widecheck{e_4(\tx^b)}}+O(r^{-1})\widecheck{\widecheck{e_4(\tx^c)}}+O(r^{-1})\widecheck{\widecheck{\nab(\tx^b)}}\\
&&+O(r^{-1})\widecheck{\widecheck{\nab(\tx^c)}}+\widetilde{\Ga}_{g}\c\widetilde{\Ga}_{b},\quad b,c=1,2,
\eeaa
and hence, using again Definition \ref{definition.Ga_gGa_b:widetildecase}, 
\beaa
\widecheck{\widecheck{\g}}^{\rr\rr}=r\widetilde{\Ga}_{b}, \quad \widecheck{\widecheck{\g}}^{\rr\ttt}=r\widetilde{\Ga}_{g},\quad \widecheck{\widecheck{\g}}^{\ttt\ttt}=r\widetilde{\Ga}_{g}, \quad \widecheck{\widecheck{\g}}^{\rr a}=\widetilde{\Ga}_{b}, \quad \widecheck{\widecheck{\g}}^{\ttt a}=\widetilde{\Ga}_{g},\quad \widecheck{\widecheck{\g}}^{ab}=r^{-1}\widetilde{\Ga}_{g},
\eeaa 
as stated in \eqref{eq:controloflinearizedinversemetriccoefficients}. Also, using again the above identity for $\widecheck{\widecheck{\g}}^{\ttt\ttt}$, we have
\beaa
\widecheck{\widecheck{\g}}^{\ttt\ttt} &=& O(r^{-2})\widecheck{\widecheck{e_3(\ttt)}}+O(1)\widecheck{\widecheck{e_4(\ttt)}}+O(r^{-1})\widecheck{\widecheck{\nab(\ttt)}}+(r\widetilde{\Ga}_{g})^2\\
&=& O(1)\widecheck{\widecheck{e_4(\ttt)}}+\widetilde{\Ga}_{g}
\eeaa
which together with \eqref{eq:additionestimatewidechecke4tttbyeprminus2:withdoublewidecheck} and \eqref{eq:decaypropertiesofwidehatGab1widehatGag1} implies \eqref{eq:controloflinearizedinversemetriccoefficients:inversegtautau}. This concludes the proof of Lemma \ref{lemma:controloflinearizedinversemetriccoefficients}.
\end{proof}

%%%%%%%%%%%%%%%%%%%%%%%%%%%%%%%%%%%%%%

\subsection{Regular scalarization of tensorial wave equations}
\lab{sec:regularscalarization}
  
%%%%%%%%%%%%%%%%%%%%%%%%%%%%%%%%%%%%%%

In this section, we start by recalling the regular scalarization procedure of tensorial wave equations introduced in Section 3 of \cite{MaSz26}, which relies on the notion of regular triplets $\Om_i$, $i=1,2,3$, see Definition \ref{def:definitionofregulartripletOmii=123}. The goal of the section in then to construct a suitable regular triplet $\Om_i$, $i=1,2,3$ in $\MM$ in Section \ref{sec:defandmainpropofregtripletOmiinKerrpert}.

%%%%%%%%%%%%%%%%%%%%%%%%%%%%%%%%%%%%%%%%%%%%

\subsubsection{Regular scalarization of tensorial wave equations in \cite{MaSz26}}
\lab{sec:regularscalarizationinMaSz26}

%%%%%%%%%%%%%%%%%%%%%%%%%%%%%%%%%%%%%%%%%%%%

In this section, we recall the regular scalarization procedure for tensorial wave equations introduced in \cite{MaSz26}.

%%%%%%%%%%%%%%%%%%%%%%%%%%%%%%%%%%%%%%%%%%

\paragraph{\textit{Scalarization using a regular triplet $\Om_i$, $i=1,2,3$}.}

%%%%%%%%%%%%%%%%%%%%%%%%%%%%%%%%%%%%%%%%%%

In order to scalarize horizontal tensors, we will rely on the following definition, see Definition 3.1 in \cite{MaSz26}. 

\begin{definition}[Regular triplet]
\lab{def:definitionofregulartripletOmii=123}
Let $(\MM, \g)$ be a spacetime, $(e_3, e_4)$ be a null pair, and consider the corresponding horizontal structure $\O(\MM)$ introduced in Section \ref{subsection:review-horiz.structures}. We say that vectorfields $\Om_i$, $i=1,2,3$, identified with elements of $\sk_1$, form a regular triplet if they are regular and satisfy the following identities 
\bea\lab{eq:fundamentalpropertiesof1formsOmi}
x^i\Om_i=0, \qquad (\Om^i)_a(\Om_i)_b=\de_{ab}, \qquad \Om_i\c\Om_j=\de_{ij}-x^ix^j, \qquad \Om_i\c\dual\Om_j=\in_{ijk}x^k,
\eea
where, by convention, we denote $\Om^i=\Om_i$, i.e., the $i$-index is lowered or raised using $\de_{ij}$ or $\de^{ij}$.
\end{definition}

\begin{remark}\lab{rmk:generalcontructionofregulartripletsingivenspacetime}
For a specific choice of a regular triplet in Kerr, see Definition \ref{def:regulartripletinKerrOmii=123}, and for the construction of a regular triplet in $\MM$, see Section \ref{sec:defandmainpropofregtripletOmiinKerrpert}.
\end{remark}

Next, we introduce the following 1-forms on $\MM$.
\begin{definition}
\lab{def:Mialphaj:Kerr}
Let $(\MM, \g)$ be a spacetime, $(e_3, e_4)$ be a null pair, and consider the corresponding horizontal structure $\O(\MM)$ introduced in Section \ref{subsection:review-horiz.structures}. Let $\Om_i$, $1,2,3$ be a regular triplet in the sense of Definition \ref{def:definitionofregulartripletOmii=123}. We  define the following 1-forms on $\MM$ 
\bea\lab{eq:definitionofMalphaijwithoutambiguity}
M_{i\a}^j:=(\Ddot_\a\Om_i)\c\Om^j, \quad \forall \a,i,j.
\eea
Further, we define $M_{i}^{j \a}:=\g^{\a\b} M_{i\b}^j$.
\end{definition}

\begin{lemma}\lab{lemma:introductionandpropertiesoftheMalphaij}
Let $M_{i\a}^j$ be the 1-forms on $(\MM, \g)$ as defined in Definition \ref{def:Mialphaj:Kerr}.  Then we have
\bea
\label{def:Mialphaj}
\Ddot_\a\Om_i=M_{i\a}^j\Om_j.
\eea
\end{lemma}

\begin{proof}
See Lemma 3.4 in \cite{MaSz26}.
\end{proof}

%%%%%%%%%%%%%%%%%%%%%%%%%%%%%%%%

\paragraph{\textit{From tensors to regular scalars and back}.}

%%%%%%%%%%%%%%%%%%%%%%%%%%%%%%%%

The following lemma allows to pass from horizontal tensors in $\sk_2(\mathbb{C})$ to scalars and reciprocally.
\begin{lemma}\lab{lemma:backandforthbetweenhorizontaltensorsk2andscalars:complex}
Let $(\MM, \g)$ be a spacetime, $(e_3, e_4)$ be a null pair, and consider the corresponding horizontal structure $\O(\MM)$ introduced in Section \ref{subsection:review-horiz.structures}. Assume that $\Om_i$, $1,2,3$ is a regular triplet in the sense of Definition \ref{def:definitionofregulartripletOmii=123}. Then, the following holds:
\begin{enumerate}
\item\lab{item1:sk2Csatisfyconditions} Let $\pmb\psi\in\sk_2(\mathbb{C})$ and define the complex-valued scalars $\psi_{ij}:=\pmb\psi(\Om_i, \Om_j)$, $i,j=1,2,3$. Then:
\begin{itemize}
\item The complex-valued scalars $\psi_{ij}$ satisfy 
\bea\lab{eq:fundamentalidentitiestoderivefromscalarizationoftensor:complexcase}
\psi_{ij}=\psi_{ji}, \qquad x^i\psi_{ij}=0, \qquad (\de^{ij}-x^ix^j)\psi_{ij}=0, \qquad \in_{ikl}x^l\psi_{kj}+i\psi_{ij}=0.
\eea
\item We may recover the tensor $\pmb\psi$ from the scalars $\psi_{ij}$ by the formula 
\beaa
\pmb\psi_{ab}=\psi_{ij}(\Om^i)_a(\Om^j)_b.
\eeaa
\end{itemize}
\item Reciprocally, let $\psi_{ij}$ be complex-valued scalars satisfying the identities \eqref{eq:fundamentalidentitiestoderivefromscalarizationoftensor:complexcase}, and introduce the complex-valued horizontal 2-tensor $\pmb\psi$ by $\pmb\psi_{ab}:=\psi_{ij}(\Om^i)_a(\Om^j)_b$, $a,b=1,2$. Then, we have $\pmb\psi\in\sk_2(\mathbb{C})$ and $\pmb\psi(\Om_i, \Om_j)=\psi_{ij}$ for all $i,j=1,2,3$.
\end{enumerate}
\end{lemma}

\begin{proof}
See Lemma 3.8 in \cite{MaSz26}.
\end{proof}

%%%%%%%%%%%%%%%%%%%%%%%%%%%%%%%%%%%%%%%%

\paragraph{\textit{Scalarization of the tensorial wave operator $\squared_2$}.}

%%%%%%%%%%%%%%%%%%%%%%%%%%%%%%%%%%%%%%%%

The following lemma provides the scalarization of the tensorial wave operator $\squared_2$.
\begin{lemma}\lab{lemma:formoffirstordertermsinscalarazationtensorialwaveeq}
Let $(\MM, \g)$ be a spacetime, $(e_3, e_4)$ be a null pair, and consider the corresponding horizontal structure $\O(\MM)$ introduced in Section \ref{subsection:review-horiz.structures}. Assume that $\Om_i$, $1,2,3$ is a regular triplet in the sense of Definition \ref{def:definitionofregulartripletOmii=123}. Also, let $\pmb\psi\in\sk_2(\mathbb{C})$ and let $\psi_{ij}$ be the scalars associated to it in view of Lemma \ref{lemma:backandforthbetweenhorizontaltensorsk2andscalars:complex}. Then, we have
\bea
\squared_2\pmb\psi(\Om_i, \Om_j) &=& \square_\g(\psi_{ij}) -S(\psi)_{ij} - (Q\psi)_{ij}
\eea
where 
\bsub
\label{SandV}
\begin{align}
S(\psi)_{ij} ={}& 2M_{i}^{k\a}\pr_\a(\psi_{kj}) +2M_{j}^{k\a}\pr_\a(\psi_{ik}),\\
(Q\psi)_{ij} ={}& (\Ddot^\a M_{i\a}^k)\psi_{kj}+(\Ddot^\a M_{j\a}^k)\psi_{ik} -M_{i\a}^kM_k^{l\a}\psi_{lj}-2M_{i\a}^kM_{j}^{l\a}\psi_{kl}-M_{j\a}^kM_k^{l\a}\psi_{il},
\end{align}
\esub
with the 1-forms $M_{i\a}^j$ defined by \eqref{eq:definitionofMalphaijwithoutambiguity}.
\end{lemma}

\begin{proof}
See Lemma 3.9 in \cite{MaSz26}.
\end{proof}

%%%%%%%%%%%%%%%%%%%%%%%%%%%%%%%%%%%%%%%%%%%%%%%%%%%%%%%%%%%%%%%%%%%%%%%%%%%%%%%%%%%%%%%%%%%%%%%%

\paragraph{\textit{Differentiation with respect to $\pr_\tau$ and $\widehat{\pr}_{\tphi}$ preserving identities \eqref{eq:fundamentalidentitiestoderivefromscalarizationoftensor:complexcase}}.}

%%%%%%%%%%%%%%%%%%%%%%%%%%%%%%%%%%%%%%%%%%%%%%%%%%%%%%%%%%%%%%%%%%%%%%%%%%%%%%%%%%%%%%%%%%%%%%%%

We start by noticing that differentiation w.r.t. $\pr_\tau$ preserves the identities  \eqref{eq:fundamentalidentitiestoderivefromscalarizationoftensor:complexcase}. 
\begin{lemma}\lab{lemma:differentiatingwrtprtaupreservetheidentitiesscaloftensors}
Let $\psi_{ij}$ be a family of complex-valued scalars satisfying the identities \eqref{eq:fundamentalidentitiestoderivefromscalarizationoftensor:complexcase}. Then, $\pr_\tau(\psi_{ij})$ satisfies the identities \eqref{eq:fundamentalidentitiestoderivefromscalarizationoftensor:complexcase} as well.
\end{lemma}

\begin{proof}
See Lemma 3.15 in \cite{MaSz26}.
\end{proof}

While $\pr_{\tphi}(x^3)=0$, we have $\pr_{\tphi}(x^1)=-x^2$ and $\pr_{\tphi}(x^2)=x^1$. Hence, differentiation w.r.t.  $\pr_{\tphi}$ does not preserve the identities \eqref{eq:fundamentalidentitiestoderivefromscalarizationoftensor:complexcase} and we will instead use the following modification. 

\begin{definition}\lab{def:widehatprtphi}
Let $\widehat{\pr}_{\tphi}$ denote the first-order operator acting on families of complex-valued scalars $\psi_{ij}$ as follows 
\beaa
\widehat{\pr}_{\tphi}(\psi)_{ij}:=\pr_{\tphi}(\psi_{ij}) +\in_{ik3}\psi_{kj} +\in_{jk3}\psi_{ki}.
\eeaa
\end{definition}

\begin{remark}
In Kerr, if $\psi_{ij}=\pmb\psi(\Om_i, \Om_j)$ with $\pmb\psi\in\sk_2$, then we have $\widehat{\pr}_{\tphi}(\psi)_{ij}=\Lieb_{\pr_{\tphi}}\pmb\psi(\Om_i, \Om_j)$, see Lemma \ref{lemma:InKerrlinkbetweenprtauwidehatprtphiandLieb[rtauLiebprtphi}, where  the horizontal Lie derivative $\Lieb$ has been introduced in Definition \ref{definition:hor-Lie-derivative}. This motivates Definition \ref{def:widehatprtphi}.
\end{remark}

The following lemma proves that differentiation w.r.t. $\widehat{\pr}_{\tphi}$ preserves the identities \eqref{eq:fundamentalidentitiestoderivefromscalarizationoftensor:complexcase}. 
\begin{lemma}\lab{lemma:differentiatingwrtwidehatprtphipreservetheidentitiesscaloftensors}
Let $\psi_{ij}$ be a family of complex-valued scalars satisfying the identities \eqref{eq:fundamentalidentitiestoderivefromscalarizationoftensor:complexcase} and let $\widehat{\pr}_{\tphi}$ be as in Definition \ref{def:widehatprtphi}. Then, $\widehat{\pr}_{\tphi}(\psi)_{ij}$ satisfies the identities \eqref{eq:fundamentalidentitiestoderivefromscalarizationoftensor:complexcase} as well.
\end{lemma}

\begin{proof}
See Lemma 3.18 in \cite{MaSz26}.
\end{proof}
 
\begin{remark}
\lab{rem:prtauandprtphihatpreservesk2C}
In view of Lemmas \ref{lemma:backandforthbetweenhorizontaltensorsk2andscalars:complex},  \ref{lemma:differentiatingwrtprtaupreservetheidentitiesscaloftensors} and   \ref{lemma:differentiatingwrtwidehatprtphipreservetheidentitiesscaloftensors}, we immediately infer the fact that 
if $\psi_{ij}=\pmb\psi(\Om_i, \Om_j)$ for $\pmb\psi\in \sk_2(\mathbb{C})$, then for any $k,l\in\mathbb{N}$, there exists $\pmb\psi_{(k,l)}\in\sk_2(\mathbb{C})$ such that $\pr_\tau^k\widehat{\pr}_{\tphi}^l(\psi)_{ij}=\pmb\psi_{(k,l)}(\Om_i,\Om_j)$.
\end{remark}

%%%%%%%%%%%%%%%%%%%%%%%%%%%%%

\paragraph{\textit{Regular triplet in Kerr}.}

%%%%%%%%%%%%%%%%%%%%%%%%%%%%%

\begin{definition}[Regular triplet in Kerr]
\lab{def:regulartripletinKerrOmii=123}
Let 
\bea\lab{eq:definitionofxii=123usedintheformulaoftheregulartripletsinKerr}
x^1:=\cos\tphi\sin\th, \qquad x^2:=\sin\tphi\sin\th, \qquad x^3:=\cos\th.
\eea
Then, we define the following horizontal vectorfields $\Om^i$ in Kerr by 
\bea
\Om^i :=|q|\dual\nab(x^i), \quad i=1,2,3.
\eea
\end{definition}

\begin{lemma}\lab{lemma:fundamentalpropertiesof1formsOmi}
The horizontal vectorfields $\Om_i$ in Kerr introduced in Definition \ref{def:regulartripletinKerrOmii=123} satisfy \eqref{eq:fundamentalpropertiesof1formsOmi}. In particular, they form a regular triplet in Kerr in the sense of Definition \ref{def:definitionofregulartripletOmii=123}.
\end{lemma}

\begin{proof}
See Lemma 3.21 in \cite{MaSz26}.
\end{proof}

We also derive the following properties of the 1-forms $M_{i\a}^j$ in Kerr.
\begin{lemma}\lab{lemma:computationoftheMialphajinKerr}
Let $\Om_i$, $i=1,2,3$, be the regular triplet in Kerr of Definition \ref{def:regulartripletinKerrOmii=123}, and let $M_{i\a}^j$ be the corresponding 1-forms in Kerr given by \eqref{eq:definitionofMalphaijwithoutambiguity}. Then, we have, for $i,j=1,2,3$,  
\beaa
M_{i3}^j &=& \phimod'(r)\in_{ki3}(\de^{kj}-x^kx^j)+\frac{a\cos\th}{|q|^2}\in_{ijk}x^k,\\
M_{i4}^j &=& \frac{2a -\De\phimod'(r)}{|q|^2}\in_{ki3}(\de^{kj}-x^kx^j)+\frac{a\cos\th\De}{|q|^4}\in_{ijk}x^k,\\
M_{i\a}^j(\pr_\tau)^\a &=& -\frac{2amr\cos\th}{|q|^4}\in_{ijk}x^k,
\eeaa
and the following asymptotic holds, for $r$ large and $i,j=1,2,3$,
\beaa
M_{ia}^j=O(r^{-1}), \quad a=1,2.
\eeaa
Also, we have for $r\in [r_+(1+2\dbl),  12m]$
\beaa
M_{i\a}^j(\pr_r)^\a=0,\qquad \g^{r\a}M_{i\a}^j=0, \quad i,j=1,2,3.
\eeaa
\end{lemma}

\begin{proof}
See Lemma 3.22 in \cite{MaSz26}.
\end{proof}

\begin{lemma}\lab{lemma:InKerrlinkbetweenprtauwidehatprtphiandLieb[rtauLiebprtphi}
In Kerr, if $\psi_{ij}=\pmb\psi(\Om_i, \Om_j)$ with $\pmb\psi\in\sk_2$, then 
\beaa
\Lieb_{\pr_\tau}\Om_i=0, \qquad \Lieb_{\pr_{\tphi}}\Om_i = -\in_{ij3}\Om_j,
\eeaa
and
\beaa
\pr_\tau(\psi_{ij})=\Lieb_{\pr_\tau}\pmb\psi(\Om_i, \Om_j), \qquad \widehat{\pr}_{\tphi}(\psi)_{ij}=\Lieb_{\pr_{\tphi}}\pmb\psi(\Om_i, \Om_j),
\eeaa
where the horizontal Lie derivative $\Lieb$ has been introduced in Definition \ref{definition:hor-Lie-derivative}, and where $\widehat{\pr}_{\tphi}$ has been introduced in Definition \ref{def:widehatprtphi}.
\end{lemma}

\begin{proof}
See Lemma 3.23 in \cite{MaSz26}.
\end{proof}

Finally, we provide the explicit formula of $M_{ia}^j$, $i,j=1,2,3$, $a=1,2$, in Kerr.
\begin{lemma}\lab{lemma:explicitcomputationofMiaja=1or2}
Let $\Om_i$, $i=1,2,3$, be the regular triplet in Kerr of Definition \ref{def:regulartripletinKerrOmii=123}, and let $M_{i\a}^j$ be the corresponding 1-forms in Kerr given by \eqref{eq:definitionofMalphaijwithoutambiguity}. Then, we have, for $a=1,2$, and $j=1,2,3$,  
\beaa
M_{1a}^j &=& -\frac{a^2x^3}{|q|^3}(\de_{1j}-x^1x^j)(\dual\Om_3)_a +\frac{r^2}{|q|^3}x^1(\dual\Om_j)_{a}-\frac{2a^2x^3}{|q|^3}x^2(\Om_j)_{a},\\
M_{2a}^j &=& -\frac{a^2x^3}{|q|^3}(\de_{2j}-x^2x^j)(\dual\Om_3)_a +\frac{r^2}{|q|^3}x^2(\dual\Om_j)_{a}+\frac{2a^2x^3}{|q|^3}x^1(\Om_j)_{a},\\
M_{3a}^j &=& -\frac{a^2x^3}{|q|^3}(\de_{3j}-x^3x^j)(\dual\Om_3)_a +\frac{r^2+a^2}{|q|^3}x^3(\dual\Om_j)_{a}.
\eeaa
\end{lemma}

\begin{proof}
See Lemma 5.21 in \cite{Sze}.
\end{proof}

%%%%%%%%%%%%%%%%%%%%%%%%%%%%%%%%%%%%%%%%%%%%%%%%%%

\subsubsection{Definition and main properties of a regular triplet $\Om_i$, $i=1,2,3$, in $\MM$}
\lab{sec:defandmainpropofregtripletOmiinKerrpert}

%%%%%%%%%%%%%%%%%%%%%%%%%%%%%%%%%%%%%%%%%%%%%%%%%%

In the following lemma, we exhibit a regular triplet $\Om_i$, $i=1,2,3$, in $\MM$.
\begin{lemma}\lab{lemma:constructionofhorizontal1formsOmiinperturbationsofKerr}
On $\MM$, there exist horizontal 1-forms $\Om_i$, $i=1,2,3$, such that 
\begin{enumerate}
\item The horizontal 1-forms $\Om_i$, $i=1,2,3$ satisfy on $\MM$
\bea\lab{eq:fundamentalpropertiesof1formsOmi:txversiononMM}
\tx^i\Om_i=0, \qquad (\Om^i)_a(\Om_i)_b=\de_{ab}, \qquad \Om_i\c\Om_j=\de_{ij}-\tx^i\tx^j, \qquad \Om_i\c\dual\Om_j=\in_{ijk}\tx^k,
\eea
with  $(\tx^1, \tx^2, \tx^3)=(\tx^1_p, \tx^2_p, \cos(\widetilde{\th}))$, where $(\widetilde{\th}, \tx^1, \tx^2)$ are defined in  \eqref{definitionofthenewcoordinatetildethetandtphi}. In particular, $\Om_i$, $i=1,2,3$ forms a regular triplet on $\MM$.

\item We have 
\bea\lab{eq:estimatefor15weightedderivativesofdoublewidechechnabalphaOmionMM}
\bsplit
|\dk^{\leq 15}\widecheck{\widecheck{M_{i4}^j}}|&\les \ep\min\Big\{r^{-2}\tau^{-\frac{1+\dec}{2}}, \, r^{-1}\tau^{-1-\frac{3\dec}{4}}\Big\},\quad i,j=1,2,3,\quad\textrm{on}\quad\MM,\\
|\dk^{\leq 15}\big(\widecheck{\widecheck{M_{i3}^j}}, \widecheck{\widecheck{M_{ia}^j}}\big)|&\les \ep r^{-1}\tau^{-1-\frac{3\dec}{4}},\quad i,j=1,2,3,\quad a=1,2,\quad\textrm{on}\quad\MM,
\end{split}
\eea
where the notation $M_{i\a}^j$, $\a=1,2,3,4$, has been introduced in \eqref{eq:definitionofMalphaijwithoutambiguity}, where $\widecheck{\widecheck{M}}_{i\a}^j$, $\a=3,4$, are given by
 \beaa
 \widecheck{\widecheck{M_{i3}^j}} &:=& M_{i3}^j - \left(\phimod'(\rr)\in_{ki3}(\de^{kj}-\tx^k\tx^j)+\frac{a\cos(\widetilde{\th})}{|\tq|^2}\in_{ijk}\tx^k\right),\\
 \widecheck{\widecheck{M_{i4}^j}} &:=& M_{i4}^j - \left(\frac{2a -\De(\rr)\phimod'(\rr)}{|\tq|^2}\in_{ki3}(\de^{kj}-\tx^k\tx^j)+\frac{a\cos(\widetilde{\th})\De(\rr)}{|\tq|^4}\in_{ijk}\tx^k\right),
 \eeaa
and where $\widecheck{\widecheck{M}}_{ia}^j$, $a=1,2$ are given by 
\beaa
\widecheck{\widecheck{M_{1a}^j}} &:=& M_{1a}^j +\frac{a^2\tx^3}{|\tq|^3}(\de_{1j}-\tx^1\tx^j)(\dual\Om_3)_a -\frac{\rr^2}{|\tq|^3}\tx^1(\dual\Om_j)_{a}+\frac{2a^2\tx^3}{|\tq|^3}\tx^2(\Om_j)_{a},\\
\widecheck{\widecheck{M}}_{2a}^j &:=& M_{2a}^j +\frac{a^2\tx^3}{|\tq|^3}(\de_{2j}-\tx^2\tx^j)(\dual\Om_3)_a -\frac{\rr^2}{|\tq|^3}\tx^2(\dual\Om_j)_{a}-\frac{2a^2\tx^3}{|\tq|^3}\tx^1(\Om_j)_{a},\\
\widecheck{\widecheck{M}}_{3a}^j &:=& M_{2a}^j +\frac{a^2\tx^3}{|\tq|^3}(\de_{3j}-\tx^3\tx^j)(\dual\Om_3)_a -\frac{\rr^2+a^2}{|\tq|^3}\tx^3(\dual\Om_j)_{a},
\eeaa
with $(\rr, \tx^1, \tx^2)$ the coordinates defined in \eqref{definitionofthenewcoordinatetilder} \eqref{definitionofthenewcoordinatetildethetandtphi}, and $\tJk$ the complex 1-form defined in \eqref{definitionofthenewcomplex1formstJktJkpm}.

\item We have
\bea\lab{eq:estimatefor15weightedderivativesofrenormalizedLiebtphiOmi}
|\dk^{\leq 15}\big(\Lieb_{\pr_{\widetilde{\tphi}}}\Om_i +\in_{ij3}\Om^j\big)|&\les \ep\tau^{-1-\frac{3\dec}{4}},\quad i=1,2,3,\quad a=1,2,\quad\textrm{on}\quad\MM.
\eea
\end{enumerate}
\end{lemma}

\begin{remark}
In view of \eqref{eq:estimatefor15weightedderivativesofdoublewidechechnabalphaOmionMM}, \eqref{eq:estimatefor15weightedderivativesofrenormalizedLiebtphiOmi}, and the estimates \eqref{eq:decaypropertiesofwidehatGab1widehatGag1} satisfied by $(\widetilde{\Ga}_{g}, \widetilde{\Ga}_{b})$, we have
\bsub\lab{eq:assumptionsonregulartripletinperturbationsofKerr}
\bea
&&\lab{eq:assumptionsonregulartripletinperturbationsofKerr:0}
\widecheck{\widecheck{M_{i4}^j}}=\widetilde{\Ga}_g, \qquad \widecheck{\widecheck{M_{i3}^j}}=\widetilde{\Ga}_b,  \qquad \widecheck{\widecheck{M_{ia}^j}}=\widetilde{\Ga}_b, \quad \forall\, i,j,a,\\
&&\lab{eq:assumptionforLiebprtphiOmiinKerrperturbation}
\Lieb_{\pr_{\widetilde{\tphi}}}\Om_i +\in_{ij3}\Om^j=r\widetilde{\Ga}_b, \quad\textrm{for}\quad i=1,2,3.
\eea
\esub 
\end{remark}

\begin{proof}
The proof proceeds in the following step. 

\noindent {\bf Step 1.} First, we define 1-forms $(\Om_*)_i$, $i=1,2,3$, on $\Si_*$ by 
\bea\lab{eq:definition1formsOmstaronSigmastar}
(\Om_*)_1=\dual f_+, \qquad (\Om_*)_2=\dual f_-, \qquad (\Om_*)_3=f_0,\quad\textrm{on}\quad\Si_*,
\eea
where the 1-forms $f_0$, $f_+$ and $f_-$ on $\Si_*$ have been introduced in Definition \ref{def:definitionoff0fplusfminus}. In view of \eqref{eq:definition1formsOmstaronSigmastar} and in view of the transport equations \eqref{eq:extendionofthereal1formsf0fpmfromS*toSigma*} for $f_0$, $f_+$ and $f_-$, we have 
\bea\lab{eq:transporteqationinnustaralongSigmastarforOmii=123}
\nab_{\nu_*}(\Om_*)_i =0, \quad i=1,2,3,\quad\textrm{on}\quad\Si_*.
\eea
Also, in view of \eqref{eq:definition1formsOmstaronSigmastar}, as well as the definition of  $f_0$, $f_+$ and $f_-$ on $S_*$, see \eqref{eq:definitionofthereal1formsf0fpmonS*}, we have in particular the following identities on $S_*$
\bea\lab{eq:fundamentalpropertiesof1formsOmi*:onS*}
\bsplit
& x^i(\Om_*)_i=0, \quad (\Om^i_*)_a((\Om_*)_i)_b=\de_{ab}, \qquad (\Om_*)_i\c(\Om_*)_j=\de_{ij}-x^ix^j, \\
& (\Om_*)_i\c\dual(\Om_*)_j=\in_{ijk}x^k,\quad\textrm{on}\quad S_*, \qquad (x^1, x^2, x^3):=(x^1_p, x^2_p, \cos(\th)).
 \end{split}
\eea
Then, using \eqref{eq:transporteqationinnustaralongSigmastarforOmii=123} together with the fact that $\nu_*(x^i)=0$, $i=1,2,3$, along $\Si_*$ in view of \eqref{eq:corrdinatesthandvarphiarepropagatedalongSi*bynu*}, we have 
\beaa
&&\nab_{\nu_*}\big(x^i(\Om_*)_i\big)=0, \qquad \nab_{\nu_*}\big((\Om^i_*)_a((\Om_*)_i)_b-\de_{ab}\big)=0, \qquad \nu_*\big((\Om_*)_i\c(\Om_*)_j-\de_{ij}+x^ix^j\big)=0,\\
&& \nu_*\big((\Om_*)_i\c\dual(\Om_*)_j-\in_{ijk}x^k\big)=0, \quad\textrm{on}\quad\Si_*,
\eeaa
which allows to extend \eqref{eq:fundamentalpropertiesof1formsOmi*:onS*} from $S_*$ to $\Si_*$, i.e., 
\bea\lab{eq:fundamentalpropertiesof1formsOmi*:onSigma*}
\bsplit
& x^i(\Om_*)_i=0, \quad (\Om^i_*)_a((\Om_*)_i)_b=\de_{ab}, \qquad (\Om_*)_i\c(\Om_*)_j=\de_{ij}-x^ix^j, \\
& (\Om_*)_i\c\dual(\Om_*)_j=\in_{ijk}x^k,\quad\textrm{on}\quad \Si_*.
 \end{split}
\eea
so that $(\Om_*)_i$, $i=1,2,3$ satisfy \eqref{eq:fundamentalpropertiesof1formsOmi} on $\Si_*$.

\noindent {\bf Step 2.} Next, we estimate angular derivatives of $(\Om_*)_i$, $i=1,2,3$ on $\Si_*$. We introduce the following linearized quantities 
\bea\lab{eq:defofwidechecknabstarOmstarii=123onSigma*}
\widecheck{\nab_*(\Om_*)_i} := \nab_*(\Om_*)_i - \frac{1}{r}x^i\in, \quad i=1,2,3,
\eea
which together with \eqref{eq:definition1formsOmstaronSigmastar} and \eqref{eq:assumptionsonSigmastarforpartII} implies, using also \eqref{eq:dominantconditionforrstarcomparedtotaustartonS*},
\bea\lab{eq:controlofwidechecknabstarOmstarii=123onSigma*}
|\dk_*^{\leq 15}\widecheck{\nab_*(\Om_*)_i}|\les \frac{\ep}{r^2\tau^{\frac{1}{2}+\dec}}\les \frac{\ep}{r\tau^{1+\dec}}\quad\textrm{on}\quad\Si_*.
\eea

Next, we control $\Lieb_{\pr_{\tphi}}(\Om_*)_i$ on $S_*$. First, in view of \eqref{eq:specialorthonormalbasisofSstar}, we have
\beaa
\Lieb_{\pr_{\tphi}}(e_*)_1 &=& \frac{1}{re^{\phi_*}}[\pr_{\tphi}, \pr_\th] -\pr_{\tphi}(\phi_*)(e_*)_1=-\pr_{\tphi}(\phi_*)(e_*)_1,\quad\textrm{on}\quad S_*,\\
\Lieb_{\pr_{\tphi}}(e_*)_2 &=& \frac{1}{r\sin\th e^{\phi_*}}[\pr_{\tphi}, \pr_{\tphi}] -\pr_{\tphi}(\phi_*)(e_*)_2=-\pr_{\tphi}(\phi_*)(e_*)_2,\quad\textrm{on}\quad S_*,
\eeaa
which together with \eqref{eq:definitionofthereal1formsf0fpmonS*} yields 
\beaa
\Lieb_{\pr_{\tphi}}f_0 = -\pr_{\tphi}(\phi_*)f_0, \qquad \Lieb_{\pr_{\tphi}}f_+ =-f_- -\pr_{\tphi}(\phi_*)f_+, \qquad \Lieb_{\pr_{\tphi}}f_- =f_+ -\pr_{\tphi}(\phi_*)f_-, \quad\textrm{on}\quad S_*.
\eeaa
In view of \eqref{eq:definition1formsOmstaronSigmastar}, we infer
\beaa
\Lieb_{\pr_{\tphi}}(\Om_*)_i +\in_{ij3}(\Om_*)^j &=& -\pr_{\tphi}(\phi_*)(\Om_*)_i, \quad i=1,2,3,\quad\textrm{on}\quad S_*,
\eeaa
and hence, using \eqref{eq:controloftheconformaluniformizationfactorofinducedmetricS*}, as well as \eqref{eq:dominantconditionforrstarcomparedtotaustartonS*}, we obtain
\bea\lab{eq:controlof15weightedderivativesrenormalizedLiebtphiOmstarionS*}
\big|\dkb_*^{\leq 15}\big(\Lieb_{\pr_{\tphi}}(\Om_*)_i +\in_{ij3}(\Om_*)^j\big)\big|\les\frac{\ep}{r\tau^{\frac{1}{2}+\dec}}\les\frac{\ep}{\tau^{1+\dec}}, \quad i=1,2,3,\quad\textrm{on}\quad S_*.
\eea

\noindent {\bf Step 3.} Next, we extend \eqref{eq:controlof15weightedderivativesrenormalizedLiebtphiOmstarionS*} to $\Si_*$. To this end, we introduce the following approximate Killing vectorfield
\bea\lab{eq:definitionapproxKiilingvectorfieldZ*onSi*}
\Z_* := rf_0^b(e_*)_b \quad\textrm{on}\quad\Si_*.
\eea
In particular, we have $\Z_*(r)=\Z_*(\tau)=0$ on $\Si_*$ in view of \eqref{eq:identitiesonSigma*fornullframeadaptedSigma*}, and since 
\beaa
f_0\c f_+=-\sin\th\sin\tphi, \qquad f_0\c f_- =\sin\th\cos\tphi, \quad\textrm{on}\quad\Si_*,
\eeaa
in view of \eqref{eq:scalarproductreal1formsf0fpmonSi*}, together with Definition \ref{def:renormalizationforf0fpfm}, we also have 
\beaa
\Z_*(x^i)=\pr_{\tphi}(x^i), \quad i=1,2,3,\quad (x^1, x^2, x^3)=(\cos\th, x^1_p, x^2_p),
\eeaa
so that
\bea\lab{eq:infactZ*=prtphionSi*}
\Z_* = \pr_{\tphi}\quad\textrm{on}\quad\Si_*.
\eea
Also, we have from \eqref{eq:definitionapproxKiilingvectorfieldZ*onSi*}
\beaa
\g(\D_{(e_*)_a}\Z_*, (e_*)_b) &=& r((\nab_*)_af_0)_b=r\left(\nab_*f_0-\frac{1}{r}\cos\th\in\right)_{ab}+\cos\th\in_{ab}
\eeaa
so that 
\beaa
\nab_{\Z_*}(\Om_*)_i &=& \Lieb_{\Z_*}(\Om_*)_i -\cos\th\dual(\Om_*)_i -r\left(\nab_*f_0-\frac{1}{r}\cos\th\in\right)\c(\Om_*)_i.
\eeaa
Together with \eqref{eq:assumptionsonSigmastarforpartII} and \eqref{eq:dominantconditionforrstarcomparedtotaustartonS*}, we infer
\bea\lab{eq:comparisionrenomralizationofLiebtphiOmiandnabZ*OmionSi*}
\nn&&\left|\dk^{\leq 15}\left[\Lieb_{\pr_{\tphi}}(\Om_*)_i +\in_{ij3}(\Om_*)^j -\left(\nab_{\Z_*}(\Om_*)_i +\in_{ij3}(\Om_*)^j +\cos\th\dual(\Om_*)_i\right)\right]\right|\\ 
&\les& \frac{\ep}{r\tau^{\frac{1}{2}+\dec}}\les\frac{\ep}{\tau^{1+\dec}} \quad\textrm{on}\quad\Si_*.
\eea
In particular, \eqref{eq:controlof15weightedderivativesrenormalizedLiebtphiOmstarionS*} and \eqref{eq:comparisionrenomralizationofLiebtphiOmiandnabZ*OmionSi*} imply 
\bea\lab{eq:controlof15weightedderivativesrenormalizednabtphiOmstarionS*}
\big|\dkb_*^{\leq 15}\big(\nab_{\Z_*}(\Om_*)_i +\in_{ij3}(\Om_*)^j +\cos\th\dual(\Om_*)_i\big)\big|\les\frac{\ep}{\tau^{1+\dec}}, \quad i=1,2,3,\quad\textrm{on}\quad S_*.
\eea

Next, we extend \eqref{eq:controlof15weightedderivativesrenormalizednabtphiOmstarionS*} to $\Si_*$. To this end, we first fix the following transversality conditions 
\bea\lab{eq:transversailityconditionnab*e*4ofOmi=0onSi*}
\nab_{(e_*)_4}(\Om_*)_i &=& 0, \quad i=1,2,3, \quad\textrm{on}\quad\Si_*.
\eea
Then, in view of  \eqref{eq:corrdinatesthandvarphiarepropagatedalongSi*bynu*}, \eqref{eq:extendionofthereal1formsf0fpmfromS*toSigma*}, \eqref{eq:transporteqationinnustaralongSigmastarforOmii=123}, \eqref{eq:definitionapproxKiilingvectorfieldZ*onSi*} and \eqref{eq:transversailityconditionnab*e*4ofOmi=0onSi*}, we have 
\beaa
&&\nab_{\nu_*}\big(\nab_{\Z_*}(\Om_*)_i +\in_{ij3}(\Om_*)^j +\cos\th\dual(\Om_*)_i\big)\\
&=& \nab_{\nu_*}\nab_{\Z_*}(\Om_*)_i=\nab_{\nu_*}(f_0^br\nab_{(e_*)_b}(\Om_*)_i)= f_0^b\nab_{\nu_*}(r\nab_{(e_*)_b}(\Om_*)_i)\\
&=& f_0^b[\nab_{\nu_*},r\nab_{(e_*)_b}](\Om_*)_i\\
&=& f_0^b[\nab_{(e_3)_*},r\nab_{(e_*)_b}](\Om_*)_i+b_*f_0^b[\nab_{(e_4)_*},r\nab_{(e_*)_b}](\Om_*)_i -rf_0\c\nab_*(b_*)\nab_{(e_4)_*}(\Om_*)_i\\
&=& f_0^b[\nab_{(e_3)_*},r\nab_{(e_*)_b}](\Om_*)_i+b_*f_0^b[\nab_{(e_4)_*},r\nab_{(e_*)_b}](\Om_*)_i.
\eeaa
In view of the commutation formulas \eqref{commutator-3-a-u-b} and \eqref{commutator-4-a-u-b}, we infer
\beaa
\nab_{\nu_*}\big(\nab_{\Z_*}(\Om_*)_i +\in_{ij3}(\Om_*)^j +\cos\th\dual(\Om_*)_i\big) &=& r\Ga_b^*(\nab_{(e_*)_4}, \nab_{(e_*)_3}, \nab_*)\Om_i+\Ga_b^*,
\eeaa
which together with \eqref{eq:transporteqationinnustaralongSigmastarforOmii=123} and \eqref{eq:transversailityconditionnab*e*4ofOmi=0onSi*} yields
\beaa
\nab_{\nu_*}\big(\nab_{\Z_*}(\Om_*)_i +\in_{ij3}(\Om_*)^j +\cos\th\dual(\Om_*)_i\big) &=& \Ga_b^*\c r\widecheck{\nab_*\Om_i}+\Ga_b^*.
\eeaa
Commuting with the general commutation formulas of Lemma \ref{LEMMA:COMM-GEN-B}, using the values of $\B_{ab\mu\nu}$ given by Proposition \ref{proposition:componentsofB}, using also \eqref{eq:transporteqationinnustaralongSigmastarforOmii=123} and \eqref{eq:transversailityconditionnab*e*4ofOmi=0onSi*}, and then integrating the resulting transport equation in $\nu_*$ from $S_*$ where we have \eqref{eq:controlof15weightedderivativesrenormalizednabtphiOmstarionS*}, we infer
\beaa
\big|\dk_*^{\leq 15}\big(\nab_{\Z_*}(\Om_*)_i +\in_{ij3}(\Om_*)^j +\cos\th\dual(\Om_*)_i\big)\big|\les\frac{\ep}{\tau^{1+\dec}}+\int_{\tau}^{\tau_*}\frac{\ep}{r{\tau'}^{1+\dec}}d\tau'\les \frac{\ep}{\tau^{1+\dec}}+\frac{\ep}{r\tau^{\dec}},
\eeaa
which together with \eqref{eq:dominantconditionforrstarcomparedtotaustartonS*} yields
\bea\lab{eq:controlof15weightedderivativesrenormalizedLiebtZstarOmstarionSigma*}
\big|\dk_*^{\leq 15}\big(\nab_{\Z_*}(\Om_*)_i +\in_{ij3}(\Om_*)^j +\cos\th\dual(\Om_*)_i\big)\big|\les\frac{\ep}{\tau^{1+\dec}}, \quad i=1,2,3,\quad\textrm{on}\quad \Si_*.
\eea

\noindent {\bf Step 4.} Next, we initialize the horizontal 1-forms $\Om_i$, $i=1,2,3$ on $\Si_*$ as follows 
\bea\lab{eq:initializationofOmii=123onSigma*fromOm*i}
\langle \Om_i, e_a\rangle=\langle (\Om_*)_i), (e_*)_a\rangle, \qquad i=1,2,3, \quad a=1,2,\quad \textrm{on}\quad\Si_*.
\eea
We then immediately infer from \eqref{eq:fundamentalpropertiesof1formsOmi*:onSigma*} and \eqref{definitionofthenewcoordinatetildethetandtphi}
\bea\lab{eq:fundamentalpropertiesof1formsOmi:onSigma*}
\tx^i\Om_i=0, \quad (\Om^i)_a(\Om_i)_b=\de_{ab}, \quad \Om_i\c\Om_j=\de_{ij}-\tx^i\tx^j, \quad \Om_i\c\dual\Om_j=\in_{ijk}\tx^k,\quad\textrm{on}\quad\Si_*,
\eea
with $(\tx^1, \tx^2, \tx^3):=(\tx^1_p, \tx^2_p, \cos(\widetilde{\th}))$ so that $\Om_i$, $i=1,2,3$ satisfy \eqref{eq:fundamentalpropertiesof1formsOmi:txversiononMM} on $\Si_*$.

Next, we extend the horizontal 1-forms $\Om_i$, $i=1,2,3$ to $\MM$ using the following transport equations 
\bea\lab{eq:trnasrportequationinYtauforextensionOmii=123tpMM}
\bsplit
\nab_{\Ytau}'\Om_i =& [\Ytau^3]_K(\rr)\left(\in^j\!\!_{i3}\phimod'(\rr)\Om_j -\frac{a\cos(\widetilde{\th})}{|\tq|^2}\dual\Om_i\right)\\
&+[\Ytau^4|q|^{-2}]_K(\rr)\left(\in^j\!\!_{i3}(2a -\De(\rr)\phimod'(\rr))\Om_j -\frac{a \De(\rr)\cos(\widetilde{\th})}{|\tq|^2}\dual\Om_i\right)\quad\textrm{on}\quad\MM,
\end{split}
\eea
where the vectorfield $\Ytau$ has beed introduced in Definition \ref{def:definitionofthevectorfildYtau}, and where $[\Ytau^4|q|^{-2}]_K(\rr)$ and $[\Ytau^3]_K(\rr)$ are defined in \eqref{eq:defintionofKerrvaluesofYtau4modqminus2andYtau3}. 

Using \eqref{eq:trnasrportequationinYtauforextensionOmii=123tpMM}, as well as the transport equations along $\Ytau$ in \eqref{definitionofthenewcoordinatetildethetandtphi} for $(\tx^1, \tx^2, \tx^3)=(\tx^1_p, \tx^2_p, \cos(\widetilde{\th}))$, we obtain 
\beaa
\nab_{\Ytau}'(\tx^i\Om_i) &=& -\Big([\Ytau^3]_K(\rr)+\De[\Ytau^4|q|^{-2}]_K(\rr)\Big)\frac{a\cos(\widetilde{\th})}{|\tq|^2}\dual(\tx^i\Om_i),\\
\nab_{\Ytau}'((\Om^i)_a(\Om_i)_b-\de_{ab}) &=&  -\Big([\Ytau^3]_K(\rr)+\De[\Ytau^4|q|^{-2}]_K(\rr)\Big)\frac{a\cos(\widetilde{\th})}{|\tq|^2}\\
&&\qquad\qquad\times\Big(\in_{ac}\big((\Om_i)_c(\Om^i)_b-\de_{bc}\big)+\in_{bc}\big((\Om_i)_c(\Om^i)_a-\de_
{ac}\big)\Big),\\ 
\Ytau\Big(\Om_i\c\Om_j-\de_{ij}+\tx^i\tx^j\Big) &=& \Big([\Ytau^4|q|^{-2}]_K(\rr)(2a -\De(\rr)\phimod'(\rr))+[\Ytau^3]_K(\rr)\phimod'(\rr)\Big)\\
&&\times\Big(\in_{ki3}\big(\Om_k\c\Om_j-\de_{kj}+\tx^k\tx^j\big)+\in_{kj3}\big(\Om_k\c\Om_i-\de_{ki}+\tx^k\tx^i\big)\Big),\\
\Ytau\big(\Om_i\c\dual\Om_j-\in_{ijk}\tx^k\big) &=& \Big([\Ytau^4|q|^{-2}]_K(\rr)(2a -\De\phimod'(\rr))+[\Ytau^3]_K(\rr)\phimod'(\rr)\Big)\\
&&\qquad\times\Big(\in^k\!\!_{i3}\big(\Om_k\c\dual\Om_j -\in_{kjl}\tx^l\big) -\in^k\!\!_{j3}\big(\Om_k\c\dual\Om_i -\in_{kil}\tx^l\big)\Big).
\eeaa
Together with \eqref{eq:fundamentalpropertiesof1formsOmi:onSigma*}, this immediately extends \eqref{eq:fundamentalpropertiesof1formsOmi:onSigma*} from $\Si_*$ to $\MM$ so that 
so that $\Om_i$, $i=1,2,3$ satisfy \eqref{eq:fundamentalpropertiesof1formsOmi:txversiononMM} on $\MM$ with $(\tx^1, \tx^2, \tx^3)=(\tx^1_p, \tx^2_p, \cos(\widetilde{\th}))$ as stated.

\noindent {\bf Step 5.} It remains to prove \eqref{eq:estimatefor15weightedderivativesofdoublewidechechnabalphaOmionMM} \eqref{eq:estimatefor15weightedderivativesofrenormalizedLiebtphiOmi}. First, we define the renormalized quantities of derivatives of $\Om_i$, $i=1,2,3$. To this end, notice from \eqref{eq:fundamentalpropertiesof1formsOmi:txversiononMM} and \eqref{eq:definitionofMalphaijwithoutambiguity} that 
\beaa
(\Ddot_\a\Om_i)_a &=& (\Ddot_\a\Om_i)_b\de_{ab}=(\Ddot_\a\Om_i)_b(\Om_j)_a(\Om^j)_b= (\Ddot_\a\Om_i\c\Om^j)(\Om_j)_a\\
&=& M_{i\a}^j(\Om_j)_a
\eeaa
and hence
\beaa
\nab_4\Om_i=M_{i4}^j\Om_j, \qquad \nab_3\Om_i=M_{i3}^j\Om_j, \qquad \nab_a\Om_i=M_{ia}^j\Om_j, \quad i=1,2,3,\quad a=1,2.
\eeaa
We thus define the renormalized quantities of derivatives of $\Om_i$, $i=1,2,3$, as follows 
\bea\lab{eq:definitionlinearizedquantiiesnabaphaOmiwithdoublecheck}
\widecheck{\widecheck{\nab_4\Om_i}}=\widecheck{\widecheck{M_{i4}^j}}\Om_j, \qquad \widecheck{\widecheck{\nab_3\Om_i}}=\widecheck{\widecheck{M_{i3}^j}}\Om_j, \qquad \widecheck{\widecheck{\nab_a\Om_i}}=\widecheck{\widecheck{M_{ia}^j}}\Om_j, \quad i=1,2,3,\quad a=1,2.
\eea

Next, we consider the control of $\Om_i$, $i=1,2,3$, on $\Si_*$ in the global null frame $(e_3', e_4', e_1', e_2')$. In view of \eqref{eq:definitionlinearizedquantiiesnabaphaOmiwithdoublecheck} and the definition of $\widecheck{\widecheck{M_{i\a}^j}}$, we have, using also the fact that $\rr=r$ and $(\tx^1, \tx^2, \tx^3)=(x^1, x^2, x^3)$ on $\Si_*$,  
\bea\lab{eq:identityusingasymptoticofKerrvaluefornabalphaOmitocontrollinearizednabalphaOmionSi*}
\bsplit
\nab_3'\Om_i &= O(mr^{-2})+\widecheck{\widecheck{\nab_3'\Om_i}}, \qquad \nab_4'\Om_i =O(mr^{-2})+\widecheck{\widecheck{\nab_4'\Om_i}},\\ 
(\nab_a'\Om_i)_b &=\frac{x^i}{r}\in_{ab}+O(m^2r^{-3})+(\widecheck{\widecheck{\nab_a'\Om_i}})_b\quad\textrm{on}\quad\Si_*.
\end{split}
\eea
Also, we have in view of Remark \ref{rmk:formofYautforrgeqr0intermsofdoubleprimedframe}, \eqref{eq:defintionofKerrvaluesofYtau4modqminus2andYtau3} and \eqref{eq:trnasrportequationinYtauforextensionOmii=123tpMM}
\bea\lab{eq:identityusingasymptoticofKerrvaluefornabalphaOmitocontrollinearizednabalphaOmionSi*:bis}
\nab_4'\Om_i &=& O(mr^{-2})+\big(O(m^2r^{-2})+\Ga_g'\big)\nab_3'\Om_i. 
\eea
Next, using the change of frame formulas \eqref{eq:changeofframecoefffromfromframeSigmastarttodoubleprimedframe} \eqref{eq:changeofframefromfromframeSigmastarttodoubleprimedframe}, the initialization \eqref{eq:initializationofOmii=123onSigma*fromOm*i} on $\Si_*$, and the transport equation \eqref{eq:transporteqationinnustaralongSigmastarforOmii=123} on $\Si_*$, we obtain
\beaa
\bsplit
0 = \nab_{\nu_*}(\Om_*)_i =& \big(1+O(r^{-2})+O(r^{-1})\widecheck{\fb_*}+O(r^{-2})\widecheck{b}_*\big)\nab_3'\Om_i+\big(O(r^{-1})+\widecheck{\fb_*}+O(r^{-1})\widecheck{b}_*\big)\nab_a'\Om_i\\
&+\big(O(1)+\widecheck{b}_*+O(r^{-2})+O(r^{-1})\widecheck{\fb_*}\big)\nab_4'\Om_i
\end{split}
\eeaa
and 
\bea\lab{eq:auxiliaryidentitybetweennabstarOmi*andnabOmionSi*}
(\nab_*)_a(\Om_*)_i = \big(1+O(r^{-2})+O(r^{-1})\widecheck{\fb_*}\big)\nab_a'\Om_i+\big(O(r^{-1})+\widecheck{\fb_*}\big)\nab_3'\Om_i+\big(O(r^{-1})+\widecheck{\fb_*}\big)\nab_4'\Om_i
\eea
which together with \eqref{eq:identityusingasymptoticofKerrvaluefornabalphaOmitocontrollinearizednabalphaOmionSi*} and \eqref{eq:identityusingasymptoticofKerrvaluefornabalphaOmitocontrollinearizednabalphaOmionSi*:bis} yields on $\Si_*$
\beaa
\widecheck{\widecheck{\nab_3'\Om_i}} &=& O(mr^{-2})+O(r^{-1})\widecheck{\fb_*}+O(r^{-2})\widecheck{b}_*+r^{-2}\Ga_g'+O(r^{-1})\widecheck{\nab_*(\Om_*)_i},\\
\widecheck{\widecheck{\nab_4'\Om_i}} &=& O(mr^{-2})+O(r^{-3})\widecheck{\fb_*}+O(r^{-4})\widecheck{b}_*+r^{-2}\Ga_g'+O(r^{-3})\widecheck{\nab_*(\Om_*)_i},\\
\widecheck{\widecheck{\nab_a'\Om_i}} &=& O(mr^{-2})+O(r^{-2})\widecheck{\fb_*}+O(r^{-3})\widecheck{b}_*+r^{-3}\Ga_g'+O(1)\widecheck{\nab_*(\Om_*)_i}.
\eeaa
In view of \eqref{eq:dominantconditionforrstarcomparedtotaustartonS*}, \eqref{eq:assumptionsonMMforpartII:secondglobalframeauxassfor15derivatives}, \eqref{eq:assumptionsonSigmastarforpartII} and \eqref{eq:controlofwidechecknabstarOmstarii=123onSigma*}, we deduce 
\bea\lab{eq:controlofnab3primebabaprimeOmionSi*}
|\dk^{\leq 15}\widecheck{\widecheck{\nab_3'\Om_i}}|+|\dk^{\leq 15}\widecheck{\widecheck{\nab_a'\Om_i}}|\les \frac{1}{r^2}+\frac{\ep}{r\tau^{1+\dec}}\les\frac{\ep}{r\tau^{1+\dec}}\quad\textrm{on}\quad\Si_*.
\eea
and 
\bea\lab{eq:controlofnab4primeOmionSi*:nonsharp}
|\dk^{\leq 15}\widecheck{\widecheck{\nab_4'\Om_i}}|\les \frac{1}{r^2}+\frac{\ep}{r\tau^{1+\dec}}\les\frac{\ep}{r\tau^{1+\dec}}\quad\textrm{on}\quad\Si_*,
\eea
where \eqref{eq:controlofnab4primeOmionSi*:nonsharp} is non-sharp and used for convenience.

Next, we introduce the approximate Killing vectorfield $\Z'$ w.r.t. the global null frame $(e_4', e_3', e_a')$, i.e. 
\bea\lab{eq:defonZ'expressedindoubleprimedframe}
\Z '= \frac 1 2 \left(2(r^2+a^2)\Re(\Jk)^be_b' -a(\sin\th)^2 e_4' -\frac{a(\sin\th)^2\De}{ |q|^2} e_3'\right).
\eea
In view of \eqref{eq:relationbetweenJkJkpmandf0onSigmastar} and \eqref{eq:defonZ'expressedindoubleprimedframe}, we have
\beaa
\Z' &=& rf_0^be_b'+O(r^{-1})\nab'+O(1)e_4'+O(1)e_3'\quad\textrm{on}\quad\Si_*.
\eeaa
which together with \eqref{eq:definitionapproxKiilingvectorfieldZ*onSi*} and \eqref{eq:auxiliaryidentitybetweennabstarOmi*andnabOmionSi*} implies on $\Si_*$
\begin{align*}
(\nab_*)_{\Z_*}(\Om_*)_i =& rf_0^b(\nab_*)_b(\Om_*)_i\\
=& rf_0^b\big(1+O(r^{-2})+O(r^{-1})\widecheck{\fb_*}\big)\nab_b'\Om_i+\big(O(1)+r\widecheck{\fb_*}\big)\nab_3'\Om_i+\big(O(1)+r\widecheck{\fb_*}\big)\nab_4'\Om_i\\
=& \nab_{Z'}\Om_i +\big(O(r^{-1})+O(1)\widecheck{\fb_*}\big)\nab'\Om_i+\big(O(1)+r\widecheck{\fb_*}\big)\nab_4'\Om_i+\big(O(1)+r\widecheck{\fb_*}\big)\nab_3'\Om_i\\
=& \nab_{Z'}\Om_i +O(r^{-2})+\big(O(r^{-1})+O(1)\widecheck{\fb_*}\big)\widecheck{\widecheck{\nab'\Om_i}}+\big(O(1)+r\widecheck{\fb_*}\big)\widecheck{\widecheck{\nab_4'\Om_i}}+\big(O(1)+r\widecheck{\fb_*}\big)\widecheck{\widecheck{\nab_3'\Om_i}}.
\end{align*}
In view of \eqref{eq:dominantconditionforrstarcomparedtotaustartonS*}, \eqref{eq:assumptionsonSigmastarforpartII}, \eqref{eq:controlofnab3primebabaprimeOmionSi*} and \eqref{eq:controlofnab4primeOmionSi*:nonsharp}, we deduce
\beaa
|\dk^{\leq 15}(\nab_{\Z'}\Om_i-(\nab_*)_{\Z_*}(\Om_*)_i)|\les \frac{1}{r^2}+\frac{\ep}{r\tau^{1+\dec}}\les\frac{\ep}{r\tau^{1+\dec}}\quad\textrm{on}\quad\Si_*,
\eeaa
which together with \eqref{eq:controlof15weightedderivativesrenormalizedLiebtZstarOmstarionSigma*} and  \eqref{eq:initializationofOmii=123onSigma*fromOm*i} implies
\beaa
\big|\dk_*^{\leq 15}\big(\nab_{\Z'}\Om_i +\in_{ij3}\Om^j +\cos\th\dual\Om_i\big)\phi_*\big|\les\frac{\ep}{\tau^{1+\dec}}, \quad i=1,2,3,\quad\textrm{on}\quad \Si_*.
\eeaa
On the other hand, we have by Lemma \ref{lemma:basicpropertiesLiebTfasdiuhakdisug:chap9}
\beaa
\nab_{\Z'}\Om_i &=& \Lieb_{\Z'}\Om_i -\cos\th\dual\Om_i+r\Ga_b'+O(r^{-2}).
\eeaa
We deduce, using also \eqref{eq:dominantconditionforrstarcomparedtotaustartonS*} and \eqref{eq:assumptionsonMMforpartII:secondglobalframeauxassfor15derivatives}, 
\beaa
\big|\dk_*^{\leq 15}\big(\Lieb_{\Z'}\Om_i +\in_{ij3}\Om^j\big)\big|\les r^{-2}+\frac{\ep}{\tau^{1+\dec}}\les\frac{\ep}{\tau^{1+\dec}}, \quad i=1,2,3,\quad\textrm{on}\quad \Si_*.
\eeaa
Together with \eqref{eq:trnasrportequationinYtauforextensionOmii=123tpMM}, \eqref{commutatorbetweenLieTLieZandnabnab4nab3:2} and Remark \ref{rmk:formofYautforrgeqr0intermsofdoubleprimedframe}, we finally obtain
\bea\lab{eq:renormalizedLiebZprimeOmionSi*}
\big|\dk^{\leq 15}\big(\Lieb_{\Z'}\Om_i +\in_{ij3}\Om^j\big)\big|\les r^{-2}+\frac{\ep}{\tau^{1+\dec}}\les\frac{\ep}{\tau^{1+\dec}}, \quad i=1,2,3,\quad\textrm{on}\quad \Si_*.
\eea

\noindent {\bf Step 6.} Next, we extend \eqref{eq:controlofnab3primebabaprimeOmionSi*} to $\Mext\cap\{r\geq r_0+1\}$. To this end, we assume the following bootstrap assumptions 
\bea\lab{eq:bootass15weightedderivativeswidechecknabalphaOmiMMgeqr1}
\bsplit
|\dk^{\leq 15}\widecheck{\widecheck{\nab_4'\Om_i}}|&\leq\min\left\{\frac{\sqrt{\ep}}{r^{2}\tau^{\frac{1+\dec}{2}}}, \, \frac{\sqrt{\ep}}{r\tau^{1+\frac{3\dec}{4}}}\right\},\quad i=1,2,3,\,\,\,\textrm{on}\,\,\,\Mext\cap\{r\geq r_1\},\\
|\dk^{\leq 15}\big(\widecheck{\widecheck{\nab_3'\Om_i}}, \widecheck{\widecheck{\nab_a'\Om_i}}\big)|&\leq \frac{\sqrt{\ep}}{r\tau^{1+\frac{3\dec}{4}}},\quad i=1,2,3,\quad a=1,2,\quad\textrm{on}\quad\Mext\cap\{r\geq r_1\},
\end{split}
\eea
for some $r_1\geq r_0+1$. Then, we commute the transport equations \eqref{eq:trnasrportequationinYtauforextensionOmii=123tpMM} with $\Lieb_{\T'}$ using the first commutation formula in \eqref{eq:commutatorsYtauwithT'andrnab'}, where $\T'$ given by \eqref{eq:defonT'expressedindoubleprimedframe}. Together with Remark \ref{rmk:formofYautforrgeqr0intermsofdoubleprimedframe}, we obtain 
\beaa
\nab_{\Ytau}'\Lieb_{\T'}\Om_i &=& O(r^{-2})\Lieb_{\T'}\Om_j+O(r^{-3})\T'(\rr)+O(r^{-2})\T'(\cos(\widetilde{\th}))+\dk^{\leq 1}\Ga_g'\c \nab_3'\Om_i\\
&&+r^{-1}\dk^{\leq 1}(\Ga_b' \c \Om_i)\\
&=& O(r^{-2})\Lieb_{\T'}\Om_j+r^{-1}\dk^{\leq 1}\Ga_b'+\dk^{\leq 1}\Ga_g'\c\widecheck{\nab_3'\Om_i}+\Ga_b'\c\big(\widecheck{\nab_a'\Om_i}, \widecheck{\nab_4'\Om_i}\big)+r^{-2}\widetilde{\Ga}_b.
\eeaa
Commuting with $(\Lieb_{\T'}, r\nab')$ using \eqref{eq:commutatorsYtauwithT'andrnab'}, and integrating the resulting transport equation from $\Si_*$ where \eqref{eq:controlofnab3primebabaprimeOmionSi*} \eqref{eq:controlofnab4primeOmionSi*:nonsharp} hold, we infer, using also \eqref{eq:assumptionsonMMforpartII:secondglobalframeauxassfor15derivatives}, \eqref{eq:decaypropertiesofwidehatGab1widehatGag1} and \eqref{eq:bootass15weightedderivativeswidechecknabalphaOmiMMgeqr1}
\beaa
|(\Lieb_{\T'}, r\nab', r\nab_{\Ytau}')^{\leq 15}\Lieb_{\T'}\Om_i| &\les& \frac{\ep}{r\tau^{1+\dec}}+\int_r^{r_*}\frac{\ep}{{r'}^2\tau^{1+\frac{3\dec}{4}}}dr'\les\frac{\ep}{r\tau^{1+\frac{3\dec}{4}}}.
\eeaa
Using Lemma \ref{lemma:basicpropertiesLiebTfasdiuhakdisug:chap9} and the fact that $(\nab_4', \nab_3')$ are generated by $(\nab_{\T'}, \nab_{\Ytau}', \nab')$, we deduce 
\bea\lab{eq:bootass15weightedderivativeswidechecknabalphaOmiMMgeqr1:improvementLiebT'Omi}
|\dk^{\leq 15}\Lieb_{\T'}\Om_i| &\les& \frac{\ep}{r\tau^{1+\frac{3\dec}{4}}},\quad i=1,2,3,\quad\textrm{on}\quad\Mext\cap\{r\geq r_1\}.
\eea

Next, we control $\widecheck{\widecheck{\nab_a\Om_i}}$ on $\Mext\cap\{r\geq r_1\}$. To this end, we commute the transport equations \eqref{eq:trnasrportequationinYtauforextensionOmii=123tpMM} with $r\nab'$ using the second commutation formula in \eqref{eq:commutatorsYtauwithT'andrnab'} which yields
\beaa
\nab_{\Ytau}'(r\nab'\Om_i) &=& O(\rr^{-2})r\nab'\Om_i+O(\rr^{-2})r\nab'(\rr)+O(\rr^{-2})r\nab'(\cos(\widetilde{\th}))\\
&&+ \big(O(r^{-2})+\dk^{\leq 1}\Ga_g'\big)\c r\nab'\Om_i+\big(O(r^{-3})+\dk^{\leq 1}\Ga_g'\big)\c\nab_{\Ytau}'\Om_i\\
&&+\big(O(r^{-3})+\dk^{\leq 1}\Ga_g'\big)\Lieb_{\T'}\Om_i.
\eeaa
Linearizing, and using the transport equations in \eqref{definitionofthenewcoordinatetilder} \eqref{definitionofthenewcoordinatetildethetandtphi} and \eqref{eq:trnasrportequationinYtauforextensionOmii=123tpMM} along $\Ytau$, we deduce 
\beaa
\nab_{\Ytau}'(r\widecheck{\widecheck{\nab'\Om_i}}) &=& O(\rr^{-2})r\widecheck{\widecheck{\nab'\Om_j}}+\frac{\rr-r}{r^3}+\frac{\cos(\widetilde{\th})-\cos\th}{r^3}+O(r^{-1})\nab'(\rr)+r^{-1}\widetilde{\Ga}_g\\
&&+\dk^{\leq 1}\Ga_g'+ \big(O(r^{-2})+\dk^{\leq 1}\Ga_g'\big)\c\widecheck{\widecheck{r\nab'\Om_i}}+\big(O(r^{-3})+\dk^{\leq 1}\Ga_g'\big)\Lieb_{\T'}\Om_i.
\eeaa
Commuting with $(\Lieb_{\T'}, r\nab')$ using \eqref{eq:commutatorsYtauwithT'andrnab'}, and integrating the resulting transport equation from $\Si_*$ where \eqref{eq:controlofnab3primebabaprimeOmionSi*} holds, we infer, using also\footnote{Note in particular that \eqref{eq:esitmatefor15weightedderivativesoffirstorderdrivgivescheckofrr} yields
\beaa
r^{-1}|\dk^{\leq 15}\nab'(r)|\les \frac{\ep}{r^{1+\frac{\dec}{4-\dec}}\tau^{1+\frac{3\dec}{4}}}. 
\eeaa}
\eqref{eq:assumptionsonMMforpartII:secondglobalframeauxassfor15derivatives}, \eqref{eq:esitmatefor16weightedderivativesofrrminusr}, \eqref{eq:esitmatefor15weightedderivativesoffirstorderdrivgivescheckofrr}, \eqref{eq:esitmatefor16weightedderivativesoftxbminusxb}, \eqref{eq:decaypropertiesofwidehatGab1widehatGag1} and \eqref{eq:bootass15weightedderivativeswidechecknabalphaOmiMMgeqr1:improvementLiebT'Omi}
\beaa
|(\Lieb_{\T'}, r\nab', r\nab_{\Ytau}')^{\leq 15}\widecheck{\widecheck{\nab'\Om_i}}| &\les& \frac{\ep}{r\tau^{1+\dec}}+\frac{1}{r}\int_r^{r_*}\frac{\ep}{{r'}^{1+\frac{\dec}{4-\dec}}\tau^{1+\frac{3\dec}{4}}}dr'\les\frac{\ep}{r\tau^{1+\frac{3\dec}{4}}}.
\eeaa
Using Lemma \ref{lemma:basicpropertiesLiebTfasdiuhakdisug:chap9} and the fact that $(\nab_4', \nab_3')$ are generated by $(\nab_{\T'}, \nab_{\Ytau}', \nab')$, we deduce 
\bea\lab{eq:bootass15weightedderivativeswidechecknabalphaOmiMMgeqr1:improvementdoublewidechecknabaOmi}
|\dk^{\leq 15}\widecheck{\widecheck{\nab'\Om_i}}| &\les& \frac{\ep}{r\tau^{1+\dec}}+\frac{1}{r}\int_r^{r_*}\frac{\ep}{{r'}^2\tau^{1+\frac{3\dec}{4}}}dr'\les\frac{\ep}{r\tau^{1+\frac{3\dec}{4}}}.
\eea

Next, using \eqref{eq:trnasrportequationinYtauforextensionOmii=123tpMM}, Remark \ref{rmk:formofYautforrgeqr0intermsofdoubleprimedframe}, \eqref{eq:assumptionsonMMforpartII:secondglobalframeauxassfor15derivatives},  
\eqref{eq:bootass15weightedderivativeswidechecknabalphaOmiMMgeqr1:improvementLiebT'Omi},  \eqref{eq:bootass15weightedderivativeswidechecknabalphaOmiMMgeqr1:improvementdoublewidechecknabaOmi}, the decomposition of $\T'$ on the null frame $(e_4', e_3', e_1', e_2')$ and Lemma \ref{lemma:basicpropertiesLiebTfasdiuhakdisug:chap9}, we infer
\beaa
r^2\left|\dk^{\leq 15}\widecheck{\widecheck{\nab_4'\Om_i}}\right| &\les& r^2\left|\dk^{\leq 15}\Ga_g'\right|+r\left|\dk^{\leq 15}\left(\Lieb_{\T'}\Om_i, \widecheck{\widecheck{\nab_a'\Om_i}}\right)\right|\\
&\les& \frac{\ep}{\tau^{\frac{1}{2}+\dec}}\quad\textrm{on}\quad\Mext\cap\{r\geq r_1\},
\eeaa
which together with \eqref{eq:bootass15weightedderivativeswidechecknabalphaOmiMMgeqr1:improvementLiebT'Omi},  \eqref{eq:bootass15weightedderivativeswidechecknabalphaOmiMMgeqr1:improvementdoublewidechecknabaOmi}, the decomposition of $\T'$ on the null frame $(e_4', e_3', e_1', e_2')$ and Lemma \ref{lemma:basicpropertiesLiebTfasdiuhakdisug:chap9} implies the following improvement of the bootstrap assumptions \eqref{eq:bootass15weightedderivativeswidechecknabalphaOmiMMgeqr1}
\beaa
\bsplit
|\dk^{\leq 15}\widecheck{\widecheck{\nab_4'\Om_i}}|&\les\min\left\{\frac{\ep}{r^{2}\tau^{\frac{1+\dec}{2}}}, \, \frac{\ep}{r\tau^{1+\frac{3\dec}{4}}}\right\},\quad i=1,2,3,\,\,\,\textrm{on}\,\,\,\Mext\cap\{r\geq r_1\},\\
|\dk^{\leq 15}\big(\widecheck{\widecheck{\nab_3'\Om_i}}, \widecheck{\widecheck{\nab_a'\Om_i}}\big)|&\les \frac{\ep}{r\tau^{1+\frac{3\dec}{4}}},\quad i=1,2,3,\quad a=1,2,\quad\textrm{on}\quad\Mext\cap\{r\geq r_1\}.
\end{split}
\eeaa
We thus deduce $r_1=r_0+1$ and 
\beaa
\bsplit
|\dk^{\leq 15}\widecheck{\widecheck{\nab_4'\Om_i}}|&\les\min\left\{\frac{\ep}{r^{2}\tau^{\frac{1+\dec}{2}}}, \, \frac{\ep}{r\tau^{1+\frac{3\dec}{4}}}\right\},\quad i=1,2,3,\,\,\,\textrm{on}\,\,\,\Mext\cap\{r\geq r_0+1\},\\
|\dk^{\leq 15}\big(\widecheck{\widecheck{\nab_3'\Om_i}}, \widecheck{\widecheck{\nab_a'\Om_i}}\big)|&\les \frac{\ep}{r\tau^{1+\frac{3\dec}{4}}},\quad i=1,2,3,\quad a=1,2,\quad\textrm{on}\quad\Mext\cap\{r\geq r_0+1\}.
\end{split}
\eeaa
Then, we easily extend this estimate from $\{r=r_0+1\}$ to $\MM\cap\{r\leq r_0+1\}$ by proceeding similarly as above which yields 
\beaa
\bsplit
|\dk^{\leq 15}\widecheck{\widecheck{\nab_4'\Om_i}}|&\les\min\left\{\frac{\ep}{r^{2}\tau^{\frac{1+\dec}{2}}}, \, \frac{\ep}{r\tau^{1+\frac{3\dec}{4}}}\right\},\quad i=1,2,3,\,\,\,\textrm{on}\,\,\,\MM,\\
|\dk^{\leq 15}\big(\widecheck{\widecheck{\nab_3'\Om_i}}, \widecheck{\widecheck{\nab_a'\Om_i}}\big)|&\les \frac{\ep}{r\tau^{1+\frac{3\dec}{4}}},\quad i=1,2,3,\quad a=1,2,\quad\textrm{on}\quad\MM.
\end{split}
\eeaa
Together with the change of frame formulas \eqref{frametransformation:fromfirsttosecondglobalframe}, and using also the control of $(f, \fb)$ in \eqref{eq:assumptionsonMMforpartII:secondglobalframeauxassfor15derivatives}, we infer 
\bea\lab{eq:estimatefor15weightedderivativesofdoublewidechechnabalphaOmionMM:forOmifinalestimateproof}
\bsplit
|\dk^{\leq 15}\widecheck{\widecheck{\nab_4\Om_i}}|&\les\min\left\{\frac{\ep}{r^{2}\tau^{\frac{1+\dec}{2}}}, \, \frac{\ep}{r\tau^{1+\frac{3\dec}{4}}}\right\},\quad i=1,2,3,\,\,\,\textrm{on}\,\,\,\MM,\\
|\dk^{\leq 15}\big(\widecheck{\widecheck{\nab_3\Om_i}}, \widecheck{\widecheck{\nab_a\Om_i}}\big)|&\les \frac{\ep}{r\tau^{1+\frac{3\dec}{4}}},\quad i=1,2,3,\quad a=1,2,\quad\textrm{on}\quad\MM.
\end{split}
\eea
Together with \eqref{eq:definitionlinearizedquantiiesnabaphaOmiwithdoublecheck}, this implies the stated estimate \eqref{eq:estimatefor15weightedderivativesofdoublewidechechnabalphaOmionMM} for $\widecheck{\widecheck{M^j_{i\a}}}$.

\noindent {\bf Step 7.} Finally, we extend \eqref{eq:renormalizedLiebZprimeOmionSi*} to $\MM$. We have
\beaa
\nab_{\Ytau}'\big(\Lieb_{\Z'}\Om_i +\in_{ij3}\Om^j\big) &=& [\nab_{\Ytau}', \Lieb_{\Z'}]\Om_i+\Lieb_{\Z'}\nab_{\Ytau}'\Om_i+\in_{ij3}\nab_{\Ytau}'\Om^j
\eeaa
which together with \eqref{eq:trnasrportequationinYtauforextensionOmii=123tpMM} yields
\beaa
\nab_{\Ytau}'\big(\Lieb_{\Z'}\Om_i +\in_{ij3}\Om^j\big) &=& O(r^{-2})\big(\Lieb_{\Z'}\Om_i +\in_{ij3}\Om^j\big)+[\nab_{\Ytau}', \Lieb_{\Z'}]\Om_i +\rr^{-3}\Z'(\rr)+\rr^{-2}\Z'(\cos(\widetilde{\th}))\\
&=& O(r^{-2})\big(\Lieb_{\Z'}\Om_i +\in_{ij3}\Om^j\big)+[\nab_{\Ytau}', \Lieb_{\Z'}]\Om_i +r^{-1}\widetilde{\Ga}_g.
\eeaa
Hence, using \eqref{commutatorbetweenLieTLieZandnabnab4nab3:2}, \eqref{eq:additionalpropertiesseconddgolbalnullframe} and Remark \ref{rmk:formofYautforrgeqr0intermsofdoubleprimedframe}, this implies
\beaa
\nab_{\Ytau}'\big(\Lieb_{\Z'}\Om_i +\in_{ij3}\Om^j\big) &=& O(r^{-2})\big(\Lieb_{\Z'}\Om_i +\in_{ij3}\Om^j\big)+\dk^{\leq 1}\Ga_g' +r^{-1}\widetilde{\Ga}_g,
\eeaa
Commuting with $(\Lieb_{\T'}, r\nab')$ using \eqref{eq:commutatorsYtauwithT'andrnab'}, and integrating the resulting transport equation from $\Si_*$ where \eqref{eq:renormalizedLiebZprimeOmionSi*} holds, we infer, using also \eqref{eq:assumptionsonMMforpartII:secondglobalframeauxassfor15derivatives} and \eqref{eq:decaypropertiesofwidehatGab1widehatGag1}, 
\beaa
\big|(\Lieb_{\T'}, r\nab', r\nab_{\Ytau}')^{\leq 15}\big(\Lieb_{\Z'}\Om_i +\in_{ij3}\Om^j\big)\big| &\les& \frac{\ep}{\tau^{1+\dec}}+\int_r^{r_*}\frac{\ep}{{r'}^{1+\frac{\dec}{2}}\tau^{1+\frac{3\dec}{4}}}dr'\les\frac{\ep}{\tau^{1+\frac{3\dec}{4}}}.
\eeaa
Using Lemma \ref{lemma:basicpropertiesLiebTfasdiuhakdisug:chap9} and the fact that $(\nab_4', \nab_3')$ are generated by $(\nab_{\T'}, \nab_{\Ytau}', \nab')$, we deduce 
\beaa
\big|\dk^{\leq 15}\big(\Lieb_{\Z'}\Om_i +\in_{ij3}\Om^j\big)\big| &\les& \frac{\ep}{\tau^{1+\frac{3\dec}{4}}}, \quad i=1,2,3, \quad\textrm{on}\quad\MM.
\eeaa
Together with the change of frame formulas \eqref{frametransformation:fromfirsttosecondglobalframe}, and using also the control of $(f, \fb)$ in \eqref{eq:assumptionsonMMforpartII:secondglobalframeauxassfor15derivatives}, we infer 
\bea\lab{eq:alsmostdoneasestimateforrenomralizedLiebZOmiinsteadofcoordderivative}
\big|\dk^{\leq 15}\big(\Lieb_{\Z}\Om_i +\in_{ij3}\Om^j\big)\big| &\les& \frac{\ep}{\tau^{1+\frac{3\dec}{4}}}, \quad i=1,2,3, \quad\textrm{on}\quad\MM,
\eea
where the approximate Killing vectorfield $\Z$ is given by \eqref{eq:definitionofTandPhithataretheapproximateKillingvectorifeldinKerrpert}.

Next, we use the fact that 
$$\Z(\rr)=r^2\widetilde{\Ga}_g,\qquad \Z(\ttt)=r^2\widetilde{\Ga}_g, \qquad \Z(\tx^b)=\pr_{\widetilde{\tphi}}(\tx^b)+r\widetilde{\Ga}_g,$$
which implies 
$$\Z=\pr_{\widetilde{\tphi}}+r^2\widetilde{\Ga}_g\pr_{\rr}+r^2\widetilde{\Ga}_g\pr_{\ttt}+r\widetilde{\Ga}_g\pr_{\tx^b}$$
and hence
\beaa
\Lieb_{\Z}\Om_i=\Lieb_{\widetilde{\tphi}}\Om_i+r\widetilde{\Ga}_g+r^2\widetilde{\Ga}_g\c\big(\widecheck{\widecheck{\nab_3\Om_i}}, \widecheck{\widecheck{\nab_4\Om_i}}, \widecheck{\widecheck{\nab_a\Om_i}}\big)
\eeaa
which together with \eqref{eq:decaypropertiesofwidehatGab1widehatGag1},  \eqref{eq:estimatefor15weightedderivativesofdoublewidechechnabalphaOmionMM:forOmifinalestimateproof} and \eqref{eq:alsmostdoneasestimateforrenomralizedLiebZOmiinsteadofcoordderivative} implies the stated estimate \eqref{eq:estimatefor15weightedderivativesofrenormalizedLiebtphiOmi}. This concludes the proof of Lemma \ref{lemma:constructionofhorizontal1formsOmiinperturbationsofKerr}.
\end{proof}

%%%%%%%%%%%%%%%%%%%%%%%%%%%%%%%%%%%%%%%%%%%

\section{Extension of Theorems M1 and M2 to the full subextremal range}
\lab{sec:decayestimatesforAandAb:analogPartIIGKS22}

%%%%%%%%%%%%%%%%%%%%%%%%%%%%%%%%%%%%%%%%%%%

The goal of this section is to prove Theorems \ref{theoremM1:intro} and \ref{theoremM2:intro} on the extension of the decay estimates for Teukolsky equations of Theorems M1 and M2 in Section 3.7.1 of \cite{KS:Kerr} to the full subextremal range.

%%%%%%%%%%%%%%%%%%%%%%%%%%%%%
   
\subsection{Preliminaries}
        
%%%%%%%%%%%%%%%%%%%%%%%%%%%%%

%%%%%%%%%%%%%%%%%%%%%%%%%%%%%
   
    \subsubsection{Main norms for Teukolsky equations}
        \lab{subsection:basicnormsforpsi}
        
%%%%%%%%%%%%%%%%%%%%%%%%%%%%%

We introduce in this section the energy, Morawetz, flux and $r^p$ norms for solutions to Teukolsky equations. Below, $\psi$ denotes a horizontal tensor in $\sk_2(\mathbb{C})$.

   {\bf 1. Morawetz  norm.} 
\bea
\bsplit
\M_\de[\psi](\tau_1,\tau_2):=&\int_{\Mntrap(\tau_1,\tau_2)}\left(\frac{|\nab_3\psi|^2}{r^{1+\de}} +\frac{|\nabla\psi|^2}{r}\right)
+\int_{\MM(\tau_1,\tau_2)}\left(\frac{|\nab_{\Rhat}\psi|^2}{r^{1+\de}} +\frac{|\psi|^2}{r^3}\right),
\end{split}
\eea
for any given $0\leq\de\leq 1$, where $\Rhat$ denotes the following vectorfield 
\bea\lab{eq:defofRhat}
\Rhat := \frac 1 2 \left( \frac{|q|^2}{r^2+a^2} e_4-\frac{\De}{r^2+a^2}  e_3\right).
\eea

{\bf 2. Energy norm.} 
\bea
\E[\psi](\tau) :=\int_{\Si (\tau)} \Big( |\nab_4\psi|^2 +   r^{-2}|\nab_3\psi|^2 +|\nab\psi|^2 + r^{-2}|\psi|^2\Big).
\eea

{\bf 3. Flux  norm.} 
\bea
\bsplit
\F[\psi](\tau_1,\tau_2) : =& \F_{\AA}[\psi](\tau_1,\tau_2)+\F_{\Si_*}[\psi](\tau_1,\tau_2),\\
\F_{\AA}[\psi](\tau_1,\tau_2) := & \int_{\AA(\tau_1, \tau_2)}\Big( |\nab_4\psi |^2+|\nab_3\psi|^2+|\nab\psi|^2+r^{-2} | \psi |^2\Big), \\
\F_{\Si_*}[\psi](\tau_1,\tau_2) := & \int_{\Si_*(\tau_1, \tau_2)}\Big( |\nab_4\psi |^2+|\nab_3\psi|^2+|\nab\psi|^2+ r^{-2} | \psi |^2\Big).
\end{split}
\eea

{\bf 4. Norm for the inhomogeneous term.} For $N\in\sk_k(\mathbb{C})$, $k=1,2$,
\bea\lab{eq:defmathcalNpsif:00}
\widetilde{\mathcal{N}}_\de[\psi, N](\tau_1, \tau_2) &:=& \!\!\int_{\MM(\tau_1, \tau_2)}r^{1+\de}|N|^2\\
&+&\!\!\!\!\min\left[\left(\int_{\Mtrap(\tau_1, \tau_2)}|N|^2\right)^{\frac{1}{2}} \left(\int_{\Mtrap(\tau_1, \tau_2)}|\dk{^{\leq 1}}\psi|^2\right)^{\frac{1}{2}}, 
\int_{\Mtrap(\tau_1, \tau_2)}\tau^{1+\de}|N|^2\right]. \nn
\eea

{\bf 5. Weighted  bulk norm.} For $0<p<2$, we define
\bea
\B_{p}[\psi](\tau_1, \tau_2)&:=& \M_\de[\psi](\tau_1,\tau_2) +\int_{\MM_{r\geq 10m}(\tau_1,\tau_2)} r^{p-3}\Big(|\dk \psi|^2 +|\psi|^2  \Big).
\eea

{\bf 6. Weighted  energy  norm.} For $0<p<2$, we define
\bea
\E_p[\psi](\tau) := \left\{\ba{l} 
\E[\psi]+\displaystyle\int_{\Si_{r\geq 10m}(\tau)} r^{p}\Big(|\nab_4\psi|^2  + r^{-2}|\psi|^2  \Big) \quad\textrm{for }p\leq 1-\de,\\[3mm]
\E[\psi]+\displaystyle\int_{\Si_{r\geq 10m}(\tau)} r^{p}\Big(|r^{-1}\nab_4(r\psi)|^2  + r^{-p-1-\de}|\psi|^2  \Big) \quad\textrm{for }p> 1-\de.
\ea\right.
\eea

{\bf 7. Weighted flux   norm.} For $0<p<2$, we define
\bea
\F_p[\psi](\tau_1,\tau_2):= \F[\psi](\tau_1,\tau_2) +\int_{\Si_*(\tau_1, \tau_2)} r^p\Big(|\nab_4 \psi|^2 +|\nab\psi|^2+ r^{-2}|\psi|^2 \Big).
\eea

{\bf 8. Combined norms.}  We denote the combined norm
\bea
\BEF_p[\psi](\tau_1,\tau_2):= \sup_{\tau\in[\tau_1, \tau_2]}\E_p[\psi](\tau) +\B_p[\psi](\tau_1, \tau_2) +\F_p[\psi](\tau_1, \tau_2),
\eea
and define similarly the combined norms $\EF_p[\psi](\tau_1,\tau_2)$, $\BE_p[\psi](\tau_1,\tau_2)$ and $\BF_p[\psi](\tau_1,\tau_2)$.

{\bf 9. Weighted norm for the inhomogeneous term.}.} For $0<p<2$, we define
\bea
\NN_p[\psi,  N](\tau_1, \tau_2) &=& \widetilde{\mathcal{N}}_\de[\psi,  N](\tau_1, \tau_2)+\left| \int_{\MM_{r\geq 10m}}  r^{p-1}  \, \Re(\nab_4 (r\psi ) \c  \ov{N})\right|.
\eea

{\bf 10. Higher order norms.}  We define the higher  derivative norms $\M_\de^s[\psi]$, $\E^s[\psi]$, $\F^s[\psi]$, $\B_p^s[\psi]$, $\E_p^s[\psi]$,  $\F_p^s[\psi]$, 
 by the  general procedure for a norm $\Q[\psi]$, i.e. 
 \beaa
 \Q^s[\psi]=\sum_{k\le s}\Q[\dk^k\psi].
 \eeaa
Also, we define $\widetilde{\mathcal{N}}_\de^s[\psi,  N]$ and $\NN_p^s[\psi, N]$ by
\beaa
\widetilde{\mathcal{N}}_\de^s[\psi,  N]=\sum_{k\le s}\widetilde{\mathcal{N}}_\de[\dk^k\psi, \dk^k N], \qquad \NN_p^s[\psi,  N]=\sum_{k\le s}\NN_p[\dk^k\psi, \dk^k N].
\eeaa

%%%%%%%%%%%%%%%%%%%%%%%%%%%%%%%%%%%%%%%%%%%

\subsubsection{Commutators with the D'Alembertian}

%%%%%%%%%%%%%%%%%%%%%%%%%%%%%%%%%%%%%%%%%%%

\begin{proposition}\label{PROPOSITION:COMMUTATORS-SQUARED-LIED-T-Z} 
The following  commutation  formulas hold true  for any horizontal tensor $\psi$:
\beaa
\, [ {\Lieb_\T}, \squared_2]\psi = {\dk \big(\Ga_g \c \dk \psi\big)}, \qquad \, [{\Lieb_\Z}, \squared_2]\psi =    {  \dk(\Ga_g   \c \dk \psi)+r\Ga_b \c \squared_2 \psi+\Ddot_3 \big(r \xi \c \nab_3 \psi  \big). }
\eeaa
\end{proposition}

\begin{proof} 
See Proposition 4.3.4 in \cite{GKS22}.
\end{proof}

%%%%%%%%%%%%%%%%%%%%%%%%%%%%%%%%%%%%%%%%%%%

\subsubsection{Interpolated estimates for $(\Ga_b, \Ga_g)$}
\lab{sec:interpolatedestimatesforGabGag}

%%%%%%%%%%%%%%%%%%%%%%%%%%%%%%%%%%%%%%%%%%%

We denote from now on by $\kst$ an integer such that 
\bea\lab{eq:choiceofinterpolationintegerk*forTheoremM1andM2}
0<\kst - \frac{\kl}{2}\leq\frac{\kl\dec}{16},
\eea
which exists in view of the assumption $\kl\dec\gg 1$, see \eqref{eq:constraintsonthemainsmallconstantsepanddelta}. Then, interpolating between the assumptions in \eqref{eq:assumptionsonMMforpartII}, we obtain  
\beaa
|\dk^{\leq\kst}\Ga_g| &\les& |\dk^{\leq \frac{\kl}{2}}\Ga_g|^{1-\frac{2\kst -\kl}{\kl}}|\dk^{\leq \kl}\Ga_g|^{\frac{2\kst -\kl}{\kl}}\les \frac{\ep}{r^2\tau^{(\frac{1}{2}+\dec)(1-\frac{2\kst-\kl}{\kl})}},\\
|\dk^{\leq\kst}\Ga_b| &\les& |\dk^{\leq \frac{\kl}{2}}\Ga_b|^{1-\frac{2\kst -\kl}{\kl}}|\dk^{\leq \kl}\Ga_b|^{\frac{2\kst -\kl}{\kl}}\les \frac{\ep}{r\tau^{(1+\dec)(1-\frac{2\kst-\kl}{\kl})}}
\eeaa
and since we have from \eqref{eq:choiceofinterpolationintegerk*forTheoremM1andM2}
\beaa
(1+\dec)\frac{2\kst-\kl}{\kl}\leq \frac{4}{\kl}\left(\kst - \frac{\kl}{2}\right)\leq \frac{4}{\kl}\frac{\kl\dec}{16}=\frac{\dec}{4},
\eeaa
we infer
\bea\lab{eq:interpolationestimatesforGabGagfork*derivativeswithk*largerthanklover2:forproofTheoremM1andM2}
r^2|\dk^{\leq k}\Ga_g|\leq \frac{\ep}{\tau^{\frac{1}{2}+\frac{3\dec}{4}}}, \qquad r|\dk^{\leq k}(\Ga_g, \Ga_b)| \leq \frac{\ep}{\tau^{1+\frac{3\dec}{4}}}, \qquad  k\leq \kst.
\eea

%%%%%%%%%%%%%%%%%%%%%%%%%%%%%%%%%%%%%%%%%%%

\subsubsection{Control of error terms}

%%%%%%%%%%%%%%%%%%%%%%%%%%%%%%%%%%%%%%%%%%%

The following lemma allows us to control error terms arising in Morawetz estimates.
\begin{lemma}\lab{lemma:controloftheerrortermsGagpsiinwaveeq}
We have, for any $s\leq\kst-1$, with $\kst$ an integer satisfying \eqref{eq:choiceofinterpolationintegerk*forTheoremM1andM2}, 
\bea
\widetilde{\NN}^{s}_\de[\psi_1, \dk(\Ga_g\c\psi_2)](\tt_1, \tt_2) &\les& \ep\B^{s}_\de[\psi_2](\tau_1, \tau_2)+\ep\sup_{\tau\in[\tau_1, \tau_2]}\E^{s}[\psi_2](\tau).
\eea
\end{lemma}

\begin{proof}
In view of \eqref{eq:defmathcalNpsif:00}, we have
\begin{align*}
\widetilde{\NN}^{s}_\de[\psi_1, \dk(\Ga_g\c\psi_2)](\tt_1, \tt_2) \les& \int_{\MM(\tau_1, \tau_2)}r^{1+\de}|\dk^{\leq s+1}(\Ga_g\c\psi_2)|^2+\int_{\Mtrap(\tau_1, \tau_2)}\tau^{1+\de}|\dk^{\leq s+1}(\Ga_g\c\psi_2)|^2\\
\les& \int_{\Mntrap(\tau_1, \tau_2)}r^{1+\de}|\dk^{\leq s+1}(\Ga_g\c\psi_2)|^2 +\int_{\Mtrap(\tau_1, \tau_2)}\tau^{1+\de}|\dk^{\leq s+1}(\Ga_g\c\psi_2)|^2
\end{align*}
which together with \eqref{eq:interpolationestimatesforGabGagfork*derivativeswithk*largerthanklover2:forproofTheoremM1andM2} implies, for $s\leq\kst-1$,
\beaa
\widetilde{\NN}^{s}_\de[\psi_1, \dk(\Ga_g\c\psi_2)](\tt_1, \tt_2) &\les& \ep\B^{s}_\de[\psi_2](\tau_1, \tau_2)+\ep\left(\int_1^{+\infty}\tau^{-1+\de-\frac{3}{2}\dec}\right)\sup_{\tau\in[\tau_1, \tau_2]}\E^{s}[\psi_2](\tau)\\
&\les& \ep\B^{s}_\de[\psi_2](\tau_1, \tau_2)+\ep\sup_{\tau\in[\tau_1, \tau_2]}\E^{s}[\psi_2](\tau),
\eeaa
where we have used the fact that $\dec>\de$ in view of \eqref{eq:constraintsonthemainsmallconstantsepanddelta}. This concludes the proof of Lemma \ref{lemma:controloftheerrortermsGagpsiinwaveeq}.
\end{proof}

%%%%%%%%%%%%%%%%%%%%%%%%%%%%%%%%%%%%%%%%%%%%%%%%%%%%%%%

\subsection{Energy-Morawetz estimates for Teukolsky equations in $\MM$}
\lab{sec:energyMorawetzesitmatesforTeukoslkyonMM:higherorderderivatives:00}

%%%%%%%%%%%%%%%%%%%%%%%%%%%%%%%%%%%%%%%%%%%%%%%%%%%%%%%

%%%%%%%%%%%%%%%%%%%%%%%%%%%%%%%%%%%%%%%%%%%%%%%%%%

\subsubsection{Energy-Morawetz estimates for Teukolsky for up to 14 derivatives}
\lab{sec:energyMorawetzesitmatesforTeukoslkyonMM:upto15derivatives}

%%%%%%%%%%%%%%%%%%%%%%%%%%%%%%%%%%%%%%%%%%%%%%%%%%

In this section, we derive energy-Morawetz estimates for up to 14 derivatives for solutions to Teukolsky equation. To this end, we consider $\pmb\phi_s^{(p)}\in\sk_2(\mathbb{C})$, $s=\pm 2$, $p=0,1,2$, given by \eqref{eq:definitionofthephiplus2phierarchy:perturbationofKerr} \eqref{eq:definitionofthephiminus2phierarchy:perturbationofKerr}, and write the Teukolsky wave system as in Corollary \ref{cor:derivationoftheTeukolskytensorialwavesystemfors=plusminus2:kerrpert} w.r.t. the coordinates system $(\ttt, \rr, \tx^1, \tx^2)$ on $\MM$, where $\ttt$, $\rr$ and $(\tx^1, \tx^2)$ are introduced respectively in Lemma \ref{lemma:controloftildeuonMext}, and in \eqref{definitionofthenewcoordinatetilder} and 
\eqref{definitionofthenewcoordinatetildethetandtphi}, i.e.,
\bsub
\lab{eq:TensorialTeuSysandlinearterms:rescaleRHScontaine2:general:Kerrperturbation:tildeform}
\bea
\lab{eq:TensorialTeuSys:rescaleRHScontaine2:general:Kerrperturbation:tildeform}
\bigg(\squared_2 -\frac{4ia\cos\widetilde{\th}}{|\widetilde{q}|^2}\nab_{\pr_{\ttt}}- \frac{4-2\de_{p0}}{|\widetilde{q}|^2}\bigg){\phis{p}} = \L_{s}^{(p)}[\pmb\phi_{s}]+\N_{W,s}^{(p)}, \quad s=\pm2, \quad p=0,1,2,
\eea
where the linear coupling terms $\L_{s}^{(p)}[{\pmb\phi_s}]$ have the following schematic forms
\bea
\lab{eq:tensor:Lsn:onlye_2present:general:Kerrperturbation:bis:tildeform}
\bsplit
{\L_{s}^{(0)}[\pmb\phi_{s}]}={}& (2s\rr^{-3} +O(m\rr^{-4}))\phis{1}+ O(m\rr^{-3}) \nab_{\widetilde{\Xcal}_s}^{\leq 1}\phis{0},\\
{\L_{s}^{(1)}[\pmb\phi_{s}]}={}& (s\rr^{-3} +O(m\rr^{-4}))\phis{2}+ O(m\rr^{-3}) \nab_{\widetilde{\Xcal}_s}^{\leq 1}  \phis{1}+O(m\rr^{-2})\nab_{\pr_{\widetilde{\tphi}}+a\pr_{\ttt}}^{\leq 1}\phis{0},\\
{\L_{s}^{(2)}[\pmb\phi_{s}]}={}&O(m\rr^{-3})\phis{2}+O(m\rr^{-2})\nab_{\pr_{\widetilde{\tphi}}+a\pr_{\ttt}}^{\leq 1}\phis{1}+O(m^2\rr^{-2})\phis{0},
\end{split}
\eea
\esub
with $\widetilde{\Xcal}_s$, $s=\pm 2$, being the regular vectorfields given by
\bea\lab{eq:formofregularhorizontalvectorfieldmathcalXs:Kerrperturbation:cor:000}
\widetilde{\Xcal}_{s} &:=& s\big(\pr_{\widetilde{\tphi}}+a(\sin\widetilde{\th})^2\pr_{\ttt}\big) - \frac{2a\cos\widetilde{\th}}{\rr}\sin\widetilde{\th}\pr_{\widetilde{\th}}, \quad s=\pm 2, 
\eea
with all the coefficients in  \eqref{eq:tensor:Lsn:onlye_2present:general:Kerrperturbation:bis:tildeform} being independent of coordinates $\ttt$ and 
$\widetilde{\tphi}$, and with the coefficients in front of the terms $\phis{2}$ and $\nab_{\pr_{\widetilde{\tphi}}+a\pr_{\ttt}}\phis{1}$ on the RHS of equation of ${\L_{s}^{(2)}[\pmb\phi_{s}]}$ in \eqref{eq:tensor:Lsn:onlye_2present:general:Kerrperturbation:bis:tildeform} being real functions.

Moreover, the tensorial Teukolsky transport equations w.r.t. the coordinates system $(\ttt, \rr, \tx^1, \tx^2)$ on $\MM$ are given by
\bsub\lab{def:TensorialTeuScalars:wavesystem:Kerrperturbation:tildeform}
 \bea
\lab{def:TensorialTeuScalars:wavesystem:Kerrperturbation:+2:tildeform}
\nab_3 \left(\frac{\rr\ov{\widetilde{q}}}{\widetilde{q}}\left(\frac{\rr^2}{|\widetilde{q}|^2}\right)^{p-2}\pmb\phi_{+2}^{(p)}\right)=\frac{\ov{\widetilde{q}}}{\rr\widetilde{q}}\left(\frac{\rr^2}{|\widetilde{q}|^2}\right)^{p-1}\pmb\phi_{+2}^{(p+1)}+\N_{T,+2}^{(p)}, \quad p=0,1,
\eea
and
\bea
\lab{def:TensorialTeuScalars:wavesystem:Kerrperturbation:-2:tildeform}
\nab_4\left(\frac{\rr\widetilde{q}}{\ov{\widetilde{q}}}\left(\frac{{\rr^2}}{|\widetilde{q}|^2}\right)^{p-2}\pmb\phi_{-2}^{(p)}\right)=\frac{\widetilde{q}}{\rr\ov{\widetilde{q}}}\left(\frac{\rr^2}{|\widetilde{q}|^2}\right)^{p-1}\frac{\De(\rr)}{|\widetilde{q}|^2}\pmb\phi_{-2}^{(p+1)}+\N_{T,-2}^{(p)}, \,\,\,\, p=0,1.
\eea
\esub

Also, in view of \eqref{eq:definitionofthephiminus2phierarchy:perturbationofKerr}, $\pmb\phi_{-2}^{(0)}$ degenerates at $r=r_+$ and does thus not allow to recover estimates for $\Ab$ near $r=r_+$. To remedy this problem, we will rely on the following lemma.
\begin{lemma}\lab{lemma:pmbphiminus2p=0isdegenerateatr=rplus}
The horizontal tensor $\Ab\in\sk_2(\mathbb{C})$ satisfies in the redshift region $r\leq r_+(1+2\dred)$:
\bea\lab{eq:waveequationpmbphip=0sminus2nodeginredshiftregion}
\squared_2\nab_4^p\Ab = (2-p)\pr_r\left(\frac{\De}{|q|^2}\right)\nab_3\nab_4^p\Ab
+ O(1)\big(\nab_{4}\nab_4^{\leq p}\Ab,\nab\nab_4^{\leq p}\Ab, \nab_4^{\leq p}\Ab\big)+\N_{\nab_4^p\Ab}, \quad p=0,1,2,
\eea
where we have in $r\leq r_+(1+2\dred)$
\bea
\N_{\Ab}=\dk^{\leq 1}(\Ga_g\c\Ga_b), \qquad \N_{\nab_4^p\Ab}=\nab_4^p\N_{\Ab}+\dk^{\leq p+1}(\Ga_g\c\Ab), \quad p=1,2.
\eea
\end{lemma}

\begin{proof}
Recall from Lemma 5.3.3 in \cite{GKS22} that
\bea\lab{eq:squared2AbformofTeukolskyfromLemma533inGKS22}
\squared_2 \Ab&=& \left(2\ov{\tr \Xb} +4\omb \right) \nab_4 \Ab-4\om \nab_3 \Ab -\left( 4\Hb-4\ze \right)\c \nab \Ab +V \Ab+ \err
 \eea
 where 
 \beaa
 V&=& \frac 3 4 \tr X \ov{\tr \Xb} -\frac 1 4 \tr \Xb  \ov{\tr X}  - 4P-4\om\left( \frac 1 2 \tr \Xb +2\ov{\tr \Xb} \right)+2 \omb \tr X-4\nab_3\om  -8\omb\om\\
&&+  \DD\c \ov{Z}- 2Z \c \ov{Z}+2 \ze\c \left( 4\Hb+2\eta \right)- 4 \etab \c \eta  +2i\etab \wedge \eta
 \eeaa
 and the error terms are given by
 \beaa
 \err&=& r^{-1}\dk^{\leq 1}\big(  \Ga_b \c \Bb \big) + \Ga_b \c \Ga_b \c \Ga_g.
 \eeaa
The case $p=0$ in \eqref{eq:waveequationpmbphip=0sminus2nodeginredshiftregion} is then a non-sharp consequence of \eqref{eq:squared2AbformofTeukolskyfromLemma533inGKS22}.

Next, the cases $p=1,2$ follow from the case $p=0$, commutation with $\nab_4^p$, the commutation formulas of Section \ref{sec:generalcommutationformulasrealcase}, as well as the following commutation formula for $\psi\in\sk_2$ in $r\leq r_+(1+2\dred )$:
\beaa
\begin{split}
[\nab_4, \squared_2]\psi =& \pr_r\left(\frac{\De}{|q|^2}\right)\nab_3\nab_4\psi+O(|r-r_+|)\big(\squared_2\psi, \nab_3\nab_4^{\leq 1}\psi\big)\\
&+O(1)(\nab_4\nab_4^{\leq 1}\psi, \nab\nab_4^{\leq 1}\psi, \nab_4^{\leq 1}\psi)+\dk^{\leq 2}(\Ga_g\c\psi),
\end{split}
\eeaa
which is proved in Lemma 6.18 of \cite{MaSz26}. This concludes the proof of Lemma \ref{lemma:pmbphiminus2p=0isdegenerateatr=rplus}.
\end{proof}

We now state Energy-Morawetz estimates for up to 14 derivatives for solutions to Teukolsky equations.

\begin{theorem}
\label{thm:main:MaSz26}
Let $(\MM, \g)$ satisfy the assumptions of Section \ref{sec:geometricsetupfordecsayTeukolsky}, and let $(\ttt, \rr, \tx^1, \tx^2)$ be the coordinates system on $\MM$, with $\ttt$, $\rr$ and $(\tx^1, \tx^2)$ introduced respectively in Lemma \ref{lemma:controloftildeuonMext}, and in \eqref{definitionofthenewcoordinatetilder} and 
\eqref{definitionofthenewcoordinatetildethetandtphi}. Then, for $\ep>0$ small enough, we have for solutions $\pmb\phi_s^{(p)}$, $s=\pm 2$, $p=0,1,2$, to the tensorial Teukolsky wave/transport systems \eqref{eq:TensorialTeuSysandlinearterms:rescaleRHScontaine2:general:Kerrperturbation:tildeform}  \eqref{def:TensorialTeuScalars:wavesystem:Kerrperturbation:tildeform} on $\MM$ the following energy-Morawetz-flux estimates, for any $1\leq\tau_1<\tau_2\leq\tau_*$ and any $0<\de\leq \frac{1}{3}$, 
\bea
\label{MainEnerMora:psi:plus2case}
\sum_{p=0}^2\EMF^{{11}}_{\de}[\pmb\phi_{+2}^{(p)}](\tau_1, \tau_2)
&\les&  \sum_{p=0}^2\E^{{11}}[\pmb\phi_{+2}^{(p)}](\tau_1)
+\sum_{p=0}^2\widetilde{\NN}^{{11}}_\de[\pmb\phi_{+2}^{(p)}, \N_{W,+2}^{(p)}](\tt_1, \tt_2)\nn\\ 
&&
+\sum_{p=0}^1\int_{\MM(\tt_1, \tt_2)}r^{-1+\de}|\dk^{\leq 12}\N_{T,+2}^{(p)}|^2
\eea
and, assuming also that $\Ab$ satisfies \eqref{eq:waveequationpmbphip=0sminus2nodeginredshiftregion},
\bea
\label{MainEnerMora:psi:minus2case}
\nn&&\sum_{p=0}^2\EMF^{{14}}_\de[\pmb\phi_{-2}^{(p)}](\tau_1, \tau_2) + \sum_{p=0}^2\EMF^{{14}}_{r\leq r_+(1+\dred)}[\nab_4^p\Ab](\tau_1, \tau_2)\\
\nn&\les& \sum_{p=0}^2\E^{{14}}[\pmb\phi_{-2}^{(p)}](\tau_1)+\sum_{p=0}^2\E^{{14}}_{r\leq r_+(1+\dred)}[\nab_4^p\Ab](\tau_1) +\int_{\Si(\tau_1)}|\dk^{\leq 12}\N_{T,-2}^{(0)}|^2\\
\nn&&+\sum_{p=0}^1\int_{\MM(\tt_1, \tt_2)}r^{-1+\de}|\dk^{\leq {15}}\N_{T,-2}^{(p)}|^2+\sum_{p=0}^2\widetilde{\NN}_\de^{{14}}[\pmb\phi_{-2}^{(p)}, \N_{W,-2}^{(p)}](\tt_1, \tt_2)\\ 
&&+\sum_{p=0}^2\int_{\MM_{r\leq r_+(1+2\dred)}(\tau_1, \tau_2)}|\dk^{\leq {14}}\N_{\nab_4^p\Ab}|^2,
\eea
where the norms $\EMF^{\reg}_\de[\c](\tau_1,\tau_2)$, $\E^{\reg}[\c](\tau_1)$ and $\widetilde{\mathcal{N}}^{\reg}_\de[\c, \c](\tau_1, \tau_2)$ have been introduced in Section \ref{subsection:basicnormsforpsi}, where $\N_{W,s}^{(p)}$, $\N_{T,s}^{(p)}$ and $\N_{\nab_4^p\Ab}$ are introduced in equations \eqref{eq:TensorialTeuSysandlinearterms:rescaleRHScontaine2:general:Kerrperturbation:tildeform},   \eqref{def:TensorialTeuScalars:wavesystem:Kerrperturbation:tildeform} and \eqref{eq:waveequationpmbphip=0sminus2nodeginredshiftregion}, and where the implicit constant in $\lesssim$ only depends on $a$, $m$, $\de$ and $\dec$.
\end{theorem}

\begin{proof}
Recall Definition \ref{def:renormalizationofallnonsmallquantitiesinPGstructurebyKerrvalue:widetildecase} on linearized quantities with Kerr values computed w.r.t. the coordinates $(\rr, \tx^1, \tx^2)$ and the complex 1-form $(\tJk, \tJk_{\pm})$, as well as Definition \ref{definition.Ga_gGa_b:widetildecase} for the quantities $(\widetilde{\Ga}_g, \widetilde{\Ga}_b)$, where $(\tJk, \tJk_{\pm})$ are given by Lemma \ref{lemma:controlofthecomplex1formswidetildeJkandwidetildeJkplusminus}. Then, given that: 
\begin{itemize}
\item the complex 1-forms $(\Jk, \Jk_{\pm})$ satisfy the identities \eqref{eq:usefulalgebraicidentitiesinvolvingscalarproductsReJkReJkpm:Kerrpert:tilde},

\item $(\widetilde{\Ga}_{g}, \widetilde{\Ga}_{b})$ satisfy \eqref{eq:decaypropertiesofwidehatGab1widehatGag1} \eqref{eq:additionestimatewidechecke4tttbyeprminus2:withdoublewidecheck}, and the linearized components of the inverse metric in the coordinates system $(\ttt, \rr, \tx^1, \tx^2)$ satisfy \eqref{eq:controloflinearizedinversemetriccoefficients} \eqref{eq:controloflinearizedinversemetriccoefficients:inversegtautau},

\item the regular triplet $\Om_i$, $i=1,2,3$, on $\MM$ constructed in Lemma \ref{lemma:constructionofhorizontal1formsOmiinperturbationsofKerr} satisfies \eqref{eq:assumptionsonregulartripletinperturbationsofKerr},
\end{itemize} 
the assumptions (4.27), (4.30), (4.34), (4.35) and (4.42) in Section 4.2.3 of \cite{Sze} are satisfied and we are thus in position to apply Theorem 4.52 in \cite{Sze} which immediately yields \eqref{MainEnerMora:psi:plus2case} and \eqref{MainEnerMora:psi:minus2case} as stated.
\end{proof}

\begin{remark}
Theorem 4.52 in \cite{Sze}, used in the above proof of Theorem \ref{thm:main:MaSz26}, is the extension of Theorem 7.1 in \cite{MaSz26} to the case where $\MM$ only extends to the spacelike hypersuface\footnote{Note that Theorem 7.1 in \cite{MaSz26} is proved in the case where the spacetime $\MM$ extends to $\II_+$.} $\Si_*$. 
\end{remark}

%%%%%%%%%%%%%%%%%%%%%%%%%%%%%%%%%%%%%%%%%%%%%%%%%%%%%%

\subsubsection{Energy-Morawetz estimates for Teukolsky equations in $\MM$ up to $\kst-3$ derivatives}
\lab{sec:energyMorawetzesitmatesforTeukoslkyonMM:higherorderderivatives}

%%%%%%%%%%%%%%%%%%%%%%%%%%%%%%%%%%%%%%%%%%%%%%%%%%%%%%%

The goal of this section is to extend the Energy-Morawetz estimates for up to 14 derivatives for solutions to Teukolsky equations in Theorem 
\ref{thm:main:MaSz26} to an arbitrary number of derivatives. This is done in the following theorem. 
\begin{theorem}
\label{thm:main:MaSz26:extendhigherorderderivatives}
Let $(\MM, \g)$ satisfy the assumptions of Section \ref{sec:geometricsetupfordecsayTeukolsky}. Then, for $\ep>0$ small enough, we have for solutions $\pmb\phi_s^{(p)}$, $s=\pm 2$, $p=0,1,2$, to the Teukolsky wave/transport system \eqref{eq:TensorialTeuSysandlinearterms:rescaleRHScontaine2:general:Kerrperturbation:alternateformnullframeinsteadcoordvectorfield} \eqref{eq:transportequationsins=plus2andminus2caseforp=0and1} the following energy-Morawetz-flux estimates, for any $1\leq\tau_1<\tau_2\leq\tau_*$ and any $11\leq\reg\leq \kst-3$
\bea
\label{MainEnerMora:psi:plus2case:extendhigherorderderivatives}
&&\sum_{p=0}^2\EMF_\de^{\reg}[\pmb\phi_{+2}^{(p)}](\tau_1, \tau_2) \nn\\
&\les&  \sum_{p=0}^2\E^{\reg}[\pmb\phi_{+2}^{(p)}](\tau_1)+\sum_{p=0}^2\widetilde{\NN}^{\reg}_\de[\pmb\phi_{+2}^{(p)}, \widetilde{\N}_{W,+2}^{(p)}](\tt_1, \tt_2)+ \sum_{p=0}^1\int_{\MM(\tt_1, \tt_2)}r^{-1+\de}|\dk^{\leq\reg+1}\N_{T,+2}^{(p)}|^2\nn\\
&& +\sum_{p=0}^1\int_{\MM(\tt_1, \tt_2)}r^{1+\de}|\dk^{\leq\reg+1}\big(\Ga_b\c\pmb\phi_{+2}^{(p)}\big)|^2+\ep\sum_{p=0}^2\B^{\reg}_\de[\pmb\phi_{+2}^{(p)}](\tau_1, \tau_2),
\eea
and, for $14\leq\reg\leq\kst -3$, assuming also that $\Ab$ satisfies \eqref{eq:waveequationpmbphip=0sminus2nodeginredshiftregion},
\bea
\label{MainEnerMora:psi:minus2case:extendhigherorderderivatives}
\nn&&\sum_{p=0}^2\Big(\EMF_\de^{\reg}[\pmb\phi_{-2}^{(p)}](\tau_1, \tau_2)+\EMF^{\reg}_{r\leq r_+(1+\dred)}[\nab_4^p\Ab](\tau_1, \tau_2)\Big)\\
\nn&\les& \sum_{p=0}^2\E^{\reg}[\pmb\phi_{-2}^{(p)}](\tau_1)+\sum_{p=0}^2\E^{\reg}_{r\leq r_+(1+\dred)}[\nab_4^p\Ab](\tau_1)+\int_{\Si(\tau_1)}|\dk^{\leq\reg-2}\N_{T,-2}^{(0)}|^2\nn\\
&&+\sum_{p=0}^1\int_{\MM(\tt_1, \tt_2)}r^{-1+\de}|\dk^{\leq\reg+1}\N_{T,-2}^{(p)}|^2+\sum_{p=0}^2\widetilde{\NN}_\de^{\reg}[\pmb\phi_{-2}^{(p)}, \widetilde{\N}_{W,-2}^{(p)}](\tt_1, \tt_2)\nn\\
&& +\ep\sum_{p=0}^2\B^{\reg}_\de[\pmb\phi_{-2}^{(p)}](\tau_1, \tau_2)+ \frac{\ep^2_0}{\tau_1^{3+3\dec}},
\eea
where the norms $\EMF^{\reg}_\de[\c](\tau_1,\tau_2)$, $\E^{\reg}[\c](\tau_1)$, $\B^{\reg}_\de[\c](\tau_1,\tau_2)$ and $\widetilde{\mathcal{N}}^{\reg}_\de[\c, \c](\tau_1, \tau_2)$ have been introduced in Section \ref{subsection:basicnormsforpsi}, where $\widetilde{\N}_{W,s}^{(p)}$, $s=\pm2$, $p=0,1,2$, are introduced in \eqref{eq:TensorialTeuSysandlinearterms:rescaleRHScontaine2:general:Kerrperturbation:alternateformnullframeinsteadcoordvectorfield} and $\N_{T,s}^{(p)}$, $s=\pm2$, $p=0,1$, are introduced in \eqref{eq:transportequationsins=plus2andminus2caseforp=0and1}, and where the implicit constant in $\lesssim$ only depends on $a$, $m$, $\de$, $\Nmic$ and $\dec$.
\end{theorem}

\begin{proof}
Extending basic energy-Morawetz estimates, such as the ones for up to 14 in Theorem 
\ref{thm:main:MaSz26}, to an arbitrary number of derivatives is by now standard. The proof below, in the case of Teukolsky in perturbations of Kerr, is inspired by Section 3.6 in \cite{MaSz24}, and Sections 10.1 and 10.2 in \cite{MaSz26}. It proceeds in several steps.

\noindent{\bf Step 1.} Let $\reg_s$, $s=\pm 2$ be integers given by 
\bea\lab{eq:choiceofrregplus2andminus2}
\reg_{+2}:=11, \qquad \reg_{-2}:=14.
\eea  
For any $\reg_s\leq\reg\leq \kst-3$, we commute \eqref{eq:TensorialTeuSys:rescaleRHScontaine2:general:Kerrperturbation:alternateformnullframeinsteadcoordvectorfield} and 
\eqref{eq:transportequationsins=plus2andminus2caseforp=0and1} with $\Lieb_{\T}^{\reg-\reg_s}$, where the horizontal Lie derivative $\Lieb$ has been introduced in Definition \ref{definition:hor-Lie-derivative}, which yields
\bea\lab{eq:TeukwavetransportsystemcommuttedwithLiebregminusks}
\bigg(\squared_2 -\frac{4ia\cos\th}{|q|^2}\nab_{\T}- \frac{4-2\de_{p0}}{|q|^2}\bigg)\Lieb_{\T}^{\reg-\reg_s}\phis{p} &=& \widetilde{\L}_s^{(p)}[\Lieb_{\T}^{\reg-\reg_s}\pmb\phi_{s}]+\Lieb_{\T}^{\reg-\reg_s}\widetilde{\N}_{W,s}^{(p)}\nn\\
&&+\err_{W,s, \Lieb_\T}^{(p),\reg-\reg_s},\,\,\, s=\pm 2, \,\,\, p=0,1,2,
\eea
and
\bsub\lab{eq:TeuktransportsystemforLiebregminusregspmbphisp}
\bea
\nab_3\left(\frac{r\ov{q}}{q}\left(\frac{r^2}{|q|^2}\right)^{p-2}\Lieb_{\T}^{\reg-\reg_s}\pmb\phi_{+2}^{(p)}\right) &=& \frac{\ov{q}}{rq}\left(\frac{r^2}{|q|^2}\right)^{p-1}\Lieb_{\T}^{\reg-\reg_s}\pmb\phi_{+2}^{(p+1)}+\Lieb_{\T}^{\reg-\reg_s}\N_{T,+2}^{(p)}\nn\\
&&+\err_{T,+2, \Lieb_\T}^{(p),\reg-\reg_{+2}}, \quad p=0,1,
\eea
\bea
\nab_4\left(\frac{rq}{\ov{q}}\left(\frac{r^2}{|q|^2}\right)^{p-2}\Lieb_{\T}^{\reg-\reg_s}\pmb\phi_{-2}^{(p)}\right) &=& \frac{q}{r\ov{q}}\left(\frac{r^2}{|q|^2}\right)^{p-1}\frac{\De}{|q|^2}\Lieb_{\T}^{\reg-\reg_s}\pmb\phi_{-2}^{(p+1)}+\Lieb_{\T}^{\reg-\reg_s}\N_{T,-2}^{(p)}\nn\\
&&+\err_{T,-2, \Lieb_\T}^{(p),\reg-\reg_{-2}}, \quad p=0,1,
\eea
\esub
where $\err_{W,s, \Lieb_\T}^{(p),\reg-\reg_s}$, $s=\pm 2$, $p=0,1,2$, and $\err_{T,s, \Lieb_\T}^{(p),\reg-\reg_{s}}$, $s=\pm 2$, $p=0,1$, are given by
\beaa
\err_{W,s, \Lieb_\T}^{(p),\reg-\reg_s} &=& [\squared_2, \Lieb_{\T}^{\reg-\reg_s}]\phis{p}+\left[\Lieb_{\T}^{\reg-\reg_s}, \frac{4ia\cos\th}{|q|^2}\nab_{\T}\right]\pmb\phi_{-2}^{(p)}+\left[\Lieb_{\T}^{\reg-\reg_s}, \frac{4-2\de_{p0}}{|q|^2}\right]\pmb\phi_{-2}^{(p)}\\
&&+[\Lieb_{\T}^{\reg-\reg_s}, \widetilde{\L}_s^{(p)}]\pmb\phi_{-2}^{(p)}, 
\eeaa
\beaa
\err_{T,+2, \Lieb_\T}^{(p),\reg-\reg_{+2}} = \left[\nab_3\left(\frac{r\ov{q}}{q}\left(\frac{r^2}{|q|^2}\right)^{p-2}\cdot\right), \Lieb_{\T}^{\reg-\reg_s}\right]\pmb\phi_{+2}^{(p)}+\left[\Lieb_{\T}^{\reg-\reg_s}, \frac{\ov{q}}{rq}\left(\frac{r^2}{|q|^2}\right)^{p-1}\right]\pmb\phi_{+2}^{(p+1)}, 
\eeaa
and
\beaa
\err_{T,-2, \Lieb_\T}^{(p),\reg-\reg_{-2}} = \left[\nab_4\left(\frac{rq}{\ov{q}}\left(\frac{r^2}{|q|^2}\right)^{p-2}\c\right), \Lieb_{\T}^{\reg-\reg_s}\right]\pmb\phi_{-2}^{(p)}+\left[\Lieb_{\T}^{\reg-\reg_s}, \frac{q}{r\ov{q}}\left(\frac{r^2}{|q|^2}\right)^{p-1}\frac{\De}{|q|^2}\right]\pmb\phi_{-2}^{(p+1)}.
\eeaa
Together with Proposition \ref{PROPOSITION:COMMUTATORS-SQUARED-LIED-T-Z} and Lemma \ref{lemma:commutatorbetweenLieTLieZandnabnab4nab3}, we infer
\bea\lab{eq:controloftheerrortermserrWsLiebTregminusregsanderrTsLiebTregminusregs}
\bsplit
\err_{W,s, \Lieb_\T}^{(p),\reg-\reg_s} =& \dk^{\leq\reg-\reg_s}\big(\Ga_g \c \dk\phis{p}\big)+r^{-3}\dk^{\leq\reg-\reg_s}(\Ga_b\c\phis{p+1})\de_{p\leq 1}+r^{-1}\dk^{\leq\reg-\reg_s}(\Ga_b\c\phis{p-1})\de_{p\geq 1}\\
&+r^{-2}\dk^{\leq\reg-\reg_s}(\Ga_b\c\phis{0})\de_{p=2}, \quad s=\pm 2, \quad p=0,1,2,\\
\err_{T,+2, \Lieb_\T}^{(p),\reg-\reg_{+2}} =& r\dk^{\leq\reg-\reg_s}\big(\Ga_b\c\pmb\phi_{+2}^{(p)}\big)+r^{-1}\dk^{\leq\reg-\reg_s}\big(\Ga_b\c\pmb\phi_{+2}^{(p+1)}\big), \quad p=0,1,\\
\err_{T,-2, \Lieb_\T}^{(p),\reg-\reg_{-2}} =& r\dk^{\leq\reg-\reg_s}\big(\Ga_g\c\pmb\phi_{-2}^{(p)}\big)+r^{-1}\dk^{\leq\reg-\reg_s}\big(\Ga_b\c\pmb\phi_{-2}^{(p+1)}\big), \quad p=0,1.
\end{split}
\eea

Next, we consider the coordinates $(\ttt, \rr, \tx^1, \tx^2)$ with $\ttt$, $\rr$ and $(\tx^1, \tx^2)$ introduced respectively in Lemma \ref{lemma:controloftildeuonMext}, and in \eqref{definitionofthenewcoordinatetilder} and 
\eqref{definitionofthenewcoordinatetildethetandtphi}. Also, we consider the complex 1-forms $(\tJk, \tJk_{\pm})$ introduced in \eqref{definitionofthenewcomplex1formstJktJkpm}. Then, \eqref{eq:propoertiesoflinearizedfirstorderderivativesoftphiintermsofGabGag} and 
\eqref{eq:identiteReJkcReJkandReJkcReJkpmuptoGagtermswhichisenough} hold in view of Definition \ref{definition.Ga_gGa_b:widetildecase}, \eqref{eq:decaypropertiesofwidehatGab1widehatGag1} and \eqref{eq:usefulalgebraicidentitiesinvolvingscalarproductsReJkReJkpm:Kerrpert:tilde}. We may thus apply Corollary \ref{cor:derivationoftheTeukolskytensorialwavesystemfors=plusminus2:kerrpert} to \eqref{eq:TeukwavetransportsystemcommuttedwithLiebregminusks} which implies that $\Lieb_{\T}^{\reg-\reg_s}\pmb\phi_s^{(p)}$, $s=\pm 2$, $p=0,1,2$ satisfies the following wave system, for $\reg_s\leq\reg\leq \kst-3$, 
\bea\lab{eq:TeukwavetransportsystemcommuttedwithLiebregminusks:gobacktoMaSz26form}
\bigg(\squared_2 -\frac{4ia\cos\th}{|q|^2}\nab_{\pr_{\ttt}}- \frac{4-2\de_{p0}}{|q|^2}\bigg)\Lieb_{\T}^{\reg-\reg_s}\phis{p} = \L_s^{(p)}[\Lieb_{\T}^{\reg-\reg_s}\pmb\phi_{s}]+\N_{W,s, \Lieb_\T}^{(p),\reg-\reg_s},
\eea
where $\L_{s}^{(p)}[\pmb\phi_{s}]$ is given by \eqref{eq:tensor:Lsn:onlye_2present:general:Kerrperturbation:bis:tildeform}, i.e., w.r.t. the coordinates $(\ttt, \rr, \tx^1, \tx^2)$, and where, in view of \eqref{eq:schematicformofNpWs:comparisionNWandwidetildeWW}, $\N_{W,s, \Lieb_\T}^{(p),\reg-\reg_s}$, $s=\pm 2$, $p=0,1,2$ satisfies
\bea\lab{eq:structureofRHSNWspTeukwavessytemLiebTregminusregspmbphisspinformofMaSz26}
\bsplit
\N_{W,s, \Lieb_\T}^{(0),\reg-\reg_s} =& \Lieb_{\T}^{\reg-\reg_s}\widetilde{\N}_{W,s}^{(0)}+\err_{W,s, \Lieb_\T}^{(0),\reg-\reg_s} + r^{-1}\widetilde{\Ga}_b\c\dk\Lieb_{\T}^{\reg-\reg_s}\pmb\phi_s^{(0)}\\
&+r^{-3}\Big(O(\rr-r)+O(\tx^b-\tx^b)\Big)\Lieb_{\T}^{\reg-\reg_s}\big(\pmb\phi_s^{(0)}, \pmb\phi_s^{(1)}\big),\\
\N_{W,s, \Lieb_\T}^{(1),\reg-\reg_s} =& \Lieb_{\T}^{\reg-\reg_s}\widetilde{\N}_{W,s}^{(1)}+\err_{W,s, \Lieb_\T}^{(1),\reg-\reg_s} + r^{-1}\widetilde{\Ga}_b\c\dk\Lieb_{\T}^{\reg-\reg_s}\big(\pmb\phi_s^{(0)}, \pmb\phi_s^{(1)}\big)\\
&+r^{-3}\Big(O(\rr-r)+O(\tx^b-\tx^b)\Big)\Lieb_{\T}^{\reg-\reg_s}\big(r\dk^{\leq 1}\pmb\phi_s^{(0)}, \dk^{\leq 1}\pmb\phi_s^{(1)}, \pmb\phi_s^{(2)}\big),\\
\N_{W,s, \Lieb_\T}^{(2),\reg-\reg_s} =& \Lieb_{\T}^{\reg-\reg_s}\widetilde{\N}_{W,s}^{(2)}+\err_{W,s, \Lieb_\T}^{(2),\reg-\reg_s} + r^{-1}\widetilde{\Ga}_b\c\dk\Lieb_{\T}^{\reg-\reg_s}\pmb\phi_s^{(2)}+\widetilde{\Ga}_g\c\dk\Lieb_{\T}^{\reg-\reg_s}\pmb\phi_s^{(1)},\\
&+r^{-3}\Big(O(\rr-r)+O(\tx^b-\tx^b)\Big)\Lieb_{\T}^{\reg-\reg_s}\big(r\pmb\phi_s^{(0)}, r\dk^{\leq 1}\pmb\phi_s^{(1)}, \pmb\phi_s^{(2)}\big).
\end{split}
\eea
Also, we have, in view of \eqref{eq:TeuktransportsystemforLiebregminusregspmbphisp},
\bsub\lab{eq:TeuktransportsystemforLiebregminusregspmbphisp:modifiedform}
\bea
\nab_3\left(\frac{\rr\ov{\tq}}{\tq}\left(\frac{\rr^2}{|\tq|^2}\right)^{p-2}\Lieb_{\T}^{\reg-\reg_{+2}}\pmb\phi_{+2}^{(p)}\right) = \frac{\ov{\tq}}{\rr\tq}\left(\frac{\rr^2}{|\tq|^2}\right)^{p-1}\Lieb_{\T}^{\reg-\reg_{+2}}\pmb\phi_{+2}^{(p+1)}+\N_{T,+2, \Lieb_\T}^{(p),\reg-\reg_{+2}}, \quad p=0,1,
\eea
\bea
\nab_4\left(\frac{\rr\tq}{\ov{\tq}}\left(\frac{\rr^2}{|\tq|^2}\right)^{p-2}\Lieb_{\T}^{\reg-\reg_{-2}}\pmb\phi_{-2}^{(p)}\right) &=& \frac{\tq}{\rr\ov{\tq}}\left(\frac{\rr^2}{|\tq|^2}\right)^{p-1}\frac{\De(\rr)}{|\tq|^2}\Lieb_{\T}^{\reg-\reg_{-2}}\pmb\phi_{-2}^{(p+1)}+\N_{T,-2, \Lieb_\T}^{(p),\reg-\reg_{-2}}\nn\\
&&+\Lieb_{\T}^{\reg-\reg_{-2}}\N_{T,-2}^{(p)}+\err_{T,-2, \Lieb_\T}^{(p),\reg-\reg_{-2}}, \quad p=0,1,
\eea
\esub
where, $\N_{T,s, \Lieb_\T}^{(p),\reg-\reg_s}$, $s=\pm 2$, $p=0,1$, satisfies
\bsub\lab{eq:structureofRHSNTspTeuktransportsytemLiebTregminusregspmbphisspinformofMaSz26}
\bea
\N_{T,+2, \Lieb_\T}^{(p),\reg-\reg_{+2}} &=& \Lieb_{\T}^{\reg-\reg_{+2}}\N_{T,+2}^{(p)}+\err_{T,+2, \Lieb_\T}^{(p),\reg-\reg_{+2}}\nn\\
&& +\Big(O(\rr-r)+O(\tx^b-\tx^b)\Big)\Lieb_{\T}^{\reg-\reg_{+2}}\big(\dk^{\leq 1}\pmb\phi_{+2}^{(p)}, r^{-2}\pmb\phi_{+2}^{(p+1)}\big),\\
\N_{T,-2, \Lieb_\T}^{(p),\reg-\reg_{-2}} &=& \Lieb_{\T}^{\reg-\reg_{-2}}\N_{T,-2}^{(p)}+\err_{T,-2, \Lieb_\T}^{(p),\reg-\reg_{-2}}\nn\\
&& +\Big(O(\rr-r)+O(\tx^b-\tx^b)\Big)\Lieb_{\T}^{\reg-\reg_{-2}}\big(r^{-1}\dk^{\leq 1}\pmb\phi_{-2}^{(p)}, r^{-2}\pmb\phi_{-2}^{(p+1)}\big).
\eea
\esub
Moreover, commuting \eqref{eq:waveequationpmbphip=0sminus2nodeginredshiftregion} with $\Lieb_{\T}^{\reg-\reg_{-2}}$, for $\reg_s\leq\reg\leq \kst-3$, we immediately obtain in the redshift region $r\leq r_+(1+2\dred)$:
\bea\lab{eq:waveequationpmbphip=0sminus2nodeginredshiftregion:commutationLiebTregminusregs}
\nn\squared_2\nab_4^p\Lieb_{\T}^{\reg-\reg_{-2}}\Ab &=& (2-p)\pr_r\left(\frac{\De}{|q|^2}\right)\nab_3\nab_4^p\Lieb_{\T}^{\reg-\reg_{-2}}\Ab\\
&&+ O(1)\big(\nab_{4}\nab_4^{\leq p}\Lieb_{\T}^{\reg-\reg_{-2}}\Ab,\nab\nab_4^{\leq p}\Lieb_{\T}^{\reg-\reg_{-2}}\Ab, \nab_4^{\leq p}\Lieb_{\T}^{\reg-\reg_{-2}}\Ab\big)\nn\\
&&+\N_{\nab_4^p\Lieb_{\T}^{\reg-\reg_{-2}}\Ab}, \qquad p=0,1,2,
\eea
where we have in $r\leq r_+(1+2\dred)$
\bea\lab{eq:strctureofNnab4pLiebTregminusregminus2ofAbinredshiftregion}
\N_{\nab_4^p\Lieb_{\T}^{\reg-\reg_{-2}}\Ab}=\dk^{\leq \reg-\reg_{-2}+p+1}(\Ga_g\c\Ga_b), \quad p=0,1,2.
\eea

In view of \eqref{eq:TeukwavetransportsystemcommuttedwithLiebregminusks:gobacktoMaSz26form} \eqref{eq:TeuktransportsystemforLiebregminusregspmbphisp:modifiedform}, $\Lieb_{\T}^{\reg-\reg_s}\phis{p}$, $s=\pm 2$, $p=0,1,2$, satisfies the Teukolsky wave/transport systems \eqref{eq:TensorialTeuSysandlinearterms:rescaleRHScontaine2:general:Kerrperturbation:tildeform}  \eqref{def:TensorialTeuScalars:wavesystem:Kerrperturbation:tildeform} with RHS $(\N_{W,s}^{(p)}, \N_{T,s}^{(p)})\to (\N_{W,s, \Lieb_\T}^{(p),\reg-\reg_s}, \N_{T,s, \Lieb_\T}^{(p),\reg-\reg_s})$. Also, $\nab_4^p\Lieb_{\T}^{\reg-\reg_{-2}}\Ab$ satisfies \eqref{eq:waveequationpmbphip=0sminus2nodeginredshiftregion} with RHS $\N_{\nab_4^p\Ab}\to\N_{\nab_4^p\Lieb_{\T}^{\reg-\reg_{-2}}\Ab}$. We may thus apply Theorem \ref{thm:main:MaSz26} 
which implies, taking also \eqref{eq:choiceofrregplus2andminus2} into account, for $11\leq\reg\leq\kst -3$, 
\bea\lab{MainEnerMora:psi:plus2case:commutationLiebTregminusregplus2case}
\sum_{p=0}^2\EMF^{({11})}_{\de}[\Lieb_{\T}^{\reg-11}\pmb\phi_{+2}^{(p)}](\tau_1, \tau_2)
&\les&  \sum_{p=0}^2\E^{\reg}[\pmb\phi_{+2}^{(p)}](\tau_1)
+\sum_{p=0}^2\widetilde{\NN}^{({11})}_\de[\Lieb_{\T}^{\reg-11}\pmb\phi_{+2}^{(p)}, \N_{W,+2, \Lieb_\T}^{(p),\reg-11}](\tt_1, \tt_2)\nn\\ 
&&
+\sum_{p=0}^1\int_{\MM(\tt_1, \tt_2)}r^{-1+\de}|\dk^{\leq 12}\N_{T,+2, \Lieb_\T}^{(p),\reg-11}|^2
\eea
and, for $14\leq\reg\leq\kst -3$, 
\bea\lab{MainEnerMora:psi:minus2case:commutationLiebTregminusregplus2case}
\nn&&\sum_{p=0}^2\EMF^{({14})}_\de[\Lieb_{\T}^{\reg-14}\pmb\phi_{-2}^{(p)}](\tau_1, \tau_2) + \sum_{p=0}^2\EMF^{({14})}_{r\leq r_+(1+\dred)}[\nab_4^p\Lieb_{\T}^{\reg-14}\Ab](\tau_1, \tau_2)\\
\nn&\les& \sum_{p=0}^2\E^{\reg}[\pmb\phi_{-2}^{(p)}](\tau_1)+\sum_{p=0}^2\E^{\reg}_{r\leq r_+(1+\dred)}[\nab_4^p\Ab](\tau_1) +\int_{\Si(\tau_1)}|\dk^{\leq 12}\N_{T,-2, \Lieb_\T}^{(0),\reg-14}|^2\\
\nn&&+\sum_{p=0}^1\int_{\MM(\tt_1, \tt_2)}r^{-1+\de}|\dk^{\leq {15}}\N_{T,-2, \Lieb_\T}^{(p),\reg-14}|^2+\sum_{p=0}^2\widetilde{\NN}_\de^{({14})}[\Lieb_{\T}^{\reg-14}\pmb\phi_{-2}^{(p)}, \N_{W,-2, \Lieb_\T}^{(p),\reg-14}](\tt_1, \tt_2)\\ 
&&+\sum_{p=0}^2\int_{\MM_{r\leq r_+(1+2\dred)}(\tau_1, \tau_2)}|\dk^{\leq {14}}\N_{\nab_4^p\Lieb_{\T}^{\reg-14}\Ab}|^2.
\eea

Next, we estimate the terms involving $\N_{W,s, \Lieb_\T}^{(p),\reg-\reg_s}$, $\N_{T,s, \Lieb_\T}^{(p),\reg-\reg_s}$ and $\N_{\nab_4^p\Lieb_{\T}^{\reg-\reg_{-2}}\Ab}$ on the RHS of \eqref{MainEnerMora:psi:plus2case:commutationLiebTregminusregplus2case} \eqref{MainEnerMora:psi:minus2case:commutationLiebTregminusregplus2case}. In view of \eqref{eq:structureofRHSNWspTeukwavessytemLiebTregminusregspmbphisspinformofMaSz26}, \eqref{eq:controloftheerrortermserrWsLiebTregminusregsanderrTsLiebTregminusregs} and Lemma \ref{lemma:controloftheerrortermsGagpsiinwaveeq}, we have
\bea\lab{eq:controlwidetildeNN11deLiebregminus11pmbphiplus2NWplus2LiebTpregminus11}
&&\sum_{p=0}^2\widetilde{\NN}^{({11})}_\de[\Lieb_{\T}^{\reg-11}\pmb\phi_{+2}^{(p)}, \N_{W,+2, \Lieb_\T}^{(p),\reg-11}](\tt_1, \tt_2)\nn\\
&\les& \sum_{p=0}^2\widetilde{\NN}^{({11})}_\de[\Lieb_{\T}^{\reg-11}\pmb\phi_{+2}^{(p)},  \Lieb_{\T}^{\reg-11}\widetilde{\N}_{W,+2}^{(p)}](\tt_1, \tt_2)+\ep\sum_{p=0}^2\bigg(\B^{\reg}_\de[\pmb\phi_{+2}^{(p)}](\tau_1, \tau_2)+\sup_{\tau\in[\tau_1, \tau_2]}\E^{\reg}[\pmb\phi_{+2}^{(p)}](\tau)\bigg)\nn\\
&\les& \sum_{p=0}^2\widetilde{\NN}^{\reg}_\de[\pmb\phi_{+2}^{(p)}, \widetilde{\N}_{W,+2}^{(p)}](\tt_1, \tt_2)+\ep\sum_{p=0}^2\bigg(\B^{\reg}_\de[\pmb\phi_{+2}^{(p)}](\tau_1, \tau_2)+\sup_{\tau\in[\tau_1, \tau_2]}\E^{\reg}[\pmb\phi_{+2}^{(p)}](\tau)\bigg)
\eea
and
\bea\lab{eq:controlwidetildeNN14deLiebregminus14pmbphiminus2NWminus2LiebTpregminus14}
&&\sum_{p=0}^2\widetilde{\NN}_\de^{({14})}[\Lieb_{\T}^{\reg-14}\pmb\phi_{-2}^{(p)}, \N_{W,-2, \Lieb_\T}^{(p),\reg-14}](\tt_1, \tt_2)\nn\\
&\les& \sum_{p=0}^2\widetilde{\NN}_\de^{({14})}[\Lieb_{\T}^{\reg-14}\pmb\phi_{-2}^{(p)},  \Lieb_{\T}^{\reg-14}\widetilde{\N}_{W,-2}^{(p)}](\tt_1, \tt_2)+\ep\sum_{p=0}^2\bigg(\B^{\reg}_\de[\pmb\phi_{-2}^{(p)}](\tau_1, \tau_2)+\sup_{\tau\in[\tau_1, \tau_2]}\E^{\reg}[\pmb\phi_{-2}^{(p)}](\tau)\bigg)\nn\\
&\les& \sum_{p=0}^2\widetilde{\NN}_\de^{\reg}[\pmb\phi_{-2}^{(p)}, \widetilde{\N}_{W,-2}^{(p)}](\tt_1, \tt_2)+\ep\sum_{p=0}^2\bigg(\B^{\reg}_\de[\pmb\phi_{-2}^{(p)}](\tau_1, \tau_2)+\sup_{\tau\in[\tau_1, \tau_2]}\E^{\reg}[\pmb\phi_{-2}^{(p)}](\tau)\bigg),
\eea
where we have also used \eqref{eq:interpolationestimatesforGabGagfork*derivativeswithk*largerthanklover2:forproofTheoremM1andM2}, \eqref{eq:decaypropertiesofwidehatGab1widehatGag1} and the control of $\rr-r$ and $\tx^b-x^b$ provided respectively by \eqref{eq:esitmatefor16weightedderivativesofrrminusr} and \eqref{eq:esitmatefor16weightedderivativesoftxbminusxb}. Also, in view of \eqref{eq:structureofRHSNTspTeuktransportsytemLiebTregminusregspmbphisspinformofMaSz26}, \eqref{eq:controloftheerrortermserrWsLiebTregminusregsanderrTsLiebTregminusregs} and Lemma \ref{lemma:controloftheerrortermsGagpsiinwaveeq}, we have
\begin{align}\lab{eq:controlwidetildeNN11deLiebregminus11pmbphiplus2NWplus2LiebTpregminus11:NTcase}
&\sum_{p=0}^1\int_{\MM(\tt_1, \tt_2)}r^{-1+\de}|\dk^{\leq 12}\N_{T,+2, \Lieb_\T}^{(p),\reg-11}|^2\nn\\
\les& \sum_{p=0}^1\int_{\MM(\tt_1, \tt_2)}r^{-1+\de}|\dk^{\leq 12}\Lieb_{\T}^{\reg-11}\N_{T,+2}^{(p)}|^2 +\sum_{p=0}^1\int_{\MM(\tt_1, \tt_2)}r^{1+\de}|\dk^{\leq\reg+1}\big(\Ga_b\c\pmb\phi_{+2}^{(p)}\big)|^2\nn\\
& +\ep\sum_{p=0}^2\bigg(\B^{\reg}_\de[\pmb\phi_{+2}^{(p)}](\tau_1, \tau_2)+\sup_{\tau\in[\tau_1, \tau_2]}\E^{\reg}[\pmb\phi_{+2}^{(p)}](\tau)\bigg)\nn\\
\les&  \sum_{p=0}^1\int_{\MM(\tt_1, \tt_2)}r^{-1+\de}|\dk^{\leq\reg+1}\N_{T,+2}^{(p)}|^2 +\sum_{p=0}^1\int_{\MM(\tt_1, \tt_2)}r^{1+\de}|\dk^{\leq\reg+1}\big(\Ga_b\c\pmb\phi_{+2}^{(p)}\big)|^2\nn\\
& +\ep\sum_{p=0}^2\bigg(\B^{\reg}_\de[\pmb\phi_{+2}^{(p)}](\tau_1, \tau_2)+\sup_{\tau\in[\tau_1, \tau_2]}\E^{\reg}[\pmb\phi_{+2}^{(p)}](\tau)\bigg) 
\end{align}
and
\bea\lab{eq:controlwidetildeNN14deLiebregminus14pmbphiminus2NWminus2LiebTpregminus14:NTcase}
&&\int_{\Si(\tau_1)}|\dk^{\leq 12}\N_{T,-2, \Lieb_\T}^{(0),\reg-14}|^2+\sum_{p=0}^1\int_{\MM(\tt_1, \tt_2)}r^{-1+\de}|\dk^{\leq {15}}\N_{T,-2, \Lieb_\T}^{(p),\reg-14}|^2\nn\\
&\les& \int_{\Si(\tau_1)}|\dk^{\leq 12}\Lieb_{\T}^{\reg-14}\N_{T,-2}^{(p)}|^2+\sum_{p=0}^1\int_{\MM(\tt_1, \tt_2)}r^{-1+\de}|\dk^{\leq {15}}\Lieb_{\T}^{\reg-14}\N_{T,-2}^{(p)}|^2\nn\\
&&+\ep\sum_{p=0}^2\bigg(\B^{\reg}_\de[\pmb\phi_{-2}^{(p)}](\tau_1, \tau_2)+\sup_{\tau\in[\tau_1, \tau_2]}\E^{\reg}[\pmb\phi_{-2}^{(p)}](\tau)\bigg)\nn\\
&\les& \int_{\Si(\tau_1)}|\dk^{\leq\reg-2}\N_{T,-2}^{(0)}|^2+\sum_{p=0}^1\int_{\MM(\tt_1, \tt_2)}r^{-1+\de}|\dk^{\leq\reg+1}\N_{T,-2}^{(p)}|^2\nn\\
&&+\ep\sum_{p=0}^2\bigg(\B^{\reg}_\de[\pmb\phi_{-2}^{(p)}](\tau_1, \tau_2)+\sup_{\tau\in[\tau_1, \tau_2]}\E^{\reg}[\pmb\phi_{-2}^{(p)}](\tau)\bigg),
\eea
where we have used \eqref{eq:interpolationestimatesforGabGagfork*derivativeswithk*largerthanklover2:forproofTheoremM1andM2} and the control of $\rr-r$ and $\tx^b-x^b$ provided respectively by \eqref{eq:esitmatefor16weightedderivativesofrrminusr} and \eqref{eq:esitmatefor16weightedderivativesoftxbminusxb}. Next, we have in view of \eqref{eq:strctureofNnab4pLiebTregminusregminus2ofAbinredshiftregion} and \eqref{eq:interpolationestimatesforGabGagfork*derivativeswithk*largerthanklover2:forproofTheoremM1andM2}
\beaa
\sum_{p=0}^2\int_{\MM_{r\leq r_+(1+2\dred)}(\tau_1, \tau_2)}|\dk^{\leq {14}}\N_{\nab_4^p\Lieb_{\T}^{\reg-14}\Ab}|^2 &\les& \int_{\MM_{r\leq r_+(1+2\dred)}(\tau_1, \tau_2)}|\dk^{\leq\reg+3}(\Ga_g\c\Ga_b)|^2\\
&\les& \frac{\ep^4}{\tau_1^{3+3\dec}}\les \frac{\ep^2_0}{\tau_1^{3+3\dec}}.
\eeaa
Together with \eqref{MainEnerMora:psi:plus2case:commutationLiebTregminusregplus2case}, \eqref{MainEnerMora:psi:minus2case:commutationLiebTregminusregplus2case}, \eqref{eq:controlwidetildeNN11deLiebregminus11pmbphiplus2NWplus2LiebTpregminus11}, \eqref{eq:controlwidetildeNN14deLiebregminus14pmbphiminus2NWminus2LiebTpregminus14}, \eqref{eq:controlwidetildeNN11deLiebregminus11pmbphiplus2NWplus2LiebTpregminus11:NTcase} and \eqref{eq:controlwidetildeNN14deLiebregminus14pmbphiminus2NWminus2LiebTpregminus14:NTcase}, we infer, for $11\leq\reg\leq\kst -3$,
\beaa
&&\sum_{p=0}^2\EMF^{({11})}_{\de}[\Lieb_{\T}^{\reg-11}\pmb\phi_{+2}^{(p)}](\tau_1, \tau_2)\nn\\
&\les&  \sum_{p=0}^2\E^{\reg}[\pmb\phi_{+2}^{(p)}](\tau_1)+\sum_{p=0}^2\widetilde{\NN}^{\reg}_\de[\pmb\phi_{+2}^{(p)}, \widetilde{\N}_{W,+2}^{(p)}](\tt_1, \tt_2)+ \sum_{p=0}^1\int_{\MM(\tt_1, \tt_2)}r^{-1+\de}|\dk^{\leq\reg+1}\N_{T,+2}^{(p)}|^2\nn\\
&& +\sum_{p=0}^1\int_{\MM(\tt_1, \tt_2)}r^{1+\de}|\dk^{\leq\reg+1}\big(\Ga_b\c\pmb\phi_{+2}^{(p)}\big)|^2 +\ep\sum_{p=0}^2\bigg(\B^{\reg}_\de[\pmb\phi_{+2}^{(p)}](\tau_1, \tau_2)+\sup_{\tau\in[\tau_1, \tau_2]}\E^{\reg}[\pmb\phi_{+2}^{(p)}](\tau)\bigg) 
\eeaa
and, for $14\leq\reg\leq\kst -3$,
\beaa
\nn&&\sum_{p=0}^2\EMF^{({14})}_\de[\Lieb_{\T}^{\reg-14}\pmb\phi_{-2}^{(p)}](\tau_1, \tau_2) + \sum_{p=0}^2\EMF^{({14})}_{r\leq r_+(1+\dred)}[\nab_4^p\Lieb_{\T}^{\reg-14}\Ab](\tau_1, \tau_2)\\
\nn&\les& \sum_{p=0}^2\E^{\reg}[\pmb\phi_{-2}^{(p)}](\tau_1)+\sum_{p=0}^2\E^{\reg}_{r\leq r_+(1+\dred)}[\nab_4^p\Ab](\tau_1)+\int_{\Si(\tau_1)}|\dk^{\leq\reg-2}\N_{T,-2}^{(0)}|^2\nn\\
&&+\sum_{p=0}^1\int_{\MM(\tt_1, \tt_2)}r^{-1+\de}|\dk^{\leq\reg+1}\N_{T,-2}^{(p)}|^2+\sum_{p=0}^2\widetilde{\NN}_\de^{\reg}[\pmb\phi_{-2}^{(p)}, \widetilde{\N}_{W,-2}^{(p)}](\tt_1, \tt_2)\nn\\
&&+\ep\sum_{p=0}^2\bigg(\B^{\reg}_\de[\pmb\phi_{-2}^{(p)}](\tau_1, \tau_2)+\sup_{\tau\in[\tau_1, \tau_2]}\E^{\reg}[\pmb\phi_{-2}^{(p)}](\tau)\bigg)+ \frac{\ep^2_0}{\tau_1^{3+3\dec}}.
\eeaa
This yields, for $11\leq\reg\leq\kst -3$, 
\begin{align}\lab{MainEnerMora:psi:plus2case:commutationLiebTregminusregplus2case:bis}
&\sum_{p=0}^2\EMF_{\de}[\Lieb_{\T}^{\reg}\pmb\phi_{+2}^{(p)}](\tau_1, \tau_2)\nn\\
\les&  \sum_{p=0}^2\E^{\reg}[\pmb\phi_{+2}^{(p)}](\tau_1)+\sum_{p=0}^2\widetilde{\NN}^{\reg}_\de[\pmb\phi_{+2}^{(p)}, \widetilde{\N}_{W,+2}^{(p)}](\tt_1, \tt_2)+ \sum_{p=0}^1\int_{\MM(\tt_1, \tt_2)}r^{-1+\de}|\dk^{\leq\reg+1}\N_{T,+2}^{(p)}|^2\nn\\
& +\sum_{p=0}^1\int_{\MM(\tt_1, \tt_2)}r^{1+\de}|\dk^{\leq\reg+1}\big(\Ga_b\c\pmb\phi_{+2}^{(p)}\big)|^2 +\ep\sum_{p=0}^2\bigg(\B^{\reg}_\de[\pmb\phi_{+2}^{(p)}](\tau_1, \tau_2)+\sup_{\tau\in[\tau_1, \tau_2]}\E^{\reg}[\pmb\phi_{+2}^{(p)}](\tau)\bigg) 
\end{align}
and, for $14\leq\reg\leq\kst -3$, 
\bea\lab{MainEnerMora:psi:minus2case:commutationLiebTregminusregplus2case:bis}
\nn&&\sum_{p=0}^2\EMF_\de[\Lieb_{\T}^{\reg}\pmb\phi_{-2}^{(p)}](\tau_1, \tau_2) + \sum_{p=0}^2\EMF_{r\leq r_+(1+\dred)}[\nab_4^p\Lieb_{\T}^{\reg}\Ab](\tau_1, \tau_2)\\
\nn&\les& \sum_{p=0}^2\E^{\reg}[\pmb\phi_{-2}^{(p)}](\tau_1)+\sum_{p=0}^2\E^{\reg}_{r\leq r_+(1+\dred)}[\nab_4^p\Ab](\tau_1)+\int_{\Si(\tau_1)}|\dk^{\leq\reg-2}\N_{T,-2}^{(0)}|^2\nn\\
&&+\sum_{p=0}^1\int_{\MM(\tt_1, \tt_2)}r^{-1+\de}|\dk^{\leq\reg+1}\N_{T,-2}^{(p)}|^2+\sum_{p=0}^2\widetilde{\NN}_\de^{\reg}[\pmb\phi_{-2}^{(p)}, \widetilde{\N}_{W,-2}^{(p)}](\tt_1, \tt_2)\nn\\
&&+\ep\sum_{p=0}^2\bigg(\B^{\reg}_\de[\pmb\phi_{-2}^{(p)}](\tau_1, \tau_2)+\sup_{\tau\in[\tau_1, \tau_2]}\E^{\reg}[\pmb\phi_{-2}^{(p)}](\tau)\bigg)+ \frac{\ep^2_0}{\tau_1^{3+3\dec}}.
\eea

\noindent{\bf Step 2.} In this step, we rely on \eqref{MainEnerMora:psi:plus2case:commutationLiebTregminusregplus2case:bis} \eqref{MainEnerMora:psi:minus2case:commutationLiebTregminusregplus2case:bis} to obtain the control of $\M_{R,2R}^{\reg}[\pmb\phi_s^{(p)}]$, $s=\pm 2$, $p=0,1,2$, for any $R\geq 22m$ and any $\reg\leq\kst-3$. For $1\leq\reg\leq\kst-3$ and $0\leq j\leq \reg-1$, we introduce a smooth cut-off function $\chi_{j,\reg}=\chi_{j,\reg}(r)$ such that 
\bea\lab{eq:propertyofthecutofffunctionchijregofr}
\chi_{j,\reg}(r)=1\quad\textrm{for}\quad r\in [R_{j,\reg}^{(2)}, R_{j,\reg}^{(3)}],\qquad \chi_{j,\reg}(r)=0\quad\textrm{on}\quad r\in [R_{j,\reg}^{(1)}, R_{j,\reg}^{(4)}],
\eea
where the real numbers $R_{j,\reg}^{(l)}$, $l=1,2,3,4$, are chosen as follows  
\beaa
11m\leq\frac{R}{2}<R_{j,\reg}^{(1)}<R_{j,\reg}^{(2)}<R<2R<R_{j,\reg}^{(3)}<R_{j,\reg}^{(4)}<3R,
\eeaa
which implies in particular 
\beaa
\chi_{j,\reg}(r)=1\quad\textrm{for}\quad r\in [R, 2R], \qquad \textrm{supp}(\chi_{j,\reg})\subset\left(\frac{R}{2}, 3R\right)\subset(11m, 3R).
\eeaa
We also introduce smooth cut-off functions $\widetilde{\chi}_q(\th)$, $q=1,2$, such that 
\bea\lab{eq:propertiesofsmoothcutoffswidetildechiqoftheta}
\widetilde{\chi}_1^2(\th)+\widetilde{\chi}_2^2(\th)=1\quad\textrm{on}\quad[0,\pi], \quad \textrm{supp}(\widetilde{\chi}_1)\subset\left[\frac{\pi}{4}, \frac{3\pi}{4}\right], \quad \textrm{supp}(\widetilde{\chi}_2)\subset\left[0, \frac{\pi}{3}\right]\cup\left[\frac{2\pi}{3}, \pi\right].
\eea

Next, we rely on the coordinates system $(\tau, r, x^1, x^2)$ of $\MM$ and we denote
\bea\lab{eq:defprunweightedderivativeincoordinatesvectorfields}
\pr:=(\nab_{\pr_\tau}, \nab_{\pr_r}, r^{-1}\nab_{\pr_{x^1}}, r^{-1}\nab_{\pr_{x^2}}). 
\eea
Commuting \eqref{eq:TensorialTeuSys:rescaleRHScontaine2:general:Kerrperturbation:alternateformnullframeinsteadcoordvectorfield} with $\widetilde{\chi}_q\chi_{j,\reg}\pr^j\Lieb_{\T}^{\reg-j-1}$, and using also \eqref{eq:tensor:Lsn:onlye_2present:general:Kerrperturbation:alternateformnullframeinsteadcoordvectorfield}, we easily obtain the following schematic tensorial wave equations, for $p=0,1,2$, $s=\pm 2$, $q=1,2$, $1\leq\reg\leq\kst-3$ and $0\leq j\leq \reg-1$,
\beaa
\squared_2(\widetilde{\chi}_q\chi_{j,\reg}\pr^j\Lieb_{\T}^{\reg-j-1}\pmb\phi_s^{(p)}) &=& O(R)\pr^{\leq 2}\big(\widetilde{\chi}_q\chi_{j,\reg}\big)\pr^{\leq\reg}\pmb\phi_s+\chi_{j,\reg}\pr^j\Lieb_{\T}^{\reg-j-1}\widetilde{\N}_{W,s}^{(p)}.
\eeaa
Then, using the regular triplet $\Om_i$, $i=1,2,3$, on $\MM$ exhibited in Lemma \ref{lemma:constructionofhorizontal1formsOmiinperturbationsofKerr}, we define the scalar functions $\phi_{s,kl}^{(p),j,\reg}$ from $\pr^j\Lieb_{\T}^{\reg-j-1}\pmb\phi_s^{(p)}$ as follows, for $p=0,1,2$, $s=\pm 2$,
\bea\lab{eq:definitionofsk2Ctensorsphisklpjreg}
\phi_{s,kl}^{(p),j,\reg}:=\pr^j\Lieb_{\T}^{\reg-j-1}\pmb\phi_s^{(p)}(\Om_k, \Om_l), \quad k,l=1,2,3, \quad 1\leq\reg\leq\kst-3, \quad 0\leq j\leq \reg-1, 
\eea
and scalarize the above system of tensorial wave equations using Lemma \ref{lemma:formoffirstordertermsinscalarazationtensorialwaveeq} which yields the following schematic system of scalar wave equations 
\beaa
\square_\g(\widetilde{\chi}_q\chi_{j,\reg}\phi_{s,kl}^{(p),j,\reg}) = O(R)\pr^{\leq 2}\big(\widetilde{\chi}_q\chi_{j,\reg}\big)\pr^{\leq\reg}\phi_s+\widetilde{\chi}_q\chi_{j,\reg}\pr^j\Lieb_{\T}^{\reg-j-1}\widetilde{\N}_{W,s}^{(p)}(\Om_k, \Om_l).
\eeaa 
Proceeding as in the proof of Lemma 3.13 in \cite{MaSz24}, we decompose the scalar wave operator as $\square_\g=\square_{\gam}+\square_\g-\square_{\gam}$ and obtain, on the support of $\chi_{j,\reg}$,   
\bea\lab{eq:gamalphabetapralphaprbetaofwidechiqchijregphiskljregisequalto}
\nn\gam^{\a\b}\pr_\a\pr_\b(\widetilde{\chi}_q\chi_{j,\reg}\phi_{s,kl}^{(p),j,\reg}) &=& O(R)\pr^{\leq 2}\big(\widetilde{\chi}_q\chi_{j,\reg}\big)\pr^{\leq\reg}\phi_s+O(\ep)\widetilde{\chi}_q\chi_{j,\reg}\pr^{\reg+1}\phi_s\\
&&+\widetilde{\chi}_q\chi_{j,\reg}\pr^j\Lieb_{\T}^{\reg-j-1}\widetilde{\N}_{W,s}^{(p)}(\Om_k, \Om_l),
\eea
and then
\beaa
\gam^{in}\pr_i\pr_n(\widetilde{\chi}_q\chi_{j,\reg}\phi_{s,kl}^{(p),j,\reg}) &=& O(R)\pr^{\leq 2}\big(\widetilde{\chi}_q\chi_{j,\reg}\big)\pr^{\leq\reg}\phi_s+O(R)\widetilde{\chi}_q\chi_{j,\reg}\pr\pr_\tau\phi_{s,kl}^{(p),j,\reg}\\
&&+O(\ep)\widetilde{\chi}_q\chi_{j,\reg}\pr^{\reg+1}\phi_s+\widetilde{\chi}_q\chi_{j,\reg}\pr^j\Lieb_{\T}^{\reg-j-1}\widetilde{\N}_{W,s}^{(p)}(\Om_k, \Om_l).
\eeaa
Squaring and integrating on $\MM(\tau_1, \tau_2)$, and summing over $p=0,1,2$, we infer, since, in view of \eqref{eq:propertyofthecutofffunctionchijregofr}, $\chi_{j,\reg}(r)$ is supported in $r\in [R_{j,\reg}^{(1)}, R_{j,\reg}^{(4)}]$ and in particular in $\Mntrap$, 
\beaa
&&\sum_{p=0}^2\int_{\MM(\tau_1, \tau_2)}|\gam^{in}\pr_i\pr_n(\widetilde{\chi}_q\chi_{j,\reg}\phi_{s,kl}^{(p),j,\reg})|^2\\ 
&\les_R& \sum_{p=0}^2\Big[\M_{R_{j,\reg}^{(1)}, R_{j,\reg}^{(4)}}[\widetilde{\chi}_q\pr_\tau\phi_{s,kl}^{(p),j,\reg}](\tau_1, \tau_2)
+\M_{R_{j,\reg}^{(1)}, R_{j,\reg}^{(4)}}^{(\reg-1)}[\pmb\phi_s^{(p)}](\tau_1, \tau_2)+\ep\M_{11m, 3R}^{\reg}[\pmb\phi_s^{(p)}](\tau_1, \tau_2)\Big]\\
&&+\sum_{p=0}^2\int_{\MM_{11m,3R}(\tau_1, \tau_2)}|\dk^{\leq\reg-1}\widetilde{\N}_{W,s}^{(p)}|^2.
\eeaa
Now, noticing from \eqref{eq:inverse:hypercoord} that, for a scalar function $\psi$,  
\beaa
\gam^{ij}\pr_i\psi\pr_j\psi\gtrsim_R |(\pr_r, \pr_{x^1}, \pr_{x^2})\psi|^2\quad\textrm{for}\,\, r\in(11m,3R), 
\eeaa
and using integration by parts, we infer, noticing that no boundary terms are generated since $\pr_i$, $=1,2,3$ is tangent to $\Si(\tau)$, 
\beaa
&&\sum_{p=0}^2\int_{\MM(\tau_1, \tau_2)}\widetilde{\chi}_q^2\chi_{j,\reg}^2|(\pr_r, \pr_{x^1}, \pr_{x^2})^2\phi_{s,kl}^{(p),j,\reg}|^2\\ 
&\les_R& \sum_{p=0}^2\left(\int_{\MM(\tau_1, \tau_2)}\widetilde{\chi}_q^2\chi_{j,\reg}^2|(\pr_r, \pr_{x^1}, \pr_{x^2})^2\phi_{s,kl}^{(p),j,\reg}|^2\right)^{\frac{1}{2}}\Big(\M_{R_{j,\reg}^{(1)}, R_{j,\reg}^{(4)}}^{(\reg-1)}[\pmb\phi_s^{(p)}](\tau_1, \tau_2)\Big)^{\frac{1}{2}}\\
&&+\sum_{p=0}^2\Big[\M_{R_{j,\reg}^{(1)}, R_{j,\reg}^{(4)}}[\widetilde{\chi}_q\pr_\tau\phi_{s,kl}^{(p),j,\reg}](\tau_1, \tau_2)
+\M_{R_{j,\reg}^{(1)}, R_{j,\reg}^{(4)}}^{(\reg-1)}[\pmb\phi_s^{(p)}](\tau_1, \tau_2)+\ep\M_{11m, 3R}^{\reg}[\pmb\phi_s^{(p)}](\tau_1, \tau_2)\Big]\\
&&+\sum_{p=0}^2\int_{\MM_{11m,3R}(\tau_1, \tau_2)}|\dk^{\leq\reg-1}\widetilde{\N}_{W,s}^{(p)}|^2,
\eeaa
and hence
\beaa
&&\sum_{p=0}^2\int_{\MM(\tau_1, \tau_2)}\widetilde{\chi}_q^2\chi_{j,\reg}^2|(\pr_r, \pr_{x^1}, \pr_{x^2})^2\phi_{s,kl}^{(p),j,\reg}|^2\\ 
&\les_R& \sum_{p=0}^2\Big[\M_{R_{j,\reg}^{(1)}, R_{j,\reg}^{(4)}}[\widetilde{\chi}_q\pr_\tau\phi_{s,kl}^{(p),j,\reg}](\tau_1, \tau_2)
+\M_{R_{j,\reg}^{(1)}, R_{j,\reg}^{(4)}}^{(\reg-1)}[\pmb\phi_s^{(p)}](\tau_1, \tau_2)+\ep\M_{11m, 3R}^{\reg}[\pmb\phi_s^{(p)}](\tau_1, \tau_2)\Big]\\
&&+\sum_{p=0}^2\int_{\MM_{11m,3R}(\tau_1, \tau_2)}|\dk^{\leq\reg-1}\widetilde{\N}_{W,s}^{(p)}|^2.
\eeaa
Since $\chi_{j,\reg}(r)=1$ for $r\in [R_{j,\reg}^{(2)}, R_{j,\reg}^{(3)}]$ in view of \eqref{eq:propertyofthecutofffunctionchijregofr}, we have, 
\beaa
&& \sum_{p=0}^2\M_{R_{j,\reg}^{(2)}, R_{j,\reg}^{(3)}}[\widetilde{\chi}_q\pr\phi_{s,kl}^{(p),j,\reg}](\tau_1, \tau_2)\\
&\les& \sum_{p=0}^2\int_{\MM(\tau_1, \tau_2)}\widetilde{\chi}_q^2\chi_{j,\reg}^2|(\pr_r, \pr_{x^1}, \pr_{x^2})^2\phi_{s,kl}^{(p),j,\reg})|^2\\
&& + \sum_{p=0}^2\Big[\M_{R_{j,\reg}^{(1)}, R_{j,\reg}^{(4)}}[\widetilde{\chi}_q\pr_\tau\phi_{s,kl}^{(p),j,\reg}](\tau_1, \tau_2)+\M_{R_{j,\reg}^{(1)}, R_{j,\reg}^{(4)}}^{(\reg-1)}[\pmb\phi_s^{(p)}](\tau_1, \tau_2)\Big]
\eeaa
and hence
\beaa
&&\sum_{p=0}^2\M_{R_{j,\reg}^{(2)}, R_{j,\reg}^{(3)}}[\widetilde{\chi}_q\pr\phi_{s,kl}^{(p),j,\reg}](\tau_1, \tau_2)\\ 
&\les_R& \sum_{p=0}^2\Big[\M_{R_{j,\reg}^{(1)}, R_{j,\reg}^{(4)}}[\widetilde{\chi}_q\pr_\tau\phi_{s,kl}^{(p),j,\reg}](\tau_1, \tau_2)
+\M_{R_{j,\reg}^{(1)}, R_{j,\reg}^{(4)}}^{(\reg-1)}[\pmb\phi_s^{(p)}](\tau_1, \tau_2)+\ep\M_{11m, 3R}^{\reg}[\pmb\phi_s^{(p)}](\tau_1, \tau_2)\Big]\\
&&+\sum_{p=0}^2\int_{\MM_{11m,3R}(\tau_1, \tau_2)}|\dk^{\leq\reg-1}\widetilde{\N}_{W,s}^{(p)}|^2.
\eeaa
In view of the definition \eqref{eq:definitionofsk2Ctensorsphisklpjreg} of $\phi_{s,kl}^{(p),j,\reg}$, as well as the fact that $\pr_\tau=\T+r\Ga_b\dk$, using also the identity relating $\pr_\tau$ and $\Lieb_{\pr_\tau}$ in Kerr in Lemma \ref{lemma:InKerrlinkbetweenprtauwidehatprtphiandLieb[rtauLiebprtphi}, we deduce
\beaa
\nn&&\sum_{p=0}^2\M_{R_{j,\reg}^{(2)}, R_{j,\reg}^{(3)}}[\widetilde{\chi}_q\pr^{j+1}\Lieb_{\T}^{\reg-j-1}\pmb\phi_s^{(p)}](\tau_1, \tau_2)\\ 
\nn&\les_R& \sum_{p=0}^2\Big[\M_{R_{j,\reg}^{(1)}, R_{j,\reg}^{(4)}}[\widetilde{\chi}_q\pr^j\Lieb_{\T}^{\reg-j}\pmb\phi_s^{(p)}](\tau_1, \tau_2)
+\M_{R_{j,\reg}^{(1)}, R_{j,\reg}^{(4)}}^{(\reg-1)}[\pmb\phi_s^{(p)}](\tau_1, \tau_2)+\ep\M_{11m, 3R}^{\reg}[\pmb\phi_s^{(p)}](\tau_1, \tau_2)\Big]\\
&&+\sum_{p=0}^2\int_{\MM_{11m,3R}(\tau_1, \tau_2)}|\dk^{\leq\reg-1}\widetilde{\N}_{W,s}^{(p)}|^2.
\eeaa
Summing over $q=1,2$, and using \eqref{eq:propertiesofsmoothcutoffswidetildechiqoftheta}, we infer, for $s=\pm 2$, $1\leq\reg\leq\kst-3$ and $0\leq j\leq \reg-1$,
\bea\lab{eq:controlofMrjreg2Rjreg3ofprjplus1LiebTregminusjminus1pmbphisp}
\nn&&\sum_{p=0}^2\M_{R_{j,\reg}^{(2)}, R_{j,\reg}^{(3)}}[\dk^{j+1}\Lieb_{\T}^{\reg-j-1}\pmb\phi_s^{(p)}](\tau_1, \tau_2)\\ 
\nn&\les_R& \sum_{p=0}^2\Big[\M_{R_{j,\reg}^{(1)}, R_{j,\reg}^{(4)}}[\dk^j\Lieb_{\T}^{\reg-j}\pmb\phi_s^{(p)}](\tau_1, \tau_2)
+\M_{R_{j,\reg}^{(1)}, R_{j,\reg}^{(4)}}^{(\reg-1)}[\pmb\phi_s^{(p)}](\tau_1, \tau_2)+\ep\M_{11m, 3R}^{\reg}[\pmb\phi_s^{(p)}](\tau_1, \tau_2)\Big]\\
&&+\sum_{p=0}^2\int_{\MM_{11m,3R}(\tau_1, \tau_2)}|\dk^{\leq\reg-1}\widetilde{\N}_{W,s}^{(p)}|^2.
\eea

Then, relying on \eqref{eq:controlofMrjreg2Rjreg3ofprjplus1LiebTregminusjminus1pmbphisp}, we argue by iteration on $\reg$ and on $j$ to infer, for any $R\geq 22m$, $s=\pm 2$ and $\reg\leq\kst-3$,
\beaa
\sum_{p=0}^2\M_{R, 2R}^{\reg}[\pmb\phi_s^{(p)}](\tau_1, \tau_2) &\les_R& \sum_{p=0}^2\Big[\M_{11m, 3R}[\Lieb_{\T}^{\leq\reg}\pmb\phi_s^{(p)}](\tau_1, \tau_2)+\ep\M_{11m, 3R}^{\reg}[\pmb\phi_s^{(p)}](\tau_1, \tau_2)\Big]\\
&& +\sum_{p=0}^2\int_{\MM_{11m,3R}(\tau_1, \tau_2)}|\dk^{\leq\reg-1}\widetilde{\N}_{W,s}^{(p)}|^2.
\eeaa
Together with the control of $\M_{11m, 3R}[\Lieb_{\T}^{\leq\reg}\pmb\phi_s^{(p)}](\tau_1, \tau_2)$ provided by \eqref{MainEnerMora:psi:plus2case:commutationLiebTregminusregplus2case:bis} and \eqref{MainEnerMora:psi:minus2case:commutationLiebTregminusregplus2case:bis}, we deduce, for any $R\geq 22m$ and $11\leq\reg\leq\kst-3$,
\begin{align}\lab{MainEnerMora:psi:plus2case:commutationLiebTregminusregplus2case:bis:allderivativesonRto2RforRgeq22m}
&\sum_{p=0}^2\M_{R, 2R}^{\reg}[\pmb\phi_{+2}^{(p)}](\tau_1, \tau_2)\nn\\
\les_R&  \sum_{p=0}^2\E^{\reg}[\pmb\phi_{+2}^{(p)}](\tau_1)+\sum_{p=0}^2\widetilde{\NN}^{\reg}_\de[\pmb\phi_{+2}^{(p)}, \widetilde{\N}_{W,+2}^{(p)}](\tt_1, \tt_2)+ \sum_{p=0}^1\int_{\MM(\tt_1, \tt_2)}r^{-1+\de}|\dk^{\leq\reg+1}\N_{T,+2}^{(p)}|^2\nn\\
& +\sum_{p=0}^1\int_{\MM(\tt_1, \tt_2)}r^{1+\de}|\dk^{\leq\reg+1}\big(\Ga_b\c\pmb\phi_{+2}^{(p)}\big)|^2 +\ep\sum_{p=0}^2\bigg(\B^{\reg}_\de[\pmb\phi_{+2}^{(p)}](\tau_1, \tau_2)+\sup_{\tau\in[\tau_1, \tau_2]}\E^{\reg}[\pmb\phi_{+2}^{(p)}](\tau)\bigg) 
\end{align}
and, for $14\leq\reg\leq\kst -3$, 
\bea\lab{MainEnerMora:psi:minus2case:commutationLiebTregminusregplus2case:bis:allderivativesonRto2RforRgeq22m}
\nn&&\sum_{p=0}^2\M_{R, 2R}^{\reg}[\pmb\phi_{-2}^{(p)}](\tau_1, \tau_2)\\
\nn&\les_R& \sum_{p=0}^2\E^{\reg}[\pmb\phi_{-2}^{(p)}](\tau_1)+\sum_{p=0}^2\E^{\reg}_{r\leq r_+(1+\dred)}[\nab_4^p\Ab](\tau_1)+\int_{\Si(\tau_1)}|\dk^{\leq\reg-2}\N_{T,-2}^{(0)}|^2\nn\\
&&+\sum_{p=0}^1\int_{\MM(\tt_1, \tt_2)}r^{-1+\de}|\dk^{\leq\reg+1}\N_{T,-2}^{(p)}|^2+\sum_{p=0}^2\widetilde{\NN}_\de^{\reg}[\pmb\phi_{-2}^{(p)}, \widetilde{\N}_{W,-2}^{(p)}](\tt_1, \tt_2)\nn\\
&&+\ep\sum_{p=0}^2\bigg(\B^{\reg}_\de[\pmb\phi_{-2}^{(p)}](\tau_1, \tau_2)+\sup_{\tau\in[\tau_1, \tau_2]}\E^{\reg}[\pmb\phi_{-2}^{(p)}](\tau)\bigg)+ \frac{\ep^2_0}{\tau_1^{3+3\dec}}.
\eea

\noindent{\bf Step 3.} Next, assuming that $R\geq 22m$ is large enough, we have, for $\reg\leq\kst-3$,
\bea\lab{eq:EMnearinfinity:highorderweightedderivatives:Teu:pm2:analogforhigherderivativesinglobalunprimedframeMM}
\nn&&\sum_{p=0}^2\EMF^{\reg}_{\de, \geq 2R}[\pmb\phi_{s}^{(p)}](\tau_1, \tau_2)\\
\nn&\les_R& \sum_{p=0}^2\M^{\reg}_{R,2R}[\pmb\phi_{s}^{(p)}](\tau_1, \tau_2)+\sum_{p=0}^2\E^{\reg}_{r\geq R}[\pmb\phi_{s}^{(p)}](\tau_1)+\sum_{p=0}^2\int_{\MM_{r\geq R}(\tt_1, \tt_2)}r^{1+\de}|\dk^{\leq \reg}\widetilde{\N}_{W,s}^{(p)}|^2\\
&&+\sum_{p=0}^1\int_{\MM_{r\geq R}(\tt_1, \tt_2)}r^{-1+\de}|\dk^{\leq \reg+1} \N^{(p)}_{T,s}|^2.
\eea
Indeed, \eqref{eq:EMnearinfinity:highorderweightedderivatives:Teu:pm2:analogforhigherderivativesinglobalunprimedframeMM} is the non-sharp analog of Proposition 7.8 in \cite{MaSz26} whose proof in Section 11 of \cite{MaSz26} immediately extends to the case where $\MM$ has the spacelike hypersurface $\Si_*$, instead of $\II_+$, as part of its future boundary.

Next, we fix $R\geq 22m$ large enough so that \eqref{eq:EMnearinfinity:highorderweightedderivatives:Teu:pm2:analogforhigherderivativesinglobalunprimedframeMM} holds, and we may thus forget about the dependence in $R$ from now on. Then, \eqref{eq:EMnearinfinity:highorderweightedderivatives:Teu:pm2:analogforhigherderivativesinglobalunprimedframeMM} together with \eqref{MainEnerMora:psi:plus2case:commutationLiebTregminusregplus2case:bis:allderivativesonRto2RforRgeq22m} \eqref{MainEnerMora:psi:minus2case:commutationLiebTregminusregplus2case:bis:allderivativesonRto2RforRgeq22m} implies, for $11\leq\reg\leq\kst-3$,
\begin{align}\lab{MainEnerMora:psi:plus2case:commutationLiebTregminusregplus2case:bis::allderivativesonrgeqRforRgeq22mlargeenough}
&\sum_{p=0}^2\EMF^{\reg}_{\de, \geq 2R}[\pmb\phi_{+2}^{(p)}](\tau_1, \tau_2)\nn\\
\les&  \sum_{p=0}^2\E^{\reg}[\pmb\phi_{+2}^{(p)}](\tau_1)+\sum_{p=0}^2\widetilde{\NN}^{\reg}_\de[\pmb\phi_{+2}^{(p)}, \widetilde{\N}_{W,+2}^{(p)}](\tt_1, \tt_2)+ \sum_{p=0}^1\int_{\MM(\tt_1, \tt_2)}r^{-1+\de}|\dk^{\leq\reg+1}\N_{T,+2}^{(p)}|^2\nn\\
& +\sum_{p=0}^1\int_{\MM(\tt_1, \tt_2)}r^{1+\de}|\dk^{\leq\reg+1}\big(\Ga_b\c\pmb\phi_{+2}^{(p)}\big)|^2 +\ep\sum_{p=0}^2\bigg(\B^{\reg}_\de[\pmb\phi_{+2}^{(p)}](\tau_1, \tau_2)+\sup_{\tau\in[\tau_1, \tau_2]}\E^{\reg}[\pmb\phi_{+2}^{(p)}](\tau)\bigg) 
\end{align}
and, for $14\leq\reg\leq\kst -3$, 
\bea\lab{MainEnerMora:psi:minus2case:commutationLiebTregminusregplus2case:bis:allderivativesonrgeqRforRgeq22mlargeenough}
\nn&&\sum_{p=0}^2\EMF^{\reg}_{\de, \geq 2R}[\pmb\phi_{-2}^{(p)}](\tau_1, \tau_2)\\
\nn&\les& \sum_{p=0}^2\E^{\reg}[\pmb\phi_{-2}^{(p)}](\tau_1)+\sum_{p=0}^2\E^{\reg}_{r\leq r_+(1+\dred)}[\nab_4^p\Ab](\tau_1)+\int_{\Si(\tau_1)}|\dk^{\leq\reg-2}\N_{T,-2}^{(0)}|^2\nn\\
&&+\sum_{p=0}^1\int_{\MM(\tt_1, \tt_2)}r^{-1+\de}|\dk^{\leq\reg+1}\N_{T,-2}^{(p)}|^2+\sum_{p=0}^2\widetilde{\NN}_\de^{\reg}[\pmb\phi_{-2}^{(p)}, \widetilde{\N}_{W,-2}^{(p)}](\tt_1, \tt_2)\nn\\
&&+\ep\sum_{p=0}^2\bigg(\B^{\reg}_\de[\pmb\phi_{-2}^{(p)}](\tau_1, \tau_2)+\sup_{\tau\in[\tau_1, \tau_2]}\E^{\reg}[\pmb\phi_{-2}^{(p)}](\tau)\bigg)+ \frac{\ep^2_0}{\tau_1^{3+3\dec}}.
\eea

\noindent{\bf Step 4.} Next, let $R\geq 22m$ be the large enough constant fixed at the end Step 3 such that \eqref{MainEnerMora:psi:plus2case:commutationLiebTregminusregplus2case:bis::allderivativesonrgeqRforRgeq22mlargeenough} \eqref{MainEnerMora:psi:minus2case:commutationLiebTregminusregplus2case:bis:allderivativesonrgeqRforRgeq22mlargeenough} hold, and let $\chi_R=\chi_R(r)$ be a smooth cut-off function such that
\bea
\lab{def:cutoffcuntionchi0:00}
0\leq\chi_R\leq 1, \qquad \chi_R=1 \qquad \text{for } \,\, r\leq \frac{5R}{2}, \quad \chi_R=0 \quad \text{for }\,\, r\geq 3R.
\eea
For any $\reg_0+\reg_1\leq \kst-3$ with $\reg_0\geq 1$ and $\reg_1\geq\reg_s$, we commute \eqref{eq:TensorialTeuSys:rescaleRHScontaine2:general:Kerrperturbation:alternateformnullframeinsteadcoordvectorfield} and 
\eqref{eq:transportequationsins=plus2andminus2caseforp=0and1} with $(\chi_R\Lieb_{\Z})^{\reg_0}\Lieb_{\T}^{\reg_1-\reg_s}$. Proceeding as in Step 1, noticing that the proof is significantly simpler since $(\chi_R\Lieb_{\Z})^{\reg_0}\Lieb_{\T}^{\reg_1-\reg_s}$ is localized in $r\leq 3R$ in view of \eqref{def:cutoffcuntionchi0:00} and the fact that $\reg_0\geq 1$, we obtain the following analog of \eqref{MainEnerMora:psi:plus2case:commutationLiebTregminusregplus2case:bis} \eqref{MainEnerMora:psi:minus2case:commutationLiebTregminusregplus2case:bis}, $\reg_0+\reg_1\leq \kst-3$ with $\reg_0\geq 1$, 
\beaa
&&\sum_{p=0}^2\EMF_{\de}[(\chi_R\Lieb_{\Z})^{\reg_0}\Lieb_{\T}^{\reg_1}\pmb\phi_{+2}^{(p)}](\tau_1, \tau_2)\nn\\
&\les&  \sum_{p=0}^2\E^{\reg}[\pmb\phi_{+2}^{(p)}](\tau_1)+\sum_{p=0}^2\widetilde{\NN}^{\reg}_\de[\pmb\phi_{+2}^{(p)}, \widetilde{\N}_{W,+2}^{(p)}](\tt_1, \tt_2)+ \sum_{p=0}^1\int_{\MM(\tt_1, \tt_2)}r^{-1+\de}|\dk^{\leq\reg+1}\N_{T,+2}^{(p)}|^2\nn\\
&&+\sum_{p=0}^2\M^{\reg}_{2R,3R}[\pmb\phi_{+2}^{(p)}](\tau_1, \tau_2) +\sum_{p=0}^1\int_{\MM(\tt_1, \tt_2)}r^{1+\de}|\dk^{\leq\reg+1}\big(\Ga_b\c\pmb\phi_{+2}^{(p)}\big)|^2\\
&& +\ep\sum_{p=0}^2\bigg(\B^{\reg}_\de[\pmb\phi_{+2}^{(p)}](\tau_1, \tau_2)+\sup_{\tau\in[\tau_1, \tau_2]}\E^{\reg}[\pmb\phi_{+2}^{(p)}](\tau)\bigg) 
\eeaa
and
\beaa
\nn&&\sum_{p=0}^2\EMF_\de[(\chi_R\Lieb_{\Z})^{\reg_0}\Lieb_{\T}^{\reg_1}\pmb\phi_{-2}^{(p)}](\tau_1, \tau_2) + \sum_{p=0}^2\EMF_{r\leq r_+(1+\dred)}[\nab_4^p\Lieb_{\T}^{\reg}\Ab](\tau_1, \tau_2)\\
\nn&\les& \sum_{p=0}^2\E^{\reg}[\pmb\phi_{-2}^{(p)}](\tau_1)+\sum_{p=0}^2\E^{\reg}_{r\leq r_+(1+\dred)}[\nab_4^p\Ab](\tau_1)+\int_{\Si(\tau_1)}|\dk^{\leq\reg-2}\N_{T,-2}^{(0)}|^2\nn\\
&&+\sum_{p=0}^1\int_{\MM(\tt_1, \tt_2)}r^{-1+\de}|\dk^{\leq\reg+1}\N_{T,-2}^{(p)}|^2+\sum_{p=0}^2\widetilde{\NN}_\de^{\reg}[\pmb\phi_{-2}^{(p)}, \widetilde{\N}_{W,-2}^{(p)}](\tt_1, \tt_2)\nn\\
&&+\sum_{p=0}^2\M^{\reg}_{2R,3R}[\pmb\phi_{-2}^{(p)}](\tau_1, \tau_2)+\ep\sum_{p=0}^2\bigg(\B^{\reg}_\de[\pmb\phi_{-2}^{(p)}](\tau_1, \tau_2)+\sup_{\tau\in[\tau_1, \tau_2]}\E^{\reg}[\pmb\phi_{-2}^{(p)}](\tau)\bigg)+ \frac{\ep^2_0}{\tau_1^{3+3\dec}}.
\eeaa
Together with \eqref{MainEnerMora:psi:plus2case:commutationLiebTregminusregplus2case:bis} \eqref{MainEnerMora:psi:minus2case:commutationLiebTregminusregplus2case:bis}, we deduce, for $\reg\leq\kst -3$, 
\beaa
&&\sum_{p=0}^2\EMF_{\de}[(\chi_R\Lieb_{\Z}, \Lieb_{\T})^{\reg}\pmb\phi_{+2}^{(p)}](\tau_1, \tau_2)\nn\\
&\les&  \sum_{p=0}^2\E^{\reg}[\pmb\phi_{+2}^{(p)}](\tau_1)+\sum_{p=0}^2\widetilde{\NN}^{\reg}_\de[\pmb\phi_{+2}^{(p)}, \widetilde{\N}_{W,+2}^{(p)}](\tt_1, \tt_2)+ \sum_{p=0}^1\int_{\MM(\tt_1, \tt_2)}r^{-1+\de}|\dk^{\leq\reg+1}\N_{T,+2}^{(p)}|^2\nn\\
&& +\sum_{p=0}^2\M^{\reg}_{2R,3R}[\pmb\phi_{+2}^{(p)}](\tau_1, \tau_2) +\sum_{p=0}^1\int_{\MM(\tt_1, \tt_2)}r^{1+\de}|\dk^{\leq\reg+1}\big(\Ga_b\c\pmb\phi_{+2}^{(p)}\big)|^2 \nn\\
&&+\ep\sum_{p=0}^2\bigg(\B^{\reg}_\de[\pmb\phi_{+2}^{(p)}](\tau_1, \tau_2)+\sup_{\tau\in[\tau_1, \tau_2]}\E^{\reg}[\pmb\phi_{+2}^{(p)}](\tau)\bigg) 
\eeaa
and
\beaa
\nn&&\sum_{p=0}^2\EMF_\de[(\chi_R\Lieb_{\Z}, \Lieb_{\T})^{\reg}\pmb\phi_{-2}^{(p)}](\tau_1, \tau_2) + \sum_{p=0}^2\EMF_{r\leq r_+(1+\dred)}[\nab_4^p(\chi_R\Lieb_{\Z}, \Lieb_{\T})^{\reg}\Ab](\tau_1, \tau_2)\\
\nn&\les& \sum_{p=0}^2\E^{\reg}[\pmb\phi_{-2}^{(p)}](\tau_1)+\sum_{p=0}^2\E^{\reg}_{r\leq r_+(1+\dred)}[\nab_4^p\Ab](\tau_1)+\int_{\Si(\tau_1)}|\dk^{\leq\reg-2}\N_{T,-2}^{(0)}|^2\nn\\
&&+\sum_{p=0}^1\int_{\MM(\tt_1, \tt_2)}r^{-1+\de}|\dk^{\leq\reg+1}\N_{T,-2}^{(p)}|^2+\sum_{p=0}^2\widetilde{\NN}_\de^{\reg}[\pmb\phi_{-2}^{(p)}, \widetilde{\N}_{W,-2}^{(p)}](\tt_1, \tt_2)\nn\\
&&+\sum_{p=0}^2\M^{\reg}_{2R,3R}[\pmb\phi_{-2}^{(p)}](\tau_1, \tau_2) +\ep\sum_{p=0}^2\bigg(\B^{\reg}_\de[\pmb\phi_{-2}^{(p)}](\tau_1, \tau_2)+\sup_{\tau\in[\tau_1, \tau_2]}\E^{\reg}[\pmb\phi_{-2}^{(p)}](\tau)\bigg)+ \frac{\ep^2_0}{\tau_1^{3+3\dec}}.
\eeaa
Finally, relying on \eqref{MainEnerMora:psi:plus2case:commutationLiebTregminusregplus2case:bis::allderivativesonrgeqRforRgeq22mlargeenough} \eqref{MainEnerMora:psi:minus2case:commutationLiebTregminusregplus2case:bis:allderivativesonrgeqRforRgeq22mlargeenough} to control $\M^{\reg}_{2R,3R}[\pmb\phi_s^{(p)}](\tau_1, \tau_2)$ for $s=\pm 2$ and $p=0,1,2$, we infer, for $11\leq\reg\leq\kst -3$,   
\begin{align}\lab{MainEnerMora:psi:plus2case:commutationLiebTregminusregplus2case:ter}
&\sum_{p=0}^2\EMF_{\de}[(\chi_R\Lieb_{\Z}, \Lieb_{\T})^{\reg}\pmb\phi_{+2}^{(p)}](\tau_1, \tau_2)\nn\\
\les&  \sum_{p=0}^2\E^{\reg}[\pmb\phi_{+2}^{(p)}](\tau_1)+\sum_{p=0}^2\widetilde{\NN}^{\reg}_\de[\pmb\phi_{+2}^{(p)}, \widetilde{\N}_{W,+2}^{(p)}](\tt_1, \tt_2)+ \sum_{p=0}^1\int_{\MM(\tt_1, \tt_2)}r^{-1+\de}|\dk^{\leq\reg+1}\N_{T,+2}^{(p)}|^2\nn\\
& +\sum_{p=0}^1\int_{\MM(\tt_1, \tt_2)}r^{1+\de}|\dk^{\leq\reg+1}\big(\Ga_b\c\pmb\phi_{+2}^{(p)}\big)|^2+\ep\sum_{p=0}^2\bigg(\B^{\reg}_\de[\pmb\phi_{+2}^{(p)}](\tau_1, \tau_2)+\sup_{\tau\in[\tau_1, \tau_2]}\E^{\reg}[\pmb\phi_{+2}^{(p)}](\tau)\bigg) 
\end{align}
and, for $14\leq\reg\leq\kst -3$, 
\bea\lab{MainEnerMora:psi:minus2case:commutationLiebTregminusregplus2case:ter}
\nn&&\sum_{p=0}^2\EMF_\de[(\chi_R\Lieb_{\Z}, \Lieb_{\T})^{\reg}\pmb\phi_{-2}^{(p)}](\tau_1, \tau_2) + \sum_{p=0}^2\EMF_{r\leq r_+(1+\dred)}[\nab_4^p(\chi_R\Lieb_{\Z}, \Lieb_{\T})^{\reg}\Ab](\tau_1, \tau_2)\\
\nn&\les& \sum_{p=0}^2\E^{\reg}[\pmb\phi_{-2}^{(p)}](\tau_1)+\sum_{p=0}^2\E^{\reg}_{r\leq r_+(1+\dred)}[\nab_4^p\Ab](\tau_1)+\int_{\Si(\tau_1)}|\dk^{\leq\reg-2}\N_{T,-2}^{(0)}|^2\nn\\
&&+\sum_{p=0}^1\int_{\MM(\tt_1, \tt_2)}r^{-1+\de}|\dk^{\leq\reg+1}\N_{T,-2}^{(p)}|^2+\sum_{p=0}^2\widetilde{\NN}_\de^{\reg}[\pmb\phi_{-2}^{(p)}, \widetilde{\N}_{W,-2}^{(p)}](\tt_1, \tt_2)\nn\\
&& +\ep\sum_{p=0}^2\bigg(\B^{\reg}_\de[\pmb\phi_{-2}^{(p)}](\tau_1, \tau_2)+\sup_{\tau\in[\tau_1, \tau_2]}\E^{\reg}[\pmb\phi_{-2}^{(p)}](\tau)\bigg)+ \frac{\ep^2_0}{\tau_1^{3+3\dec}}.
\eea

\noindent{\bf Step 5.} Next, we rely on \eqref{MainEnerMora:psi:plus2case:commutationLiebTregminusregplus2case:ter} \eqref{MainEnerMora:psi:minus2case:commutationLiebTregminusregplus2case:ter} to recover the flux on $\AA$ and Morawetz estimates for $\pmb\phi_s^{(p)}$, $s=\pm 2$, $p=0,1,2$, in $r\leq 2R$ for higher order derivatives, with $R\geq 22m$ being the large enough constant fixed at the end Step 3 such that \eqref{MainEnerMora:psi:plus2case:commutationLiebTregminusregplus2case:bis::allderivativesonrgeqRforRgeq22mlargeenough} \eqref{MainEnerMora:psi:minus2case:commutationLiebTregminusregplus2case:bis:allderivativesonrgeqRforRgeq22mlargeenough} hold. To this end, we introduce a smooth cut-off function $\chi_{\dred,R}$ satisfying 
\bea
\lab{def:cutoffcuntionchi0:00:dredRversion}
\chi_{\dred,R}=1 \quad \text{for } \,\, r_+(1+\dred)\leq r\leq 2R, \quad \chi_{\dred,R}=0 \quad \text{on }\,\, \mathbb{R}\setminus\left[r_+\left(1+\frac{\dred}{2}\right), \frac{5R}{2}\right].
\eea
Using also the notation $\pr$ in \eqref{eq:defprunweightedderivativeincoordinatesvectorfields}, and the regular triplet $\Om_i$, $i=1,2,3$, on $\MM$ exhibited in Lemma \ref{lemma:constructionofhorizontal1formsOmiinperturbationsofKerr}, we define the scalar functions $\widetilde{\phi}_{s,kl}^{(p),j,\reg}$ from $\pr^j(\Lieb_{\T}, \chi_R\Lieb_{\Z})^{\reg-j-1}\pmb\phi_s^{(p)}$ as follows, for $p=0,1,2$, $s=\pm 2$,
\bea\lab{eq:definitionofsk2Ctensorsphisklpjreg:recoveringderivativesinrleq2R}
\widetilde{\phi}_{s,kl}^{(p),j,\reg}:=\pr^j(\Lieb_{\T}, \chi_R\Lieb_{\Z})^{\reg-j-1}\pmb\phi_s^{(p)}(\Om_k, \Om_l), \,\,\, k,l=1,2,3, \,\,\, 1\leq\reg\leq\kst-3, \,\, 0\leq j\leq \reg-1.
\eea
Then, commuting \eqref{eq:TensorialTeuSys:rescaleRHScontaine2:general:Kerrperturbation:alternateformnullframeinsteadcoordvectorfield} with $\widetilde{\chi}_q\chi_{\dred,R}\pr^j(\Lieb_{\T}, \chi_R\Lieb_{\Z})^{\reg-j-1}$, with the cut-off functions $\widetilde{\chi}_q(\th)$, $q=1,2$, being introduced in \eqref{eq:propertiesofsmoothcutoffswidetildechiqoftheta}, and then scalarizing w.r.t. the regular triplet $\Om_i$, $i=1,2,3$, we obtain the following analog of \eqref{eq:gamalphabetapralphaprbetaofwidechiqchijregphiskljregisequalto}
\bea\lab{eq:gamalphabetapralphaprbetaofwidechiqchijregphiskljregisequalto:recoveringderivativesinrleq2R}
\gam^{\a\b}\pr_\a\pr_\b(\widetilde{\chi}_q\chi_{\dred,R}\widetilde{\phi}_{s,kl}^{(p),j,\reg}) &=& O(1)\pr^{\leq 2}\big(\widetilde{\chi}_q\chi_{\dred,R}\big)\pr^{\leq\reg}\phi_s+O(\ep)\widetilde{\chi}_q\chi_{\dred,R}\pr^{\reg+1}\phi_s\nn\\
&&+\widetilde{\chi}_q\chi_{\dred,R}\pr^j(\Lieb_{\T}, \chi_R\Lieb_{\Z})^{\reg-j-1}\widetilde{\N}_{W,s}^{(p)}(\Om_k, \Om_l),
\eea
where we do not track dependance of constants on $R$ since it have been fixed in Step 3. Then, using (3.33) in \cite{MaSz24}, we infer
\bea\lab{eq:Depr2plusmathringgaabpraprbofwidetildechiqchidrerRwidetildephisklpjregisequaltoRHSforcontrolonrleq2R}
\nn&&\Big(\De\pr_r^2+\mathring{\ga}^{ab}\pr_{x^a}\pr_{x^b}\Big)(\widetilde{\chi}_q\chi_{\dred,R}\widetilde{\phi}_{s,kl}^{(p),j,\reg})\\ 
\nn&=& O(1)(\pr_\tau, \pr_{\tphi})(\pr_\tau, \pr_{\tphi},\pr_r)(\widetilde{\chi}_q\chi_{\dred,R}\widetilde{\phi}_{s,kl}^{(p),j,\reg})+O(1)\pr^{\leq 2}\big(\widetilde{\chi}_q\chi_{\dred,R}\big)\pr^{\leq\reg}\phi_s\\
&&+O(\ep)\widetilde{\chi}_q\chi_{\dred,R}\pr^{\reg+1}\phi_s+\widetilde{\chi}_q\chi_{\dred,R}\pr^j(\Lieb_{\T}, \chi_R\Lieb_{\Z})^{\reg-j-1}\widetilde{\N}_{W,s}^{(p)}(\Om_k, \Om_l).
\eea
Now, we proceed as in Step 5 of \cite{MaSz24} by first multiplying by $\pr_r^2(\widetilde{\chi}_q\chi_{\dred,R}\widetilde{\phi}_{s,kl}^{(p),j,\reg})$ and integrating on $\MM(\tau_1, \tau_2)$ to obtain, after integration by parts, 
\beaa
&&\int_{\MM(\tau_1, \tau_2)}\widetilde{\chi}_q^2\chi_{\dred,R}^2|\pr_r\pr\widetilde{\phi}_{s,kl}^{(p),j,\reg}|^2\\
&\les& \left(\int_{\MM(\tau_1, \tau_2)}\widetilde{\chi}_q^2\chi_{\dred,R}^2|\pr_r\pr\widetilde{\phi}_{s,kl}^{(p),j,\reg}|^2\right)^{\frac{1}{2}}\\
&&\times\Big(\M[\widetilde{\chi}_q\chi_{\dred,R}(\pr_\tau, \pr_{\tphi})\widetilde{\phi}_{s,kl}^{(p),j,\reg}](\tau_1, \tau_2)+\M^{(\reg-1)}_{r\leq 5R/2}[\pmb\phi_s^{(p)}](\tau_1, \tau_2)\Big)^{\frac{1}{2}}\\
&&+\M[\widetilde{\chi}_q\chi_{\dred,R}(\pr_\tau, \pr_{\tphi})\widetilde{\phi}_{s,kl}^{(p),j,\reg}](\tau_1, \tau_2)+\M^{(\reg-1)}_{r\leq 5R/2}[\pmb\phi_s^{(p)}](\tau_1, \tau_2)+\ep\M^{\reg}_{r\leq 5R/2}[\pmb\phi_s^{(p)}](\tau_1, \tau_2)\\
&&+\int_{\MM_{r\leq 5R/2}(\tau_1, \tau_2)}|\dk^{\leq\reg-1}\widetilde{\N}_{W,s}^{(p)}|^2
\eeaa
and hence
\beaa
&&\int_{\MM(\tau_1, \tau_2)}\widetilde{\chi}_q^2\chi_{\dred,R}^2|\pr_r\pr\widetilde{\phi}_{s,kl}^{(p),j,\reg}|^2\\
&\les& \M[\widetilde{\chi}_q\chi_{\dred,R}(\pr_\tau, \widehat{\pr}_{\tphi})\widetilde{\phi}_{s,kl}^{(p),j,\reg}](\tau_1, \tau_2)+\M^{(\reg-1)}_{r\leq 5R/2}[\pmb\phi_s^{(p)}](\tau_1, \tau_2)+\ep\M^{\reg}_{r\leq 5R/2}[\pmb\phi_s^{(p)}](\tau_1, \tau_2)\\
&&+\int_{\MM_{r\leq 5R/2}(\tau_1, \tau_2)}|\dk^{\leq\reg-1}\widetilde{\N}_{W,s}^{(p)}|^2,
\eeaa
where $\widehat{\pr}_{\tphi}$ has been introduced in Definition \ref{def:widehatprtphi}. Then, relying on the definition of $\widetilde{\phi}_{s,kl}^{(p),j,\reg}$ in \eqref{eq:definitionofsk2Ctensorsphisklpjreg:recoveringderivativesinrleq2R}, using also the identity relating $\pr_\tau$ and $\Lieb_{\pr_\tau}$, as well as $\widehat{\pr}_{\tphi}$ and $\Lieb_{\pr_{\tphi}}$, in Kerr in Lemma \ref{lemma:InKerrlinkbetweenprtauwidehatprtphiandLieb[rtauLiebprtphi}, and the fact that $\T=\pr_\tau+r\Ga_b\dk$ and $\Z=\pr_{\tphi}+r^2\Ga_g\dk$, we deduce, for $s=\pm 2$, $p=0,1,2$, $1\leq\reg\leq\kst-3$ and $0\leq j\leq \reg-1$, 
\bea\lab{eq:controlofprrderivativesforMorwetzforwidetildechiqchidredRprjplus1LiebTandchiLiebZofpmbphisp}
\nn&&\int_{\MM(\tau_1, \tau_2)}\widetilde{\chi}_q^2\chi_{\dred,R}^2|\nab_{\pr_r}\pr^{j+1}(\Lieb_{\T}, \chi_R\Lieb_{\Z})^{\reg-j-1}\pmb\phi_s^{(p)}|^2\\
\nn&\les& \M[\chi_{\dred,R}\pr^j(\Lieb_{\T}, \chi_R\Lieb_{\Z})^{\reg-j}\pmb\phi_s^{(p)}](\tau_1, \tau_2)+\M^{(\reg-1)}_{r\leq 5R/2}[\pmb\phi_s^{(p)}](\tau_1, \tau_2)\\
&&+\ep\M^{\reg}_{r\leq 5R/2}[\pmb\phi_s^{(p)}](\tau_1, \tau_2)+\int_{\MM_{r\leq 5R/2}(\tau_1, \tau_2)}|\dk^{\leq\reg-1}\widetilde{\N}_{W,s}^{(p)}|^2.
\eea

Also, rewriting \eqref{eq:gamalphabetapralphaprbetaofwidechiqchijregphiskljregisequalto:recoveringderivativesinrleq2R} using the equation below (3.35) in \cite{MaSz24}, we have
\beaa
\mathring{\ga}^{ab}\pr_{x^a}\pr_{x^b}(\widetilde{\chi}_q\chi_{\dred,R}\widetilde{\phi}_{s,kl}^{(p),j,\reg}) &=& O(1)(\pr_\tau, \pr_{\tphi},\pr_r)\pr(\widetilde{\chi}_q\chi_{\dred,R}\widetilde{\phi}_{s,kl}^{(p),j,\reg})+O(1)\pr^{\leq 2}\big(\widetilde{\chi}_q\chi_{\dred,R}\big)\pr^{\leq\reg}\phi_s\\
&&+O(\ep)\widetilde{\chi}_q\chi_{\dred,R}\pr^{\reg+1}\phi_s+\widetilde{\chi}_q\chi_{\dred,R}\pr^j(\Lieb_{\T}, \chi_R\Lieb_{\Z})^{\reg-j-1}\widetilde{\N}_{W,s}^{(p)}(\Om_k, \Om_l).
\eeaa
Squaring and integrating on $\Mntrap(\tau_1, \tau_2)$, proceeding as above, and adding to \eqref{eq:controlofprrderivativesforMorwetzforwidetildechiqchidredRprjplus1LiebTandchiLiebZofpmbphisp}, we finally obtain, for $s=\pm 2$, $p=0,1,2$, $1\leq\reg\leq\kst-3$ and $0\leq j\leq \reg-1$,  
\beaa
\nn&& \M[\chi_{\dred,R}\dk^{j+1}(\Lieb_{\T}, \chi_R\Lieb_{\Z})^{\reg-j-1}\pmb\phi_s^{(p)}](\tau_1, \tau_2)\\
\nn&\les& \M[\chi_{\dred,R}\dk^j(\Lieb_{\T}, \chi_R\Lieb_{\Z})^{\reg-j}\pmb\phi_s^{(p)}](\tau_1, \tau_2)+\M^{(\reg-1)}_{r\leq 5R/2}[\pmb\phi_s^{(p)}](\tau_1, \tau_2)\\
&&+\ep\M^{\reg}_{r\leq 5R/2}[\pmb\phi_s^{(p)}](\tau_1, \tau_2)+\int_{\MM_{r\leq 5R/2}(\tau_1, \tau_2)}|\dk^{\leq\reg-1}\widetilde{\N}_{W,s}^{(p)}|^2.
\eeaa
Arguing by iteration on $j$, this implies, for $s=\pm 2$, $p=0,1,2$, and $1\leq\reg\leq\kst-3$,
\beaa
\M^{\reg}[\chi_{\dred,R}\pmb\phi_s^{(p)}](\tau_1, \tau_2) &\les& \M[(\Lieb_{\T}, \chi_R\Lieb_{\Z})^{\reg}\pmb\phi_s^{(p)}](\tau_1, \tau_2)+\M^{(\reg-1)}_{r\leq 5R/2}[\pmb\phi_s^{(p)}](\tau_1, \tau_2)\\
&&+\ep\M^{\reg}_{r\leq 5R/2}[\pmb\phi_s^{(p)}](\tau_1, \tau_2)+\int_{\MM_{r\leq 5R/2}(\tau_1, \tau_2)}|\dk^{\leq\reg-1}\widetilde{\N}_{W,s}^{(p)}|^2.
\eeaa
Together with \eqref{MainEnerMora:psi:plus2case:commutationLiebTregminusregplus2case:ter} \eqref{MainEnerMora:psi:minus2case:commutationLiebTregminusregplus2case:ter}, we deduce, for  $11\leq\reg\leq\kst-3$,
\beaa
&&\sum_{p=0}^2\M^{\reg}[\chi_{\dred,R}\pmb\phi_{+2}^{(p)}](\tau_1, \tau_2) \nn\\
&\les&  \sum_{p=0}^2\E^{\reg}[\pmb\phi_{+2}^{(p)}](\tau_1)+\sum_{p=0}^2\widetilde{\NN}^{\reg}_\de[\pmb\phi_{+2}^{(p)}, \widetilde{\N}_{W,+2}^{(p)}](\tt_1, \tt_2)+ \sum_{p=0}^1\int_{\MM(\tt_1, \tt_2)}r^{-1+\de}|\dk^{\leq\reg+1}\N_{T,+2}^{(p)}|^2\nn\\
&& +\sum_{p=0}^1\int_{\MM(\tt_1, \tt_2)}r^{1+\de}|\dk^{\leq\reg+1}\big(\Ga_b\c\pmb\phi_{+2}^{(p)}\big)|^2+\ep\sum_{p=0}^2\bigg(\B^{\reg}_\de[\pmb\phi_{+2}^{(p)}](\tau_1, \tau_2)+\sup_{\tau\in[\tau_1, \tau_2]}\E^{\reg}[\pmb\phi_{+2}^{(p)}](\tau)\bigg)\\
&&+\M^{(\reg-1)}_{r\leq 5R/2}[\pmb\phi_{+2}^{(p)}](\tau_1, \tau_2) 
\eeaa
and, for $14\leq\reg\leq\kst -3$, 
\beaa
\nn&&\sum_{p=0}^2\M^{\reg}[\chi_{\dred,R}\pmb\phi_{-2}^{(p)}](\tau_1, \tau_2)\\
\nn&\les& \sum_{p=0}^2\E^{\reg}[\pmb\phi_{-2}^{(p)}](\tau_1)+\sum_{p=0}^2\E^{\reg}_{r\leq r_+(1+\dred)}[\nab_4^p\Ab](\tau_1)+\int_{\Si(\tau_1)}|\dk^{\leq\reg-2}\N_{T,-2}^{(0)}|^2\nn\\
&&+\sum_{p=0}^1\int_{\MM(\tt_1, \tt_2)}r^{-1+\de}|\dk^{\leq\reg+1}\N_{T,-2}^{(p)}|^2+\sum_{p=0}^2\widetilde{\NN}_\de^{\reg}[\pmb\phi_{-2}^{(p)}, \widetilde{\N}_{W,-2}^{(p)}](\tt_1, \tt_2)\nn\\
&& +\ep\sum_{p=0}^2\bigg(\B^{\reg}_\de[\pmb\phi_{-2}^{(p)}](\tau_1, \tau_2)+\sup_{\tau\in[\tau_1, \tau_2]}\E^{\reg}[\pmb\phi_{-2}^{(p)}](\tau)\bigg)+ \frac{\ep^2_0}{\tau_1^{3+3\dec}}+\M^{(\reg-1)}_{r\leq 5R/2}[\pmb\phi_{-2}^{(p)}](\tau_1, \tau_2).
\eeaa
Then, applying standard redshift estimates, see in particular the general redshift estimates of Section 9.4 in \cite{GKS22}, to the tensorial wave equations \eqref{eq:TensorialTeuSys:rescaleRHScontaine2:general:Kerrperturbation:alternateformnullframeinsteadcoordvectorfield} and \eqref{eq:waveequationpmbphip=0sminus2nodeginredshiftregion}, we obtain, using also \eqref{def:cutoffcuntionchi0:00:dredRversion}, for  $11\leq\reg\leq\kst-3$,
\beaa
&&\sum_{p=0}^2\Big(\EF_{r\leq r_+(1+\dred)}^{\reg}[\pmb\phi_{+2}^{(p)}](\tau_1, \tau_2)+\M^{\reg}_{r\leq 2R}[\pmb\phi_{+2}^{(p)}](\tau_1, \tau_2)\Big) \nn\\
&\les&  \sum_{p=0}^2\E^{\reg}[\pmb\phi_{+2}^{(p)}](\tau_1)+\sum_{p=0}^2\widetilde{\NN}^{\reg}_\de[\pmb\phi_{+2}^{(p)}, \widetilde{\N}_{W,+2}^{(p)}](\tt_1, \tt_2)+ \sum_{p=0}^1\int_{\MM(\tt_1, \tt_2)}r^{-1+\de}|\dk^{\leq\reg+1}\N_{T,+2}^{(p)}|^2\nn\\
&& +\sum_{p=0}^1\int_{\MM(\tt_1, \tt_2)}r^{1+\de}|\dk^{\leq\reg+1}\big(\Ga_b\c\pmb\phi_{+2}^{(p)}\big)|^2+\ep\sum_{p=0}^2\bigg(\B^{\reg}_\de[\pmb\phi_{+2}^{(p)}](\tau_1, \tau_2)+\sup_{\tau\in[\tau_1, \tau_2]}\E^{\reg}[\pmb\phi_{+2}^{(p)}](\tau)\bigg)\\
&&+\M^{(\reg-1)}_{r\leq 5R/2}[\pmb\phi_{+2}^{(p)}](\tau_1, \tau_2) 
\eeaa
and, for $14\leq\reg\leq\kst -3$, 
\beaa
\nn&&\sum_{p=0}^2\Big(\EF_{r\leq r_+(1+\dred)}^{\reg}[\pmb\phi_{-2}^{(p)}](\tau_1, \tau_2)+\EMF^{\reg}_{r\leq r_+(1+\dred)}[\nab_4^p\Ab](\tau_1, \tau_2)+\M^{\reg}_{r\leq 2R}[\pmb\phi_{-2}^{(p)}](\tau_1, \tau_2)\Big)\\
\nn&\les& \sum_{p=0}^2\E^{\reg}[\pmb\phi_{-2}^{(p)}](\tau_1)+\sum_{p=0}^2\E^{\reg}_{r\leq r_+(1+\dred)}[\nab_4^p\Ab](\tau_1)+\int_{\Si(\tau_1)}|\dk^{\leq\reg-2}\N_{T,-2}^{(0)}|^2\nn\\
&&+\sum_{p=0}^1\int_{\MM(\tt_1, \tt_2)}r^{-1+\de}|\dk^{\leq\reg+1}\N_{T,-2}^{(p)}|^2+\sum_{p=0}^2\widetilde{\NN}_\de^{\reg}[\pmb\phi_{-2}^{(p)}, \widetilde{\N}_{W,-2}^{(p)}](\tt_1, \tt_2)\nn\\
&& +\ep\sum_{p=0}^2\bigg(\B^{\reg}_\de[\pmb\phi_{-2}^{(p)}](\tau_1, \tau_2)+\sup_{\tau\in[\tau_1, \tau_2]}\E^{\reg}[\pmb\phi_{-2}^{(p)}](\tau)\bigg)+ \frac{\ep^2_0}{\tau_1^{3+3\dec}}+\M^{(\reg-1)}_{r\leq 5R/2}[\pmb\phi_{-2}^{(p)}](\tau_1, \tau_2).
\eeaa
Together with \eqref{MainEnerMora:psi:plus2case:commutationLiebTregminusregplus2case:bis::allderivativesonrgeqRforRgeq22mlargeenough} and \eqref{MainEnerMora:psi:minus2case:commutationLiebTregminusregplus2case:bis:allderivativesonrgeqRforRgeq22mlargeenough}, we infer, for  $11\leq\reg\leq\kst-3$,
\beaa
&&\sum_{p=0}^2\Big(\EF_{r\leq r_+(1+\dred)}^{\reg}[\pmb\phi_{+2}^{(p)}](\tau_1, \tau_2)+\M^{\reg}_\de[\pmb\phi_{+2}^{(p)}](\tau_1, \tau_2)+\EF^{\reg}_{r\geq 2R}[\pmb\phi_{+2}^{(p)}](\tau_1, \tau_2)\Big) \nn\\
&\les&  \sum_{p=0}^2\E^{\reg}[\pmb\phi_{+2}^{(p)}](\tau_1)+\sum_{p=0}^2\widetilde{\NN}^{\reg}_\de[\pmb\phi_{+2}^{(p)}, \widetilde{\N}_{W,+2}^{(p)}](\tt_1, \tt_2)+ \sum_{p=0}^1\int_{\MM(\tt_1, \tt_2)}r^{-1+\de}|\dk^{\leq\reg+1}\N_{T,+2}^{(p)}|^2\nn\\
&& +\sum_{p=0}^1\int_{\MM(\tt_1, \tt_2)}r^{1+\de}|\dk^{\leq\reg+1}\big(\Ga_b\c\pmb\phi_{+2}^{(p)}\big)|^2+\ep\sum_{p=0}^2\bigg(\B^{\reg}_\de[\pmb\phi_{+2}^{(p)}](\tau_1, \tau_2)+\sup_{\tau\in[\tau_1, \tau_2]}\E^{\reg}[\pmb\phi_{+2}^{(p)}](\tau)\bigg)\\
&&+\M^{(\reg-1)}_{r\leq 2R}[\pmb\phi_{+2}^{(p)}](\tau_1, \tau_2) 
\eeaa
and, for $14\leq\reg\leq\kst -3$, 
\beaa
\nn&&\sum_{p=0}^2\Big(\EF_{r\leq r_+(1+\dred)}^{\reg}[\pmb\phi_{-2}^{(p)}](\tau_1, \tau_2)+\M^{\reg}_\de[\pmb\phi_{-2}^{(p)}](\tau_1, \tau_2)+\EF^{\reg}_{r\geq 2R}[\pmb\phi_{-2}^{(p)}](\tau_1, \tau_2)\Big) \nn\\
&&\sum_{p=0}^2\EMF^{\reg}_{r\leq r_+(1+\dred)}[\nab_4^p\Ab](\tau_1, \tau_2)\\
\nn&\les& \sum_{p=0}^2\E^{\reg}[\pmb\phi_{-2}^{(p)}](\tau_1)+\sum_{p=0}^2\E^{\reg}_{r\leq r_+(1+\dred)}[\nab_4^p\Ab](\tau_1)+\int_{\Si(\tau_1)}|\dk^{\leq\reg-2}\N_{T,-2}^{(0)}|^2\nn\\
&&+\sum_{p=0}^1\int_{\MM(\tt_1, \tt_2)}r^{-1+\de}|\dk^{\leq\reg+1}\N_{T,-2}^{(p)}|^2+\sum_{p=0}^2\widetilde{\NN}_\de^{\reg}[\pmb\phi_{-2}^{(p)}, \widetilde{\N}_{W,-2}^{(p)}](\tt_1, \tt_2)\nn\\
&& +\ep\sum_{p=0}^2\bigg(\B^{\reg}_\de[\pmb\phi_{-2}^{(p)}](\tau_1, \tau_2)+\sup_{\tau\in[\tau_1, \tau_2]}\E^{\reg}[\pmb\phi_{-2}^{(p)}](\tau)\bigg)+ \frac{\ep^2_0}{\tau_1^{3+3\dec}}+\M^{(\reg-1)}_{r\leq 2R}[\pmb\phi_{-2}^{(p)}](\tau_1, \tau_2).
\eeaa
Finally, arguing by iteration on $\reg$ and using Theorem \ref{thm:main:MaSz26}, we deduce, for  $11\leq\reg\leq\kst-3$,
\begin{align}\lab{eq:EnergyMorawetzFlusdelahigherorderderivtaiavespmbphispalmostcompleandonlylackinenergyonrplus1plusdredto2R:plus2}
&\sum_{p=0}^2\Big(\EF_{r\leq r_+(1+\dred)}^{\reg}[\pmb\phi_{+2}^{(p)}](\tau_1, \tau_2)+\M^{\reg}_\de[\pmb\phi_{+2}^{(p)}](\tau_1, \tau_2)+\EF^{\reg}_{r\geq 2R}[\pmb\phi_{+2}^{(p)}](\tau_1, \tau_2)\Big) \nn\\
\les&  \sum_{p=0}^2\E^{\reg}[\pmb\phi_{+2}^{(p)}](\tau_1)+\sum_{p=0}^2\widetilde{\NN}^{\reg}_\de[\pmb\phi_{+2}^{(p)}, \widetilde{\N}_{W,+2}^{(p)}](\tt_1, \tt_2)+ \sum_{p=0}^1\int_{\MM(\tt_1, \tt_2)}r^{-1+\de}|\dk^{\leq\reg+1}\N_{T,+2}^{(p)}|^2\nn\\
& +\sum_{p=0}^1\int_{\MM(\tt_1, \tt_2)}r^{1+\de}|\dk^{\leq\reg+1}\big(\Ga_b\c\pmb\phi_{+2}^{(p)}\big)|^2+\ep\sum_{p=0}^2\bigg(\B^{\reg}_\de[\pmb\phi_{+2}^{(p)}](\tau_1, \tau_2)+\sup_{\tau\in[\tau_1, \tau_2]}\E^{\reg}[\pmb\phi_{+2}^{(p)}](\tau)\bigg)
\end{align}
and, for $14\leq\reg\leq\kst -3$, 
\bea\lab{eq:EnergyMorawetzFlusdelahigherorderderivtaiavespmbphispalmostcompleandonlylackinenergyonrplus1plusdredto2R:minus2}
\nn&&\sum_{p=0}^2\Big(\EF_{r\leq r_+(1+\dred)}^{\reg}[\pmb\phi_{-2}^{(p)}](\tau_1, \tau_2)+\M^{\reg}_\de[\pmb\phi_{-2}^{(p)}](\tau_1, \tau_2)+\EF^{\reg}_{r\geq 2R}[\pmb\phi_{-2}^{(p)}](\tau_1, \tau_2)\Big) \nn\\
\nn&&\sum_{p=0}^2\EMF^{\reg}_{r\leq r_+(1+\dred)}[\nab_4^p\Ab](\tau_1, \tau_2)\\
\nn&\les& \sum_{p=0}^2\E^{\reg}[\pmb\phi_{-2}^{(p)}](\tau_1)+\sum_{p=0}^2\E^{\reg}_{r\leq r_+(1+\dred)}[\nab_4^p\Ab](\tau_1)+\int_{\Si(\tau_1)}|\dk^{\leq\reg-2}\N_{T,-2}^{(0)}|^2\nn\\
&&+\sum_{p=0}^1\int_{\MM(\tt_1, \tt_2)}r^{-1+\de}|\dk^{\leq\reg+1}\N_{T,-2}^{(p)}|^2+\sum_{p=0}^2\widetilde{\NN}_\de^{\reg}[\pmb\phi_{-2}^{(p)}, \widetilde{\N}_{W,-2}^{(p)}](\tt_1, \tt_2)\nn\\
&& +\ep\sum_{p=0}^2\bigg(\B^{\reg}_\de[\pmb\phi_{-2}^{(p)}](\tau_1, \tau_2)+\sup_{\tau\in[\tau_1, \tau_2]}\E^{\reg}[\pmb\phi_{-2}^{(p)}](\tau)\bigg)+ \frac{\ep^2_0}{\tau_1^{3+3\dec}}.
\eea

\noindent{\bf Step 6.} In view of \eqref{eq:EnergyMorawetzFlusdelahigherorderderivtaiavespmbphispalmostcompleandonlylackinenergyonrplus1plusdredto2R:plus2} \eqref{eq:EnergyMorawetzFlusdelahigherorderderivtaiavespmbphispalmostcompleandonlylackinenergyonrplus1plusdredto2R:minus2}, it remains to recover the energy for higher order derivatives of $\pmb\phi_s^{(p)}$, $s=\pm 2$, $p=0,1,2$, for $r\in[r_+(1+\dred), 2R]$. To this end, we square \eqref{eq:Depr2plusmathringgaabpraprbofwidetildechiqchidrerRwidetildephisklpjregisequaltoRHSforcontrolonrleq2R} and integrate it over $\Si(\tau)$ for $\tau\in[\tau_1, \tau_2]$. We then integrate by parts the LHS, noticing that all derivatives involved are tangent to $\Si(\tau)$. Proceeding similarly as in Step 5, we obtain, for $s=\pm 2$, $p=0,1,2$, $1\leq\reg\leq\kst-3$, $0\leq j\leq \reg-1$ and $\tau\in[\tau_1, \tau_2]$, 
\begin{align*}
\nn& \E[\chi_{\dred,R}\dk^{j+1}(\Lieb_{\T}, \chi_R\Lieb_{\Z})^{\reg-j-1}\pmb\phi_s^{(p)}](\tau)\\
\nn\les& \E[\chi_{\dred,R}\dk^j(\Lieb_{\T}, \chi_R\Lieb_{\Z})^{\reg-j}\pmb\phi_s^{(p)}](\tau)+\E^{(\reg-1)}_{r\leq 5R/2}[\pmb\phi_s^{(p)}](\tau)+\ep\E^{\reg}_{r\leq 5R/2}[\pmb\phi_s^{(p)}](\tau)+\int_{\Si_{r\leq 5R/2}(\tau)}|\dk^{\leq\reg-1}\widetilde{\N}_{W,s}^{(p)}|^2.
\end{align*}
Arguing by iteration on $j$, and using a trace estimate for the last term, we obtain, for $s=\pm 2$, $p=0,1,2$, $1\leq\reg\leq\kst-3$ and $\tau\in[\tau_1, \tau_2]$,
\beaa
 \E^{\reg}[\chi_{\dred,R}\pmb\phi_s^{(p)}](\tau) &\les& \E[\chi_{\dred,R}(\Lieb_{\T}, \chi_R\Lieb_{\Z})^{\reg}\pmb\phi_s^{(p)}](\tau)+\E^{(\reg-1)}_{r\leq 5R/2}[\pmb\phi_s^{(p)}](\tau)+\ep\E^{\reg}_{r\leq 5R/2}[\pmb\phi_s^{(p)}](\tau)\\
&&+\int_{\MM_{r\leq 5R/2}(\tau_1, \tau_2)}|\dk^{\leq\reg}\widetilde{\N}_{W,s}^{(p)}|^2.
\eeaa
In view of \eqref{def:cutoffcuntionchi0:00:dredRversion}, arguing by iteration on $\reg$, and taking the supremum in $\tau\in[\tau_1, \tau_2]$, we infer
\beaa
 &&\sup_{\tau\in[\tau_1, \tau_2]}\E^{\reg}_{r_+(1+\dred),2R}[\pmb\phi_s^{(p)}](\tau)\\ 
 &\les& \sup_{\tau\in[\tau_1, \tau_2]}\E[(\Lieb_{\T}, \chi_R\Lieb_{\Z})^{\reg}\pmb\phi_s^{(p)}](\tau)+\sup_{\tau\in[\tau_1, \tau_2]}\E_{\leq 2R}[\pmb\phi_s^{(p)}](\tau)\\
 &&+\sup_{\tau\in[\tau_1, \tau_2]}\E^{(\reg-1)}_{2R,5R/2}[\pmb\phi_s^{(p)}](\tau)+\ep\sup_{\tau\in[\tau_1, \tau_2]}\E^{\reg}_{r\leq 5R/2}[\pmb\phi_s^{(p)}](\tau)+\int_{\MM_{r\leq 5R/2}(\tau_1, \tau_2)}|\dk^{\leq\reg}\widetilde{\N}_{W,s}^{(p)}|^2.
\eeaa
Finally, relying on \eqref{MainEnerMora:psi:plus2case:commutationLiebTregminusregplus2case:ter} \eqref{MainEnerMora:psi:minus2case:commutationLiebTregminusregplus2case:ter} to control the first term on the RHS, on Theorem \ref{thm:main:MaSz26} to control the second term on the RHS, and on \eqref{eq:EnergyMorawetzFlusdelahigherorderderivtaiavespmbphispalmostcompleandonlylackinenergyonrplus1plusdredto2R:plus2} \eqref{eq:EnergyMorawetzFlusdelahigherorderderivtaiavespmbphispalmostcompleandonlylackinenergyonrplus1plusdredto2R:minus2} to control the third term on the RHS, and using $\ep>0$ small enough, we obtain, for  $11\leq\reg\leq\kst-3$,
\beaa
&&\sum_{p=0}^2\EMF_\de^{\reg}[\pmb\phi_{+2}^{(p)}](\tau_1, \tau_2) \nn\\
&\les&  \sum_{p=0}^2\E^{\reg}[\pmb\phi_{+2}^{(p)}](\tau_1)+\sum_{p=0}^2\widetilde{\NN}^{\reg}_\de[\pmb\phi_{+2}^{(p)}, \widetilde{\N}_{W,+2}^{(p)}](\tt_1, \tt_2)+ \sum_{p=0}^1\int_{\MM(\tt_1, \tt_2)}r^{-1+\de}|\dk^{\leq\reg+1}\N_{T,+2}^{(p)}|^2\nn\\
&& +\sum_{p=0}^1\int_{\MM(\tt_1, \tt_2)}r^{1+\de}|\dk^{\leq\reg+1}\big(\Ga_b\c\pmb\phi_{+2}^{(p)}\big)|^2+\ep\sum_{p=0}^2\B^{\reg}_\de[\pmb\phi_{+2}^{(p)}](\tau_1, \tau_2)
\eeaa
and, for $14\leq\reg\leq\kst -3$, 
\beaa
\nn&&\sum_{p=0}^2\Big(\EMF_\de^{\reg}[\pmb\phi_{-2}^{(p)}](\tau_1, \tau_2)+\EMF^{\reg}_{r\leq r_+(1+\dred)}[\nab_4^p\Ab](\tau_1, \tau_2)\Big)\\
\nn&\les& \sum_{p=0}^2\E^{\reg}[\pmb\phi_{-2}^{(p)}](\tau_1)+\sum_{p=0}^2\E^{\reg}_{r\leq r_+(1+\dred)}[\nab_4^p\Ab](\tau_1)+\int_{\Si(\tau_1)}|\dk^{\leq\reg-2}\N_{T,-2}^{(0)}|^2\nn\\
&&+\sum_{p=0}^1\int_{\MM(\tt_1, \tt_2)}r^{-1+\de}|\dk^{\leq\reg+1}\N_{T,-2}^{(p)}|^2+\sum_{p=0}^2\widetilde{\NN}_\de^{\reg}[\pmb\phi_{-2}^{(p)}, \widetilde{\N}_{W,-2}^{(p)}](\tt_1, \tt_2)\nn\\
&& +\ep\sum_{p=0}^2\B^{\reg}_\de[\pmb\phi_{-2}^{(p)}](\tau_1, \tau_2)+ \frac{\ep^2_0}{\tau_1^{3+3\dec}}.
\eeaa
This concludes the proof of Theorem \ref{thm:main:MaSz26:extendhigherorderderivatives}.
\end{proof}

%%%%%%%%%%%%%%%%%%%%%%%%%%%%%%%%%%%%%%%%%%%

\subsection{$r^p$-weighted estimates for Teukolsky equations in $\MM$}
\lab{sec:rpweightedestimatesforwaveequationsandTeukolskyinMM}

%%%%%%%%%%%%%%%%%%%%%%%%%%%%%%%%%%%%%%%%%%%

Now that Morawetz estimates for Teukolsky equations have been derived in Theorem \ref{thm:main:MaSz26:extendhigherorderderivatives}, we turn to $r^p$-weighted estimates. 
Given that we only need to focus on the region $r\geq R$ with $R\gg m$ large enough, we may rely on the results of Sections 10 to 12 in \cite{GKS22} as they are valid for all $|a|<m$ on such large $r$ regions.

%%%%%%%%%%%%%%%%%%%%%%%%%%%%%%%%%%%%%%%%%%%

\subsubsection{$r^p$-weighted estimates for a basic tensorial wave equation}

%%%%%%%%%%%%%%%%%%%%%%%%%%%%%%%%%%%%%%%%%%%

We first consider the solution $\psi\in\sk_2(\mathbb{C})$ satisfying on $\MM$ the following tensorial wave equation 
\bea\lab{eq:Gen.RW-chap10}
\bigg(\squared_2 -\frac{4ia\cos\th}{|q|^2}\nab_{\T}- \frac{4}{|q|^2}\bigg)\psi = N,
\eea
for some RHS $N$. Then, the following basic $r^p$-weighted estimates hold for solutions to \eqref{eq:Gen.RW-chap10}. 

\begin{proposition}[Basic $r^p$-weighted estimates for \eqref{eq:Gen.RW-chap10}]
\lab{Proposition:Step3-Chap10}
Let $R\gg m$ large enough. We have, for $\de\leq p\leq 2-\de$,  $0\le s\le k_L$,
\bea
\lab{eq:Proposition-Step3-Chap10}
\bsplit
\BEF^s_{p, \geq R}[\psi](\tau_1, \tau_2) \les& \E^s_{p,\geq \frac{R}{2}}[\psi](\tau_1)+\NN^s_{p, \ge R/2}[\psi, N]{(\tau_1, \tau_2)} +\M^s_{\frac{R}{2}, R}[\psi](\tau_1, \tau_2).
\end{split}
\eea
\end{proposition}

\begin{proof}
Note that for $\psi\in\sk_2(\mathbb{C})$ satisfying \eqref{eq:Gen.RW-chap10}, $\Re(\psi)\in\sk_2$ satisfies
\beaa
&&\bigg(\squared_2 -\frac{4a\cos\th}{|q|^2}\dual\nab_{\T}- V\bigg)\Re(\psi) = \widetilde{N}, \\
&&\widetilde{N}:=\Re(N)+V_1\Re(\psi),\qquad  V:= \frac{4\De}{ (r^2+a^2) |q|^2}, \qquad V_1=\frac{4}{|q|^2}-V=O(mr^{-3}),
\eeaa
which is  (10.0.1) in \cite{GKS22}. We may thus directly apply the $r^p$-weighted estimates of Proposition 10.1.2 in \cite{GKS22} which yield, for $r\geq R$ with $R\gg m$ large enough, 
\beaa
\BEF^s_{p, \geq R}[\psi](\tau_1, \tau_2) &\les& \E^s_{p,\geq \frac{R}{2}}[\psi](\tau_1)+\NN^s_{p, \ge R/2}[\Re(\psi), \widetilde{N}]{(\tau_1, \tau_2)} +\M^s_{\frac{R}{2}, R}[\psi](\tau_1, \tau_2).
\eeaa
Also, since $V_1=O(mr^{-3})$, we have
\beaa
\NN^s_{p, \ge R/2}[\Re(\psi), \widetilde{N}]{(\tau_1, \tau_2)} &\les& \NN^s_{p, \ge R/2}[\psi, N]{(\tau_1, \tau_2)}+\NN^s_{p, \ge R/2}[\Re(\psi), O(r^{-3})\Re(\psi)]{(\tau_1, \tau_2)}
\eeaa
where 
\beaa
&&\NN^s_{p, \ge R/2}[\Re(\psi), O(r^{-3})\Re(\psi)]{(\tau_1, \tau_2)} \\
&\les& \int_{\MM_{r\geq R/2}(\tau_1, \tau_2)}r^{1+\de-6}|\dk^{\leq s}\psi|^2+\int_{\MM_{r\geq R/2}(\tau_1, \tau_2)}r^{p-1}|\nab_4(r\dk^{\leq s}\psi)|r^{-3}|\dk^{\leq s}\psi|\\
&\les& R^{-1}\B^s_{p, \geq R/2}[\psi](\tau_1, \tau_2)
\eeaa
can be absorbed for $R$ large enough. We thus deduce 
\beaa
\BEF^s_{p, \geq R}[\psi](\tau_1, \tau_2) \les \E^s_{p,\geq \frac{R}{2}}[\psi](\tau_1)+\NN^s_{p, \ge R/2}[\psi, N]{(\tau_1, \tau_2)} +\M^s_{\frac{R}{2}, R}[\psi](\tau_1, \tau_2)
\eeaa
as stated. This concludes the proof of Proposition \ref{Proposition:Step3-Chap10}.
\end{proof}

While we have relied on the full wave-transport system derived in Theorem \ref{thm:derivationoftheTeukolskytensorialwavesystemfors=plusminus2:kerrpert:alternateformnullframeinsteadcoordvectorfield} and Proposition \ref{prop:transportequationphis=plusminus2p=0and1} for the energy-Morawetz estimates of Theorem \ref{thm:main:MaSz26:extendhigherorderderivatives}, following \cite{GKS22}, we will only use the tensorial wave equation for $\pmb\phi_s^{(2)}$ with $s=\pm 2$ and the transport equations of Proposition \ref{prop:transportequationphis=plusminus2p=0and1}  to derive $r^p$-weighted estimates for Teukolsky. In the next section we derive the needed transport estimates.

%%%%%%%%%%%%%%%%%%%%%%%%%%%%%%%%%%%%%%%%%%%%%%

\subsubsection{Weighted estimates for transport equations}

%%%%%%%%%%%%%%%%%%%%%%%%%%%%%%%%%%%%%%%%%%%%%%

We start with transport equations along $e_4$.

\begin{lemma}\lab{lemma:transportlemmaforspinminus2:largerregion}
Assume that $\pmb\phi, \,\pmb\psi\in\sk_2(\mathbb{C})$ are such that 
\bea
\lab{eq:nab4phi=psi:transportestis}
\nab_4\pmb\phi=\pmb\psi.
\eea
 Then, for $\de_0>0$, $\reg\leq\kl$,  and for $R$ large enough, we have for any $1\leq\tau_1<\tau_2\leq\tau_*$
\begin{align}
\lab{eq:transporte4:general:anyordersummary:lemma}
&\int_{\MM_{r\geq R}(\tau_1, \tau_2)}r^{-\de_0-3}|\dk^{\leq \reg}\pmb\phi|^2+\int_{\Si_{r\geq R}(\tau_2)}r^{-\de_0-4}|\dk^{\leq \reg}\pmb\phi|^2+\int_{\Si_*(\tau_1, \tau_2)}r^{-\de_0-2}|\dk^{\leq \reg}\pmb\phi|^2\nn\\
\les_{\de_0}&\int_{\Si_{r\geq R/2}(\tau_1)}r^{-\de_0-4}|\dk^{\leq \reg}\pmb\phi|^2+\int_{\MM_{R/2,R}(\tau_1, \tau_2)}r^{-\de_0-3}|\dk^{\leq \reg}\pmb\phi|^2+ \int_{\MM_{r\geq R/2}(\tau_1, \tau_2)}r^{-\de_0-1}|\dk^{\leq \reg}\pmb\psi|^2.
\end{align}
\end{lemma}

\begin{proof}
We compute 
\beaa
\Div ( e_4) &=& \g^{43}\g(\D_4e_4, e_3) +\g^{43}\g(\D_3e_4, e_4)+\g^{bc}\g(\D_be_4, e_c)\\
&=& -\frac{1}{2}4\om +\trch = \frac{2}{r}+O(mr^{-2}) -2\omc+\trchc\\
&=& \frac{2}{r}+O(mr^{-2}) +\Ga_g= \frac{2}{r}+O(mr^{-2}) 
\eeaa
and 
\beaa
e_4(r)=1+O(mr^{-1})+\widecheck{e_4(r)}=1+O(mr^{-1})+\Ga_g=1+O(mr^{-1}),
\eeaa
and hence, for $\vartheta_R(r)=\vartheta(r/R)$ a cut-off function such that $\vartheta(r)=1$ for $r\geq 1$ and $\vartheta(r)=0$ for $r\leq 1/2$, we infer
\beaa
\Div (r^{-\de_0-2}|\vartheta_R\pmb\phi|^2e_4) &=& 2r^{-\de_0-2}\vartheta_R^2\Re(\pmb\phi\c\ov{\nab_4\pmb\phi})+2r^{-\de_0-2}R^{-1}\vartheta_R\vartheta'(r/R)e_4(r)|\pmb\phi|^2\\
&&+r^{-\de_0-3}\Big((-\de_0-2)e_4(r)+r\Div ( e_4)\Big)\vartheta_R^2|\pmb\phi|^2\\
&=& 2r^{-\de_0-2}\vartheta_R^2\Re(\pmb\phi\c\ov{\pmb\psi})+r^{-\de_0-3}\big(-\de_0+O(mr^{-1})\big)\vartheta_R^2|\pmb\phi|^2\\
&& +2r^{-\de_0-2}R^{-1}\big(1+O(mr^{-1})\big)\vartheta_R\vartheta'(r/R)|\pmb\phi|^2.
\eeaa
Then, for $\de_0>0$ and $R\geq R(\de_0)$ large enough, we deduce 
\beaa
\frac{1}{2} r^{-\de_0-3}\vartheta_R^2|\pmb\phi|^2+{\frac{1}{\de_0}}\Div (r^{-\de_0-2}|\vartheta_R\pmb\phi|^2e_4) \lesssim \frac{1}{\de_0^2}r^{-\de_0-1}\vartheta_R^2|\pmb\psi|^2+{\frac{2}{\de_0}R^{-\de_0-3}|\vartheta'(r/R)| |\pmb\phi|^2}.
\eeaa
Integrating on $\MM(\tau_1,\tau_2)$ and using the support properties of $\vartheta_R$ and $\vartheta'$, we infer
\bea
\lab{eq:transporte4:general:zeroorder:lemmaprof}
&&\int_{\MM_{r\geq R}(\tau_1, \tau_2)}r^{-\de_0-3}|\pmb\phi|^2+\int_{\Si_{r\geq R}(\tau_2)}r^{-\de_0-4}|\pmb\phi|^2+\int_{\Si_*(\tau_1, \tau_2)}r^{-\de_0-2}|\pmb\phi|^2\nn\\
&\les_{\de_0}& \int_{\MM_{r\geq R/2}(\tau_1, \tau_2)}r^{-\de_0-1}|\pmb\psi|^2 +\int_{\Si_{r\geq R/2}(\tau_1)}r^{-\de_0-4}|\pmb\phi|^2+\int_{\MM_{R/2,R}(\tau_1, \tau_2)}r^{-\de_0-3}|\pmb\phi|^2,
\eea
 where we also used the fact that for $r\geq 13m$, we have
\beaa
\g(e_4, N_\Sigma)=-e_4(\tau)=-\frac{m^2}{r^2}+\Ga_g=-\frac{m^2}{r^2}(1+O(\ep))\simeq -\frac{m^2}{r^2},
\eeaa
and the fact that $\g(N_{\Si_*}, e_4)=-1+O(mr^{-1})$ on $\Si_*$. 

Next, commuting \eqref{eq:nab4phi=psi:transportestis} with $r\nab$ and using the commutator formulas in Lemma \ref{LEMMA:COMM-GEN-B}, we obtain 
\beaa
\nab_4 (r\nab\pmb\phi)=r\nab\pmb\psi +r\xi\nab_3\pmb\phi +\big(O(r^{-2})+\Ga_g\big)\dk^{\leq 1}\pmb\phi.
\eeaa
Hence, applying estimate \eqref{eq:transporte4:general:zeroorder:lemmaprof}  to this equation, we infer
\bea
\lab{eq:transporte4:general:firstorderangular:lemmaprof}
&&\int_{\MM_{r\geq R}(\tau_1, \tau_2)}r^{-\de_0-3}|r\nab\pmb\phi|^2+\int_{\Si_{r\geq R}(\tau_2)}r^{-\de_0-4}|r\nab\pmb\phi|^2+\int_{\Si_*(\tau_1, \tau_2)}r^{-\de_0-2}|r\nab\pmb\phi|^2\nn\\
&\les_{\de_0}& \int_{\MM_{r\geq R/2}(\tau_1, \tau_2)}r^{-\de_0-1}\Big|r\nab\pmb\psi +r\xi\nab_3\pmb\phi +\big(O(r^{-2})+\Ga_g\big)\dk^{\leq 1}\pmb\phi\Big|^2\nn\\
&& +\int_{\Si_{r\geq R/2}(\tau_1)}r^{-\de_0-4}|r\nab\pmb\phi|^2+\int_{\MM_{R/2,R}(\tau_1, \tau_2)}r^{-\de_0-3}|r\nab\pmb\phi|^2\nn\\
&\les_{\de_0}&\int_{\Si_{r\geq R/2}(\tau_1)}r^{-\de_0-4}|r\nab\pmb\phi|^2+\int_{\MM_{R/2,R}(\tau_1, \tau_2)}r^{-\de_0-3}|r\nab\pmb\phi|^2\nn\\
&&+ \int_{\MM_{r\geq R/2}(\tau_1, \tau_2)}r^{-\de_0-1}|r\nab\pmb\psi|^2 +R^{-2}\int_{\MM_{r\geq R/2}(\tau_1, \tau_2)}r^{-\de_0-3}|\dk^{\leq 1}\pmb\phi|^2.
\eea

Next, commuting \eqref{eq:nab4phi=psi:transportestis} with $\nab_3$ and using the commutator formulas in Lemma \ref{LEMMA:COMM-GEN-B}, we obtain 
\beaa
\nab_4(\nab_3\pmb\phi)&=& \nab_3\pmb\psi +O(r^{-2})\nab_3\pmb\phi +\big(O(r^{-3})+\Ga_g\big)\dk^{\leq 1}\pmb\phi.
\eeaa
Hence, applying estimate \eqref{eq:transporte4:general:zeroorder:lemmaprof}  to this equation, we infer
\bea
\lab{eq:transporte4:general:firstordere3:lemmaprof}
&&\int_{\MM_{r\geq R}(\tau_1, \tau_2)}r^{-\de_0-3}|\nab_{3}\pmb\phi|^2+\int_{\Si_{r\geq R}(\tau_2)}r^{-\de_0-4}|\nab_{3}\pmb\phi|^2+\int_{\Si_*(\tau_1, \tau_2)}r^{-\de_0-2}|\nab_{3}\pmb\phi|^2\nn\\
&\les_{\de_0}& \int_{\MM_{r\geq R/2}(\tau_1, \tau_2)}r^{-\de_0-1}\big|\nab_3\pmb\psi +O(r^{-2})\nab_3\pmb\phi +\big(O(r^{-3})+\Ga_g\big)\dk^{\leq 1}\pmb\phi\big|^2\nn\\
&& +\int_{\Si_{r\geq R/2}(\tau_1)}r^{-\de_0-4}|\nab_{3}\pmb\phi|^2 +\int_{\MM_{R/2,R}(\tau_1, \tau_2)}r^{-\de_0-3}|\nab_{3}\pmb\phi|^2\nn\\
&\les_{\de_0}&\int_{\Si_{r\geq R/2}(\tau_1)}r^{-\de_0-4}|\nab_{3}\pmb\phi|^2+\int_{\MM_{R/2,R}(\tau_1, \tau_2)}r^{-\de_0-3}|\nab_{3}\pmb\phi|^2\nn\\
&&+ \int_{\MM_{r\geq R/2}(\tau_1, \tau_2)}r^{-\de_0-1}|\nab_{3}\pmb\psi|^2 +R^{-2}\int_{\MM_{r\geq R/2}(\tau_1, \tau_2)}r^{-\de_0-3}|\dk^{\leq 1}\pmb\phi|^2.
\eea

Next, applying $\nab_4( r\c)$ to both sides of  the transport equation \eqref{eq:nab4phi=psi:transportestis}, we obtain  the following transport equation for $r\nab_4\pmb\phi$:
\bea
\nab_4 (r\nab_4\pmb\phi)&=&\nab_4(r\pmb\psi). 
\eea
Hence, applying estimate \eqref{eq:transporte4:general:zeroorder:lemmaprof}  to this equation, we infer
\begin{align}
\lab{eq:transporte4:general:firstordere4:lemmaprof}
&\int_{\MM_{r\geq R}(\tau_1, \tau_2)}r^{-\de_0-3}|r\nab_4\pmb\phi|^2+\int_{\Si_{r\geq R}(\tau_2)}r^{-\de_0-4}|r\nab_4\pmb\phi|^2+\int_{\Si_*(\tau_1, \tau_2)}r^{-\de_0-2}|r\nab_4\pmb\phi|^2\nn\\
\les_{\de_0}&\int_{\Si_{r\geq R/2}(\tau_1)}r^{-\de_0-4}|r\nab_4\pmb\phi|^2+\int_{\MM_{R/2,R}(\tau_1, \tau_2)}r^{-\de_0-3}|r\nab_4\pmb\phi|^2+ \int_{\MM_{r\geq R/2}(\tau_1, \tau_2)}r^{-\de_0-1}|\nab_4(r\pmb\psi)|^2.
\end{align}

Combining the estimates \eqref{eq:transporte4:general:zeroorder:lemmaprof}, \eqref{eq:transporte4:general:firstorderangular:lemmaprof}, \eqref{eq:transporte4:general:firstordere3:lemmaprof} and \eqref{eq:transporte4:general:firstordere4:lemmaprof}, and taking $R$ large enough so that the integral with an $R^{-2}$ prefactor is absorbed, we infer
\beaa
&&\int_{\MM_{r\geq R}(\tau_1, \tau_2)}r^{-\de_0-3}|\dk^{\leq 1}\pmb\phi|^2+\int_{\Si_{r\geq R}(\tau_2)}r^{-\de_0-4}|\dk^{\leq 1}\pmb\phi|^2+\int_{\Si_*(\tau_1, \tau_2)}r^{-\de_0-2}|\dk^{\leq 1}\pmb\phi|^2\nn\\
&\les_{\de_0}&\int_{\Si_{r\geq R/2}(\tau_1)}r^{-\de_0-4}|\dk^{\leq 1}\pmb\phi|^2+\int_{\MM_{R/2,R}(\tau_1, \tau_2)}r^{-\de_0-3}|\dk^{\leq 1}\pmb\phi|^2+ \int_{\MM_{r\geq R/2}(\tau_1, \tau_2)}r^{-\de_0-1}|\dk^{\leq 1}\pmb\psi|^2
\eeaa
which proves the desired estimate \eqref{eq:transporte4:general:anyordersummary:lemma} in the case $\reg=1$. The general $\reg\leq\kl$ case follows inductively by the same method.  This concludes the proof of Lemma \ref{lemma:transportlemmaforspinminus2:largerregion}.
\end{proof}

Next, we consider transport equations along $e_3$. 
\begin{lemma}\lab{lemma:transportlemmaforspinplus2:largerregion}
Assume that $\pmb\phi, \,\pmb\psi\in\sk_2(\mathbb{C})$ are such that 
\bea
\lab{eq:nab3phi=psi:transportestis}
\nab_3\pmb\phi=\pmb\psi.
\eea
 Then, for $\de_0>0$, $\reg\leq\kl$,  and for $R$ large enough, we have for any $1\leq\tau_1<\tau_2\leq\tau_*$
\bea
\lab{eq:transporte3:general:anyordersummary:lemma}
&&\int_{\MM_{r\geq R}(\tau_1, \tau_2)}r^{\de_0-3}|\dk^{\leq \reg}\pmb\phi|^2+\int_{\Si_{r\geq R}(\tau_2)}r^{\de_0-2}|\dk^{\leq \reg}\pmb\phi|^2+\int_{\Si_*(\tau_1, \tau_2)}r^{\de_0-2}|\dk^{\leq \reg}\pmb\phi|^2\nn\\
&\les_{\de_0}&\int_{\Si_{r\geq R/2}(\tau_1)}r^{\de_0-2}|\dk^{\leq \reg}\pmb\phi|^2+\int_{\MM_{R/2,R}(\tau_1, \tau_2)}r^{\de_0-3}|\dk^{\leq \reg}\pmb\phi|^2+ \int_{\MM_{r\geq R/2}(\tau_1, \tau_2)}r^{\de_0-1}|\dk^{\leq \reg}\pmb\psi|^2.
\eea
\end{lemma}

\begin{proof}
We compute 
\beaa
\Div(e_3) &=& \g^{43}\g(\D_4e_3, e_3) +\g^{43}\g(\D_3e_3, e_4)+\g^{bc}\g(\D_be_3, e_c)\\
&=& -\frac{1}{2}4\omb +\trchb = -\frac{2}{r} +O(m^2r^{-3})  -2\omb+\trchbc\\
&=& -\frac{2}{r}+O(mr^{-2}) +\Ga_b= -\frac{2}{r}+O(mr^{-2})+O(\ep r^{-1}) 
\eeaa
and 
\beaa
e_3(r)=-1+\widecheck{e_3(r)}=-1+r\Ga_b=-1+O(\ep),
\eeaa
and hence, for $\vartheta_R(r)=\vartheta(r/R)$ a cut-off function such that $\vartheta(r)=1$ for $r\geq 1$ and $\vartheta(r)=0$ for $r\leq 1/2$, we infer
\beaa
\Div (r^{\de_0-2}|\vartheta_R\pmb\phi|^2e_3) &=& 2r^{\de_0-2}\vartheta_R^2\Re(\pmb\phi\c\ov{\nab_3\pmb\phi})+2r^{\de_0-2}R^{-1}\vartheta_R\vartheta'(r/R)e_3(r)|\pmb\phi|^2\\
&&+r^{\de_0-3}\Big((\de_0-2)e_3(r)+r\Div ( e_3)\Big)\vartheta_R^2|\pmb\phi|^2\\
&=& 2r^{\de_0-2}\vartheta_R^2\Re(\pmb\phi\c\ov{\pmb\psi})+r^{\de_0-3}\big(-\de_0+O(mr^{-1})+O(\ep)\big)\vartheta_R^2|\pmb\phi|^2\\
&& +2r^{-\de_0-2}R^{-1}\big(-1+O(\ep)\big)\vartheta_R\vartheta'(r/R)|\pmb\phi|^2.
\eeaa
Then, for $\de_0>0$, $R\geq R(\de_0)$ large enough, and $0<\ep\ll\de_0$, we deduce 
\beaa
\frac{1}{2} r^{\de_0-3}\vartheta_R^2|\pmb\phi|^2+{\frac{1}{\de_0}}\Div (r^{\de_0-2}|\vartheta_R\pmb\phi|^2e_3) \lesssim \frac{1}{\de_0^2}r^{\de_0-1}\vartheta_R^2|\pmb\psi|^2+{\frac{2}{\de_0}R^{\de_0-3}|\vartheta'(r/R)| |\pmb\phi|^2}.
\eeaa
Integrating on $\MM(\tau_1,\tau_2)$ and using the support properties of $\vartheta_R$ and $\vartheta'$, we infer
\bea
\lab{eq:transporte3:general:zeroorder:lemmaprof}
&&\int_{\MM_{r\geq R}(\tau_1, \tau_2)}r^{\de_0-3}|\pmb\phi|^2+\int_{\Si_{r\geq R}(\tau_2)}r^{\de_0-2}|\pmb\phi|^2+\int_{\Si_*(\tau_1, \tau_2)}r^{\de_0-2}|\pmb\phi|^2\nn\\
&\les_{\de_0}& \int_{\MM_{r\geq R/2}(\tau_1, \tau_2)}r^{\de_0-1}|\pmb\psi|^2 +\int_{\Si_{r\geq R/2}(\tau_1)}r^{\de_0-2}|\pmb\phi|^2+\int_{\MM_{R/2,R}(\tau_1, \tau_2)}r^{\de_0-3}|\pmb\phi|^2,
\eea
 where we also used the fact that for $r\geq 13m$, we have
\beaa
\g(e_3, N_\Sigma)=-e_3(\tau)=-2+O(mr^{-1})+r\Ga_b=-2+O(mr^{-1})+O(\ep)\simeq -1,
\eeaa
and the fact that $\g(N_{\Si_*}, e_3)=-1+O(mr^{-1})$ on $\Si_*$. 

Next, commuting \eqref{eq:nab3phi=psi:transportestis} with $r\nab$ and using the commutator formulas in Lemma \ref{LEMMA:COMM-GEN-B}, we obtain 
\beaa
\nab_3(r\nab\pmb\phi) &=& r\nab\pmb\psi +r(\eta-\ze)\nab_3\pmb\phi +\big(O(r^{-2})+\Ga_b\big)\dk^{\leq 1}\pmb\phi\\
&=& r\nab\pmb\psi +\big(O(r^{-1})+r\Ga_b\big)\pmb\psi +\big(O(r^{-2})+\Ga_b\big)\dk^{\leq 1}\pmb\phi.
\eeaa
Hence, applying estimate \eqref{eq:transporte3:general:zeroorder:lemmaprof}  to this equation, we infer
\bea
\lab{eq:transporte3:general:firstorderangular:lemmaprof}
&&\int_{\MM_{r\geq R}(\tau_1, \tau_2)}r^{\de_0-3}|r\nab\pmb\phi|^2+\int_{\Si_{r\geq R}(\tau_2)}r^{\de_0-2}|r\nab\pmb\phi|^2+\int_{\Si_*(\tau_1, \tau_2)}r^{-\de_0-2}|r\nab\pmb\phi|^2\nn\\
&\les_{\de_0}& \int_{\MM_{r\geq R/2}(\tau_1, \tau_2)}r^{\de_0-1}\Big|r\nab\pmb\psi +\big(O(r^{-1})+r\Ga_b\big)\pmb\psi +\big(O(r^{-2})+\Ga_b\big)\dk^{\leq 1}\pmb\phi\Big|^2\nn\\
&& +\int_{\Si_{r\geq R/2}(\tau_1)}r^{\de_0-2}|r\nab\pmb\phi|^2+\int_{\MM_{R/2,R}(\tau_1, \tau_2)}r^{\de_0-3}|r\nab\pmb\phi|^2\nn\\
&\les_{\de_0}&\int_{\Si_{r\geq R/2}(\tau_1)}r^{\de_0-2}|r\nab\pmb\phi|^2+\int_{\MM_{R/2,R}(\tau_1, \tau_2)}r^{\de_0-3}|r\nab\pmb\phi|^2\nn\\
&&+ \int_{\MM_{r\geq R/2}(\tau_1, \tau_2)}r^{-\de_0-1}|(r\nab)^{\leq 1}\pmb\psi|^2 +(R^{-2}+\ep^2)\int_{\MM_{r\geq R/2}(\tau_1, \tau_2)}r^{-\de_0-3}|\dk^{\leq 1}\pmb\phi|^2.
\eea

Next, commuting \eqref{eq:nab3phi=psi:transportestis} with $\nab_4$ and using the commutator formulas in Lemma \ref{LEMMA:COMM-GEN-B}, we obtain 
\beaa
\nab_3(\nab_4\pmb\phi)&=& \nab_4\pmb\psi +O(r^{-2})\nab_3\pmb\phi +\big(O(r^{-3})+\Ga_g\big)\dk^{\leq 1}\pmb\phi\\
&=& \nab_4\pmb\psi +O(r^{-2})\pmb\psi +\big(O(r^{-3})+\Ga_g\big)\dk^{\leq 1}\pmb\phi.
\eeaa
Hence, applying estimate \eqref{eq:transporte3:general:zeroorder:lemmaprof}  to this equation with $\de_0\to \de_0+2$, we infer
\bea
\lab{eq:transporte3:general:firstordere4:lemmaprof}
&&\int_{\MM_{r\geq R}(\tau_1, \tau_2)}r^{\de_0-3}|r\nab_{4}\pmb\phi|^2+\int_{\Si_{r\geq R}(\tau_2)}r^{\de_0-2}|r\nab_{4}\pmb\phi|^2+\int_{\Si_*(\tau_1, \tau_2)}r^{\de_0-2}|r\nab_{4}\pmb\phi|^2\nn\\
&\les_{\de_0}& \int_{\MM_{r\geq R/2}(\tau_1, \tau_2)}r^{\de_0+1}\big|\nab_4\pmb\psi +O(r^{-2})\pmb\psi +\big(O(r^{-3})+\Ga_g\big)\dk^{\leq 1}\pmb\phi\big|^2\nn\\
&& +\int_{\Si_{r\geq R/2}(\tau_1)}r^{\de_0-2}|r\nab_{4}\pmb\phi|^2 +\int_{\MM_{R/2,R}(\tau_1, \tau_2)}r^{\de_0-3}|r\nab_{4}\pmb\phi|^2\nn\\
&\les_{\de_0}&\int_{\Si_{r\geq R/2}(\tau_1)}r^{\de_0-2}|r\nab_{4}\pmb\phi|^2+\int_{\MM_{R/2,R}(\tau_1, \tau_2)}r^{\de_0-3}|r\nab_{4}\pmb\phi|^2\nn\\
&&+ \int_{\MM_{r\geq R/2}(\tau_1, \tau_2)}r^{\de_0-1}|(r\nab_{4})^{\leq 1}\pmb\psi|^2 +(R^{-2}+\ep^2)\int_{\MM_{r\geq R/2}(\tau_1, \tau_2)}r^{\de_0-3}|\dk^{\leq 1}\pmb\phi|^2.
\eea

Next, applying $\nab_3$ to both sides of  the transport equation \eqref{eq:nab3phi=psi:transportestis}, we obtain $\nab_3(\nab_3\pmb\phi)=\nab_3\pmb\psi$ and applying estimate \eqref{eq:transporte3:general:zeroorder:lemmaprof}  to this equation, we infer
\bea
\lab{eq:transporte3:general:firstordere3:lemmaprof}
&&\int_{\MM_{r\geq R}(\tau_1, \tau_2)}r^{\de_0-3}|\nab_3\pmb\phi|^2+\int_{\Si_{r\geq R}(\tau_2)}r^{\de_0-2}|\nab_3\pmb\phi|^2+\int_{\Si_*(\tau_1, \tau_2)}r^{\de_0-2}|\nab_3\pmb\phi|^2\nn\\
&\les_{\de_0}&\int_{\Si_{r\geq R/2}(\tau_1)}r^{-\de_0-4}|r\nab_4\pmb\phi|^2+\int_{\MM_{R/2,R}(\tau_1, \tau_2)}r^{-\de_0-3}|r\nab_4\pmb\phi|^2\nn\\
&&+ \int_{\MM_{r\geq R/2}(\tau_1, \tau_2)}r^{\de_0-1}|\nab_3\pmb\psi|^2.
\eea

Combining the estimates \eqref{eq:transporte3:general:zeroorder:lemmaprof}, \eqref{eq:transporte3:general:firstorderangular:lemmaprof}, \eqref{eq:transporte3:general:firstordere4:lemmaprof} and \eqref{eq:transporte3:general:firstordere3:lemmaprof}, and taking $R$ large enough and $\ep>0$ small enough so that the integral with an $R^{-2}+\ep^2$ prefactor is absorbed, we infer
\beaa
&&\int_{\MM_{r\geq R}(\tau_1, \tau_2)}r^{\de_0-3}|\dk^{\leq 1}\pmb\phi|^2+\int_{\Si_{r\geq R}(\tau_2)}r^{\de_0-2}|\dk^{\leq 1}\pmb\phi|^2+\int_{\Si_*(\tau_1, \tau_2)}r^{\de_0-2}|\dk^{\leq 1}\pmb\phi|^2\nn\\
&\les_{\de_0}&\int_{\Si_{r\geq R/2}(\tau_1)}r^{\de_0-2}|\dk^{\leq 1}\pmb\phi|^2+\int_{\MM_{R/2,R}(\tau_1, \tau_2)}r^{\de_0-3}|\dk^{\leq 1}\pmb\phi|^2+ \int_{\MM_{r\geq R/2}(\tau_1, \tau_2)}r^{\de_0-1}|\dk^{\leq 1}\pmb\psi|^2
\eeaa
which proves the desired estimate \eqref{eq:transporte3:general:anyordersummary:lemma} in the case $\reg=1$. The general $\reg\leq\kl$ case follows inductively by the same method.  This concludes the proof of Lemma \ref{lemma:transportlemmaforspinplus2:largerregion}.
\end{proof}

%%%%%%%%%%%%%%%%%%%%%%%%%%%%%%%%%%%%%%%%%%%%%%

\subsubsection{$r^p$-weighted estimates for $\pmb\phi_{+2}^{(q)}$, $q=0,1,2$}

%%%%%%%%%%%%%%%%%%%%%%%%%%%%%%%%%%%%%%%%%%%%%%

We start by defining combined $r^p$-norms for the Teukolsky system in the case $s=+2$. 
\begin{definition}\lab{def:combinednormEpandBEFpTeukolskywavetransportsystem:plus2case}
Let $p$ be a real number such that $\de\leq p\leq 2-\de$ and $s\leq\kl-1$. We define the combined norms $\E_p^s[\pmb\phi_{+2}](\tau)$ and $\BEF_p^s[\pmb\phi_{+2}](\tau_1, \tau_2)$ as follows 
\beaa
\E_p^s[\pmb\phi_{+2}](\tau) &:=& \E_p^s[\pmb\phi_{+2}^{(2)}](\tau)+\int_{\Si(\tau)}r^p|\dk^{\leq s+1}\pmb\phi_{+2}^{(1)}|^2+\int_{\Si(\tau)}r^{p+2}|\dk^{\leq s+1}\pmb\phi_{+2}^{(0)}|^2,
\eeaa
and
\begin{align*}
\BEF_p^s[\pmb\phi_{+2}](\tau_1, \tau_2) :=&\sum_{q=0}^1\EMF^s_\de[\pmb\phi_{+2}^{(q)}](\tau_1, \tau_2)+\BEF_p^s[\pmb\phi_{+2}^{(2)}](\tau_1, \tau_2)+\int_{\MM_{r\geq 10m}(\tau_1, \tau_2)}r^{p-1}|\dk^{\leq s+1}\pmb\phi_{+2}^{(1)}|^2\\
&+\sup_{\tau\in[\tau_1, \tau_2]}\int_{\Si(\tau)}r^p|\dk^{\leq s+1}\pmb\phi_{+2}^{(1)}|^2+\int_{\Si_*(\tau_1, \tau_2)}r^p|\dk^{\leq s+1}\pmb\phi_{+2}^{(1)}|^2\\
&+\int_{\MM_{r\geq 10m}(\tau_1, \tau_2)}r^{p+1}|\dk^{\leq s+1}\pmb\phi_{+2}^{(0)}|^2+\sup_{\tau\in[\tau_1, \tau_2]}\int_{\Si(\tau)}r^{p+2}|\dk^{\leq s+1}\pmb\phi_{+2}^{(0)}|^2\\
&+\int_{\Si_*(\tau_1, \tau_2)}r^{p+2}|\dk^{\leq s+1}\pmb\phi_{+2}^{(0)}|^2.
\end{align*}
\end{definition}

\begin{proposition}\lab{prop:rpweightedestimatesforcombinednomrBEFsppmbphiplus2forallp}
Assume that $\phi_{+2}^{(q)}$, $q=0,1,2$, satisfies \eqref{eq:TensorialTeuSys:rescaleRHScontaine2:general:Kerrperturbation:alternateformnullframeinsteadcoordvectorfield} for $q=2$ and \eqref{eq:transportequationsins=+2caseforp=0and1}. Then, for $R=R(\de)\gg m$ larger enough, and for $s\leq\kl-1$ and $\de\leq p\leq 2-\de$, we have
\beaa
\nn\BEF^s_{p, \geq R}[\pmb\phi_{+2}](\tau_1, \tau_2) &\les& \E^s_{p,\geq \frac{R}{2}}[\pmb\phi_{+2}](\tau_1)+\NN^s_{p, \ge R/2}[\pmb\phi_{+2}^{(2)}, \widetilde{\N}_{W,+2}^{(2)}]{(\tau_1, \tau_2)} +\sum_{q=0}^2\M^s_{\frac{R}{2}, R}[\pmb\phi_{+2}^{(q)}](\tau_1, \tau_2)\\
&&+\int_{\MM_{r\geq R/2}(\tau_1, \tau_2)}r^{p-1}|\dk^{\leq s+1}\N_{T,+2}^{(1)}|^2+\int_{\MM_{r\geq R/2}(\tau_1, \tau_2)}r^{p+1}|\dk^{\leq s+1}\N_{T,+2}^{(0)}|^2,
\eeaa
where the combined norms $\BEF^s[\pmb\phi_{+2}](\tau_1, \tau_2)$ and $\E^s_{p}[\pmb\phi_{+2}](\tau_1)$ have are given by Definition \ref{def:combinednormEpandBEFpTeukolskywavetransportsystem:plus2case}.
\end{proposition}

\begin{proof}
In view of \eqref{eq:TensorialTeuSys:rescaleRHScontaine2:general:Kerrperturbation:alternateformnullframeinsteadcoordvectorfield}, we have
\beaa
\bigg(\squared_2 -\frac{4ia\cos\th}{|q|^2}\nab_{\T}- \frac{4}{|q|^2}\bigg)\pmb\phi_{+2}^{(2)} = \widetilde{\L}_{+2}^{(2)}[\pmb\phi_{+2}]+\widetilde{\N}_{W,+2}^{(2)}, 
\eeaa
which is of the form \eqref{eq:Gen.RW-chap10} with $N=\widetilde{\L}_{+2}^{(2)}[\pmb\phi_{+2}]+\widetilde{\N}_{W,+2}^{(2)}$. We may thus apply Proposition \ref{Proposition:Step3-Chap10} which yields, for $R\gg m$ large enough, $\de\leq p\leq 2-\de$ and $0\le s\le k_L-1$,
\bea\lab{eq:rpweigthedestimateRlargeforpmbphiplus2p=2:intermediarydirectapplication}
\nn\BEF^s_{p, \geq R}[\pmb\phi_{+2}^{(2)}](\tau_1, \tau_2) &\les& \E^s_{p,\geq \frac{R}{2}}[\pmb\phi_{+2}^{(2)}](\tau_1)+\NN^s_{p, \ge R/2}[\pmb\phi_{+2}^{(2)}, \widetilde{\N}_{W,+2}^{(2)}]{(\tau_1, \tau_2)} +\M^s_{\frac{R}{2}, R}[\pmb\phi_{+2}^{(2)}](\tau_1, \tau_2)\\
&&+\NN^s_{p, \ge R/2}[\pmb\phi_{+2}^{(2)}, \widetilde{\L}_{+2}^{(2)}[\pmb\phi_{+2}]](\tau_1, \tau_2).
\eea

Next, we estimate the last term on the RHS of \eqref{eq:rpweigthedestimateRlargeforpmbphiplus2p=2:intermediarydirectapplication}. Recall from \eqref{eq:tensor:Lsn:onlye_2present:general:Kerrperturbation:alternateformnullframeinsteadcoordvectorfield} that 
\beaa
{\widetilde{\L}_{+2}^{(2)}[\pmb\phi_{+2}]}&=&O(mr^{-3})\pmb\phi_{+2}^{(2)}+O(mr^{-2})\nab_{\Z+a\T}^{\leq 1}\pmb\phi_{+2}^{(1)}+O(m^2 r^{-2})\pmb\phi_{+2}^{(0)},
\eeaa
so that 
\beaa
&&\NN^s_{p, \ge R/2}[\pmb\phi_{+2}^{(2)}, \widetilde{\L}_{+2}^{(2)}[\pmb\phi_{+2}]](\tau_1, \tau_2)\\ 
&\les& \int_{\MM_{r\geq R/2}(\tau_1, \tau_2)}r^{p-1}|\nab_4(r\dk^{\leq s}\pmb\phi_{+2}^{(2)})|\Big(r^{-3}|\dk^{\leq s}\pmb\phi_{+2}^{(2)}|+r^{-2}|\dk^{\leq s+1}\pmb\phi_{+2}^{(1)}|+r^{-2}|\dk^{\leq s}\pmb\phi_{+2}^{(0)}|\Big)\\
&\les& R^{-1}\B^s_{p, \geq R}[\pmb\phi_{+2}^{(2)}](\tau_1, \tau_2)\\
&&+\Big(\B^s_{p, \geq R}[\pmb\phi_{+2}^{(2)}](\tau_1, \tau_2)\Big)^{\frac{1}{2}}\left(\int_{\MM_{r\geq R/2}(\tau_1, \tau_2)}r^{p-3}\Big(|\dk^{\leq s+1}\pmb\phi_{+2}^{(1)}|^2+|\dk^{\leq s}\pmb\phi_{+2}^{(0)}|^2\Big)\right)^{\frac{1}{2}}.
\eeaa
Hence, plugging into \eqref{eq:rpweigthedestimateRlargeforpmbphiplus2p=2:intermediarydirectapplication}, we infer, for $R$ large enough, 
\beaa
\nn\BEF^s_{p, \geq R}[\pmb\phi_{+2}^{(2)}](\tau_1, \tau_2) &\les& \E^s_{p,\geq \frac{R}{2}}[\pmb\phi_{+2}^{(2)}](\tau_1)+\NN^s_{p, \ge R/2}[\pmb\phi_{+2}^{(2)}, \widetilde{\N}_{W,+2}^{(2)}]{(\tau_1, \tau_2)} +\M^s_{\frac{R}{2}, R}[\pmb\phi_{+2}^{(2)}](\tau_1, \tau_2)\\
&&+\int_{\MM_{r\geq R/2}(\tau_1, \tau_2)}r^{p-3}\Big(|\dk^{\leq s+1}\pmb\phi_{+2}^{(1)}|^2+|\dk^{\leq s}\pmb\phi_{+2}^{(0)}|^2\Big).
\eeaa
In view of Definition \ref{def:combinednormEpandBEFpTeukolskywavetransportsystem:plus2case}, we deduce, for $R$ large enough, 
\bea\lab{eq:rpweigthedestimateRlargeforpmbphiplus2p=2:intermediarydirectapplication:bis}
\nn\BEF^s_{p, \geq R}[\pmb\phi_{+2}^{(2)}](\tau_1, \tau_2) &\les& \E^s_{p,\geq \frac{R}{2}}[\pmb\phi_{+2}^{(2)}](\tau_1)+\NN^s_{p, \ge R/2}[\pmb\phi_{+2}^{(2)}, \widetilde{\N}_{W,+2}^{(2)}]{(\tau_1, \tau_2)} +\M^s_{\frac{R}{2}, R}[\pmb\phi_{+2}^{(2)}](\tau_1, \tau_2)\\
&&+R^{-2}\BEF^s_{p, \geq R/2}[\pmb\phi_{+2}](\tau_1, \tau_2).
\eea

Next, we estimate $\pmb\phi_{+2}^{(q)}$, $q=0,1$. Recalling the transport equations \eqref{eq:transportequationsins=+2caseforp=0and1}, we apply Lemma \ref{lemma:transportlemmaforspinplus2:largerregion} for $q=1$ with $\de_0=p\geq\de$ and for $q=0$ with $\de_0=p+2$ which yields, for $R$ large enough,  
\beaa
&&\int_{\MM_{r\geq R}(\tau_1, \tau_2)}r^{p-1}|\dk^{\leq s+1}\pmb\phi_{+2}^{(1)}|^2+\int_{\Si_{r\geq R}(\tau_2)}r^{p}|\dk^{\leq s+1}\pmb\phi_{+2}^{(1)}|^2+\int_{\Si_*(\tau_1, \tau_2)}r^{p}|\dk^{\leq s+1}\pmb\phi_{+2}^{(1)}|^2\nn\\
&\les_{\de}&\int_{\Si_{r\geq R/2}(\tau_1)}r^{p}|\dk^{\leq s+1}\pmb\phi_{+2}^{(1)}|^2+\int_{\MM_{R/2,R}(\tau_1, \tau_2)}r^{p-1}|\dk^{\leq s+1}\pmb\phi_{+2}^{(1)}|^2\nn\\
&&+ \int_{\MM_{r\geq R/2}(\tau_1, \tau_2)}r^{p-3}|\dk^{\leq s+1}\pmb\phi_{+2}^{(2)}|^2+\int_{\MM_{r\geq R/2}(\tau_1, \tau_2)}r^{p-1}|\dk^{\leq s+1}\N_{T,+2}^{(1)}|^2,
\eeaa
and 
\beaa
&&\int_{\MM_{r\geq R}(\tau_1, \tau_2)}r^{p+1}|\dk^{\leq s+1}\pmb\phi_{+2}^{(0)}|^2+\int_{\Si_{r\geq R}(\tau_2)}r^{p+2}|\dk^{\leq s+1}\pmb\phi_{+2}^{(0)}|^2+\int_{\Si_*(\tau_1, \tau_2)}r^{p+2}|\dk^{\leq s+1}\pmb\phi_{+2}^{(0)}|^2\nn\\
&\les&\int_{\Si_{r\geq R/2}(\tau_1)}r^{p+2}|\dk^{\leq s+1}\pmb\phi_{+2}^{(0)}|^2+\int_{\MM_{R/2,R}(\tau_1, \tau_2)}r^{p+1}|\dk^{\leq s+1}\pmb\phi_{+2}^{(0)}|^2\nn\\
&&+ \int_{\MM_{r\geq R/2}(\tau_1, \tau_2)}r^{p-1}|\dk^{\leq s+1}\pmb\phi_{+2}^{(1)}|^2+\int_{\MM_{r\geq R/2}(\tau_1, \tau_2)}r^{p+1}|\dk^{\leq s+1}\N_{T,+2}^{(0)}|^2.
\eeaa
Adding these estimates to \eqref{eq:rpweigthedestimateRlargeforpmbphiplus2p=2:intermediarydirectapplication:bis}, we infer, in view of Definition \ref{def:combinednormEpandBEFpTeukolskywavetransportsystem:plus2case},
\beaa
\nn\BEF^s_{p, \geq R}[\pmb\phi_{+2}](\tau_1, \tau_2) &\les_\de& \E^s_{p,\geq \frac{R}{2}}[\pmb\phi_{+2}](\tau_1)+\NN^s_{p, \ge R/2}[\pmb\phi_{+2}^{(2)}, \widetilde{\N}_{W,+2}^{(2)}]{(\tau_1, \tau_2)} +\sum_{q=0}^2\M^s_{\frac{R}{2}, R}[\pmb\phi_{+2}^{(q)}](\tau_1, \tau_2)\\
&&+\int_{\MM_{r\geq R/2}(\tau_1, \tau_2)}r^{p-1}|\dk^{\leq s+1}\N_{T,+2}^{(1)}|^2+\int_{\MM_{r\geq R/2}(\tau_1, \tau_2)}r^{p+1}|\dk^{\leq s+1}\N_{T,+2}^{(0)}|^2\\
&&+R^{-2}\BEF^s_{p, \geq R/2}[\pmb\phi_{+2}](\tau_1, \tau_2),
\eeaa
and hence, for $R=R(\de)$ large enough, and for $s\leq\kl-1$ and $\de\leq p\leq 2-\de$,
\beaa
\nn\BEF^s_{p, \geq R}[\pmb\phi_{+2}](\tau_1, \tau_2) &\les& \E^s_{p,\geq \frac{R}{2}}[\pmb\phi_{+2}](\tau_1)+\NN^s_{p, \ge R/2}[\pmb\phi_{+2}^{(2)}, \widetilde{\N}_{W,+2}^{(2)}]{(\tau_1, \tau_2)} +\sum_{q=0}^2\M^s_{\frac{R}{2}, R}[\pmb\phi_{+2}^{(q)}](\tau_1, \tau_2)\\
&&+\int_{\MM_{r\geq R/2}(\tau_1, \tau_2)}r^{p-1}|\dk^{\leq s+1}\N_{T,+2}^{(1)}|^2+\int_{\MM_{r\geq R/2}(\tau_1, \tau_2)}r^{p+1}|\dk^{\leq s+1}\N_{T,+2}^{(0)}|^2
\eeaa
as stated. This concludes the proof of Proposition \ref{prop:rpweightedestimatesforcombinednomrBEFsppmbphiplus2forallp}.
\end{proof}

We now have the following corollary of Proposition \ref{prop:rpweightedestimatesforcombinednomrBEFsppmbphiplus2forallp} and Theorem \ref{thm:main:MaSz26:extendhigherorderderivatives}.
\begin{corollary}\lab{cor:rpweightedestimatesforcombinednomrBEFsppmbphiplus2forallp}
Let $(\MM, \g)$ satisfy the assumptions of Section \ref{sec:geometricsetupfordecsayTeukolsky}. Then, for $\ep>0$ small enough, we have for solutions $\pmb\phi_{+2}^{(q)}$, $q=0,1,2$, to the Teukolsky wave/transport system \eqref{eq:TensorialTeuSysandlinearterms:rescaleRHScontaine2:general:Kerrperturbation:alternateformnullframeinsteadcoordvectorfield} \eqref{eq:transportequationsins=plus2andminus2caseforp=0and1} the following $r^p$-weighted estimates, for $11\leq s\leq \kst-3$ and $\de\leq p\leq 2-\de$, 
\beaa
\nn\BEF^s_{p}[\pmb\phi_{+2}](\tau_1, \tau_2) &\les& \E^s_{p}[\pmb\phi_{+2}](\tau_1)+\NN^s_{p}[\pmb\phi_{+2}^{(2)}, \widetilde{\N}_{W,+2}^{(2)}]{(\tau_1, \tau_2)} +\widetilde{\NN}^s_\de[\pmb\phi_{+2}^{(0)}, \widetilde{\N}_{W,+2}^{(0)}](\tt_1, \tt_2)\\
&&+\widetilde{\NN}^s_\de[\pmb\phi_{+2}^{(1)}, \widetilde{\N}_{W,+2}^{(1)}](\tt_1, \tt_2).
\eeaa
\end{corollary}

\begin{proof}
In view of Proposition \ref{prop:rpweightedestimatesforcombinednomrBEFsppmbphiplus2forallp}, \eqref{MainEnerMora:psi:plus2case:extendhigherorderderivatives} and Definition \ref{def:combinednormEpandBEFpTeukolskywavetransportsystem:plus2case}, we have, for $11\leq s\leq \kst-3$ and $\de\leq p\leq 2-\de$, 
\beaa
\nn\BEF^s_{p}[\pmb\phi_{+2}](\tau_1, \tau_2) &\les& \E^s_{p}[\pmb\phi_{+2}](\tau_1)+\NN^s_{p}[\pmb\phi_{+2}^{(2)}, \widetilde{\N}_{W,+2}^{(2)}]{(\tau_1, \tau_2)} +\widetilde{\NN}^s_\de[\pmb\phi_{+2}^{(0)}, \widetilde{\N}_{W,+2}^{(0)}](\tt_1, \tt_2)\\
&&+\widetilde{\NN}^{s}_\de[\pmb\phi_{+2}^{(1)}, \widetilde{\N}_{W,+2}^{(1)}](\tt_1, \tt_2)+\int_{\MM(\tau_1, \tau_2)}r^{p-1}|\dk^{\leq s+1}\N_{T,+2}^{(1)}|^2\\
&&+\int_{\MM(\tau_1, \tau_2)}r^{p+1}|\dk^{\leq s+1}\N_{T,+2}^{(0)}|^2 +\sum_{q=0}^1\int_{\MM(\tt_1, \tt_2)}r^{1+\de}|\dk^{\leq s+1}\big(\Ga_b\c\pmb\phi_{+2}^{(q)}\big)|^2\\
&&+\ep\B^s_{\de}[\pmb\phi_{+2}](\tau_1, \tau_2).
\eeaa
In particular, since $p\geq\de$, we infer, for $\ep>0$ small enough, for $11\leq s\leq \kst-3$ and $\de\leq p\leq 2-\de$, 
\begin{align}\lab{eq:intermediarycontrolofBEFsppmbphiplus2combinednormwhere3termsonRHSstillneedtobeestimated}
\nn\BEF^s_{p}[\pmb\phi_{+2}](\tau_1, \tau_2) \les& \E^s_{p}[\pmb\phi_{+2}](\tau_1)+\NN^s_{p}[\pmb\phi_{+2}^{(2)}, \widetilde{\N}_{W,+2}^{(2)}]{(\tau_1, \tau_2)} +\widetilde{\NN}^s_\de[\pmb\phi_{+2}^{(0)}, \widetilde{\N}_{W,+2}^{(0)}](\tt_1, \tt_2)\\
\nn&+\widetilde{\NN}^s_\de[\pmb\phi_{+2}^{(1)}, \widetilde{\N}_{W,+2}^{(1)}](\tt_1, \tt_2)+\int_{\MM(\tau_1, \tau_2)}r^{p-1}|\dk^{\leq s+1}\N_{T,+2}^{(1)}|^2\\
&+\int_{\MM(\tau_1, \tau_2)}r^{p+1}|\dk^{\leq s+1}\N_{T,+2}^{(0)}|^2 +\sum_{q=0}^1\int_{\MM(\tt_1, \tt_2)}r^{1+\de}|\dk^{\leq s+1}\big(\Ga_b\c\pmb\phi_{+2}^{(q)}\big)|^2.
\end{align}

Next, we estimate the last three terms on the RHS of \eqref{eq:intermediarycontrolofBEFsppmbphiplus2combinednormwhere3termsonRHSstillneedtobeestimated}. Note from \eqref{eq:transportequationsins=+2caseforp=0and1:RHSNT+2p=0and1} and \eqref{eq:definitionofthephiplus2phierarchy:perturbationofKerr} that 
\beaa
\N_{T,+2}^{(0)}=r\Ga_b\pmb\phi_{+2}^{(0)}, \qquad \N_{T,+2}^{(1)}=r\Ga_b\pmb\phi_{+2}^{(1)}+r^2\dk^{\leq 1}(\Ga_b)\pmb\phi_{+2}^{(0)},
\eeaa
and hence, for $s\leq \kst-3$ and $\de\leq p\leq 2-\de$, 
\begin{align*}
\nn& \int_{\MM(\tau_1, \tau_2)}r^{p-1}|\dk^{\leq s+1}\N_{T,+2}^{(1)}|^2+\int_{\MM(\tau_1, \tau_2)}r^{p+1}|\dk^{\leq s+1}\N_{T,+2}^{(0)}|^2 +\sum_{q=0}^1\int_{\MM(\tt_1, \tt_2)}r^{1+\de}|\dk^{\leq s+1}\big(\Ga_b\c\pmb\phi_{+2}^{(q)}\big)|^2\\
\nn\les& \ep^2\int_{\MM(\tau_1, \tau_2)}r^{p-1}|\dk^{\leq s+1}\pmb\phi_{+2}^{(1)}|^2+\ep^2\int_{\MM(\tau_1, \tau_2)}r^{p+1}|\dk^{\leq s+1}\pmb\phi_{+2}^{(0)}|^2\\
\les& \ep^2\BEF^s_{p}[\pmb\phi_{+2}](\tau_1, \tau_2),
\end{align*}
where we used Definition \ref{def:combinednormEpandBEFpTeukolskywavetransportsystem:plus2case} in the last inequality. Plugging in \eqref{eq:intermediarycontrolofBEFsppmbphiplus2combinednormwhere3termsonRHSstillneedtobeestimated}, we infer, for $\ep>0$ small enough, and for $11\leq s\leq \kst-3$ and $\de\leq p\leq 2-\de$, 
\beaa
\nn\BEF^s_{p}[\pmb\phi_{+2}](\tau_1, \tau_2) &\les& \E^s_{p}[\pmb\phi_{+2}](\tau_1)+\NN^s_{p}[\pmb\phi_{+2}^{(2)}, \widetilde{\N}_{W,+2}^{(2)}]{(\tau_1, \tau_2)} +\widetilde{\NN}^s_\de[\pmb\phi_{+2}^{(0)}, \widetilde{\N}_{W,+2}^{(0)}](\tt_1, \tt_2)\\
&&+\widetilde{\NN}^s_\de[\pmb\phi_{+2}^{(1)}, \widetilde{\N}_{W,+2}^{(1)}](\tt_1, \tt_2),
\eeaa
as stated. This concludes the proof of Corollary \ref{cor:rpweightedestimatesforcombinednomrBEFsppmbphiplus2forallp}.
\end{proof}

%%%%%%%%%%%%%%%%%%%%%%%%%%%%%%%%%%%%%%%%%%%%%%

\subsubsection{$r^p$-weighted estimates for $\pmb\phi_{-2}^{(q)}$, $q=0,1,2$}

%%%%%%%%%%%%%%%%%%%%%%%%%%%%%%%%%%%%%%%%%%%%%%

We start by defining combined $r^p$-norms for the Teukolsky system in the case $s=-2$. 
\begin{definition}\lab{def:combinednormEpandBEFpTeukolskywavetransportsystem:minus2case}
Let $p$ be a real number such that $\de\leq p\leq 2-\de$ and $s\leq\kl-1$. We define the combined norms $\E_p^s[\pmb\phi_{-2}](\tau)$ and $\BEF_p^s[\pmb\phi_{-2}](\tau_1, \tau_2)$ as follows 
\beaa
\E_p^s[\pmb\phi_{-2}](\tau) &:=& \E_p^s[\pmb\phi_{-2}^{(2)}](\tau)+\int_{\Si(\tau)}r^p|\dk^{\leq s+1}\pmb\phi_{+2}^{(1)}|^2+\int_{\Si(\tau)}r^{p+2}|\dk^{\leq s+1}\pmb\phi_{+2}^{(0)}|^2+\sum_{q=0}^2\E^s_{r\leq r_+(1+\dred)}[\nab_4^q\Ab](\tau)\\
&&+\sum_{q=0,1}\int_{\Si(\tau)}r^{p-4}|\dk^{\leq s+1}\pmb\phi_{-2}^{(q)}|^2,
\eeaa
and
\begin{align*}
\BEF_p^s[\pmb\phi_{-2}](\tau_1, \tau_2) :=&\sum_{q=0}^1\EMF^s_\de[\pmb\phi_{-2}^{(q)}](\tau_1, \tau_2)+\BEF_p^s[\pmb\phi_{-2}^{(2)}](\tau_1, \tau_2)+\sum_{q=0}^2\EMF^s_{r\leq r_+(1+\dred)}[\nab_4^q\Ab](\tau_1, \tau_2)\\
&+\sum_{q=0,1}\int_{\MM_{r\geq 10m}(\tau_1, \tau_2)}r^{p-3}|\dk^{\leq s+1}\pmb\phi_{-2}^{(q)}|^2+\sum_{q=0,1}\sup_{\tau\in[\tau_1, \tau_2]}\int_{\Si_{r\geq 10m}(\tau)}r^{p-4}|\dk^{\leq s+1}\pmb\phi_{-2}^{(q)}|^2\\
&+\sum_{q=0,1}\int_{\Si_*(\tau_1, \tau_2)}r^{p-2}|\dk^{\leq s+1}\pmb\phi_{-2}^{(q)}|^2.
\end{align*}
\end{definition} 

\begin{proposition}\lab{prop:rpweightedestimatesforcombinednomrBEFsppmbphiminus2forallp}
Assume that $\phi_{-2}^{(q)}$, $q=0,1,2$, satisfies \eqref{eq:TensorialTeuSys:rescaleRHScontaine2:general:Kerrperturbation:alternateformnullframeinsteadcoordvectorfield} for $q=2$ and \eqref{eq:transportequationsins=-2caseforp=0and1}. Then, for $R=R(\de)\gg m$ larger enough, and for $s\leq\kl-1$ and $\de\leq p\leq 2-\de$, we have
\beaa
\nn\BEF^s_{p, \geq R}[\pmb\phi_{-2}](\tau_1, \tau_2) &\les& \E^s_{p,\geq \frac{R}{2}}[\pmb\phi_{-2}](\tau_1)+\NN^s_{p, \ge R/2}[\pmb\phi_{-2}^{(2)}, \widetilde{\N}_{W,-2}^{(2)}]{(\tau_1, \tau_2)} +\sum_{q=0}^2\M^s_{\frac{R}{2}, R}[\pmb\phi_{-2}^{(q)}](\tau_1, \tau_2)\\
&&+\sum_{q=0}^1\int_{\MM_{r\geq R/2}(\tau_1, \tau_2)}r^{p-3}|\dk^{\leq s+1}\N_{T,-2}^{(q)}|^2,
\eeaa
where the combined norms $\BEF^s[\pmb\phi_{-2}](\tau_1, \tau_2)$ and $\E^s_{p}[\pmb\phi_{-2}](\tau_1)$ have are given by Definition \ref{def:combinednormEpandBEFpTeukolskywavetransportsystem:minus2case}.
\end{proposition}

\begin{proof}
In view of \eqref{eq:TensorialTeuSys:rescaleRHScontaine2:general:Kerrperturbation:alternateformnullframeinsteadcoordvectorfield}, we have
\beaa
\bigg(\squared_2 -\frac{4ia\cos\th}{|q|^2}\nab_{\T}- \frac{4}{|q|^2}\bigg)\pmb\phi_{-2}^{(2)} = \widetilde{\L}_{-2}^{(2)}[\pmb\phi_{-2}]+\widetilde{\N}_{W,-2}^{(2)}, 
\eeaa
which is of the form \eqref{eq:Gen.RW-chap10} with $N=\widetilde{\L}_{-2}^{(2)}[\pmb\phi_{-2}]+\widetilde{\N}_{W,-2}^{(2)}$. We may thus apply Proposition \ref{Proposition:Step3-Chap10} which yields, for $R\gg m$ large enough, $\de\leq p\leq 2-\de$ and $0\le s\le k_L-1$,
\bea\lab{eq:rpweigthedestimateRlargeforpmbphiminus2p=2:intermediarydirectapplication}
\nn\BEF^s_{p, \geq R}[\pmb\phi_{-2}^{(2)}](\tau_1, \tau_2) &\les& \E^s_{p,\geq \frac{R}{2}}[\pmb\phi_{-2}^{(2)}](\tau_1)+\NN^s_{p, \ge R/2}[\pmb\phi_{-2}^{(2)}, \widetilde{\N}_{W,-2}^{(2)}]{(\tau_1, \tau_2)} +\M^s_{\frac{R}{2}, R}[\pmb\phi_{-2}^{(2)}](\tau_1, \tau_2)\\
&&+\NN^s_{p, \ge R/2}[\pmb\phi_{-2}^{(2)}, \widetilde{\L}_{-2}^{(2)}[\pmb\phi_{-2}]](\tau_1, \tau_2).
\eea

Next, we estimate the last term on the RHS of \eqref{eq:rpweigthedestimateRlargeforpmbphiminus2p=2:intermediarydirectapplication}. Recall from \eqref{eq:tensor:Lsn:onlye_2present:general:Kerrperturbation:alternateformnullframeinsteadcoordvectorfield} that 
\beaa
{\widetilde{\L}_{-2}^{(2)}[\pmb\phi_{-2}]}&=&O(mr^{-3})\pmb\phi_{-2}^{(2)}+O(mr^{-2})\nab_{\Z+a\T}^{\leq 1}\pmb\phi_{-2}^{(1)}+O(m^2 r^{-2})\pmb\phi_{-2}^{(0)},
\eeaa
so that 
\begin{align*}
&\NN^s_{p, \ge R/2}[\pmb\phi_{-2}^{(2)}, \widetilde{\L}_{-2}^{(2)}[\pmb\phi_{-2}]](\tau_1, \tau_2)\\ 
\les& \int_{\MM_{r\geq R/2}(\tau_1, \tau_2)}r^{p-1}|\nab_4(r\dk^{\leq s}\pmb\phi_{-2}^{(2)})|\Big(r^{-3}|\dk^{\leq s}\pmb\phi_{-2}^{(2)}|+r^{-2}|\dk^{\leq s+1}\pmb\phi_{-2}^{(1)}|+r^{-2}|\dk^{\leq s}\pmb\phi_{-2}^{(0)}|\Big)\\
\les& R^{-1}\B^s_{p, \geq R}[\pmb\phi_{-2}^{(2)}](\tau_1, \tau_2)+\Big(\B^s_{p, \geq R}[\pmb\phi_{-2}^{(2)}](\tau_1, \tau_2)\Big)^{\frac{1}{2}}\left(\int_{\MM_{r\geq R/2}(\tau_1, \tau_2)}r^{p-3}\Big(|\dk^{\leq s+1}\pmb\phi_{-2}^{(1)}|^2+|\dk^{\leq s}\pmb\phi_{-2}^{(0)}|^2\Big)\right)^{\frac{1}{2}}.
\end{align*}
Hence, plugging into \eqref{eq:rpweigthedestimateRlargeforpmbphiminus2p=2:intermediarydirectapplication}, we infer, for $R$ large enough, 
\bea\lab{eq:rpweigthedestimateRlargeforpmbphiminus2p=2:intermediarydirectapplication:bis}
\nn\BEF^s_{p, \geq R}[\pmb\phi_{-2}^{(2)}](\tau_1, \tau_2) &\les& \E^s_{p,\geq \frac{R}{2}}[\pmb\phi_{-2}^{(2)}](\tau_1)+\NN^s_{p, \ge R/2}[\pmb\phi_{-2}^{(2)}, \widetilde{\N}_{W,-2}^{(2)}]{(\tau_1, \tau_2)} +\M^s_{\frac{R}{2}, R}[\pmb\phi_{-2}^{(2)}](\tau_1, \tau_2)\\
&&+\int_{\MM_{r\geq R/2}(\tau_1, \tau_2)}r^{p-3}\Big(|\dk^{\leq s+1}\pmb\phi_{-2}^{(1)}|^2+|\dk^{\leq s}\pmb\phi_{-2}^{(0)}|^2\Big).
\eea

Next, we estimate $\pmb\phi_{-2}^{(q)}$, $q=0,1$. Recalling the transport equations \eqref{eq:transportequationsins=-2caseforp=0and1}, we apply Lemma \ref{lemma:transportlemmaforspinminus2:largerregion} with $\de_0=2-p\geq\de$ which yields, for $R$ large enough and $q=0,1$, 
\bea\lab{eq:transportestimateforpmbphiminus2p=0and1:usefullhereaswellaslater}
&&\int_{\MM_{r\geq R}(\tau_1, \tau_2)}r^{p-3}|\dk^{\leq s+1}\pmb\phi_{-2}^{(q)}|^2+\int_{\Si_{r\geq R}(\tau_2)}r^{p-4}|\dk^{\leq s+1}\pmb\phi_{-2}^{(q)}|^2+\int_{\Si_*(\tau_1, \tau_2)}r^{p-2}|\dk^{\leq s+1}\pmb\phi_{-2}^{(q)}|^2\nn\\
&\les_{\de}&\int_{\Si_{r\geq R/2}(\tau_1)}r^{p-4}|\dk^{\leq s+1}\pmb\phi_{-2}^{(q)}|^2+\int_{\MM_{R/2,R}(\tau_1, \tau_2)}r^{p-3}|\dk^{\leq s+1}\pmb\phi_{-2}^{(q)}|^2\nn\\
&&+ \int_{\MM_{r\geq R/2}(\tau_1, \tau_2)}r^{p-5}|\dk^{\leq s+1}\pmb\phi_{-2}^{(q+1)}|^2+\int_{\MM_{r\geq R/2}(\tau_1, \tau_2)}r^{p-3}|\dk^{\leq s+1}\N_{T,-2}^{(q)}|^2.
\eea
Adding these estimates to \eqref{eq:rpweigthedestimateRlargeforpmbphiminus2p=2:intermediarydirectapplication:bis}, we infer, in view of Definition \ref{def:combinednormEpandBEFpTeukolskywavetransportsystem:minus2case},
\beaa
\nn\BEF^s_{p, \geq R}[\pmb\phi_{-2}](\tau_1, \tau_2) &\les_\de& \E^s_{p,\geq \frac{R}{2}}[\pmb\phi_{-2}](\tau_1)+\NN^s_{p, \ge R/2}[\pmb\phi_{-2}^{(2)}, \widetilde{\N}_{W,-2}^{(2)}]{(\tau_1, \tau_2)} +\sum_{q=0}^2\M^s_{\frac{R}{2}, R}[\pmb\phi_{-2}^{(q)}](\tau_1, \tau_2)\\
&&+\sum_{q=0}^1\int_{\MM_{r\geq R/2}(\tau_1, \tau_2)}r^{p-3}|\dk^{\leq s+1}\N_{T,-2}^{(q)}|^2\\
&&+\sum_{q=1}^2R^{-2}\int_{\MM_{r\geq R/2}(\tau_1, \tau_2)}r^{p-3}|\dk^{\leq s+1}\pmb\phi_{-2}^{(q)}|^2\\
&\les_\de&  \E^s_{p,\geq \frac{R}{2}}[\pmb\phi_{-2}](\tau_1)+\NN^s_{p, \ge R/2}[\pmb\phi_{-2}^{(2)}, \widetilde{\N}_{W,-2}^{(2)}]{(\tau_1, \tau_2)} +\sum_{q=0}^2\M^s_{\frac{R}{2}, R}[\pmb\phi_{-2}^{(q)}](\tau_1, \tau_2)\\
&&+\sum_{q=0}^1\int_{\MM_{r\geq R/2}(\tau_1, \tau_2)}r^{p-3}|\dk^{\leq s+1}\N_{T,-2}^{(q)}|^2+R^{-2}\BEF^s_{p, \geq R}[\pmb\phi_{-2}](\tau_1, \tau_2)
\eeaa
and hence, for $R=R(\de)$ large enough, and for $s\leq\kl-1$ and $\de\leq p\leq 2-\de$,
\beaa
\nn\BEF^s_{p, \geq R}[\pmb\phi_{-2}](\tau_1, \tau_2) &\les& \E^s_{p,\geq \frac{R}{2}}[\pmb\phi_{-2}](\tau_1)+\NN^s_{p, \ge R/2}[\pmb\phi_{-2}^{(2)}, \widetilde{\N}_{W,-2}^{(2)}]{(\tau_1, \tau_2)} +\sum_{q=0}^2\M^s_{\frac{R}{2}, R}[\pmb\phi_{-2}^{(q)}](\tau_1, \tau_2)\\
&&+\sum_{q=0}^1\int_{\MM_{r\geq R/2}(\tau_1, \tau_2)}r^{p-3}|\dk^{\leq s+1}\N_{T,-2}^{(q)}|^2
\eeaa
as stated. This concludes the proof of Proposition \ref{prop:rpweightedestimatesforcombinednomrBEFsppmbphiminus2forallp}.
\end{proof}

We now have the following corollary of Proposition \ref{prop:rpweightedestimatesforcombinednomrBEFsppmbphiminus2forallp} and Theorem \ref{thm:main:MaSz26:extendhigherorderderivatives}.
\begin{corollary}\lab{cor:rpweightedestimatesforcombinednomrBEFsppmbphiminus2forallp}
Let $(\MM, \g)$ satisfy the assumptions of Section \ref{sec:geometricsetupfordecsayTeukolsky}. Then, for $\ep>0$ small enough, we have for solutions $\pmb\phi_{-2}^{(q)}$, $q=0,1,2$, to the Teukolsky wave/transport system \eqref{eq:TensorialTeuSysandlinearterms:rescaleRHScontaine2:general:Kerrperturbation:alternateformnullframeinsteadcoordvectorfield} \eqref{eq:transportequationsins=plus2andminus2caseforp=0and1} the following $r^p$-weighted estimates, for $14\leq s\leq \kst-3$ and $\de\leq p\leq 2-\de$, 
\beaa
\nn\BEF^s_{p}[\pmb\phi_{-2}](\tau_1, \tau_2) &\les& \E^s_{p}[\pmb\phi_{-2}](\tau_1)+\NN^s_{p}[\pmb\phi_{-2}^{(2)}, \widetilde{\N}_{W,-2}^{(2)}]{(\tau_1, \tau_2)} +\widetilde{\NN}^s_\de[\pmb\phi_{-2}^{(0)}, \widetilde{\N}_{W,-2}^{(0)}](\tt_1, \tt_2)\\
\nn&&+\widetilde{\NN}^{s}_\de[\pmb\phi_{-2}^{(1)}, \widetilde{\N}_{W,-2}^{(1)}](\tt_1, \tt_2)+\frac{\ep^2}{\tau_1^{2+\frac{3}{2}\dec}}\int_{\MM(\tau_1, \tau_2)}r^{1+\de}|\dk^{\leq s+1}\Xi|^2+\frac{\ep^2_0}{\tau_1^{2+3\dec}}.
\eeaa
\end{corollary}

\begin{proof}
In view of Proposition \ref{prop:rpweightedestimatesforcombinednomrBEFsppmbphiminus2forallp}, \eqref{MainEnerMora:psi:minus2case:extendhigherorderderivatives} and Definition \ref{def:combinednormEpandBEFpTeukolskywavetransportsystem:minus2case}, we have, for $14\leq s\leq \kst-3$ and $\de\leq p\leq 2-\de$, noticing that $p-3\leq -1-\de<-1+\de$,
\beaa
\nn\BEF^s_{p}[\pmb\phi_{-2}](\tau_1, \tau_2) &\les& \E^s_{p}[\pmb\phi_{-2}](\tau_1)+\NN^s_{p}[\pmb\phi_{-2}^{(2)}, \widetilde{\N}_{W,-2}^{(2)}]{(\tau_1, \tau_2)} +\widetilde{\NN}^s_\de[\pmb\phi_{-2}^{(0)}, \widetilde{\N}_{W,-2}^{(0)}](\tt_1, \tt_2)\\
&&+\widetilde{\NN}^{s}_\de[\pmb\phi_{-2}^{(1)}, \widetilde{\N}_{W,-2}^{(1)}](\tt_1, \tt_2)+\sum_{q=0}^1\int_{\MM(\tau_1, \tau_2)}r^{-1+\de}|\dk^{\leq s+1}\N_{T,-2}^{(q)}|^2\\
&& +\int_{\Si(\tau_1)}|\dk^{\leq s-2}\N_{T,-2}^{(0)}|^2+\ep\B^s_{\de}[\pmb\phi_{-2}](\tau_1, \tau_2)+ \frac{\ep^2_0}{\tau_1^{3+3\dec}}.
\eeaa
In particular, since $p\geq\de$, we infer, for $\ep>0$ small enough, for $14\leq s\leq \kst-3$ and $\de\leq p\leq 2-\de$, 
\bea\lab{eq:intermediarycontrolofBEFsppmbphiminus2combinednormwhere3termsonRHSstillneedtobeestimated}
\nn\BEF^s_{p}[\pmb\phi_{-2}](\tau_1, \tau_2) &\les& \E^s_{p}[\pmb\phi_{-2}](\tau_1)+\NN^s_{p}[\pmb\phi_{-2}^{(2)}, \widetilde{\N}_{W,-2}^{(2)}]{(\tau_1, \tau_2)} +\widetilde{\NN}^s_\de[\pmb\phi_{-2}^{(0)}, \widetilde{\N}_{W,-2}^{(0)}](\tt_1, \tt_2)\\
\nn&&+\widetilde{\NN}^{s}_\de[\pmb\phi_{-2}^{(1)}, \widetilde{\N}_{W,-2}^{(1)}](\tt_1, \tt_2)+\sum_{q=0}^1\int_{\MM(\tau_1, \tau_2)}r^{-1+\de}|\dk^{\leq s+1}\N_{T,-2}^{(q)}|^2\\
&& +\int_{\Si(\tau_1)}|\dk^{\leq s-2}\N_{T,-2}^{(0)}|^2+ \frac{\ep^2_0}{\tau_1^{3+3\dec}}.
\eea

Next, we estimate the terms involving $\N_{T,-2}^{(q)}$, $q=0,1$, on the RHS of \eqref{eq:intermediarycontrolofBEFsppmbphiminus2combinednormwhere3termsonRHSstillneedtobeestimated}. Note from \eqref{eq:transportequationsins=-2caseforp=0and1:RHSNT-2p=0and1} that, for $s\leq \kst-3$,
\bea\lab{eq:estimatesforNTminus2qq=01appearinginrpweightedestimates:usefullforlater}
\nn&& \sum_{q=0}^1\int_{\MM(\tau_1, \tau_2)}r^{-1+\de}|\dk^{\leq s+1}\N_{T,-2}^{(q)}|^2 +\int_{\Si(\tau_1)}|\dk^{\leq s-2}\N_{T,-2}^{(0)}|^2\\
&\les& \frac{\ep^2}{\tau_1^{2+\frac{3}{2}\dec}}\int_{\MM(\tau_1, \tau_2)}r^{1+\de}|\dk^{\leq s+1}\Xi|^2+\frac{\ep^2_0}{\tau_1^{2+3\dec}}
\eea
and hence, plugging in \eqref{eq:intermediarycontrolofBEFsppmbphiminus2combinednormwhere3termsonRHSstillneedtobeestimated}, we infer, for $\ep>0$ small enough, and for $14\leq s\leq \kst-3$ and $\de\leq p\leq 2-\de$, 
\beaa
\nn\BEF^s_{p}[\pmb\phi_{-2}](\tau_1, \tau_2) &\les& \E^s_{p}[\pmb\phi_{-2}](\tau_1)+\NN^s_{p}[\pmb\phi_{-2}^{(2)}, \widetilde{\N}_{W,-2}^{(2)}]{(\tau_1, \tau_2)} +\widetilde{\NN}^s_\de[\pmb\phi_{-2}^{(0)}, \widetilde{\N}_{W,-2}^{(0)}](\tt_1, \tt_2)\\
\nn&&+\widetilde{\NN}^{s}_\de[\pmb\phi_{-2}^{(1)}, \widetilde{\N}_{W,-2}^{(1)}](\tt_1, \tt_2)+\frac{\ep^2}{\tau_1^{2+\frac{3}{2}\dec}}\int_{\MM(\tau_1, \tau_2)}r^{1+\de}|\dk^{\leq s+1}\Xi|^2+\frac{\ep^2_0}{\tau_1^{2+3\dec}}.
\eeaa
as stated. This concludes the proof of Corollary \ref{cor:rpweightedestimatesforcombinednomrBEFsppmbphiminus2forallp}.
\end{proof}

%%%%%%%%%%%%%%%%%%%%%%%%%%%%%%%%%%%%%%%%%%%

\subsection{Proof of Theorem \ref{theoremM1:intro}}
\lab{sec:proofofThmM1}

%%%%%%%%%%%%%%%%%%%%%%%%%%%%%%%%%%%%%%%%%%%

In this section we  make use of Corollary \ref{cor:rpweightedestimatesforcombinednomrBEFsppmbphiplus2forallp} to prove Theorem \ref{theoremM1:intro} which we state in precise form below.
 
\begin{theorem}[Extension of Theorem M1 in \cite{KS:Kerr} to the full subextremal range]
\lab{theoremM1:Chap11}
Assume that the spacetime $\MM$  verifies the assumptions {stated in Section \ref{section:assumptionsneededforstatementofTheoremM1}}, as well as  the assumption \eqref{eq:controlofinitialdataforThM8-intro} on the  initial data. Then, if $\ep_0>0$ is sufficiently small, there exists $\dee>\dec$ such that we have the following   estimates for {$\pmb\phi_{+2}^{(2)}$ on  $\Si_*$, for all  $k\le \kst -13$,
 \bea\lab{eq:goal1}
\sup_{\Si_*}\Big(\big(\tau^{1+\dee}+r\tau^{\frac{1}{2}+\dee}\big)|\dk^k\pmb\phi_{+2}^{(2)}|+r\tau^{1+\dee}|\dk^{k-1}\nab_3\pmb\phi_{+2}^{(2)}|\Big) &\les& \ep_0
\eea
and 
\bea\lab{eq:goal2}
\int_{\Si_*(\geq\tau)}|\dk^{k-1}\nab_3\pmb\phi_{+2}^{(2)}|^2 &\les& \ep_0^2\tau^{-2-2\dec}.
\eea
Also, $A$ satisfies in $\MM$, for all  $k \le \kst -13$,
\bea\lab{eq:estimateThM1forAonMext}\lab{eq:goal3} 
\bsplit
\sup_{\MM}\Big(r^2\tau^{1+\dee}+r^3(2r+\tau)^{\frac{1}{2}+\dee}\Big)|\dk^kA|\\
+\sup_{\MM}\Big(r^3\tau^{1+\dee}+r^4\tau^{\frac{1}{2}+\dee}+r^{\frac{9}{2}+\dec}\Big)|\dk^{k-1}\nab_3A|\\
+\sup_{\MM}\Big(r^4\tau^{1+\dee}+r^{\frac{9}{2}+\dec}\Big)|\dk^{k-2}\nab_3^2A| &\les \ep_0.
\end{split}
\eea}
\end{theorem}

%%%%%%%%%%%%%%%%%%%%%%%%%%%%%%%%%

\subsubsection{Initial data assumptions}
\lab{sec:assumptioninitialdatanorm:allresults}

%%%%%%%%%%%%%%%%%%%%%%%%%%%%%%%%%

We now introduce the initial data assumptions used in the rest of the paper. To this end, we recall from Definition 1.5.1 the initial data norm denoted by $\Ik_k$ which measures the size of the perturbation from Kerr at $\tau=1$,   for $k$ derivatives of the  curvature tensor.
\begin{definition}\lab{def:initialdatanorm}
We define the following initial data  norms on $\Si_1:=\Si(\tau=1)$
\bea\lab{eq:defintitionoftheinitialdatanormJkforThM1and2andcurvestM8}
\bsplit
\Ik_{k} := &\sup_{S\subset\Si_1 }r^{\frac{5}{2}  +\de_B} \Big( { \big\| \dk^k\, A \big\|_{L^2(S)}}  + \big\| \dk^k\, B\big\|_{L^2(S)}\Big)\\
&+ \sup_{S\subset\Si_1 }\Big(  r^2 \big\| \dk^k\,  \Pc  \big\|_{L^2(S)}+  r  \big\| \dk^k\, \Bb\big\|_{L^2(S)} + \big\| \dk^k\, \Ab\big\|_{L^2(S)}\Big),
\end{split}
\eea
where $S=S(\tau=1, r)$ denote the level surfaces of $r$ in $\Si_1$, and where $(A, B, \Pc, \Bb, \Ab)$ in \eqref{eq:defintitionoftheinitialdatanormJkforThM1and2andcurvestM8} are defined w.r.t. the choice of global null frame $(e_3, e_4, e_1, e_2)$ fixed by Remark \ref{rmk:choiceofmaingloblaframeinthewholepaper}, i.e., 
\begin{itemize}
\item In Section \ref{sec:proofofThmM1}, the global null frame $(e_3, e_4, e_1, e_2)$ in \eqref{eq:defintitionoftheinitialdatanormJkforThM1and2andcurvestM8}  is the one in Section 3.6.4 in \cite{KS:Kerr}.

\item In Section \ref{sec:proofofThmM2}, the global null frame $(e_3, e_4, e_1, e_2)$ in \eqref{eq:defintitionoftheinitialdatanormJkforThM1and2andcurvestM8} is the one in Section 3.6.5 in \cite{KS:Kerr}. 

\item In Sections \ref{sec:higherordercurvatureestimatesforprovingThM8largea} \ref{sec:energyMorawetzforPc} \ref{sec:energyMorawetzforABBbAb}, the global null frame $(e_3, e_4, e_1, e_2)$ in \eqref{eq:defintitionoftheinitialdatanormJkforThM1and2andcurvestM8} is the one in Section 9.6.1 in \cite{KS:Kerr}.
\end{itemize}
\end{definition}

In the rest of the paper, we make the following assumption on the control of the initial data norm\footnote{The original assumption on initial data in \cite{KS:Kerr} is stated for $k_{large}+10$ derivatives, see {(3.52)} in that paper, in a given frame of an initial data layer $\LL(a_0, m_0)$. The control \eqref{eq:controlofinitialdataforThM8-intro} in the frames  used in this paper  are obtained in Theorem M0 of Section 3.7.1 in \cite{KS:Kerr}, and in Theorem {9.46} in  \cite{KS:Kerr} for $\kl+7$ derivatives with $\kl=k_{large}$.}
\bea\lab{eq:controlofinitialdataforThM8-intro}
\Ik_{\kl+7}\leq \ep_0.
\eea

%%%%%%%%%%%%%%%%%%%%%%%%%%%%%%%%%%%%%%%%%%%

\subsubsection{Spacetime assumptions for Theorem \ref{theoremM1:Chap11}}
\lab{section:assumptionsneededforstatementofTheoremM1}

%%%%%%%%%%%%%%%%%%%%%%%%%%%%%%%%%%%%%%%%%%%
 
The assumptions for Theorem \ref{theoremM1:Chap11} are as follows, see also Section 11.7.1 in \cite{GKS22}:
\begin{enumerate}
\item Assumptions on $\MM(\tau\leq\tau_*-2)$:
\begin{enumerate}
\item The setting on $(\MM, \g)$ is the one introduced in Sections \ref{sec:smallnesconstants}--\ref{sec:geometricsetuponSigma*}. Also, recall from Remarks \ref{rmk:choiceofmaingloblaframeinthewholepaper} and \ref{rmk:choiceofsecondgloblaframeinthewholepaper} that the global null frame $(e_3, e_4, e_1, e_2)$ introduced in Section \ref{section:SpacetimeMM-chap6} is the one in Section 3.6.4 in \cite{KS:Kerr}, while the second global null frame $(e_3', e_4', e_1', e_2')$ introduced in \ref{sec:secondglobalnullframeonMM} is the one of Section 3.6.5 in \cite{KS:Kerr}. 

\item The quantitative assumptions in Section \ref{sec:controlofGabandGagfromBA} hold on $\MM(\tau\leq\tau_*-2)$.

\item We choose the integer $\kst$ introduced in Section \ref{sec:interpolatedestimatesforGabGag} such that\footnote{Recall that $\kst$ satisfies $0<\kst - \frac{\kl}{2}\leq\frac{\kl\dec}{16}$ in view of \eqref{eq:choiceofinterpolationintegerk*forTheoremM1andM2}. Since $\kl\dec\gg 1$, see \eqref{eq:constraintsonthemainsmallconstantsepanddelta}, we may indeed enforce \eqref{eq:choiceofinterpolationintegerk*forTheoremM1andM2:moreprecise} provided $\kl\dec$ is chosen large enough. Also, the relation between $\kl$ and the integers appearing in the statement of Theorem M1 in Section 3.7.1 of \cite{KS:Kerr} is $\kl=k_{large}$ and $\frac{\kl}{2}=\frac{k_{large}}{2}=k_{small}$ so that \eqref{eq:choiceofinterpolationintegerk*forTheoremM1andM2:moreprecise} implies $\kst\geq k_{small}+120$.}
\bea\lab{eq:choiceofinterpolationintegerk*forTheoremM1andM2:moreprecise}
\kst\geq \frac{\kl}{2}+120.
\eea
Then, \eqref{eq:interpolationestimatesforGabGagfork*derivativeswithk*largerthanklover2:forproofTheoremM1andM2} implies that $(\Ga_g, \Ga_b)$ w.r.t. the global null frame $(e_3, e_4, e_1, e_2)$ of $\MM$ introduced in Section \ref{section:SpacetimeMM-chap6} satisfies on $\MM(\tau\leq\tau_*-2)$
\bea\lab{eq:assumptionsonMextforpartII-1}
r^2|\dk^{\leq k}\Ga_g|\leq \frac{\ep}{\tau^{\frac{1}{2}+\frac{3\dec}{4}}}, \qquad r|\dk^{\leq k}(\Ga_g, \Ga_b)| \leq \frac{\ep}{\tau^{1+\frac{3\dec}{4}}}, \qquad  k\leq \frac{\kl}{2}+120.
\eea
We also assume that the curvature components   $A, B$ verify,  for  $ k\leq k_L$, on $\MM(\tau\leq\tau_*-2)$
 \bea
 \lab{eq:assumptionsonMextforpartII-2}
  r^{7/2+{\dt}} |\dk^{\leq k} (A, B)| &\leq \ep,
 \eea
for some constant $\dt$ such that $\dt>2\dec$.

\item We make also the following additional assumptions on $\MM(\tau\leq\tau_*-2)$
 \bea
 \lab{eq:GlobalFrame-HcinGa_g1}
 \Hc \in \Ga_g
 \eea
 and 
 \bea
 \lab{eq:GlobalFrame-HcinGa_g2}
 r^3|\dk^{\le k} \Xi|\le \ep, \quad k\leq\kl, \qquad \nab_3\Xi\in r^{-1}\dk^{\leq 1}\Ga_g.
 \eea
\end{enumerate} 
\item Assumptions on $\MM(\tau\geq\tau_*-2)$:
\begin{enumerate}
\item  The region  $\MM(\tau\ge \tau_*-2)$ is covered by a second frame   $( e''_1, e''_2, e''_3, e''_4 )$ for which  all assumptions {\eqref{eq:assumptionsonMextforpartII-1}--\eqref{eq:GlobalFrame-HcinGa_g2}} are verified. 
  
 \item   The  change of frame coefficients $(f', \fb', \la')$ from  $(e_1, e_2, e_3, e_4)$  to  $(e_1'', e_2'', e_3'', e_4'')$  {verifies}, in $\MM(\tau_*-2, \tau_*)$. 
 \bea
 \lab{eq:frame-chpater11}
 |\dk^{\le k+1}(f', \fb', \la'-1)|\les \frac{\ep}{{r\tau^{1+\dec}}}, \qquad {  k\le\kst}.
 \eea
 
 \item The frames  $(e_1'', e_2'', e_3'', e_4'')$ and  $(e_1, e_2, e_3, e_4)$ coincide  in  a neighborhood  of  $\{\tau=\tau_*\}$. 
\end{enumerate}
\end{enumerate}

\begin{remark}
The additional assumptions \eqref{eq:GlobalFrame-HcinGa_g1} \eqref{eq:GlobalFrame-HcinGa_g2} hold true provided the global null frame $(e_3, e_4, e_1, e_2)$ is the one in Section 3.6.4 in \cite{KS:Kerr}, which is indeed our choice here as announced in Remarks \ref{rmk:choiceofmaingloblaframeinthewholepaper}. They are crucial\footnote{This fact also played an important role  in \cite{KS}, see {Sections} 3.4.6  and  3.5.3  and the application to the proof of Theorem M1 in chapter 5, as well as Theorem 5.2.9 in \cite{GKS22}.} in  deriving the  correct structure of the nonlinear correction terms $\widetilde{\N}_{W,+2}^{(p)}$ in \eqref{eq:schematicformofNpWsplus2}, see Theorem \ref{thm:derivationoftheTeukolskytensorialwavesystemfors=plusminus2:kerrpert:alternateformnullframeinsteadcoordvectorfield}.
\end{remark}

\begin{remark}\lab{remark:wheretofindtheframesofTheoremM1inKS2021}
The global null frame $( e_1'', e_2'', e_3'', e_4'')$ on $\MM(\tau\ge \tau_*-2)$ is the one constructed\footnote{{In fact, the frame $( e_1'', e_2'', e_3'', e_4'')$ in Proposition   4.4.3  of \cite{KS:Kerr} satisfies $\Xi'\in r^{-1}\Ga_g$ which is stronger than \eqref{eq:GlobalFrame-HcinGa_g2}, but we will not need this stronger control for $\Xi'$ here.}} in  Proposition   {4.8}  of \cite{KS:Kerr}. Also, \eqref{eq:frame-chpater11} corresponds to Property (d) of Lemma {4.13} in \cite{KS:Kerr}.  
\end{remark}

\begin{remark}
In fact, for the global null frame $(e_3, e_4, e_1, e_2)$ in Section 3.6.4 in \cite{KS:Kerr}, the quantitative assumptions for $(\Ga_g, \Ga_b)$ in \eqref{eq:assumptionsonMMforpartII} hold on $\MM$, while the assumption for $\xi$ in \eqref{eq:assumptionsonMMforpartII} only holds on $\MM(1, \tau_*-2)$. In particular, \eqref{eq:assumptionsonMextforpartII-1} holds on $\MM$ as well, and actually, so does \eqref{eq:assumptionsonMextforpartII-2}. Thus, the only quantitive assumptions that hold only on $\MM(1, \tau_*-2)$ are the ones in \eqref{eq:GlobalFrame-HcinGa_g1} \eqref{eq:GlobalFrame-HcinGa_g2}.
\end{remark}

%%%%%%%%%%%%%%%%%%%%%%%%%%%%%%%%%%%%%%%%%%

\subsubsection{Proof of Theorem \ref{theoremM1:Chap11}}
\lab{sec:finallyproofofThmM1}

%%%%%%%%%%%%%%%%%%%%%%%%%%%%%%%%%%%%%%%%%%
 
We now prove Theorem \ref{theoremM1:Chap11}. As in Section 11.7 of \cite{GKS22}, the proof is done first on $\MM(1, \tau_*-2)$ and then on $\MM(\tau_*-2, \tau_*)$.

%%%%%%%%%%%%%%%%%%%%%%%%%%%%%%%%%%%%%%%%%%%%

\paragraph{\textit{Proof of Theorem \ref{theoremM1:Chap11} on $\MM(1, \tau_*-2)$}.}

%%%%%%%%%%%%%%%%%%%%%%%%%%%%%%%%%%%%%%%%%%%% 

The strategy is a slight adaptation of the one in Section 11.7.3 of \cite{GKS22}. \\

\noindent{\bf Step 1.} In view of Corollary \ref{cor:rpweightedestimatesforcombinednomrBEFsppmbphiplus2forallp}, we have, for $11\leq s\leq \kst-3$ and $\de\leq p\leq 2-\de$, 
\bea\lab{eq:BEFsppmbphiplus2combinednormwithstilltheerrortermstobeestimatedalsmostfinal}
\nn\BEF^s_{p}[\pmb\phi_{+2}](\tau_1, \tau_2) &\les& \E^s_{p}[\pmb\phi_{+2}](\tau_1)+\NN^s_{p}[\pmb\phi_{+2}^{(2)}, \widetilde{\N}_{W,+2}^{(2)}]{(\tau_1, \tau_2)} +\widetilde{\NN}^s_\de[\pmb\phi_{+2}^{(0)}, \widetilde{\N}_{W,+2}^{(0)}](\tt_1, \tt_2)\\
&&+\widetilde{\NN}^s_\de[\pmb\phi_{+2}^{(1)}, \widetilde{\N}_{W,+2}^{(1)}](\tt_1, \tt_2),
\eea
which is the analog of Theorem 11.2.3 in \cite{GKS22}. We then need to estimate the terms involving $\widetilde{\N}_{W,+2}^{(p)}$ on the RHS of \eqref{eq:BEFsppmbphiplus2combinednormwithstilltheerrortermstobeestimatedalsmostfinal}. 

First, noticing that the term $\NN^s_{p}[\pmb\phi_{+2}^{(2)}, \widetilde{\N}_{W,+2}^{(2)}]{(\tau_1, \tau_2)}$ on the RHS of \eqref{eq:BEFsppmbphiplus2combinednormwithstilltheerrortermstobeestimatedalsmostfinal} coincides with the term\footnote{Indeed, we have $\qf=\pmb\phi_{+2}^{(2)}$ in view of Remark \ref{rmk:compasisionphiplus2p=0withMaSz26andphiplus2p=2withqfGKS22}.} $\NN^s_{p}[\qf, N_{\err}]{(\tau_1, \tau_2)}$ on the RHS of (11.2.1) in \cite{GKS22}, we may directly apply the estimate for $\NN^s_{p}[\qf, N_{\err}]{(\tau_1, \tau_2)}$ in the proof of Theorem 11.6.1, which does not rely on the smallness of $a$, to obtain, for $11\leq s\leq \kst-3$ and $\de\leq p\leq 2-\de$,  
\beaa
\NN^s_{p}[\pmb\phi_{+2}^{(2)}, \widetilde{\N}_{W,+2}^{(2)}]{(\tau_1, \tau_2)} \les \ep_0^2\tau_1^{p-2-3\dec}.
\eeaa
Then, notice that $\widetilde{\NN}^s_\de[\pmb\phi_{+2}^{(q)}, \widetilde{\N}_{W,+2}^{(q)}](\tt_1, \tt_2)$, $q=0,1$, clearly satisfy the corresponding estimate for $p=\de$, i.e. 
\beaa
\sum_{q=0}^1\NN^s_{\de}[\pmb\phi_{+2}^{(q)}, \widetilde{\N}_{W,+2}^{(q)}]{(\tau_1, \tau_2)} \les \ep_0^2\tau_1^{\de-2-3\dec},
\eeaa
since $\pmb\phi_{+2}^{(q)}$, $q=0,1$, satisfy better decay estimates than $\pmb\phi_{+2}^{(2)}$ and $\widetilde{\N}^{(q)}_{W,+2}$, $q=0,1$, display a more favorable structure than $\widetilde{\N}^{(2)}_{W,+2}$ in view of \eqref{eq:schematicformofNpWsplus2}. Plugging in \eqref{eq:BEFsppmbphiplus2combinednormwithstilltheerrortermstobeestimatedalsmostfinal}, we infer, for $11\leq s\leq \kst-3$ and $\de\leq p\leq 2-\de$, 
\bea\lab{eq:BEFsppmbphiplus2combinednormwithstilltheerrortermstobeestimatedalsmostfinal:evenmorefinal}
\BEF^s_{p}[\pmb\phi_{+2}](\tau_1, \tau_2) &\les& \E^s_{p}[\pmb\phi_{+2}](\tau_1)+\ep_0^2\tau_1^{p-2-3\dec},
\eea
which is the analog of (11.6.1) in \cite{GKS22}.

\noindent{\bf Step 2.}  Relying on \eqref{eq:BEFsppmbphiplus2combinednormwithstilltheerrortermstobeestimatedalsmostfinal:evenmorefinal}, 
we  run the basic mean value argument for ${\tau}\le \tau_*-2$ and obtain, for $s\le \kst-5$,
\bea
\lab{eq:theoremM1-Chap11-3}
\BEF^s_p[\pmb\phi_{+2}](\tau, \tau_*-2) &\les& \ep_0^2\tau^{-(2-p-\de)}.
\eea

\noindent{\bf Step 3.} Next, we note that the improved $r^p$-weighted estimates for the tensorial wave equation for $\pmb\phi_{+2}^{(2)}$ in Theorem 11.6.2 in \cite{GKS22} only involve estimates in $r\ge R$ for $R$ large and apply immediately without any restriction on $a$, recalling again that $\qf=\pmb\phi_{+2}^{(2)}$ in view of Remark \ref{rmk:compasisionphiplus2p=0withMaSz26andphiplus2p=2withqfGKS22}. We can then immediately follow the procedure in Steps 2 and 3 of Section 11.7.3 in \cite{GKS22}, which relies on \eqref{eq:theoremM1-Chap11-3}, the improved $r^p$-weighted estimates for $\pmb\phi_{+2}^{(2)}$ in Theorem 11.6.2 in \cite{GKS22}, and a mean value argument. This yields, for ${\tau}\le \tau_*-2$ and for $s\le \kst-6$,
\beaa
\E^s_{2-\de}[\pmb\phi_{+2}^{(2)}](\tau) &\les& \ep_0^2\tau^{-3\dec-\de}.
\eeaa

\noindent{\bf Step 4.} We now rely on the following slight variation of \eqref{eq:rpweigthedestimateRlargeforpmbphiplus2p=2:intermediarydirectapplication:bis}, for ${\tau}\le \tau_*-2$, $\de\leq p\leq 2-\de$ and $s\leq\kl-1$, 
\beaa
\nn\BEF^s_{p, \geq R}[\pmb\phi_{+2}^{(2)}](\tau, \tau_*-2) &\les& \E^s_{p, \geq \frac{R}{2}}[\pmb\phi_{+2}^{(2)}](\tau)+\NN^s_{p, \ge R/2}[\pmb\phi_{+2}^{(2)}, \widetilde{\N}_{W,+2}^{(2)}]{(\tau, \tau_*-2)} +\M^s_{\frac{R}{2}, R}[\pmb\phi_{+2}^{(2)}](\tau, \tau_*-2)\\
&&+R^{-2\de}\BEF^s_{\de, \geq R/2}[\pmb\phi_{+2}](\tau, \tau_*-2)\\
&\les& E^s_{p}[\pmb\phi_{+2}^{(2)}](\tau)+\NN^s_{p}[\pmb\phi_{+2}^{(2)}, \widetilde{\N}_{W,+2}^{(2)}]{(\tau, \tau_*-2)}+\BEF^s_{\de}[\pmb\phi_{+2}](\tau, \tau_*-2).
\eeaa
Together with the estimate for $\BEF^s_\de[\pmb\phi_{+2}](\tau, \tau_*-2)$ given by \eqref{eq:theoremM1-Chap11-3} and the estimate in Step 1 for $\NN^s_{2-\de}[\pmb\phi_{+2}^{(2)}, \widetilde{\N}_{W,+2}^{(2)}]{(\tau, \tau_*-2)}$, we obtain, for ${\tau}\le \tau_*-2$, $\de\leq p\leq 2-\de$ and $s\le \kst-5$,
\beaa
\nn\BEF^s_{p}[\pmb\phi_{+2}^{(2)}](\tau, \tau_*-2) &\les& \E^s_{p}[\pmb\phi_{+2}^{(2)}](\tau)+\ep_0^2\tau^{p-2-3\dec}+\ep_0^2\tau^{-(2-2\de)}.
\eeaa
In particular, we infer for $s\le \kst-5$, {$\tau_1<\tau_2\le \tau_*-2$ and $3\dec +2\de \leq p\leq 2-\de$},
 \beaa
  \BEF^s_{p}[\pmb\phi_{+2}^{(2)}](\tau_1, \tau_2) \les
       \E^s_{p}[\pmb\phi_{+2}^{(2)}](\tau_1)  +\ep_0^2\tau_1^{-(2+3\dec -p)}.
  \eeaa
Together with Step 3 and the usual mean value argument, we deduce  for $s\le \kst-7$,
 \bea\lab{eq:Flux-psiStep4:beforetolast}
  \BEF^s_{p}[\pmb\phi_{+2}^{(2)}](\tau_1, \tau_2) \les
        \ep_0^2\tau_1^{-(2+3\dec -p)}, \qquad 1-\de \leq p\leq 2-\de.
  \eea
In particular,  for $s\le \kst -7$, we have 
\bea
\lab{eq:Flux-psiStep4}
\BEF^s_{1-\de}[\pmb\phi_{+2}^{(2)}](\tau_1, \tau_2)   \les
        \ep_0^2\tau_1^{-(3\dec+1+\de)}.
  \eea

\noindent{\bf Step 5.} As in Proposition 5.12 in \cite{KS}, {we use the trace theorem and} we interpolate the control on $\E_p^s[\pmb\phi_{+2}^{(2)}]$ provided \eqref{eq:Flux-psiStep4:beforetolast} between $p=1+\de$ and $p=1-\de$. This yields 
\beaa
\tau^{1+3\dec}\int_{S_r}|\dk^{\le s}\pmb\phi_{+2}^{(2)}|^2 &\les& \ep_0^2.  
\eeaa
Using Sobolev, we infer the following pointwise decay estimate for $\pmb\phi_{+2}^{(2)}$, for all $s\le \kst-9$,
\bea\lab{eq:pointwisedecayforpsiatrminus1plustaudecauStep5}
|\dk^{\le s }\pmb\phi_{+2}^{(2)}| &\les& \ep_0r^{-1}\tau^{-\frac{1+3\dec}{2}}.
\eea

\noindent{\bf Step 6.}   We now integrate the transport equations \eqref{eq:transportequationsins=+2caseforp=0and1} {from the initial data layer} and using   the pointwise decay for $\pmb\phi_{+2}^{(2)}$ in \eqref{eq:pointwisedecayforpsiatrminus1plustaudecauStep5}, 
 we  infer the following pointwise decay estimates, for all  $s\le \kst-9$,
\bea
|\dk^{\le s}\pmb\phi_{+2}^{(1)}| \les \ep_0r^{-2}\tau^{-\frac{1+3\dec}{2}},\qquad |\dk^{\le s}\pmb\phi_{+2}^{(0)}| \les \ep_0r^{-3}\tau^{-\frac{1+3\dec}{2}}.
\eea
Integrating these estimates on $\Si(\tau)$, we  deduce, for all $s\le \kst-10$ and $\tau\leq\tau_*-2$,
\bea
\lab{eq:Flux-psiStep5}
\int_{\Si(\tau)}r^{1-\de}|\dk^{\le s+1}\pmb\phi_{+2}^{(1)}|^2+\int_{\Si(\tau)}r^{3-\de}|\dk^{\le s+1}\pmb\phi_{+2}^{(0)}|^2\les \ep_0^2\tau^{-(1+3\dec)}.
\eea

\noindent{\bf Step 7.} In view of Definition \ref{def:combinednormEpandBEFpTeukolskywavetransportsystem:plus2case}, the estimates  \eqref{eq:Flux-psiStep4}  and \eqref{eq:Flux-psiStep5}  imply, for all $s\le \kst -10$,
\beaa
\E^s_{1-\de}[\pmb\phi_{+2}](\tau) &\les& \ep_0^2\tau^{-(1+3\dec)}.
\eeaa
Together with \eqref{eq:BEFsppmbphiplus2combinednormwithstilltheerrortermstobeestimatedalsmostfinal:evenmorefinal}, we un again  the standard mean value argument and deduce, for all $s\le \kst-11$ and $\tau\leq\tau_*-2$,
\beaa
\BEF^s_{\de}[\pmb\phi_{+2}](\tau, \tau_*-2) &\les& \ep_0^2\tau^{-(2+3\dec-2\de)}.
\eeaa
In particular,  for $s\le \kst -11$,
\bea
\lab{eq:Flux-psiStep7}
\BEF^s_{\de}[\pmb\phi_{+2}^{(2)}](\tau, \tau_*-2)  &\les& \ep_0^2\tau^{-(2+3\dec- 2\de)}
\eea
and, in view of the definition  of the flux norms $\F_\de^s[\psi]$  in {Section} \ref{subsection:basicnormsforpsi}, we deduce{, on $\Si_*$,} 
\beaa
 \int_{\Si_*{(\tau, \tau_*-2)}}|\nab_3 \dk^{s-1}\pmb\phi_{+2}^{(2)}|^2 &\les& \ep_0^2\tau^{-(2+3\dec-3 \de)}.
\eeaa
Hence, choosing $\de_{extra}= \frac{3\dec-2\de}{2}>\dec$, for $s\le \kst-10$,
\beaa
 \int_{\Si_*{(\tau, \tau_*-2)}}|\nab_3 \dk^{s-1} \psi|^2 &\les& \ep_0^2 \tau^{-2-2\de_{extra} }
 \eeaa
which establishes  the   desired  estimate \eqref{eq:goal2} {for $\tau\leq \tau_*-2$}.

\noindent{\bf Step 8.}  Making  use of  the estimate \eqref{eq:Flux-psiStep7} 
 and proceeding as in the  derivation of  the estimate  (5.2.7) in   Proposition 5.12 of \cite{KS},  we derive the estimate, for any $S_r\subset\Si(\tau)$, {$\tau\leq \tau_*-2$}, $s\le \kst-11$, 
\bea
\lab{eq:Flux-psiStep8-1}
 r^{-1} \int_{S_r} | \dk^{\le s}\pmb\phi_{+2}^{(2)}|^2 &\les &\ep_0^2  \tau^{-(2+3\dec-2\de)}.
\eea
Similarly, proceeding as in the estimate (5.2.9)  in Proposition 5.13 in \cite{KS}  we  derive
 for any $S_r\subset\Si(\tau)$, {$\tau\leq \tau_*-2$}, $s\le \kst-11$,
\bea
\lab{eq:Flux-psiStep8-2}
\int_{S_r} | \dk^{\le s-1} \nab_3\pmb\phi_{+2}^{(2)}|^2 &\les &\ep_0^2  \tau^{-(2+3\dec-2\de)}.
\eea

\noindent{\bf Step 9.} In view of \eqref{eq:pointwisedecayforpsiatrminus1plustaudecauStep5}, 
 \eqref{eq:Flux-psiStep8-1} and \eqref{eq:Flux-psiStep8-2}, we  deduce the following estimate  for $\qf$, for all $s\le \kst-13$,
\bea
\lab{eq:theoremM1-Chap11-psi}
\sup_{\MM{(\tau\leq\tau_*-2)}}\Big(r\tau^{\frac{1}{2}+\dee}+\tau^{1+\dee}\Big)|\dk^{\leq s}\pmb\phi_{+2}^{(2)}|+\sup_{\MM{(\tau\leq\tau_*-2)}}r\tau^{1+\dee}|\dk^{\leq s-1}\nab_3\pmb\phi_{+2}^{(2)}| \les \ep_0
\eea
for $\de_{extra}= \frac{3\dec-2\de}{2}>\dec$.

\noindent{\bf Step 10.}    Using  \eqref{eq:theoremM1-Chap11-psi}, integrating the transport equations \eqref{eq:transportequationsins=+2caseforp=0and1} {from the initial data layer}, and using the definition of $\pmb\phi_{+2}^{(q)}$, $q=0,1$, in terms of $A$ and $\nab_3A$ in \eqref{eq:definitionofthephiplus2phierarchy:perturbationofKerr}, we derive the following pointwise estimate for $A$, for all $s\le \kst -13$,
 \beaa
\sup_{\MM{(\tau\leq\tau_*-2)}}\Big(r^2\tau^{1+\dee}+r^3(2r+\tau)^{\frac{1}{2}+\dee}\Big)\Big(|\dk^{\le s} A|+r|\dk^{\le s-1}\nab_3 A|\Big)&&\\
{+\sup_{\MM(\tau\leq\tau_*-2)}\Big(r^4\tau^{1+\dee}+r^5\tau^{\frac{1}{2}+\dee}\Big)|\dk^{\le s-2}\nab_3^2A|} &\les& \ep_0
\eeaa
as stated in \eqref{eq:goal3}. Together with the proof of \eqref{eq:goal2} in Step 7, this concludes the proof of Theorem \ref{theoremM1:Chap11} {for $\tau\leq \tau_*-2$}.

%%%%%%%%%%%%%%%%%%%%%%%%%%%%%%%%%%%%%%%%%%%

\paragraph{\textit{End of the proof of Theorem \ref{theoremM1:Chap11}}.}

%%%%%%%%%%%%%%%%%%%%%%%%%%%%%%%%%%%%%%%%%%% 

Given that $\MM(\tau_*-2,\tau_*)$ is a local existence type region, the arguments in Sections 11.7.4 to 11.7.7 in \cite{GKS22} never rely on the size of $a$ and are thus immediately valid in the full subextremal range. 
In particular, this immediately yields the proof of Theorem \ref{theoremM1:Chap11} on $\MM(\tau_*-2,\tau_*)$ in the full subextremal range. Together with the proof of Theorem \ref{theoremM1:Chap11} on $\tau\leq \tau_*-2$ in the previous paragraph, this concludes the proof of Theorem \ref{theoremM1:Chap11}.

%%%%%%%%%%%%%%%%%%%%%%%%%%%%%%%%%%%%%%%%%%%

\subsection{Proof of Theorem \ref{theoremM2:intro}}
\lab{sec:proofofThmM2}

%%%%%%%%%%%%%%%%%%%%%%%%%%%%%%%%%%%%%%%%%%%

In this section we  make use of Corollary \ref{cor:rpweightedestimatesforcombinednomrBEFsppmbphiminus2forallp} to prove Theorem \ref{theoremM2:intro} which we state in precise form below.
 
\begin{theorem}[Extension of Theorem M2 in \cite{KS:Kerr} to the full subextremal range]
\lab{thm:restatementofTheoremM2}
Assume that the global frame of $\MM$ satisfies the assumptions of {Section} \ref{sec:assumptionsontheframe:Chapter12}, as well as  the assumption \eqref{eq:controlofinitialdataforThM8-intro} on the  initial data. Also, assume that the control of the flux of $\pmb\phi_{+2}^{(2)}$ provided by \eqref{eq:goal2} holds, i.e., for $k\leq\kst-14$, 
\bea\lab{eq:assumptionsonqfforTheoremM2compatiblewithTheoremM1:chap12}
\int_{\Si_*(\geq\tau)}|\nab_3\dk^k\pmb\phi_{+2}^{(2)}|^2 &\les& \ep_0^2\tau^{-2-2\dec}.
\eea
Then, $\aa$ satisfies the following estimate on $\Si_*$, for all $1\leq \tau\leq \tau_*$ and $k\leq \kst -14$, 
\beaa
\int_{\Si_*(\geq \tau)}|\dk^k\aa|^2 &\les& \ep_0^2\tau^{-2-2\dec}.
\eeaa
\end{theorem}

%%%%%%%%%%%%%%%%%%%%%%%%%%%%%%%%%%%%%%%%%%%%%

\subsubsection{Spacetime assumptions for Theorem \ref{thm:restatementofTheoremM2}}
\lab{sec:assumptionsontheframe:Chapter12}

%%%%%%%%%%%%%%%%%%%%%%%%%%%%%%%%%%%%%%%%%%%%%

The assumptions for Theorem \ref{thm:restatementofTheoremM2} are as follows, see also Section 12.4.1 in \cite{GKS22}:
\begin{enumerate}
\item Assumptions on $\MM$:
\begin{enumerate}
\item The setting on $(\MM, \g)$ is the one introduced in Sections \ref{sec:smallnesconstants}--\ref{sec:geometricsetuponSigma*}. Also, recall from Remarks \ref{rmk:choiceofmaingloblaframeinthewholepaper} and \ref{rmk:choiceofsecondgloblaframeinthewholepaper} that the global null frame $(e_3, e_4, e_1, e_2)$ introduced in Section \ref{section:SpacetimeMM-chap6} and second global null frame $(e_3', e_4', e_1', e_2')$ introduced in \ref{sec:secondglobalnullframeonMM} both coincidence with the one of Section 3.6.5 in \cite{KS:Kerr}. 

\item We choose the integer $\kst$ introduced in Section \ref{sec:interpolatedestimatesforGabGag} such that \eqref{eq:choiceofinterpolationintegerk*forTheoremM1andM2:moreprecise} holds. Then, as in \eqref{eq:assumptionsonMextforpartII-1}, $(\Ga_g, \Ga_b)$ w.r.t. the global null frame $(e_3, e_4, e_1, e_2)$ of $\MM$ satisfies on $\MM$
\bea\lab{eq:assumptionsonMextforpartIIchap12}
r^2|\dk^{\leq k}\Ga_g|\leq \frac{\ep}{\tau^{\frac{1}{2}+\frac{3\dec}{4}}}, \qquad r|\dk^{\leq k}(\Ga_g, \Ga_b)| \leq \frac{\ep}{\tau^{1+\frac{3\dec}{4}}}, \qquad  k\leq \frac{\kl}{2}+120.
\eea
We also assume that the curvature components   $A, B$ verify,  for  $ k\leq k_L$, on $\MM$
 \bea
 \lab{eq:assumptionsonMextforpartIIchap12:moredecayinrAB}
  r^{7/2+{\dec}} |\dk^{\leq k} (A, B)| &\leq \ep.
 \eea
\end{enumerate} 
\item We make also the following additional assumptions on $\MM(\tau\leq\tau_*-2)$
 \bea
\lab{eq:Xi-Hb-chapter12}
\Xi\in r^{-2}\Ga_g, \qquad \Hbc\in r^{-1}\Ga_g.
\eea
\end{enumerate}

\begin{remark}
The additional assumptions \eqref{eq:Xi-Hb-chapter12} hold true provided the global null frame $(e_3, e_4, e_1, e_2)$ is the one in Section 3.6.5 in \cite{KS:Kerr}, which is indeed our choice here as announced in Remarks \ref{rmk:choiceofmaingloblaframeinthewholepaper}. They are crucial\footnote{This fact already played an important role  in Theorem 5.3.7 in \cite{GKS22}.} in  deriving the  correct structure of the nonlinear correction terms $\widetilde{\N}_{W,-2}^{(p)}$ in \eqref{eq:schematicformofNpWsminus2}, see Theorem \ref{thm:derivationoftheTeukolskytensorialwavesystemfors=plusminus2:kerrpert:alternateformnullframeinsteadcoordvectorfield}.
\end{remark}

%%%%%%%%%%%%%%%%%%%%%%%%%%%%%%%%%%%%%%%%%%%

\subsubsection{Proof of Theorem \ref{thm:restatementofTheoremM2}}

%%%%%%%%%%%%%%%%%%%%%%%%%%%%%%%%%%%%%%%%%%%

We now prove Theorem \ref{thm:restatementofTheoremM2}. As in Section 12.4 of \cite{GKS22}, the proof is done first on $\MM(1, \tau_*-2)$ and then on $\MM(\tau_*-2, \tau_*)$.

%%%%%%%%%%%%%%%%%%%%%%%%%%%%%%%%%%%%%%%%%%%%%%%%%%%%%%%%%%

\paragraph{\textit{Estimates on $\MM(1, \tau_*-2)$ for the proof of Theorem \ref{thm:restatementofTheoremM2}}.}

%%%%%%%%%%%%%%%%%%%%%%%%%%%%%%%%%%%%%%%%%%%%%%%%%%%%%%%%%%

We start with estimates on $\MM(1, \tau_*-2)$. 

\noindent{\bf Step 1.} We start with proving an analog of Theorem 12.2.4 in \cite{GKS22}. In view of Corollary \ref{cor:rpweightedestimatesforcombinednomrBEFsppmbphiminus2forallp}, we have, for $14\leq s\leq \kst-3$ and $\de\leq p\leq 2-\de$, 
\beaa
\nn\BEF^s_{p}[\pmb\phi_{-2}](\tau_1, \tau_2) &\les& \E^s_{p}[\pmb\phi_{-2}](\tau_1)+\NN^s_{p}[\pmb\phi_{-2}^{(2)}, \widetilde{\N}_{W,-2}^{(2)}]{(\tau_1, \tau_2)} +\widetilde{\NN}^s_\de[\pmb\phi_{-2}^{(0)}, \widetilde{\N}_{W,-2}^{(0)}](\tt_1, \tt_2)\\
\nn&&+\widetilde{\NN}^{s}_\de[\pmb\phi_{-2}^{(1)}, \widetilde{\N}_{W,-2}^{(1)}](\tt_1, \tt_2)+\frac{\ep^2}{\tau_1^{2+\frac{3}{2}\dec}}\int_{\MM(\tau_1, \tau_2)}r^{1+\de}|\dk^{\leq s+1}\Xi|^2+\frac{\ep^2_0}{\tau_1^{2+3\dec}}.
\eeaa
Together with \eqref{eq:Xi-Hb-chapter12} and \eqref{eq:assumptionsonMextforpartIIchap12}, this yields, for $14\leq s\leq \kst-3$, $\tau_1<\tau_2\leq\tau_*-2$ and $\de\leq p\leq 2-\de$,
\bea\lab{eq:BEFsppmbphiminus2combinednormwithstilltheerrortermstobeestimatedalsmostfinal}
\nn\BEF^s_{p}[\pmb\phi_{-2}](\tau_1, \tau_2) &\les& \E^s_{p}[\pmb\phi_{-2}](\tau_1)+\NN^s_{p}[\pmb\phi_{-2}^{(2)}, \widetilde{\N}_{W,-2}^{(2)}]{(\tau_1, \tau_2)} +\widetilde{\NN}^s_\de[\pmb\phi_{-2}^{(0)}, \widetilde{\N}_{W,-2}^{(0)}](\tt_1, \tt_2)\\
&&+\widetilde{\NN}^{s}_\de[\pmb\phi_{-2}^{(1)}, \widetilde{\N}_{W,-2}^{(1)}](\tt_1, \tt_2)+\frac{\ep^2_0}{\tau_1^{2+3\dec}}.
\eea
Then, in view of \eqref{eq:schematicformofNpWsminus2} and \eqref{eq:Xi-Hb-chapter12}, we have, on $\tau_1<\tau_2\leq\tau_*-2$,
\bea\lab{eq:structureofNpWminus2forp=012matchingtheoneinGKS22chap12}
\bsplit
\widetilde{\N}^{(0)}_{W,-2} &=  \dk^{\leq 1}\big(\Ga_b\c \Ga_g\big),\qquad \widetilde{\N}^{(1)}_{W,-2} = \dk^{\leq 2}\big(\Ga_b\c \Ga_g\big),\\ 
\widetilde{\N}^{(2)}_{W,-2} &=  \widetilde{N}_{\err}+ \dk^{\leq 3} (\Ga_g \c \Ga_b),\qquad \widetilde{N}_{\err}=r^2 \dk^{\leq 2}(\Ga_b \c (A, B)),
\end{split}
 \eea
 and the second estimate in (12.2.9) of \cite{GKS22}, whose proof does not depend on the size of $a$, yields, for $s\leq \kst-3$, $\tau_1<\tau_2\leq\tau_*-2$ and $\de\leq p\leq 2-\de$,
 \beaa
 \widetilde{\NN}^{s}_\de[\pmb\phi_{-2}^{(0)}, \widetilde{\N}_{W,-2}^{(0)}]{(\tau_1, \tau_2)} +\widetilde{\NN}^s_\de[\pmb\phi_{-2}^{(1)}, \widetilde{\N}_{W,-2}^{(1)}](\tt_1, \tt_2)\\
+\NN^s_{p}[\pmb\phi_{-2}^{(2)}, \widetilde{\N}_{W,-2}^{(2)}-\widetilde{N}_{\err}](\tt_1, \tt_2) &\les& \ep_0\tau_1^{-1-\dec}\Big(\BEF^s_p[\pmb\phi_{-2}](\tau_1, \tau_2)\Big)^{\frac{1}{2}} {+\ep_0^2\tau_1^{-2-2\dec}}.
 \eeaa
Plugging in \eqref{eq:BEFsppmbphiminus2combinednormwithstilltheerrortermstobeestimatedalsmostfinal}, we deduce, for $s\leq \kst-3$, $\tau_1<\tau_2\leq\tau_*-2$ and $\de\leq p\leq 2-\de$,
\bea\lab{eq:BEFsppmbphiminus2combinednormwithstilltheerrortermstobeestimatedalsmostfinal:almosttherenow}
\BEF^s_{p}[\pmb\phi_{-2}](\tau_1, \tau_2) &\les& \E^s_{p}[\pmb\phi_{-2}](\tau_1)+\NN^s_{p}[\pmb\phi_{-2}^{(2)}, \widetilde{N}_{\err}]{(\tau_1, \tau_2)} +\ep_0^2\tau_1^{-2-2\dec},
\eea
which is the analog of Theorem 12.2.4 in \cite{GKS22}. 

\noindent{\bf Step 2.} Next, we derive an analog of Corollary 12.2.5 in \cite{GKS22}. To this end, we then need to estimate the terms involving $\widetilde{N}_{\err}$ on the RHS of \eqref{eq:BEFsppmbphiminus2combinednormwithstilltheerrortermstobeestimatedalsmostfinal:almosttherenow}. 
This is done in Lemma 12.2.6 in \cite{GKS22} whose proof does not depend on the size of $a$ and yields, for $s\leq \kst-3$ and $\de\leq p\leq 1-\de$,
\beaa
   \NN^s_{p}[\pmb\phi_{-2}^{(2)}, \widetilde{N}_{\err}]{(\tau_1, \tau_2)}
      &\les& \ep_0\tau_1^{-1-\frac{3}{2}\dec+\frac{p+\de}{2}}\Big(\BEF^s_p[\pmb\phi_{-2}^{(2)}](\tau_1, \tau_2)\Big)^{\frac{1}{2}}.
\eeaa    
Plugging in \eqref{eq:BEFsppmbphiminus2combinednormwithstilltheerrortermstobeestimatedalsmostfinal:almosttherenow}, we infer, , for $s\leq \kst-3$ and $\tau_1<\tau_2\leq\tau_*-2$,
\bea\lab{eq:step1fluxSigamstarLiebTAbchap12:1'} 
\BEF^s_{p}[\pmb\phi_{-2}](\tau_1, \tau_2) &\les&  \E^s_{p}[\pmb\phi_{-2}](\tau_1) +\ep_0^2\tau_1^{-2-3\dec+p+\de} +\ep_0^2\tau_1^{-2-2\dec}, \quad  \de\leq p\leq 1-\de,
\eea
which is the analog of Corollary 12.2.5 in \cite{GKS22}. 

\noindent{\bf Step 3.} Next, we derive the following analog\footnote{\lab{footnote:onelossofderivativeforeq:analogofequationparenthesis12dot2dot7parenthesisinGKS22comparedtoGKS22}In fact, \eqref{eq:analogofequationparenthesis12dot2dot7parenthesisinGKS22} loses a derivative compared to (12.2.7) in \cite{GKS22}. In turn, this induces a loss of few derivatives in the statement of Theorem \ref{thm:restatementofTheoremM2} compared to the one of Theorem M2 in \cite{GKS22}.} of Proposition 12.2.9 in \cite{GKS22}. 

\begin{proposition}\lab{prop:analogofProposition12dot2dot9inGKS22}
We have, for $1\leq\tau_1<\tau_2\leq\tau_*-2$, $s\leq \kst-3$ and $\de\leq p\leq 2-\de$,
\begin{align}\lab{eq:analogofequationparenthesis12dot2dot7parenthesisinGKS22}
&\int_{\MM(\tau_1, \tau_2)}r^{p-3}|\dk^{\leq s}(\Ab, \underline{\Psi})|^2+\sup_{\tau\in[\tau_1, \tau_2]}\int_{\Si(\tau)}r^{p-4}|\dk^{\leq s}(\Ab, \underline{\Psi})|^2+\int_{\Si_*(\tau_1, \tau_2)}r^{p-2}|\dk^{\leq s}(\Ab, \underline{\Psi})|^2\nn\\
\les & \int_{\Si(\tau_1)}r^{p-4}|\dk^{\leq s}(\Ab, \underline{\Psi})|^2 + \B^{s}_\de[\pmb\phi_{-2}^{(2)}](\tau_1, \tau_2)+\ep_0^2\tau_1^{-2-2\dec},
\end{align}
where $\underline{\Psi}\in\sk_2(\mathbb{C})$ is given by Definition \ref{def:PsibforintegrationofAbfromqfb:chap12}.
\end{proposition}

\begin{proof}
Suppose $\Phi_1, \Phi_2 \in \sk_2(\CCC)$ with signature $s\leq -1$ satisfy the differential relation 
\bea\label{eq:relation-nab4Phi1Phi2}
\nabc_4\Phi_1=\Phi_2.
\eea
Then, for every {$\de_0\geq\de$}, we have
 \bea\label{eq:general-integrated-estimate-e4}
 \bsplit
&\int_{\MM(\tau_1, \tau_2)}  r^{-\de_0-3} |\Phi_1|^2
+ \sup_{\tau\in[\tau_1,\tau_2]}\int_{\Si(\tau)} r^{-\de_0 -4}|\Phi_1|^2+  \int_{\AA\cup \Si_*(\tau_1,\tau_2)} r^{-\de_0 -2}|\Phi_1|^2  \\
  &\les \int_{\MM(\tau_1, \tau_2)} r^{-\de_0 -1}|\Phi_2|^2+ \int_{\Si(\tau_1)} r^{-\de_0 -4}|\Phi_1|^2.
  \end{split}
  \eea
Indeed, \eqref{eq:general-integrated-estimate-e4} corresponds to (12.3.3) of \cite{GKS22}, whose proof in Section 12.3.1 of \cite{GKS22} does not depend on the size of $a$ and hence applies here.  

Next, recall from Proposition \ref{cor:systemoftransportequationsforPsibandAbfromqfb:chap12} that $(\Ab, \underline{\Psi})$ satisfies the following system of transport equations  
\beaa
\nabc_4(r\Psib)=\frac{q}{r\ov{q}}\pmb\phi_{-2}^{(2)}  +r\dk^{\leq 1}(\Ga_g\c\Ga_b), \qquad \nabc_4\left( \frac{q^4}{r^3}\Ab\right)= \frac{1}{r}\Psib+r\Ga_g\c\Ga_b.
\eeaa
Applying \eqref{eq:general-integrated-estimate-e4} with $\de_0=2-p\geq\de$, we immediately obtain, for $\de\leq p\leq 2-\de$,  
 \bea\label{eq:general-integrated-estimate-e4:applicationtosystemfoAbandundPsi}
 \bsplit
&\int_{\MM(\tau_1, \tau_2)}  r^{p-3} |(\Ab, \underline{\Psi})|^2 + \sup_{\tau\in[\tau_1,\tau_2]}\int_{\Si(\tau)} r^{p-4}|(\Ab, \underline{\Psi})|^2+  \int_{\AA\cup \Si_*(\tau_1,\tau_2)} r^{p-2}|(\Ab, \underline{\Psi})|^2  \\
  &\les \int_{\MM(\tau_1, \tau_2)} r^{p-5}|\pmb\phi_{-2}^{(2)}|^2+ \int_{\Si(\tau_1)} r^{p -4}|(\Ab, \underline{\Psi})|^2+\ep_0^2\tau_1^{-2-2\dec},
  \end{split}
  \eea
  which implies \eqref{eq:analogofequationparenthesis12dot2dot7parenthesisinGKS22} for $s=0$. 

Next, we commute the above system of transport equations satisfied by $(\Ab, \underline{\Psi})$ with $\Lieb_\T$, $q\DD$ and $r\nab_4$ which yields, in view of \eqref{commutatorbetweenLieTLieZandnabnab4nab3:1}
\beaa
\nabc_4(r\Lieb_T\Psib)=\frac{q}{r\ov{q}}\Lieb_T\pmb\phi_{-2}^{(2)}  +r\dk^{\leq 2}(\Ga_g\c\Ga_b), \qquad \nabc_4\left( \frac{q^4}{r^3}\Lieb_T\Ab\right)= \frac{1}{r}\Lieb_T\Psib+r\dk^{\leq 1}(\Ga_g\c\Ga_b),
\eeaa
in view of Lemma \ref{LEMMA:COMM-GEN-B} and \eqref{eq:Xi-Hb-chapter12}, for $\tau\leq\tau_*-2$, 
\beaa
\nabc_4(r(q\DD)\Psib) &=& \frac{q}{r\ov{q}}(q\DD)\pmb\phi_{-2}^{(2)} +O(r^{-2})\pmb\phi_{-2}^{(2)} +O(r^{-1})\Psib+r\dk^{\leq 2}(\Ga_g\c\Ga_b), \\ 
\nabc_4\left( \frac{q^4}{r^3}(q\DD)\Ab\right) &=& \frac{1}{r}(q\DD)\Psib  +O(r^{-2})\Psib +O(r^{-1})\Ab +r\dk^{\leq 1}(\Ga_g\c\Ga_b),
\eeaa
and
\beaa
\nabc_4(r(r\nab_4)\Psib) &=& \frac{q}{r\ov{q}}(r\nab_4)\pmb\phi_{-2}^{(2)} +O(r^{-1})\pmb\phi_{-2}^{(2)} +O(1)\Psib+r\dk^{\leq 2}(\Ga_g\c\Ga_b), \\ 
\nabc_4\left( \frac{q^4}{r^3}(r\nab_4)\Ab\right) &=& \frac{1}{r}(r\nab_4)\Psib  +O(r^{-1})\Psib +O(1)\Ab +r\dk^{\leq 1}(\Ga_g\c\Ga_b).
\eeaa
Applying \eqref{eq:general-integrated-estimate-e4} with $\de_0=2-p\geq\de$ to these system of transport equations, we infer, using also \eqref{eq:general-integrated-estimate-e4:applicationtosystemfoAbandundPsi} and Lemma \ref{lemma:basicpropertiesLiebTfasdiuhakdisug:chap9}, 
 \beaa
&&\int_{\MM(\tau_1, \tau_2)}  r^{p-3} |\dk(\Ab, \underline{\Psi})|^2 + \sup_{\tau\in[\tau_1,\tau_2]}\int_{\Si(\tau)} r^{p-4}|\dk(\Ab, \underline{\Psi})|^2+  \int_{\AA\cup \Si_*(\tau_1,\tau_2)} r^{p-2}|\dk(\Ab, \underline{\Psi})|^2  \\
  &\les& \int_{\MM(\tau_1, \tau_2)} r^{p-5}|\dk^{\leq 1}\pmb\phi_{-2}^{(2)}|^2+ \int_{\Si(\tau_1)} r^{p -4}|\dk^{\leq 1}(\Ab, \underline{\Psi})|^2+\ep_0^2\tau_1^{-2-2\dec}.
  \eeaa
Together with \eqref{eq:general-integrated-estimate-e4:applicationtosystemfoAbandundPsi}, this implies \eqref{eq:analogofequationparenthesis12dot2dot7parenthesisinGKS22} for $s=1$. The general case then follows by iteration. This concludes the proof of Proposition \ref{prop:analogofProposition12dot2dot9inGKS22}.
\end{proof}

\noindent{\bf Step 4.} We now claim the following proposition which is the analog of Proposition 12.4.6 in \cite{GKS22}.
\begin{proposition}\lab{prop:decayofpr2tqfb:onMMfirst}
The following decay estimate holds, for $1\leq\tau_1<\tau_2\leq\tau_*-2$ and $s\leq\kst -16$,
\beaa
\int_{\MM(\tau_1, \tau_2)}r^{-1-\de}|\dk^{\leq s}\Lieb_\T^2(\Ab, \underline{\Psi})|^2 &\les& \ep_0^2 \tau_1^{-2-2\dec}.
 \eeaa
 Furthermore, there exists  a sequence of times $\tau^{(j)}$ such that, for $s\leq\kst-16$, 
\beaa
\int_{\Si(\tau^{(j)})}r^{-2}|\dk^{\leq s}\Lieb_\T^2(\Ab, \underline{\Psi})|^2 &\les& \ep_0^2(\tau^{(j)})^{-2-2\dec}, \qquad \tau^{(j)}\sim 2^j.
\eeaa
\end{proposition}

\begin{proof}
The proof of Proposition 12.4.6 in \cite{GKS22} in Section 14.4.2 of \cite{GKS22} consists of a sequence of standard mean value arguments based on Corollary 12.2.5 in \cite{GKS22} and Proposition 12.2.9 in \cite{GKS22}, and does not use the smallness of $a$. Given that \eqref{eq:step1fluxSigamstarLiebTAbchap12:1'} is the analog of Corollary 12.2.5 in \cite{GKS22} and Proposition \ref{prop:analogofProposition12dot2dot9inGKS22} is the analog of Proposition 12.2.9 in \cite{GKS22}, the proof immediately extends to our setting of Proposition \ref{prop:decayofpr2tqfb:onMMfirst}. Finally, recall from Footnote \ref{footnote:onelossofderivativeforeq:analogofequationparenthesis12dot2dot7parenthesisinGKS22comparedtoGKS22} that \eqref{eq:analogofequationparenthesis12dot2dot7parenthesisinGKS22} loses a derivative compared to (12.2.7) in \cite{GKS22}. Given that (12.2.7) in \cite{GKS22} is used four times in the proof of Proposition 12.4.6 in \cite{GKS22}, this results in an additional loss of four derivatives for Proposition \ref{prop:decayofpr2tqfb:onMMfirst} compared to Proposition 12.4.6 in \cite{GKS22}.
\end{proof}

%%%%%%%%%%%%%%%%%%%%%%%%%%%%%%%%%%%%%%%%%%%%%%%%%%%%%%%

\paragraph{\textit{End of the proof of Theorem \ref{thm:restatementofTheoremM2}}.}

%%%%%%%%%%%%%%%%%%%%%%%%%%%%%%%%%%%%%%%%%%%%%%%%%%%%%%%

We now claim the following proposition which is the analog of Proposition 12.4.6 in \cite{GKS22}.
\begin{proposition}\lab{prop:decayofpr2tqfb}
The following decay estimate holds for $\Lieb_\T^2\Ab$, for $s\leq\kst -17$,
\bea
\F^s_{\Si_*}[\Lieb_\T^2\Ab](\tau, \tau_*) &\les& \ep_0^2\tau^{-2-2\dec},
\eea
where $\F^s_{\Si_*}$ denotes the flux on $\Si_*$. 
\end{proposition}

\begin{proof}
The proof of Proposition 12.4.6 in Section 12.4.3 \cite{GKS22} is done first on $\MM(1,\tau_*-2)$ and then on $\MM(\tau_*-2, \tau_*)$. The proof on $\MM(1,\tau_*-2)$ relies on Proposition 12.4.6 in \cite{GKS22} and  basic estimates for the transport system of Proposition \ref{cor:systemoftransportequationsforPsibandAbfromqfb:chap12} for $(\Ab, \underline{\Psi})$ that do not depend on the size of $a$. The proof thus immediately extends to the one of Proposition \ref{prop:decayofpr2tqfb:onMMfirst} on $\MM(1,\tau_*-2)$. Then, given that $\MM(\tau_*-2,\tau_*)$ is a local existence type region, the proof on $\MM(\tau_*-2, \tau_*)$ does not rely on the size of $a$ and is thus immediately valid in the full subextremal range. 
\end{proof}

The  rest of the proof of Theorem M2 in Sections 12.4.4 and 12.4.5 in \cite{GKS22} then rely on Proposition 12.4.6 in \cite{GKS22} and estimates on the far away hypersurface $\Si_*$ where the size of $a$ is never used. Thus, starting from Proposition \ref{prop:decayofpr2tqfb}, it immediately yields the conclusion of Theorem \ref{thm:restatementofTheoremM2}.

%%%%%%%%%%%%%%%%%%%%%%%%%%%%%%%%%%%%%%%%%%%%%%%%%

\section{Extension of the curvature estimates in Theorems M8 to the full subextremal range}
\lab{sec:higherordercurvatureestimatesforprovingThM8largea:00}

%%%%%%%%%%%%%%%%%%%%%%%%%%%%%%%%%%%%%%%%%%%%%%%%%

The goal of this section is to prove Theorem \ref{THEOREMM8:INTRO} on the extension of the curvature estimates of Theorem M8 in Section 3.7.1 of \cite{KS:Kerr}, stated in Theorem 9.65 of \cite{KS:Kerr}, to the full subextremal range.

%%%%%%%%%%%%%%%%%%%%%%%%%%%%%%%%%%%%%%%%%%%%%%%%%

\subsection{Proof of Theorem \ref{THEOREMM8:INTRO}}
\lab{sec:higherordercurvatureestimatesforprovingThM8largea}

%%%%%%%%%%%%%%%%%%%%%%%%%%%%%%%%%%%%%%%%%%%%%%%%%

%%%%%%%%%%%%%%%%%%%%%%%%%%%%%%%%%%%%%%%%%%%%%%%%%%%%%%%%%%%%%%%%%%%%%%%%%%     
   
\subsubsection{Geometric set-up for Section \ref{sec:higherordercurvatureestimatesforprovingThM8largea:00}}
\lab{sec:geometricsetupspacetimeMMforsections789curvatureestimatesTheoremM8}
  
%%%%%%%%%%%%%%%%%%%%%%%%%%%%%%%%%%%%%%%%%%%%%%%%%%%%%%%%%%%%%%%%%%%%%%%%%%

The geometric set-up for Section \ref{sec:higherordercurvatureestimatesforprovingThM8largea:00} is the one in Sections \ref{sec:smallnesconstants}--\ref{sec:geometricsetuponSigma*}. In particular, in view of Section \ref{section:SpacetimeMM-chap6}:
\begin{itemize}
\item the vacuum spacetime $(\MM, \g)$ in endowed with a global null pair $(e_3, e_4)$, a pair of constants $(a, m)$, scalar functions $(\tau, r, \th, \tphi)$ and complex horizontal 1-forms $\Jk$, $\Jk_{\pm}$,

\item the level sets $\Si(\tau)$ of $\tau$ are spacelike, and $\tau\in[1,\tau_*]$ on $\MM$ for some arbitrary large constant $\tau_*$,

\item the boundary of $\MM$ is given by \eqref{eq:prMM=AAcupSistarcupSi1cupSitaustar}, 

\item the spacelike future boundaries $\AA$ and $\Si_*$ are defined respectively by \eqref{eq:defAA=r=rplusoneminusedeh} and \eqref{eq:rangetauandronSigmastar}, and the minimum $r_*$ of $r$ on $\Si_*$ satisfies the dominance condition \eqref{eq:dominantconditionforrstarcomparedtotaustartonS*},

\item the spacetime $\MM$ is decomposed in $\Mint$ and $\Mext$ defined in \eqref{eq:definitionofMextandMint},

\item the subregions $\MM(\tt_1,\tt_2)$, $\MM_{red}$, $\Mtrap$ and $\Mntrap$ of $\MM$ are defined as in \eqref{eq:defofsubregionsofMM}. 
\end{itemize}
While we also assume the existence and the properties of the second global null frame $(e_3', e_4', e_1', e_2')$ in Section \ref{sec:secondglobalnullframeonMM}, as well as the specific geometric set-up of Section \ref{sec:geometricsetuponSigma*} on $\Si_*$, these assumptions play only an auxiliary role in Section \ref{sec:higherordercurvatureestimatesforprovingThM8largea:00} as they are solely needed in order for the energy-Morawetz estimates for up to 14 derivatives for solutions to Teukolsky equations in Theorem \ref{thm:main:MaSz26} to hold.

\begin{remark}\lab{rmk:choiceofmaingloblaframeinthewholepaper:sections789}
Recall from Remark \ref{rmk:choiceofmaingloblaframeinthewholepaper} that the global null frame $(e_3, e_4, e_1, e_2)$ used in Section \ref{sec:higherordercurvatureestimatesforprovingThM8largea:00} is the one in Section 9.6.1 in \cite{KS:Kerr}. Also, recall from Remark \ref{rmk:choiceofsecondgloblaframeinthewholepaper} that second global null frame $(e_3', e_4', e_1', e_2')$ 
coincides with $(e_3, e_4, e_1, e_2)$ in these sections.
\end{remark}

Finally we  provide a definition\footnote{The definition here differs from that in  {Section} \ref{sec:definitionofGabandGagfirsttime} in that we separate  the linearized curvature quantities,  denoted $\Rc_g, \Rc_b$, from the Ricci and metric coefficients, denoted by $\Ga_g, \Ga_b$.} for $\Ga_b$, $\Ga_g$, $\Rc_b$ and $\Rc_g$ used throughout Section \ref{sec:higherordercurvatureestimatesforprovingThM8largea:00} for the Ricci coefficients and curvature components associated to the  global null frame $(e_3, e_4, e_1, e_2)$. 
 
\begin{definition}\lab{def:GabGagRcbRcg:P3}
We define the quantities $\Ga_g, \Ga_b$
\beaa
\Ga_g&=&  \Big\{\Xi, \,  \omc, \, \trXc,\,   \Xh,\,  \Zc,\, \Hbc, \,  \trXbc,  \,  {\widecheck{e_4(r)}}, \,  r^{-1}\nab(r), \, \widecheck{e_4(\tau)}, \, r^{-1}\widecheck{\DD(\tau)},\, {e_4(\cos\th)}, \,  {r\widecheck{\nab_4\Jk}}\Big\},\\
\Ga_b&=&  \Big\{\Hc, \,  \Xbh, \,  {\omb}, \,  \Xib,\, r^{-1}\widecheck{e_3(r)}, \, r^{-1}\widecheck{e_3(\tau)}, \,    \widecheck{\DD(\cos\th)}, \,  e_3(\cos\th),  \,   r\,\widecheck{\ov{\DD}\c\Jk}, \,  r\,\DD\hot\Jk, \,  r\,\widecheck{\nab_3\Jk}\Big\}.
\eeaa
We also define the   following sets of curvature quantities
\beaa
\Rc_g =\Big\{ \Pc, \,   B, \,   A\Big\}, \qquad 
\Rc_b =\Big\{ r\Bb,  \Ab\Big\}. 
\eeaa
Moreover, we denote by $\dk^{\le k}\Ga_g$, $\dk^{\le k}\Ga_b$, $\dk^{\le k}\Rc_g$, $\dk^{\le k}\Rc_b$ all derivatives  up to order $k$ with respect to the weighted derivatives $\dk=\{\nab_3, r\nab_4, \dkb=r\nab \}$.
\end{definition}

\begin{remark}
 \lab{Remark:Ga_g=r^{-1}Ga_b}
In the norms $\Sk_k$ for Ricci and metric coefficients introduced in {Section} \ref{subsection:MainNormsM8}, $\Ga_g$  behaves precisely as $r^{-1} \Ga_b$. Also, in the norms $\Rk_k$ for curvature components introduced in {Section} \ref{subsection:MainNormsM8}, $\Rc_g$ behaves precisely as  $r^{-2}\Rc_b$. Since, in Section \ref{sec:higherordercurvatureestimatesforprovingThM8largea:00}, we are   only interested in deriving  estimates for the  top derivatives of  curvature using the norms $\Sk_k$ and $\Rk_k$, we will often identify\footnote{The sole properties for which one might have to distinguish are the low derivatives decay rates in $\tau$ that are only needed when  estimating   quadratic error terms. But, even in that case, we in fact only need $\tau^{-1-\dec}$ decay in $\MM_{trap}$ which is again consistent with these identifications.} $\Ga_g$ with $r^{-1} \Ga_b$ and $\Rc_g$ with  $r^{-2}\Rc_b$.
 \end{remark}

%%%%%%%%%%%%%%%%%%%%%%%%%%%%%%%%%%%%%%%%%%%%%%%%%%%%%%%%%%%%%%%%%%%%%%%%%

\subsubsection{Main norms for Section \ref{sec:higherordercurvatureestimatesforprovingThM8largea:00}}
\lab{subsection:MainNormsM8}

%%%%%%%%%%%%%%%%%%%%%%%%%%%%%%%%%%%%%%%%%%%%%%%%%%%%%%%%%%%%%%%%%%%%%%%%%

The norms introduced in this section are norms for the Ricci coefficients and curvature components associated to the  global null frame $(e_3, e_4, e_1, e_2)$ of $\MM$, i.e., the one in Section 9.6.1 in \cite{KS:Kerr}, see Remark \ref{rmk:choiceofmaingloblaframeinthewholepaper:sections789}. These are the norms  appearing  in Section 9.4.1 of \cite{KS:Kerr} in connection with the higher order curvatures estimates of Theorem {9.65} in \cite{KS:Kerr}.  

We start with the norms for Ricci coefficients and curvature components on $\Mext$.

\begin{definition}
\lab{def:exteriorGac.norms}
We define the following norms  for the Ricci coefficients in $\Mext$ 
\beaa
\Skext_k^2 := \sup_{\la \ge r_0}  \int_{r=\la}\Big(r^2\big|\dk^{\le k}(\Ga_g\setminus\{\trXbc\}) |^2+ \big| \dk^{\le k}(\Ga_b\setminus\{\Xib\})\big|^2 +r^{2-\dt}|\dk^{\le k}\trXbc |^2+ r^{-\dt}|\dk^{\le k}\Xib |^2\Big).
\eeaa
For convenience we introduce the notation 
 \bea\lab{eq:defintionofGabprimeandGagprime:partIII}
 \Ga_b':=\Ga_b\setminus \{\Xib\}, \qquad \Ga_g':=\Ga_g\setminus \{\trXbc\}. 
 \eea
\end{definition}

\begin{definition}
\lab{Definiition:Rkext}
We define  the following norms for the  curvature coefficients in $\Mext$
\beaa
\Rkext^2_k := \int_{\Mext} r^{3+\dt}|\dk^{\le k}(A, B)|^2  +r^{3-\dt}\big(|\dk^{\le k}\Pc |^2 +r^{-2} |\dk^{\le k}\Bb|^2 +r^{-4} |\dk^{\le k}\Ab|^2\big).
\eeaa
\end{definition}

Next, we continue with the norms for Ricci coefficients and curvature components on $\Mint$.
    
\begin{definition}
\lab{definit:norms-SkMint'}
We define the following norms\footnote{{As $r$ is bounded on $\Mint$, we do not need to include $r$-weights in the norms $\Skint_k$ and $\Rkint_k$ introduced below.}}  for the Ricci coefficients  in   $\Mint$
\beaa
\Skint_k^2 := \int_{\Mint}\big|\dk^{\le k}\Gac\big|^2, 
\eeaa
 where $\Gac$ denotes the set of all  linearized Ricci and metric coefficients, i.e.,   
  \beaa
\Gac &:=& \Big\{\Xib,\, \omb,\,  \trXbc,\,  \Xbh,\,  \widecheck{e_3(r)},\, \widecheck{e_3(\tau)}, \, e_3(\cos\th),\, \widecheck{\nab_3\Jk},\, \Zc,\,  \Hbc,\, \Hc,\,  \DD r,\, \widecheck{ \DD\tau},\, \widecheck{ \DD\cos \th},\, \widecheck{ \ov{\DD}\c \Jk},\,  \DD\hot\Jk,\,  \trXc,\, \Xh, \\
&&\,\,\,\,\omc,\,  \widecheck{e_4(r)},\, \widecheck{e_4(\tau)}, \, e_4(\cos\th),\, \widecheck{\nab_4 \Jk},\, \Xi\Big\}.
  \eeaa
 \end{definition}

\begin{definition}
\lab{definit:norms-RkMint'}
We define the following norms for the null  curvature  components  in $\Mint$
\beaa
\Rkint_k^2 =  \int_{\Mint}\Big( \big| \nab_{\widetilde{R}} \dk^{\le k-1}\Rc\big|^2
+|\dk^{\le k-1}\Rc|^2\Big) +\int_{\Mint_{\ntrap}}  \big| \dk^{\le k}\Rc\big|^2   +\sup_\tau\int_{\Mint\cap\Si(\tau)}  |\dk^{\le k}  \Rc|^2,
\eeaa
where the vectorfield $\widetilde{R}$ is defined by
\bea\lab{def:Rhat}
\widetilde{R} :=  \frac 1 2  \left( e_4-\frac{\De}{|q|^2}  e_3\right), 
\eea
and  where $\Rc$ denotes the set of all  linearized curvature components, i.e.,   
  \beaa
\Rc &:=& \Big\{\Ab,\,\, \Bb,\,\,  \Pc,\,\,  B,\,\,  A\Big\}.
  \eeaa
\end{definition}

Finally, we define global  norms  of $\MM$ as follows
\bea\lab{eq:defintionofglobalnorms}
\Sk_k:=\Skext_k+\Skint_k, \qquad \Rk_k:=\Rkext_k+\Rkint_k.
\eea

We conclude with the following lemma which provides basic properties of the above norms. 
\begin{lemma}
\lab{lemma:auxilliarynormsforGa_b}
The following estimates hold true:
\begin{enumerate}
\item We have, in  $\Mext$,  for $p> 1+\dt$.
\bea
\lab{eq:auxilliarynormsforGa_b1}
\int_{\Mext}\Big(r^{-p}  | \dk^{\le k}\Ga_b|^2 +r^{-p+2}  | \dk^{\le k}\Ga_g|^2\Big) &\les  &  r_0^{1-p+\dt}    \Skext_k^2.
\eea
\item We also have the stronger  estimate for the notations $\Ga_b'$ and $\Ga_g'$ introduced in \eqref{eq:defintionofGabprimeandGagprime:partIII}, for $p>1$,
\bea
\lab{eq:auxilliarynormsforGa_b1-strong}
\int_{\Mext} \Big(r^{-p }| \dk^{\le k}\Ga' _b|^2 +r^{-p+2 }| \dk^{\le k}\Ga' _g|^2\Big) &\les  &  r_0^{1-p}    \Skext_k^2.
\eea

\item  We have
\bea
\lab{eq:auxilliarynormsforRb:boundary}
\nn\sup_{\tau\in[1, \tau_*]}\int_{\Si(\tau)} r^{3-\dt}\Big(| \dk^{\le k}\Rc_g|^2 +r^{-4}  | \dk^{\le k}\Rc_b|^2\Big)\\
+\int_{\AA\cup\Si_*} r^{3-\dt}\Big(| \dk^{\le k}\Rc_g|^2 +r^{-4}  | \dk^{\le k}\Rc_b|^2\Big) &\les& \Rk_{k+1}\Rk_k, 
\eea
and, for  $ p>1$, 
\bea
\lab{eq:auxilliarynormsforGa_b2}
\nn\sup_{\tau\in[1, \tau_*]}\int_{\Si(\tau)} \Big(r^{-p}  | \dk^{\le k}\Ga_b'|^2 +r^{-p+2}  | \dk^{\le k}\Ga_g'|^2\Big)\\
\nn +\int_{\AA\cup\Si_*} \Big(r^{-p}  | \dk^{\le k}\Ga_b'|^2 +r^{-p+2}  | \dk^{\le k}\Ga_g'|^2\Big)\\
 \nn +\sup_{\tau\in[1, \tau_*]}\int_{\Si(\tau)}  \Big(r^{-p-\dt}  | \dk^{\le k}\Xib|^2+r^{-p-\dt+2}  | \dk^{\le k}\trXbc|^2\Big)\\
  +\int_{\AA\cup\Si_*}  \Big(r^{-p-\dt}  | \dk^{\le k}\Xib|^2+r^{-p-\dt+2}  | \dk^{\le k}\trXbc|^2\Big) &\les& \Sk_{k+1}\Sk_k. 
\eea
\end{enumerate}
\end{lemma}

\begin{proof}
See Lemma 13.5.5 in \cite{GKS22}.
\end{proof}

%%%%%%%%%%%%%%%%%%%%%

\subsubsection{Main assumptions} 
\lab{section:MainAss.ThmM8}

%%%%%%%%%%%%%%%%%%%%%

The goal of Section \ref{sec:higherordercurvatureestimatesforprovingThM8largea:00} is to prove the control of high order curvature estimates in Theorem M8 of \cite{KS:Kerr}. In this section, we recall the assumptions on which this proof rests.

%%%%%%%%%%%%%%%%%%%%%%%%

\paragraph{\textit{Control of the initial data}.} 

%%%%%%%%%%%%%%%%%%%%%%%%

We have the following control of the initial data norm of Definition \ref{def:initialdatanorm}
\bea\lab{eq:controlofinitialdataforThM8}
\Ik_{k_L}\les \ep_0.
\eea
This {result} has been proved in Theorem {9.46} of \cite{KS:Kerr}.

%%%%%%%%%%%%%%%%%%%%%

\paragraph{\textit{Iteration assumption}.} 

%%%%%%%%%%%%%%%%%%%%%

We make  the  iteration assumption for $J$  in the range  $\frac{k_L}{2}\leq J\leq k_L-1$,
\bea\lab{eq:iterationassumptiondiscussionThM8:bis}
\Sk_J+\Rk_J  & \les&\ep_J,
\eea
see the iteration assumption {(9.52)} in \cite{KS:Kerr}.

\begin{remark}
We refer to {(9.53)} in \cite{KS:Kerr} for the specific choice of $\ep_J$. We do not recall it here\footnote{{In fact, $\ep_J:=\ep_0+L_*(J)$ with $L_*(J)$ defined in {(9.53)} in \cite{KS:Kerr}.}} since it is irrelevant for the statements and proofs of Sections \ref{sec:higherordercurvatureestimatesforprovingThM8largea}--\ref{sec:energyMorawetzforABBbAb}.
\end{remark}

%%%%%%%%%%%%%%%%%%%%%%%

\paragraph{\textit{Bootstrap assumptions}.} 

%%%%%%%%%%%%%%%%%%%%%%%

{Relative to the global norms \eqref{eq:defintionofglobalnorms}, we make a bootstrap assumption on decay for low derivatives} to deal with trapping. Recall  the scalar function 
$\tau_{trap}$ defined by
\beaa
\tau_{trap} := \left\{\ba{lll}
1+\tau & \textrm{on} & \MM_{trap},\\
1& \textrm{on} & \Mntrap.
\ea\right.
\eeaa
Then, we assume that $(\Ga_g, \Ga_b)$ and $(A, B, \Pc, \Rc_b)$ satisfy the following estimates on $\MM$
\bea\lab{eq:auxlowderivativebootassdecayforchapte13}
\bsplit
r^{\frac{7}{2}+\dec}|\dk^{\leq k}(A, B)|+r^3|\dk^{\leq k}\Pc|+r^2|\dk^{\leq k}\Ga_g|+r|\dk^{\leq k}(\Rc_b, \Ga_b)| &\leq \frac{\ep}{\tau_{trap}^{1+\dec}}, \quad k\leq \frac{\kl}{2},
\end{split}
\eea
 {see\footnote{{Note that \eqref{eq:auxlowderivativebootassdecayforchapte13} is in fact weaker than the corresponding control in {(9.43)} of \cite{KS:Kerr}.}} {(9.43)} in \cite{KS:Kerr}. Also, we assume   
\bea\lab{eq:mainbootassforchapte13:onlyforkLover2derivativesfirst}
\Sk_k+\Rk_k &\leq& \ep, \quad k\le \frac{k_L}{2},
\eea
see {(9.42)} in \cite{KS:Kerr}.}

%%%%%%%%%%%%%%%%%%%%%%%%%%%%%%%%%%%%%%%%%%%%%%%%%

\paragraph{\textit{Identities satisfied by the global frame of $\MM$}.} 

%%%%%%%%%%%%%%%%%%%%%%%%%%%%%%%%%%%%%%%%%%%%%%%%%

We assume that the following properties hold for  the global null frame of $\MM$
{\bea\lab{eq:specialidentityforthegloablframeofMMinpartIII}
\Xi\in r^{-1}\Ga_g
\eea
and
\bea\lab{eq:specialidentityforthegloablframeofMMinpartIII:bis}
\Xi=0, \qquad \Hbc=0, \quad\textrm{on}\quad \{r\geq r_0+1\}\cap\{\tau\leq \tau_*-3\}.
\eea}

{\begin{remark}\lab{rmk:specialidentityforthegloablframeofMMinpartIII}
Recall from Remark \ref{rmk:choiceofmaingloblaframeinthewholepaper:sections789} that the global null frame used in Section \ref{sec:higherordercurvatureestimatesforprovingThM8largea:00} is constructed in {Section} 9.6.1 of \cite{KS:Kerr} and indeed satisfies \eqref{eq:specialidentityforthegloablframeofMMinpartIII} and  \eqref{eq:specialidentityforthegloablframeofMMinpartIII:bis}. The estimate \eqref{eq:specialidentityforthegloablframeofMMinpartIII} is used to control the error terms of some commutators appearing in Sections \ref{sec:higherordercurvatureestimatesforprovingThM8largea:00}, while the identities \eqref{eq:specialidentityforthegloablframeofMMinpartIII:bis} are used for the control of $\pmb\phi_{-2}^{(2)}$ in Theorem \ref{Thm:Estimates-forqf}.
\end{remark}}

%%%%%%%%%%%%%%%%%%%%%%%%%%%%%%%%%%%%%%%%%%%%%%

\subsubsection{Precise statement of Theorem \ref{THEOREMM8:INTRO}}

%%%%%%%%%%%%%%%%%%%%%%%%%%%%%%%%%%%%%%%%%%%%%%

We provide below a precise statement of Theorem \ref{THEOREMM8:INTRO}.

\begin{theorem}[Extension of Theorem 13.6.3 in \cite{GKS22} to the full subextremal range]
\lab{prop:rpweightedestimatesiterationassupmtionThM8}
Let $J$ {be} such that $\frac{k_L}{2}\leq J\leq k_L-1$. Assume that the spacetime $\MM$  as defined in Section \ref{sec:geometricsetupspacetimeMMforsections789curvatureestimatesTheoremM8}  
verifies the assumption \eqref{eq:controlofinitialdataforThM8} on initial data, the quantitative assumptions \eqref{eq:auxlowderivativebootassdecayforchapte13} \eqref{eq:mainbootassforchapte13:onlyforkLover2derivativesfirst}, and the iteration assumption \eqref{eq:iterationassumptiondiscussionThM8:bis}. Assume also that the assumptions for up to 16 derivatives in \eqref{eq:assumptionsonMMforpartII:secondglobalframeauxassfor15derivatives} \eqref{eq:assumptionsonSigmastarforpartII} hold\footnote{These auxiliary assumptions for up to 16 derivatives are solely needed to ensure that the energy-Morawetz estimates for solutions to Teukolsky equations of Theorem \ref{thm:main:MaSz26} hold. Also, recall from Remark \ref{rmk:choiceofmaingloblaframeinthewholepaper:sections789} that the global null frame appearing in \eqref{eq:assumptionsonMMforpartII:secondglobalframeauxassfor15derivatives} is the global null frame $(e_3, e_4, e_1, e_2)$ used in this section, i.e., the one in Section 9.6.1 in \cite{KS:Kerr}.} as well. Then,  we the following control of the curvature components holds  in $\MM$ 
\bea
\lab{eq:InteriorcurvEstimatesThmM8}
   \Rkint_{J+1}^2&\les&  r_0^{18}\Big(\ep_J^{\frac{1}{6}}(\Sk_{J+1}+\Rk_{J+1})^{\frac{11}{6}}+\ep_0^2\Big),\\
\lab{eq:ExteriorcurvEstimatesThmM8}
\Rkext_{J+1}^2 &\les& r_0^{3+\de_B}\Rkint^2_{J+1}+ r_0^{-\de_B} \Skext^2_{J+1} +\ep_J^2 +\ep_0^2 {+\ep(\Rk_{J+1}+\Sk_{J+1})^2},
\eea
where the constant in $\les$ is   independent of $r_0$.
\end{theorem}

\begin{remark}
The estimates \eqref{eq:InteriorcurvEstimatesThmM8} and \eqref{eq:ExteriorcurvEstimatesThmM8} imply 
\bea\lab{eq:TotalIntplusExtcurvEstimatesThmM8}
\Rk_{J+1} &\les& r_0^{\frac{21}{2}+\frac{\dt}{2}}\Big(\ep_J^{\frac{1}{12}}(\Sk_{J+1}+\Rk_{J+1})^{\frac{11}{12}}+\ep_0+\ep_J\Big) + r_0^{-\frac{\de_B}{2}} \Skext_{J+1} +\sqrt{\ep}\Sk_{J+1},
\eea
where the constant in $\les$ is independent of $r_0$. The estimates \eqref{eq:TotalIntplusExtcurvEstimatesThmM8} with the choice $\kl=k_{large}+7$ proves the extension of 
Theorem {9.65} in \cite{KS:Kerr} to the full extremal range\footnote{Note that \eqref{eq:TotalIntplusExtcurvEstimatesThmM8} differs slightly from the main estimate claimed in Theorem {9.65} of \cite{KS:Kerr}. This  discrepancy does not matter for the proof of Theorem M8 in \cite{KS:Kerr} since the only relevant property of \eqref{eq:TotalIntplusExtcurvEstimatesThmM8} is that $\Skext_{J+1}$ is multiplied by $r_0^{-\frac{\de_B}{2}}$ and that the power of $\Sk_{J+1}+\Rk_{J+1}$ is smaller than one.}. The control of $\Rk_{J+1}$ provided by \eqref{eq:TotalIntplusExtcurvEstimatesThmM8}, together with the control for $\Sk_{J+1}$ provided by {Propositions 9.50--9.53} in \cite{KS:Kerr}, allows to obtain the iteration assumption \eqref{eq:iterationassumptiondiscussionThM8:bis} with $J$ replaced by $J+1$ for $\ep_J$ given by {(9.52) (9.53)} in \cite{KS:Kerr}. This iteration procedure then concludes the proof of Theorem M8 of \cite{KS:Kerr}, see {Section} 9.4.7 in \cite{KS:Kerr} for the iteration procedure, and {Section} 9.4.3 in \cite{KS:Kerr} for the proof of Theorem M8.
\end{remark}

\begin{remark}\lab{rmk:theinteriorestimatesareinfactglobal}
The estimate \eqref{eq:InteriorcurvEstimatesThmM8} results in fact from  global energy-Morawetz estimates on $\MM$ which are then restricted to $\Mint$. 
\end{remark}

\begin{remark}\label{remark:whydoweputBAuptokLbyabuselagnguage!!}
The nonlinear error terms appearing in the proof of \eqref{eq:InteriorcurvEstimatesThmM8} and \eqref{eq:ExteriorcurvEstimatesThmM8} may all be controlled by\footnote{{The term $\ep\int_1^{+\infty}\frac{d\tau}{\tau^{1+\dec}}$ on the first line of \eqref{eq:nonlinearestimate:chap9addedcomment} comes from the control of $\err_{J+1}$ in $\Mtrap$ where we rely on  \eqref{eq:auxlowderivativebootassdecayforchapte13}.}} 
\bea\lab{eq:nonlinearestimate:chap9addedcomment}
\err_{J+1} &\les& \left(\ep\int_1^{+\infty}\frac{d\tau}{\tau^{1+\dec}}+\Rk_{\frac{k_L}{2}}+\Sk_{\frac{k_L}{2}}\right)(\Rk_{J+1}+\Sk_{J+1})^2 \les \ep(\Rk_{J+1}+\Sk_{J+1})^2,
\eea
where we used \eqref{eq:mainbootassforchapte13:onlyforkLover2derivativesfirst} in the second estimate. The systematic use of \eqref{eq:nonlinearestimate:chap9addedcomment} would result in carrying the term $\ep(\Rk_{J+1}+\Sk_{J+1})^2$ on the RHS of all estimates throughout Sections \ref{sec:energyMorawetzforPc} \ref{sec:energyMorawetzforABBbAb}. Thus, to simplify the exposition, we will ignore these type of error terms in all subsequent estimates involved in the proof of Theorem \ref{prop:rpweightedestimatesiterationassupmtionThM8}.
\end{remark}

%%%%%%%%%%%%%%%%%%%%%%%%%%%%%%%%%%%%%%%%%%%%%%%%%%%%%%%%%

\subsubsection{Proof of Theorem \ref{prop:rpweightedestimatesiterationassupmtionThM8}}

%%%%%%%%%%%%%%%%%%%%%%%%%%%%%%%%%%%%%%%%%%%%%%%%%%%%%%%%%

The proof of Theorem \ref{prop:rpweightedestimatesiterationassupmtionThM8} relies on the following four intermediary results. 

\begin{theorem}
\lab{theorem:Morawetz-EnergyPc}
The following estimates hold true in $\MM=\MM(1, \tau_*)$
\bea\lab{eq:conclusionchapter14forlargea}
\BEF_\de^J[r^2\Pc]   &\les& r_0^{15}\Big( \Sk_{J+1} \Sk_J +\Rk_{J+1}\Rk_J +\ep_J^2+\ep_0^2\Big),
\eea
where the norm $\BEF^k_p$ has been introduced in {Section} \ref{subsection:basicnormsforpsi}.
\end{theorem}

  \begin{proposition}
     \lab{proposition:EstimatesBBb-interior}
    The following estimates hold true in 
     $\MM=\MM(1, \tau_*) $
         \bea
 \lab{eq:mainestimateBBb-M8}
 \BEF^J_\de[r^2B]+\BEF^J_\de[\Bb] &\les&   \de_{J+1}[\Pc] +\Sk_J \Sk_{J+1} +\Rk_{J} \Rk_{J+1} +\ep_0^2 +\Sk_{J+1}^{\frac{4}{3}}\Rk_{J+1}^{\frac{1}{3}}\Rk_{J}^{\frac{1}{3}}+\Sk_{J+1}\sqrt{ \de_{J+1}[\Pc]}\nn\\
 &&+\ep_J^2+\Sk_{J+1}\left(\de_{J+1}[\Pc] +\Sk_J \Sk_{J+1} +\Rk_{J} \Rk_{J+1} +\ep_0^2\right)^{\frac{1}{2}},
 \eea
   where
   \beaa
   \de_{J+1}[\Pc] :=\BEF_\de^{J}[r^2 \Pc ].
   \eeaa
   \end{proposition}
   
    \begin{proposition}
     \lab{proposition:EstimatesAAb-interior}
    The following estimates hold true in  $\MM=\MM(1, \tau_*)$
     \bea
 \lab{eq:proposition-EstimatesAAb-interior}
  \bsplit
   \BEF^J_\de[ r^2A]\les& \de_{J+1}[B] +\Sk_J \Sk_{J+1} +\Rk_{J} \Rk_{J+1} +\ep_0^2 +\Sk_{J+1}^{\frac{4}{3}}\Rk_{J+1}^{\frac{1}{3}}\Rk_{J}^{\frac{1}{3}}+\Sk_{J+1}\sqrt{ \de_{J+1}[B]}\\
 &+\ep_J^2+\Sk_{J+1}\left(\de_{J+1}[B] +\Sk_J \Sk_{J+1} +\Rk_{J} \Rk_{J+1} +\ep_0^2\right)^{\frac{1}{2}},\\
    \BEF^J_\de[\Ab]\les& \de_{J+1}[\Bb] +\Sk_J \Sk_{J+1} +\Rk_{J} \Rk_{J+1} +\ep_0^2 +\Sk_{J+1}^{\frac{4}{3}}\Rk_{J+1}^{\frac{1}{3}}\Rk_{J}^{\frac{1}{3}}+\Sk_{J+1}\sqrt{ \de_{J+1}[\Bb]}\\
 &+\ep_J^2+\Sk_{J+1}\left(\de_{J+1}[\Bb] +\Sk_J \Sk_{J+1} +\Rk_{J} \Rk_{J+1} +\ep_0^2\right)^{\frac{1}{2}},
   \end{split}
   \eea
   where 
   \beaa
\de_{J+1}[ B] :=\BEF_\de^J[ r^2 B],\qquad \de_{J+1}[ \Bb] := \BEF_\de^J[\Bb].
\eeaa
    \end{proposition}  
   
\begin{remark}
Theorem \ref{theorem:Morawetz-EnergyPc} and Propositions \ref{proposition:EstimatesBBb-interior} and \ref{proposition:EstimatesAAb-interior} correspond respectively to the extension of Theorem 14.1.3, and Propositions 15.1.1 and 15.1.2 in \cite{GKS22} to all subextremal angular momenta.
\end{remark}  
   
 \begin{theorem}
\lab{THM:MAINRESULTMEXT}
The following estimate holds  
\bea
\Rkext_{J+1}^2 &\les& r_0^{3+\de_B}\Rkint^2_{J+1} + r_0^{-\de_B}\Skext^2_{J+1} +\ep_J^2 +\ep_0^2,
\eea
with $\Rkext_k$, $\Rkint_k$ and $\Skext_k$  defined as in  {Section} \ref{subsection:MainNormsM8}.
\end{theorem}

\begin{proof}
Theorem \ref{THM:MAINRESULTMEXT} corresponds to Theorem 16.1.1 in \cite{GKS22} whose proof, based on $r^p$ weighted estimates for the Bianchi system for $r\geq r_0$ with $r_0\gg m$ large enough, is readily valid for all angular momenta and hence yields the proof of Theorem \ref{THM:MAINRESULTMEXT}.
\end{proof}  

We are now ready to prove Theorem \ref{prop:rpweightedestimatesiterationassupmtionThM8}.
\begin{proof}[Proof of Theorem \ref{prop:rpweightedestimatesiterationassupmtionThM8}]
Combining Theorem \ref{theorem:Morawetz-EnergyPc} and Propositions \ref{proposition:EstimatesBBb-interior} \ref{eq:proposition-EstimatesAAb-interior}, we obtain
\beaa
\BEF^J_\de[A, B, \Pc, \Bb, \Ab] \les r_0^{15}\Big(\ep_J^{\frac{1}{6}}(\Sk_{J+1}+\Rk_{J+1})^{\frac{11}{6}}+\ep_0^2\Big).
\eeaa
Also, in view of the definition of $\Rkint_k$ in {Section} \ref{subsection:MainNormsM8}, and the one of $\B^k_\de$ in {Section} \ref{subsection:basicnormsforpsi}, we have
\beaa
\Rkint_{J+1}^2 &\les& r_0^3 \BEF^J_\de[A, B, \Pc, \Bb, \Ab].
\eeaa
We deduce 
\beaa
\Rkint_{J+1}^2 \les r_0^{18}\Big(\ep_J^{\frac{1}{6}}(\Sk_{J+1}+\Rk_{J+1})^{\frac{11}{6}}+\ep_0^2\Big)
\eeaa
which concludes the proof of  \eqref{eq:InteriorcurvEstimatesThmM8}. Since Theorem \ref{THM:MAINRESULTMEXT} implies \eqref{eq:ExteriorcurvEstimatesThmM8}, this concludes the proof of Theorem \ref{prop:rpweightedestimatesiterationassupmtionThM8}. 
\end{proof}

The rest of this section is dedicated to the proof of Theorem \ref{theorem:Morawetz-EnergyPc} and Propositions \ref{proposition:EstimatesBBb-interior} and \ref{proposition:EstimatesAAb-interior}. Theorem \ref{theorem:Morawetz-EnergyPc} is proved in Section \ref{sec:energyMorawetzforPc}, Proposition \ref{proposition:EstimatesBBb-interior} is proved in Section \ref{sec:proofof:proposition:EstimatesBBb-interior}, 
and Proposition \ref{proposition:EstimatesAAb-interior} is proved in Section \ref{sec:proofof:proposition:EstimatesAAb-interior}.

%%%%%%%%%%%%%%%%%%%%%%%%%%%%%%%%%%%%%%%%%%%%%%%%%%%%%%%

\subsection{Energy-Morawetz estimates for $\Pc$ (Proof of Theorem \ref{theorem:Morawetz-EnergyPc})}
\lab{sec:energyMorawetzforPc}

%%%%%%%%%%%%%%%%%%%%%%%%%%%%%%%%%%%%%%%%%%%%%%%%%%%%%%%

%%%%%%%%%%%%%%%%%%%%%%%%%%%%%%%%%

\subsubsection{Preliminaries for the control of $\Pc$}

%%%%%%%%%%%%%%%%%%%%%%%%%%%%%%%%%

In this section, we collect several ingredients needed for the proof of Theorem \ref{theorem:Morawetz-EnergyPc}. We start with the following weighted  estimates for the inhomogeneous scalar wave equation on $\MM$. 
 \begin{theorem}
\lab{Prop:scalarwavePsi-M8}
Let $\psi$ be a solution to the inhomogeneous scalar wave equation on $\MM$
\bea\lab{eq:inhomogeouswaveequationusedforcontrolofPcinsection8}
\square_\g\psi&=& N.
\eea 
Then, the following weighted estimates hold true in $\MM$, for $1\leq\tau_1<\tau_2\leq\tau_*$ and $1\leq s\leq\kl-1$, 
\bea\lab{eq:Estimatesforpsi-M8-2}
\BEF^s_\de[\psi](\tau_1, \tau_2) &\les& \E^s_\de[\psi](\tau_1)+\NN_\de^s[\psi, N](\tau_1, \tau_2),
\eea
where the norms $\BEF^s_\de[\psi](\tau_1, \tau_2)$, $\E^s_\de[\psi](\tau_1)$ and $\NN_\de^s[\psi, N](\tau_1, \tau_2)$ are defined as in Section \ref{subsection:basicnormsforpsi}. 
\end{theorem}

\begin{proof}
First, we derive the following energy-Morawetz estimates for the  inhomogeneous scalar wave equation \eqref{eq:inhomogeouswaveequationusedforcontrolofPcinsection8} in $\MM$, for $1\leq\tau_1<\tau_2\leq\tau_*$ and $1\leq s\leq\kl-1$,
\bea\lab{eq:Estimatesforpsi-M8-2:energyMorawetz}
\EMF^s_\de[\psi](\tau_1, \tau_2) &\les& \E^s[\psi](\tau_1)+\widetilde{\NN}_\de^s[\psi, N](\tau_1, \tau_2)+\ep\B_\de^s[\psi](\tau_1, \tau_2),
\eea
where the norms $\EMF^s_\de[\psi](\tau_1, \tau_2)$, $\E^s[\psi](\tau_1)$, $\widetilde{\NN}_\de^s[\psi, N](\tau_1, \tau_2)$ and $\B_\de^s[\psi](\tau_1, \tau_2)$ are as in Section \ref{subsection:basicnormsforpsi} except for the Morawetz norm $\M_\de[\psi](\tau_1, \tau_2)$ which is defined for a scalar field $\psi$ as follows 
\beaa
\M_\de[\psi](\tau_1, \tau_2) &=& \int_{\MM(\tau_1, \tau_2)}\left(\frac{|\nab_3\psi|^2}{r^{1+\de}} +\frac{|\nabla\psi|^2}{r}\right)+\int_{\MM(\tau_1,\tau_2)}\left(\frac{|\nab_{\Rhat}\psi|^2}{r^{1+\de}} +\frac{|\psi|^2}{r^{3+\de}}\right).
\eeaa
The proof of \eqref{eq:Estimatesforpsi-M8-2:energyMorawetz} for $s=1$, which is a simple extension of the main result of \cite{MaSz24} to the case where $\MM$ only extends to the spacelike hypersurface $\Si_*$, is given in Theorem 4.76 of \cite{Sze}. The extension to higher order derivatives $1\leq s\leq\kl-1$ proceeds similarly as in the proof of Theorem \ref{thm:main:MaSz26:extendhigherorderderivatives} (see also Section 3.6 in \cite{MaSz24}), which is the corresponding statement for Teukolsky, and is in fact much simpler. We then combine it with the analog of Proposition \ref{Proposition:Step3-Chap10} with $p=\de$ for \eqref{eq:inhomogeouswaveequationusedforcontrolofPcinsection8} and conclude the proof of \eqref{eq:Estimatesforpsi-M8-2} using $\ep>0$ small enough to absorb the term $\ep\B_\de[\psi](\tau_1, \tau_2)$ on the RHS of \eqref{eq:Estimatesforpsi-M8-2:energyMorawetz}.
\end{proof}

We  will also  make use of the  following control of $\pmb\phi_{-2}^{(2)}$.
\begin{theorem}
\lab{Thm:Estimates-forqf}
We have,  for $k\le k_L-3$,
\bea
\lab{eq:Estimatesforqf-M8}
\BEF^{k} _{\de}[\pmb\phi_{-2}^{(2)}]{(1,\tau_*-3)} &\les \ep_0^2. 
\eea
\end{theorem} 

\begin{proof}
First we easily obtain the following analog of \eqref{MainEnerMora:psi:minus2case:extendhigherorderderivatives}, for $1\leq\tau_1<\tau_2\leq\tau_*-3$ and $14\leq\reg\leq\kl -3$, 
\begin{align}\lab{MainEnerMora:psi:minus2case:extendhigherorderderivatives:analogvaliduptotoporder}
\nn&\sum_{p=0}^2\Big(\EMF_\de^{\reg}[\pmb\phi_{-2}^{(p)}](\tau_1, \tau_2)+\EMF^{\reg}_{r\leq r_+(1+\dred)}[\nab_4^p\Ab](\tau_1, \tau_2)\Big)\\
\nn\les& \sum_{p=0}^2\E^{\reg}[\pmb\phi_{-2}^{(p)}](\tau_1)+\sum_{p=0}^2\E^{\reg}_{r\leq r_+(1+\dred)}[\nab_4^p\Ab](\tau_1)+\int_{\Si(\tau_1)}|\dk^{\leq\reg-2}\N_{T,-2}^{(0)}|^2\nn\\
&+\sum_{p=0}^1\int_{\MM(\tt_1, \tt_2)}r^{-1+\de}|\dk^{\leq\reg+1}\N_{T,-2}^{(p)}|^2+\sum_{p=0}^2\widetilde{\NN}_\de^{\reg}[\pmb\phi_{-2}^{(p)}, \widetilde{\N}_{W,-2}^{(p)}](\tt_1, \tt_2) +\ep\sum_{p=0}^2\B^{\reg}_\de[\pmb\phi_{-2}^{(p)}](\tau_1, \tau_2)+ \ep_0^2,
\end{align}
where the norms $\EMF^{\reg}_\de[\c](\tau_1,\tau_2)$, $\E^{\reg}[\c](\tau_1)$, $\B^{\reg}_\de[\c](\tau_1,\tau_2)$ and $\widetilde{\mathcal{N}}^{\reg}_\de[\c, \c](\tau_1, \tau_2)$ have been introduced in Section \ref{subsection:basicnormsforpsi}, where $\widetilde{\N}_{W,s}^{(p)}$, $s=\pm2$, $p=0,1,2$, are introduced in \eqref{eq:TensorialTeuSysandlinearterms:rescaleRHScontaine2:general:Kerrperturbation:alternateformnullframeinsteadcoordvectorfield} and $\N_{T,s}^{(p)}$, $s=\pm2$, $p=0,1$, are introduced in \eqref{eq:transportequationsins=plus2andminus2caseforp=0and1}. We combine \eqref{MainEnerMora:psi:minus2case:extendhigherorderderivatives:analogvaliduptotoporder} with Proposition \ref{prop:rpweightedestimatesforcombinednomrBEFsppmbphiminus2forallp} applied with $p=\de$, as in the proof of Corollary \ref{cor:rpweightedestimatesforcombinednomrBEFsppmbphiminus2forallp}, which yields, using also the fact that $\Xi=0$ on $\{r\geq r_0+1\}\cap\{\tau\leq \tau_*-3\}$ in view of \eqref{eq:specialidentityforthegloablframeofMMinpartIII:bis}, for $1\leq\tau_1<\tau_2\leq\tau_*-3$ and $14\leq s\leq \kl-3$,
\bea\lab{eq:analogofcor:rpweightedestimatesforcombinednomrBEFsppmbphiminus2forallp:sections789}
\nn\BEF^s_{\de}[\pmb\phi_{-2}](\tau_1, \tau_2) &\les& \E^s_{p}[\pmb\phi_{-2}](\tau_1)+\NN^s_{\de}[\pmb\phi_{-2}^{(2)}, \widetilde{\N}_{W,-2}^{(2)}]{(\tau_1, \tau_2)} +\widetilde{\NN}^s_\de[\pmb\phi_{-2}^{(0)}, \widetilde{\N}_{W,-2}^{(0)}](\tt_1, \tt_2)\\
&&+\widetilde{\NN}^{s}_\de[\pmb\phi_{-2}^{(1)}, \widetilde{\N}_{W,-2}^{(1)}](\tt_1, \tt_2)+\ep^2_0.
\eea

Next, we have, in view of the definition of $\widetilde{\NN}^{s}_\de[\c,\c]$ and $\NN^{s}_\de[\c,\c]$, 
\beaa
\nn&&\NN^s_{\de}[\pmb\phi_{-2}^{(2)}, \widetilde{\N}_{W,-2}^{(2)}]{(\tau_1, \tau_2)} +\widetilde{\NN}^s_\de[\pmb\phi_{-2}^{(0)}, \widetilde{\N}_{W,-2}^{(0)}](\tt_1, \tt_2)+\widetilde{\NN}^{s}_\de[\pmb\phi_{-2}^{(1)}, \widetilde{\N}_{W,-2}^{(1)}](\tt_1, \tt_2)\\
\nn&\les& \Big(\BEF^s_\de[\pmb\phi_{-2}^{(2)}](\tau_1, \tau_2)\Big)^{\frac{1}{2}}\left(\int_{\tau_1}^{\tau_2}\|\dk^{\leq s}\widetilde{\N}_{W,-2}^{(2)}\|_{L^2(\Si_{trap}(\tau))}+\left(\int_{\MM(\tau_1, \tau_2)}r^{\de+1}|\dk^{\leq s}\widetilde{\N}_{W,-2}^{(2)}|^2\right)^{\frac{1}{2}}\right)\\
&&+\sum_{p=0}^2\left(\left(\int_{\tau_1}^{\tau_2}\|\dk^{\leq s}\widetilde{\N}_{W,-2}^{(p)}\|_{L^2(\Si_{trap}(\tau))}\right)^2+\int_{\MM(\tau_1, \tau_2)}r^{\de+1}|\dk^{\leq s}\widetilde{\N}_{W,-2}^{(p)}|^2\right)
\eeaa
Together with the structure of $\widetilde{\N}_{W,-2}^{(p)}$, $p=0,1,2$, in \eqref{eq:structureofNpWminus2forp=012matchingtheoneinGKS22chap12}, and the control provided by the bootstrap assumptions {\eqref{eq:mainbootassforchapte13:onlyforkLover2derivativesfirst}} {and}  \eqref{eq:auxlowderivativebootassdecayforchapte13}, we infer, for $1\leq\tau_1<\tau_2\leq\tau_*-3$ and $s\le k_L-3$,
\beaa
\nn&&\NN^s_{\de}[\pmb\phi_{-2}^{(2)}, \widetilde{\N}_{W,-2}^{(2)}]{(\tau_1, \tau_2)} +\widetilde{\NN}^s_\de[\pmb\phi_{-2}^{(0)}, \widetilde{\N}_{W,-2}^{(0)}](\tt_1, \tt_2)+\widetilde{\NN}^{s}_\de[\pmb\phi_{-2}^{(1)}, \widetilde{\N}_{W,-2}^{(1)}](\tt_1, \tt_2)\\
\nn&\les& \Big(\BEF^s_\de[\pmb\phi_{-2}^{(2)}](\tau_1, \tau_2)\Big)^{\frac{1}{2}}\ep(\Sk_{s+3}+\Rk_{s+3})+\ep^2(\Sk_{s+3}+\Rk_{s+3})^2,
\eeaa
where the estimate for $\widetilde{\N}_{W,-2}^{(2)}$ is proved in Theorem 14.1.6 in \cite{GKS22} (with $N_{\err}=\widetilde{\N}_{W,-2}^{(2)}$ in the notations of the proof of Theorem 14.1.6 in \cite{GKS22}), and where that proof also applies to the estimates for $\widetilde{\N}_{W,-2}^{(p)}$, $p=0,1$, in view of their similar and simpler structure in \eqref{eq:structureofNpWminus2forp=012matchingtheoneinGKS22chap12}. Plugging this estimate in \eqref{eq:analogofcor:rpweightedestimatesforcombinednomrBEFsppmbphiminus2forallp:sections789}, we infer, for $1\leq\tau_1<\tau_2\leq\tau_*-3$ and $14\leq s\leq \kl-3$,
\beaa
\nn\BEF^s_{\de}[\pmb\phi_{-2}](\tau_1, \tau_2) &\les& \E^s_{p}[\pmb\phi_{-2}](\tau_1)+\ep^2(\Sk_{s+3}+\Rk_{s+3})^2+\ep^2_0.
\eeaa
In view of our convention in Remark \ref{remark:whydoweputBAuptokLbyabuselagnguage!!}, we may ignore the terms of the type $\ep(\Sk_{s+3}+\Rk_{s+3})^2$, and we thus deduce, for $14\leq s\leq \kl-3$, 
\beaa
        \BEF_\de^s[\pmb\phi_{-2}]{(1, \tau_*-3)} \les  \ep_0^2
 \eeaa
which implies the stated estimate. This concludes the proof of Theorem \ref{Thm:Estimates-forqf}.
\end{proof}

We will in addition rely on the following wave equation for $P$ which follows from Bianchi identities.
\begin{lemma}
The curvature component $P$ satisfies the following scalar wave equation
\bea\lab{eq:P-WaveEq-M8}
\nn\square_\g P &=&  \tr X\nab_3P +\ov{\tr\underline{X}}\nab_4P  -\ov{H}\c\DD P  - \Hb\c\ov{\DD}P + V P  \\
&& + r^{-3}\mathfrak{d}^{\leq 1}(\Gamma_b \c\widecheck{R}_b)-\underline{A}\c \ov{A}, \qquad V:=\frac{3}{2}\Big[ \ov{\tr\Xb} \tr X+2P  -2 \Hb \c  \ov{H}  \Big].
\eea
\end{lemma}

\begin{proof}
See Lemma 5.5.1 and (14.2.1) in \cite{GKS22}.
\end{proof}

The following lemma will also be useful.
\begin{lemma}
\lab{Le:squareq^2Psi}
Let $\Psi$ {be} a scalar function solution to the following wave equation 
\bea
\lab{eq:Le-squareq^2Psi}
\square_\g\Psi &=&  \tr X\nab_3\Psi +\ov{\tr\Xb}\nab_4\Psi  -\ov{H}\c\DD\Psi  - \Hb\c\ov{\DD}\Psi +V\Psi+F, 
\eea
where $V$ is a potential and $F$ a scalar function. Then, we have
\beaa
\square_\g(q^2\Psi) &=& \Big[V +q^{-2}\square_\g(q^2)\Big]q^2\Psi +r\Ga_b\c \dk\Psi+q^2F.
\eeaa
\end{lemma}

\begin{proof}
See Lemma 5.5.3 in \cite{GKS22}.
\end{proof}

The following lemma provides an identity relating  $\DDc\hot \DDc\Pc$ to  $\pmb\phi_{-2}^{(2)}$.
  \begin{lemma}
  \lab{Lemma:Formula-qf-DDhotDDPc}
  The following relation between $\pmb\phi_{-2}^{(2)}$ and $\Pc$ holds true.
\bea
 \lab{eq-Formula-qf-DDhotDDPc:P3}
 \bsplit
 \pmb\phi_{-2}^{(2)} =& \frac{1}{2}\ov{q}q^3 \DDc\hot \DDc\Pc +  \dk^{\leq 1}\Ga_b'  +O(a)\dk^{\leq 1}\Rc_b + O(ar) \dk^{\le 1}\Pc +r\dk^{\leq 1}(\Ga_b\c \Rc_b),
 \end{split}
 \eea
 where we recall that $\Ga_b'=\Ga_b\setminus\{\Xib\}$.
  \end{lemma}
  
\begin{proof}
See Lemma 4.11.7 in \cite{GKS22}, recalling from Remark \ref{rmk:compasisionphiminus2p=0withMaSz26andphiminus2p=2withqfGKS22} that $\pmb\phi_{-2}^{(2)}$ is denoted by $\qfb$ in \cite{GKS22}.
\end{proof}

We end this section with the following lemma comparing the  norms $\BEF_\de$ with the norms $\Rk_k$ and $\Sk_k$ for curvature components and Ricci coefficients.
\begin{lemma}\lab{lemma:basicrelationBEFdenormandnormsRkandSikforcurvatureandRiccicoeff:section789}
We have
\bea
 \lab{induction:hypoth-Pc}
 \B^{k-1}_\de[ r^2\Rc_g] +\B^{k-1}_\de[\Rc_b] \les \Rk_k^2, \qquad   \EF^{k-1}_\de[ r^2 \Rc_g] +\EF^{k-1}_\de[\Rc_b] \les \Rk_k\Rk_{k+1},
 \eea
 and
 \bea
\B^{k-1}_\de[\Ga_b'] + \B^{k-1}_\de[r^{-\de_B} \Ga_b] \les \mathfrak{S}_k^2, \qquad \EF^{k-1}_\de[\Ga_b'] +\EF^{k-1}_\de[r^{-\de_B} \Ga_b] \les \mathfrak{S}_k  \mathfrak{S}_{k+1}. 
\eea
\end{lemma}

\begin{proof}
See Lemmas 14.1.1 and  14.1.2 in \cite{GKS22}. 
\end{proof}

%%%%%%%%%%%%%%%%%%%%%%%%%%%%%%%%%%%%

\subsubsection{Proof of Theorem \ref{theorem:Morawetz-EnergyPc}}

%%%%%%%%%%%%%%%%%%%%%%%%%%%%%%%%%%%%

We assume that all the bootstrap, initial data assumptions and  the induction hypothesis made in {Section} \ref{section:MainAss.ThmM8} hold true. The goal of this section is to prove Theorem \ref{theorem:Morawetz-EnergyPc}. To this end, we rely in particular on the non-integrable Hodge estimates of Section \ref{sec:nonintegrableHodgeestimatesforcontrolPc} and the commutators of Section \ref{sec:commutatorsforcontrolPc}, and we proceed in the following steps. 

.

\noindent{\bf Step 1.} In view of the control for $\pmb\phi_{-2}^{(2)}$ provided by Theorem \ref{Thm:Estimates-forqf} and the identity \eqref{eq-Formula-qf-DDhotDDPc:P3} relating  $\DDc\hot \DDc\Pc$ to  $\pmb\phi_{-2}^{(2)}$, we immediately obtain, using also Lemma \ref{lemma:basicrelationBEFdenormandnormsRkandSikforcurvatureandRiccicoeff:section789}, 
\beaa
\BEF^{J-2}_\de[r^4{}^{(c)}\DD\hot {}^{(c)}\DD\widecheck{P}](1, \tau_*-3) &\les& \BEF^{J-1}_\de[\Ga_b']+\BEF^{J-1}_\de[\widecheck{R}_b]+\BEF^{J-1}_\de[r\widecheck{P}]+\ep_0^2\\
&\les& \mathfrak{S}_{J+1}\mathfrak{S}_J+\mathfrak{R}_{J+1}\mathfrak{R}_J+\ep_0^2.
\eeaa 
 Since the difference between ${}^{(c)}\DD$ and $\DD$ yields lower order terms, we infer
\beaa
\BEF^{J-2}_\de[r^4\DD\hot\DD\widecheck{P}](1, \tau_*-3) &\les& \mathfrak{S}_{J+1}\mathfrak{S}_J+\mathfrak{R}_{J+1}\mathfrak{R}_J+\ep_J^2+\ep_0^2.
\eeaa  
Relying on the above estimate, Steps 2--9 below yield the proof of Theorem \ref{theorem:Morawetz-EnergyPc} on $\MM(1, \tau_*-3)$, and we complete the proof in Step 10 by dealing with the local existence type region $\MM(\tau_*-3, \tau_*)$. For convenience, we may thus drop $(1, \tau_*-3)$ on all LHS in Steps 2--9, and reintroduce at the end of Step 9. In particular, we rewrite the above estimate simply as
\bea\lab{eq:proofcontrolPc:controlBEFdenormr4DDdhotDDdPc}
\BEF^{J-2}_\de[r^4\DD\hot\DD\widecheck{P}] &\les& \mathfrak{S}_{J+1}\mathfrak{S}_J+\mathfrak{R}_{J+1}\mathfrak{R}_J+\ep_J^2+\ep_0^2.
\eea

\noindent{\bf Step 2.} Next, we apply the non-integrable Hodge estimates of Lemma \ref{lemma:howtogofromDDhotDDtonab2} and infer from \eqref{eq:proofcontrolPc:controlBEFdenormr4DDdhotDDdPc}
\bea\lab{eq:proofcontrolPc:consaeuenceHodgeestimatecontrolBdenormr4nab2Pc}
\nn\B^{J-2}_\de[r^4\nab^2\widecheck{P}] &\les& \mathfrak{S}_{J+1}\mathfrak{S}_J+\mathfrak{R}_{J+1}\mathfrak{R}_J+\ep_J^2+\ep_0^2+\sqrt{\BEF^{J-1}_\de[r^2\widecheck{P}]\BEF^J_\de[r^2\widecheck{P}]}\\
&\les& \mathfrak{S}_{J+1}\mathfrak{S}_J+\mathfrak{R}_{J+1}\mathfrak{R}_J+\ep_J^2+\ep_0^2+\sqrt{\mathfrak{R}_{J+1}\mathfrak{R}_J\EF^J_\de[r^2\widecheck{P}]},
\eea  
where we have also used Lemma \ref{lemma:basicrelationBEFdenormandnormsRkandSikforcurvatureandRiccicoeff:section789}.

Also, in view of Lemma \ref{lemma:howtocontrolcrosstermsnabnabTandnabnabRpsi} and Lemma \ref{lemma:basicrelationBEFdenormandnormsRkandSikforcurvatureandRiccicoeff:section789}, we have 
\begin{align}\lab{eq:proofcontrolPc:controlBdenormr2rnabrThatandRhatPc}
\nn&\B^{J-2}_\de[r^2(r\nab)(\nab_{\widehat{T}}, \nab_{\widehat{R}})\widecheck{P}]\\
\nn\les& \sqrt{\B^{J-2}_\de[r^4\nab^2\widecheck{P}]\B^{J-2}_\de[r^2(\nab_{\widehat{T}}^2, \nab_{\widehat{R}}^2)\widecheck{P}]}+\sqrt{\EF^{J-1}_\de[r^2\widecheck{P}]\EF^J_\de[r^2\widecheck{P}]}+ \B^{J-2}_\de[r^4\nab^2\widecheck{P}]+\B^{J-1}_\de[r^2\widecheck{P}]\\
\les& \sqrt{\B^{J-2}_\de[r^4\nab^2\widecheck{P}]\B^{J-2}_\de[r^2(\nab_{\widehat{T}}^2, \nab_{\widehat{R}}^2)\widecheck{P}]}+\sqrt{\mathfrak{R}_{J+1}\mathfrak{R}_J\EF^J_\de[r^2\widecheck{P}]}+ \B^{J-2}_\de[r^4\nab^2\widecheck{P}]+\B^{J-1}_\de[r^2\widecheck{P}].
\end{align}
and 
\bea\lab{eq:proofcontrolPc:controlBdenormr2rnabThatnabRhatPc}
\nn\B^{J-2}_\de[r^2\nab_{\widehat{T}}\nab_{\widehat{R}}\widecheck{P}] &\les& \B^{J-2}_\de[r^2\nab_{\widehat{T}}^2\widecheck{P}]+\B^{J-2}_\de[r^2\nab_{\widehat{R}}^2\widecheck{P}] +\sqrt{\EF^{J-1}_\de[r^2\widecheck{P}]\EF^J_\de[r^2\widecheck{P}]}+\B^{J-1}_\de[r^2\widecheck{P}]\\
&\les& \B^{J-2}_\de[r^2\nab_{\widehat{T}}^2\widecheck{P}]+\B^{J-2}_\de[r^2\nab_{\widehat{R}}^2\widecheck{P}] +\sqrt{\mathfrak{R}_{J+1}\mathfrak{R}_J\EF^J_\de[r^2\widecheck{P}]}+\mathfrak{R}_{J+1}\mathfrak{R}_J.
\eea

\noindent{\bf Step 3.} Next, recall from Lemma 4.7.6 in \cite{GKS22} that we have for a scalar function $\psi$
\beaa
|q|^2 \square\psi =\frac{(r^2+a^2)^2}{\De} \big( -  \nab_{\widehat{T}}^2\psi+  \nab_{\widehat{R}}^2\psi \big) +2r \nab_{\widehat{R}}\psi+  |q|^2 \De\psi   + |q|^2  (\eta+\etab) \c \nab \psi  + r^2 \Ga_g \c \mathfrak{d} \psi,
\eeaa
where the vectorfields $\That$ and $\Rhat$ are defined by 
 \bea
   \lab{eq:ThatRhat-e_3e_4-Kerr}
 \That :=\frac 1 2 \left( \frac{|q|^2}{r^2+a^2} e_4+\frac{\De}{r^2+a^2}  e_3\right), \qquad \Rhat := \frac 1 2 \left( \frac{|q|^2}{r^2+a^2} e_4-\frac{\De}{r^2+a^2} e_3\right).
 \eea

Also, recall from (14.2.1) in \cite{GKS22} that $P$ satisfies 
\bea
\lab{eq:P-WaveEq-linearized}
\square_\g \Pc =  \tr X\nab_3\Pc  +\ov{\tr\Xb}\nab_4\Pc  -\ov{H}\c\DD \Pc   - \Hb\c\ov{\DD}\Pc + r^{-2}\Rc_g    + r^{-4}  \dk^{\le 1} \Ga_b  + r^{-3}\dk^{\leq 1}(\Ga_b \c\Rc_b) -\Ab\c \ov{A}.
\eea
We infer 
\beaa
 \B^{J-2}_\de[r^2(\nab_{\widehat{T}}^2\widecheck{P} - \nab_{\widehat{R}}^2\widecheck{P})] &\lesssim& \B^{J-2}_\de[r^2\De\widecheck{P}]+\B^{J-1}_\de[r\widecheck{P}]+\B^{J-2}_\de[\Rc_g]+\B^{J-2}_\de[r^{-2}  \dk^{\le 1} \Ga_b]+\ep_0^2
\eeaa
and hence, using \eqref{eq:proofcontrolPc:consaeuenceHodgeestimatecontrolBdenormr4nab2Pc} and Lemma \ref{lemma:basicrelationBEFdenormandnormsRkandSikforcurvatureandRiccicoeff:section789}, we deduce 
\bea\lab{eq:proofcontrolPc:controlBdenormr2rnabThat2minusnabRhat2Pc}
\B^{J-2}_\de[r^2(\nab_{\widehat{T}}^2\widecheck{P} - \nab_{\widehat{R}}^2\widecheck{P})] &\lesssim& \mathfrak{S}_{J+1}\mathfrak{S}_J+\mathfrak{R}_{J+1}\mathfrak{R}_J+\ep_J^2+\ep_0^2+\sqrt{\mathfrak{R}_{J+1}\mathfrak{R}_J\EF^J_\de[r^2\widecheck{P}]}.
\eea

\noindent{\bf Step 4.} Next, we derive a wave equation for $\Rhat^2P$, with the vectorfield $\Rhat$ defined in \eqref{eq:ThatRhat-e_3e_4-Kerr}. Commuting the wave equation \eqref{eq:P-WaveEq-M8} for $P$ with $\widehat{R}^2$, we have
\beaa
\square_\g(\widehat{R}^2P) &=& [\square_\g, \widehat{R}^2]P+\tr X\nab_3\nab_{\widehat{R}}^2P +\ov{\tr\underline{X}}\nab_4\nab_{\widehat{R}}^2P  -\ov{H}\c\DD\nab_{\widehat{R}}^2P  - \Hb\c\ov{\DD}\nab_{\widehat{R}}^2P + V\nab_{\widehat{R}}^2P\\
&& +O(r^{-2})\mathfrak{d}^{\leq 1}\nab_{\widehat{R}}^{\leq 1}P  + r^{-3}\mathfrak{d}^{\leq 3}(\Gamma_b \c\widecheck{R}_b)-\mathfrak{d}^{\leq 2}(\underline{A}\c \ov{A}).
\eeaa 
Now, in view of Corollary \ref{cor:commutatorbetweenRhatsquareandmod|q|2square}, we have
\beaa
[\Rhat^2, \square_\g]P &=& O(mr^{-2})\Rhat^{\leq 1}\square P +O(r^{-1})\Rhat^{\leq 1}\Delta P  +O(mr^{-2})\nab_{\Rhat}^{\leq 1}\nab\That P   +O(m^2r^{-3})\Rhat^3P  \\
&&+O(m^2r^{-4})\nab^2P +O(r^{-2})\Rhat^{\leq 1}\dk P +\dk^{\leq 2}(\Ga_g\dk P).
\eeaa
Plugging the wave equation \eqref{eq:P-WaveEq-M8} for $P$, we deduce 
 \beaa
  [\square_\g, \widehat{R}^2]P &=& O(r^{-1})\Rhat^{\leq 1}\Delta P  +O(mr^{-2})\nab_{\Rhat}^{\leq 1}\nab\That P   +O(m^2r^{-3})\Rhat^3P  \\
&&+O(m^2r^{-4})\nab^2P +O(r^{-2})\Rhat^{\leq 1}\dk P +\dk^{\leq 2}(\Ga_g\dk P)+ r^{-3}\mathfrak{d}^{\leq 2}(\Gamma_b \c\widecheck{R}_b).
 \eeaa 
  Coming back to the equation for $\widehat{R}^2P$ and plugging the commutator $[\square_\g, \widehat{R}^2]P$, we obtain
  \beaa
\square_\g(\widehat{R}^2P) &=& O(mr^{-2})\nab_{\Rhat}^{\leq 1}\nab\That P + O(r^{-1})\Rhat^{\leq 1}\Delta P     +O(m^2r^{-3})\Rhat^3P  +\tr X\nab_3\nab_{\widehat{R}}^2P +\ov{\tr\underline{X}}\nab_4\nab_{\widehat{R}}^2P\\
&&  -\ov{H}\c\DD\nab_{\widehat{R}}^2P  - \Hb\c\ov{\DD}\nab_{\widehat{R}}^2P + V\nab_{\widehat{R}}^2P  +O(m^2r^{-4})\nab^2P +O(r^{-2})\Rhat^{\leq 1}\dk P \\
&& +\dk^{\leq 2}(\Ga_g\dk P) + r^{-3}\mathfrak{d}^{\leq 3}(\Gamma_b \c\widecheck{R}_b)-\mathfrak{d}^{\leq 2}(\underline{A}\c \ov{A}).
\eeaa
Renormalizing by $q^2$ using Lemma \ref{Le:squareq^2Psi}, we infer
 \beaa
\square_\g(q^2\widehat{R}^2P) &=& O(m)\nab_{\Rhat}^{\leq 1}\nab\That P + O(r)\Rhat^{\leq 1}\Delta P     +O(m^2r^{-3})q^2\Rhat^3P  + Wq^2\widehat{R}^2P  +O(m^2r^{-2})\nab^2P\\
&& +O(1)\Rhat^{\leq 1}\dk P  +r^2\dk^{\leq 2}(\Ga_g\dk P) + r^{-1}\mathfrak{d}^{\leq 3}(\Gamma_b \c\widecheck{R}_b)-r^2\mathfrak{d}^{\leq 2}(\underline{A}\c \ov{A}),
\eeaa
where $W:=V+q^{-2}\square_\g(q^2)$ satisfies $W=O(mr^{-3})+r^{-1}\Ga_g$ so that 
\beaa
\square_\g(q^2\widehat{R}^2P) &=& O(m)\nab_{\Rhat}^{\leq 1}\nab\That P + O(r)\Rhat^{\leq 1}\Delta P     +O(m^2r^{-3})q^2\Rhat^3P    +O(m^2r^{-2})\nab^2P +O(1)\Rhat^{\leq 1}\dk P \\
&& +r^2\dk^{\leq 2}(\Ga_g\dk P) + r^{-1}\mathfrak{d}^{\leq 3}(\Gamma_b \c\widecheck{R}_b)-r^2\mathfrak{d}^{\leq 2}(\underline{A}\c \ov{A}).
\eeaa

\noindent{\bf Step 5.} We now linearize the wave equation for $q^2\widehat{R}^2P$ of Step 4 and obtain 
 \beaa
\square_\g\big(\widecheck{q^2\widehat{R}^2P}\big) &=& O(m^2r^{-3})\Rhat\big(\widecheck{q^2\Rhat^2\Pc}\big) + O(m)\nab_{\Rhat}^{\leq 1}\nab\That\Pc + O(r)\Rhat^{\leq 1}\Delta\Pc   +O(m^2r^{-2})\nab^2\Pc\\
&& +O(1)\Rhat^{\leq 1}\dk\Pc  +r^{-2}\mathfrak{d}^{\leq 2}\Ga_b  + r^{-1}\mathfrak{d}^{\leq 3}(\Gamma_b \c\widecheck{R}_b)-r^2\mathfrak{d}^{\leq 2}(\underline{A}\c \ov{A}).
\eeaa
We denote by $h=h(r, \cos\th)$ the regular $O(mr^{-3})$ scalar function in front of the first term on the RHS, i.e., we rewrite this wave equation as 
\bea\lab{eq:proofcontrolPc:waveeqwidecheckq2nabwidehatR2P}
\nn\square_\g(\widecheck{q^2\nab_{\widehat{R}}^2P}) &=& h(r, \cos\th)\Rhat\big(\widecheck{q^2\Rhat^2\Pc}\big) + O(m)\nab_{\Rhat}^{\leq 1}\nab\That\Pc + O(r)\Rhat^{\leq 1}\Delta\Pc   +O(m^2r^{-2})\nab^2\Pc\\ &&+O(1)\Rhat^{\leq 1}\dk\Pc  +r^{-2}\mathfrak{d}^{\leq 2}\Ga_b  + r^{-1}\mathfrak{d}^{\leq 3}(\Gamma_b \c\widecheck{R}_b)-r^2\mathfrak{d}^{\leq 2}(\underline{A}\c \ov{A}).
\eea 

Now, for a scalar function $f=f(r, \cos\th)$ such that $f=1+O(r^{-2})$ and a scalar function $\psi$, we have
\beaa
\square_\g(f\psi) &=& f\square_\g\psi -e_4(f)e_3(\psi) - e_3(f)e_4(\psi)+2\nab(f)\c\nab\psi+\square_\g(f)\psi\\
&=& f\square_\g\psi -e_4(r)\pr_rf e_3(\psi) -e_3(r)\pr_rf e_4(\psi) +O(r^{-3})\nab\psi +O(r^{-4})\psi +r^{-3}\dk^{\leq 1}(\Ga_b\psi)\\
&=&  f\square_\g\psi -\frac{\Delta}{|q|^2}\pr_rf e_3(\psi) +\pr_rf e_4(\psi) +O(r^{-3})\nab\psi +O(r^{-4})\psi +r^{-3}\dk^{\leq 1}(\Ga_b\psi) \\
&=& f\left(\square_\g\psi  +2\pr_r\log(f)\frac{r^2+a^2}{|q|^2}\widehat{R}\psi\right) + O(r^{-3})\nab\psi +O(r^{-4})\psi+r^{-3}\dk^{\leq 1}(\Ga_b\psi). 
\eeaa
where we have used the definition of $\Rhat$ in \eqref{eq:ThatRhat-e_3e_4-Kerr}. We thus choose 
\beaa
f(r, \cos\th) &=& \exp\left(\frac{1}{2}\int_r^{+\infty}\frac{(r')^2+a^2(\cos\th)^2}{(r')^2+a^2}h(r', \cos\th)dr'\right)
\eeaa
so that 
\beaa
2\pr_r\log(f)\frac{r^2+a^2}{|q|^2}=-h, \qquad f=1+O(r^{-2}).
\eeaa
With this choice of $f$, we infer
\beaa
\square_\g(f\psi) &=& f\left(\square_\g\psi  -h\widehat{R}\psi\right) + O(r^{-3})\nab\psi +O(r^{-4})\psi +r^{-3}\dk^{\leq 1}(\Ga_b\psi). 
\eeaa
Applying this to $\psi=\widecheck{q^2\widehat{R}^2P}$ and using the wave equation \eqref{eq:proofcontrolPc:waveeqwidecheckq2nabwidehatR2P} for $\widecheck{q^2\widehat{R}^2P}$, we infer
\bea\lab{eq:proofcontrolPc:waveeqwidecheckq2nabwidehatR2P:renormalizedwithf}
\nn\square_\g\big(f\widecheck{q^2\widehat{R}^2P}\big) &=&  O(m)\nab_{\Rhat}^{\leq 1}\nab\That\Pc +O(r^{-1})\nab_{\widehat{R}}\nabla\Rhat\widecheck{P} + O(r)\Rhat^{\leq 1}\Delta\Pc   +O(m^2r^{-2})\nab^2\Pc\\ &&+O(1)\Rhat^{\leq 1}\dk\Pc  +r^{-2}\mathfrak{d}^{\leq 2}\Ga_b  + r^{-1}\mathfrak{d}^{\leq 3}(\Gamma_b \c\widecheck{R}_b)-r^2\mathfrak{d}^{\leq 2}(\underline{A}\c \ov{A}).
\eea

\noindent{\bf Step 6.} We rewrite the wave equation \eqref{eq:proofcontrolPc:waveeqwidecheckq2nabwidehatR2P:renormalizedwithf} as follows 
\beaa
\bsplit
\square_\g(f\widecheck{q^2\widehat{R}^2P}) =& h,\\
h:=& O(m)\nab_{\Rhat}^{\leq 1}\nab\That\Pc +O(r^{-1})\nab_{\widehat{R}}\nabla\Rhat\widecheck{P} + O(r)\Rhat^{\leq 1}\Delta\Pc   +O(m^2r^{-2})\nab^2\Pc +O(1)\Rhat^{\leq 1}\dk\Pc \\
& +r^{-2}\mathfrak{d}^{\leq 2}\Ga_b  + r^{-1}\mathfrak{d}^{\leq 3}(\Gamma_b \c\widecheck{R}_b)-r^2\mathfrak{d}^{\leq 2}(\underline{A}\c \ov{A}).
\end{split}
\eeaa
We then apply the energy Morawetz estimate of Theorem \ref{Prop:scalarwavePsi-M8}  to the  above wave equation and obtain
\beaa
\BEF^{J-2}_{\de}[f\widecheck{q^2\widehat{R}^2P}] &\les& \ep_0^2+\NN_\de^{J-2}[f\widecheck{q^2\widehat{R}^2P},h].
\eeaa
Now, in view of the definition of $\NN_\de^s[\c, \c]$ in Section \ref{subsection:basicnormsforpsi}, we have 
\beaa
\NN_\de^{J-2}[f\widecheck{q^2\widehat{R}^2P},h] &\les& \left(\int_{\Mtrap}|\dk^{\leq J-2}h|^2\right)^{\frac{1}{2}}\left(\int_{\Mtrap}|\dk^{\leq J-1}f\widecheck{q^2\widehat{R}^2P}|^2\right)^{\frac{1}{2}}+\int_{\MM(\tau_1, \tau_2)}r^{1+\de}|\dk^{\leq J-2}h|^2\\
&&+\left(\int_{\MM(\tau_1, \tau_2)}r^{1+\de}|\dk^{\leq J-2}h|^2\right)^{\frac{1}{2}}\left(\int_{\Mntrap(\tau_1, \tau_2)}r^{\de-3}|\nab_4(r\dk^{\leq J-2}f\widecheck{q^2\widehat{R}^2P})|^2\right)^{\frac{1}{2}}\\
&\les& \int_{\MM(\tau_1, \tau_2)}r^{1+\de}|\dk^{\leq J-2}h|^2+\left(\int_{\MM(\tau_1, \tau_2)}r^{1+\de}|\dk^{\leq J-2}h|^2\right)^{\frac{1}{2}}\left(\B_\de^J[r^2\Pc]\right)^{\frac{1}{2}}.
\eeaa
Also, in view of the above definition of $h$, we have
\beaa
\int_{\MM}r^{1+\de}|\mathfrak{d}^{\leq J-2}h|^2 &\les& \B^{J-2}_{\de}[r^4\nab^2\widecheck{P}]+\B^{J-2}_{\de}[r^2(r\nab)(\nab_{\widehat{T}}, \nab_{\widehat{R}})\widecheck{P}]+\ep_J^2+\ep_0^2.
\eeaa
Plugging these estimates in the above control for $\BEF^{J-2}_{\de}[f\widecheck{q^2\widehat{R}^2P}]$ yields
\beaa
\BEF^{J-2}_{\de}[f\widecheck{q^2\nab_{\widehat{R}}^2P}] &\les& \B^{J-2}_{\de}[r^4\nab^2\widecheck{P}]+\B^{J-2}_{\de}[r^2(r\nab)(\nab_{\widehat{T}}, \nab_{\widehat{R}})\widecheck{P}]+\ep_J^2+\ep_0^2\\
&&+\left(\B^{J-2}_{\de}[r^4\nab^2\widecheck{P}]+\B^{J-2}_{\de}[r^2(r\nab)(\nab_{\widehat{T}}, \nab_{\widehat{R}})\widecheck{P}]+\ep_J^2+\ep_0^2\right)^{\frac{1}{2}}\sqrt{\B^J_{\de}[r^2\widecheck{P}]}
\eeaa
and hence, using also $f\widecheck{q^2\nab_{\widehat{R}}^2P}=fq^2\nab_{\widehat{R}}^2\Pc+r^{-1}\dk^{\leq 1}\Ga_b$ and Lemma \ref{lemma:basicrelationBEFdenormandnormsRkandSikforcurvatureandRiccicoeff:section789}, we deduce
\begin{align}\lab{eq:proofcontrolPc:BEFdeJminusesimtateforr2nabR2Pc}
\BEF^{J-2}_{\de}[r^2\nab_{\widehat{R}}^2\widecheck{P}] \les& \B^{J-2}_{\de}[r^4\nab^2\widecheck{P}]+\B^{J-2}_{\de}[r^2(r\nab)(\nab_{\widehat{T}}, \nab_{\widehat{R}})\widecheck{P}]+\mathfrak{S}_J\mathfrak{S}_{J+1}+\ep_J^2+\ep_0^2\nn\\
&+\left(\B^{J-2}_{\de}[r^4\nab^2\widecheck{P}]+\B^{J-2}_{\de}[r^2(r\nab)(\nab_{\widehat{T}}, \nab_{\widehat{R}})\widecheck{P}]+\ep_J^2+\ep_0^2\right)^{\frac{1}{2}}\sqrt{\B^J_{\de}[r^2\widecheck{P}]}.
\end{align}

In view of \eqref{eq:proofcontrolPc:consaeuenceHodgeestimatecontrolBdenormr4nab2Pc}, \eqref{eq:proofcontrolPc:controlBdenormr2rnabrThatandRhatPc},  \eqref{eq:proofcontrolPc:controlBdenormr2rnabThatnabRhatPc}, \eqref{eq:proofcontrolPc:controlBdenormr2rnabThat2minusnabRhat2Pc} and \eqref{eq:proofcontrolPc:BEFdeJminusesimtateforr2nabR2Pc}, we obtain
\begin{align}\lab{eq:proofcontrolPc:controlBdenormr2ofThatRhatnabsquareofPc}
\B^{J-2}_\de[r^2(r\nab, \nab_{\widehat{T}}, \nab_{\widehat{R}})^2\widecheck{P}] \les&  \mathfrak{S}_{J+1}\mathfrak{S}_J+\mathfrak{R}_{J+1}\mathfrak{R}_J+\ep_J^2+\ep_0^2+\sqrt{\mathfrak{R}_{J+1}\mathfrak{R}_J\EF^J_\de[r^2\widecheck{P}]}\\
\nn&+\Bigg(\mathfrak{S}_{J+1}\mathfrak{S}_J+\mathfrak{R}_{J+1}\mathfrak{R}_J+\ep_J^2+\ep_0^2+\sqrt{\mathfrak{R}_{J+1}\mathfrak{R}_J\EF^J_\de[r^2\widecheck{P}]}
\Bigg)^{\frac{1}{3}}\Big(\B^J_{\de}[r^2\widecheck{P}]\Big)^{\frac{2}{3}}.
\end{align}

\noindent{\bf Step 7.} Next, we recover estimates for $r\nab_4$ derivatives in $r\geq r_0$. To this end, we first apply Lemma \ref{Le:squareq^2Psi} to the wave equation \eqref{eq:P-WaveEq-M8} for $P$ which yields 
\beaa
\nn\square_\g(q^2P) &=&  Wq^2P  + r^{-1}\mathfrak{d}^{\leq 1}(\Gamma_b \c\widecheck{R}_b)-r^2\underline{A}\c \ov{A}, 
\eeaa
where $W=O(mr^{-3})+r^{-1}\Ga_g$. Then, we rely on the following commutator 
\beaa
\begin{split}
[r \nab_4, \square_\g]\psi &=-\nab_4\nab_4\psi -r\left(\frac 1 2 \trch -2\om\right) \square_\g\psi -r\left(\frac 1 2 \trch +2\om\right) \lap\psi  \\
&+O(r^{-2}) \dk^{\leq 1} \psi+O(r^{-3}) \dk^{\leq 2} \psi+r\Ddot_3 \big( \xi \c \Ddot_a \psi \big)+ r^{-1}\dk \big( \Ga_g \c \dk \psi),
\end{split}
\eeaa
see (4.7.14) in \cite{GKS22}. Together with the above wave equation for $q^2P$ and the assumption \eqref{eq:specialidentityforthegloablframeofMMinpartIII}, we infer
\beaa
\square_\g(r\nab_4(q^2P)) &=& \nab_4\nab_4(q^2P)  +r\left(\frac 1 2 \trch +2\om\right) \lap(q^2P) +O(r^{-2}) \dk^{\leq 1}(q^2P)+O(r^{-3}) \dk^{\leq 2}(q^2P)\\
&&+ r^{-1}\dk^{\leq 1}\big( \Ga_g \c \dk (q^2P))+ r\nab_4(Wq^2P)+r^{-1}\mathfrak{d}^{\leq 2}(\Gamma_b \c\widecheck{R}_b)-r^2\dk^{\leq 1}(\underline{A}\c \ov{A}),
\eeaa
and linearizing, we infer
\bea\lab{eq:waveeqrwidecheckq2Pforrpweightedsetimates}
\square_\g(r\widecheck{\nab_4(q^2P)}) &=& \nab_4\widecheck{\nab_4(q^2P)}  +r\left(\frac 1 2 \trch +2\om\right) \widecheck{\lap(q^2P)} +O(r^{-2}) \dk^{\leq 1}(q^2\Pc)+O(r^{-3}) \dk^{\leq 2}(q^2\Pc)\nn\\
&&+ r^{-1}\dk^{\leq 1}\big( \Ga_g \c \dk (q^2P))+r^{-1}\mathfrak{d}^{\leq 2}(\Gamma_b \c\widecheck{R}_b)-r^2\dk^{\leq 1}(\underline{A}\c \ov{A})+r^{-2}\dk^{\leq 1}\Ga_b.
\eea
Using a smooth cut-off in $r$ supported in $r\geq r_0/2$, and equal to $1$ for $r\geq r_0$, and running a standard $r^p$-weighted estimate with $p=\de$, we easily obtain, noticing that the first term on the RHS of \eqref{eq:waveeqrwidecheckq2Pforrpweightedsetimates} has the good sign, 
\beaa
\dot{\BEF}_{\de; r\geq r_0}^{J-1}[r\widecheck{\nab_4(q^2P)}] &\les& \ep_0^2+\ep_J^2+r_0\B_{\de; r_0/2\leq r\leq r_0}^{J-1}[r\widecheck{\nab_4(q^2P)}] +(\ep_0+\ep_J)\sqrt{\dot{\B}_{\de; r\geq r_0}^{J-1}[r\widecheck{\nab_4(q^2P)}]}\\
&&+\sqrt{\dot{\B}_{\de; r\geq r_0}^{J-2}[\widecheck{\De(q^2P)}]}\sqrt{\dot{\B}_{\de; r\geq r_0}^{J-1}[r\widecheck{\nab_4(q^2P)}]}+r_0^{-1}\dot{\B}_{\de; r\geq r_0}^J[r^2\Pc],
\eeaa
where we denote by $\dot{\B}_{p; r\geq r_0}$ and $\dot{\BEF}_{p; r\geq r_0}$ the pure $r^p$-part of $\B_{p; r\geq r_0}$ and $\BEF_{p; r\geq r_0}$ (i.e., without the contribution respectively from $\M_{\de, r\geq r_0}$ or $\EMF_{\de, r\geq r_0}$). Together with Lemma \ref{lemma:basicrelationBEFdenormandnormsRkandSikforcurvatureandRiccicoeff:section789}, and for $r_0$ large enough, we infer
\bea\lab{eq:proofcontrolPc:dotBEFrgeqr1forr3nab4Pc}
\dot{\BEF}_{\de; r\geq r_0}^{J-1}[r^3\nab_4\Pc] \les \ep_0^2+\ep_J^2+\Sk_J\Sk_{J+1}+r_0\B_{\de; r_0/2\leq r\leq r_0}^{J-1}[r^3\nab_4\Pc] +\B_{\de; r\geq r_0}^{J-2}[r^2(r\nab, \That)^2\Pc].
\eea

In addition to the control for $\dot{\BEF}_{\de; r\geq r_0}^{J-1}[r^3\nab_4\Pc]$ in \eqref{eq:proofcontrolPc:dotBEFrgeqr1forr3nab4Pc}, we also need to recover a control for $\T$ derivatives of $r^3\nab_4\Pc$. To this end, we rely on the following commutator in Lemma \ref{lemma:commutationof|q|^2squregwithnab3andrnab4forlocalexistenceregiontau*minus3tau*}
\beaa
\frac{1}{|q|^2}[r\nab_4, |q|^2\square_\g]P &=& -\frac{4}{r}\T(r\nab_4P)+O(r^{-2})\dk(r\nab_4P) +\dk^{\leq 1}(\Ga_g\dk P)+O(r^{-2})\dk P,
\eeaa
which together with the wave equation \eqref{eq:P-WaveEq-M8} for $P$ implies 
\beaa
\square_\g(r\nab_4P) &=& \frac{1}{|q|^2}r\nab_4(|q|^2\square_\g P) -\frac{1}{|q|^2}[r\nab_4, |q|^2\square_\g]P\\
&=& \frac{4}{r}\T(r\nab_4P) +\tr X\nab_3(r\nab_4P) +\ov{\tr\underline{X}}(r\nab_4P)  -\ov{H}\c\DD(r\nab_4P)  - \Hb\c\ov{\DD}(r\nab_4P) + Vr\nab_4P  \\
&& + r^{-3}\mathfrak{d}^{\leq 2}(\Gamma_b \c\widecheck{R}_b)-\dk^{\leq 1}(\underline{A}\c \ov{A})+O(r^{-2})\dk(r\nab_4P) +\dk^{\leq 1}(\Ga_g\dk P)+O(r^{-2})\dk P.
\eeaa
We then linearize this wave equation which yields
\begin{align*}
\square_\g(r\widecheck{\nab_4P}) =& \frac{4}{r}\T(r\widecheck{\nab_4P}) +\tr X\nab_3(r\widecheck{\nab_4P}) +\ov{\tr\underline{X}}(r\widecheck{\nab_4P})  -\ov{H}\c\DD(r\widecheck{\nab_4P})  - \Hb\c\ov{\DD}(r\widecheck{\nab_4P}) + Vr\widecheck{\nab_4P}  \\
& + r^{-3}\mathfrak{d}^{\leq 2}(\Gamma_b \c\widecheck{R}_b)-\dk^{\leq 1}(\underline{A}\c \ov{A})+O(r^{-2})\dk(r\widecheck{\nab_4P}) +\dk^{\leq 1}(\Ga_g\dk P)+O(r^{-2})\dk\Pc+r^{-4}\dk^{\leq 1}\Ga_b.
\end{align*}
Together with Lemma \ref{Le:squareq^2Psi}, we infer
\bea\lab{eq:waveequationforq2rwidecheckna4P:usefulenergyestimates}
\square_\g(q^2r\widecheck{\nab_4P}) &=& \frac{4}{r}\T(q^2r\widecheck{\nab_4P}) +Wq^2r\widecheck{\nab_4P} +O(r^{-2})\dk^{\leq 1}(q^2r\widecheck{\nab_4P})+O(r^{-2})q^2\dk\Pc\nn\\
&&+r^{-2}\dk^{\leq 1}\Ga_b +r^2\dk^{\leq 1}(\Ga_g\dk P) + r^{-1}\mathfrak{d}^{\leq 2}(\Gamma_b \c\widecheck{R}_b)-r^2\dk^{\leq 1}(\underline{A}\c \ov{A}),
\eea
where $W=O(mr^{-3})$. Using a smooth cut-off in $r$ supported in $r\geq r_0/2$, and equal to $1$ for $r\geq r_0$, and running an energy estimates with the vectorfield $\T$ which is timelike in $r\geq r_0/2$, we easily obtain, noticing that the first term on the RHS of \eqref{eq:waveequationforq2rwidecheckna4P:usefulenergyestimates} has the good sign, 
\beaa
&&\EF_{r\geq r_0}^{J-1}[q^2r\widecheck{\nab_4P}]+\int_{\MM_{r\geq r_0}}\frac{|\T(q^2r\dk^{J-1}\widecheck{\nab_4P}|^2}{r}\\ 
&\les&  \ep_0^2+\ep_J^2+r_0\B_{\de; r_0/2\leq r\leq r_0}^{J-1}[r\widecheck{\nab_4(q^2P)}] \\
&&+\left(\int_{\MM}\frac{|\T(q^2r\dk^{J-1}\widecheck{\nab_4P}|^2}{r}\right)^{\frac{1}{2}}\left(\ep_0^2+\ep_J^2+\dot{\B}_{\de; r\geq r_0}^{J-1}[r\widecheck{\nab_4(q^2P)}]+\B_{\de; r\geq r_0}^{J-1}[r^2\Pc]\right)^{\frac{1}{2}}
\eeaa
and hence, using in particular Lemma \ref{lemma:basicrelationBEFdenormandnormsRkandSikforcurvatureandRiccicoeff:section789}, 
\beaa
\EF_{r\geq r_0}^{J-1}[r^3\nab_4\Pc]+\int_{\MM_{r\geq r_0}}\frac{|r^3\T\dk^{J-1}\nab_4\Pc|^2}{r} &\les&  \ep_0^2+\ep_J^2+\Sk_J\Sk_{J+1}+r_0\B_{\de; r_0/2\leq r\leq r_0}^{J-1}[r\widecheck{\nab_4(q^2P)}] \\
&&+\dot{\B}_{\de; r\geq r_0}^{J-1}[r^3\nab_4\Pc].
\eeaa
Together with \eqref{eq:proofcontrolPc:dotBEFrgeqr1forr3nab4Pc}, we infer
\bea\lab{eq:proofcontrolPc:controlBEFdenormr2ofrnab4Pc:largerregion}
\BEF_{\de; r\geq r_0}^{J-1}[r^3\nab_4\Pc] \les \ep_0^2+\ep_J^2+\Sk_J\Sk_{J+1}+r_0\B_{\de; r_0/2\leq r\leq r_0}^{J-1}[r^3\nab_4\Pc] +\B_{\de; r\geq r_0}^{J-2}[r^2(r\nab, \That)^2\Pc].
\eea
Combining with \eqref{eq:proofcontrolPc:controlBdenormr2ofThatRhatnabsquareofPc}, this yields 
\begin{align}\lab{eq:proofcontrolPc:controlBdenormr2ofThatrnab4nabofPc}
\B^{J-1}_\de[r^2(r\nab, \nab_{\widehat{T}}, r\nab_4)\widecheck{P}] \les&  r_0\Big(\mathfrak{S}_{J+1}\mathfrak{S}_J+\mathfrak{R}_{J+1}\mathfrak{R}_J+\ep_J^2+\ep_0^2+\sqrt{\mathfrak{R}_{J+1}\mathfrak{R}_J\EF^J_\de[r^2\widecheck{P}]}\Big)\\
\nn+&r_0\Bigg(\mathfrak{S}_{J+1}\mathfrak{S}_J+\mathfrak{R}_{J+1}\mathfrak{R}_J+\ep_J^2+\ep_0^2+\sqrt{\mathfrak{R}_{J+1}\mathfrak{R}_J\EF^J_\de[r^2\widecheck{P}]}
\Bigg)^{\frac{1}{3}}\Big(\B^J_{\de}[r^2\widecheck{P}]\Big)^{\frac{2}{3}}.
\end{align}

Finally, we recover $\nab_3$ derivatives in the redshift region. To this end, we proceed as in Step 9 of Section 14.2.1 in \cite{GKS22}, noticing that the smallness of $a$ was not used in this particular step. We obtain 
\bea\lab{eq:proofcontrolPc:controlBdenormr2ofnab3ofPcinredshiftregion}
\BEF^{J-1}_\de[r^2\chi_{red}\nab_3\widecheck{P}]\lesssim \de_{red}^{-3}\B^{J-1}_\de[r^2(\nab_{\widehat{T}}, r\nab, r\nab_4)^2\widecheck{P}]+\ep_J^2+\ep_0^2.
\eea
Together with \eqref{eq:proofcontrolPc:controlBdenormr2ofThatrnab4nabofPc}, we obtain\footnote{Recall that $\les$ is allowed to depend implicitly on $\dred$ by convention, see Section \ref{sec:smallnesconstants}.}
\begin{align*}
\B^J_\de[r^2\Pc] \les&  r_0\Big(\mathfrak{S}_{J+1}\mathfrak{S}_J+\mathfrak{R}_{J+1}\mathfrak{R}_J+\ep_J^2+\ep_0^2+\sqrt{\mathfrak{R}_{J+1}\mathfrak{R}_J\EF^J_\de[r^2\widecheck{P}]}\Big)\\
\nn&+r_0\Bigg(\mathfrak{S}_{J+1}\mathfrak{S}_J+\mathfrak{R}_{J+1}\mathfrak{R}_J+\ep_J^2+\ep_0^2+\sqrt{\mathfrak{R}_{J+1}\mathfrak{R}_J\EF^J_\de[r^2\widecheck{P}]}
\Bigg)^{\frac{1}{3}}\Big(\B^J_{\de}[r^2\widecheck{P}]\Big)^{\frac{2}{3}}.
\end{align*}
and hence
\begin{align}\lab{eq:proofcontrolPc:controlBdenormr2ofPcforallweightedderivatives}
\B^J_\de[r^2\Pc] \les&  r_0^2\Big(\mathfrak{S}_{J+1}\mathfrak{S}_J+\mathfrak{R}_{J+1}\mathfrak{R}_J+\ep_J^2+\ep_0^2+\sqrt{\mathfrak{R}_{J+1}\mathfrak{R}_J\EF^J_\de[r^2\widecheck{P}]}\Big).
\end{align}

\noindent{\bf Step 8.} In view of \eqref{eq:proofcontrolPc:controlBdenormr2ofPcforallweightedderivatives}, we still need to recover energy and flux norms. To this end, recall from Lemma 14.2.3 in \cite{GKS22} that $\psi= q^2(\T P, \Z P)$ satisfies a wave equation of the type
\bea\lab{eq:waveequationforq2TPandZPfromGKS22}
\square_\g(q^2(\T P, \Z P)) = Wq^2(\T P, \Z P) +   r^{-2} \mathfrak{d}^{\leq 1}\Ga_b + r^{-1}\mathfrak{d}^{\leq 2}(\Ga_b \c\widecheck{R}_b)-r^2\mathfrak{d}^{\leq 1}(\underline{A}\c \ov{A}), \,\,\,\, W=O( m r^{-3}).
\eea
We then multiply $q^2(\T P, \Z P)$ by a smooth cut-off function in $r$ which vanishes on $\Mtrap$ and equals to $1$ on $r\leq r_+(1+\dred)$ and $r\geq 11m$ and apply Theorem \ref{Prop:scalarwavePsi-M8} to the resulting wave equation, using also the fact that $\T P=\T\Pc+r^{-3}\Ga_b$ and Lemma \ref{lemma:basicrelationBEFdenormandnormsRkandSikforcurvatureandRiccicoeff:section789}. We obtain 
\bea\lab{eq:proofcontrolPc:controlBEFdenormr2ofTandZPc:awayfromtrapping}
\BEF^{J-1}_{\de, r\leq r_+(1+\dred)}[r^2(\T\Pc, \Z\Pc)]+ \BEF^{J-1}_{\de, r\geq 11m}[r^2(\T\Pc, \Z\Pc)] \les \mathfrak{S}_{J+1}\mathfrak{S}_J+\mathfrak{R}_{J+1}\mathfrak{R}_J+\ep_J^2+\ep_0^2.
\eea
Together with \eqref{eq:proofcontrolPc:controlBEFdenormr4DDdhotDDdPc} and the non-integrable Hodge estimates in Lemma  \ref{lemma:improvednonintegrableHoadgeestimatefortheenergyonSigmatauusingcontroTandZderivatives}, we infer
\bea\lab{eq:proofcontrolPc:controlBEFdenormr4nabsquaredPc:awayfromtrapping:0000}
\BEF^{J-2}_{\de, r\leq r_+(1+\dred)}[r^4\nab^2\widecheck{P}] +\BEF^{J-2}_{\de, r\geq 11m}[r^4\nab^2\widecheck{P}] &\les& \mathfrak{S}_{J+1}\mathfrak{S}_J+\mathfrak{R}_{J+1}\mathfrak{R}_J+\ep_J^2+\ep_0^2.
\eea

Combining \eqref{eq:proofcontrolPc:BEFdeJminusesimtateforr2nabR2Pc}, \eqref{eq:proofcontrolPc:controlBEFdenormr2ofrnab4Pc:largerregion}, \eqref{eq:proofcontrolPc:controlBdenormr2ofnab3ofPcinredshiftregion}, \eqref{eq:proofcontrolPc:controlBdenormr2ofPcforallweightedderivatives}, \eqref{eq:proofcontrolPc:controlBEFdenormr2ofTandZPc:awayfromtrapping} and \eqref{eq:proofcontrolPc:controlBEFdenormr4nabsquaredPc:awayfromtrapping:0000}, we obtain 
\bea\lab{eq:proofcontrolPc:controlBEFdenormr4nabsquaredPc:awayfromtrapping}
\B^J_\de[r^2\Pc] +\EF^J_{\de, r\leq r_+(1+\dred)}[r^2\widecheck{P}] +\EF^J_{\de, r\geq 11m}[r^2\widecheck{P}] \les r_0^2\Big(\mathfrak{S}_{J+1}\mathfrak{S}_J+\mathfrak{R}_{J+1}\mathfrak{R}_J+\ep_J^2+\ep_0^2\Big).
\eea

\noindent{\bf Step 9.} In view of \eqref{eq:proofcontrolPc:controlBEFdenormr4nabsquaredPc:awayfromtrapping}, we still need to recover the energy in the region $r_+(1+\dred)\leq r\leq 11m$. To this end, we consider an integer $n$ such that $1\leq n<n+1\leq\tau_*-3$ and combine \eqref{eq:proofcontrolPc:controlBEFdenormr4DDdhotDDdPc}, Lemma \ref{lemma:basicrelationBEFdenormandnormsRkandSikforcurvatureandRiccicoeff:section789} and Lemma \ref{lemma:howtogofromDDhotDDtonab2} which yields
\bea\lab{eq:proofcontrolPc:consaeuenceHodgeestimatecontrolEnormintegratednnplus1rplusto11mr4nab2Pc}
\int_{n}^{n+1}\E^{J-2}_{r_+(1+\dred)\leq r\leq 11m}[r^2\nab^2\Pc](\tau)d\tau \les \sqrt{\Rk_j\Rk_{J+1}\sup_{\tau}\E^J[\Pc](\tau)}+\mathfrak{S}_{J+1}\mathfrak{S}_J+\mathfrak{R}_{J+1}\mathfrak{R}_J+\ep_J^2+\ep_0^2.
\eea
Also, arguing as in Step 3, we have
\beaa
&&\int_{n}^{n+1}\E^{J-2}_{r_+(1+\dred)\leq r\leq 11m}[r^2(\nab_{\widehat{T}}^2\widecheck{P} - \nab_{\widehat{R}}^2\widecheck{P})](\tau)d\tau\\ 
&\lesssim& \int_{n}^{n+1}\Big(\E^{J-2}_{r_+(1+\dred)\leq r\leq 11m}[r^2\De\widecheck{P}](\tau)d\tau+\E^{J-1}_{r_+(1+\dred)\leq r\leq 11m}[r\widecheck{P}](\tau)+\E^{J-2}_{r_+(1+\dred)\leq r\leq 11m}[\Rc_g](\tau)\\
&&+\E^{J-2}_{r_+(1+\dred)\leq r\leq 11m}[r^{-2}  \dk^{\le 1} \Ga_b](\tau)\Big)d\tau+\ep_0^2
\eeaa
and hence, using \eqref{eq:proofcontrolPc:consaeuenceHodgeestimatecontrolEnormintegratednnplus1rplusto11mr4nab2Pc} and Lemma \ref{lemma:basicrelationBEFdenormandnormsRkandSikforcurvatureandRiccicoeff:section789}, we deduce 
\bea\lab{eq:proofcontrolPc:controlBdenormr2rnabThat2minusnabRhat2Pc:energy}
\int_{n}^{n+1}\E^{J-2}_{r_+(1+\dred)\leq r\leq 11m}[r^2(\nab_{\widehat{T}}^2\widecheck{P} - \nab_{\widehat{R}}^2\widecheck{P})](\tau)d\tau &\lesssim& \mathfrak{S}_{J+1}\mathfrak{S}_J+\mathfrak{R}_{J+1}\mathfrak{R}_J+\ep_J^2+\ep_0^2\nn\\
&&+\sqrt{\Rk_j\Rk_{J+1}\sup_{\tau}\E^J[\Pc](\tau)}.
\eea

Then, \eqref{eq:proofcontrolPc:BEFdeJminusesimtateforr2nabR2Pc}, \eqref{eq:proofcontrolPc:controlBEFdenormr4nabsquaredPc:awayfromtrapping}, \eqref{eq:proofcontrolPc:consaeuenceHodgeestimatecontrolEnormintegratednnplus1rplusto11mr4nab2Pc}, \eqref{eq:proofcontrolPc:controlBdenormr2rnabThat2minusnabRhat2Pc:energy}, Lemma \ref{lemma:howtocontrolcrosstermsnabnabTandnabnabRpsi} and Lemma \ref{lemma:basicrelationBEFdenormandnormsRkandSikforcurvatureandRiccicoeff:section789}, we infer, for any integer $n$ such that $1\leq n<n+1\leq\tau_*-3$,
\beaa
\int_{n}^{n+1}\E^{J}[r^2\Pc](\tau)d\tau &\lesssim& \mathfrak{S}_{J+1}\mathfrak{S}_J+\mathfrak{R}_{J+1}\mathfrak{R}_J+\ep_J^2+\ep_0^2+\sqrt{\Rk_j\Rk_{J+1}\sup_{\tau}\E^J[\Pc](\tau)}.
\eeaa
Thus, for any integer $n$ such that $1\leq n<n+1\leq\tau_*-3$, we deduce the existence of $\tau^{(n)}\in[n,n+1]$ such that 
\bea\lab{eq:proofcontrolPc:controlEnormrallJweightedderivativesPcattau=taun}
\E^{J}[r^2\Pc](\tau^{(n)}) &\lesssim& \mathfrak{S}_{J+1}\mathfrak{S}_J+\mathfrak{R}_{J+1}\mathfrak{R}_J+\ep_J^2+\ep_0^2+\sqrt{\Rk_j\Rk_{J+1}\sup_{\tau}\E^J[\Pc](\tau)}.
\eea
Next, using local energy estimates for the wave equation \eqref{eq:waveequationforq2TPandZPfromGKS22} on $\tau_{(n)}\leq\tau\leq\min(\tau_{(n)}+1, \tau_*)$, and relying on \eqref{eq:proofcontrolPc:controlEnormrallJweightedderivativesPcattau=taun}, we obtain 
\bea\lab{eq:proofcontrolPc:controlEnormJ-1weightedderivativesTandZPc}
\sup_{\tau}\E^{J-1}[r^2(\T\Pc, \Z\Pc)](\tau) &\lesssim& \mathfrak{S}_{J+1}\mathfrak{S}_J+\mathfrak{R}_{J+1}\mathfrak{R}_J+\ep_J^2+\ep_0^2+\sqrt{\Rk_j\Rk_{J+1}\sup_{\tau}\E^J[\Pc](\tau)}.
\eea
We then, deduce from \eqref{eq:proofcontrolPc:controlBEFdenormr4DDdhotDDdPc}, \eqref{eq:proofcontrolPc:controlEnormJ-1weightedderivativesTandZPc} and from the non-integrable Hodge estimates of Lemma \ref{lemma:improvednonintegrableHoadgeestimatefortheenergyonSigmatauusingcontroTandZderivatives}
\beaa
\sup_{\tau}\E^{J-2}[r^4\nab^2\widecheck{P}](\tau) &\les& \mathfrak{S}_{J+1}\mathfrak{S}_J+\mathfrak{R}_{J+1}\mathfrak{R}_J+\ep_J^2+\ep_0^2+\sqrt{\Rk_j\Rk_{J+1}\sup_{\tau}\E^J[\Pc](\tau)},
\eeaa
which together with \eqref{eq:proofcontrolPc:BEFdeJminusesimtateforr2nabR2Pc} and \eqref{eq:proofcontrolPc:controlEnormJ-1weightedderivativesTandZPc} implies
\beaa
\sup_{\tau}\E^J[\widecheck{P}](\tau) &\les& \mathfrak{S}_{J+1}\mathfrak{S}_J+\mathfrak{R}_{J+1}\mathfrak{R}_J+\ep_J^2+\ep_0^2+\sqrt{\Rk_j\Rk_{J+1}\sup_{\tau}\E^J[\Pc](\tau)},
\eeaa
and hence
\beaa
\sup_{\tau}\E^J[\widecheck{P}](\tau) &\les& \mathfrak{S}_{J+1}\mathfrak{S}_J+\mathfrak{R}_{J+1}\mathfrak{R}_J+\ep_J^2+\ep_0^2.
\eeaa
Together with \eqref{eq:proofcontrolPc:controlBEFdenormr4nabsquaredPc:awayfromtrapping}, we deduce
\beaa
\BEF^J_\de[r^2\Pc] &\les& r_0^2\Big(\mathfrak{S}_{J+1}\mathfrak{S}_J+\mathfrak{R}_{J+1}\mathfrak{R}_J+\ep_J^2+\ep_0^2\Big).
\eeaa
Finally, recalling that we have dropped for convenience $(1, \tau_*-3)$ on all LHS in Steps 2--9, see the end of Step 1, we have in fact obtained so far
\bea\lab{eq:proofcontrolPc:controlBEFdenormrJderivativesPc:onlyontauleqtaustarminus3}
\BEF^J_\de[r^2\Pc](1, \tau_*-3) &\les& r_0^2\Big(\mathfrak{S}_{J+1}\mathfrak{S}_J+\mathfrak{R}_{J+1}\mathfrak{R}_J+\ep_J^2+\ep_0^2\Big),
\eea
which corresponds to the proof of \eqref{eq:conclusionchapter14forlargea} on $\MM(1, \tau_*-3)$.

\noindent{\bf Step 10.} We now complete the proof of Theorem \ref{theorem:Morawetz-EnergyPc} by extending \eqref{eq:proofcontrolPc:controlBEFdenormrJderivativesPc:onlyontauleqtaustarminus3} to $\MM$. To this end, recall that \eqref{eq:proofcontrolPc:controlEnormJ-1weightedderivativesTandZPc} follows form applying  local energy estimates for the wave equation \eqref{eq:waveequationforq2TPandZPfromGKS22} on $\tau_{(n)}\leq\tau\leq\min(\tau_{(n)}+1, \tau_*)$ and relying on \eqref{eq:proofcontrolPc:controlEnormrallJweightedderivativesPcattau=taun}, so that \eqref{eq:proofcontrolPc:controlEnormJ-1weightedderivativesTandZPc} holds in fact on $\MM$. In particular, we have
\bea\lab{eq:proofcontrolPc:controlEnormJ-1weightedderivativesTandZPc:evenbetyondtau*minus3}
\sup_{\tau\in[\tau_*-3, \tau_*]}\EF^{J-1}[r^2(\T\Pc, \Z\Pc)](\tau) \lesssim r_0^2\Big(\mathfrak{S}_{J+1}\mathfrak{S}_J+\mathfrak{R}_{J+1}\mathfrak{R}_J+\ep_J^2+\ep_0^2\Big).
\eea

Next, we commute the wave equation \eqref{eq:P-WaveEq-M8} for $P$ with $\nab_3$ and rely on Lemma \ref{lemma:commutationof|q|^2squregwithnab3andrnab4forlocalexistenceregiontau*minus3tau*} which implies
\beaa
\nn\square_\g(\nab_3P) &=&  \tr X\nab_3(\nab_3P) +\ov{\tr\underline{X}}\nab_4(\nab_3P)  -\ov{H}\c\DD(\nab_3P)  - \Hb\c\ov{\DD}(\nab_3P) + V\nab_3P  \\
&& + r^{-3}\mathfrak{d}^{\leq 2}(\Gamma_b \c\widecheck{R}_b)-\dk(\underline{A}\c \ov{A})+ O(r^{-2})\dk\nab_3P +\dk^{\leq 1}(\Ga_g\dk P)+O(r^{-2})\dk P.
\eeaa
We linearize this wave equation which yields
\beaa
\nn\square_\g(\widecheck{\nab_3P}) &=&  \tr X\nab_3(\widecheck{\nab_3P}) +\ov{\tr\underline{X}}\nab_4(\widecheck{\nab_3P})  -\ov{H}\c\DD(\widecheck{\nab_3P})  - \Hb\c\ov{\DD}(\widecheck{\nab_3P}) + V\widecheck{\nab_3P}  \\
&& + r^{-3}\mathfrak{d}^{\leq 2}(\Gamma_b \c\widecheck{R}_b)-\dk(\underline{A}\c \ov{A})+ O(r^{-2})\dk\widecheck{\nab_3P} +O(r^{-2})\dk\Pc+r^{-4}\dk^{\leq 1}\Ga_b.
\eeaa
Then, applying Lemma \ref{Le:squareq^2Psi}, we infer
\beaa
\nn\square_\g(q^2\widecheck{\nab_3P}) &=&  Wq^2\widecheck{\nab_3P} + O(1)\dk\widecheck{\nab_3P} +O(1)\dk\Pc+r^{-2}\dk^{\leq 1}\Ga_b  + r^{-1}\mathfrak{d}^{\leq 2}(\Gamma_b \c\widecheck{R}_b)-r^2\dk(\underline{A}\c \ov{A}),\\
\nn\square_\g(q^2r\widecheck{\nab_4P}) &=&  Wq^2r\widecheck{\nab_4P} +\frac{4}{r}\T(q^2r\widecheck{\nab_4P}) + O(1)\dk(r\widecheck{\nab_4P}) +O(1)\dk P+r^{-2}\dk^{\leq 1}\Ga_b \\
&& + r^{-1}\mathfrak{d}^{\leq 2}(\Gamma_b \c\widecheck{R}_b)-r^2\dk(\underline{A}\c \ov{A}),
\eeaa
for a scalar function $W$ satisfying $W=O(mr^{-3})$, where the second wave equation has been derived in \eqref{eq:waveequationforq2rwidecheckna4P:usefulenergyestimates}. Then:
\begin{itemize}
\item commuting with $\dk^k$ for $k\leq J-1$, which preserves the structure of these wave equations, 

\item relying on local energy estimates for the corresponding wave equations using a vectorfield $\widetilde{T}$ which is globally timelike on $\MM(\tau_*-3, \tau_*)$ and agrees with $\T$ for $r\geq 3m$, 

\item noticing that the term $r^{-1}\T(q^2r\widecheck{\nab_4P})$ on the RHS of the wave equation for $q^2r\widecheck{\nab_4P}$ has the good sign in $r\geq 3m$,

\item using \eqref{eq:auxilliarynormsforGa_b1} to control the contribution of the terms $r^{-2}\dk^{\leq 1}\Ga_b$ on the RHS of the wave equations, 

\item and using \eqref{eq:proofcontrolPc:controlBEFdenormrJderivativesPc:onlyontauleqtaustarminus3} to control the initial data on $\Si(\tau_*-3)$, 
\end{itemize}
we easily obtain, after applying Gronwall to the resulting inequalities, the following estimates 
\beaa
\EF^{J-1}[r^2(\widecheck{\nab_3P}, r\widecheck{\nab_4P})](\tau_*-3, \tau_*) &\lesssim& r_0^2\Big(\mathfrak{S}_{J+1}\mathfrak{S}_J+\mathfrak{R}_{J+1}\mathfrak{R}_J+\ep_J^2+\ep_0^2\Big).
\eeaa
Since $r^2\widecheck{\nab_3P}=r^2\nab_3\Pc+r^{-1}\Ga_b$ and $rr^3\widecheck{\nab_4P}=r^3\nab_4\Pc+r^{-1}\Ga_g$, we infer, in view of Lemma \ref{lemma:basicrelationBEFdenormandnormsRkandSikforcurvatureandRiccicoeff:section789}, 
\beaa
\EF^{J-1}[r^2(\nab_3\Pc, r\nab_4\Pc)](\tau_*-3, \tau_*) &\lesssim& r_0^2\Big(\mathfrak{S}_{J+1}\mathfrak{S}_J+\mathfrak{R}_{J+1}\mathfrak{R}_J+\ep_J^2+\ep_0^2\Big).
\eeaa
Together with \eqref{eq:proofcontrolPc:controlEnormJ-1weightedderivativesTandZPc:evenbetyondtau*minus3}, this yields 
\bea\lab{eq:conclusionchapter14forlargea:sofaronlyontauleqtaustarminus3:extendedbutonlyfornab3Pandrnab4P}
\EF^{J-1}[r^2(\T\Pc, \Z\Pc, \nab_3\Pc, r\nab_4\Pc)](\tau_*-3, \tau_*) &\lesssim& r_0^2\Big(\mathfrak{S}_{J+1}\mathfrak{S}_J+\mathfrak{R}_{J+1}\mathfrak{R}_J+\ep_J^2+\ep_0^2\Big).
\eea

Next, using in particular \eqref{eq:specialidentityforthegloablframeofMMinpartIII}, we have the following commutator, see Step 5 of Section 14.2.1 in \cite{GKS22},
\beaa
\square_\g (|q|^2\lap \psi) &=& \lap(|q|^2\square_\g\psi) + O(m^2r^{-3}) \dk^{\leq 2}\psi+ \dk^2 \big( \Ga_g \c \dk \psi).
\eeaa
Together with the wave equation \eqref{eq:P-WaveEq-M8} for $P$, we infer 
\beaa
\square_\g (|q|^2\lap P) &=& \lap\Big(|q|^2\Big(\tr X\nab_3P +\ov{\tr\underline{X}}\nab_4P  -\ov{H}\c\DD P  - \Hb\c\ov{\DD}P + V P\Big)\Big) + O(m^2r^{-3}) \dk^{\leq 2}P\\
&& +r^{-3}\mathfrak{d}^{\leq 3}(\Gamma_b \c\widecheck{R}_b)-\dk^{\leq 2}(\underline{A}\c \ov{A})+ \dk^2 \big( \Ga_g \c \dk P)\\
&=& \tr X\nab_3(|q|^2\lap P) +\ov{\tr\underline{X}}\nab_4(|q|^2\lap P)  -\ov{H}\c\DD(|q|^2\lap P)  - \Hb\c\ov{\DD}(|q|^2\lap P) + V|q|^2\lap P \\
&& + O(mr^{-2}) \dk^{\leq 2}P +r^{-3}\mathfrak{d}^{\leq 3}(\Gamma_b \c\widecheck{R}_b)-\dk^{\leq 2}(\underline{A}\c \ov{A})+ \dk^2 \big( \Ga_g \c \dk P).
\eeaa
Then, applying Lemma \ref{Le:squareq^2Psi}, we infer
\beaa
\square_\g (q^2|q|^2\lap P) &=& Wq^2|q|^2\lap P  + O(1) \dk^{\leq 2}P +r^{-1}\mathfrak{d}^{\leq 3}(\Gamma_b \c\widecheck{R}_b)-r^2\dk^{\leq 2}(\underline{A}\c \ov{A})+ r^2\dk^2 \big( \Ga_g \c \dk P),
\eeaa
where $W=O(mr^{-3})$. Then, linearizing and relying on Lemma \ref{lemma:basicrelationBEFdenormandnormsRkandSikforcurvatureandRiccicoeff:section789}, we deduce 
\beaa
\square_\g (\widecheck{q^2|q|^2\lap P}) &=& W\widecheck{q^2|q|^2\lap P}  + O(1) \dk^{\leq 2}\Pc +O(r^{-2} ) \dk^{\le 2} \Ga_b \\
&&+r^{-1}\mathfrak{d}^{\leq 3}(\Gamma_b \c\widecheck{R}_b)-r^2\dk^{\leq 2}(\underline{A}\c \ov{A})+ r^2\dk^2 \big( \Ga_g \c \dk\Pc),
\eeaa
Running a local energy estimate on $\tau_*-3\leq\tau\leq\tau_*$ as above, we infer 
\beaa
\EF^{J-2}[r^2\widecheck{q^2 \lap  P}](\tau_*-3, \tau_*) &\les& r_0^2\Big(\mathfrak{S}_{J+1}\mathfrak{S}_J+\mathfrak{R}_{J+1}\mathfrak{R}_J+\ep_J^2+\ep_0^2\Big).
\eeaa
In view of Lemma \ref{lemma:basicrelationBEFdenormandnormsRkandSikforcurvatureandRiccicoeff:section789}, we deduce
\beaa
\EF^{J-2}[r^4\lap\Pc](\tau_*-3, \tau_*) &\les& r_0^2\Big(\mathfrak{S}_{J+1}\mathfrak{S}_J+\mathfrak{R}_{J+1}\mathfrak{R}_J+\ep_J^2+\ep_0^2\Big),
\eeaa
which together with \eqref{eq:proofcontrolPc:controlEnormJ-1weightedderivativesTandZPc:evenbetyondtau*minus3} and the non-integrable Hodge estimates of Lemma \ref{lemma:improvednonintegrableHoadgeestimatefortheenergyonSigmatauusingcontroTandZderivatives} yields
\beaa
\EF^{J-2}[r^4\nab^2\Pc](\tau_*-3, \tau_*) &\les& r_0^2\Big(\mathfrak{S}_{J+1}\mathfrak{S}_J+\mathfrak{R}_{J+1}\mathfrak{R}_J+\ep_J^2+\ep_0^2\Big).
\eeaa
Combining with \eqref{eq:conclusionchapter14forlargea:sofaronlyontauleqtaustarminus3:extendedbutonlyfornab3Pandrnab4P}, we obtain
\bea\lab{eq:conclusionchapter14forlargea:sofaronlyontauleqtaustarminus3:extendedbutonlyforEFnormsomissesBEFde}
\EF^{J}[r^2\Pc](\tau_*-3, \tau_*) &\lesssim& r_0^2\Big(\mathfrak{S}_{J+1}\mathfrak{S}_J+\mathfrak{R}_{J+1}\mathfrak{R}_J+\ep_J^2+\ep_0^2\Big).
\eea

Next, using the wave equation \eqref{eq:waveequationforq2TPandZPfromGKS22} and the above wave equations for $\widecheck{\nab_3P}$ and $\widecheck{q^2|q|^2\lap P}$, and applying a Morawetz  estimate  in $r\geq r_0$, we easily infer, using also \eqref{eq:conclusionchapter14forlargea:sofaronlyontauleqtaustarminus3:extendedbutonlyforEFnormsomissesBEFde},  
\beaa
&&\M^{J-1}_{\de, \geq r_0}[r^2(\T\Pc, \Z\Pc, \nab_3\Pc)](\tau_*-3, \tau_*) +\M^{J-2}_{\de, \geq r_0}[r^4\De\Pc](\tau_*-3, \tau_*) \\  
&\les& \EF^{J}[r^2\Pc](\tau_*-3, \tau_*)+\mathfrak{S}_{J+1}\mathfrak{S}_J+\mathfrak{R}_{J+1}\mathfrak{R}_J+\ep_J^2+\ep_0^2\\
&\les& r_0^2\Big(\mathfrak{S}_{J+1}\mathfrak{S}_J+\mathfrak{R}_{J+1}\mathfrak{R}_J+\ep_J^2+\ep_0^2\Big),
\eeaa 
which together with the non-integrable Hodge estimates of Lemma \ref{lemma:improvednonintegrableHoadgeestimatefortheenergyonSigmatauusingcontroTandZderivatives} yields
\beaa
\M^{J-1}_{\de, \geq r_0}[r^2(\T\Pc, \Z\Pc, \nab_3\Pc)](\tau_*-3, \tau_*) +\M^{J-2}_{\de, \geq r_0}[r^4\nab^2\Pc](\tau_*-3, \tau_*) \les r_0^2\Big(\mathfrak{S}_{J+1}\mathfrak{S}_J+\mathfrak{R}_{J+1}\mathfrak{R}_J+\ep_J^2+\ep_0^2\Big).
\eeaa 
Also, the above energy estimate for the wave equation \eqref{eq:waveequationforq2rwidecheckna4P:usefulenergyestimates} satisfied by $q^2r\widecheck{\nab_4P}$ already yields 
\beaa
\M^{J-1}_{\de, \geq r_0}[r^3\nab_4\Pc)](\tau_*-3, \tau_*)  \les r_0^2\Big(\mathfrak{S}_{J+1}\mathfrak{S}_J+\mathfrak{R}_{J+1}\mathfrak{R}_J+\ep_J^2+\ep_0^2\Big),
\eeaa
thanks to the fact that the term $r^{-1}\T(q^2r\widecheck{\nab_4P})$ on the RHS of the wave equation for $q^2r\widecheck{\nab_4P}$ has the good sign in $r\geq 3m$. We thus deduce, together with \eqref{eq:conclusionchapter14forlargea:sofaronlyontauleqtaustarminus3:extendedbutonlyforEFnormsomissesBEFde}, 
\beaa
\EMF^{J}_\de[r^2\Pc](\tau_*-3, \tau_*) &\lesssim& r_0^2\Big(\mathfrak{S}_{J+1}\mathfrak{S}_J+\mathfrak{R}_{J+1}\mathfrak{R}_J+\ep_J^2+\ep_0^2\Big).
\eeaa
Noticing that on $\MM(\tau_*-3, \tau_*)$, we have
\beaa
\B^{J}_\de[r^2\Pc](\tau_*-3, \tau_*) \les \EMF^{J}_\de[r^2\Pc](\tau_*-3, \tau_*), 
\eeaa
we deduce 
\beaa
\B^{J}_\de[r^2\Pc](\tau_*-3, \tau_*)+\EF^{J}[r^2\Pc](\tau_*-3, \tau_*) &\lesssim& r_0^2\Big(\mathfrak{S}_{J+1}\mathfrak{S}_J+\mathfrak{R}_{J+1}\mathfrak{R}_J+\ep_J^2+\ep_0^2\Big).
\eeaa
Finally, to recover $\EF^{J}_\de[r^2\Pc](\tau_*-3, \tau_*)$, we run $r^p$ weighted estimates in $r\geq r_0$ for the wave equation \eqref{eq:waveequationforq2TPandZPfromGKS22}, for the above wave equations for $\widecheck{\nab_3P}$ and $\widecheck{q^2|q|^2\lap P}$, and for the wave equation for $r\widecheck{\nab_4(q^2P)}$ in \eqref{eq:waveeqrwidecheckq2Pforrpweightedsetimates} which easily yields
\begin{align*}
\EF_\de^{J}[r^2\Pc](\tau_*-3, \tau_*) \les& \B^{J}_\de[r^2\Pc](\tau_*-3, \tau_*)+\EF^{J}[r^2\Pc](\tau_*-3, \tau_*) + r_0^2\Big(\mathfrak{S}_{J+1}\mathfrak{S}_J+\mathfrak{R}_{J+1}\mathfrak{R}_J+\ep_J^2+\ep_0^2\Big)\\
\les& r_0^2\Big(\mathfrak{S}_{J+1}\mathfrak{S}_J+\mathfrak{R}_{J+1}\mathfrak{R}_J+\ep_J^2+\ep_0^2\Big)
\end{align*}
and hence, combining with the above, 
\beaa
\BEF^{J}_\de[r^2\Pc](\tau_*-3, \tau_*) &\lesssim& r_0^2\Big(\mathfrak{S}_{J+1}\mathfrak{S}_J+\mathfrak{R}_{J+1}\mathfrak{R}_J+\ep_J^2+\ep_0^2\Big).
\eeaa
Together with \eqref{eq:proofcontrolPc:controlBEFdenormrJderivativesPc:onlyontauleqtaustarminus3}, this implies 
\beaa
\BEF^J_\de[r^2\Pc](1, \tau_*) &\les& \mathfrak{S}_{J+1}\mathfrak{S}_J+\mathfrak{R}_{J+1}\mathfrak{R}_J+\ep_J^2+\ep_0^2,
\eeaa
as stated in \eqref{eq:conclusionchapter14forlargea}. This concludes the proof of Theorem \ref{theorem:Morawetz-EnergyPc}.

%%%%%%%%%%%%%%%%%%%%%%%%%%%%%%%%%%%%%%

\subsubsection{Non-integrable Hodge estimates for the control of $\Pc$}
\lab{sec:nonintegrableHodgeestimatesforcontrolPc}

%%%%%%%%%%%%%%%%%%%%%%%%%%%%%%%%%%%%%%

In this section, we derive non-integrable Hodge estimates used in the proof of Theorem \ref{theorem:Morawetz-EnergyPc}.
\begin{lemma}\lab{lemma:howtogofromDDhotDDtonab2}
Let $\psi$ be a scalar function. Then, we have
\beaa
\B_\de^{J-2}[r^2\nab^2\psi] \les \B_\de^{J-2}[r^2\DD\hot\DD\psi]+\sqrt{\BEF^{J-1}_\de[\psi]\BEF^J_\de[\psi]}+\ep\M_1^{J}[\psi].
\eeaa
Also, for any $[n, n+1]\subset[1,\tau_*]$, we have 
\beaa
&&\int_{n}^{n+1}\E^{J-2}_{r_+(1+\dred)\leq r\leq 11m}[r^2\nab^2\psi](\tau)d\tau\\
&\les& \sup_{\tau\in[n,n+1]}\Big(\E^{J-2}[r^2\DD\hot\DD\psi](\tau)+\sqrt{\E^{J-1}[\psi](\tau)\E^J[\psi](\tau)}+\ep\E^J[\psi](\tau)\Big).
\eeaa
\end{lemma}

\begin{proof}
As in (15.3.1) in \cite{GKS22}, we introduce the following complex non-integrable Hodge  operators
\bea
\lab{def:complexadjointperators.M8}
\DDs_2 &=&-\frac 1 2 \DD\hot, \qquad \DDd_2 =\frac 1  2 \DDov \c, \qquad \DDs_1= -\DD, \qquad \DDd_1=\frac 1 2 \DDov\c,
\eea
which act respectively on $\sk_1(\mathbb{C})$, $\sk_2(\mathbb{C})$, $\sk_0(\mathbb{C})$ and $\sk_1(\mathbb{C})$. Then, in view of Lemma 15.3.3 in \cite{GKS22}, these operators satisfy the following identities
\bea\lab{eq:identitiesrelatingDDdDDsorDDsDDdwithDeltap=012}
\bsplit
\DDd_1\DDs_1 &= -\lap_0     - \frac 1 2 i    (\atrch\nab_3+\atrchb \nab_4),\\
\DDs_1\, \DDd_1 &=  -\lap_1  + \frac 1 2 i   (\atrch\nab_3+\atrchb \nab_4)  + \Kh,\\
\DDd_2\DDs_2 &= - \lap_1 - \frac 1 2 i    (\atrch\nab_3+\atrchb \nab_4)  - \Kh,\\
\DDs_2 \DDd_2 &=-\lap_2   + \frac 1 2 i    (\atrch\nab_3+\atrchb \nab_4) + 2 \Kh,
\end{split}
\eea
with the scalar $\Kh$ given by the formula
\bea\label{eq:definition-K-in}
\Kh := - \frac 14  \trch \trchb-\frac 1 4 \atrch \atrchb+\frac 1 2 \chih \c \chibh-  \rho.  
\eea
We also recall Lemma 2.1.40 in \cite{GKS22} which allows to write horizontal divergences in terms of spacetime divergences, for $f\in\sk_1$, as follows 
\bea\lab{eq:rewritinghorizontaldivergenceasaspacetimedivergence}
\D^\a f_\a&=& \div(f) + (\eta +\etab) \c f.
\eea
Then, using \eqref{eq:rewritinghorizontaldivergenceasaspacetimedivergence}, integration by parts, $\nab(r)\in r\Ga_g$, and $\eta+\etab=O(r^{-2})+\Ga_b$, we have, for a scalar function $h$ and a real number $q$,
\beaa
\int_{\MM}r^q|\DD\hot\DD(h)|^2 &=&  \int_{\MM}r^q|\DDs_2\DDs_1(h)|^2\\
 &=& -\int_{\MM}r^q\DDs_1(h)\c\DDd_2\DDs_2\DDs_1(h)+O(1)\int_{\MM}r^{q-1}|\nab^2(h)||\nab(h)|\\
 &&+O(1)\int_{\pr\MM}r^{q-1}|\nab^2(h)||\nab(h)|,
\eeaa
which together with the third and the second identities in \eqref{eq:identitiesrelatingDDdDDsorDDsDDdwithDeltap=012} yields, using again \eqref{eq:rewritinghorizontaldivergenceasaspacetimedivergence} and integration by parts,
\beaa
\int_{\MM}r^q|\DD\hot\DD(h)|^2 &=& -\int_{\MM}r^q\DDs_1(h)\c\DDs_1\DDd_1\DDs_1(h)+O(1)\int_{\MM}r^{q-2}|\nab^2(h)||\dk^{\leq 1}h|\\
&&+O(1)\int_{\pr\MM}r^{q-1}|\nab^2(h)||\nab(h)|\\
&=& \int_{\MM}r^q|\DDd_1\DDs_1(h)|^2+O(1)\int_{\MM}r^{q-2}|\nab^2(h)||\dk^{\leq 1}h|+O(1)\int_{\pr\MM}r^{q-1}|\nab^2(h)||\nab(h)|.
\eeaa
Using the first identity in \eqref{eq:identitiesrelatingDDdDDsorDDsDDdwithDeltap=012}, we infer
\beaa
\int_{\MM}r^q|\DD\hot\DD(h)|^2 &=& \int_{\MM}r^q|\De_0(h)|^2+O(1)\int_{\MM}r^{q-2}|\nab^2(h)||\dk^{\leq 1}h|+O(1)\int_{\pr\MM}r^{q-1}|\nab^2(h)||\nab(h)|.
\eeaa

Next, recall from Proposition 2.1.48 in \cite{GKS22} the following identity for a scalar function $h$
\beaa
\big|\Delta_0(h)|^2 &= \big|\nab^2(h)\big|^2 +{}^{(h)}K |\nab(h)|^2 +\err_0[\Delta(h)]+\div_0[\Delta(h)],
\eeaa
with
\beaa
\err_0[\Delta(h)]&:=& -\frac 1 2 \nab(h)\c   \big(\atrch\nab_3+\atrchb \nab_4\big)\dual \nab(h), \\
\div_0[\Delta\psi]&:=& \nab_a\Big(\nab^a(h)\c  \Delta_0(h)  - \frac 1 2 \nab^a|\nab(h)|^2 \Big).
\eeaa
In view of the above, we deduce 
\bea\lab{eq:basicidentityyieldingthecontrolofrqnab2h2byrqDDhotDDhsquare}
\int_{\MM}r^q|\nab^2(h)|^2 \les \int_{\MM}r^q|\DD\hot\DD(h)|^2 +\int_{\MM}r^{q-3}|\nab\dk^{\leq 1}(h)||\dk^{\leq 1}h|+\int_{\pr\MM}r^{q-1}|\nab^2(h)||\nab(h)|.
\eea

Applying the identity \eqref{eq:basicidentityyieldingthecontrolofrqnab2h2byrqDDhotDDhsquare} respectively with:  
\begin{itemize}
\item $q=1+\de$ and $h=\dk^{J-2}\psi$,

\item $q=3-\de$ and $h=\nab_{\widehat{R}}\dk^{J-2}\psi$,

\item $q=3+\de$ and $h=\chi(r)(\nab, \nab_4)\dk^{J-2}\psi$, where the smooth cut-off function $\chi=\chi(r)$ vanishes on $\MM_{trap}$, and $\chi=1$ on $r\geq 11m$ and on $r\leq r_+(1+\dred)$,

\item $q=0$ and $\chi_{red}(r)\nab_3\dk^{J-2}\psi$ where the smooth cut-off function $\chi_{red}$ is supported in $r\leq r_+(1+2\dred)$, and $\chi_{red}=1$ on $r\leq r_+(1+\dred)$,
\end{itemize}
we obtain
\bea\lab{eq:intermediaryestimiateforHodgeestimateBdenormJminus2r2nab2whichdoesnotcoverallvaluesofr}
&&\B_{\de,r\leq r_+(1+\dred)}^{J-2}[r^2\nab^2\psi]+\B_{\de,r_+(1+2\dbl)\leq r\leq 10m}^{J-2}[r^2\nab^2\psi]+\B_{\de, r\geq 11m}^{J-2}[r^2\nab^2\psi] \nn\\
&\les& \B_\de^{J-2}[r^2\DD\hot\DD\psi]+\sqrt{\BEF^{J-1}_\de[\psi]\BEF^J_\de[\psi]}
\eea
where we recall from \eqref{eq:defofsubregionsofMM} that $\Mtrap=\MM_{r_+(1+2\dbl), 10m}$. 

Next, we need extend \eqref{eq:intermediaryestimiateforHodgeestimateBdenormJminus2r2nab2whichdoesnotcoverallvaluesofr} to $r\in[r_+(1+\dred), r_+(1+2\dbl)]\cup[10m, 11m]$. To this end, we consider the following modification of $\nab$
$$\breve{\nab}:=\nab-\frac{1}{e_3(r)}\nab(r)e_3, \qquad \breve{\nab}(r)=0, \qquad \breve{\nab}=\nab+r\Ga_g\dk.$$
This allow to derive the following analog of \eqref{eq:basicidentityyieldingthecontrolofrqnab2h2byrqDDhotDDhsquare}, for $r_1<r_2\leq 11m$, 
\bea\lab{eq:basicidentityyieldingthecontrolofrqnab2h2byrqDDhotDDhsquare:usingbrevenabasauxstep}
\int_{\MM_{r_1, r_2}(\tau_1, \tau_2)}|\nab^2(h)|^2 &\les& \int_{\MM_{r_1, r_2}(\tau_1, \tau_2)}|\DD\hot\DD(h)|^2 +\int_{\MM_{r_1, r_2}(\tau_1, \tau_2)}|\nab\dk^{\leq 1}(h)||\dk^{\leq 1}h|\nn\\
&&+\int_{\Si(\tau_1)\cup\Si(\tau_2)}|\nab^2(h)||\nab(h)|+\ep^2\int_{\MM_{r_1, r_2}(\tau_1, \tau_2)}|\dk^{\leq 2}(h)|^2\nn\\
&&+\ep^2\int_{\Si(\tau_1)\cup\Si(\tau_2)}|\dk^{\leq 2}(h)||\dk^{\leq 1}(h)|,
\eea
where we used $\breve{\nab}=\nab+r\Ga_g\dk$ to go back and forth between $\nab$ and $\breve{\nab}$ in order to do integrations by parts only using $\breve{\nab}$ which does not generate boundary terms on $r=r_1$ and $r=r_2$ thanks to $\breve{\nab}(r)=0$. Proceeding as in the proof of \eqref{eq:intermediaryestimiateforHodgeestimateBdenormJminus2r2nab2whichdoesnotcoverallvaluesofr}, and using \eqref{eq:basicidentityyieldingthecontrolofrqnab2h2byrqDDhotDDhsquare:usingbrevenabasauxstep} instead of \eqref{eq:basicidentityyieldingthecontrolofrqnab2h2byrqDDhotDDhsquare}, respectively with the choices $r_1=r_+(1+\dred)$, $r_2=r_+(1+2\dbl)$ and $r_1=10m$, $r_2=11m$, we obtain 
\beaa
&&\B_{\de,r_+(1+\dred)\leq r\leq r_+(1+2\dbl)}^{J-2}[r^2\nab^2\psi]+\B_{\de, 10m\leq r\leq 11m}^{J-2}[r^2\nab^2\psi]\\
&\les& \B_\de^{J-2}[r^2\DD\hot\DD\psi]+\sqrt{\BEF^{J-1}_\de[\psi]\BEF^J_\de[\psi]}+\ep^2\M_1^J[\psi].
\eeaa
Combining with \eqref{eq:intermediaryestimiateforHodgeestimateBdenormJminus2r2nab2whichdoesnotcoverallvaluesofr}, we infer 
\beaa
\B_{\de}^{J-2}[r^2\nab^2\psi] &\les& \B_\de^{J-2}[r^2\DD\hot\DD\psi]+\sqrt{\BEF^{J-1}_\de[\psi]\BEF^J_\de[\psi]}+\ep^2\M_1^J[\psi]
\eeaa
as stated in the first estimate. 

The second estimate follows immediately from the following bound
\beaa
\int_{n}^{n+1}\E^{J-2}[\chi(r)r^2\nab^2\psi](\tau)d\tau \les \sup_{\tau\in[n,n+1]}\Big(\E^{J-2}[r^2\DD\hot\DD\psi](\tau)+\sqrt{\E^{J-1}[\psi](\tau)\E^J[\psi](\tau)}+\ep\E^J[\psi](\tau)\Big),
\eeaa
where $\chi=\chi(r)$ a smooth cut-off function which is supported in $r_+(1+\dred/2)\leq r\leq 12m$ and equal to $1$ on $r_+(1+\dred)\leq r\leq 11m$. This bound, which is proved similarly as above and is in fact simpler, is left to the reader. This concludes the proof of Lemma \ref{lemma:howtogofromDDhotDDtonab2}. 
\end{proof}

The following lemma will allow us to obtain non-integrable Hodge estimates on the spheres $S=S(r, \tau)$.
\begin{lemma}\lab{lemma:pairoftangentvectofieldstoSdecomposedasebpluslinearcomboTandPhipluserror}
Given an orthonormal basis $e_b$, $b=1,2$, of horizontal vectorfields, there exists and orthonormal basis of vectorfields $e_b^S$, $b=1,2$, tangent to $S$ satisfying 
\beaa
e_b^S=e_b - e_b(\tau)\T -a\mu e_b(\tau)\Z+\Ga_b\dk, \quad b=1,2, \qquad \mu:=\frac{|q|^2+r^2+a^2-\De}{\Sigma(|q|^2+\Sigma)},
\eeaa
where we recall that $\Si=\sqrt{(r^2+a^2)|q|^2+2mra^2(\sin\th)^2}=\sqrt{(r^2+a^2)^2 - a^2(\sin\th)^2\De}$.
\end{lemma}

\begin{proof}
We define the following pair of vectorfields 
\beaa
X_b^S:= e_b - e_b(\tau)\T -a\mu e_b(\tau)\Z, \qquad \mu:=\frac{|q|^2+r^2+a^2-\De}{\Sigma(|q|^2+\Sigma)}, \quad b=1,2.
\eeaa
Since 
\beaa
&&\nab(\tau)=a\Re(\Jk)+\Ga_b, \quad \nab(r)=r\Ga_g, \quad \T(r)=r\Ga_b, \quad  \T(\tau)=1+r\Ga_b, \\
&&\Z(r)=r^2\Ga_g, \quad \Z(\tau)=r\Ga_b, 
\eeaa
and since $\mu=O(r^{-2})$, we infer
\beaa
X_b^S(r)=r\Ga_g, \qquad X_b^S(\tau)=\Ga_b, \quad b=1,2.
\eeaa
Also, we have
\beaa
\g(X_b^S, X_c^S) &=& \de_{bc} - e_b(\tau)\g(\T, e_c) - e_c(\tau)\g(\T, e_b) -a\mu e_b(\tau)\g(\Z, e_c) -a\mu e_c(\tau)\g(\Z, e_b)\\
&& +e_b(\tau)e_c(\tau)\Big(\g(\T, \T)+2a\mu\g(\T, \Z)+a^2\mu^2\g(\Z, \Z)\Big).
\eeaa
Using again $\nab(\tau)=a\Re(\Jk)+\Ga_b$ and $\mu=O(r^{-2})$, and since 
\beaa
&&\g(e_b, \T)=-a\Re(\Jk)_b, \quad \g(\T, \T)=-\frac{\De-a^2(\sin\th)^2}{|q|^2}, \quad \g(e_b, \Z)=(r^2+a^2)\Re(\Jk)_b,\\
&&\g(\Z, \T)=-\frac{2amr(\sin\th)^2}{|q|^2}, \quad \g(\Z, \Z)=\frac{(r^2+a^2)^2-a^2(\sin\th)^2\De}{|q|^2}(\sin\th)^2,
\eeaa
we infer
\beaa
\g(X_b^S, X_c^S) &=& \de_{bc} +a^2\Big(2 +\g(\T, \T) +2\big(-(r^2+a^2)+a\g(\T, \Z)\big)\mu\\
&& +a^2\mu^2\g(\Z, \Z)\Big)\Re(\Jk)_b\Re(\Jk)_c+r^{-1}\Ga_b\\
&=& \de_{bc} +a^2\bigg\{2 -\frac{\De-a^2(\sin\th)^2}{|q|^2} +2\bigg(-(r^2+a^2) -\frac{2a^2mr(\sin\th)^2}{|q|^2}\bigg)\mu\\
&& +a^2\mu^2\frac{\Si^2}{|q|^2}(\sin\th)^2\bigg\}\Re(\Jk)_b\Re(\Jk)_c+r^{-1}\Ga_b.
\eeaa
Now, in view of the above choice of $\mu$, we have
\beaa
&&  2-\frac{\De-a^2(\sin\th)^2}{|q|^2} +2\bigg(-(r^2+a^2) -\frac{2a^2mr(\sin\th)^2}{|q|^2}\bigg)\mu +a^2\mu^2\frac{\Si^2}{|q|^2}(\sin\th)^2\\
&=& 2-\frac{\De-a^2(\sin\th)^2}{|q|^2} -\frac{2\Si^2}{|q|^2}\mu +a^2\mu^2\frac{\Si^2}{|q|^2}(\sin\th)^2\\
&=& 2-\frac{\De-a^2(\sin\th)^2}{|q|^2} -\frac{2\Si}{|q|^2}\frac{|q|^2+r^2+a^2-\De}{|q|^2+\Sigma} +\frac{(|q|^2+r^2+a^2-\De)^2}{|q|^2(|q|^2+\Sigma)^2}a^2(\sin\th)^2\\
&=& \frac{(|q|^2+r^2+a^2+\De)^2}{|q|^2(|q|^2+\Sigma)^2}\bigg(2|q|^2(|q|^2+\Sigma)^2 -(\De - a^2(\sin\th)^2)(|q|^2+\Sigma)^2\\
&& -2\Si(|q|^2+r^2+a^2-\De)(|q|^2+\Sigma)+(|q|^2+r^2+a^2-\De)^2a^2(\sin\th)^2\bigg)\\
&=& 0
\eeaa
and hence
\beaa
\g(X_b^S, X_c^S) &=&  \de_{bc} +r^{-1}\Ga_b.
\eeaa
Thus, the vectorfields $X_b^S$ satisfy 
\beaa
X_b^S(r)=r\Ga_g, \qquad X_b^S(\tau)=\Ga_b, \quad b=1,2, \qquad \g(X_b^S, X_c^S) =  \de_{bc} +r^{-1}\Ga_b, \quad b,c=1,2.
\eeaa
Then, we define
\beaa
\widetilde{X}_b^S =X_b^S+\la_be_3+\und{\la}_be_4, \quad b=1,2,
\eeaa
where we would like to chose $\la_b$ and $\und{\la}_b$, $b=1,2$, such that $\widetilde{X}_b^S(r)=\widetilde{X}_b^S(\tau)=0$, $b=1,2$, which is equivalent to 
\beaa
e_3(r)\la_b + e_4(r)\und{\la}_b = - X_b^S(r)=r\Ga_g, \qquad e_3(\tau)\la_b + e_4(\tau)\und{\la}_b = - X_b^S(\tau)=\Ga_b,\quad b=1,2.
\eeaa
Since 
\beaa
e_3(r)=-1+r\Ga_b, \quad e_4(r)=\frac{\De}{|q|^2}+\Ga_g, \quad e_3(\tau)\gtrsim 1, \qquad e_4(\tau)\gtrsim m^2r^{-2}, 
\eeaa
we have
\beaa
e_3(r)e_4(\tau)-e_4(r)e_3(\tau) &=& -|e_4(\tau)| -\frac{\De}{|q|^2}|e_3(\tau)|+r\Ga_b\les -1.
\eeaa
Hence, the above system in $(\la_b, \und{\la}_b)$, $b=1,2$, is uniquely solvable and we have
\beaa
\la_b=\frac{-X^S(r)e_4(\tau)+X^S(\tau)e_4(r)}{e_3(r)e_4(\tau)-e_4(r)e_3(\tau)}=\Ga_b, \qquad \und{\la}_b=\frac{-e_3(r)X_b^S(\tau)+X_b^S(r)e_3(\tau)}{e_3(r)e_4(\tau)-e_4(r)e_3(\tau)}=r\Ga_g,
\eeaa
so that there indeed exists  $\la_b$ and $\und{\la}_b$, $b=1,2$, with $\la_b=\Ga_b$ and $\und{\la}_b=r\Ga_g$, $b=1,2$, such that 
\beaa
\widetilde{X}_b^S(r)=\widetilde{X}_b^S(\tau)=0, \quad b=1,2.
\eeaa
In particular, the vectorfields $\widetilde{X}_b^S$ are tangent to $S(r, \tau)$ and we have in addition 
\beaa
\g(\widetilde{X}_b^S, \widetilde{X}_c^S) &=& \g(X_b^S+\la_be_3+\und{\la}_be_4, X_c^S+\la_ce_3+\und{\la}_ce_4)\\
&=& \g(X_b^S, X_c^S)+ \la_b\g(X_c^S, e_3) + \la_c\g(X_b^S, e_3) + \und{\la}_b\g(X_c^S, e_4) + \und{\la}_c\g(X_b^S, e_4)\\
&& -2\la_b\und{\la}_c-2\la_c\und{\la}_b\\ 
&=& \de_{bc}+\Ga_g,
\eeaa
where we used the fact that 
\beaa
\g(X_b^S, X_c^S) =  \de_{bc} +r^{-1}\Ga_b, \quad \la_b=\Ga_b, \quad \und{\la}_b=r\Ga_g, \quad \g(X_b^S, e_3),\, \g(X_b^S, e_4)=O(ar^{-1})+\Ga_b.
\eeaa
Finally, we use Gram-Schmidt and define $e_b^S$, $b=1,2$ as follows
\beaa
e_1^S:=\frac{\widetilde{X}_1^S}{|\widetilde{X}_1^S|}=(1+\Ga_g)\widetilde{X}_1^S, \quad e_2^S:=\frac{\widetilde{X}_2^S-\g(\widetilde{X}_2^S, e_1^S)e_1^S}{|\widetilde{X}_2^S-\g(\widetilde{X}_2^S, e_1^S)e_1^S|}=(1+\Ga_g)\widetilde{X}_2^S+\Ga_g\widetilde{X}_1^S
\eeaa
so that $e_b^S$, $b=1,2$, is an orthonormal basis of tangent vectorfields to $S$. Furthermore, since 
\beaa
\widetilde{X}_b^S =X_b^S+\la_be_3+\und{\la}_be_4=X_b^S+\Ga_b\dk, \quad b=1,2, 
\eeaa
and in view of the definition of $X_b^S$, $b=1,2$, we infer
\beaa
e_b^S=e_b - e_b(\tau)\T -a\mu e_b(\tau)\Z+\Ga_b\dk, \quad b=1,2,
\eeaa
as stated. This concludes the proof of Lemma \ref{lemma:pairoftangentvectofieldstoSdecomposedasebpluslinearcomboTandPhipluserror}.
\end{proof}

The following lemma provides non-integrable Hodge estimates on the spheres $S=S(r, \tau)$.
\begin{lemma}\lab{lemma:improvednonintegrableHoadgeestimatefortheenergyonSigmatauusingcontroTandZderivatives}
Let $\psi$ be a scalar function. Then, we have, for any sphere $S=S(r, \tau)$ of $\MM$ and $k\leq\kl-2$,
\beaa
\int_S|r^2\nab^2\dk^k\psi|^2 &\les& \int_S|r^2\dk^k\DD\hot\DD\psi|^2 + \int_S|\dk^{k+1}\T\psi|^2 + \int_S|\dk^{k+1}\Z\psi|^2 +\int_S|\dk^{\leq k+1}\psi|^2+\ep^2\int_S|\dk^{k+2}\psi|^2,\\ 
\int_S|r^2\nab^2\dk^k\psi|^2 &\les& \int_S|r^2\dk^k\Delta\psi|^2 + \int_S|\dk^{k+1}\T\psi|^2 + \int_S|\dk^{k+1}\Z\psi|^2 +\int_S|\dk^{\leq k+1}\psi|^2+\ep^2\int_S|\dk^{k+2}\psi|^2.
\eeaa
\end{lemma}

\begin{proof}
In view of Lemma \ref{lemma:pairoftangentvectofieldstoSdecomposedasebpluslinearcomboTandPhipluserror} and basic commutator estimates, we have
\beaa
\int_S|r^2\nab^2\psi|^2 &\les& \int_S|r^2(\nab^S)^2\psi|^2 + \int_S|\dk\T\psi|^2 + \int_S|\dk\Z\psi|^2 +\int_S|\dk^{\leq 1}\psi|^2+\ep^2\int_S|\dk^2\psi|^2,
\eeaa
and 
\beaa
\int_S|r^2\DD^S\hot\DD^S\psi|^2 &\les& \int_S|r^2\DD\,\hot\DD\psi|^2 + \int_S|\dk\T\psi|^2 + \int_S|\dk\Z\psi|^2 +\int_S|\dk^{\leq 1}\psi|^2+\ep^2\int_S|\dk^2\psi|^2,
\eeaa
so that the first stated identity for $k=0$ follows from 
\beaa
\int_S|r^2(\nab^S)^2\psi|^2 &\les& \int_S|r^2\DD^S\hot\DD^S\psi|^2 +\int_S|\psi|^2
\eeaa
which is an immediate consequence of standard integrable Hodge estimates on the spheres $S(r, \tau)$. The first identity for higher order derivatives then follows by a similar reduction, this time to higher order integrable Hodge estimates on $S$. The proof of the second stated estimate follows along the same lines which concludes the proof of Lemma \ref{lemma:improvednonintegrableHoadgeestimatefortheenergyonSigmatauusingcontroTandZderivatives}.
\end{proof}

\begin{lemma}\lab{lemma:howtocontrolcrosstermsnabnabTandnabnabRpsi}
Let $\psi$ be a scalar function. Then, we have
\beaa
\B^{J-2}_\de[(r\nab)(\nab_{\widehat{T}}, \nab_{\widehat{R}})\psi] &\les& \sqrt{\B^{J-2}_\de[r^2\nab^2\psi]\B^{J-2}_\de[(\nab_{\widehat{T}}^2, \nab_{\widehat{R}}^2)\psi]}\\
&&+\sqrt{\EF^{J-1}_\de[\psi]\EF^J_\de[\psi]}+ \B^{J-2}_\de[r^2\nab^2\psi]+\B^{J-1}_\de[\psi],\\
\B^{J-2}_\de[\nab_{\widehat{T}}\nab_{\widehat{R}}\psi] &\les& \B^{J-2}_\de[\nab_{\widehat{T}}^2\psi]+\B^{J-2}_\de[\nab_{\widehat{R}}^2\psi]+\sqrt{\EF^{J-1}_\de[\psi]\EF^J_\de[\psi]}+\B^{J-1}_\de[\psi].
\eeaa
Also, for any $[n, n+1]\subset[1,\tau_*]$, we have
\beaa
\int_n^{n+1}\E_{r_+(1+\dred)\leq r\leq 11m}^{J-2}[(r\nab)(\nab_{\widehat{T}}, \nab_{\widehat{R}})\psi](\tau)d\tau &\les& \sup_{\tau\in[n,n+1]}\Big(\sqrt{\E^{J-2}[r^2\nab^2\psi](\tau)\E^{J-2}[(\nab_{\widehat{T}}^2, \nab_{\widehat{R}}^2)\psi](\tau)}\\
&&+\sqrt{\E^{J-1}[\psi](\tau)\E^J_\de[\psi](\tau)}+ \E^{J-1}[\psi](\tau)\Big),\\
\int_n^{n+1}\E^{J-2}_{r_+(1+\dred)\leq r\leq 11m}[\nab_{\widehat{T}}\nab_{\widehat{R}}\psi](\tau)d\tau &\les& \sup_{\tau\in[n,n+1]}\Big(\E^{J-2}_\de[\nab_{\widehat{T}}^2\psi](\tau)+\E^{J-2}_\de[\nab_{\widehat{R}}^2\psi](\tau)\\
&&+\sqrt{\E^{J-1}[\psi](\tau)\E^J[\psi](\tau)}+\E^{J-1}[\psi](\tau)d\tau.
\eeaa
\end{lemma}

\begin{proof}
The proof of the estimates for $\nab_{\widehat{T}}\nab_{\widehat{R}}\psi$ being similar, and in fact simpler, we focus on the ones for $(r\nab)(\nab_{\widehat{T}}, \nab_{\widehat{R}})\psi$. We start with the bulk term. Letting $\psi^{(J-2)}:=\mathfrak{d}^{\leq J-2}\psi$, we have
\beaa
\B^{J-2}_\de[(r\nab)(\nab_{\widehat{T}}, \nab_{\widehat{R}})\psi] &=& \int_{\MM}\Big(|\nab_{\widehat{R}}(r\nab)(\nab_{\widehat{T}}, \nab_{\widehat{R}})\psi^{(J-2)}|^2+|(r\nab)(\nab_{\widehat{T}}, \nab_{\widehat{R}})\psi^{(J-2)}|^2\Big)\\
&&+\int_{\MM_{trap}}|(\nab_{\widehat{T}}, \nab)(r\nab)(\nab_{\widehat{T}}, \nab_{\widehat{R}})\psi^{(J-2)}|^2\\
&\les& \int_{\MM}\Big(|\nab_{\widehat{R}}(r\nab)(\nab_{\widehat{T}}, \nab_{\widehat{R}})\psi^{(J-2)}|^2+|(r\nab)\nab_{\widehat{T}}\psi^{(J-2)}|^2\Big)\\
&&+\int_{\MM_{trap}}|\nab_{\widehat{T}}(r\nab)(\nab_{\widehat{T}}, \nab_{\widehat{R}})\psi^{(J-2)}|^2\\
&&+ B^{J-2}_\de[r^2\nab^2\psi]+B^{J-1}_\de[\psi].
\eeaa
Since $\nab$ and $\widehat{T}$ are tangent to the boundary of $\MM_{trap}$, we may integrate by parts by  $\nab$ and $\widehat{T}$ on $\MM_{trap}$, and by any derivative on $\MM$. We obtain 
\beaa
\B^{J-2}_\de[(r\nab)(\nab_{\widehat{T}}, \nab_{\widehat{R}})\psi] &\les& \sqrt{\B^{J-2}_\de[r^2\nab^2\psi]\B^{J-2}_\de[(\nab_{\widehat{T}}^2, \nab_{\widehat{R}}^2)\psi]}\\
&&+\sqrt{\EF^{J-1}_\de[\psi]\EF^J_\de[\psi]}+ \B^{J-2}_\de[r^2\nab^2\psi]+\B^{J-1}_\de[\psi].
\eeaa

The corresponding estimate for the energy norm follows immediately from the following bound
\beaa
\int_n^{n+1}\E^{J-2}[\chi(r)(r\nab)(\nab_{\widehat{T}}, \nab_{\widehat{R}})\psi](\tau)d\tau &\les& \sup_{\tau\in[n,n+1]}\Big(\sqrt{\E^{J-2}[r^2\nab^2\psi](\tau)\E^{J-2}[(\nab_{\widehat{T}}^2, \nab_{\widehat{R}}^2)\psi](\tau)}\\
&&+\sqrt{\E^{J-1}[\psi](\tau)\E^J_\de[\psi](\tau)}+ \E^{J-1}[\psi](\tau)\Big),
\eeaa
where $\chi=\chi(r)$ a smooth cut-off function which is supported in $r_+(1+\dred/2)\leq r\leq 12m$ and equal to $1$ on $r_+(1+\dred)\leq r\leq 11m$. This bound, which is proved similarly as above and is in fact simpler, is left to the reader. This concludes the proof of Lemma \ref{lemma:howtocontrolcrosstermsnabnabTandnabnabRpsi}.
\end{proof}

%%%%%%%%%%%%%%%%%%%%%%%%%%%%%%%%%%%%%%

\subsubsection{Commutators for the control of $\Pc$}
\lab{sec:commutatorsforcontrolPc}

%%%%%%%%%%%%%%%%%%%%%%%%%%%%%%%%%%%%%%

We consider the commutator $[\Rhat, \square]$. 
\begin{lemma}\lab{lemma:commutatorbetweenRhatandmod|q|2square}
For a scalar function $\psi$, we have the following commutator identity 
\beaa
[\Rhat, \square]\psi &=& O(mr^{-2})\square\psi+O(r^{-1})\Delta\psi +O(mr^{-2})\nab\That\psi+O(m^2r^{-3})\Rhat^2\psi +O(r^{-2})\dk\psi\\
&& +\dk^{\leq 1}(\Ga_g\dk\psi).
\eeaa
\end{lemma}

\begin{proof}
For a scalar function $\psi$, we have
\beaa
\square\psi &=& -\nab_3\nab_4\psi -\frac{1}{2}\trchb\nab_4\psi +\left(2\om -\frac{1}{2}\trch\right)\nab_3\psi+\De\psi+2\etab\c\nab\psi. 
\eeaa
and hence
\beaa
|q|^2\square\psi &=& -|q|^2\nab_3\nab_4\psi -\frac{1}{2}|q|^2\trchb\nab_4\psi +|q|^2\left(2\om -\frac{1}{2}\trch\right)\nab_3\psi+|q|^2\De\psi+2|q|^2\etab\c\nab\psi. 
\eeaa
Then, rewriting 
\beaa
|q|^2 \De\psi  &=& |q|\nab\c(|q|\nab\psi)  -\frac{1}{2}\nab(|q|^2)\c \nab \psi,
\eeaa
we obtain
\beaa
|q|^2\square\psi &=& -|q|^2\nab_3\nab_4\psi -\frac{1}{2}|q|^2\trchb\nab_4\psi +|q|^2\left(2\om -\frac{1}{2}\trch\right)\nab_3\psi+|q|\nab\c(|q|\nab\psi)\\
&& -\frac{1}{2}\nab(|q|^2)\c \nab \psi +2|q|^2\etab\c\nab\psi,
\eeaa
and hence
\beaa
[\Rhat, |q|^2\square]\psi &=& -\Rhat(|q|^2)\nab_3\nab_4\psi -|q|^2[\Rhat, \nab_3\nab_4]\psi +[\Rhat, |q|\nab]\c(|q|\nab\psi) +|q|\nab\c([\Rhat, |q|\nab]\psi)\\
&&+\left[\Rhat, -\frac{1}{2}|q|^2\trchb\nab_4 +|q|^2\left(2\om -\frac{1}{2}\trch\right)\nab_3 -\frac{1}{2}\nab(|q|^2)\c \nab +2|q|^2\etab\c\nab\right]\psi\\
&=& -\Rhat(|q|^2)\nab_3\nab_4\psi -|q|^2[\Rhat, \nab_3\nab_4]\psi +[\Rhat, |q|\nab]\c(|q|\nab\psi) +|q|\nab\c([\Rhat, |q|\nab]\psi)\\
&&+O(1)\dk\psi+r\dk^{\leq 1}(\Ga_b)\dk\psi.
\eeaa
Next, we use the the fact that 
\bea\lab{eq:commutationofRhatwith|q|nab}
[\nab_{\Rhat}, |q|\nab]\psi=r\Ga_b\dk\psi, \quad \psi\in\sk_0, \qquad [\nab_{\Rhat}, |q|\nab]U=O(ar^{-2})U+r\Ga_b\dk^{\leq 1}U, \quad U\in\sk_1,
\eea
where the second commutator formula follows from Corollary 9.2.3 page 325 in \cite{GKS22}, and the first can be obtained by the same proof and is in fact simpler. This yields 
\beaa
[\Rhat, |q|^2\square]\psi &=& -\Rhat(|q|^2)\nab_3\nab_4\psi -|q|^2[\Rhat, \nab_3\nab_4]\psi +O(1)\dk\psi+r\dk^{\leq 1}(\Ga_b\dk\psi).
\eeaa

It remains to deal with the commutator $[\Rhat, \nab_3\nab_4]\psi$. To this end, we use the following 
\bea
2\Rhat=Y^4e_4 - Y^3e_3, \qquad Y^4:=\frac{|q|^2}{r^2+a^2}, \qquad Y^3:=\frac{\De}{r^2+a^2},
\eea
which yields 
\beaa
2[\Rhat, \nab_3\nab_4]\psi &=& [Y^4e_4-Y^3e_3, e_3e_4]\psi\\
&=& Y^4[e_4, e_3]e_4\psi -Y^3e_3[e_3, e_4]\psi +[Y^4,e_3e_4]e_4\psi -[Y^3,e_3e_4]e_3\psi \\
&=& Y^4\big(2(\eta-\etab)\c\nab +2\om\nab_3+\Ga_b\nab_4\big)e_4\psi  \\
&&+Y^3e_3\big(2(\eta-\etab)\c\nab +2\om\nab_3+\Ga_b\nab_4\big)\psi  -e_3(Y^4)e_4^2\psi  -e_4(Y^4)e_3e_4\psi \\
&& +e_3(Y^3)e_4e_3\psi  +e_4(Y^3)e_3^2\psi -e_3e_4(Y^4)e_4\psi +e_3e_4(Y^3)e_3\psi\\
&=& \big(2(\eta-\etab)\c\nab +2\om\nab_3\big)(Y^4e_4+Y^3e_3)\psi -\big(e_4(Y^4) - e_3(Y^3)\big)e_3e_4\psi\\
&& +\pr_r(Y^4)e_4^2\psi +\frac{\De}{|q|^2}\pr_r(Y^3)e_3^2\psi +O(mr^{-3})\dk\psi+r^{-1}\dk^{\leq 1}(\Ga_b\dk\psi).
\eeaa
Now, we have
\beaa
2\That=Y^4e_4+Y^3e_3
\eeaa
and hence
\beaa
\frac{\De}{|q|^2}\pr_r(Y^3)e_3^2 &=& \pr_r(Y^3)\frac{r^2+a^2}{|q|^2}Y^3e_3^2 = \pr_r(Y^3)\frac{r^2+a^2}{|q|^2}(2\That-Y^4e_4)e_3
\eeaa
as well as 
\beaa
\pr_r(Y^4)e_4^2 &=& \pr_r(Y^4)\frac{r^2+a^2}{|q|^2}Y^4e_4^2= \pr_r(Y^4)\frac{r^2+a^2}{|q|^2}(2\That-Y^3e_3)e_4.
\eeaa
We infer
\beaa
2[\Rhat, \nab_3\nab_4]\psi  &=& \big(2(\eta-\etab)\c\nab +2\om\nab_3\big)(Y^4e_4+Y^3e_3)\psi -\big(e_4(Y^4) -e_3(Y^3)\big)e_3e_4\psi\\
&& +\pr_r(Y^4)e_4^2\psi +\frac{\De}{|q|^2}\pr_r(Y^3)e_3^2\psi+O(mr^{-3})\dk\psi+r^{-1}\dk^{\leq 1}(\Ga_b\dk\psi)\\
&=& 2\big(2(\eta-\etab)\c\nab +2\om\nab_3\big)\That -\big(e_4(Y^4) - e_3(Y^3)\big)e_3e_4\\
&& +\pr_r(Y^4)\frac{r^2+a^2}{|q|^2}(2\That-Y^3e_3)e_4+\pr_r(Y^3)\frac{r^2+a^2}{|q|^2}(2\That-Y^4e_4)e_3\\
&& +O(mr^{-3})\dk\psi+r^{-1}\dk^{\leq 1}(\Ga_b\dk\psi)\\
&=& \left(4(\eta-\etab)\c\nab +4\om\nab_3+2\pr_r(Y^4)\frac{r^2+a^2}{|q|^2}e_4 + 2\pr_r(Y^3)\frac{r^2+a^2}{|q|^2}e_3\right)\That\psi\\
&& +\left(-\big(e_4(Y^4) - e_3(Y^3)\big) -\pr_r(Y^4)\frac{r^2+a^2}{|q|^2}Y^3 - \pr_r(Y^3)\frac{r^2+a^2}{|q|^2}Y^4\right)e_3e_4\psi\\
&& +O(mr^{-3})\dk\psi+r^{-1}\dk^{\leq 1}(\Ga_b\dk\psi)\\
&=& \left(4(\eta-\etab)\c\nab +4\om\nab_3+2\pr_r(Y^4)\frac{r^2+a^2}{|q|^2}e_4 + 2\pr_r(Y^3)\frac{r^2+a^2}{|q|^2}e_3\right)\That\psi\\
&& +O(mr^{-2})e_3e_4\psi +O(mr^{-3})\dk\psi+r^{-1}\dk^{\leq 1}(\Ga_b\dk\psi).
\eeaa
Now, as
\beaa
&&2\om\nab_3+\pr_r(Y^4)\frac{r^2+a^2}{|q|^2}e_4 +\pr_r(Y^3)\frac{r^2+a^2}{|q|^2}e_3\\
&=& \left(-\pr_r\left(\frac{\De}{|q|^2}\right) +\pr_r\left(\frac{\De}{r^2+a^2}\right)\frac{r^2+a^2}{|q|^2}\right)e_3 +\pr_r\left(\frac{|q|^2}{r^2+a^2}\right)\frac{r^2+a^2}{|q|^2}e_4+\Ga_g\dk\psi\\
&=& \left(-\pr_r\left(\frac{1}{|q|^2}\right) +\pr_r\left(\frac{1}{r^2+a^2}\right)\frac{r^2+a^2}{|q|^2}\right)\De e_3 +\pr_r\left(\frac{|q|^2}{r^2+a^2}\right)\frac{r^2+a^2}{|q|^2}e_4+\Ga_g\dk\psi\\
&=& \pr_r\left(\frac{|q|^2}{r^2+a^2}\right)\frac{r^2+a^2}{|q|^2}\left(e_4+\frac{\De}{|q|^2}e_3\right)+\Ga_g\dk\psi\\
&=& \frac{1}{2}\pr_r\left(\frac{|q|^2}{r^2+a^2}\right)\left(\frac{r^2+a^2}{|q|^2}\right)^2\That+\Ga_g\dk\psi,
\eeaa
we deduce
\beaa
2[\Rhat, \nab_3\nab_4]\psi  &=& \left(4(\eta-\etab)\c\nab +\pr_r\left(\frac{|q|^2}{r^2+a^2}\right)\left(\frac{r^2+a^2}{|q|^2}\right)^2\That\right)\That\psi\\
&& +O(mr^{-2})e_3e_4\psi +O(mr^{-3})\dk\psi+\dk^{\leq 1}(\Ga_g\dk\psi)\\
&=& O(mr^{-2})\nab\That\psi+O(m^2r^{-3})\That^2\psi +O(mr^{-2})e_3e_4\psi +O(mr^{-3})\dk\psi+\dk^{\leq 1}(\Ga_g\dk\psi).
\eeaa
Plugging in the above, we deduce
\beaa
[\Rhat, |q|^2\square]\psi &=& -\Rhat(|q|^2)\nab_3\nab_4\psi -|q|^2[\Rhat, \nab_3\nab_4]\psi +O(1)\dk\psi+r\dk^{\leq 1}(\Ga_b\dk\psi)\\
&=& \big(-\Rhat(|q|^2)+O(m)\big)e_3e_4\psi +O(m)\nab\That\psi+O(m^2r^{-1})\That^2\psi +O(mr^{-1})\dk\psi\\
&&+r^2\dk^{\leq 1}(\Ga_g\dk\psi).
\eeaa

Next, we use the expression of $\square$ in a null frame to replace $e_3e_4\psi$ which yields 
\beaa
[\Rhat, |q|^2\square]\psi &=& \big(\Rhat(|q|^2)+O(m)\big)\square\psi+O(r)\Delta\psi +O(m)\nab\That\psi+O(m^2r^{-1})\That^2\psi +O(1)\dk\psi\\
&&+r^2\dk^{\leq 1}(\Ga_g\dk\psi).
\eeaa
Then, recalling from Lemma 4.7.6 in \cite{GKS22} that we have for a scalar function $\psi$
\beaa
\begin{split}
|q|^2 \square\psi &=\frac{(r^2+a^2)^2}{\De} \big( -  \nab_{\widehat{T}}^2\psi+  \nab_{\widehat{R}}^2\psi \big) +2r \nab_{\widehat{R}}\psi\\
&+  |q|^2 \De\psi   + |q|^2  (\eta+\etab) \c \nab \psi  + r^2 \Ga_g \c \mathfrak{d} \psi,
\end{split}
\eeaa
we infer
\beaa
\That^2\psi &=& \Rhat^2\psi +O(1)\square\psi +O(1)\De\psi+O(r^{-1})\dk\psi+\Ga_g\dk\psi,
\eeaa
and hence
\beaa
[\Rhat, |q|^2\square]\psi &=& \big(\Rhat(|q|^2)+O(m)\big)\square\psi+O(r)\Delta\psi +O(m)\nab\That\psi+O(m^2r^{-1})\Rhat^2\psi +O(1)\dk\psi\\
&&+r^2\dk^{\leq 1}(\Ga_g\dk\psi).
\eeaa
Finally, as 
\beaa
[\Rhat, \square]\psi =\frac{1}{|q|^2}\left([\Rhat, |q|^2\square]\psi - \Rhat(|q|^2)\square\psi\right), 
\eeaa
we deduce
\beaa
[\Rhat, \square]\psi = O(mr^{-2})\square\psi+O(r^{-1})\Delta\psi +O(mr^{-2})\nab\That\psi+O(m^2r^{-3})\Rhat^2\psi +O(r^{-2})\dk\psi +\dk^{\leq 1}(\Ga_g\dk\psi),
\eeaa
which concludes the proof of Lemma \ref{lemma:commutatorbetweenRhatandmod|q|2square}.
\end{proof}

Next, we have the following corollary of Lemma \ref{lemma:commutatorbetweenRhatandmod|q|2square}. 
\begin{corollary}\lab{cor:commutatorbetweenRhatsquareandmod|q|2square}
For a scalar function $\psi$, we have the following commutator identity 
\beaa
[\Rhat^2, \square]\psi &=& O(mr^{-2})\Rhat^{\leq 1}\square\psi +O(r^{-1})\Rhat^{\leq 1}\Delta\psi  +O(mr^{-2})\nab_{\Rhat}^{\leq 1}\nab\That\psi   +O(m^2r^{-3})\Rhat^3\psi  \\
&&+O(m^2r^{-4})\nab^2\psi+O(r^{-2})\Rhat^{\leq 1}\dk\psi +\dk^{\leq 2}(\Ga_g\dk\psi).
\eeaa
\end{corollary}

\begin{proof}
Recall from Lemma \ref{lemma:commutatorbetweenRhatandmod|q|2square} that the following commutator identity holds 
\beaa
[\Rhat, \square]\psi &=& O(mr^{-2})\square\psi+O(r^{-1})\Delta\psi +O(mr^{-2})\nab\That\psi+O(m^2r^{-3})\Rhat^2\psi +O(r^{-2})\dk\psi\\
&& +\dk^{\leq 1}(\Ga_g\dk\psi).
\eeaa
We infer
\beaa
[\Rhat^2, \square]\psi &=& \Rhat[\Rhat, \square]\psi +[\Rhat, \square]\Rhat\psi\\
&=& O(mr^{-2})\Rhat^{\leq 1}\square\psi +O(mr^{-2})\square\Rhat\psi +O(r^{-1})\Rhat^{\leq 1}\Delta\psi +O(r^{-1})\Delta\Rhat\psi\\
&& +O(mr^{-2})\nab_{\Rhat}^{\leq 1}\nab\That\psi  +O(mr^{-2})\nab\That\Rhat\psi +O(m^2r^{-3})\Rhat^3\psi  \\
&&+O(r^{-2})\Rhat^{\leq 1}\dk\psi +\dk^{\leq 2}(\Ga_g\dk\psi)
\eeaa
and hence, using again Lemma \ref{lemma:commutatorbetweenRhatandmod|q|2square} to evaluate the second term on the RHS, we obtain 
\beaa
[\Rhat^2, \square]\psi &=& O(mr^{-2})\Rhat^{\leq 1}\square\psi +O(r^{-1})\Rhat^{\leq 1}\Delta\psi +O(r^{-1})\Delta\Rhat\psi\\
&& +O(mr^{-2})\nab_{\Rhat}^{\leq 1}\nab\That\psi  +O(mr^{-2})\nab\That\Rhat\psi +O(m^2r^{-3})\Rhat^3\psi  \\
&&+O(r^{-2})\Rhat^{\leq 1}\dk\psi +\dk^{\leq 2}(\Ga_g\dk\psi). 
\eeaa

Next, we have in view of \eqref{eq:commutationofRhatwith|q|nab}
\beaa
O(r^{-1})\Delta\Rhat\psi &=& O(r^{-3})|q|^2\Delta\Rhat\psi=O(r^{-3})|q|\nab\c(|q|\nab\Rhat\psi) + O(r^{-3})\nab(|q|)|q|\nab\Rhat\psi\\
&=&  O(r^{-1})\Rhat^{\leq 1}\Delta\psi +O(mr^{-2})\Rhat^{\leq 1}\dk\psi +r^{-2}\dk^{\leq 1}(\Ga_b\dk\psi).
\eeaa
Also, we have for a scalar function $\psi$
\beaa
[\Rhat, \That]\psi &=& \left[\frac{|q|^2}{r^2+a^2}e_4 - \frac{\De}{r^2+a^2}e_3, \frac{|q|^2}{r^2+a^2}e_4 + \frac{\De}{r^2+a^2}e_3\right]\psi\\
&=&  2\left[\frac{|q|^2}{r^2+a^2}e_4,  \frac{\De}{r^2+a^2}e_3\right]\psi\\
&=& \frac{2|q|^2}{r^2+a^2}e_4\left(\frac{\De}{r^2+a^2}\right)e_3\psi - \frac{2\De}{r^2+a^2}e_3\left(\frac{|q|^2}{r^2+a^2}\right)e_4\psi +\frac{2|q|^2\De}{(r^2+a^2)^2}[e_4, e_3]\psi\\
&=& \frac{2\De}{r^2+a^2}\pr_r\left(\frac{\De}{r^2+a^2}\right)e_3\psi +\frac{2\De}{r^2+a^2}\pr_r\left(\frac{|q|^2}{r^2+a^2}\right)e_4\psi \\
&& +\frac{2|q|^2\De}{(r^2+a^2)^2}\Big(2(\etab-\eta)\c\nab\psi+2\om e_3\psi\Big)+r^{-1}\Ga_b\dk\psi,
\eeaa
or
\beaa
[\Rhat, \That]\psi &=& \frac{4\De}{|q|^2}\pr_r\left(\frac{|q|^2}{r^2+a^2}\right)\That\psi  +\frac{4|q|^2\De}{(r^2+a^2)^2}(\etab-\eta)\c\nab\psi+\Ga_g\dk\psi.
\eeaa
We deduce, using also \eqref{eq:commutationofRhatwith|q|nab}, 
\beaa
O(mr^{-2})\nab\That\Rhat\psi &=& O(mr^{-2})\nab_{\Rhat}\nab\That\psi +O(mr^{-3})\nab\That\psi+O(m^2r^{-4})\nab^2\psi +O(m^2r^{-5})\nab\psi \\
&&+r^{-2}\dk(\Ga_b\dk\psi).
\eeaa
Plugging in the above, this yields 
\beaa
[\Rhat^2, \square]\psi &=& O(mr^{-2})\Rhat^{\leq 1}\square\psi +O(r^{-1})\Rhat^{\leq 1}\Delta\psi  +O(mr^{-2})\nab_{\Rhat}^{\leq 1}\nab\That\psi   +O(m^2r^{-3})\Rhat^3\psi  \\
&&+O(m^2r^{-4})\nab^2\psi+O(r^{-2})\Rhat^{\leq 1}\dk\psi +\dk^{\leq 2}(\Ga_g\dk\psi)
\eeaa
as stated. This concludes the proof of Corollary \ref{cor:commutatorbetweenRhatsquareandmod|q|2square}.
\end{proof}
 
 We conclude this section with the following lemma on commutators between $\square_\g$ and $(\nab_3, r\nab_4)$.  
\begin{lemma}\lab{lemma:commutationof|q|^2squregwithnab3andrnab4forlocalexistenceregiontau*minus3tau*}
We have
\beaa
\frac{1}{|q|^2}[\nab_3, |q|^2\square_\g]\psi &=& O(r^{-2})\dk\nab_3\psi +\dk^{\leq 1}(\Ga_g\dk\psi)+O(r^{-2})\dk\psi,\\
\frac{1}{|q|^2}[r\nab_4, |q|^2\square_\g]\psi &=& -\frac{4}{r}\T(r\nab_4\psi)+O(r^{-2})\dk(r\nab_4\psi) +\dk^{\leq 1}(\Ga_g\dk\psi)+O(r^{-2})\dk\psi.
\eeaa
\end{lemma} 

\begin{proof}
For a scalar function $\psi$, we have, in view of Lemma \ref{lemma:expression-wave-operator},
\beaa
\square_\g\psi &=& -\nab_3\nab_4\psi +\left(2\omb -\frac 1 2 \trchb\right)\nab_4\psi -\frac 1 2 \trch \nab_3\psi+\lap\psi+2\eta \c\nab \psi. 
\eeaa
We then multiply with $|q|^2$ and commute respectively with $\nab_3$ and $r\nab_4$ which yields, in view of the commutation formulas \eqref{eq:comm-nab3-nab4-naba-f-general}, 
\beaa
[\nab_3, |q|^2\square_\g]\psi &=& -[\nab_3, |q|^2\nab_3\nab_4]\psi +\left[\nab_3, |q|^2\left(2\omb -\frac 1 2 \trchb\right)\nab_4\right]\psi -\frac 1 2 [\nab_3, \trch |q|^2\nab_3]\psi\\
&& +[\nab_3, |q|^2\lap]\psi+2[\nab_3, |q|^2\eta \c\nab]\psi\\
&=& O(1)\dk\nab_3\psi +r^2\dk^{\leq 1}(\Ga_g\dk\psi)+O(1)\dk\psi,
\eeaa
and
\beaa
[r\nab_4, |q|^2\square_\g]\psi &=& -[r\nab_4, |q|^2\nab_3\nab_4]\psi +\left[r\nab_4, |q|^2\left(2\omb -\frac 1 2 \trchb\right)\nab_4\right]\psi -\frac 1 2 [r\nab_4, \trch |q|^2\nab_3]\psi\\
&& +[r\nab_4, |q|^2\lap]\psi+2[r\nab_4, |q|^2\eta \c\nab]\psi\\
&=& -2r\nab_3(r\nab_4\psi)+O(1)\dk(r\nab_4\psi) +r^2\dk^{\leq 1}(\Ga_g\dk\psi)+O(1)\dk\psi\\
&=& -4r\T(r\nab_4\psi)+O(1)\dk(r\nab_4\psi) +r^2\dk^{\leq 1}(\Ga_g\dk\psi)+O(1)\dk\psi,
\eeaa
where we have used in particular the following non-sharp consequences of Lemma 4.7.10 in \cite{GKS22}
\beaa
\,[\nab_3, |q|^2\De]\psi &=& O(1)\dk\nab_3\psi +r^2\dk^{\leq 1}(\Ga_g\dk\psi)+O(1)\dk\psi, \\
\,[r\nab_4, |q|^2\De]\psi &=& O(1)\dk(r\nab_4\psi) +r^2\dk^{\leq 1}(\Ga_g\dk\psi)+O(1)\dk\psi.
\eeaa
We infer
\beaa
\frac{1}{|q|^2}[\nab_3, |q|^2\square_\g]\psi &=& O(r^{-2})\dk\nab_3\psi +\dk^{\leq 1}(\Ga_g\dk\psi)+O(r^{-2})\dk\psi,\\
\frac{1}{|q|^2}[r\nab_4, |q|^2\square_\g]\psi &=& -\frac{4}{r}\T(r\nab_4\psi)+O(r^{-2})\dk(r\nab_4\psi) +\dk^{\leq 1}(\Ga_g\dk\psi)+O(r^{-2})\dk\psi,
\eeaa
as stated. This concludes the proof of Lemma \ref{lemma:commutationof|q|^2squregwithnab3andrnab4forlocalexistenceregiontau*minus3tau*}.
\end{proof}

%%%%%%%%%%%%%%%%%%%%%%%%%%%%%%%%%%%%

\subsection{Energy-Morawetz estimates for $A, B, \Bb, \Ab$}
\lab{sec:energyMorawetzforABBbAb}

%%%%%%%%%%%%%%%%%%%%%%%%%%%%%%%%%%%%

The goal of this section is to prove Propositions \ref{proposition:EstimatesBBb-interior} and \ref{proposition:EstimatesAAb-interior}. We start by deriving useful non-integrable Hodge estimates.

%%%%%%%%%%%%%%%%%%%%%%%%%%%%%%%%%%%%%%%%%%%%%%%%%%%%%%%%%%%%%%%%%%%%%%%%%%%%%%%%%%

\subsubsection{Non-integrable Hodge estimates for the proof of Propositions \ref{proposition:EstimatesBBb-interior} and \ref{proposition:EstimatesAAb-interior}}
\lab{sec:proofof:proposition:EstimatesBBb-interior:nonintegrableHodgeestimates}

%%%%%%%%%%%%%%%%%%%%%%%%%%%%%%%%%%%%%%%%%%%%%%%%%%%%%%%%%%%%%%%%%%%%%%%%%%%%%%%%%%

Our first non-integrable Hodge estimates are derived in the region $r\geq r_1$ where $r_1=r_1(m)\geq 20m$ will be chosen large enough.
\begin{lemma}\lab{lemma:HodgeestimatesforthecontrolofBBbandAAb:largerregionrgeqr1}
Let $\psi$ be a horizontal tensor in $\sk_j(\mathbb{C})$, $j=1,2$. Then, we have, for any sphere $S=S(r, \tau)$ of $\MM$, and for $r\geq r_1$ and $k\leq\kl-1$,
\beaa
\int_S|r\nab\dk^k\psi|^2 &\les& \int_S|r\dk^k\ov{\DD}\c\psi|^2 + \int_S|\dk^{k}\nab_{\T}\psi|^2  +\int_S|\dk^{\leq k}\psi|^2+(r_1^{-2}+\ep^2)\int_S|\dk^{k+1}\psi|^2,\\ 
\max_{k_1+k_2=k}\int_S|(r\nab)^{k_1+1}\nab_{\T}^{k_2}\psi|^2 &\les& \int_S|r\dk^k\ov{\DD}\c\psi|^2 + \int_S|\nab_{\T}^{k+1}\psi|^2  +\int_S|\dk^{\leq k}\psi|^2+(r_1^{-2}+\ep^2)\int_S|\dk^{k+1}\psi|^2.
\eeaa
\end{lemma}

\begin{proof}
Recall from Lemma 9.3.1 in \cite{GKS22} that given an orthonormal basis $e_b$, $b=1,2$, of horizontal vectorfields, there exists and orthonormal basis of vectorfields $e_b^S$, $b=1,2$, tangent to $S$ satisfying 
\beaa
\eS_b&=& e_b+O(mr^{-1} )\T +O(mr^{-2})\dk+r\Ga_g( e_3, e_4). 
\eeaa
This implies 
\beaa
\int_S|r\nab\psi|^2 &\les& \int_S|r\nab^S\psi|^2 + \int_S|\nab_{\T}\psi|^2  +\int_S|\psi|^2+(r_1^{-2}+\ep^2)\int_S|\dk\psi|^2
\eeaa
and
\beaa
\int_S|r\ov{\DD}^S\c\psi|^2 &\les& \int_S|r\ov{\DD}\c\psi|^2 + \int_S|\nab_{\T}\psi|^2  +\int_S|\psi|^2+(r_1^{-2}+\ep^2)\int_S|\dk\psi|^2
\eeaa
so that the two stated identity for $k=0$ follow from 
\beaa
\int_S|r\nab^S\psi|^2 &\les& \int_S|r\ov{\DD}^S\c\psi|^2 +\int_S|\psi|^2
\eeaa
which is an immediate consequence of standard integrable Hodge estimates on the spheres $S(r, \tau)$. These identities for higher order derivatives then follow by a similar reduction, this time to higher order integrable Hodge estimates on $S$. This concludes the proof of Lemma \ref{lemma:HodgeestimatesforthecontrolofBBbandAAb:largerregionrgeqr1}.
\end{proof}

Then, we consider non-integrable Hodge estimates on $\Si_{r_+\leq r\leq r_1}(\tau)$.
\begin{lemma}\lab{lemma:HodgeestimatesforthecontrolofBBbandAAb:regionrplusleqrleqr1}
Let $\psi$ be a horizontal tensor in $\sk_j(\mathbb{C})$, $j=1,2$. Then, we have, for $p\in\mathbb{R}$ and $k\leq\kl-1$,
\beaa
\int_{\Si_{r_+\leq r\leq r_1}(\tau)}r^p|r\nab\dk^k\psi|^2 &\les& \int_{\Si(\tau)}r^p\Big(|r\dk^k\ov{\DD}\c\psi|^2 + |\dk^k\nab_3\psi|^2 + r_1|\dk^{\leq k+1}\psi||\dk^{\leq k}\psi|\Big),\\
\int_{\Si_{r_+\leq r\leq r_1}(\tau)}r^p|r\nab\dk^k\psi|^2 &\les& \int_{\Si(\tau)}r^p\Big(|r\dk^k\ov{\DD}\c\psi|^2 + r_1^4|\dk^k\nab_4\psi|^2 + r_1|\dk^{\leq k+1}\psi||\dk^{\leq k}\psi|\Big).
\eeaa
\end{lemma}

\begin{proof}
To prove the first estimate, we introduce the following vectorfields 
\beaa
(e_{\Si(\tau)})_a=e_a -\frac{e_a(\tau)}{e_3(\tau)}e_3, \quad a=1,2, \qquad \g((e_{\Si(\tau)})_a, (e_{\Si(\tau)})_b)=\de_{ab}, \quad a,b=1,2,
\eeaa
so that $(e_{\Si(\tau)})_a$, $a=1,2$ are tangent to $\Si(\tau)$ and orthonormal. We also introduce   
\beaa
\bsplit
(e_{\Si(\tau)})_4:=&e_4 -\frac{2}{e_3(\tau)}\nab\tau+\frac{|\nab\tau|^2}{(e_3(\tau))^2}e_3, \\ 
N_{\Si(\tau)}:=&\la\left(e_4-\frac{e_4(\tau)}{e_3(\tau)}e_3-\frac{2\nab^b(\tau)}{e_3(\tau)}(e_{\Si(\tau)})_b\right),\qquad \la:=\frac{e_3(\tau)}{2\sqrt{|\g(\D\tau, \D\tau)|}},
\end{split}
\eeaa 
which verify
\beaa
&&\g((e_{\Si(\tau)})_4, (e_{\Si(\tau)})_a)=0, \quad a=1,2, \qquad \g((e_{\Si(\tau)})_4, (e_{\Si(\tau)})_4)=0, \qquad \g(e_3, (e_{\Si(\tau)})_4)=-2,\\
&&\g(N_{\Si(\tau)}, (e_{\Si(\tau)})_a)=0, \quad a=1,2,\qquad \g(N_{\Si(\tau)}, N_{\Si(\tau)})=1.
\eeaa
In particular, $(e_{\Si(\tau)})_a$, $a=1,2$, is an orthonormal basis of the horizontal structure $\HH_{\Si(\tau)}:=\{e_3, (e_{\Si(\tau)})_4\}^\perp$ and $((e_{\Si(\tau)})_1, (e_{\Si(\tau)})_2, N_{\Si(\tau)})$ is an orthonormal basis of the tangent space of $\Si(\tau)$. We then denote by $\nab_{\Si(\tau)}$ the horizontal covariant operator associated to the horizontal structure $\HH_{\Si(\tau)}$ and $\widetilde{\nab}_{\Si(\tau)}$ the covariant derivative on $\Si(\tau)$. 

Next, we rely on a smooth cut-off $\chi=\chi(r)$ such that $\chi$ is supported in $r\in[r_+(1-\dhor/2), 2r_1]$ and $\chi=1$ in $r\in[r_+, r_1]$, and we estimate 
\beaa
\int_{\Si_{r_+\leq r\leq r_1}(\tau)}r^p|r\nab\psi|^2 &\les& \int_{\Si(\tau)}r^p\chi(\tau)|r\nab\psi|^2\\
&\les& \int_{\Si(\tau)}r^p\chi(\tau)|r\nab_{\Si(\tau)}\psi|^2 + \int_{\Si(\tau)}r^p|\nab_3\psi|^2
\eeaa
Now, the following identity holds in view of Proposition 2.1.37 and Proposition 2.1.43 in \cite{GKS22}, for $\psi\in\sk_j(\mathbb{C})$, $j=1,2$,
\beaa
|\nab_{\Si(\tau)}\psi|^2 &=& \frac{j}{4}|\ov{\DD}_{\Si(\tau)}\c\psi|^2 + \frac{i}{j}\big(\atrch(\nab_{\Si(\tau)})_3+\atrchb(\nab_{\Si(\tau)})_4\big)\psi\c\psi \\
&&+(\nab_{\Si(\tau)})_a \left( \nab_{\Si(\tau)}^a\psi \c\psi - \frac{j}{2}\ov{\nab_{\Si(\tau)}}\psi\c\psi \right).
\eeaa
Also, denoting  by $\textrm{Div}_{\Si(\tau)}$ the divergence on $\Si(\tau)$, we have for a horizontal 1-form $f$
\beaa
\textrm{Div}_{\Si(\tau)}(f) =\div_{\Si(\tau)}(f) -f_a\g(\D_{N_{\Si(\tau)}}N_{\Si(\tau)}, e_a)
\eeaa
which together with the above formula for $N_{\Si(\tau)}$ implies 
\beaa
\textrm{Div}_{\Si(\tau)}(f) =\div_{\Si(\tau)}(f)+(O(1)+r^2\dk^{\leq 1}\Ga_g)f
\eeaa
We then multiply the above identity for $|\nab_{\Si(\tau)}\psi|^2$ by $r^{p+2}\chi(\tau)$, integrate on $\Si(\tau)$ and use integration by parts to obtain
\beaa
\int_{\Si(\tau)}r^p\chi(\tau)|r\nab_{\Si(\tau)}\psi|^2 &\les& \int_{\Si(\tau)}r^p\Big(|r\ov{\DD}_{\Si(\tau)}\c\psi|^2 + r_1|\dk\psi||\psi|\Big)\\
&\les& \int_{\Si(\tau)}r^p\Big(|r\ov{\DD}\c\psi|^2 + |\nab_3\psi|^2 + r_1|\dk\psi||\psi|\Big),
\eeaa
where we used the decomposition of $(e_{\Si(\tau)})_a$. Plugging in the above, we infer
\beaa
\int_{\Si_{r_+\leq r\leq r_1}(\tau)}r^p|r\nab\psi|^2 &\les& \int_{\Si(\tau)}r^p\Big(|r\ov{\DD}\c\psi|^2 + |\nab_3\psi|^2 + r_1|\dk\psi||\psi|\Big)
\eeaa
which is the first stated estimate for $k=0$. The proof of the first estimate for $1\leq k\leq \kl-1$ proceeds similarly. The proof of the second stated estimate then proceeds as above by interchanging $e_3$ and $e_4$ in the definition of $(e_{\Si(\tau)})_a$, $a=1,2$, and $N_{\Si(\tau)}$, and noticing that $e_3(\tau)=2+O(m^2r^{-2})+r\Ga_b$ while $e_4(\tau)=\frac{m^2}{r^2}(1+O(mr^{-1}))+\Ga_g$ which accounts for the different weights in $r_1$ in the two stated estimates. This concludes the proof of Lemma \ref{lemma:HodgeestimatesforthecontrolofBBbandAAb:regionrplusleqrleqr1}.
\end{proof}

Finally, we consider non-integrable Hodge estimates in $r\leq r_+$.
\begin{lemma}\lab{lemma:HodgeestimatesforthecontrolofBBbandAAb:regionrleqrplus}
Let $\psi$ be a horizontal tensor in $\sk_j(\mathbb{C})$, $j=1,2$. Then, we have, for any sphere $S=S(r, \tau)$ of $\MM$, and for $r\leq r_+$ and $k\leq\kl-1$,
\beaa
\max_{k_1+k_2=k}\int_S|(r\nab)\dk^{k_1+1}\nab_4^{k_2}\psi|^2 &\les& \int_S\Big(|r\dk^k\ov{\DD}\c\psi|^2 + |\nab_4^{k+1}\psi|^2 + |\dk^k\nab_3\psi|^2+|\dk^{\leq k}\psi|^2\Big),\\
\max_{k_1+k_2=k}\int_S|(r\nab)\dk^{k_1+1}\nab_3^{k_2}\psi|^2 &\les& \int_S\Big(|r\dk^k\ov{\DD}\c\psi|^2 + |\nab_3^{k+1}\psi|^2 + |\dk^k\nab_4\psi|^2+|\dk^{\leq k}\psi|^2\Big).
\eeaa
\end{lemma}

\begin{proof}
On $r_+(1-\dhor)\leq r\leq r_+$, we have
\beaa
\left(e_2-\frac{a\sin\th|q|}{2(r^2+a^2)}e_4\right)\tau = O(\dhor)+O(\ep), \qquad \left(e_2-\frac{a\sin\th|q|}{2(r^2+a^2)}e_4\right)r = O(\dhor)+O(\ep),
\eeaa
so that, given an orthonormal basis $e_b$, $b=1,2$, of horizontal vectorfields, there exists and orthonormal basis of vectorfields $e_b^S$, $b=1,2$, tangent to $S$ satisfying 
\beaa
\eS_b&=& e_b+O(1)e_4 +O(\dhor)\dk+O(\ep)\dk.
\eeaa
This implies 
\beaa
\int_S|\nab\psi|^2\les (1+O(\dhor)+O(\ep))\int_S|\nab^S\psi|^2+\int_S|\nab_4\psi|^2+(O(\dhor)+O(\ep))\int_S|\nab_3\psi|^2
\eeaa
and
\beaa
\int_S|\ov{\DD}_S\c\psi|^2\les (1+O(\dhor)+O(\ep))\int_S|\ov{\DD}\c\psi|^2+\int_S|\nab_4\psi|^2+(O(\dhor)+O(\ep))\int_S|\nab_3\psi|^2
\eeaa
so that the two stated identity for $k=0$ follow, for $\ep>0$ and $\dhor>0$ small enough, 
\beaa
\int_S|r\nab^S\psi|^2 &\les& \int_S|r\ov{\DD}^S\c\psi|^2 +\int_S|\psi|^2
\eeaa
which is an immediate consequence of standard integrable Hodge estimates on the spheres $S(r, \tau)$. These identities for higher order derivatives then follow by a similar reduction, this time to higher order integrable Hodge estimates on $S$. This concludes the proof of Lemma \ref{lemma:HodgeestimatesforthecontrolofBBbandAAb:regionrleqrplus}.
\end{proof}

%%%%%%%%%%%%%%%%%%%%%%%%%%%%%%%%%%%%%%%%%%%%%%%%%%%%%%%%%%%%%

\subsubsection{Energy-Morawetz estimates for $B, \Bb$ (proof of Proposition \ref{proposition:EstimatesBBb-interior})}
\lab{sec:proofof:proposition:EstimatesBBb-interior}

%%%%%%%%%%%%%%%%%%%%%%%%%%%%%%%%%%%%%%%%%%%%%%%%%%%%%%%%%%%%%

The proof of Proposition \ref{proposition:EstimatesBBb-interior} adapts the derivation of energy-Morawetz estimates for $B, \Bb$ in Section 15.5 of  \cite{GKS22}. 
The restriction on the size of $a$ in that proof is essentially due to the use of non-integrable Hodge estimates to recover bounds for higher order energies and fluxes. To remove this restriction, unlike the proof in Section 15.5 in \cite{GKS22}, we first control the bulk norms, and then recover flux and energy norms relying in particular on the non-integrable Hodge estimates of Lemmas \ref{lemma:HodgeestimatesforthecontrolofBBbandAAb:largerregionrgeqr1}, \ref{lemma:HodgeestimatesforthecontrolofBBbandAAb:regionrplusleqrleqr1} and \ref{lemma:HodgeestimatesforthecontrolofBBbandAAb:regionrleqrplus}. We proceed in the following steps. 

\noindent{\bf Step 1.} We start with the following immediate consequence of Bianchi identities for the pair $(B, \Pc)$, see (15.5.4) in \cite{GKS22},
\bea
\lab{eq:EstmatesforB-1}
\bsplit
\BEF^{J-1}_\de[ r^2 \nabc_3B]+ \BEF^{J-1}_\de[ r^3 \DDov\c B  ] &\les \de_{J+1}[\Pc] +\Sk_J \Sk_{J+1} +\Rk_{J} \Rk_{J+1} +\ep_0^2.
\end{split}
\eea
Then, we control of $\B^{J-1}_\de[ r^3\nab B]$. To this end, we rely on the Bochner identity (2.1.22) in \cite{GKS22} and Lemma 2.1.40 in \cite{GKS22} which yields after integration by parts 
\beaa
\B^{J-1}_\de[r^3\nab B] &\les& \B^{J-1}_\de[r^3\DDov\c B] + \Big(\BEF^J_\de[r^2B]\Big)^{\frac{1}{2}}\Big(\BEF^{J-1}_\de[r^2B]\Big)^{\frac{1}{2}}. 
\eeaa
Together with \eqref{eq:EstmatesforB-1} and \eqref{induction:hypoth-Pc}, this yields
\bea\lab{eq:BEFder2nab3BBEFr3DDovBBr3nabB}
\nn&&\BEF^{J-1}_\de[ r^2 \nabc_3B]+ \EF^{J-1}_\de[ r^3 \DDov\c B  ] + \B^{J-1}_\de[r^3\nab B]\\
&\les& \de_{J+1}[\Pc] +\Sk_J \Sk_{J+1} +\Rk_{J} \Rk_{J+1} +\Big(\BEF^J_\de[r^2B]\Big)^{\frac{1}{2}}\Big(\Rk_{J} \Rk_{J+1}\Big)^{\frac{1}{2}}+\ep_0^2.
\eea
Also, proceeding similarly for the Bianchi pair $(\Pc, \Bb)$, we obtain the following analog of \eqref{eq:BEFder2nab3BBEFr3DDovBBr3nabB}
\bea
\lab{eq:BEFdernab4BbBEFrDDovBbBrnabBb}
\nn&&\BEF^{J-1}_\de[r\nabc_4\Bb] + \EF^{J-1}_\de[r\DDov\c\Bb]  + \B^{J-1}_\de[r\nab\Bb]   \\
&\les& \de_{J+1}[\Pc] +\Sk_J \Sk_{J+1} +\Rk_{J} \Rk_{J+1} +\Big(\BEF^J_\de[\Bb]\Big)^{\frac{1}{2}}\Big(\Rk_{J} \Rk_{J+1}\Big)^{\frac{1}{2}} +\ep_0^2.
\eea

\noindent{\bf Step 2.} Next, as in (15.5.8) of \cite{GKS22}, we introduce the following quantities 
  \bea\lab{eq:definitionofBtPtpPtmandBbttorecovermissingderivativeforBandBb}
    \Bt:=(\ov{q}\nabc_4 )^{\le J-2}  \nabc^2_{\widetilde{R}}\Bdot,  \qquad      \Bbt:=(\nabc_3 )^{\le J-2}\nabc^2_{\widetilde{R}}\Bbdot, \qquad \Bdot:=\Lieb_{\T}B, \qquad \Bbdot:=\Lieb_{\T}\Bb,
   \eea  
   where the vectorfield $\widetilde{R}$ is given by \eqref{def:Rhat}. We first estimate $\Bt$ by revisiting the proof of Proposition 15.5.7 in \cite{GKS22}. As in that proof, we have, see (15.5.15) in \cite{GKS22}, for $b=2+\de$,
 \bea
 \lab{eq:EstimateBt:1}
   \int_{\MM } r^{b-1} |\Bt |^2+\int_{\pr^+\MM} r^{b}|\Bt |^2 &\les&  \de_{J+1}[\Pc] +  \int_{\MM} r^{b+1}   \big|\Ft_{(1, J-2 )}'\big|^2\nn\\
   &&+ \sqrt{ \de_{J+1}[\Pc]} \left(\int_{\MM} r^{b+1} \big|\Ft_{(2, J-2 )}\big|^2\right)^{\frac{1}{2}}
   +\big|I\big| +\ep_0^2,
  \eea
  where $\Ft_{(1, J-2)}'$, $\Ft_{(2, J-2)}$ and $I$ are given by
   \beaa
 \Ft_{(1, J-2)}' &=&   O(r^{-1} )  \dk^{J-2} \nab_{\widetilde{R}} \big(\nabc_3  \Bdot, \dkb \Bdot\big)  + O(r^{-1})  \dk^{J-1} \nab_{\widetilde{R}} \Pdot +   O(r^{-1})  \dk^{\le J-1}\Pdot\\
 && + O(r^{-3} )\dk^{\le J}  \Ga_b  + r^{-2}   \dk^{\le J+1}(\Ga_b\c \Rc_b) + r^{-2} \dk^{\le J} \big( \Ga_b\c\Ga_b\big),\\
   \Ft_{(2, J-2)} &=& O(r^{-1})  \dk^{J-1} \nab_{\widetilde{R}} (\Pdot, \Bdot)+ O(r^{-3} )   \dk^{\le J+1}   \Ga_b + O(r^{-1})  \dk^{\le J-1}(\Pdot, \Bdot) + r^{-2}   \dk^{\le J+1}(\Ga_b\c \Rc_b),\\
I &=& \int_{\MM}    O(r^{-3} )    |q|^{b}\Re(    (\nab_3, \nab)\dk^{\le J}\Ga_b \c\ov{\Bt}) +\int_{\MM}O(r^{-1})|q|^{b}\Re(\dk^{\le J}B\c\ov{\Bt}),
\eeaa
with $\Pdot:=\T(P)$ and $\Bdot:=\Lieb_{\T}B$. Then, still following the proof of Proposition 15.5.7 in \cite{GKS22}, we have, 
  \beaa
&&\int_{\MM} r^{b+1}   \big|\Ft_{(1, J-2 )}'\big|^2 + \sqrt{ \de_{J+1}[\Pc]} \left(\int_{\MM} r^{b+1} \big|\Ft_{(2, J-2 )}\big|^2\right)^{\frac{1}{2}}\\
&\les& \B^{J-1}_\de[ r^2(\nab_3B, r\nab B)] +\de_{J+1}[\Pc]  +\ep_J^2   +\ep_0^2 +\sqrt{ \de_{J+1}[\Pc]} \left(\B^J_\de[r^2B]+\Sk_{J+1}^2\right)^{\frac{1}{2}}.
\eeaa
Departing from the proof of Proposition 15.5.7 in \cite{GKS22}, we now rely on \eqref{eq:BEFder2nab3BBEFr3DDovBBr3nabB} to infer
\begin{align*}
&\int_{\MM} r^{b+1}   \big|\Ft_{(1, J-2 )}'\big|^2 + \sqrt{ \de_{J+1}[\Pc]} \left(\int_{\MM} r^{b+1} \big|\Ft_{(2, J-2 )}\big|^2\right)^{\frac{1}{2}}\\
\les& \de_{J+1}[\Pc] +\Sk_J \Sk_{J+1} +\Rk_{J} \Rk_{J+1} +\ep_0^2  +\sqrt{ \de_{J+1}[\Pc]} \left(\B^J_\de[r^2B]+\Sk_{J+1}^2\right)^{\frac{1}{2}}+\Big(\BEF^J_\de[r^2B]\Big)^{\frac{1}{2}}\Big(\Rk_{J} \Rk_{J+1}\Big)^{\frac{1}{2}}.
\end{align*}
 Back to \eqref{eq:EstimateBt:1} we infer the following improvement of (15.5.16) in \cite{GKS22}
  \begin{align}
 \lab{eq:EstimateBt:2}
   \int_{\MM } r^{b-1} |\Bt |^2+\sup_{\tau\leq\tau_*}\int_{\pr^+\MM(1,\tau)}r^{b}|\Bt |^2 \les& \big|I\big|+ \de_{J+1}[\Pc] +\Sk_J \Sk_{J+1} +\Rk_{J} \Rk_{J+1} +\ep_0^2 \\
   &+\Big(\BEF^J_\de[r^2B]\Big)^{\frac{1}{2}}\Big(\Rk_{J} \Rk_{J+1}\Big)^{\frac{1}{2}}+\sqrt{ \de_{J+1}[\Pc]} \left(\B^J_\de[r^2B]+\Sk_{J+1}^2\right)^{\frac{1}{2}}. \nn
  \end{align}
  Then, as in \cite{GKS22}, we decompose 
   \beaa
  I = I_1+I_2,\qquad  I_1= \int_{\MM}    O(r^{-3} )    |q|^{b}\Re(  (\nab_3, \nab)\dk^{\le J}\Ga_b \c\ov{\Bt}),\qquad  I_2 = \int_{\MM}O(r^{-1})|q|^{b}\Re(\dk^{\le J}B\c\ov{\Bt}).
  \eeaa
  For $I_1$, we do not use its final estimate in \cite{GKS22}, but instead the estimate immediately before, i.e.,
   \beaa
     |I_1|  &\les& \Sk_{J+1}\Big(\B^{J-1}_\de[ r^2(\nab_3B, r\nab B)]\Big)^{\frac{1}{2}}+\ep_J\Big(\B^J_\de[ r^2B]\Big)^{\frac{1}{2}}+\ep_0^2.
     \eeaa
  Also, we use the following estimate for $I_2$ in \cite{GKS22}
  \beaa
|I_2| &\les& \ep_J\Big(\B_\de^{J}[r^2B]\Big)^{\frac{1}{2}}+\ep_0^2.
\eeaa
We infer 
  \beaa
 |I| &\les& |I_1|+|I_2| \les  \Sk_{J+1}\Big(\B^{J-1}_\de[ r^2(\nab_3B, r\nab B)]\Big)^{\frac{1}{2}}+\ep_J\Big(\B^J_\de[ r^2B]\Big)^{\frac{1}{2}}+\ep_0^2.
     \eeaa
 Together with \eqref{eq:BEFder2nab3BBEFr3DDovBBr3nabB}, this yields the following improvement of the control of $I$ in \cite{GKS22}
   \beaa
 |I|  &\les&  \Sk_{J+1}\left(\de_{J+1}[\Pc] +\Sk_J \Sk_{J+1} +\Rk_{J} \Rk_{J+1} +\ep_0^2+\Big(\BEF^J_\de[r^2B]\Big)^{\frac{1}{2}}\Big(\Rk_{J} \Rk_{J+1}\Big)^{\frac{1}{2}}\right)^{\frac{1}{2}}\\
 &&+\ep_J\Big(\B^J_\de[ r^2B]\Big)^{\frac{1}{2}}+\ep_0^2.
     \eeaa
 Finally, plugging the above estimate for $I$ in \eqref{eq:EstimateBt:2}, and recalling that $b=2+\de$, we infer the following improvement of (15.5.13) in \cite{GKS22}
  \begin{align}\lab{eq:lemma-EstimatesBtPtp3}
   &\int_{\MM } r^{1+\de} |\Bt |^2 +\sup_{\tau\leq\tau_*}\int_{\pr^+\MM(1,\tau)}r^{2+\de}|\Bt |^2\nn\\
    \les&  \de_{J+1}[\Pc] +\Sk_J \Sk_{J+1} +\Rk_{J} \Rk_{J+1} +\ep_0^2+\sqrt{ \de_{J+1}[\Pc]} \left(\B^J_\de[r^2B]+\Sk_{J+1}^2\right)^{\frac{1}{2}}+\ep_J\Big(\B^J_\de[ r^2B]\Big)^{\frac{1}{2}}\nn\\
   &+\Sk_{J+1}\left(\de_{J+1}[\Pc] +\Sk_J \Sk_{J+1} +\Rk_{J} \Rk_{J+1} +\ep_0^2+\Big(\BEF^J_\de[r^2B]\Big)^{\frac{1}{2}}\Big(\Rk_{J} \Rk_{J+1}\Big)^{\frac{1}{2}}\right)^{\frac{1}{2}}.
  \end{align}

\noindent{\bf Step 3.}  Next, we estimate $\Bbt$ introduced in \eqref{eq:definitionofBtPtpPtmandBbttorecovermissingderivativeforBandBb} by revisiting the proof of Proposition 15.5.8 in \cite{GKS22}. As in that proof, we have, for $b=-\de$,
\beaa
 \int_{\MM } r^{b-1}   | \Bbt |^2 +\sup_{\tau\leq\tau_*}\int_{\pr^+\MM(1,\tau)} r^{ b-2 } |\Bbt |^2 &\les &  \big|I_2 \big| + \de_{J+1}[\Pc]   + \sqrt{\de_{J+1}[\Pc]}\left(\int_{\MM } r^{b+1} |\Ft_{(1,J-2)}  |^2\right)^{\frac{1}{2}}\nn\\
 && + \int_{\MM } r^{b+1} |\Ft_{(2,J-2)}' |^2 + \ep_0^2,
\eeaa
where $\Ft_{(1,J-2)}$, $\Ft'_{(2,J-2)}$ and $I_2$ are given by
\beaa
\Ft_{(1,J-2)} &=& O(r^{-1})\dk^{J-1} \nab_{\widetilde{R}} (\Pdot, \Bbdot)+ O(r^{-3} )   \dk^{\le J+1}   \Ga_b + O(r^{- 1 })  \dk^{\le J-1}(\Pdot, \Bbdot) + r^{-1}   \dk^{\le J+1}(\Ga_b\c \Rc_b),\\
\Ft'_{(2,J-2)}   &=& O(r^{-1} )\dk^{\le J-1} \nab_{\widetilde{R}}\Pdot +O(r^{-2} ) \dk^{\le J-2}\dkb  \nab_{\widetilde{R}}\Bbdot  + O(r^{-1})  \dk^{\le J-1}\Pdot+ O(r^{-1})  \dk^{\le J}\Bb \\
   &&+ O(r^{-3}) \dk^{\le J}\Ga_b+  \dk^{\le J+1 } \Big( \Ga_b  \c \Rc_b\Big),\\
I_2 &=& \int_{\MM}|q|^{b}\Re\Big(   O(r^{-3}) \dk^{\le J}(\nab_4, \nab)\Ga_b\c  \ov{\Bbt}\Big),
\eeaa
with $\Pdot:=\T(P)$ and $\Bbdot:=\Lieb_{\T}\Bb$. We then estimate  the integrals $\int_{\MM } r^{b+1} |\Ft_{(1,J-2)}  |^2 $ and  $\int_{\MM } r^{b+1} |\Ft'_{(2,J-2)} |^2$. Taking advantage of the structure of $\Ft_{(1,J-2)}$ and $\Ft_{(2,J-2)}'$, proceeding as for the corresponding estimate in Step 2 above, and, unlike in \cite{GKS22}, using in particular  \eqref{eq:BEFdernab4BbBEFrDDovBbBrnabBb} to control $\B^{J-1}_\de[ r(\nab_4\Bb, \nab\Bb)]$, we obtain 
\beaa
&&\sqrt{\de_{J+1}[\Pc]}\left(\int_{\MM } r^{b+1} | F_{(1,J-2)}  |^2\right)^{\frac{1}{2}} + \int_{\MM } r^{b+1} |  F_{(2,J-2)}' |^2 \\
&\les&  \de_{J+1}[\Pc] +\Sk_J \Sk_{J+1} +\Rk_{J} \Rk_{J+1} +\ep_0^2 +\sqrt{ \de_{J+1}[\Pc]} \left(\B^J_\de[\Bb]+\Sk_{J+1}^2\right)^{\frac{1}{2}} +\Big(\BEF^J_\de[\Bb]\Big)^{\frac{1}{2}}\Big(\Rk_{J} \Rk_{J+1}\Big)^{\frac{1}{2}}.
\eeaa
We deduce 
\begin{align}\lab{eq:EstimateBbt:2}
 \int_{\MM } r^{b-1}   | \Bbt |^2 +\sup_{\tau\leq\tau_*}\int_{\pr^+\MM(1,\tau)} r^{ b-2 } |\Bbt |^2 \les &  \big|I_2 \big|  +\de_{J+1}[\Pc] +\Sk_J \Sk_{J+1} +\Rk_{J} \Rk_{J+1} +\ep_0^2\\
 & +\sqrt{ \de_{J+1}[\Pc]} \left(\B^J_\de[\Bb]+\Sk_{J+1}^2\right)^{\frac{1}{2}}+\Big(\BEF^J_\de[\Bb]\Big)^{\frac{1}{2}}\Big(\Rk_{J} \Rk_{J+1}\Big)^{\frac{1}{2}}. \nn
\end{align}
Then, the term $I_2$ can be integrated by parts twice,  as the corresponding estimate in Step 2 above, to obtain 
 \beaa
     |I_2|      &\les&  \Sk_{J+1}\Big(\B^{J-1}_\de[ r(\nab_4\Bb, \nab\Bb)]\Big)^{\frac{1}{2}}+\ep_J\Big(\B^J_\de[\Bb]\Big)^{\frac{1}{2}}+\ep_0^2\\
   &\les&  \Sk_{J+1}\left(\de_{J+1}[\Pc] +\Sk_J \Sk_{J+1} +\Rk_{J} \Rk_{J+1} +\ep_0^2+\Big(\BEF^J_\de[\Bb]\Big)^{\frac{1}{2}}\Big(\Rk_{J} \Rk_{J+1}\Big)^{\frac{1}{2}}\right)^{\frac{1}{2}}\\
 &&+\ep_J\Big(\B^J_\de[ r^2B]\Big)^{\frac{1}{2}}+\ep_0^2,
     \eeaa
     where we used \eqref{eq:BEFdernab4BbBEFrDDovBbBrnabBb} in the last inequality. Finally, plugging the above estimate for $I_2$ in \eqref{eq:EstimateBbt:2}, and recalling that $b=-\de$, we deduce
\begin{align}\lab{eq:EstimateBbt:3}
 &\int_{\MM } r^{-1-\de}   | \Bbt |^2 +\sup_{\tau\leq\tau_*}\int_{\pr^+\MM(1,\tau)} r^{-2-\de} |\Bbt |^2\nn\\ 
 \les &  \de_{J+1}[\Pc] +\Sk_J \Sk_{J+1} +\Rk_{J} \Rk_{J+1} +\ep_0^2 +\sqrt{ \de_{J+1}[\Pc]} \left(\B^J_\de[\Bb]+\Sk_{J+1}^2\right)^{\frac{1}{2}} +\Big(\BEF^J_\de[\Bb]\Big)^{\frac{1}{2}}\Big(\Rk_{J} \Rk_{J+1}\Big)^{\frac{1}{2}}\nn\\
+&\ep_J\Big(\B_\de^{J}[\Bb]\Big)^{\frac{1}{2}} + \Sk_{J+1}\left(\de_{J+1}[\Pc] +\Sk_J \Sk_{J+1} +\Rk_{J} \Rk_{J+1} +\ep_0^2+\Big(\BEF^J_\de[\Bb]\Big)^{\frac{1}{2}}\Big(\Rk_{J} \Rk_{J+1}\Big)^{\frac{1}{2}}\right)^{\frac{1}{2}}.
\end{align}

\noindent{\bf Step 4.} In this step, we recover $\B^{J}_\de[r^2B]$. First, since 
  \beaa
  2\widetilde{R}=e_4+O(1)e_3, \qquad 2\T=e_4+O(1)e_3+O(1)\nab, 
  \eeaa
  and since $\Bt =(\ov{q}\nabc_4 )^{\le J-2}  \nabc^2_{\widetilde{R}}\Bdot$, we have
  \beaa
  \B_{\de}[r^{J-1}\nab_4^JB]  &\les& \int_{\MM}r^{\de+1}|\Bt|^2 + \B^{J-1}_\de[r^2(\nab_3, r\nab)B]  +\B^{J-1}_\de[r^2B],
  \eeaa
which together with \eqref{eq:BEFder2nab3BBEFr3DDovBBr3nabB} and \eqref{eq:lemma-EstimatesBtPtp3} implies, using also  \eqref{induction:hypoth-Pc},
 \bea\lab{eq:BdeofrJminus1ofnab4JB}
  \B_{\de}[r^{J-1}\nab_4^JB]  &\les&  \de_{J+1}[\Pc] +\Sk_J \Sk_{J+1} +\Rk_{J} \Rk_{J+1} +\ep_0^2+\sqrt{ \de_{J+1}[\Pc]} \left(\B^J_\de[r^2B]+\Sk_{J+1}^2\right)^{\frac{1}{2}}\nn\\
   &&+\Sk_{J+1}\left(\de_{J+1}[\Pc] +\Sk_J \Sk_{J+1} +\Rk_{J} \Rk_{J+1} +\ep_0^2+\Big(\BEF^J_\de[r^2B]\Big)^{\frac{1}{2}}\Big(\Rk_{J} \Rk_{J+1}\Big)^{\frac{1}{2}}\right)^{\frac{1}{2}}\nn\\
   &&+\ep_J\Big(\B^J_\de[ r^2B]\Big)^{\frac{1}{2}}+\Big(\BEF^J_\de[r^2B]\Big)^{\frac{1}{2}}\Big(\Rk_{J} \Rk_{J+1}\Big)^{\frac{1}{2}}.
  \eea
In particular, restricting to $r\leq 12m$, we infer, in view of \eqref{eq:BEFder2nab3BBEFr3DDovBBr3nabB} and \eqref{eq:BdeofrJminus1ofnab4JB}, using also  \eqref{induction:hypoth-Pc},
\bea\lab{eq:BdeofrJr2B:sofaronlyrleq12m}
  \B_{\de, r\leq 12m}^J[r^2B]  &\les&  \de_{J+1}[\Pc] +\Sk_J \Sk_{J+1} +\Rk_{J} \Rk_{J+1} +\ep_0^2+\sqrt{ \de_{J+1}[\Pc]} \left(\B^J_\de[r^2B]+\Sk_{J+1}^2\right)^{\frac{1}{2}}\nn\\
   &&+\Sk_{J+1}\left(\de_{J+1}[\Pc] +\Sk_J \Sk_{J+1} +\Rk_{J} \Rk_{J+1} +\ep_0^2+\Big(\BEF^J_\de[r^2B]\Big)^{\frac{1}{2}}\Big(\Rk_{J} \Rk_{J+1}\Big)^{\frac{1}{2}}\right)^{\frac{1}{2}}\nn\\
   &&+\ep_J\Big(\B^J_\de[ r^2B]\Big)^{\frac{1}{2}}+\Big(\BEF^J_\de[r^2B]\Big)^{\frac{1}{2}}\Big(\Rk_{J} \Rk_{J+1}\Big)^{\frac{1}{2}}.
  \eea

Next, we recover $\B_{\de, r\geq 12m}^J[r^2B]$. To this end, as in Step 1 of Section 15.5.5 in \cite{GKS22}, we rely on
 \bea
\begin{split}
\lab{eq:third-pair-e_4P-e_4B-new}
\nabc_3 \Psi_{(1) }+ \left(1-\frac{J+1}{2}\right)\tr\Xb \Psi_{(1)} &= - \DDs_1\,\Psi_{(2)}  +F_{(1)}, \\
\nabc_4\Psi_{(2)}+ \left(\frac{3}{2} -\frac{J+1}{2}\right)\ov{\tr X}\,\Psi_{(2)}  &=\DDd_1\Psi_{(1)}+F_{(2)},
\end{split}
\eea
where $\Psi_{(1)}=\chi_{far}(\ov{q}\nabc_4)^{J+1}B$, $\Psi_{(2)}=\chi_{far}(\ov{q}\nabc_4)^J(\ov{q}\widecheck{\nab_4\ov{P}})$, with, changing a little bit from \cite{GKS22}, $\chi_{far}(r)=0$ on $r\leq 11m$ and  $\chi_{far}(r)=1$ for $r\geq 12m$, and with $F_{(1)}, F_{(2)}$  of the form
\beaa
\bsplit
 F_{(1)} =& \chi_{far}(r)\Big[O(r^{-1} )  \dk^J\big(\nabc_3B, \dkb B\big)  + O(r^{-1})  \dk^{\leq J+1}\Pc +  O(r^{-3} )(\nab_3, \nab)\dk^{\le J}\Ga_b  \\
 &+  O(r^{-1})  \dk^{\le J}B+ O(r^{-3} )\dk^{\le J}  \Ga_b +r^{-2}\dk^{\le J+1}(\Ga_b\c \Rc_b)\Big]  + \chi_{far}'(r)\dk^{\leq J+1}(\Pc, B),\\
  F_{(2)} =&  \chi_{far}(r)\Big[O(r^{-1})+ \dk^{\leq J+1}(\Pc, B)+ O(r^{-3} )   \dk^{\le J+1}   \Ga_b +   r^{-2}\dk^{\le J+1}(\Ga_b\c \Rc_b)\Big] + \chi_{far}'(r)\dk^{\leq J+1}(\Pc, B).
\end{split}
\eeaa 
Then, as in Step 2 of Section 15.5.5 in \cite{GKS22}, we obtain, with $b=2+\de$, 
\beaa
\begin{split}
& \int_{\MM } r^{b-1} |\Psi_{(1)}|^2+\int_{\pr^+\MM} r^{b}|\Psi_{(1)}|^2\\
\les &\int_{\MM} r^{b-1} |\Psi_{(2)}|^2  +\left|\int_{\MM}|q|^{b}\Re(F_{(1)}\c\ov{\Psi_{(1)}})\right|+\int_{\MM}r^b|F_{(2)}||\Psi_{(2)}|   + \int_{\pr^-\MM} \big( r^{b}|\Psi_{(1)}|^2+  r^{b-2}|\Psi_{(2)}|^2\big).
\end{split}
\eeaa
Now, the terms on the RHS involving $\chi_{far}(r)$ are controlled using \eqref{eq:BEFder2nab3BBEFr3DDovBBr3nabB} while terms on the RHS involving $\chi_{far}'(r)$ are controlled using \eqref{eq:BdeofrJr2B:sofaronlyrleq12m}. This yields, since $b=2+\de$ and $\chi_{far}(r)=1$ for $r\geq 12m$, 
 \bea\lab{eq:intermediaryalmostdoneesitmateforBaisdhfalif:chap15bis}
  \nn  &&\int_{\MM_{r\geq 12m}} r^{1+\de}|(\ov{q}\nab_4)^{J+1}B|^2+\sup_{\tau\leq \tau_*}\int_{\pr^+\MM_{r\geq 12m}(1,\tau)} r^{2+\de}|(\ov{q}\nab_4)^{J+1}B|^2 \\
  \nn  &\les& \de_{J+1}[\Pc] +\Sk_J \Sk_{J+1} +\Rk_{J} \Rk_{J+1} +\Big(\BEF^J_\de[r^2B]\Big)^{\frac{1}{2}}\Big(\Rk_{J} \Rk_{J+1}\Big)^{\frac{1}{2}}+\ep_0^2\\
\nn  &&+\Sk_{J+1}\left(\de_{J+1}[\Pc] +\Sk_J \Sk_{J+1} +\Rk_{J} \Rk_{J+1} +\ep_0^2+\Big(\BEF^J_\de[r^2B]\Big)^{\frac{1}{2}}\Big(\Rk_{J} \Rk_{J+1}\Big)^{\frac{1}{2}}\right)^{\frac{1}{2}}\\
   &&+\sqrt{ \de_{J+1}[\Pc]} \left(\B^J_\de[r^2B]+\Sk_{J+1}^2\right)^{\frac{1}{2}}+\ep_J\Big(\B^J_\de[ r^2B]\Big)^{\frac{1}{2}}.
  \eea
Together with \eqref{eq:BEFder2nab3BBEFr3DDovBBr3nabB} and \eqref{eq:BdeofrJr2B:sofaronlyrleq12m}, we deduce, using also the fact that $\chi_{far}(r)=1$ for $r\geq 12m$,
\beaa
 \B^J_{\de}[r^2B] &\les& \de_{J+1}[\Pc] +\Sk_J \Sk_{J+1} +\Rk_{J} \Rk_{J+1} +\Big(\BEF^J_\de[r^2B]\Big)^{\frac{1}{2}}\Big(\Rk_{J} \Rk_{J+1}\Big)^{\frac{1}{2}}+\ep_0^2\\
\nn  &&+\Sk_{J+1}\left(\de_{J+1}[\Pc] +\Sk_J \Sk_{J+1} +\Rk_{J} \Rk_{J+1} +\ep_0^2+\Big(\BEF^J_\de[r^2B]\Big)^{\frac{1}{2}}\Big(\Rk_{J} \Rk_{J+1}\Big)^{\frac{1}{2}}\right)^{\frac{1}{2}}\\
   &&+\sqrt{ \de_{J+1}[\Pc]} \left(B^J_\de[r^2B]+\Sk_{J+1}^2\right)^{\frac{1}{2}}+\ep_J\Big(B^J_\de[ r^2B]\Big)^{\frac{1}{2}},
  \eeaa
  and hence
\begin{align}\lab{eq:intermediaryalmostdoneesitmateforBaisdhfalif:chap15:stillmissingsomeEFnorms}
\B^J_{\de}[r^2B] \les& \de_{J+1}[\Pc] +\Sk_J \Sk_{J+1} +\Rk_{J} \Rk_{J+1} +\Big(\EF^J[r^2B]\Big)^{\frac{1}{2}}\Big(\Rk_{J} \Rk_{J+1}\Big)^{\frac{1}{2}}+\ep_0^2+\Sk_{J+1}^{\frac{4}{3}}\Rk_{J+1}^{\frac{1}{3}}\Rk_{J}^{\frac{1}{3}}+\ep_J^2\nn\\
  +&\Sk_{J+1}\left(\de_{J+1}[\Pc] +\Sk_J \Sk_{J+1} +\Rk_{J} \Rk_{J+1} +\ep_0^2+\Big(\EF^J[r^2B]\Big)^{\frac{1}{2}}\Big(\Rk_{J} \Rk_{J+1}\Big)^{\frac{1}{2}}\right)^{\frac{1}{2}}.
  \end{align}

\noindent{\bf Step 5.}  In this step, we recover $\B^{J}_\de[\Bb]$. First, recall the following identity, see Lemma 15.5.9 in \cite{GKS22},
 \beaa
     \left( \frac{\De}{|q|^2 } \right)^2 \nab_{\widetilde{R}} \nab_3^{ J }\Bb = -2\Bbt  + O(1 )  \nab_{\widetilde{R}}  \dk^{\le J- 1 }(\nab_4, \nab)\Bb + O(r^{-1} )\nab_{\widetilde{R}}\dk^{\le J-1}\Bb+ O(r^{-1} )  \dk^{\le J-1}\Bb.
\eeaa
Together with \eqref{eq:BEFdernab4BbBEFrDDovBbBrnabBb}, \eqref{eq:EstimateBbt:3} and  \eqref{induction:hypoth-Pc}, we infer
    \begin{align}
    \lab{estimate:degenerate-Bb1}
   &  \int_{\MM } r^{-1-\de}  \left( \frac{\De}{|q|^2 } \right)^4       \big|\nab_{\widetilde{R}} \nab_3 ^{ J }  \Bb\big|^2    \nn\\
    \les &  \de_{J+1}[\Pc] +\Sk_J \Sk_{J+1} +\Rk_{J} \Rk_{J+1} +\ep_0^2 +\sqrt{ \de_{J+1}[\Pc]} \left(\B^J_\de[\Bb]+\Sk_{J+1}^2\right)^{\frac{1}{2}} +\Big(\BEF^J_\de[\Bb]\Big)^{\frac{1}{2}}\Big(\Rk_{J} \Rk_{J+1}\Big)^{\frac{1}{2}}\nn\\
+&\ep_J\Big(\B_\de^{J}[\Bb]\Big)^{\frac{1}{2}} + \Sk_{J+1}\left(\de_{J+1}[\Pc] +\Sk_J \Sk_{J+1} +\Rk_{J} \Rk_{J+1} +\ep_0^2+\Big(\BEF^J_\de[\Bb]\Big)^{\frac{1}{2}}\Big(\Rk_{J} \Rk_{J+1}\Big)^{\frac{1}{2}}\right)^{\frac{1}{2}}.
    \end{align}

Note that \eqref{estimate:degenerate-Bb1} is degenerate in the redshift region $r\leq r_+(1+\deh)$. To get rid of this degeneracy, we introduce a smooth cut-off function $\chi_{red}(r)$ such that $\chi_{red}(r)=1$ for $r\le r_+(1+\deh)$ and $\chi_{red}(r)=0$ for $r\geq r_+(1+2\deh)$ and define 
\beaa
\Psi_{(1)}:=\chi_{red}\nabc_3^{ J }\widecheck{ \nabc_3  P}, \qquad \Psi_{(2)}:=\chi_{red} \nabc_3^{ J+1 } \Bb,
\eeaa
which satisfies, see (15.5.25) in \cite{GKS22},
  \bea
\begin{split}
\lab{eq:third-pair-e_3P-e_3Bb-new}
\nabc_3 \Psi_{(1) }+ \frac{3}{2}\ov{\tr\Xb} \Psi_{(1)} &= - \DDd_1\Psi_{(2)}  +F_{(1)}, \\
\nabc_4\Psi_{(2)}+ \tr X\Psi_{(2)}  &=\DDs_1\, \Psi_{(1)}+F_{(2)},
\end{split}
\eea
with $F_{(1)}, F_{(2)}$ supported in the region  $ r\leq  r_+( 1+2\de_{red} )$ and of the form
\beaa
 F_{(1)} &=&  \chi_{red}(r)\Big[O(1) \dk^{\le J+1} (\Pc,  \Bb)  + O(1)\dk^{\le J+1}\Ga_b  +\dk^{\le J+1} \big( \Ga_b\c \Rc_b\big)\Big]+\chi_{red}'(r)\dk^{\leq J+1}(\Pc, \Bb),\\
  F_{(2)}&=&  \chi_{red}(r)\Big[O(1) \dk^{\le J}\nab B +O(1)\dk^{\le J}(\nab_4, \nab)\Ga_b+O(1) \dk^{\le J+1}\Pc\\
  && +O(1) \dk^{\le J}\Bb+O(1) \dk^{\le J}\Ga_b  +\dk^{\le J+1} \big( \Ga_b\c \Rc_b\big)\Big]+\chi_{red}'(r)\dk^{\leq J+1}(\Pc, \Bb).
\eeaa
Next, as in Step 3 of the proof of Lemma 15.5.9 in \cite{GKS22}, we obtain, with $b=-\de$,
\beaa
\begin{split}
& \int_{\MM} r^{b-1}\left(1+|2k+1|\frac{m}{r}\right) | \Psi_{(2)}|^2 +  \int_{\pr^+\MM} r^{b-2} |\Psi_{(2)}|^2\\ 
\les &\int_{\MM} r^{b-1} |\Psi_{(1)}|^2    +\left|\int_{\MM}|q|^{b}\Re(F_{(2)}\c\ov{\Psi_{(2)}})\right|    +\int_{\MM}r^{b}|F_{(1)}||\Psi_{(1)}|  +  \int_{\pr^-\MM} \big( r^{b}|\Psi_{(1)}|^2+  r^{b-2}|\Psi_{(2)}|^2\big).
 \end{split}
\eeaa 
Plugging the above formula for $F_{(1)}, F_{(2)}$, using the fact that $F_{(1)}, F_{(2)}$ are supported in $ r\leq  r_+( 1+2\de_{red} )$ and hence away from trapping, and relying on \eqref{eq:BEFdernab4BbBEFrDDovBbBrnabBb}, \eqref{estimate:degenerate-Bb1} and  \eqref{induction:hypoth-Pc}, we infer
    \begin{align}
    \lab{estimate:nondegenerate-Bb2}
   &  \int_{\MM }\big|\chi_{red}\nab_3^{J+1}\Bb\big|^2        +\sup_{\tau\leq\tau_*}\int_{\pr^+ \MM(1,\tau)}\big|\chi_{red}\nab_3^{J+1}  \Bb\big|^2  \nn\\
 \les&  \de_{J+1}[\Pc] +\Sk_J \Sk_{J+1} +\Rk_{J} \Rk_{J+1} +\ep_0^2 +\sqrt{ \de_{J+1}[\Pc]} \left(\B^J_\de[\Bb]+\Sk_{J+1}^2\right)^{\frac{1}{2}} +\Big(\BEF^J_\de[\Bb]\Big)^{\frac{1}{2}}\Big(\Rk_{J} \Rk_{J+1}\Big)^{\frac{1}{2}}\nn\\
+&\ep_J\Big(\B_\de^{J}[\Bb]\Big)^{\frac{1}{2}} + \Sk_{J+1}\left(\de_{J+1}[\Pc] +\Sk_J \Sk_{J+1} +\Rk_{J} \Rk_{J+1} +\ep_0^2+\Big(\BEF^J_\de[\Bb]\Big)^{\frac{1}{2}}\Big(\Rk_{J} \Rk_{J+1}\Big)^{\frac{1}{2}}\right)^{\frac{1}{2}}.
    \end{align}
Then, \eqref{eq:BEFdernab4BbBEFrDDovBbBrnabBb}, \eqref{estimate:degenerate-Bb1}, \eqref{estimate:nondegenerate-Bb2} and  \eqref{induction:hypoth-Pc} yield
    \begin{align*}
     \B_\de^J[\Bb] \les&  \de_{J+1}[\Pc] +\Sk_J \Sk_{J+1} +\Rk_{J} \Rk_{J+1} +\ep_0^2 +\sqrt{ \de_{J+1}[\Pc]} \left(\B^J_\de[\Bb]+\Sk_{J+1}^2\right)^{\frac{1}{2}} +\Big(\BEF^J_\de[\Bb]\Big)^{\frac{1}{2}}\Big(\Rk_{J} \Rk_{J+1}\Big)^{\frac{1}{2}}\nn\\
&+\ep_J\Big(\B_\de^{J}[\Bb]\Big)^{\frac{1}{2}} + \Sk_{J+1}\left(\de_{J+1}[\Pc] +\Sk_J \Sk_{J+1} +\Rk_{J} \Rk_{J+1} +\ep_0^2+\Big(\BEF^J_\de[\Bb]\Big)^{\frac{1}{2}}\Big(\Rk_{J} \Rk_{J+1}\Big)^{\frac{1}{2}}\right)^{\frac{1}{2}},
    \end{align*}
and hence
    \bea\lab{eq:intermediaryalmostdoneesitmateforBbaisdhfalif:chap15:stillmissingsomeEFnorms}
     \B_\de^J[\Bb] &\les&  \de_{J+1}[\Pc] +\Sk_J \Sk_{J+1} +\Rk_{J} \Rk_{J+1} +\ep_0^2 +\Big(\EF^J_\de[\Bb]\Big)^{\frac{1}{2}}\Big(\Rk_{J} \Rk_{J+1}\Big)^{\frac{1}{2}}+\Sk_{J+1}^{\frac{4}{3}}\Rk_{J+1}^{\frac{1}{3}}\Rk_{J}^{\frac{1}{3}}\nn\\
&&+ \Sk_{J+1}\left(\de_{J+1}[\Pc] +\Sk_J \Sk_{J+1} +\Rk_{J} \Rk_{J+1} +\ep_0^2+\Big(\EF^J_\de[\Bb]\Big)^{\frac{1}{2}}\Big(\Rk_{J} \Rk_{J+1}\Big)^{\frac{1}{2}}\right)^{\frac{1}{2}}.
    \eea

\noindent{\bf Step 6.} In view of \eqref{eq:intermediaryalmostdoneesitmateforBaisdhfalif:chap15:stillmissingsomeEFnorms} and \eqref{eq:intermediaryalmostdoneesitmateforBbaisdhfalif:chap15:stillmissingsomeEFnorms}, we still need to recover $\EF^J_\de[r^2B]$ and $\EF^J_\de[\Bb]$. First, we consider $B$ in the region $r\geq r_1$ for $r_1\gg m$ chosen large enough below. Relying on the first estimate in Lemma \ref{lemma:HodgeestimatesforthecontrolofBBbandAAb:largerregionrgeqr1}, we have
\beaa
\EF^{J-1}_{\de, r\geq r_1}[r^3\nab B] &\les& \EF^{J-1}_\de[r^3\DDov\c B] + \EF^{J-1}_\de[ r^2 \nabc_3B]+\EF^{J-1}_\de[r^2B]+(r_1^{-2}+\ep^2)\EF^{J}_\de[r^2B],
\eeaa
which together with \eqref{eq:BEFder2nab3BBEFr3DDovBBr3nabB}, \eqref{eq:intermediaryalmostdoneesitmateforBaisdhfalif:chap15bis} and  \eqref{induction:hypoth-Pc} implies, taking also \eqref{eq:intermediaryalmostdoneesitmateforBaisdhfalif:chap15:stillmissingsomeEFnorms} into account, 
\beaa
\EF^{J}_{\de, r\geq r_1}[r^2B] &\les& r_1^{-2}\EF^{J}_\de[r^2B] +\de_{J+1}[\Pc] +\Sk_J \Sk_{J+1} +\Rk_{J} \Rk_{J+1} +\Big(\EF^J_\de[r^2B]\Big)^{\frac{1}{2}}\Big(\Rk_{J} \Rk_{J+1}\Big)^{\frac{1}{2}}+\ep_0^2\\
\nn  &&+\Sk_{J+1}\left(\de_{J+1}[\Pc] +\Sk_J \Sk_{J+1} +\Rk_{J} \Rk_{J+1} +\ep_0^2+\Big(\EF^J_\de[r^2B]\Big)^{\frac{1}{2}}\Big(\Rk_{J} \Rk_{J+1}\Big)^{\frac{1}{2}}\right)^{\frac{1}{2}},
\eeaa
and hence, for $r_1\gg m$ large enough, we infer
\begin{align}\lab{eq:EFJnormofr2Bintheregionrgeqr1}
\EF^{J}_{\de, r\geq r_1}[r^2B] \les& \de_{J+1}[\Pc] +\Sk_J \Sk_{J+1} +\Rk_{J} \Rk_{J+1} +\Big(\EF^J_\de[r^2B]\Big)^{\frac{1}{2}}\Big(\Rk_{J} \Rk_{J+1}\Big)^{\frac{1}{2}}+\ep_0^2\nn\\
  &+\Sk_{J+1}\left(\de_{J+1}[\Pc] +\Sk_J \Sk_{J+1} +\Rk_{J} \Rk_{J+1} +\ep_0^2+\Big(\EF^J_\de[r^2B]\Big)^{\frac{1}{2}}\Big(\Rk_{J} \Rk_{J+1}\Big)^{\frac{1}{2}}\right)^{\frac{1}{2}}.
\end{align}
We now fix $r_1$ such that \eqref{eq:EFJnormofr2Bintheregionrgeqr1} holds and we apply the first estimate of Lemma \ref{lemma:HodgeestimatesforthecontrolofBBbandAAb:regionrplusleqrleqr1} which yields
\beaa
\EF^{J-1}_{\de, r_+\leq r\leq r_1}[r^3\nab B] &\les& \EF^{J-1}_\de[r^3\DDov\c B] + \EF^{J-1}_\de[ r^2 \nabc_3B]+\Big(\EF^{J-1}_\de[r^2B]\Big)^{\frac{1}{2}}\Big(\EF^{J}_\de[r^2B]\Big)^{\frac{1}{2}}.
\eeaa
Together with \eqref{eq:BEFder2nab3BBEFr3DDovBBr3nabB}, \eqref{eq:lemma-EstimatesBtPtp3}, \eqref{eq:intermediaryalmostdoneesitmateforBaisdhfalif:chap15bis}, \eqref{eq:EFJnormofr2Bintheregionrgeqr1} and  \eqref{induction:hypoth-Pc}, we infer
\begin{align}\lab{eq:EFJnormofr2Bintheregionrgeqrplus}
\EF^{J}_{\de, r\geq r_+}[r^2B] \les& \de_{J+1}[\Pc] +\Sk_J \Sk_{J+1} +\Rk_{J} \Rk_{J+1} +\Big(\EF^J_\de[r^2B]\Big)^{\frac{1}{2}}\Big(\Rk_{J} \Rk_{J+1}\Big)^{\frac{1}{2}}+\ep_0^2\\
  &+\Sk_{J+1}\left(\de_{J+1}[\Pc] +\Sk_J \Sk_{J+1} +\Rk_{J} \Rk_{J+1} +\ep_0^2+\Big(\EF^J_\de[r^2B]\Big)^{\frac{1}{2}}\Big(\Rk_{J} \Rk_{J+1}\Big)^{\frac{1}{2}}\right)^{\frac{1}{2}}. \nn
\end{align}

Next, we consider $\Bb$ in the region $r\geq r_1$ for $r_1\gg m$ chosen large enough below. Relying on the second estimate in Lemma \ref{lemma:HodgeestimatesforthecontrolofBBbandAAb:largerregionrgeqr1}, the fact that $2\T=e_3+O(r^{-1})\dk$ and the definition of $\Bbt$ in \eqref{eq:definitionofBtPtpPtmandBbttorecovermissingderivativeforBandBb}, we have
\beaa
\EF^{J-1}_{\de, r\geq r_1}[r\nab\Bb] &\les& \EF^{J-1}_\de[r\DDov\c\Bb] + \EF^{J-1}_\de[r\nabc_4\Bb] +\sup_{\tau\leq\tau_*}\int_{\pr^+\MM(1,\tau)} r^{-2-\de} |\Bbt |^2 +\EF^{J-1}_\de[\Bb]\\
&&+(r_1^{-2}+\ep^2)\EF^{J}_\de[\Bb],
\eeaa
which together with \eqref{eq:BEFdernab4BbBEFrDDovBbBrnabBb}, \eqref{eq:EstimateBbt:3} and  \eqref{induction:hypoth-Pc} implies, taking also \eqref{eq:intermediaryalmostdoneesitmateforBbaisdhfalif:chap15:stillmissingsomeEFnorms} into account, 
\beaa
\EF^{J}_{\de, r\geq r_1}[\Bb] &\les& r_1^{-2}\EF^{J}_\de[\Bb] +\de_{J+1}[\Pc] +\Sk_J \Sk_{J+1} +\Rk_{J} \Rk_{J+1} +\Big(\EF^J_\de[\Bb]\Big)^{\frac{1}{2}}\Big(\Rk_{J} \Rk_{J+1}\Big)^{\frac{1}{2}}+\ep_0^2\\
\nn  &&+\Sk_{J+1}\left(\de_{J+1}[\Pc] +\Sk_J \Sk_{J+1} +\Rk_{J} \Rk_{J+1} +\ep_0^2+\Big(\EF^J_\de[\Bb]\Big)^{\frac{1}{2}}\Big(\Rk_{J} \Rk_{J+1}\Big)^{\frac{1}{2}}\right)^{\frac{1}{2}},
\eeaa
and hence, for $r_1\gg m$ large enough, we infer
\begin{align}\lab{eq:EFJnormofBbintheregionrgeqr1}
\EF^{J}_{\de, r\geq r_1}[\Bb] \les& \de_{J+1}[\Pc] +\Sk_J \Sk_{J+1} +\Rk_{J} \Rk_{J+1} +\Big(\EF^J_\de[\Bb]\Big)^{\frac{1}{2}}\Big(\Rk_{J} \Rk_{J+1}\Big)^{\frac{1}{2}}+\ep_0^2\nn\\
  &+\Sk_{J+1}\left(\de_{J+1}[\Pc] +\Sk_J \Sk_{J+1} +\Rk_{J} \Rk_{J+1} +\ep_0^2+\Big(\EF^J_\de[\Bb]\Big)^{\frac{1}{2}}\Big(\Rk_{J} \Rk_{J+1}\Big)^{\frac{1}{2}}\right)^{\frac{1}{2}}.
\end{align}
We now fix $r_1$ such that \eqref{eq:EFJnormofBbintheregionrgeqr1} holds and we apply the second estimate of Lemma \ref{lemma:HodgeestimatesforthecontrolofBBbandAAb:regionrplusleqrleqr1} which yields
\beaa
\EF^{J-1}_{\de, r_+\leq r\leq r_1}[r\nab\Bb] &\les& \EF^{J-1}_\de[r\DDov\c\Bb] + \EF^{J-1}_\de[r\nabc_4\Bb]+\Big(\EF^{J-1}_\de[\Bb]\Big)^{\frac{1}{2}}\Big(\EF^{J}_\de[\Bb]\Big)^{\frac{1}{2}}.
\eeaa
Together with \eqref{eq:BEFdernab4BbBEFrDDovBbBrnabBb}, \eqref{eq:EstimateBbt:3}, \eqref{estimate:nondegenerate-Bb2}, \eqref{eq:EFJnormofBbintheregionrgeqr1} and  \eqref{induction:hypoth-Pc}, we infer
\begin{align}\lab{eq:EFJnormofBbintheregionrgeqrplus}
\EF^{J}_{\de, r\geq r_+}[\Bb] \les& \de_{J+1}[\Pc] +\Sk_J \Sk_{J+1} +\Rk_{J} \Rk_{J+1} +\Big(\EF^J_\de[\Bb]\Big)^{\frac{1}{2}}\Big(\Rk_{J} \Rk_{J+1}\Big)^{\frac{1}{2}}+\ep_0^2\\
  &+\Sk_{J+1}\left(\de_{J+1}[\Pc] +\Sk_J \Sk_{J+1} +\Rk_{J} \Rk_{J+1} +\ep_0^2+\Big(\EF^J_\de[\Bb]\Big)^{\frac{1}{2}}\Big(\Rk_{J} \Rk_{J+1}\Big)^{\frac{1}{2}}\right)^{\frac{1}{2}}. \nn
\end{align}

\noindent{\bf Step 7.} In view of \eqref{eq:EFJnormofr2Bintheregionrgeqrplus} and \eqref{eq:EFJnormofBbintheregionrgeqrplus}, we still need to recover $\EF^J_{\de, r\leq r_+}[r^2B]$ and $\EF^J_{\de, \leq r_+}[\Bb]$.  We have in view of the first estimate in Lemma \ref{lemma:HodgeestimatesforthecontrolofBBbandAAb:regionrleqrplus} 
\beaa
\EF_{r\leq r_+}^{J-1}[r^3\nab B] \les \EF^{J-1}_\de[ r^3 \DDov\c B]+\EF^{J-1}_\de[ r^2 \nabc_3B]+\EF^{J-1}_\de[r^2B]+\sup_{\tau\leq\tau_*}\int_{\pr^+\MM_{r\leq r_+}(1,\tau)}|\nab_4^{J+1}B|^2,
\eeaa
which together with \eqref{eq:BEFder2nab3BBEFr3DDovBBr3nabB} yields, using also \eqref{induction:hypoth-Pc},
\beaa
\EF_{r\leq r_+}^{J}[r^2B] &\les& \de_{J+1}[\Pc] +\Sk_J \Sk_{J+1} +\Rk_{J} \Rk_{J+1} +\Big(\BEF^J_\de[r^2B]\Big)^{\frac{1}{2}}\Big(\Rk_{J} \Rk_{J+1}\Big)^{\frac{1}{2}}+\ep_0^2\\
&&+ \sup_{\tau\leq\tau_*}\int_{\pr^+\MM_{r\leq r_+}(1,\tau)}|\nab_4^{J+1}B|^2.
\eeaa
Now, we have
\beaa
\int_{\pr^+\MM_{r\leq r_+}(1,\tau)}|\nab_4^{J+1}B|^2\les \int_{\pr^+\MM_{r\leq r_+}(1,\tau)}|\Bt|^2+\int_{\pr^+\MM_{r\leq r_+}(1,\tau)}|\widetilde{\Bt}|^2+\EF^{J-1}_\de[ r^2 \nabc_3B],
\eeaa
where $\Bt$ is given by \eqref{eq:definitionofBtPtpPtmandBbttorecovermissingderivativeforBandBb}, where $\widetilde{\Bt}$ is given by
 \beaa
    \widetilde{\Bt}:=(\ov{q}\nabc_4 )^{\le J-2}  \nabc^2_{\widetilde{R}}\Lieb_\Z B,
    \eeaa
    and we we used the fact that $e_4$ is in the span of $(\T, \Z)$. Thus, in view of \eqref{eq:lemma-EstimatesBtPtp3} and its analog\footnote{Note that $\widetilde{\Bt}$ corresponds to replacing $\Lieb_\T$ by $\Lieb_\Z$ in the definition of $\Bt$. One easily checks that $\widetilde{\Bt}$ can be estimated exactly as $\Bt$ and hence satisfies the analog of \eqref{eq:lemma-EstimatesBtPtp3}.} for $\widetilde{\Bt}$, we obtain, using also \eqref{eq:BEFder2nab3BBEFr3DDovBBr3nabB} and \eqref{eq:intermediaryalmostdoneesitmateforBaisdhfalif:chap15:stillmissingsomeEFnorms},  
\beaa
\EF_{r\leq r_+}^{J}[r^2B] &\les&  \de_{J+1}[\Pc] +\Sk_J \Sk_{J+1} +\Rk_{J} \Rk_{J+1} +\Big(\EF^J[r^2B]\Big)^{\frac{1}{2}}\Big(\Rk_{J} \Rk_{J+1}\Big)^{\frac{1}{2}}+\ep_0^2\nn\\
\nn  &+&\Sk_{J+1}\left(\de_{J+1}[\Pc] +\Sk_J \Sk_{J+1} +\Rk_{J} \Rk_{J+1} +\ep_0^2+\Big(\EF^J[r^2B]\Big)^{\frac{1}{2}}\Big(\Rk_{J} \Rk_{J+1}\Big)^{\frac{1}{2}}\right)^{\frac{1}{2}}\\
   &+&\Sk_{J+1}^{\frac{4}{3}}\Rk_{J+1}^{\frac{1}{3}}\Rk_{J}^{\frac{1}{3}}+\Sk_{J+1}\sqrt{ \de_{J+1}[\Pc]}+\ep_J^2.
  \eeaa
Together with \eqref{eq:EFJnormofr2Bintheregionrgeqrplus} and \eqref{eq:intermediaryalmostdoneesitmateforBaisdhfalif:chap15:stillmissingsomeEFnorms}, we deduce
\begin{align*}
\BEF^{J}[r^2B] \les&  \de_{J+1}[\Pc] +\Sk_J \Sk_{J+1} +\Rk_{J} \Rk_{J+1} +\ep_0^2 +\Sk_{J+1}^{\frac{4}{3}}\Rk_{J+1}^{\frac{1}{3}}\Rk_{J}^{\frac{1}{3}}+\Sk_{J+1}\sqrt{ \de_{J+1}[\Pc]}+\ep_J^2\\
\nn  &+\Sk_{J+1}\left(\de_{J+1}[\Pc] +\Sk_J \Sk_{J+1} +\Rk_{J} \Rk_{J+1} +\ep_0^2\right)^{\frac{1}{2}},
\end{align*}
which is the stated estimate for $\BEF^{J}[r^2B]$ in \eqref{eq:mainestimateBBb-M8}.

Next, we have in view of the second estimate in Lemma \ref{lemma:HodgeestimatesforthecontrolofBBbandAAb:regionrleqrplus} 
\beaa
\EF_{r\leq r_+}^{J-1}[r\nab\Bb] \les \EF^{J-1}_\de[r\DDov\c\Bb]+\EF^{J-1}_\de[\nabc_4\Bb]+\EF^{J-1}_\de[\Bb]+\sup_{\tau\leq\tau_*}\int_{\pr^+\MM_{r\leq r_+}(1,\tau)}|\nab_3^{J+1}\Bb|^2,
\eeaa
which together with \eqref{eq:BEFdernab4BbBEFrDDovBbBrnabBb} and \eqref{estimate:nondegenerate-Bb2} yields, using also \eqref{induction:hypoth-Pc},
\beaa
\EF_{r\leq r_+}^{J}[\Bb] &\les&   \de_{J+1}[\Pc] +\Sk_J \Sk_{J+1} +\Rk_{J} \Rk_{J+1} +\ep_0^2 +\sqrt{ \de_{J+1}[\Pc]} \left(\B^J_\de[\Bb]+\Sk_{J+1}^2\right)^{\frac{1}{2}} \\
&&+\Big(\BEF^J_\de[\Bb]\Big)^{\frac{1}{2}}\Big(\Rk_{J} \Rk_{J+1}\Big)^{\frac{1}{2}}+\ep_J\Big(\B_\de^{J}[\Bb]\Big)^{\frac{1}{2}}\\
&& + \Sk_{J+1}\left(\de_{J+1}[\Pc] +\Sk_J \Sk_{J+1} +\Rk_{J} \Rk_{J+1} +\ep_0^2+\Big(\BEF^J_\de[\Bb]\Big)^{\frac{1}{2}}\Big(\Rk_{J} \Rk_{J+1}\Big)^{\frac{1}{2}}\right)^{\frac{1}{2}}.
    \eeaa
Combining with \eqref{eq:EFJnormofBbintheregionrgeqrplus} and \eqref{eq:intermediaryalmostdoneesitmateforBbaisdhfalif:chap15:stillmissingsomeEFnorms}, we infer
\begin{align*}
\BEF^{J}[\Bb] \les&  \de_{J+1}[\Pc] +\Sk_J \Sk_{J+1} +\Rk_{J} \Rk_{J+1} +\ep_0^2 +\Sk_{J+1}^{\frac{4}{3}}\Rk_{J+1}^{\frac{1}{3}}\Rk_{J}^{\frac{1}{3}}+\Sk_{J+1}\sqrt{ \de_{J+1}[\Pc]}+\ep_J^2\\
\nn  &+\Sk_{J+1}\left(\de_{J+1}[\Pc] +\Sk_J \Sk_{J+1} +\Rk_{J} \Rk_{J+1} +\ep_0^2\right)^{\frac{1}{2}},
\end{align*}
which is the stated estimate for $\BEF^{J}[\Bb]$ in \eqref{eq:mainestimateBBb-M8}. This concludes the proof of Proposition \ref{proposition:EstimatesBBb-interior}.

%%%%%%%%%%%%%%%%%%%%%%%%%%%%%%%%%%%%%%%%%%%%%%%%%%%%%%%%%%%%%%%%

\subsubsection{Energy-Morawetz estimates for $A$, $\Ab$ (proof of Proposition \ref{proposition:EstimatesAAb-interior})}
\lab{sec:proofof:proposition:EstimatesAAb-interior}

%%%%%%%%%%%%%%%%%%%%%%%%%%%%%%%%%%%%%%%%%%%%%%%%%%%%%%%%%%%%%%%%

The proof of Proposition \ref{proposition:EstimatesAAb-interior} is very similar to the one of Proposition \ref{proposition:EstimatesBBb-interior} and we only sketch the details. First, proceeding as in Step 1 of Section \ref{sec:proofof:proposition:EstimatesBBb-interior} and considering this time the Bianchi pairs $(A, B)$ and $(\Bb, \Ab)$, we obtain the following analogs respectively of \eqref{eq:BEFder2nab3BBEFr3DDovBBr3nabB} and \eqref{eq:BEFdernab4BbBEFrDDovBbBrnabBb}
\bea\lab{eq:BEFder2nab3BBEFr3DDovBBr3nabB:Acase}
\nn&&\BEF^{J-1}_\de[ r^2 \nabc_3A]+ \EF^{J-1}_\de[ r^3 \DDov\c A  ] + \B^{J-1}_\de[r^3\nab A]\\
&\les& \de_{J+1}[B] +\Sk_J \Sk_{J+1} +\Rk_{J} \Rk_{J+1} +\Big(\BEF^J_\de[r^2A]\Big)^{\frac{1}{2}}\Big(\Rk_{J} \Rk_{J+1}\Big)^{\frac{1}{2}}+\ep_0^2,
\eea
\bea
\lab{eq:BEFdernab4BbBEFrDDovBbBrnabBb:Abcase}
\nn&&\BEF^{J-1}_\de[r\nabc_4\Ab] + \EF^{J-1}_\de[r\DDov\c\Ab]  + \B^{J-1}_\de[r\nab\Ab]   \\
&\les& \de_{J+1}[\Pc] +\Sk_J \Sk_{J+1} +\Rk_{J} \Rk_{J+1} +\Big(\BEF^J_\de[\Ab]\Big)^{\frac{1}{2}}\Big(\Rk_{J} \Rk_{J+1}\Big)^{\frac{1}{2}} +\ep_0^2.
\eea

Next, we introduce the following qualities, see (15.6.3) in \cite{GKS22},
 \bea
 \Abt:=\nabc_3^{J-1}  \nabc_{\widetilde{R}}^2 \Ab,\qquad  \At:= (\ov{q}\nabc_4)^{J-1}   \nabc_{\widetilde{R}}^2A,
   \eea
   which satisfy the Bianchi pair structure in (15.6.4) in \cite{GKS22}. We then proceed as in Step 2 and 3 of Section \ref{sec:proofof:proposition:EstimatesBBb-interior} and obtain the following analogs of \eqref{eq:lemma-EstimatesBtPtp3} and \eqref{eq:EstimateBbt:3}
\begin{align}\lab{eq:lemma-EstimatesBtPtp3:Acase}
   &\int_{\MM } r^{1+\de} |\At |^2 +\sup_{\tau\leq\tau_*}\int_{\pr^+\MM(1,\tau)}r^{2+\de}|\At |^2\nn\\
    \les&  \de_{J+1}[B] +\Sk_J \Sk_{J+1} +\Rk_{J} \Rk_{J+1} +\ep_0^2+\sqrt{ \de_{J+1}[B]} \left(\B^J_\de[r^2A]+\Sk_{J+1}^2\right)^{\frac{1}{2}}+\ep_J\Big(\B^J_\de[ r^2A]\Big)^{\frac{1}{2}}\nn\\
   &+\Sk_{J+1}\left(\de_{J+1}[B] +\Sk_J \Sk_{J+1} +\Rk_{J} \Rk_{J+1} +\ep_0^2+\Big(\BEF^J_\de[r^2A]\Big)^{\frac{1}{2}}\Big(\Rk_{J} \Rk_{J+1}\Big)^{\frac{1}{2}}\right)^{\frac{1}{2}},
  \end{align}  
 \begin{align}\lab{eq:EstimateBbt:3:Abcase}
 &\int_{\MM } r^{-1-\de}   | \Abt |^2 +\sup_{\tau\leq\tau_*}\int_{\pr^+\MM(1,\tau)} r^{-2-\de} |\Abt |^2\nn\\ 
 \les &  \de_{J+1}[\Bb] +\Sk_J \Sk_{J+1} +\Rk_{J} \Rk_{J+1} +\ep_0^2 +\sqrt{ \de_{J+1}[\Bb]} \left(\B^J_\de[\Bb]+\Sk_{J+1}^2\right)^{\frac{1}{2}} +\Big(\BEF^J_\de[\Ab]\Big)^{\frac{1}{2}}\Big(\Rk_{J} \Rk_{J+1}\Big)^{\frac{1}{2}}\nn\\
+&\ep_J\Big(\B_\de^{J}[\Ab]\Big)^{\frac{1}{2}} + \Sk_{J+1}\left(\de_{J+1}[\Bb] +\Sk_J \Sk_{J+1} +\Rk_{J} \Rk_{J+1} +\ep_0^2+\Big(\BEF^J_\de[\Ab]\Big)^{\frac{1}{2}}\Big(\Rk_{J} \Rk_{J+1}\Big)^{\frac{1}{2}}\right)^{\frac{1}{2}}.
\end{align}

Then, the rest of the proof is completely analogous to Steps 4 to 7 in Section \ref{sec:proofof:proposition:EstimatesBBb-interior}, with the respective correspondance $(\Pc, B, \Bt)\to (B, A, \At)$ and $(\Pc, \Bb, \Bbt)\to (\Bb, \Ab, \Abt)$, using also the main equations in Section 15.6 of \cite{GKS22}, details are left to the reader. This concludes the proof of Proposition \ref{proposition:EstimatesAAb-interior}.

\end{document}